\documentclass[10pt, oneside]{book}

\usepackage{amssymb,amsmath,amsthm,bbm}
\usepackage{verbatim,float,url,dsfont}
\usepackage{graphicx,subcaption,psfrag}
\usepackage{algorithm, algpseudocode}
\usepackage{algorithmicx}
\usepackage{bibentry}
\algtext*{EndWhile}
\algtext*{EndIf}
\algtext*{EndFor}
\algtext*{EndProcedure}
\usepackage{mathtools,enumitem}
\usepackage{multirow}
\usepackage{ragged2e}
\usepackage{xr-hyper}
\usepackage{array}
\usepackage{stmaryrd,scalerel}
\usepackage[dvipsnames]{xcolor}
\usepackage{tikz}
\usepackage[utf8]{inputenc} 

\usepackage{longtable}

\usepackage[colorlinks=true,citecolor=blue,urlcolor=blue,linkcolor=blue]{hyperref}
\usepackage[margin=1in]{geometry}
 \usepackage{natbib}

\usepackage[T1]{fontenc}    
\usepackage{booktabs}       
\usepackage{nicefrac}         
\usepackage{microtype}      
\usepackage{fancyhdr}

\ifdefined\TimesFont 
\usepackage{times} 
\fi

\ifdefined\ParSkip 
\usepackage{parskip} 
\fi

 \newtheorem{theorem}{Theorem}
\numberwithin{theorem}{chapter}
 \newtheorem{lemma}[theorem]{Lemma}
\newtheorem{corollary}[theorem]{Corollary}
\newtheorem{proposition}[theorem]{Proposition}
\newtheorem{definition}[theorem]{Definition}

\theoremstyle{remark}
\newtheorem{remark}[theorem]{Remark}

\newtheorem{example}[theorem]{Example}
\newtheorem{fact}[theorem]{Fact}

\renewcommand{\epsilon}{\varepsilon}

\usepackage{mathtools}

\usepackage{mathtools}
\mathtoolsset{showonlyrefs}

\usepackage{afterpage}

\makeatletter
\def\blfootnote{\xdef\@thefnmark{}\@footnotetext}
\makeatother

\usepackage{mdframed}
\usepackage{tcolorbox}
\tcbuselibrary{skins}
\tcbuselibrary{breakable}

\definecolor{theoremcolor}{HTML}{EBB900}

\newtcolorbox{proofframe}[2][]{
    enhanced,
    breakable=true,
    sharp corners=all,
    top=0mm, 
    left=4mm,
    enlarge top by=2em,
    enlarge bottom by=0.5em,
    left skip=1.5em,
    right skip=1.5em,
    borderline={0.2em}{0pt}{gray!20}, 
    boxrule=0pt,
    colback=gray!5,
    coltitle=black,
    colbacktitle=gray!5,
    colframe=gray!5,
    parbox=false,
    fonttitle=\bfseries,
    before upper={\vskip -0cm},
    title={{\begin{tikzpicture}[baseline=(title.base)]
            \node[rectangle, fill=gray!20, sharp corners, inner sep=1.5mm] (title) {#2};
          \end{tikzpicture}}
          },
    toptitle=-0.4cm,
}

\renewenvironment{proof}[1][]{%
    \begin{proofframe}{#1}\relax%
    }{\end{proofframe}}

\renewenvironment{fact}[1][]{%
   \refstepcounter{theorem}
    \ifstrempty{#1}%
    {\mdfsetup{%
    frametitle={%
        \tikz
        \node[anchor=east,rectangle,fill=theoremcolor!70]{Fact~\thetheorem};}} 
    }%
    {\mdfsetup{%
    frametitle={%
        \tikz
        \node[anchor=east,rectangle,fill=theoremcolor!70] {Fact~\thetheorem:~#1};}}%
    }%
    \mdfsetup{innertopmargin=0.2em,linecolor=theoremcolor!70, linewidth=0.2em, nobreak=true,
    topline=true, frametitleaboveskip=-0.9em,}
    \begin{mdframed}[]\relax%
    }{\end{mdframed}
}

\renewenvironment{theorem}[1][]{%
   \refstepcounter{theorem}
    \ifstrempty{#1}%
    {\mdfsetup{%
    frametitle={%
        \tikz
        \node[anchor=east,rectangle,fill=theoremcolor!70]{Theorem~\thetheorem};}} 
    }%
    {\mdfsetup{%
    frametitle={%
        \tikz
        \node[anchor=east,rectangle,fill=theoremcolor!70] {Theorem~\thetheorem:~#1};}}%
    }%
    \mdfsetup{innertopmargin=0.2em,linecolor=theoremcolor!70, linewidth=0.2em, nobreak=true,
    topline=true, frametitleaboveskip=-0.9em,}
    \begin{mdframed}[]\relax%
    }{\end{mdframed}
}

\renewenvironment{lemma}[1][]{%
    \refstepcounter{theorem}
    \ifstrempty{#1}%
    {\mdfsetup{%
    frametitle={%
        \tikz
        \node[anchor=east,rectangle,fill=theoremcolor!70]{Lemma~\thetheorem};}} 
    }%
    {\mdfsetup{%
    frametitle={%
        \tikz
        \node[anchor=east,rectangle,fill=theoremcolor!70] {Lemma~\thetheorem:~#1};}}%
    }%
    \mdfsetup{innertopmargin=0.2em,linecolor=theoremcolor!70, linewidth=0.2em, nobreak=true,
    topline=true, frametitleaboveskip=-0.9em,}
    \begin{mdframed}[]\relax%
    }{\end{mdframed}
}

\renewenvironment{example}[1][]{%
    \refstepcounter{theorem}
    \ifstrempty{#1}%
    {\mdfsetup{%
    frametitle={%
        \tikz
        \node[anchor=east,rectangle,fill=theoremcolor!70]{Example~\thetheorem};}} 
    }%
    {\mdfsetup{%
    frametitle={%
        \tikz
        \node[anchor=east,rectangle,fill=theoremcolor!70] {Example~\thetheorem:~#1};}}%
    }%
    \mdfsetup{innertopmargin=0.2em,linecolor=theoremcolor!70, linewidth=0.2em, nobreak=true,
    topline=true, frametitleaboveskip=-0.9em,}
    \begin{mdframed}[]\relax%
    }{\end{mdframed}
}

\renewenvironment{definition}[1][]{%
    \refstepcounter{theorem}
    \ifstrempty{#1}%
    {\mdfsetup{%
    frametitle={%
        \tikz
        \node[anchor=east,rectangle,fill=theoremcolor!70]{Definition~\thetheorem};}} 
    }%
    {\mdfsetup{%
    frametitle={%
        \tikz
        \node[anchor=east,rectangle,fill=theoremcolor!70] {Definition~\thetheorem:~#1};}}%
    }%
    \mdfsetup{innertopmargin=0.2em,linecolor=theoremcolor!70, linewidth=0.2em, nobreak=true,
    topline=true, frametitleaboveskip=-0.9em,}
    \begin{mdframed}[]\relax%
    }{\end{mdframed}
}

\renewenvironment{corollary}[1][]{%
    \refstepcounter{theorem}
\ifstrempty{#1}%
    {\mdfsetup{%
    frametitle={%
        \tikz
        \node[anchor=east,rectangle,fill=theoremcolor!70]{Corollary~\thetheorem};}} 
    }%
    {\mdfsetup{%
    frametitle={%
        \tikz
        \node[anchor=east,rectangle,fill=theoremcolor!70] {Corollary~\thetheorem:~#1};}}%
    }%
    \mdfsetup{innertopmargin=0.2em,linecolor=theoremcolor!70, linewidth=0.2em, nobreak=true,
    topline=true, frametitleaboveskip=-0.9em,}
    \begin{mdframed}[]\relax%
    }{\end{mdframed}
}

\renewenvironment{proposition}[1][]{%
    \refstepcounter{theorem}
    \ifstrempty{#1}%
    {\mdfsetup{%
    frametitle={%
        \tikz
        \node[anchor=east,rectangle,fill=theoremcolor!70]{Proposition~\thetheorem};}} 
    }%
    {\mdfsetup{%
    frametitle={%
        \tikz
        \node[anchor=east,rectangle,fill=theoremcolor!70] {Proposition~\thetheorem:~#1};}}%
    }%
    \mdfsetup{innertopmargin=0.2em,linecolor=theoremcolor!70, linewidth=0.2em, nobreak=true,
    topline=true, frametitleaboveskip=-0.9em,}
    \begin{mdframed}[]\relax%
    }{\end{mdframed}
}

\newenvironment{axiom}[1][]{%
    \refstepcounter{theorem}
    \ifstrempty{#1}%
    {\mdfsetup{%
    frametitle={%
        \tikz
        \node[anchor=east,rectangle,fill=theoremcolor!70]{Axiom~\thetheorem};}} 
    }%
    {\mdfsetup{%
    frametitle={%
        \tikz
        \node[anchor=east,rectangle,fill=theoremcolor!70] {Axiom~\thetheorem:~#1};}}%
    }%
    \mdfsetup{innertopmargin=0.2em,linecolor=theoremcolor!70, linewidth=0.2em, nobreak=true,
    topline=true, frametitleaboveskip=-0.9em,}
    \begin{mdframed}[]\relax%
    }{\end{mdframed}
}

\newcommand{\e}{\mathtt{e}}
\newcommand{\p}{\mathbb{P}}
\newcommand{\q}{\mathbb{Q}}
\newcommand{\s}{\mathbb{S}}
\renewcommand{\t}{\mathbb{T}}
\newcommand{\T}{\t_+}
\newcommand{\cS}{\mathcal{S}}
\newcommand{\leb}{\mathbb{L}}
\newcommand{\cP}{\mathcal{P}}
\newcommand{\cF}{\mathcal{F}}
\newcommand{\cD}{\mathcal{D}}
\newcommand{\cE}{\mathcal{E}}
\newcommand{\fE}{\mathfrak{E}}
\newcommand{\fU}{\mathfrak{U}}

\newcommand{\cU}{\mathcal{U}}
\newcommand{\cG}{\mathcal{G}}
\newcommand{\cK}{\mathcal{K}}
\newcommand{\cQ}{\mathcal{Q}}
\newcommand{\cM}{\mathcal{M}} 
\newcommand{\cN}{\mathcal{N}} 
\newcommand{\cX}{\mathcal{X}} 
\newcommand{\cL}{\mathcal{L}}
\newcommand{\fM}{\mathbb{M}}
\newcommand{\cZ}{\mathcal{Z}}
\newcommand{\cT}{\mathcal{T}}
\newcommand{\kl}{\mathrm{KL}}
\newcommand{\conv}{\mathrm{Conv}}
\newcommand{\epo}{\Lambda}
\newcommand{\bmu}{\boldsymbol{\mu}}

\renewcommand{\emptyset}{\varnothing}

\renewcommand{\d}{\mathrm{d}}  
\newcommand{\E}{\mathbb{E}}    
\newcommand{\R}{\mathbb{R}}  
\newcommand{\cR}{\mathcal{R}}   
\newcommand{\N}{\mathbb{N}}    
\usepackage{dsfont}
\newcommand{\id}{\mathds{1}} 
\newcommand{\lawis}{\buildrel \mathrm{d} \over \sim}
\newcommand{\laweq}{\buildrel \mathrm{d} \over =}
\newcommand{\var}{\mathrm{var}}
\newcommand{\mean}{\mathrm{mean}}
\newcommand{\cov}{\mathrm{cov}}

\DeclareMathOperator*{\argmax}{arg\,max}
\DeclareMathOperator*{\argmin}{arg\,min}

\newcommand{\FCR}{\mathrm{FCR}}
\newcommand{\FCP}{\mathrm{FCP}}
\newcommand{\FDR}{\mathrm{FDR}}
\newcommand{\FDP}{\mathrm{FDP}}

\usepackage{stmaryrd}
\SetSymbolFont{stmry}{bold}{U}{stmry}{m}{n}

\usepackage{booktabs,dcolumn}
\usepackage{diagbox}
\usepackage{pdflscape}
\usepackage{pgfplots}
\usepackage{multirow}  
\usepgfplotslibrary{fillbetween}
\pgfplotsset{compat=1.18}

\newcommand{\VaR}{\mathrm{VaR}} 
\newcommand{\ES}{\mathrm{ES}}

\newcommand{\simiid}{\,{\buildrel \text{iid} \over \sim\,}}

 \usetikzlibrary{arrows.meta,positioning}

\renewcommand{\thefigure}{\thetheorem}
\renewcommand{\thetable}{\thetheorem}

\makeatletter
\renewcommand{\fnum@figure}{Figure \thefigure}
\renewcommand{\fnum@table}{Table \thetable}
\makeatother

\AtBeginEnvironment{figure}{\refstepcounter{theorem}}
\AtBeginEnvironment{table}{\refstepcounter{theorem}}

\renewcommand{\thealgorithm}{\thetheorem}
\makeatletter
\renewcommand{\fnum@algorithm}{Algorithm \thealgorithm}
\makeatother

\AtBeginEnvironment{algorithm}{\refstepcounter{theorem}}

\usepackage{makeidx}
\makeindex

\fancypagestyle{plain}{%
    \fancyhf{} 
    \fancyfoot[C]{\textcolor{red}{
    } \\ \thepage} 
}
\title{Hypothesis Testing with E-values}
\author{Aaditya Ramdas\thanks{Department of Statistics and Data Science, and Machine Learning Department, Carnegie Mellon University} and Ruodu Wang\thanks{Department of Statistics and Actuarial Science, University of Waterloo}\\
\smallskip
\texttt{aramdas@cmu.edu, wang@uwaterloo.ca}
}
\date{\today}
\begin{document}

\frontmatter
\maketitle
\setcounter{tocdepth}{1}
\tableofcontents

\newpage

\vspace*{0.20\textheight}
\begin{center}
 \Large To Leila and Fangda \smallskip \\ for their unwavering support and patience\medskip 
    \\ and to Zoya and Xuelu \smallskip \\ for the joy and wonder they brought into our lives
\end{center}

\vspace*{0.54\textheight}
\begin{center}
The book is published as the inaugural issue of \smallskip \\ Foundations and Trends\textsuperscript{\textregistered}~in Statistics, Vol. 1: No. 1-2, pp 1-390 (2025)\smallskip \\
\texttt{http://dx.doi.org/10.1561/3600000002}
\end{center}

\chapter*{Preface}
\label{chap:preface}
\addcontentsline{toc}{chapter}{\nameref{chap:preface}}
\markboth{\sffamily\slshape Preface}
{\sffamily\slshape Preface}

 An e-value is a nonnegative test statistic whose expected value is at most one under the null hypothesis. 
This book is written to offer a humble, but unified, treatment of e-values for education and research in statistics and its applications. 
In our opinion, the need for such a book at this time can be explained by at least 
four reasons: 
(a) e-values have been named, utilized, and studied as a stand-alone concept only in the last few years, and  a large body of its potential users do not know what they are; 
(b) e-values are fundamental objects at the core of hypothesis testing and estimation, and they are both under-studied and under-utilized;
(c) several application domains in the natural and social sciences would benefit from knowing and adopting methodologies based on e-values in certain contexts to improve statistical efficiency and scientific reproducibility; 
(d) there has been an explosion of exciting research over the past few years,  and we think the time is ripe to collate resources in a self-contained and concise manner. 
We expand on these reasons below.

(a) As we discuss in more detail in the introduction, e-values (and their sequential extensions: nonnegative supermartingales and e-processes) have been around, either explicitly or implicitly, in the statistics and probability literature for the better part of a century. However, they have not been studied under a unified umbrella or terminology until recent years --- indeed, even the terms e-value and e-process are only roughly five years old. As a result, the e-value is not a `household concept' and statisticians often do not know much beyond its definition (if that).

(b) E-values are a fundamental concept/tool in hypothesis testing. We make five points to justify this, which will be clear from the content of the book. First, we will show that  a \emph{powered} test exists if and only if a \emph{powered} e-value exists. 
Second, when moving beyond binary accept/reject decisions as encountered in testing, if we have more than two actions, or the losses are data-dependent, then all decisions that guarantee risk control must be based on thresholding e-values.
Third, for many simple-to-state and fundamental composite testing problems (for example, testing if data come from a mixture of Gaussians), the only test we are aware of proceeds by constructing an e-value; this technique is called \emph{universal inference}. Further, there is a unique and well-defined notion of a \emph{log-optimal} e-value --- called the \emph{numeraire} --- even when testing any arbitrary composite null hypothesis, with no restriction.  
Fourth, when it comes to anytime-valid inference (in a sequential testing context) and 
combining dependent test statistics (in a multiple testing context), methods based on e-values and e-processes are the \emph{only admissible ones} to maintain statistical validity. 
Fifth, treating e-values as abstract objects leads to a deeper understanding and improvement of procedures that do not concern e-values by themselves, such as false discovery rate controlling procedures and methods to combine dependent p-values. 

(c) It is well known the sequential procedures can yield improved statistical efficiency over ones that need to fix the sample size in advance of seeing the data.
However, the recent crisis of scientific reproducibility is largely related to the use and misuse of p-values in such contexts (especially peeking at p-values and uncorrected optional stopping and continuation of experiments). This calls for methodologies that are statistically justifiable under various new and complicated environments. E-values are one tool (though certainly not the only one) to address this challenge, because they benefit from their simple definition, natural connections to game-theoretic probability and statistics, flexibility and robustness in multiple testing under dependence, and their central role in anytime-valid statistical inference. Equipping applied statisticians with the knowledge of e-values has visible benefits to the sciences and for information technology companies, and for that purpose an accessible technical book is useful. 

(d) Within the first half of this decade (2020s), research on e-values has been published at all the major venues in statistics, machine learning, information theory, and related fields --- the Annals of Statistics, the Journal of the American Statistical Association, the Journal of the Royal Statistical Society, Biometrika, Proceedings of the National Academy of Sciences, IEEE Transactions on Information Theory, Bernoulli, Statistical Science, Management Science, Operations Research, Journal of Machine Learning Research, Neural Information Processing Systems, International Conference on Machine Learning, the International Conference on Artificial Intelligence and Statistics, and many more. There have been 3 discussion papers on e-values in the Journal of the Royal Statistical Society within this period. 
These high-quality modern materials are scattered in different venues, sometimes with different notations and terminologies, due to the fast development. 
Thus, there can be significant value  in collating a good portion of them in one place.

\bigskip

We hope that, by putting the materials together in this book,  the concept of e-values becomes more accessible for  educational, research, and practical use. It is thus also the hope to  stimulate both young and established researchers, especially the next generation of statisticians, to contribute to this exciting and rapidly developing field, and to benefit the application domains of sequential anytime-valid statistical inference in the natural and social sciences.
From an educational perspective, the book should be teachable  to graduate students in statistics, probability, or computer science programs with good mathematical background, or strong undergraduate students who have a good knowledge of probability theory and statistics.

\bigskip

We do not anticipate this book to be the last word on the topic. In fact, we made the intentional choice to focus on hypothesis testing, and not on estimation. There is a rich and historied literature on \emph{confidence sequences} that we only briefly touch on in this book (Chapter~\ref{chap:eci}), but which is intimately connected to e-values. Further, except for Chapter~\ref{chap:eprocess}, we also do not get into \emph{sequential} hypothesis testing. This is for many reasons, part of which is that the optimality theory there is less well understood, and many fundamental questions related to filtrations remain open. We anticipate that books dedicated to these topics will be written in future years.

\subsection*{Acknowledgments}

We acknowledge our many wonderful coauthors and other scholars over the years on these topics, and our family members, especially our spouses, for their unwavering support.
Some parts of the book are based on joint papers with our coauthors, which will be acknowledged in the bibliographic notes in each chapter.

We have received  many useful comments and corrections  from our readers and colleagues.
We would in particular like to thank  Christopher Blier-Wong,
 Nick Koning, Peter Grünwald, Romain P\'erier, 
Vladimir Vovk, 
Will Fithian, and Zhenyuan Zhang for comprehensive comments on various parts of the book. 

We also acknowledge the detailed feedback by Aytijhya Saha, Ben Chugg, David Bickel, Diego Martinez Taboada, Hongjian Wang, Isaac Gibbs, Johannes Ruf, Kory Johnson, Lasse Fischer, Matteo Gasparin, Michael Wieck-Sosa, Neil Xu, Nikos Ignatiadis, Pratik Patil, Shubhada Agrawal, Tomas Gonzalez Lara, Yo Joong Choe, and  Yusuf Sale, who helped eliminate many typos and inconsistencies in earlier versions of the book. Last, we thank Anastasios Angelopoulos for the LaTeX formatting used for preprints of this book.

Finally, we thank the editors and editorial managers of Foundations and Trends\textsuperscript{\textregistered}~in Statistics and  the publisher for excellent editorial assistance. 

\bigskip

\bigskip

\noindent July 2025 
  \hfill 
Aaditya Ramdas

\noindent 
Pittsburgh and Waterloo
  \hfill   
Ruodu Wang

\chapter*{Lists of notation, conventions, and examples}
\label{chap:lists}
\addcontentsline{toc}{chapter}{\nameref{chap:lists}}

\section*{Standard Probability Notation}

\begin{table}[h!]
\centering\renewcommand{\arraystretch}{1.2}
\resizebox{\textwidth}{!}{
\begin{tabular}{ccl} 
\textbf{Symbol} & \textbf{for ...} & \textbf{Meaning} \\
$\Omega$ &&  sample space \\
$\cF$ &&  a $\sigma$-algebra over $\Omega$, or a filtration over $\Omega$\\ $\cX$ &&  the set of all  random variables on $(\Omega,\mathcal F)$  \\
$\cM_1$ &&  the set of all probability measures on $(\Omega,\cF)$\\
$\cM_+$ &&  the set of all nonnegative measures on $(\Omega,\cF)$\\
$\cM_1(\R)$ &   & the set of all Borel probability measures on $\R$  \\ 
$\p^n$ &$\p\in \cM_1$, $n\in \N$& the $n$-fold product measure of $\p$ on $\Omega\times \dots \times \Omega$ \\  
$\q \ll \p$ & $\p,\q\in \cM_+$ & $\q$ is absolutely continuous with respect to $\p$\\
$\E^{\p}[X]$ & $\p\in \cM_+$,  $X\in \cX$   &  $\int X\d \p$; it is the expectation under $\p$ if $\p\in \cM_1$\\ 
$\var^\p(X)$ & $\p\in \cM_1$,  $X\in \cX$   & variance of $X$ under $\p$  \\ 
$\sigma(X)$ &   $X\in \cX$   & the $\sigma$-algebra generated by $X$  \\ 
$\cP$, $\cQ$ && sets of probability measures on $(\Omega,\cF)$\\ 
$\conv(\cP)$ &  $\cP \subseteq \cM_1$ & convex hull of $\cP$ 
\\
$X \lawis \mu$ &   $\mu\in \mathcal M_1(\R)$, $X\in \cX$ & $X$ is distributed as $\mu$ (under some $\p$ from context) 
\\ 
$\mathrm{U}[a,b]$  & $a<b$  &  uniform distribution on $[a,b]$  \\
$\mathrm{N}(\mu,\sigma^2)$  & $\mu\in \R$, $\sigma\in \R_+$  &  normal distribution with parameter $(\mu,\sigma^2)$  \\ 
$\id_A $ &  $A\in \cF$   & binary-valued indicator function of event $A$ \\
$\delta_x$ & $x \in \R$ or $x\in \R^n$ & Dirac delta distribution at $x$\\
$\Phi$ & & standard normal cumulative distribution function
\end{tabular}
}
\label{tab:prob-notation}
\end{table}

\newpage
\section*{Standard Mathematics Notation}

\begin{table}[h!]
\centering\renewcommand{\arraystretch}{1.2}

\resizebox{1\textwidth}{!}{\begin{tabular}{ccl}
\textbf{Symbol} & \textbf{for ...} & \textbf{Meaning} \\
$\R$ &&   set of real numbers $(-\infty,\infty)$\\ 
$\R_+$ &&   set of nonnegative real numbers $[0,\infty)$\\
$\N$ &&    set of  positive integers $\{1, 2, \dots \}$ \\
$\N_0$ &&    set of  nonnegative integers $\{0,1, 2, \dots \}$ \\
$[K]$ & $K \in \N$  & $\{1,\dots,K\}$ \\
$|S|$ &   countable set $S$ &   cardinality of $S$ \\
$a \vee b$ & $a,b\in \R$ & $\max(a,b)$ \\
$a \wedge b$ &$a,b\in \R$  & $\min(a,b)$ \\
$a _+$ &$a\in \R$  & $a\vee 0$ \\
$a _-$ &$a\in \R$  & $(-a)\vee 0$ \\
$\e$ && $\exp(1)$\\
$\log$ & & natural logarithm to the base $\e$ \\
$\Delta_K$ & $K \in \N$ &  $\{(\lambda_1,\dots,\lambda_K) \in [0,1]^K: \sum_{k=1}^K \lambda_k = 1 \}$ \\
$\ell_K $ & $K\in \mathbb N$ &  $\sum_{k=1}^K k^{-1}$\\
$2^S$ & set $S$ & powerset of $S$ (the set of all subsets of $S$)
\\
$x_{(k)}$ & $x_1,\dots,x_K \in \R$, $k\in [K]$ &
  $k$-th ascending order statistic of $x_1,\dots,x_K$
\\
$x_{[k]}$ & $x_1,\dots,x_K \in \R$, $k\in [K]$ &
  $k$-th descending order statistic of $x_1,\dots,x_K$ \\
$\lceil x \rceil $ & $x \in \R$  &
smallest integer $n$ with $ n \ge x$
\\
$\lfloor x \rfloor $ & $x \in \R$  &
largest integer $n$ with  $n\le x$ \\$\fM_{K}$ &   $K\in \N$  & the arithmetic average function \\   
$\bigoplus_{k=1}^K \phi_k $ &  functions $\phi_1,\dots,\phi_K$ & the function $(x_1,\dots,x_K)\mapsto \sum_{k=1}^K \phi_k(x_k)$\\  
$\binom{K}{n}$  & $K,n\in\N, ~n\le K$ & $K$ choose $n$, equal to $(K!)/(n!(K-n)!)$
\\  
$\forall$  &  &  for all
\\  
$\exists$  &  & there exists\\  
$\mathbf 0$  &  & a vector of zeros  (dimension from context)
\\  
$\mathbf 1$  &  & a vector of ones  (dimension from context)
\end{tabular}}
\label{tab:real-notation}
\end{table}

\newpage

\section*{Notation Specific to the Book}

\begin{table}[h!]
\centering\renewcommand{\arraystretch}{1.2}

\resizebox{0.85\textwidth}{!}{
\begin{tabular}{ccl} 
\textbf{Symbol} & \textbf{for ...} & \textbf{Meaning} \\
$\fE$, $\fE (\cP)$ &  $\cP \subseteq \cM_1$ & the set of all e-variables for $\cP$ \\
$\fU$, $\fU (\cP)$ &  $\cP \subseteq \cM_1$ & the set of all p-variables for $\cP$ \\
$\cE$ & $\cE \subseteq \fE$ & a generic set of e-variables \\
$\cU$ & $\cU \subseteq \fU$ & a generic set of p-variables \\ 
$E$ & $E \in \fE$ & a generic e-variable \\
$P$ & $P \in \fU$ & a generic p-variable
\\
$\cK$ & & the set $\{1,\dots,K\}$, indexing multiple hypotheses 
\end{tabular}}
\label{tab:e-notation}
\end{table}

\section*{Abbreviations}

\begin{table}[h!]
\centering\renewcommand{\arraystretch}{1.2}

\resizebox{0.8\textwidth}{!}{
\begin{tabular}{ccl} 
\textbf{Abbreviation} &   & \textbf{Meaning} \\ 

iid & & independent and identically distributed\\
pdf& & probability density function\\
cdf& & cumulative distribution function\\
pmf& & probability mass function\\
KL  & & Kullback--Leibler (divergence), relative entropy\\
BH & & Benjamini--Hochberg (procedure)\\
CI & & confidence interval\\
CS && confidence sequence\\
BY && Benjamini--Yekutieli (procedure)\\
SAVI && sequential anytime-valid inference\\
\end{tabular} 
}
\label{tab:notation-abb}
\end{table}

\section*{Conventions}

\begin{enumerate}
    \item 
Scalar constants and functions are usually represented by small letters, such as $a$, $p$, $e$, $f$, or $x$. 
\item The constants $e$ (often a realized e-value) and $\e$ (exponential of $1$) are different, and we avoid using them in the same place.
\item
Boldface letters like $\bf x$ refer to vectors, and  its $k$-th entry is $x_k$ (not $\mathbf x_k$). All vectors are by default column vectors. 

    \item  Capital letters such as $X$, $P$ and $E$ are usually reserved for random variables and functions, but there are exceptions made for important constants (such as the number of hypotheses $K$ in multiple testing). Random vectors are denoted like $\bf X$, $\bf P$ and $\bf E$. 
    We sometimes also use $X$ to represent the entire data, which could be one-dimensional or multi-dimensional.

    \item  Terms like ``increasing'' and ``decreasing'' are in the non-strict (weak) sense, but  ``positive'' and ``negative'' are in the strict sense.

    \item  Inequalities between vectors and functions are always understood componentwise.

    \item The summation $\sum$ over the empty set is $0$, and the product $\prod$ over the empty set is $1$. 
    \item The set inclusion is written as $\subseteq$ and the strict set inclusion is written as $\subsetneq$.

    \item When $\E^\p[F(X)]$ appears, the suitable measurability of $F$ is  tacitly assumed.

    \item We sometimes omit the domain of integrals  and it is understood as the natural domain (e.g.,~$\R$ or $\Omega$), which should be clear from the context.  

    \item When we do not provide the proof of a nontrival result, we will point to its proof in the literature in the bibliographical note section of the corresponding chapter. 

    \item We generally rely on standard  conventions such as $0 \cdot x = 0$ for $x \in [-\infty,\infty]$ and $x \cdot \infty = \infty$ for $x \in (0,\infty]$, as well as $\log(0) = -\infty$ and $\log(\infty) = \infty$. When evaluating ratios of $[0,\infty]$-valued random variables we often use the following convention:
 $x/\infty = 0$, $\infty/x = \infty$ for $x \in [0,\infty)$, and $\infty/\infty = 1$.

\item As standard in probability theory, we remind the reader that conditional expectations like $\E^\p[Y|X]$ for random variables $X,Y$
and 
$\E^\p[Y|\mathcal F]$ for a $\sigma$-algebra $\mathcal F$  
are defined only in the $\p$-almost sure sense.  
\item If $\E[\cdot]$ or $\var(\cdot)$ appears without the probability specified, it is with respect to the true data-generating probability measure.

\end{enumerate}

\newpage
\section*{List of examples of e-values}

\begin{table}[h!]
\centering\renewcommand{\arraystretch}{1.2}
\resizebox{\textwidth}{!}{
\begin{tabular}{llll} 
\textbf{Context} & \textbf{Null} & \textbf{Alternative}  & \textbf{Section} \\ 
Likelihood ratio & $\p$ & $\q$ & \ref{sec:LR-e-variable}\\
Soft-rank e-value & exchangeable data & otherwise & \ref{sec:c1-ex2} \\
Normal two-sided & $\mathrm{N}(0,1)$& $\mathrm{N}(\mu,1)$, $\mu\ne0$  
& \ref{sec:c2-markov}\\ 
Monotone likelihood ratio & $p_\theta: \theta \le \theta_0$ & $p_\theta: \theta \ge \theta_1$ & \ref{sec:mlr-evariable}; \ref{sec:mlr-numeraire}  \\ 
Symmetry & symmetry around zero & otherwise & \ref{sec:c2-symmetry}; \ref{sec:symmetry}\\
Mean with nonnegative data &     $\E[X]\leq \mu$, $X\ge 0$ &   $\E[X] >\mu$ & \ref{sec:c2-bounded} \\
Mean and variance &     $\E[X]= \mu$, $\var(X)\leq \sigma^2$,   &    $\var(X)> \sigma^2$ & \ref{sec:c2-bounded}; \ref{sec:c9-ex-asym} \\
Sub-Gaussian distributions & sub-Gaussian, mean $\leq 0$ & otherwise & \ref{sec:c2-subG-mean}; \ref{sec:subG-mean} \\
Sub-Gaussian distributions & sub-Gaussian, mean $= 0$ & otherwise & \ref{sec:c2-subG-mean} \\
Likelihood ratio bounds & $\d\p/\d\p_0 \leq \gamma$ for given $\p_0$ & $\q$ & \ref{sec:c2-LRB}; \ref{sec:c5-ex}\\
Exchangeability &   exchangeable data & $\q$ & \ref{subsec:c2-numeraire-exch}; \ref{subsec:numeraire-exch}\\ 
Mixture method & $\p$  & $\cQ$&
\ref{sec:3-mix-plugin}\\ 
Plug-in method & $\p$  &$\cQ$&
\ref{sec:3-mix-plugin}\\
T-test (unknown $\sigma$) & $\mathrm{N}(0,\sigma^2)$, $\sigma$ unknown  & $\mathrm{N}(\mu,\sigma^2),~ \mu \neq 0$ &
\ref{sec:3-mix-plugin}\\ 
Changepoint (unknown location) & $\p^n$ & $\{\p^{n-k}\q^{k}\}_{k\in [n]}$ & \ref{sec:3-mix-plugin}\\
Universal inference e-value& $\cP$  & $\cQ$ & \ref{sec:4-split}\\ 
Numeraire & $\cP$  & $\q$ & \ref{sec:c5-numeraire}\\ 
Exponential family  & $p_\theta: \theta \le \theta_0$ & $p_{\theta_1}$ with $\theta_1>\theta_0$ & \ref{sec:exp-fam-numeraire} \\  
Mean with support $[0,1]$ &   $\E[X]\leq \mu$ & $\E[X]>\mu$ & \ref{sec:bounded-mean} \\  
Likelihood ratio process 
& $\p$  & $\q$ &\ref{sec:c6-LRP}\\
Universal inference e-process
& $\cP$  & $\cQ$ &\ref{sec:univ-eprocess}\\
Empirically adaptive e-process
& $\cP$  & $\cQ$ &\ref{sec:EAEP}\\
Time-mixture e-processes
& $\cP$  & $\cQ$ &\ref{sec:c6-time-mix}\\
Multi-armed bandit testing&      $\E[X_{k,j}|\mathcal F_{j-1} ]\le 1 $  $\forall j$ & $\E[X_{k,j}|\mathcal F_{j-1} ]>1$  $\exists j$ &  \ref{sec:c8-com-e}\\ 
 Separable compound e-values
& $\{\s \,:\, \s_k = \p_{k}\}, ~k\in [K]$ & $\q_k, ~k\in [K]$ & \ref{subsec:best_simple_separable}\\  Asymptotic zero-mean & $\E[X]=0$, $\var(X)\in(0,\infty)$ & otherwise & \ref{sec:c9-ex-asym}\\
 Compound t-test (unknown $\sigma_k$) & $\mathrm{N}(0,\sigma_k^2), ~k\in [K]$ & $\mathrm{N}(\mu_k,\sigma_k^2), ~k\in [K]$ & \ref{subsec:t_tests}\\   
Testing the $\beta$-quantile $(\VaR_\beta)$ &
 $\VaR_\beta(X)\le r$ &  $\VaR_\beta(X)> r$  & \ref{sec:c13-mean}\\
Testing an expected loss 
& $\E[\ell(X)]\le r$ & $\E[\ell(X)]> r$  & \ref{sec:c13-mean}\\
Testing the risk measure ES 
&$\ES_\beta(X)\le r$, $\VaR_\beta(X)=z$ & $\ES_\beta(X)> r$  & \ref{sec:c13-ES}\\
\end{tabular} 
}
\label{tab:examples}
\end{table}

\mainmatter

\part{Fundamental Concepts}
\label{part:background}

\chapter{Introduction}
\label{chap:introduction}

An e-value is a nonnegative test statistic whose expected value is at most one under the null hypothesis.

 We begin this chapter by offering high-level discussions on the concept of e-values and its comparison to p-values in Sections~\ref{sec:roles}
 and~\ref{sec:c1-versus}. 
Formal mathematical notation and conventions are introduced in Section~\ref{sec:notation}, and  formal definitions of e-values, p-values, and tests are given in Section~\ref{sec:1-def-E-P-test}. The first examples are then discussed in Sections~\ref{sec:LR-e-variable} and~\ref{sec:c1-ex2}. 
Historical context on e-values is provided in Section~\ref{sec:c1-history}, and a road map for the rest of the book is described in Section~\ref{sec:c1-roadmap}.

 Below are some simple messages, all formally treated in the book.
E-values, or e-variables as we usually refer to the underlying random variables, are a fundamental tool for hypothesis testing. 
 E-values are directly interpreted as evidence against the null: the larger the e-value, the more evidence that the null is false. 
The ``e'' in e-value could stand  for \emph{expectation} (their defining property) or for \emph{evidence} (their interpretation). 

The most straightforward way to view e-values is that they are composite generalizations of likelihood ratios. For testing a simple null hypothesis, all e-values are likelihood ratios (or are dominated by them). For testing a composite null hypothesis, building such ``generalized likelihood ratios'' (which are e-values) is always possible but not simple, and we will spend several chapters of the book discussing methodology and theory in this case.

E-values are named as such in order to evoke a comparison to p-values, which are a much more classical statistical tool within hypothesis testing. 
The naming analogy is that under the null hypothesis, 
\begin{align*}
&\mbox{e-values are \emph{expectation} constrained}\\&~~~~~~~~~~~~~~~~~~\mbox{versus}\\
&
\mbox{p-values are \emph{probability} constrained},
\end{align*}
since p-values must be uniformly distributed (or stochastically larger) under the null. 
Another analogy is that e-values can be seen as certain \emph{conditional expectations} on  data,
and p-values can be seen as certain \emph{conditional probabilities} on  data. 
One important interpretational difference is that for p-values, smaller values suggest that the null hypothesis is false, but it is the opposite for e-values. 
 
\begin{proof}[Writing conventions.]
\renewcommand{\qedsymbol}{}

Similarly to the situation of ``p-value'', which has several variants, it is possible that different authors might write ``e-value'' differently.

Our default choice of writing convention, also our recommendation to other authors, is to write ``e-value'' (as well as related terms) with the small letter ``e'', and without making it italic.
 When it is the initial letter of a title or a sentence, ``E'' (instead of writing ``e-Value'') should be  capitalized just like in a usual English word; otherwise we do not  capitalize it.
    This convention is adapted by the authors of this book and most of their co-authors.
    
    Some  journals have specific styles. For instance, some may require both ``e'' in e-values and ``p'' in p-values to be italic, so it is possible to see ``\textit{e}-value'' in  journal publications. 
\end{proof}

\section{The three roles of e-values} 
\label{sec:roles}

Our book will attempt to elaborate on three somewhat distinct roles that e-values play in statistical inference: as methodological tools with advantages in particular applications, as technical tools to prove other results of independent interest, and as fundamental tools lying at the very heart of statistics. We elaborate on these now, pointing forward to which angles are elaborated in which chapters.

\begin{enumerate}
    \item \textbf{Methodological.} E-values easily yield the more common inferential tools of p-values and hypothesis tests: indeed, the inverse of an e-values is a p-value (e-to-p calibration; Chapter~\ref{chap:markov}), and rejecting the null when the e-value exceeds $1/\alpha$ is a valid level-$\alpha$ test (Markov's inequality; Chapter~\ref{chap:markov}). In fact, sometimes e-values yield the only known path to performing statistical inference:  Chapter~\ref{chap:ui} presents a general technique called universal inference to derive e-values for (possibly) irregular testing problems, for which no other inferential methods are known.  Chapter~\ref{chap:numeraire} then shows that for a fixed alternative, universal inference can be improved to a unique ``log-optimal'' e-value, which always exists without regularity assumptions. Chapter~\ref{chap:eprocess} demonstrates how e-values enable sequential anytime-valid inference. E-values play a key role in selective inference by providing simple methods to control the false discovery rate (Chapter~\ref{chap:compound}) and false coverage rate (Chapter~\ref{chap:eci}), and they also yield nonparametric tests for risk measure forecasts (Chapter~\ref{chap:risk-measures}).  
    
    \item \textbf{Technical.} E-values appear as technical tools to solve problems that may appear to have nothing to do with e-values in their problem statement. For example, we will show that every procedure that controls the false discovery rate is actually an instance of the e-Benjamini--Hochberg procedure applied to certain \emph{compound} e-values (Chapter~\ref{chap:compound}). There are no such results available for p-values or the original Benjamini--Hochberg procedure. These (compound) e-values  may not be of direct interest in and of themselves, but appear as intermediate technical tools towards some other aim, such as combining inferences across many multiple testing procedures or  derandomizing multiple testing procedures that utilize external randomness. 
    Similarly, when combining arbitrarily dependent p-values into a single combined p-value, all admissible merging functions must employ e-values under the hood (Chapter~\ref{chap:pmerge}). 
    Methods based on e-values are intimately connected to other theories in modern mathematics, notably the 
    reverse information projection (Chapter~\ref{chap:numeraire}) and nonnegative (super)martingales  (Chapter~\ref{chap:eprocess}); once more, e-values arise centrally in both settings, whether or not one is interested in them.
    
    \item \textbf{Fundamental.} Chapter~\ref{chap:alternative} points out that for testing any composite null $\cP$ against any composite alternative $\cQ$, without restriction, there exists a powered test (a test with a prespecified level under any null, and power greater than that level against any alternative) if and only if there exists a bounded powered e-variable (a variable with expected value at most one under any null, and expected value larger than one under any alternative). Chapter~\ref{chap:posthoc} shows that for general decision-making problems with more than two decisions, or with data-dependent loss functions, any decision that guarantees risk control must be (implicitly) based on e-values. 
\end{enumerate}

 Overall, we view the current book as providing the building blocks, namely e-values, of \emph{a theory of evidence}. 
 This is contrast to, say, \emph{a theory of decision making} as developed by Neyman, Pearson, Wald and Savage. While the two are not at complete odds, they are also not in complete agreement. 

Given that one might argue that p-values also play similar roles in the literature as written above, we next address the relationship of e-values to p-values.

\section{P-values versus e-values}
\label{sec:c1-versus}

The book does not seek to displace or eliminate the use of p-values. Instead, we argue that e-values are complementary to p-values, and indeed a critical part of the modern statistician's toolkit. We also argue that e-values and p-values should be studied on an equal footing, and that they may each have advantages in different settings.

\subsection{E-values and p-values are statistical cousins}
\label{sec:cousins}

P-values have a long and rich history in statistics, long considered a standard tool in theory, methodology and applications. E-values, by contrast, have been relatively late to enter the scene, and are not currently  considered by common scientists to be as fundamental on all three aforementioned fronts. 
However, we argue that e-values and p-values should in fact be mentioned in the same breath --- they are closely related, equally fundamental, and we would like both to occupy an equal footing in the literature.

This is a tall ask, but e-values are up to the challenge. It is immediately apparent from their definitions that p-values and e-values are rather different objects: they are both test statistics of the data but with entirely different properties. Thus, we will focus here on the commonalities between e-values and p-values, in order to make the case that one should view them as statistical ``cousins'':
\begin{itemize}
    \item 
First, as classes of objects, p-variables, e-variables and tests exist under exactly the same conditions. Theorem~\ref{th:existence} proves that there exists a bounded powered e-variable if and only if there exists a powered level-$\alpha$ test for some (or every) $\alpha$. This holds without restriction on the (possibly composite) null hypothesis being tested, or on the composite alternative. Further, under a very mild assumption~\eqref{assm:U} that there exists a random variable that is uniformly distributed under every null and alternative, independent of all other random variables, the preceding existence statements are also equivalent to the claim that a powered p-variable exists.
\item
Second, as mentioned earlier, e-variables and p-variables can be transformed to each other by ``calibration'' (Section~\ref{sec:2-cali}). While e-values and p-values are in general just real-valued measures of evidence, they can be thresholded to yield rejection decisions of the null hypothesis with guaranteed type-I error control. While e-to-p calibration with $e \mapsto 1/e$ preserves the final decision (for any $\alpha$, the calibrated p-value is below $\alpha$ if and only if the original e-value is above $1/\alpha$), this is in general not true for all p-to-e calibrators. Nevertheless, the particular p-to-e calibrator $p \mapsto \id_{\{p \leq \alpha\}}/\alpha$ results in the all-or-nothing e-value that preserves the final decision at a fixed level $\alpha$. 
\end{itemize}

Despite the above similarities,
some key differences are as follows: 
\begin{itemize}
    \item \textbf{E-values allow for post-hoc decision making.} For testing, the level $\alpha$ can be chosen to depend on the data while maintaining a nontrivial guarantee (Section~\ref{sec:2-posthoc}). More generally, if one needs to make decisions beyond just accept/reject, e-values allow handling of data-dependent loss functions (Section~\ref{sec:post-hoc-decisions}). Critically, this is true for \emph{all} e-values, and while this may be true for \emph{some} p-values (like those derived from e-values), it is provably not true for all p-values.
    \item \textbf{E-values allow for optional continuation of experiments.}     
    Suppose a scientist collects some data that appear promising but inconclusive to make a clear decision of rejecting a null hypothesis.
Seeing this result, another scientist may then be inspired to run their own follow-up experiment, and their e-values can be merged by a simple multiplication. 
    On the other hand, sequentially obtained p-values are difficult to merge, because the very existence of the latter studies may depend on the p-value obtained from the first experiment, creating a subtle form of dependence that could lead to incorrect meta-analyses. 
These issues are discussed in detail in  Sections~\ref{sec:c4-optional}. This is also related to ``p-hacking'' (Section~\ref{sec:c7-optional}), for which e-processes provide a compelling solution.

\item \textbf{E-values are suitable for unknown sampling scheme.} 
The construction and validity of p-values rely crucially on a completely specified design of the experiment, in particular the data collecting procedure. If a dataset is simply handed to us without specifying how it was collected, one may not be able to correctly calculate a p-value, but we are often able to calculate an e-value, for a wide range of reasonable sampling schemes. This is illustrated in Example~\ref{ex:sampling}.

    \item \textbf{Compound e-values are central to multiple hypothesis testing.}  Chapter~\ref{chap:compound} shows that every multiple testing procedure that controls the false discovery rate can be written an instance of the e-Benjamini--Hochberg procedure with certain compound e-values, while Chapter~\ref{chap:combining} shows that every admissible way for combining dependent p-values must proceed via converting them to e-values and averaging those e-values.
    Chapter~\ref{chap:eci} illustrates how e-values are useful in constructing and combining  dependent confidence intervals. 
    
\end{itemize}

There are other differences one can point out in more specific settings. For example, when testing whether a Gaussian mean is nonpositive, there exists a uniformly most powerful p-value which does not depend on the (positive) alternative mean. This is not true for e-values, whose \emph{e-power} typically depends on the alternative, even in special situations such as the one considered. In contrast, for a fixed alternative, a log-optimal (``uniformly most e-powerful'') e-value always exists for any composite null (Chapter~\ref{chap:numeraire}), but uniformly most powerful p-values do not exist outside of very special cases.

Thus, p-values and e-values are very different test statistics with several contrasting properties; Chapter~\ref{chap:theoretic} contains more such contrasts. 
At a high level, p-values may be better choices in certain rigid and constrained settings, but e-values may be better choices in more open-ended or flexible settings.
Yet, they are strongly related as families of objects, and can be converted to each other. This justifies our calling them ``statistical cousins''.

\subsection{Power and e-power}
\label{c1:p-vs-e}

A question that arises very often, perhaps the first question of many users, is 
\begin{center}
\emph{Are e-values more or less powerful than p-values?}
\end{center}
This question requires some pause to carefully interpret and answer. E-values and p-values are summary statistics and do not have ``power'' by themselves.
Power refers to the probability of rejection under some specific alternative hypothesis.
The questioner often means \emph{if I construct a level-$\alpha$ test by rejecting an e-value when it exceeds $1/\alpha$, does it have lower power than the level-$\alpha$ test that rejects when a p-value is smaller than $\alpha$}. Or, \emph{if I embed an e-value in a decision making (e.g., multiple testing) procedure, will I make better decisions (e.g., more discoveries)?} This question does not specify which p-values and e-values one is using. So our answer can be split into four cases:
\begin{enumerate}
    \item If the p-value $P$ was itself derived from the e-value $E$ (since $1/E$ is a p-value; see Chapter~\ref{chap:markov}), then clearly the two tests are identical and there is no difference in power.
    \item If one starts with a p-value $P$, and uses it to define the e-value
    $E = \id_{\{P \leq \alpha\}}/\alpha$, then clearly the resulting test has the same power at level-$\alpha$ as that based on the original p-value. 
    \item If a given p-value $P$ is converted to an e-value $E$ using a different p-to-e calibrator than above (Chapter~\ref{chap:markov}), then the probability that $E$ exceeds $1/\alpha$ will be smaller than the probability that $P$ is below $\alpha$ (see also Section~\ref{sec:c2-jeffrey}).
    \item When the p-value and e-value are not derived from each other in the above ways, no general comparisons can be made. 
\end{enumerate}
The second point above can be generalized to an important, though pedantic, observation --- any level-$\alpha$ test (even those that may or may not be based on p-values) can be recovered by thresholding an appropriate e-variable at $1/\alpha$ (Chapter~\ref{chap:alternative}). Indeed, if $A$ is the rejection region of some level-$\alpha$ test, meaning that we reject when the data lies in $A$, then $E = \id_{A}/\alpha$ is an ``all-or-nothing'' e-value that equals or exceeds $\alpha$ exactly when the original test rejects. Thus, clearly, using e-values does not by itself lose expressibility or power in any way.

It is worth emphasizing that while we repeatedly use $1/\alpha$ as the ``default'' threshold for converting e-values to level-$\alpha$ tests, one must be somewhat more careful, keeping the following points in mind: 
\begin{itemize}
    \item \textbf{Full information, possible randomization.} 
    If the distribution of the e-value under the null hypothesis is known (and it could be different under each element of the null), then it is   possible to design a threshold smaller than $1/\alpha$ (Section~\ref{sec:dist-info}).
    If the null hypothesis is simple or
    if the e-value is \emph{pivotal}, meaning that its distribution is invariant under the composite null hypothesis, then this type-I error can be made exactly $\alpha$ (possibly using randomization, as in the classical Neyman--Pearson lemma).

    \item \textbf{No information, randomization.} It is always possible to reduce the threshold for rejection (from $1/\alpha$) by a random amount using the randomized Markov's inequality (Section~\ref{sec:2-rand}). One could equivalently view this as converting the e-value to a smaller randomized p-value, or equivalently stochastically rounding the e-value to an `all-or-nothing' e-value that has larger probability of exceeding $1/\alpha$ (power) than the original (Section~\ref{sec:c8-stoch}). 

    \item \textbf{Partial information, no randomization.} If partial information is available about the e-value, for example if the null distribution of the e-value  satisfies certain shape constraints, then one can utilize thresholds that are below $1/\alpha$ (Chapter~\ref{chap:shape}). 

    \item \textbf{Multiple e-values, possible randomization.} If multiple exchangeable e-values must be combined (which could occur when the e-values are calculated using a randomized procedure, such as sample splitting), then simply averaging them and thresholding at $1/\alpha$ is suboptimal, and one can improve on it by using an exchangeable improvement of Markov's inequality (Chapter~\ref{chap:ui}). If the e-values are arbitrarily dependent, they can be made exchangeable by applying a uniformly random permutation, and then the previous sentence applies.
\end{itemize}

The above points may not  alleviate all concerns. Indeed, users who directly plug in specific e-values from this book into their favorite decision procedures (that are currently based on p-values) may indeed find that the e-value procedures are less powerful. Why is that? Here is the main reason.
    
    \textbf{The typical e-values in this book are optimizing a different criterion than power.}  Indeed, the preferred criterion in this book is called ``e-power'' (Chapter~\ref{chap:alternative}), and corresponds to the expected logarithm of the e-value under the alternative $\q$, that is, $\E^\q[\log(E)]$. On that note, we make the following two observations:
    \begin{itemize}
        \item 
    If we were to optimize for power at level $\alpha$, we would simply get an all-or-nothing e-value $\id_A/\alpha$ for some set $A$. These have e-power equal to minus infinity, because they take on the value zero with non-negligible probability, which is infinitely penalized by the logarithm.
    \item 
    Log-optimal e-values, in contrast, are those which maximize the e-power against $\q$ (these always exist, even for composite nulls; Chapter~\ref{chap:numeraire}). Log-optimal e-values obviously avoid the value zero: they are $\q$-almost surely strictly positive.
    \end{itemize}
    Since our ``preferred'' e-values are log-optimal, they are rather different from the ones that maximize power. Thus, directly plugging them into the reader's favorite procedures (which may be judging the results by more traditional metrics like power), may lead to disappointment. This is not a fault of e-values, but of \emph{which e-values one chooses to use}. Carefully aligning the optimality criterion of the e-value with that of the procedure within which they are used is likely to return more favorable results. Indeed, we already pointed out that you certainly cannot \emph{lose} anything by thinking in terms of e-values instead of p-values.
     
    The reader may then wonder: what makes us prefer e-power (and log-optimality) as the criterion of choice?  Indeed, one could take any increasing function $f$ and ask for the e-value $E$ that optimizes $\E^\q[f(E)]$; power at level $\alpha$ and the logarithm are only two such choices. Nevertheless, the logarithm is a \emph{natural} choice if you take into account that e-values are interpreted as measures of evidence against the null, and one that is particularly easy to update when more evidence is collected. 
    
    For example, consider the optional continuation setting in Section~\ref{sec:cousins}. 
    Suppose a sequence of scientists collect a sequence of e-values $E_1,E_2,\dots$ against the null; then $\prod_{i =1}^n E_i$ is also an e-value, even if the e-values are sequentially dependent, meaning that the decision to collect $E_i$ was made after observing $E_1,\dots,E_{i-1}$. Now it may be immediately apparent why one wants to avoid zeros: if $E_3 = 0$, no future multiplication of large e-values will ever result in large evidence. 
    
    Further, note that $\prod_{i =1}^n E_i = \exp(n \cdot  \overline{\log E}_{[n]})$, where $\overline{\log E}_{[n]} = \sum_{i =1}^n  \log E_i / n$, so that the product increases to infinity whenever $\overline{\log E}_{[n]}$ converges to a positive constant. Further, in that case, the product converges to infinity exponentially fast, and choosing the e-values that maximize $\E^\q[\overline{\log E}_{[n]}]$ (under the alternative) also maximizes the \emph{rate} at which the product approaches infinity. Here, the logarithm was not simply a ``convenient choice of a utility function'', but arises naturally in reasoning about a sequence of products. 
    We also offer an axiomatic derivation of the logarithm in Chapter~\ref{chap:alternative}.

We are now ready to begin  with the formal mathematical framework of concepts in the book.

\section{Mathematical notation and conventions}
\label{sec:notation}

We begin with a sample space $\Omega$ equipped with a $\sigma$-algebra $\cF$, and the set $\cM_1$ of all probability measures on $(\Omega, \cF)$. 
Elements of $\cM_1$ are also called distributions. 
Formally, a random variable is a measurable function from $(\Omega,\mathcal F)$ to $\R$ or $[-\infty,\infty]$. When we allow for a random variable to take infinite values, we will write it explicitly; otherwise it takes values in $\R$.   We assume that some distribution $\p^{\dagger} \in \cM_1$ governs our data $X$ (which is a random variable). 
Note that $X$ could be a vector $(X_1,\dots,X_n)$. The variables $X_1,\dots,X_n$ could be independent and identically distributed (iid) under $\p^{\dagger}$, but they need not necessarily be so. We use $X^t$ as a shorthand for the first $t$ data points $(X_1,\dots,X_t)$.

\begin{definition}
\label{def:hypothesis}
A \emph{hypothesis} is a set of probability measures in $\cM_1$. 
A hypothesis is \emph{simple} if it is a singleton, such as $\{\p\}$  and $ \{\q\}$. Otherwise it is \emph{composite}. 
 \end{definition}

We will always use $\cP$ for the null hypothesis and $\cQ$ for the alternative hypothesis. 
We will also write $H_0$ as the null hypothesis that 
$\p^{\dagger} \in \cP$,
and $H_1$ the alternative hypothesis that $\p^{\dagger} \in \cQ$; that is, we interpret a hypothesis both as a mathematical object $\cP$ and as a statistical statement $H_0$. 
For a simple hypothesis $\cP=\{\p\}$, we say ``testing $H_0$'', ``testing $\cP$'', and ``testing $\p$'' interchangeably.
The sets $\cP$ and $\cQ$ are always assumed to be non-intersecting. Let $\conv(\cP)$ denote the convex hull of $\cP$.

As alluded above, we always use $\cP,\cQ$ for sets of distributions, and  $\p,\q$ for a single distribution.
In the context of a hypothesis $\cP$, 
statements on independence, distribution, and almost-sure equalities are meant to hold under each member.  We say that $\q$ is  absolutely continuous with respect to $\p$, denoted $\q \ll \p$, if $\p(A)=0$ implies that $\q(A)=0$ for any measurable set $A$. In this case, we use ${\d\q}/{\d\p}$ to refer to the likelihood ratio (or, more generally, the Radon--Nikodym derivative) between $\q$ and $\p$. Letting $q$ and $p$ represent the densities of $\q$ and $\p$ with respect to some reference measure $\leb$ (often the Lebesgue measure, hence the symbol), we can write $({\d\q}/{\d\p})(x) = {q(x)}/{p(x)}$.
The expectation of a random variable $X$ under $\p$ is denoted by $\E^\p[X]=\int X\d \p$.

Some notation will be used throughout.  The set  $\N=\{1,2,\dots\}$ is  set of all positive integers, and   $\R=(-\infty,\infty)$ is the set of all real numbers. 
We say that a random variable $X$ is stochastically larger (under $\p$, sometimes omitted if clear from 
 the context) than a distribution on $\R$, represented by its cumulative distribution function  (cdf) $F$, if $\p(X\le x) \le F(x)$ for all $x\in \R$. 
 We write $X\lawis F$ if the distribution of $X$ is $F$ (either a cdf or a distribution) under a probability measure $\p$ that should be clear from  context; the notation $\overset{\mathrm{iid}}\sim $ is similar. For a positive integer $K$,  we denote by $[K]=\{1,\dots,K\}$ and by $\Delta_K$ the standard simplex
 in $\R^K$, that is,
 $$\Delta_K=\left\{(v_1,\dots,v_K) \in [0,1]^K: \sum_{k=1}^K v_k = 1 \right\}.$$ 
The Euler number  $\exp(1)$ is denoted by $\e$ (please note that this is different from the variable $e$, but of course we will avoid using both in the same equation whenever possible).

All terms like ``increasing'' and ``decreasing'' are in the non-strict sense. We write $x\wedge y$ and $x\vee y$ for the minimum and the maximum of real numbers $x,y$ respectively.


We will introduce some more specific notation in different chapters. Lists of commonly used notation  and convention are collected separately in the beginning of the book, for  readers navigating non-sequentially.

\section{Definitions of e-variables,  p-variables, and tests}
\label{sec:1-def-E-P-test}

As above, let $\cP \subseteq \cM_1$ be a set of probability measures representing the null hypothesis.

\begin{definition}[E-variables, p-variables, and tests]
\label{def:e-value}
\begin{enumerate}
\item [(i)]
An \emph{e-variable} $E$ for  $\cP$ is a $[0,\infty]$-valued random variable satisfying $\E^\p[E]\le1$ for all $\p\in \cP$. An e-variable
$E$ is \emph{exact} if $\E^\p[E] = 1$ for all $\p\in \cP$. 
\item [(ii)]
A \emph{p-variable} $P$ for  $\cP$ is a  $[0,\infty)$-valued random variable satisfying $\p(P\le \alpha)\le \alpha$ for all $\alpha \in (0,1)$ and all $\p\in \cP$. A p-variable  $P$ is   \emph{exact} if $\p(P\le \alpha)= \alpha$ for all $\alpha \in (0,1)$ and all $\p\in \cP$
\item [(iii)]A  \emph{test} $\phi$  is a $[0,1]$-valued random variable.
A test is \emph{binary} if its range is $\{0,1\}$.
The \emph{type-I error (rate)} of a test $\phi$ for $\p$ is $\E^\p[\phi]$.
 A test $\phi$ has \emph{level} $\alpha\in [0,1]$  for $\cP$ if its type-I error is at  most $\alpha$ for every $\p\in \cP$. 
 \end{enumerate} 
 
 We denote by
 $\fE = \fE(\cP)$  (Fraktur E) the set of all e-variables for $\cP$  and 
by $\fU = \fU(\cP)$ (Fraktur U) the set of all p-variables for $\cP$,
with $\cP$ often omitted. 
\end{definition}

Above, we use $\fU$ (hinting at ``uniform'') instead of Fraktur P to avoid overloading the letter ``P'', which is already used for p-variables ($P$),  probability measures ($\p$), and null hypotheses ($\cP$). 

The value taken by the e-variable or p-variable upon observing the data is called an \emph{e-value} or a \emph{p-value}.
As such, we can see that e/p-values  are actually not mathematically rigorously defined objects, although they are common terms in the statistical literature.
In this book, we always use ``e/p-variables'' when making a mathematical statement or result about the corresponding random variables, and 
mention ``e/p-values'' when their statistical interpretation is more important. 
In the latter case, an e/p-value may   refer to the random variable or its realized value, and this should be clear from context.

Some immediate remarks are made below, but the reader may directly jump to Section~\ref{sec:LR-e-variable} to see a simple example of an e-variable, which is understandable without reading the following remarks.

\begin{enumerate}
\item  Our definition of a test allows for external randomization: upon observing the data, the realized value of a test $\phi$ is interpreted as the probability of rejecting the null hypothesis. More explicitly, one imagines drawing a uniform random variable $U$ independent of the data, and reject the null if $U \leq \phi$. 
With this randomization mechanism, the type-I error  for $\p$, $\E^\p[\phi]$, is precisely the unconditional probability of rejecting the null hypothesis $\p$, when the data is drawn from $\p$. A binary test rejects or does not reject without using randomization. 

    \item A p-variable defined as a random variable that is stochastically larger than a uniform random variable $U$ on $[0,1]$ under all $\p\in \mathcal P$ (an equivalent formulation is $\E^\p[f(P)]\ge \E^\p[f(U)]$ for all increasing functions $f$ and  $\p\in \cP$).
A p-variable is often truncated at $1$; this usually does not affect our discussions.  
We occasionally mention that a random variable $P$ is a p-variable even if it can take the value $\infty$ (such as one over an e-variable; see Proposition~\ref{prop:markov}), but it should be understood as saying that $P\wedge 1$ is a p-variable. 
 \item A large realized e-variable is interpreted as evidence against the null hypothesis (because their expected value is no larger than $1$ under the null). Similarly, if we observe a small realized p-variable, then we have evidence against the null hypothesis. This interpretation will be kept throughout the book. 
Two results in Chapter~\ref{chap:theoretic} will show that  e-variables are naturally represented by conditional expectations, whereas
p-variables are  naturally represented by conditional probabilities. 
 \item We sometimes use the phrase  ``an e-variable for $H_0$'' to mean ``an e-variable for $\cP$'' where the hypothesis $H_0$ is identified by $\cP$ (as in Definition~\ref{def:hypothesis}).
If $\mathcal P = \{\p\}$, then we also say ``an e-variable for $\p$''. 
\item  Most tests that we will deal with are binary, with $1$ representing rejection and $0$ representing no rejection.
 Any level-$0$ test $\phi_0$ must satisfy   $\phi_0=0$ $\p$-almost surely for every $\p \in \cP$, and 
 we typically choose $\phi_0 = 0$.
As standard in hypothesis testing, 
the general aim is that we hope to not reject the null hypothesis $H_0$ when the data are generated by some element of $\cP$,
and  we hope to reject  $H_0$ when the data are  generated by some element of $\cQ$.

\item An e-variable is allowed to take the value $\infty$;
observing $E=\infty$ for an e-variable $E$
means that we are entitled to reject the null hypothesis; this corresponds to observing $0$ for a p-variable. However, $E$ must be $\p$-almost surely finite for every $\p \in \cP$ (otherwise its expectation would not be bounded by 1 under $\p$). In particular, $E$ is allowed to be infinite only on $\cP$-nullsets (which are sets of measure zero under every $\p \in \cP$), which correspond to events that are impossible under the null (but may be possible under some alternative). At least one e-variable always exists for any $\cP$: the trivial e-variable that always equals 1 ($\p$-almost surely for every $\p \in \cP$). Of course, any constant smaller than one also yields a valid e-variable, but such an e-variable would be dominated by 1.

\end{enumerate}

We summarize some  useful properties of p-variables and e-variables below, which can be  immediately checked from Definition~\ref{def:e-value}. 

\begin{fact}
\label{fact:basic}
For a random variable $X$ and a distribution $\p\in \cM_1$, let $F_{X|\p}$ be the cdf of $X$ under $\p$.
\begin{enumerate}[label=(\roman*)]
    \item For any random variable $X$ and distribution $\p$,   $  F_{X|\p} (X)$ is a p-variable for $\p$.
\item For any random variable $X$ and hypothesis $\cP$,  $\sup_{\p\in \cP} F_{X|\p} (X)$ is a p-variable for $\cP$.
    \item If $E_1,\dots,E_K$ are   e-variables for $\cP$, then   their arithmetic average is an e-variable for $\cP$. 
    \item If $E_1,\dots,E_K$ are   e-variables for $\cP$ and they are independent under each $\p\in \cP$, then  their product is an e-variable for $\cP$.
    \item If $P$ is a p-variable for $\cP$ and $P'\ge P$, then $P'$ is a p-variable for $\cP$.
    \item If $E$ is an e-variable for $\cP$ and $0\le E'\le E$, then $E'$ is an e-variable for $\cP$.
\end{enumerate}
\end{fact}

In standard statistical practice, p-values are often computed from (i)--(ii) of Fact~\ref{fact:basic}, where $X$ is a test statistic.   
The combination properties of e-variables in (iii)--(iv) of Fact~\ref{fact:basic} are not shared by p-variables. 
The combination of e-values and that of p-values 
are  formally treated in Chapters~\ref{chap:multiple} and~\ref{chap:pmerge}, respectively.

Although we defined e-variables and p-variables together, the concept and methodologies of e-values are  not dependent on those of p-values. We will formally connect the two concepts in Chapter~\ref{chap:markov}, due to the familiarity of p-values to the statistical community, while we should keep in mind that testing with e-values without relying on p-values is the main purpose of this book. 
Moreover, in Chapters~\ref{chap:compound},~\ref{chap:combining}, and~\ref{chap:eci}, we will see how   e-values are fundamental for  multiple testing, both with and without p-values.

E-variables are also useful, or arguably even central, to sequential hypothesis testing, as we will define and discuss formally later. Informally, a finite or infinite sequence $E_1,E_2,\dots$ of e-variables (adapted to an underlying filtration) is called an \emph{e-process} if for every stopping time $\tau$ (with respect to that filtration), $E_\tau$ is also an e-value. This definition, despite its simplicity, is rather nontrivial and hides some deep connections to the theory of nonnegative (super)martingales that we discuss and explore in Chapter~\ref{chap:eprocess}.

Next, we   discuss the betting interpretation of e-values, followed by some simple examples.

\section{Betting interpretation of e-values}
\label{sec:betting-interpretation}

E-variables have a clear interpretation in terms of betting. An e-variable for $\cP$ is a valid bet against $\cP$. To elaborate, imagine two rational parties called  Forecaster and  Skeptic (here, by ``rational'' we mean that they are risk-neutral, i.e., expected value maximizers).  Forecaster claims that some $\p \in \cP$ describes the distribution of a particular (currently unseen) outcome in a satisfactory way.  Skeptic does not believe it and would like to bet against  Forecaster. The process is formalized using a statistical contract:  Skeptic pays one dollar up front, with the promise of getting $E(\omega)$ back if outcome $\omega \in \Omega$ is observed (such a function $E: \Omega \to \mathbb R_+$ must be specified by  Skeptic; this is the bet). If the contract is accepted by  Forecaster, they receive a dollar, then the outcome $\omega$ is observed, and they pay  Skeptic $E(\omega)$ dollars. When would  Forecaster accept the contract? When its terms are favorable to  Forecaster. This happens when $E$ is an e-variable, because then  Forecaster believes that they will make money by partaking in the contract. Naturally, if Skeptic believes that $\q$ describes the world better than $\cP$, they will suggest $E$ only if $\E^\q[E] > 1$. Here $E$, the e-variable, is the bet, while its realization $E(\omega)$ is an e-value, also called a betting score. 

P-values also have a betting interpretation, albeit a flawed one. The setup, continuing from above, goes as follows. For every $r \in (0,1]$, Forecaster offers a contract priced at a dollar, which pays $1/r$ dollars if some statistic $P$  of the data is smaller than $r$, and nothing otherwise. When would Forecaster offer such a contract? Only when $P$ is a p-variable for $\cP$.
Equivalently, we are really saying that if $P$ is a p-variable, then $\id_{\{P \leq r\}}/r$ is an e-variable for any $r \in (0,1]$; indeed $p \mapsto \id_{\{p \leq r\}}/r$ is a p-to-e calibrator (Section~\ref{sec:2-cali}). In any case, the allowable bets are indexed by $r$.
 Now, in this game, the realized p-value  equals the smallest value of $r$ at which  Skeptic would have received a profit, \emph{had they picked this particular, best, bet in advance}. 
 
 The above betting interpretations demonstrates why  e-values and reciprocals of p-values should not be compared on the same scale: an e-variable is genuinely a valid bet, but one over a p-variable is not---the latter is effectively ``double dipping''. Wouldn't we all like to have picked the best bet we could have in hindsight? P-values are effectively doing that, but e-values are not. In short, p-values are more akin to hindsight-optimal bets, potentially exaggerating evidence against the null. This fact helps to motivate the grids in Section~\ref{sec:c2-jeffrey}.

\section{A first example: likelihood ratio}
\label{sec:LR-e-variable}

Suppose that we are testing a simple null hypothesis $\p$ versus a simple alternative hypothesis $\q$, assuming $\q \ll \p$. For this setting, a natural e-variable is the likelihood ratio $E=(\d \q /\d \p)(X)$ where  we treat the likelihood ratio as a function of the observed data $X$. For any $\q \ll \p$, the likelihood ratio $\d \q/\d \p$ is always an e-variable for $\p$.

We specialize in two simple settings to help the reader to keep concrete examples in mind. In what follows, $X_1,\dots,X_n$ are real-valued iid data points from $\p$ or $\q$, and we denote by $S_n=\sum_{i=1}^n X_i$ their sum. 

The most standard example is  testing $\p=\mathrm{N}(0,1)$ against $\q=\mathrm{N}(\mu,1)$ for a known $\mu>0$.
Rigorously, $\p$ and $\q$ are the product measure on the data sample $X_1,\dots,X_n$ (which are iid normal), but we slightly abuse the notation here and later, which should be clear from the context.  
In this setting, the likelihood ratio e-variable is given by
\begin{equation}\label{eq:example-gaussian}
E_n= 
  \exp \left( \mu S_n - {n\mu^2}/2\right) = \prod_{i=1}^n \exp(\mu X_i - \mu^2/2).
\end{equation}
It is straightforward to check that $\E^\p[E_n]=1$, and hence the e-variable $E_n$ is exact (in fact, each of its building blocks $\exp(\mu X_i - \mu^2/2)$ is an exact e-variable). This e-variable is not an ad-hoc choice; it has certain optimality that will be formally studied in Chapter~\ref{chap:alternative}.
Write $\delta=\mu\sqrt n$, which measures the strength of the signal under the alternative hypothesis, 
and $Z=S_n/\sqrt n$, which is $\mathrm{N}(0,1)$-distributed under the null hypothesis. The e-variable in \eqref{eq:example-gaussian} can be rewritten as 
\begin{equation}\label{eq:example-gaussian-e3} E_n= 
  \exp \left( \delta Z -\delta^2/2\right).
  \end{equation}
 The Neyman--Pearson p-variable is 
  computed as 
 $
  P_n= \Phi(-S_n/ \sqrt{n}) =\Phi(-Z), $
  where $\Phi$ is the standard normal cdf.  
  Clearly,  the e-value    is  increasing in $S_n$ and the p-value   is decreasing  in   $S_n$, both confirming the intuition that a larger observed $S_n$ favours $\q$ over $\p$.

A simple calculation shows that $E_n$ is an e-variable for any choice of $\delta \in \R$ (not just for $\delta=\mu \sqrt n$ used earlier).
Suppose that we are now testing $\mathrm{N}(0,1)$ against $\mathrm{N}(\mu,1)$ for an unknown $\mu>0$.
We can choose a parameter $\delta>0$, which affects the power of the e-variable, as discussed in Chapters~\ref{chap:markov} and~\ref{chap:alternative}. 
We can also test the two-sided alternative hypothesis: $ \mathrm{N}(\mu,1)$ for an unknown $\mu\ne 0$. In this case, a simple e-variable is 
\begin{equation}\label{eq:example-gaussian-e2}
E_n= 
  \frac{\exp\left( \delta Z -\delta^2/2\right) + \exp\left( -\delta Z -\delta^2/2\right)} 2,
\end{equation}
for a parameter $\delta>0$ free to choose. We discuss how to pick $\delta$ (in a more general context) in Section~\ref{sec:3-mix-plugin}. 
The corresponding p-variable is given by 
\begin{equation}\label{eq:example-gaussian-p2}
  P_n= 2\Phi(-|Z|).\end{equation}

Another simple example is testing $\p=\mathrm{Bernoulli}(p)$ against $\q=\mathrm{Bernoulli}(q)$, where the parameters $p,q\in (0,1)$ are known. It is without loss of generality to only consider the case $q>p$.  In this setting, a natural e-variable for a sample of size $n\in \N$ is given by the likelihood ratio
\begin{equation}\label{eq:example-bernoulli}
E_n=\left( \frac{q}{p}\right)^{S_n} \left( \frac{1-q}{1-p}\right)^{n-S_n}  .
\end{equation}
 The Neyman--Pearson p-variable is 
  computed as $P_n=1-F_{n,p}(S_n)$
  where $F_{n,p}$ is the cdf of the binomial distribution with parameters $(n,p)$.   
 Again, the e-variable $E_n$ is exact, and noting that $q>p$, the e-value   is   increasing in $S_n$ and the p-value is decreasing in  $S_n$. 
The above e-variables $(E_n)_{n\in \N} $ in fact also form {e-processes},  introduced in Chapter~\ref{chap:eprocess}.

The point that e-values can be used in unknown sampling scheme, mentioned in Section~\ref{sec:cousins}, is illustrated in the following example, based on  the e-value in \eqref{eq:example-bernoulli}.
\begin{example}[Unknown sampling scheme]
    \label{ex:sampling}

    Consider a simple setting of testing $\p=\mathrm{Bernoulli}(1/2)$ against $\q=\mathrm{Bernoulli}(q)$,
     where $q>1/2$ with $n$ iid data points. 
Suppose that a scientist is presented with the dataset $(1,1,0,1,1,1,1,1)$ without explaining how the data collection was designed (e.g., from older experiments).
Consider two possibilities: 
(a) it was designed such that $8$ data points are collected;
(b) it was designed such that  data points are collected until $5$ heads in a roll are observed.
In case (a), the standard p-value is computed by $\p(N_8\ge 7)=0.035$, considered as significant in some areas of science. 
In case (b), the standard p-value is $\p(\mbox{observing 5 heads in a roll at or before $n=8$})=0.078$, often considered  as insignificant. 
The subtle point here is that if only the data are presented to the scientist but not the design (the data collector may not even have a design in mind when collecting data, and decided to stop without a plan), there is no way for the scientist to distinguish between the cases (a), (b), and other possible cases. This is a possible problematic situation in scientific practice. 
On the other hand, $E_n=q^{S_n} (1-q)^{n-S_n}/2^n$
in \eqref{eq:example-bernoulli} is a valid e-variable, regardless of how the data collecting procedure is designed as long as it collects all obtained data from the beginning. This is due to the fact that $(E_n)_{n\ge 1}$ is an e-process (formally treated in Chapter~\ref{chap:eprocess}). 

\end{example}

One may keep in mind that not all sampling schemes are allowed for e-processes either. Some obviously problematic situations include that the data were cherry-picked so that only $1$s are reported, the data collector removed the first few $0$s,
or the data collector reported the ``best'' batch of data after several trials. In these, arguably malicious, situations, any statistical methodology purely based on the collected data fails to work, including those based on e-values.

\subsection*{An illustration}

Below we illustrate in a simple case how the e-values from likelihood ratios look like and how they vary with $n$. 
Suppose that a scientist is interested in a parameter $\theta^\dagger\in\Theta$,
and iid observations $X_1,\dots,X_n$ from a distribution parameterized by $\theta^\dagger$
are available and sequentially revealed to her one by one.
She tests $H_0: \theta^\dagger  = 0$ against $H_1: \theta^\dagger \in \Theta_1$
where $0\notin\Theta_1\subseteq\Theta$. 
Let $\ell$ be the likelihood ratio function given by
\begin{equation}\label{eq:LR}
  \ell(x;\theta) = \frac{\d \q_\theta}{\d \p}(x),
\end{equation}
where $\q_{\theta}$ is the probability measure for one observation
that corresponds to $\theta\in\Theta$, and $\p$ corresponds to $\theta=0$. 
The fundamental observation, explained above, is that $\ell (X_i;\theta)$
for any $\theta\in\Theta$ and $i\in [n]$ is an e-variable for $H_0$.  
Since $H_1$ is a composite hypothesis, there is some freedom in choosing the parameter $\theta$. 
Consider the following   strategies:
\begin{enumerate} 
\item[(a)]
  Fix $\theta_1,\dots,\theta_n\in \Theta_1$,
  and define the e-variables $E_i:= \ell (X_i;\theta_i)$ for $i\in [n]$.
  One may simply choose all $\theta_i$ to be the same.
This includes the e-variable in \eqref{eq:example-gaussian}.
\item[(b)]
  Choose $\theta_1,\dots,\theta_n$ adaptively,
  where  for each $k$, $\theta_k$ is estimated from $(X_1,\dots,X_{k-1})$.
  This can be obtained, e.g., as a point estimate of $\theta$ through maximum likelihood estimation (MLE) or maximum a posteriori (MAP) estimate using some prior. 
\end{enumerate} 
The e-variables $E_i$, $i\in[n]$ are combined
to form a final e-variable $M_n:=\prod_{i=1}^n E_i$. 
All methods produce a valid final e-variable; indeed, it is straightforward to verify that the running product of $(E_i)_{i\in [n]}$ forms a nonnegative martingale $M=(M_i)_{i\in[n]}$ under $\p$ (called an e-process in Chapter~\ref{chap:eprocess}), because, with $M_0=1$, 
$$\E^\p[M_i|X_1,\dots,X_{i-1}] =M_{i-1} \E^\p[\ell(X_i;\theta_i)|X_1,\dots,X_{i-1} ] = M_{i-1}
$$ for each $i\in [n]$, as long as $\theta_i$ is computed based on $X_1,\dots,X_{i-1}$.
A further possible strategy is to take a mixture over $M$ for different strategies above (studied in Chapter~\ref{chap:alternative}).  
The example below gives a concrete illustration. 
 
\begin{figure}[t]
  \begin{center}
    \includegraphics[width=0.48\textwidth, trim={0 15 0 5}, clip]{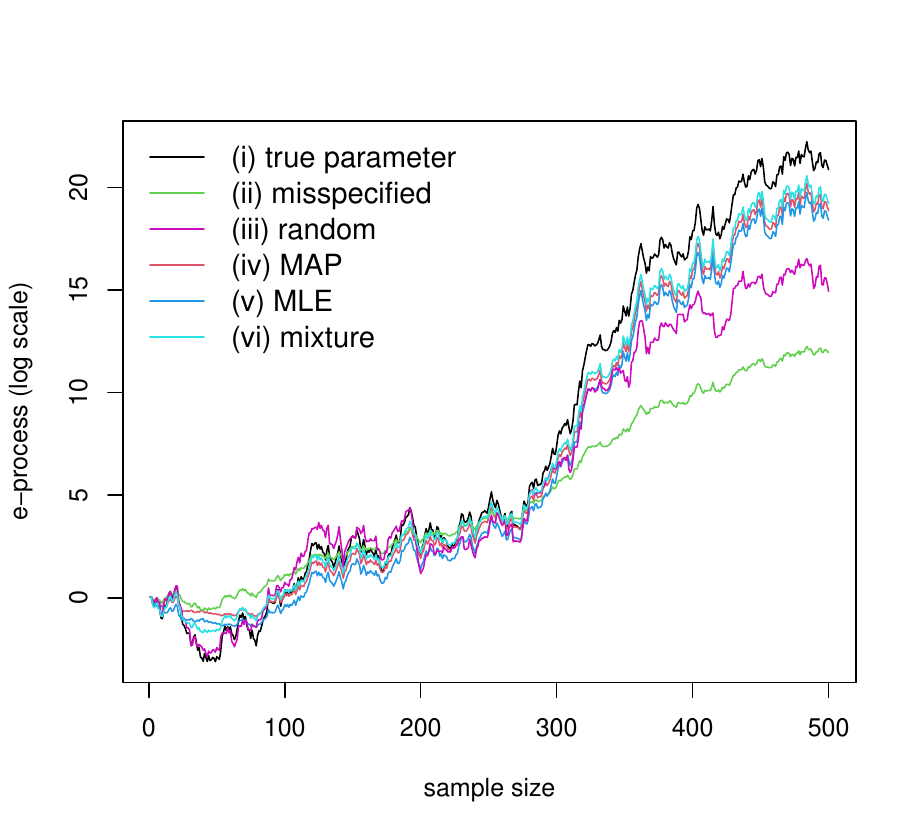}
    \includegraphics[width=0.48\textwidth, trim={0 15 0 5}, clip]{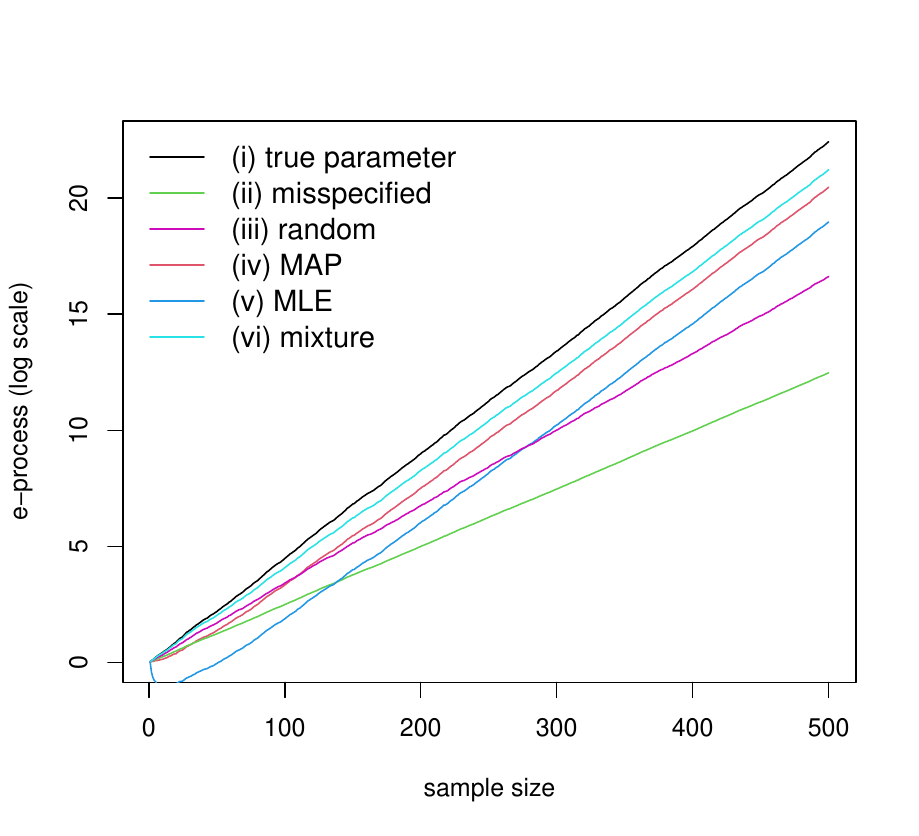}
  \end{center}
  \caption{A few ways of constructing e-values from likelihood ratio processes.
    Left: one run.
    Right: the average of 1000 runs (with average taken on the log values).}
\label{fig:c1-1}
\end{figure}

\begin{example}
\label{ex:c1-normal}
We illustrate the processes $M$ obtained from different methods for testing normal distributions.
Suppose that an iid sample  $(X_1,\dots,X_n)$ from $\mathrm{N}(\theta^\dagger,1)$ is available.
The null hypothesis is $\theta^\dagger = 0$, and the alternative is $\theta^\dagger>0$. 
We set $\theta^\dagger=0.3$, which governs the data.
We consider five different ways of constructing
$E_i := \ell(X_i;\theta_i)=\theta_i X_i-\theta_i^2/2$ for $i\in [n]$, and a further method of mixture of the processes: 
(i) choose $\theta_i := \theta^\dagger=0.3$ (knowing the true alternative),
  which grows optimally (this is called log-optimality; see Chapter~\ref{chap:alternative})
(ii) choose  $\theta_i := 0.1$, which is a misspecified alternative;
(iii) random: take iid $\theta_i$ following the uniform distribution on $[0,0.5]$;
(iv) choose $\theta_i$ by the MAP estimator with prior $\theta\sim\mathrm{N}(0.1,0.2^2)$;
(v) choose $\theta_i$ with $\theta_1:=0.1$ and $\theta_i$ the MLE of $\theta$
  based on $X_1,\dots,X_{i-1}$;
  (vi) compute a mixture of $M$ each with a fixed $\theta_i\in [0,0.6]$, uniformly weighted via a discrete grid.

The resulting $M$ curves are plotted in the logarithmic scale in Figure~\ref{fig:c1-1} with $n=500$.
For sufficiently large sample sizes, the methods (iv) and (v) based on adaptively choosing $\theta_i$ 
and the method (vi) based on the mixture 
are more powerful than (ii) and (iii) based on misspecified or random alternatives,
and they have  asymptotically the   same slope (in the log-scale)
as the best performer (i); these will be theoretically justified in Chapters~\ref{chap:alternative} and~\ref{chap:eprocess}.
The scientist may reject the null hypothesis based on a large value of $M_n$ (see Chapter~\ref{chap:markov}), but she can also reject the null hypothesis based on the maximum of $M$ or a stopping time applied to it (formally studied in Chapter~\ref{chap:eprocess}).
\end{example}

\section{A second example: the soft-rank e-value} 
\label{sec:c1-ex2}

There is a large class of classical and modern testing procedures that use some form of Monte Carlo sampling in order to produce test statistics that are exchangeable under the null, and use the rank of the original test statistic as a corresponding p-value. 
Below we give an abstract problem setup, then an explicit example---a permutation test for independence---and then describe how to build a natural e-value in this context. 

Consider testing some hypothesis $\cP$ and let $L_0$ be the test statistic calculated from the original data (a larger value of $L_0$ indicates evidence against the null). Further, suppose we can calculate  $B$ new statistics $L_1,\dots,L_B$, that are constructed to be exchangeable together with $L_0$ under the null $\cP$. Recall that random variables $L_0,\dots,L_B$ are called \emph{exchangeable} if the joint distribution of $(L_0,\dots,L_B)$ equals that of $(L_{\sigma(0)},\dots,L_{\sigma(B)})$ for any permutation $\sigma$ of $\{0,\dots, B\}$. 

For example, suppose  that our observed data consist  of
pairs $(X_i,Y_i)$, $i\in[n]$, which are assumed to be iid from some joint
distribution. Consider testing whether
$X$ is independent of $Y$ ($H_0: X\perp Y;~ H_1: X\not\perp Y$). We can reframe this as testing whether
$X_1,\dots,X_n$ are iid conditional on $Y_1,\dots,Y_n$.  Consider a simple test
statistic defined as the sample correlation between $(X_1,\dots,X_n)$ and $(Y_1,\dots,Y_n)$, denoted by $L_0$.
In order to see whether the sample correlation is sufficiently large to be statistically significant, we would compare against
the sample correlation $L_{\sigma}$ computed on the permuted data, $(X_{\sigma(1)},\dots,X_{\sigma(n)})$ and $(Y_1,\dots,Y_n)$. 
We   draw $B$ such independent random permutations $\sigma_1,\dots,\sigma_B$, and the corresponding sample correlations are denoted by $L_1,\dots,L_B$.
The key property is that under the null hypothesis, $L_0,L_{1},\dots,L_{B}$ are exchangeable.

Returning to our abstract formulation in order to build an e-value, we transform the statistics $L_0,L_1,\dots,L_B$  into some nonnegative scores $R_0,R_1,\dots,R_B$ while preserving their order and exchangeability. An explicit transformation is given in Example~\ref{ex:soft-rank-exp} below.
Define
\[
E = 
 \frac{(B+1)R_0}{\sum_{b=0}^B R_b}, 
\]
where we define $E=1$ if $\sum_{b=0}^B R_b=0$.
The fact that $E$ is an e-variable is guaranteed by $\mathbb E^\p[R_0 | \sum_{b=0}^B R_b] = \sum_{b=0}^B R_b/(B+1)$ for $\p\in \mathcal P$ due to the assumed exchangeability of $R_0,R_1,\dots,R_B$ under $\cP$. We will later see in Section~\ref{subsec:numeraire-exch} that the \emph{log-optimal} e-variable against a simple alternative (called the numeraire)  also takes the above form.

\begin{example}
    \label{ex:soft-rank-exp}
    Let $r>0$ be a prespecified constant. 
Define $L_*= \min_{b=0,\dots,B} L_b$ to be the smallest of the $(B+1)$ test statistics.
For $b=0,\dots,B$, define the transformed statistic 
\[
R_b = \frac{\exp(r L_b) - \exp(r L_*)}{r}.
\]
This transformation is performed to ensure that $R_b$ is nonnegative while the ordering amongst the test statistics is preserved.  
Note that $L_*$ can be replaced by a constant lower bound on $L_b$.
We also allow the limiting case of $r=0$, which yields $R_b=L_b-L_*\ge 0$. If  $L_1,\dots,L_B$ are assumed nonnegative, we can choose   $R_b=L_b$, i.e., using $0$ as the lower bound with $r=0$; all of these choices guarantee exchangeability and nonnegativity of $\{R_b\}_{b=0}^B$.
\end{example}


Contrasting $E$ with the usual p-variable $P$ defined by
$$
P =  \frac{\sum_{b=0}^B \id_{\{L_b  \geq L_0\}}}{B+1},$$
we find that
\[
P  =  \frac{\sum_{b=0}^B \id_{\{R_b  \geq R_0\}}}{B+1} \leq \frac{\sum_{b=0}^B R_b/R_0}{B+1} = 1/E.
\]
(If $R_0=0$ and $R_b>0$ for some $b$, then the above $1/E$ is interpreted as $\infty$.)
Since   $P$ quantifies the rank of $L_0$ amongst $L_0,\dots,L_B$, we see that $1/E$ can be seen as a smoothed notion of rank, or ``soft-rank'' for short (much like the ``soft-max'' is achieved by exponentiation in Example~\ref{ex:soft-rank-exp}). Hence, we call $E$ as the soft-rank e-value and $1/E$ as the soft-rank p-value. 
Since the direct p-value $P$ is always smaller than the soft-rank p-value $1/E$, there is apparently no advantage to using the latter for testing  a single hypothesis with the same threshold $\alpha$.  
Nevertheless, later in Chapters~\ref{chap:multiple}--\ref{chap:combining}, we will see that  e-values enjoy several advantages over p-values in the contexts of multiple testing and constructing confidence intervals.


\section{Historical context}
\label{sec:c1-history}

This section summarizes  our historical view on  hypothesis testing with e-values for the interested reader. 

\subsection{Who invented e-values?}

In some sense, the question of who invented (or, discovered?)~the concept of an e-value or e-variable is impossible to answer. As already mentioned in brief at the start of the chapter, and as we will discuss in later chapters, if $\cP = \{\p\}$ is a singleton, then all optimal e-values take the form of $\d \q/\d \p$, that is likelihood ratios of $\q$ against $\p$, for some (implicit or explicit) alternative $\q$. So, technically e-values have been around for 100 years, in the form of likelihood ratios. Similarly, for simple nulls and composite alternatives,  
Bayes factors are also likelihood ratios, and hence e-values. Likelihood ratios are central objects in parametric and nonparametric statistics, in both frequentist and Bayesian  as well as other paradigms. Therefore, in some sense, e-values have always been everywhere in modern statistics. But this is only a part of the story. 

The recent coining and usage of the term ``e-value'' is simply to recognize the importance of a more general concept that has utility much beyond point/simple nulls. Indeed, beyond the case of simple hypotheses, e-values can be viewed as nonparametric/composite generalizations of likelihood ratios to complex settings involving nonparametric and composite nulls and alternatives. Frequentists often use the generalized likelihood ratio (the ratio between the maximum likelihoods over the alternative and the null), which is not an e-value. Bayesians (or at least those who consider testing hypotheses) typically generalize Bayes factors by using the ratio of mixture likelihoods, picking a different distribution over the alternative and the null, and this is also not always an e-value. As we will later see, constructing e-values for composite nulls and alternatives is a key and central topic in the literature, and we dedicate two chapters to it in this book: Chapters~\ref{chap:ui} and~\ref{chap:numeraire}. Informally, the former chapter points out that the ratio of the mixture likelihood over alternative to the maximum likelihood over the null (a curious mixture of frequentist and Bayesian approaches) is always an e-variable. The latter chapter effectively (if only implicitly) points out that sometimes Bayes factors are indeed e-variables, if the mixture distributions are chosen in a very particular manner.

There are perhaps five central authors who we would like to centrally credit for the development of e-values, either explicitly or implicitly through their work on sequential testing, on betting, on the very related theory of nonnegative (super)martingales, or on the theory of log-optimality. These are Ville, Wald, Kelly, Robbins and Cover (in that historical order). They each had fundamental works that are intimately connected to the topics in this book, even if our motivation, treatment, and usage may substantially differ from all of them.

In its modern avatar, there have been several other central contributors to the philosophy, theory, methodology, and applications of e-values, and more generally to the subfield now called \emph{game-theoretic statistics and sequential anytime-valid inference (SAVI)}: Vovk, Shafer, Grünwald, and if we may say so, the authors of this book.   
It was in the 2018-20 period that several independent works by these authors (and their collaborators) were posted to arXiv, each with highly aligned motivation, philosophy and techniques, but using different terminology. This instigated a joint videoconferencing call in 2020 where these 5 authors jointly decided on using the term ``e-value'' whenever it is meant in a measure-theoretic context (as in this book). Of course, all of these authors have had many collaborators who have made key contributions, and they will be appropriately cited in place. 

Research on e-values has blossomed recently, often without acknowledgment of its roots. We hope that the above historical context, and the bibliographical accompaniments below are helpful for the reader. 

This particular book will focus mainly on the non-sequential context, primarily because there is already so much to say, but also because there is much more basic work to be done in the sequential context before a book can be written about it. Needless to say, such a book will be written in due course. Nevertheless, Chapter~\ref{chap:eprocess} is dedicated to the sequential setting, giving readers a glimpse of the broader framework.



\subsection{Ville, Wald, Kelly, Robbins, Cover, and Vovk}

We now discuss some of the key contributions of the aforementioned authors within the context of e-values and game-theoretic statistics (and its SAVI applications). The references are not intended to be comprehensive, but are often just exemplary of an influential line of work. We hope they serve as starting points for the reader who wishes to dig more deeply.

Ville brought martingales to the front and center of modern measure-theoretic probability theory
\citep{Ville:1939}.  He proved that for any event of measure zero, there exists a nonnegative martingale (with respect to that measure) which explodes to infinity on that event. Ville designed betting strategies that test the law of large numbers or the law of the iterated logarithm.  We will encounter some of Ville's contributions in this book. 
While our literature is certainly inspired by Ville, he was not actually motivated by statistical inference (either testing or estimation), issues around optional stopping and continuation, and so on. As the title of his thesis suggests, he was interested in a completely different goal: pointing out some basic flaws in Richard von Mises' theory of collectives. The book chapter by~\cite{shafer2022did} provides a much broader historical context for his work. 

Wald invented sequential analysis~\citep{Wald45}. Nonnegative martingales, in particular the likelihood ratio process, play a central role in Wald's methodology and theory (though he originally did not use that terminology, since the modern theory of martingales was also being developed in parallel to his work by Doob, who was aware of Ville's work). Wald emphasized his sequential test, specified in terms of a single stopping time, that optimally trades off type-I and type-II errors. In contrast, we emphasize the underlying process (the nonnegative supermartingale, or more generally, the e-process), and focus on guarantees that hold at \emph{any} stopping time, possibly not defined or anticipated in advance. 

Kelly proposed the key notion of log-optimality~\citep{Kelly56}. A few years after Shannon had developed information theory, Kelly's paper pointed out a fundamental connection between gambling and information theory. Kelly pointed out that when playing a repeated game with binary outcomes and favorable odds, it is possible to bet in such a way that your wealth grows exponentially fast, and thus it is natural to wish to optimize the rate of growth of that wealth. This immediately leads to the criterion of wanting to optimize the expected logarithm of the wealth, which we call log-optimality in this book. \cite{breiman1961optimal} made a key contribution in generalizing Kelly's work beyond the binary case, and pointing out that the same criterion also asymptotically optimizes two other related criteria: minimizing the expected time needed to reach a threshold wealth (as the threshold goes to infinity), and maximizing the expected limiting wealth  (as the time goes to infinity); the latter is often called \emph{competitive optimality}. 
While their works were not written in terms of statistical hypothesis testing, we borrow their intuition when setting up the ``testing by betting'' framework. 
There, a testing problem is transformed into a game with unfavorable odds under the null but favorable odds under the alternative, and thus a natural goal becomes to aim for log-optimality under the alternative. 

Robbins defined the fundamental concept of confidence sequences~\citep{Darling/Robbins:1967a}, the natural sequential extension of confidence intervals, and developed its duality to  ``power-one'' tests~\cite{robbins1970statistical}. He recognized and utilized that the product of (independent or sequentially dependent) e-values yields nonnegative supermartingales (despite not using the term ``e-value''), and constructed these in some nonparametric settings. Robbins' philosophical motivations were most directly aligned with ours, in that he moved away from Wald's design of a single stopping rule, to the design of anytime-valid tests and confidence sequences. However, there are still some subtle philosophical differences. Robbins focused on level-$\alpha$ concepts, whether tests or confidence sequences, while our theory emphasizes the underlying objects (e-values and e-processes), which by themselves are measures of evidence, even if they are not explicitly thresholded to Robbins' tools. 

Finally, Cover developed the framework of log-optimal portfolios for financial trading~\citep{cover1984algorithm}. Without making distributional assumptions, he developed regret bounds that compare the log-wealth of a gambler that can adaptively rebalance their portfolio in every round (reassign their current wealth to the set of available stocks) to the best constant rebalanced portfolio in hindsight (in turn which offers better returns than the best stock in hindsight). 
Cover's work does not show up as centrally as the others in this book, but this is because we only provide glimpses of the sequential setting. 

Ville, Wald, Robbins and Cover all implicitly or explicitly employed mixture or plug-in betting strategies (or likelihood ratios, or nonnegative supermartingales). The plug-in principle and the method of mixtures will be our standard go-to methods for handling composite alternatives in this work, and they can both be traced back to ideas presented by all of these authors.

Vovk spearheaded the modern development of these ideas. The core ideas underlying what is now known as \emph{game-theoretic probability} were presented in~\cite{Vovk93}. \cite{vovk1993empirical} defined and used essentially e-values in the context of algorithmic complexity and randomness (though in this context e-values had been introduced earlier by Levin and Gacs, which can then be traced back to Martin-Löf and then to Ville). \citet*{gammerman1998learning} had an explicit instance of an e-value for the composite null hypothesis of testing exchangeability in the context of conformal prediction. Following that, \cite{Shafer/Vovk:2001} wrote a book on viewing probability and finance game-theoretically. Many central ideas for game-theoretic statistics were contained in~\citet*{shafer2011test}.

\subsection{Recent progress on e-values}

As mentioned earlier, the current resurgence of interest in this topic can be attributed to a series of works that were first made public the 2018-2020 period. We discuss a dozen such works in their order of appearance (typically on arXiv) below. We also omit a dozen of our own papers that appeared during this period, as well as several papers by other authors. This has resulted in an incomplete and biased list. But we still think that this is a list that will retrospectively be viewed as containing works that developed key ideas in this space in a short time span of around two years.

\begin{itemize}
\item In Aug--Oct 2018, \citet*{howard2020time,howard2021time} identified one of the definitions of an e-process for composite nulls: a nonnegative process that is upper bounded by a family of supermartingales. They combined these with Ville's inequality to unify and improve a host of nonparametric Chernoff-style concentration inequalities. They developed several variants of the method which yielded computationally and statistically efficient  nonparametric sequential tests and confidence sequences, extending and generalizing Robbins' extensive work on the topic.

\item In Nov 2018, \citet*{kaufmann2021mixture} also used Robbins' method of mixtures applied to likelihood ratios arising from one-parameter exponential families. They then also used Ville's inequality to derive time-uniform concentration bounds for the Kullback--Leibler divergence (relative entropy) between the distributions parameterized by the empirical and true means. This yielded  sequential tests and confidence sequences with applications to parametric multi-armed bandit problems.

\item In Mar 2019, \citet*{Shafer:2021} wrote an expository article promoting the use of betting in scientific communication and arguing for the adoption of the log-optimality criterion. The paper thoroughly develops the story for testing a simple null against a simple alternative, proving that the log-optimal bet is their likelihood ratio. 

\item In Mar 2019, \citet*{Shafer/Vovk:2019} published a book on the game-theoretic foundations of probability and finance, providing a solid philosophical and mathematical basis for these fields. Many ideas in game-theoretic statistics owe their philosophical or methodological roots to ideas presented in this book. Certainly, this book predates all papers mentioned here as Mar 2019 is its publication time.

\item In Jun 2019, \citet*{GrunwaldHK19} identified what we now recognize as the simpler of the two definitions of an e-process for composite nulls: a nonnegative process that is an e-value at any stopping time. They also argued for the log-optimality criterion, pointed out the key role played by the reverse information projection in the design of log-optimal e-variables, and argued that meta-analysis (of an increasing number of studies) is a particularly sensible application area for these concepts. 

\item In Dec 2019,  \citet*{Vovk/Wang:2021} developed fundamental results related to  combining e-values, and transforming  e-values to p-values or vice versa (calibrators). For instance, they prove that the arithmetic average is essentially the only admissible symmetric function that can combine dependent e-values, and that inverting an e-value is the only admissible way to convert it to a p-value. Before this paper, authors used different names for the same concept, but following this paper, there was a joint agreement amongst authors to use the term ``e-value''.

\item In Dec 2019, \citet*{wasserman2020universal} developed particularly simple and broadly applicable e-variables and e-processes within their framework of ``universal inference''. They proved that whenever the maximum likelihood under the (composite) null hypothesis can be calculated, an e-variable can be constructed, leading to a valid (fixed-time or sequential) test by Markov's inequality.

\item In Jul 2020, \citet*{vovk2022admissible} showed that admissible ways of merging p-values under arbitrary dependence have to go through   e-values. This result is based on a classic optimal transport duality that converts probability constraints under arbitrary dependence into expectation constraints. Using this connection, many admissible methods of merging p-values are constructed from merging e-values.

\item In Sep 2020, \citet*{wang2022false} developed an analog of the Benjamini--Hochberg (BH) procedure that employed e-values instead of p-values (e-BH). They showed that unlike BH, the e-BH procedure controlled the false discovery rate under arbitrary dependence between the e-values. 

\item In Sep 2020, \citet*{ramdas2020admissible} identified  necessary and sufficient conditions for admissibility for SAVI tools (e-processes, p-processes, confidence sequences, sequential tests), and prove that all admissible SAVI methods must be based on nonnegative martingales. They also prove that the two definitions of e-processes mentioned above are actually equivalent. 

\item In Oct 2020, \citet*{waudby2020estimating} employed a nonparametric testing by betting approach to produce state-of-the-art solutions to the classic problem of estimating a bounded mean (when sampling with or without replacement). In doing so, they provide a universal representation for all composite nonnegative martingales, and a simple betting strategy that easily extends to other problems.



\end{itemize}

The pace and quality of progress did not stop here, and this book aims to collect many of the results presented in these above works, and many others that followed them.
The bibliographic notes after each chapter mention many relevant papers in addition to the above, but they are by no means exhaustive. 


\subsection{Unfortunate name clashes: other e-values in the literature}

Given the popularity of p-values, many other terms have been created in the literature by swapping the p in p-value for other letters. All the above authors had already repeatedly contemplated other options, such as the s-value (s for safe), v-value (v to honor Ville and Vovk), b-value (b for bet), m-value (m for martingale), and i-value (i for integral). In fact, the name i-value was used informally by a few individuals in the 1990s to describe a concept that is essentially an e-value in a different context, where the constraint is an integral (more general than an expectation) bounded by 1.

In the end, we as a community stuck with e-value, because it was defined by its expectation, and is directly interpreted as evidence against the null.  
There was also a consideration of its pronunciation and contrast to p-values.
However, as time has passed, we have discovered many other e-values in the literature that we were initially unaware of. We mention them below to avoid confusions for the casual reader, who may search for the term e-value and may mistakenly wander into a different literature with other goals and definitions. We apologize for any unintentional misrepresentations of their concepts, on which we are not experts.

Importantly, terms like ``e-variable'' and ``e-process'' do not have known clashes, and the reader is reassured when these terms accompany ``e-value'', so we recommend mentioning them in   papers wherever suitable.

\begin{enumerate}
    \item \textbf{BLAST E-value.} BLAST \citep{altschul1990basic} is an extremely popular software tool used in bioinformatics. 
    In the context of sequence alignment given a database of sequences, the BLAST E-value (E for Expect) is the number of expected hits of similar quality (score) that could be found just by chance. Their E-value appears to be a particular calibrator applied to a particular p-value: $E = -\log(1-P)$. This appears to keep E-values and p-values on the same scale (lower E-values point towards statistical significance). This is related to Greenland's definition of an S-value (S for surprisal or Shannon information), which is $-\log(P)$, which would be a special case of our e-value. 
    
    \item \textbf{Sensitivity E-value.} \cite{vanderweele2017sensitivity} introduced E-values in the particular context of sensitivity  analysis in causal inference.  The E-value is defined as the minimum strength of association, on the risk ratio scale, that an unmeasured confounder would need to have with both the treatment and the outcome to fully explain away a specific treatment-outcome association, conditional on the measured covariates. A large E-value implies that considerable unmeasured confounding would be needed to explain away an effect estimate. A small E-value implies little unmeasured confounding would be needed to explain away an effect estimate. This notion has become extremely popular in applied causal inference.

    \item \textbf{Bayesian E-value}. The e-value $ev(H|X)$ --- the \emph{epistemic}-value of hypothesis $H$ given observations $X$, or the
evidence-value of observations $X$ in favor (or in support) of hypothesis $H$ --- is a Bayesian statistical
significance measure introduced   together with the Full Bayesian Significance Test~\citep{pereira1999evidence}. Concise entries about the e-value and the FBST are available at the International Encyclopedia of Statistical Science and online at Wiley’s StatsRef. The recent survey by~\cite{stern2022value} points to rich and growing literature on the topic over the last 25 years. 
    
\end{enumerate}

\section{Road map}
\label{sec:c1-roadmap}

The rest of this book is organized into three parts.  
The first part, \textbf{Fundamental Concepts},  includes four chapters (including the current one) that introduce the basic concepts. \begin{itemize}
\item  
 Chapter~\ref{chap:markov} concerns the properties of e-values under the null hypothesis, including key concepts such as  Markov's inequality, calibrators, and randomized tests. It presents tools to convert e-values to p-values and tests, and vice versa. It also presents explicit e-variables in several simple examples. This chapter focuses on \emph{validity} (under the null), leaving questions of \emph{statistical efficiency} to the following chapter.
 
 \item 
 Chapter~\ref{chap:alternative} concerns 
 the properties of   e-values under the alternative  hypothesis, including
 topics of e-power, log-optimality, existence of powered e-values, and mixture and plug-in methods. 
 As an important result, the log-optimal e-variable for a simple null against a simple alternative is their likelihood ratio. It provides methods to derive ``good'' or ``useful'' e-values, and shows how e-values, p-values and tests are fundamentally connected via a unified existence result.
 
 \item Chapter~\ref{chap:posthoc} demonstrates  that e-values play a fundamental role in post-hoc testing and decision making, where some parts of the testing problem, available decisions, or data collection  depend on the current data. While the previous chapters focused on putting e-values and p-values on a similar footing, this is the first chapter to clearly separate their utility: what one can do with e-values versus what one can (or cannot) do with p-values.
\end{itemize} 

The second part, \textbf{Core Ideas}, includes five chapters. \begin{itemize}
\item   Chapter~\ref{chap:ui} introduces the general method of universal inference for constructing e-variables for general composite nulls and alternatives. This method is very simple to apply, and provides the first known test for many ``irregular'' testing problems. However, it is inadmissible in a strong sense: it is always dominated by the e-variable constructed in the next chapter. 
\item 
 Chapter~\ref{chap:numeraire} contains a key result: for a given composite null and simple alternative, there always exists a unique log-optimal e-variable (called the numeraire), and presents a strong duality result with the reverse information projection (RIPr) of that alternative onto the null. Remarkably, the numeraire is also the only e-variable that is a likelihood ratio between the alternative and some element of the (effective) null.   
 \item 
 Chapter~\ref{chap:eprocess} presents the framework of sequential anytime-valid inference (SAVI), emphasizing key concepts such as e-processes, testing by betting, optional stopping and Ville's inequality. It also discusses key differences between SAVI and Wald's sequential testing paradigm, as well as the problematic practice of peeking at p-values.
 \item 
 Chapter~\ref{chap:multiple} discusses merging e-values via averages (arbitrary dependence) or products (independence or ``sequential'' dependence), as well as the trade-offs in choosing merging methods. 
 \item 
 Chapter~\ref{chap:compound} presents the notion of  compound e-values in the context of multiple hypothesis testing,  introduces the e-Benjamini--Hochberg (e-BH) procedure and proves that every FDR controlling procedure can be seen as applying e-BH to some set of compound e-values. 
\end{itemize}
 
The third part, \textbf{Advanced Topics}, has seven  chapters.\begin{itemize}
 \item   
Chapter~\ref{chap:approximate}  studies approximate and asymptotic p-values and e-values, as well as their counterparts in multiple testing. This chapter provides a bridge between classical statistics, which often relies on asymptotics, and most  our book, which is nonasymptotic.
   \item  Chapter~\ref{chap:combining} shows that the only admissible way to combine dependent e-values is by (weighted) averaging. It also develops methods to combine a p-value  and an e-value, which give rise to an e-weighted BH procedure.   
   \item 
 Chapter~\ref{chap:pmerge} derives methods for combining multiple dependent p-values, which are proved to rely on e-values through optimal transport duality, and also ways to improve on such combination rules by using random and exchangeable improvements of Markov's inequality.   
 \item 
 Chapter~\ref{chap:eci} defines e-confidence intervals (e-CI), and presents the e-Benjamini--Yekutieli (e-BY) procedure to control the false coverage rate after selecting a subset of e-CIs to report. It also discusses how to merge dependent confidence intervals using (various versions of) majority vote.  
 \item   
Chapter~\ref{chap:theoretic} offers further contrasts between p-values and e-values, including a duality between them, a representation of e-variables and p-variables as conditional expectations and probabilities, and an existence result on powered exact e-values and p-values. 
\item  Chapter~\ref{chap:shape} shows how Markov's inequality can be improved under additional  conditions on the shape of the distribution of a single e-variable  or on the dependence between multiple e-variables.
  \item  Chapter~\ref{chap:risk-measures} studies the connection between e-values and risk measures, and offers some nonparametric methods for testing the forecast of general statistical functions (such as the mean, the variance, or a quantile) using e-values.
\end{itemize} 

In terms of the goal of each chapter, all chapters can be grouped in the following three categories. 

\begin{enumerate}
    \item Understanding the foundation of e-values: 
    Chapters~\ref{chap:markov},~\ref{chap:posthoc},~\ref{chap:approximate}, and~\ref{chap:theoretic}.
    \item  Constructing good and powerful e-values in single hypothesis testing: Chapters~\ref{chap:alternative},~\ref{chap:numeraire}, and~\ref{chap:shape}. 
    \item Designing inference methods with e-values: Chapters~\ref{chap:ui},~\ref{chap:eprocess}, 
and~\ref{chap:risk-measures}.

    \item Designing methods  for multiple testing with e-values: 
Chapters  
    ~\ref{chap:multiple}, 
~\ref{chap:compound},~\ref{chap:combining},  ~\ref{chap:pmerge}, 
and~\ref{chap:eci}.
\end{enumerate}
In the first two categories, we are interested in e-values themselves, whereas in the remaining categories, we are primarily interested in the specific statistical problem.




\chapter{Validity: E-values under the null hypothesis}

\label{chap:markov}

In this chapter, 
we study fundamental issues of e-values, in particular, on their connection to p-values, validity, and tests. These properties concern the behavior of the e-variables under the null hypothesis $\cP$. 
For this reason,
the alternative hypothesis $\cQ$  does not appear and plays no role.

\section{Testing using Markov's inequality}
\label{sec:c2-markov}
As mentioned before, e-values can be interpreted in their own right as evidence against the null hypothesis without resort to other concepts. However, the reader may naturally be curious how exactly e-values may be transformed (perhaps with some loss of efficiency) into the more classical concepts of level-$\alpha$ tests and p-values. It is in this context that a central role will be played by Markov's inequality (though it apparently was discovered by Chebyshev, Markov's advisor).


Markov's inequality guarantees that a test  rejecting an e-value larger than $1/\alpha $ yields a level-$\alpha$ test.

\begin{proposition}[Markov's inequality for e-values]\label{prop:markov}
Let $E$ be an e-variable  for $\cP$. We have
$
\p(E\ge 1/\alpha) \le \alpha
$
for all $\p\in \mathcal P$ and  $\alpha \in (0,1]$. 
Hence, $1/E$ is a p-variable, $(\alpha E) \wedge 1$ is a level-$\alpha$ test, and $\id_{\{E \geq 1/\alpha\}}$ is a level-$\alpha$ binary test.
\end{proposition}
The proof follows directly from  the simple inequality 
$\id_{\{x\ge 1 \}}\le x
$ for all $x \ge 0$, by observing that
 \begin{align*}
        \p(E\ge 1/\alpha) = \E^\p[\id_{\{\alpha E\ge 1\}} ] \le \alpha \E^\p[E]\le \alpha.
    \end{align*}

For an e-variable $E$ and a p-variable $P$, even though $1/E$ is a p-variable,
  $1/P$ is generally not an e-variable; in fact $1/P$ has infinite mean if $P$ is exact.
The p-variable $1/E$ is conservative since it cannot be  uniformly distributed on $[0,1]$, and it is typically far away from being so. 
As a consequence, the type-I error
of the test $\id_{\{E\ge 1/\alpha\}}$ for $\p$ 
 is usually smaller than $\alpha$ unless $E$ is distributed on the two points $0$ and $1/\alpha$ with mean $1$ under $\p$.

Proposition~\ref{prop:markov} can be generalized in several different directions relevant for testing with e-values; some are presented in Section~\ref{sec:2-rand} and  Chapters~\ref{chap:ui} and~\ref{chap:shape}. A sequential generalization, called Ville's inequality, is presented in Chapter~\ref{chap:eprocess}.

\begin{figure}[t]
   \begin{center}
      \includegraphics[width=0.8\textwidth]{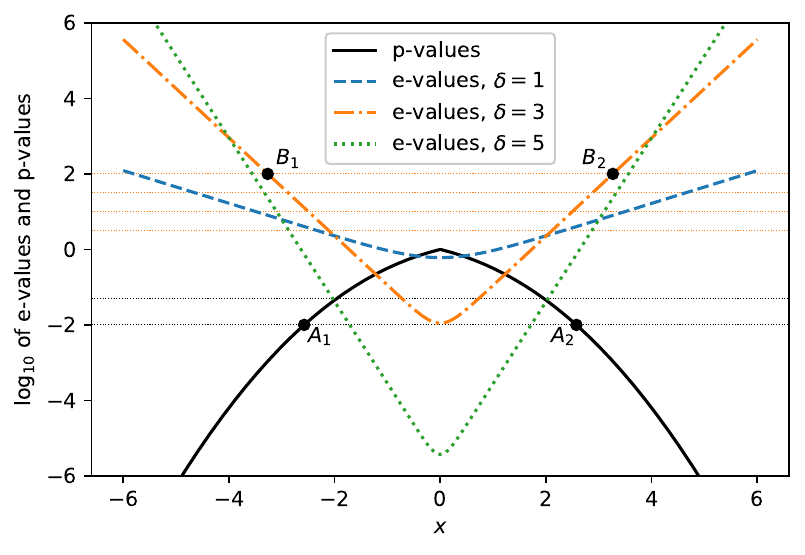}

   \end{center} 
   \caption{A comparison of e-values and p-values for the two-sided normal test in Section~\ref{sec:LR-e-variable}. The curves represent the e-values or p-values as a function of the test statistic $Z=x$. 
   The top four horizontal dotted lines correspond to $10^{\beta}$ with $\beta \in \{2, 1.5,1,0.5 \}$ for e-values.
      The bottom two horizontal dotted lines correspond to $\alpha \in \{0.05,0.01\}$ for p-values. (The reason why e-values are compared with levels $10^{\beta}$ with $\beta \in \{2, 1.5,1,0.5 \}$  is explained in Section~\ref{sec:c2-jeffrey}.)
   }
   \label{fig:c2-1}
   \end{figure}

With a p-variable $P$ and an e-variable $E$ for the same hypothesis, the tests $\id_{\{E\ge 1/\alpha\}}$ 
and $\id_{\{P\le \alpha\}}$
are in general different, except for the special case  $E = \id_{\{P \leq \alpha\}} / \alpha$
(see Section~\ref{sec:nontrivial-test} for more details on such e-variables). 
To understand how these two tests differ,  we consider the two-sided normal test in Section~\ref{sec:LR-e-variable},
where
a family of e-variables $E_n$ is given by \eqref{eq:example-gaussian-e2}, indexed by $\delta>0$,
and the p-variable $P_n$ is given by \eqref{eq:example-gaussian-p2}.
We treat both the e-variables and the p-variable as functions of the value $x$ taken by the statistic $Z$. 
Figure~\ref{fig:c2-1} illustrates 
two intervals $[A_1,A_2]$
and $[B_1,B_2]$, which correspond to $P_n(x)= 0.01$ for $x=A_1$ or $x=A_2$
and $E_n(x)=100$ for $x=B_1$ or $x=B_2$, with $\delta=3$. 
Here, $A_2=|A_1| \approx 2.58$ and $B_2=|B_1|\approx 3.27$. 
Other choices of $\delta$ lead to different intervals, as shown in Figure~\ref{fig:c2-1}. 
We can see that it is harder for an e-value to reach $100$ compared to a p-value to reach $0.01$, and this is because the Markov inequality is generally conservative. 
It is generally unfair to directly compare the values of $E$ and $1/P$. See Section~\ref{sec:c2-jeffrey} for more discussions on this point.

\section{Testing using distributional knowledge}
\label{sec:dist-info}

Proposition~\ref{prop:markov} implies that for an e-variable $E$ and $\alpha \in(0,1)$, we can always reject the null hypothesis by using the level-$\alpha$ test $\id_{\{E\ge 1/\alpha\}}$. 
If the distribution of $E$ under the null hypothesis  is known or computable (there may be multiple distributions if the null is composite), then we can find thresholds that are smaller than $1/\alpha$ to test based on $E$.
Denote by $S(\cdot;\p)$ the left-continuous survival function  of $E$ under $\p\in \cP$, that is, 
$$S(x;\p)=\p(E \ge x)$$ for $x\in \R.$
If $S(\cdot;\p)$ is the same across $\p\in \cP$, then $E$ is called pivotal (this will appear again in Section~\ref{sec:existence-exact}). 
The function $S(\cdot;\p)$ can be continuous or discontinuous, depending on both $E$ and $\p$.

For a threshold $t>0$,   we would like the test $\id_{\{E\ge t\}}$ to have level $\alpha$. Clearly, choosing $ t \ge1/ \alpha$ is sufficient, but we can do better. Indeed, by definition, as long as $\sup_{\p\in \cP} S(t|\p)\le \alpha$, $t$ is a valid threshold for the test $\id_{\{E\ge t\}}$ with level $\alpha$.
If $\sup_{\p\in \cP} S(t|\p)$ is continuous, then 
the smallest $t$ that satisfies this condition is
$$
t^*=\inf\left\{t>0: \sup_{\p\in \cP} S(t|\p)\le \alpha\right\},
$$
which is also equal to the largest  left $\alpha$-quantile of $E$ over $\p\in \cP$ (Lemma~\ref{lem:1} in Section~\ref{sec:c12-1}).
In particular, if $\cP=\{\p\}$ is simple, then we can easily see 
that 
$$
t^*=\inf \left\{t>0: \p(E\ge t)\le \alpha\right\} = \inf \left\{t>0: \p(E > t)\le \alpha\right\},
$$
which is the left  $\alpha$-quantile of $E$ under $\p$. 
The level-$\alpha$ guarantee of $\id_{\{E \ge t^*\}}$  requires 
$t \mapsto \sup_{\p\in \cP} S(t|\p)$ to be continuous  at $t^*$; otherwise we can use the test
 $\id_{\{E > t^*\}}$, or a suitable randomization between $\id_{\{E > t^*\}}$ and  $\id_{\{E \ge  t^*\}}$ (analogously to the Neyman--Pearson lemma).

    Reducing the threshold in the above manner is in fact no different from treating the e-value as any other test statistic and picking the threshold appropriately to derive a level-$\alpha$ test, as is standard in classical statistics.

    Even in situations where $S(\cdot;\p)$ is difficult to fully specify or compute, 
 bounds on $t^*$ may be obtained based on partial information of the distribution of the e-value. Several concrete cases will be studied in Chapter~\ref{chap:shape}.

\section{Relating p-values and e-values using calibrators}
 \label{sec:2-cali}
P-values and e-values can be converted between each other, via   calibrators.
Below, when we say ``any hypothesis'', we mean any specification of  the sample space $(\Omega,\cF)$ and hypothesis $\cP$ on the space.

\begin{definition}[Calibrators] \label{def:calib}
\begin{enumerate}
\item [(i)] A \emph{p-to-e calibrator} is a decreasing function $f:[0,\infty)\to[0,\infty]$ satisfying
 $f=0$ on $(1,\infty)$  such that for any hypothesis, $f(P)$ is an e-variable for any p-variable $P$.
\item [(ii)] An \emph{e-to-p calibrator} is a decreasing function $g:[0,\infty]\to [0,\infty)$ 
such that for any hypothesis, $g(E)$ is a p-variable for any e-variable $E$. 
\item [(iii)]
A calibrator $f$ is said to \emph{dominate} a calibrator $g$ if $f\ge g$ (p-to-e) or $f\le g$ (e-to-p),
and the domination is \emph{strict} if further $f\ne g$.
A calibrator is \emph{admissible} if it is not strictly dominated by any other calibrator.
\end{enumerate}
\end{definition}

We often omit ``p-to-e'' in ``p-to-e calibrators'' and simply call them calibrators, for a reason that will be made clear soon in Propositions~\ref{prop:e-to-p} and~\ref{prop:p-to-e}.

The interval $[0,\infty)$ can be safely replaced by $[0,1]$ in the above definitions since p-values with values larger than $1$ are not interesting for testing. 

A probability measure $\p$ (or the corresponding probability space) is called \emph{atomless} if there exists a continuously distributed random variable (e.g., uniform, normal) under $\p$.
Hence, when we encounter a continuously distributed random variable, implicitly the probability measure is assumed atomless. 
  In Definition~\ref{def:calib}, the property of calibrators needs to hold for any hypothesis. Nevertheless, it suffices to consider one atomless probability measure $\p$ on some $(\Omega,\cF)$, because if for the simple hypothesis $\p$, $f(P)$ is an e-variable for any p-variable $P$, then it holds also for any other, simple or composite, hypotheses.
 This observation will be useful in several places throughout. 
 A formal justification of this simplification is presented in Appendix~\ref{app:atomless}.

In the next two results we characterize all admissible calibrators in both directions.  
We first look at the simpler direction.  
The following proposition 
says that there is, essentially, only one e-to-p calibrator, $f:t \mapsto \min(1,1/t)$, that is offered by Markov's inequality. 

\begin{proposition}\label{prop:e-to-p}
  The function $f:[0,\infty]\to[0,1]$ defined by $f(t)=\min(1,1/t)$ is an e-to-p calibrator.
  It dominates every other e-to-p calibrator.
  In particular, it is the only admissible e-to-p calibrator. 
\end{proposition}

\begin{proof}[Proof.]
  The fact that $t\mapsto \min(1,1/t)$ is an e-to-p calibrator follows from Markov's inequality in Proposition~\ref{prop:markov}. 
  On the other hand, suppose that $g$ is another e-to-p calibrator.
  It suffices to check that $g$ is dominated by $t\mapsto  \min(1,1/t)$.
  Suppose $g(t)<\min(1,1/t)$ for some $t\in[0,\infty]$.
  Fix an atomless probability measure $\p$.
  Consider two cases:
  \begin{enumerate}
  \item[(a)]
    If $g(t)<\min(1,1/t)=1/t$ for some $t>1$,
    fix such $t$ and consider an e-variable $E$ that is equal to $t$ with probability $1/t$ and $0$ otherwise under $\p$.
    Then $g(E)$ is $g(t)<1/t$ with probability $1/t$,
    whereas it would have satisfied $\p (g(E)\le g(t))\le g(t)<1/t$ had it been a p-variable.
  \item[(b)]
    If $g(t)<\min(1,1/t)=1$ for some $t\in[0,1]$,
    fix such $t$ and consider the e-variable $E=t$.
    Then $g(E)=g(t)<1$,
    and so it is not a p-variable.
    \qedhere
  \end{enumerate}
\end{proof}


Instead of working with the set of all e-variables, 
we can also consider \emph{conditional}
e-to-p calibrators that are valid only on some subset $\mathcal E$ of e-variables, which is a weaker requirement,
and such calibrators depend on the choice of   $\mathcal E$. 
 This topic is treated in Chapter~\ref{chap:shape}.
  
In sharp contrast to the uniqueness of the admissible e-to-p calibrator, the class of p-to-e calibrators is much richer.

\begin{proposition}\label{prop:p-to-e}
  A decreasing function $f:[0,\infty)\to[0,\infty]$ with $f=0$ on $(1,\infty)$ is a p-to-e calibrator if and only if $\int_0^1 f (p) \d p\le 1$.
  It is admissible if and only if $f$ is upper semicontinuous, $f(0)=\infty$, and $\int_0^1 f(p) \d p = 1$.
\end{proposition}

In the context of this proposition,
being upper semicontinuous is equivalent to being left-continuous.

\begin{proof}[Proof.] 
  The first ``only if'' statement is obvious.
  To show the first ``if'' statement,
  suppose that $\int_0^1 f (p)\d p\le 1$, $P$ is a p-variable for an atomless probability measure $\p$,
  and $P'$ is uniformly distributed on $[0,1]$.
  Since $\p(P\le x) \le x \wedge 1= \p(P'\le x)$ for all $x\ge 0$ and $f$ is decreasing, 
  we have 
  \[
    \p(f(P)>x) \le \p (f(P')>x)
  \]
  for all $x\ge0$, which implies
  \[
    \E^\p[f(P)]
    \le
    \E^\p[f(P')]
    =
    \int_0^1 f(p) \d p \le 1.
  \]
  Both necessity and sufficiency in the second statement of Proposition~\ref{prop:p-to-e} are straightforward.
\end{proof}

Below we present some simple examples of admissible p-to-e calibrators. We only specify them on $(0,1]$, as they are $\infty$ at $0$ and $0$ on $(1,\infty)$. 
These calibrators 
include the power form 
\begin{equation}\label{eq:cal1}
f(p)=   \kappa p^{\kappa-1} \mbox{~~~~for some $\kappa \in (0,1)$},
\end{equation}
the mixture of \eqref{eq:cal1},  
\begin{equation}\label{eq:cal2}
  f(p)
  =
  \int_0^1 \kappa p^{\kappa-1} \d\kappa
  =
  \frac{1-p+p\log p}{p(-\log p)^2},
\end{equation}
and some special simple forms 
\begin{align}\label{eq:cal-av} f(p)&=2(1-p);\\
\label{eq:cal3} f(p)&=p^{-1/2}-1 ;\\ \label{eq:cal4} f(p)&=-\log p.
\end{align}

Converting a p-value to an e-value using a p-to-e calibrator and then back to p-value using an e-to-p calibrator generally loses quite a lot of evidence.
For instance, starting with $p=0.01$, a conversion with the p-to-e calibrator $p\mapsto p^{-1/2}-1$
gives $e=9$, and another conversion with the e-to-p calibrator $e\mapsto \min(1/e,1)$ 
yields $p'= 1/9$.
Therefore, it is not recommended in practice to calibrate iteratively in both directions.

 A possible interpretation of the above results on calibrators 
  is that e-variables and p-variables are connected via a very rough relation $1/e \sim p$.
  In one direction, the statement is precise:
  the reciprocal (truncated to 1 if needed) of an e-variable is a p-variable by Proposition~\ref{prop:e-to-p}.
  On the other hand,
  using a calibrator $f(p)=   \kappa p^{\kappa-1}$  in  \eqref{eq:cal1}  with a small $\kappa>0$ 
  and ignoring positive constant factors,
  we can see that the reciprocal of a p-variable is approximately an e-variable; this is also the case of \eqref{eq:cal2}.
  In fact, $f(p)\le 1/p$ for all $p$ when $f$ is a calibrator; otherwise $\int_0^1 f(q)\d q\ge pf(p)>1$, conflicting the condition in Proposition~\ref{prop:p-to-e}.
  However, $f(p)=1/p$ for a fixed $p$ is only possible in the extreme case $f:q\mapsto \id_{\{q\in [0,p]\}}/p$.

In the rest of the book, we omit ``p-to-e'' when mentioning calibrators, unless we want to stress the direction. 
Although all calibrators in \eqref{eq:cal1}--\eqref{eq:cal4} are admissible, 
some are more useful. In hypothesis testing, very small p-values and very large e-values are interesting.  Calibrators in \eqref{eq:cal-av}, \eqref{eq:cal3} and \eqref{eq:cal4} penalize  very small p-values more than \eqref{eq:cal2}, which can be seen from their asymptotic behaviour as $p\downarrow 0$,
and therefore they are less likely to produce very large e-values; the calibrator in  \eqref{eq:cal-av} is even bounded by $2$ for $p\ne 0$. Nevertheless, they may be useful when we take the product of many independent e-values.

We end this section with a simple property of admissible calibrators. It also explains why we choose the domain of the calibrators as $[0,\infty)$ instead of $[0,1]$, although its value is set to $0$ for input p-values larger than $1$. This result will be useful for later results in Chapter~\ref{chap:pmerge}.
\begin{proposition}\label{prop:rescale-cal}
If $f$ is a calibrator or admissible calibrator, then so is the mapping $x\mapsto a f(a x)$
for any $a\ge 1$.
    \end{proposition}
    \begin{proof}[Proof.] 
    Suppose that $f$ is an admissible calibrator. 
    Denote by $g:x\mapsto a f(a x)$.
    We need to check a few properties in  Proposition~\ref{prop:p-to-e}. 
    First,  it follows from $f(p)=0$ for $p>1$ that $g(p)=af(ap)=0$ for $p>1$, and further 
    $$
     \int_0^1 g(x) \d x =\int_0^1 a f(ax) \d x = 
 \int_0^{1/a}     a f(ax) \d x = 
    \int_0^1 f(p) \d p=1.
    $$
Upper semicontinuity and $g(0)=\infty$ follow directly from those properties of $f$.
Therefore, $g$ is an admissible calibrator. 
The same argument above shows that $g$ is also a calibrator when $f$ is a calibrator that is not necessarily admissible. 
    \end{proof}

The next result states that any   e-variable $E$ for an atomless probability measure $\p$ can be seen as some calibrator applied to  an exact p-variable. Indeed, as we see from the proof of the result below, this calibrator is 
a transform of the quantile function of $E$. 
\begin{proposition}
Let $E$ be an e-variable for an atomless probability measure $\p$. Then, there exists a calibrator $f$ and an exact p-variable $P$ for $\p$ such that $E=f(P)$, $\p$-almost surely.
Moreover, if $E$ is exact, then we can require the calibrator $f$ to be admissible. 
\end{proposition}
\begin{proof}[Proof.]  
For $p\in (0,1]$, let $f(p)$ be given by 
$
f(p)=\inf\{x\in \R: \p(E\le x)>1-p\}
$, that is,  the right quantile of $E$ at level $1-p$, 
 and further  we set $f(0)=\infty$ and $f(x)=0$ for $x>1$.
By a standard property of the quantile function (formally stated in Lemma~\ref{lem:quantile2}  in Appendix~\ref{app:quantile}),
there exists an exact p-variable $P$ such that  $E=f(P)$ almost surely. 
In case $E$ is exact, the admissibility of this calibrator  follows from  Proposition~\ref{prop:p-to-e} and the fact that the right quantile function is upper semicontinuous.  
\end{proof}


\begin{example}
In  Section~\ref{sec:LR-e-variable}, 
 for testing $\p=\mathrm{N}(0,1)$ against  $\q=\mathrm{N}(\mu,1)$ with $n$ observations, the likelihood ratio e-value  in \eqref{eq:example-gaussian-e3} is 
 $ E= 
  \exp \left( \delta Z -\delta^2/2\right),
 $ and
the Neyman--Pearson p-variable is  
 $ P=  \Phi(-Z), $
  where $Z$ is $\mathrm N(0,1)$-distributed, $\Phi$ is the standard normal cdf, and $\delta=\mu \sqrt n$.  
  Both $E$ and $P$ are exact. 
It is straightforward to see $E  =  f(P) $ 
  with the admissible calibrator $f$ given by
  $$
  f(p) =  \exp \left( -\delta  \Phi^{-1} (p)  -\delta^2/2\right) .
  $$
\end{example}

\section{Randomized tests and calibrators}
\label{sec:2-rand}

A randomized test is a test that utilizes independent external randomization.
It is possible to improve Markov's inequality using external randomization, and this has direct implications for converting e-variables into tests or p-variables. 

\begin{theorem}[Randomized Markov inequality]
\label{thm:umi}
Fix $\p\in \cM_1$.
Let $Y$ be a nonnegative random variable and $U\lawis \mathrm{U}[0,1]$ be independent of $Y$. Then, for any $\alpha > 0$, 
\begin{equation}\label{eq:UMI}
\p(Y \geq U/\alpha ) = \E^\p[\min(\alpha   Y,1)] \leq \alpha \E^\p[Y] .
\end{equation}
In particular, if $Y$ is an e-variable for $\cP$, then 
$\p(Y \geq U/\alpha ) \leq \alpha$ for any $\p\in \cP$.
\end{theorem}
   The proof is simple. Since $U$ is uniform and independent of $Y$, we get
\begin{equation}\label{eq:UMI-proof}
\p(Y \geq U/\alpha ) = \E^\p[\p( U \leq \alpha  Y \mid Y)] = \E^\p[\min(\alpha  Y,1)],
\end{equation}
yielding the claim.
The first equality in~\eqref{eq:UMI}  becomes an inequality $\leq$ if $U$ is stochastically larger than standard uniform, and in particular $U=1$ yields an alternative proof of Markov's inequality. Also, if $Y$ is bounded, that is $Y \in [0,C]$ almost surely for some $C > 0$, then the inequality in~\eqref{eq:UMI} holds with equality for any $\alpha   \leq 1/C$. In contrast, Markov's inequality holds with equality only for distributions supported on $\{0,1/\alpha  \}$.


The following corollary is a direct consequence of Theorem~\ref{thm:umi}, which offers additional interpretations of the relation between p-variables, e-variables, and randomization. 
\begin{corollary}
\label{coro:ip-comb}
For any hypothesis $\cP$, 
if  a p-variable $P$ and an e-variable $E$   are independent, then $P/E $ is a p-variable. In particular, if $E$ is an e-variable, then $U/E $ is a p-variable, where $U\lawis \mathrm{U}[0,1]$ is independent of $E$. 
\end{corollary}

Since $U < 1$ almost surely, $U/E$ is a strictly smaller p-value than $1/E$ in all but degenerate cases. Thus, the corresponding randomized test is typically strictly more powerful than the non-randomized one. 
Said differently, $e\mapsto (U/e)\wedge 1$ can be seen as a \emph{randomized e-to-p calibrator}.
In fact, this randomized e-to-p calibrator is admissible (treated in Section~\ref{sec:10-cross}) and it dominates $e\mapsto (1/e)\wedge 1$, the latter only being admissible among non-randomized e-to-p calibrators (Proposition~\ref{prop:e-to-p}). 

A more general treatment of combining a p-value and an e-value is presented in Section~\ref{sec:10-cross}.



\section{Admissible e-variables}
\label{sec:c2-admissible}

In this section, we consider admissible e-variables, 
which are, intuitively, e-variables that cannot be enlarged, formally defined below. 
  \begin{definition}
      \label{def:adm-e}
        An e-variable   $E$  for $\mathcal P$ is  \emph{admissible}
    if for any other e-variable $E'$ for $\cP$ satisfying $E'\ge E$, we have $ \p(E'=E)=1 $ for each $\p\in \cP$.
  \end{definition}

Equivalently,  we may say that  an e-variable is admissible if it is not dominated by another e-variable,
where by domination of $E'$ over $E$, we mean  $E' \geq E$  and $\p(E' > E) > 0$ for some $\p\in \cP$. 

A trivial example is that the constant $1$ is an admissible e-variable for $\cP$.

The above definition does not use any alternative $\q$, but 
if all $\p\in \cP$ and $\q$  are mutually absolutely continuous, 
then $E$ being admissible means that one cannot find another e-variable $E'$ such that $\q(E'\ge E)=1$ and 
$\q(E' > E) > 0$.  
Indeed,  
if  $\q(E'\ge E)=1$ and 
$\q(E' > E) > 0$ hold for an e-variable $E'$, then 
$E'\vee E$ is   an e-variable that dominates $E$. 

\subsection{Admissible e-variables for a simple null are likelihood ratios}
\label{sec:c2-likelihood}
We first make a straightforward but important observation: For simple nulls, every admissible e-variable can be written as a likelihood ratio.


\begin{proposition}
\label{prop:admissible-lr}
For a simple null $\p$, an  e-variable $E$ is admissible if and only if it is exact. Further, any exact e-variable for $\p$ equals a likelihood ratio ${\d\s}/{\d\p}$ $\p$-almost surely for some distribution $\s \ll \p$.
\end{proposition}
\begin{proof}[Proof.]
 If an e-variable $E$ is not exact, meaning that $\E^\p[E] < 1$,
 then $ E+(1-\E^\p[E])$ is an e-variable dominating $E$, and hence $E$ cannot be admissible.
    To show the converse, it suffices to note that if $E'$ dominates an exact e-variable $E$, then $\E^\p[E']>1$; hence, exact e-variables for $\p$ are admissible.
For the last statement, 
we can define $\s$ via the equality 
$\d \s/\d \p=E$,
and because $E $ is nonnegative and exact, we know that $\s$ is a  distribution that is absolutely continuous with respect to $\p$.
\end{proof}

  \subsection{Characterization of admissibility}

We now characterize admissible e-variables for general $\cP$.    
  For $A\in \mathcal F$ and $\delta>0$, let $\cP^\delta_A=\{\p\in \cP:\p(A)>\delta\}$
  and $\cP_A=\{\p\in \cP:\p(A)>0\}$.

The first result gives a simple necessary condition for admissible e-variables, which boils down to the exactness condition in Proposition~\ref{prop:admissible-lr} when  $\cP=\{\p\}.$
  \begin{proposition}
  \label{prop:adm-e-1}
If $E $ is an admissible e-variable for $\mathcal P$, then  for any $A\in \mathcal F$ with $\cP_A\neq \varnothing$, we have $\sup_{\p\in \cP_A} \E^\p[E]=1$. 
  \end{proposition}
  \begin{proof}[Proof.]
Suppose $\sup_{\p\in \cP_A} \E^\p[E]<1$ for some $A\in \mathcal F$. Let $\epsilon =1-\sup_{\p\in \cP_A} \E^\p[E]>0$ and  $E'= E+\epsilon \id_A$. 
Clearly $ \E^\p[E'  ]\le 1-\epsilon+\epsilon= 1$ for $\p\in \cP_A$,
  and $ \E^\p[E'  ] =  \E^\p[E   ] \le   1 $ for $\p\in \cP\setminus \cP_A$.
  Therefore, $E'$ is an e-variable for $\cP$ dominating $E$, a contradiction.
  \end{proof}

 Then next result gives a full characterization of admissible e-variables using its expectations under $\p\in \cP_A^\delta$ for $A\in \mathcal F$.   In this result, we use the convention $\inf \varnothing =0$. 
       \begin{theorem}
  \label{th:adm-e} 
A random variable $E:\Omega\to[0,\infty]$ is an admissible e-variable for $\cP$ if and only if 
\begin{align}
\label{eq:adm-cond} 
 \liminf_{\delta\downarrow 0} \inf_{\p\in \cP_A^\delta} \frac{ 1- \E^\p[E] }{\delta} =0  \mbox{ for all $A\in \mathcal F$.}
\end{align}
  \end{theorem}
  \begin{proof}[Proof.]
  We first show the ``if'' direction. 
  First,  by choosing $A=\Omega$,  \eqref{eq:adm-cond} implies that
  $(1-\sup_{\p\in \cP}\E^\p[E] )/\delta \to 0$ as $\delta\downarrow0$,
   which  
  implies that  $ \sup_{\p\in \cP}\E^\p[E] = 1$, and thus 
  $E$ is an e-variable.
  Now, 
 suppose that there exists an e-variable $E'$ with $E'\ge E$   and $\p'(E'>E)>0$ for some $\p'\in \cP$. Then, there exists  $\epsilon>0$ such that the set $A:=\{E'-E\ge \epsilon\}$ satisfies  $\p'(A)>0$. Clearly, 
  $E'\ge E+\epsilon \id_A$.
 By \eqref{eq:adm-cond}, there exists  $\delta>0$ such that  $ \sup_{\p\in \cP_A^\delta} \E^\p[E] >1-\epsilon \delta$. It follows that
 $$
 \sup_{\p\in \cP_A^\delta} \E^\p[ E'] \ge  \sup_{\p\in \cP_A^\delta} \E^\p[ E+\epsilon \id_A]
 > 1-\epsilon \delta + \epsilon \delta=1.
 $$
 This shows that $E'$ cannot be an e-variable. Hence, $E$ is an admissible e-variable.
  
  Next, we show the ``only if'' direction. Suppose that \eqref{eq:adm-cond} does not hold; that is, there exist an event $A\in \mathcal F$  and a constant $c>0$  such that 
  $$ \sup_{\p\in \cP_A^{\delta_n }} \E^\p[E] <1- c \delta_n$$ for $n$ large enough, where  we take $\delta_n=2^{-n}$. 
  Let $E'=E+c \id_A/2$.
  For $\p\in \cP_A^{\delta_n}\setminus \cP_{A}^{\delta_{n-1}}$, we have 
  $$ \E^{\p}[E']\le 1-c\delta_n +\frac{c}{2} \p(A)  \le 1-c\delta_n + c\delta_{n}\le 1.$$
  For $\p\in \cP\setminus \cP_A$, $ \E^{\p}[E']=\E^\p[E]\le 1$.
  Therefore, for all $\p\in \cP$, $\E^{\p}[E']\le 1$, and $E'$ dominates $E$. 
  This contradicts the assumption that $E$ is an admissible e-variable for $\cP$.    
  \end{proof}

 Several sufficient conditions for \eqref{eq:adm-cond} are easier to check, summarized in the next corollary. 
  \begin{corollary}
\label{coro:adm-e-1}  
  An e-variable $E$ for $\cP$ is admissible if 
  any of the following conditions holds:
  \begin{enumerate}
\item[(a)] for each $A\in \mathcal F$ with $\cP_A\neq \varnothing$, there exists $\p\in \cP_A$ such that $\E^{\p}[E]=1$;
\item[(b)]  for each $\p'\in \cP$, there exists $\p\in \cP$ such that $\p'\ll\p$
and  $\E^{\p}[E]=1$; 
\item[(c)] $E$ is exact. 
\end{enumerate}
 \end{corollary}
 \begin{proof}[Proof.]
 (a) For each $A\in \mathcal F$ with $\cP_A\neq \varnothing$, $ \inf_{\p\in \cP_A^\delta} (1- \E^\p[E])=0 $ for $\delta>0$ small enough, and hence \eqref{eq:adm-cond} holds. 
It is straightforward to check that (b) implies (a), and (c) implies (b). 
 \end{proof}

   The condition (a) in Corollary~\ref{coro:adm-e-1} deserves some further attention. 
   We will say that an e-variable is 
   \emph{strongly admissible} if condition (a) in Corollary~\ref{coro:adm-e-1} holds; note that (a) is implied by either of (b) and (c).
This condition is strictly stronger than admissibility for infinite $\cP$ as characterized in Theorem~\ref{th:adm-e} (for which there need not be any $\p\in \cP$ with $\E^\p[E]=1$).
   When $\cP$ is finite, strong admissibility is equivalent to admissibility, as shown in the next result.

      \begin{proposition}
  \label{prop:adm-e-2}
  Suppose that $\cP$ is finite. 
  Then,  $E$ is an admissible e-variable for $\cP$ if and only if $\sup_{\p\in \cP_A} \E^\p[E]=1$ for all $A\in \mathcal F$ with $\cP_A\neq \varnothing$.
  \end{proposition}
    \begin{proof}[Proof.]
  The ``only-if'' statement follows from Proposition~\ref{prop:adm-e-1}. 
    To show the ``if'' statement,
  first observe that $E$ is an e-variable by choosing $A=\Omega$. 
Then, admissibility of $E$ follows by using condition (a) in Corollary~\ref{coro:adm-e-1}. 
  \end{proof}

          

 \subsection{Admissible e-variables on product spaces}
   
  Next, we consider e-variables defined on product spaces. 
Let $\mathcal P' $ be a set of probabilities on another measurable space $(\Omega',\mathcal F')$.
We write  $\cP\times \cP'$  for the set  $\{\p\times\p': \p\in \cP,~\p'\in \cP'\}$ of product probabilities.
  For    $E$ on $\Omega$ and $E'$ on $\Omega'$, 
  let $E\times E'$ be defined by  $(\omega,\omega')\mapsto E(\omega)E'(\omega')$.
Clearly, if $E$ is an e-variable for $\cP$ and $E'$ is an e-variable for $\cP'$, then $E\times E'$ is an e-variable for $\cP\times \cP'$, because 
$$
\E^{\p\times \p'}[E\times E'] =\E^\p[E]\E^{\p'}[E']\le 1
$$
for all $ \p\in \cP$ and $\p'\in \cP'$. The next result further shows that 
$E\times E'$ is   strongly admissible when $E$ and $E'$ are strongly admissible.   \begin{proposition}
\label{prop:strong-adm}
   Let $E$ and $E'$ be strongly admissible e-variables for $\cP$ and $\cP'$, respectively.
Then
 $E\times E'$ is a strongly admissible e-variable for $\cP\times \cP'$.
  \end{proposition}
  
   \begin{proof}[Proof.]
   As we explained above, $E\times E'$ is an e-variable for $\cP\times \cP'$.
 To show its strong admissibility, take any $B$ with $(\cP\times \cP')_B\neq \varnothing$. 
It follows from standard real analysis that there exist $A\in \mathcal F$ and $A'\in \mathcal F'$
such that $A\times A'\subseteq B $ 
and $(\cP\times \cP')_{A\times A'}\neq \varnothing$. 
By strong admissibility of $E$ and $E'$, there exist $ \p\in \cP_A$ and $ \p'\in \cP'_{A'}$ such that 
 $\E^{\p} [E]= \E^{\p'}[E']=1 $.
 This gives $(\p\times \p')(B)>0$ and $ \E^{ \p\times  \p'}[E\times E']=1$, and therefore 
 $E\times E'$ is strongly admissible.
   \end{proof}
   
Combining  Propositions~\ref{prop:adm-e-2} and~\ref{prop:strong-adm}, we obtain that $E\times E'$ is an admissible e-variable for   admissible e-variables $E$ and $E'$
 when $\cP$ and $\cP'$ are finite.

\section{More examples of e-values}
\label{sec:c2-ex}

Continuing from the examples in Chapter~\ref{chap:introduction}, we now provide several simple examples of e-values for different testing problems. Some of these e-values also enjoy a notion of optimality,  introduced in Chapter~\ref{chap:alternative} and analyzed in detail in Chapter~\ref{chap:numeraire}, so these examples will be revisited later in the book.

  In most examples below, we consider a single random observation (statistic) $Z$ in a measurable space $\cZ$ and the sample space is chosen as $(\Omega, \cF)=(\cZ,\sigma(Z))$;  
 thus $\cM_1$ is the set of all probability measures on $\sigma(Z)$.

\subsection{Monotone likelihood ratio}
\label{sec:mlr-evariable}

Let $(\p_{\theta})_{\theta\in\Theta}$ be a parametric family of distributions  
with their
density functions of the test statistic $Z$ given by $(p_\theta)_{\theta\in\Theta}$,
where $\Theta\subseteq \R$ and $\cZ=\R$. 
Assume that monotone likelihood ratio 
holds: for all $\theta<\theta'$, the likelihood ratio
$
{p_{\theta'}(z)}/{p_{\theta}(z)}
$ exists and
is   increasing (or decreasing)    in 
$z$. 
\begin{proposition}
\label{prop:c2-MLR}
  Consider testing  $\cP=\{\p_\theta: \theta \le \theta_0\}$
against $\cQ= \{\p_\theta: \theta \ge  \theta_1\}
$, where $\theta_1>\theta_0$.
Under the monotone likelihood ratio assumption, for any $\theta'>\theta_0$, the random variable
$$
E= \frac{p_{\theta'}(Z)}{p_{\theta_0}(Z)}
$$
is an e-variable for $\cP$.
\end{proposition}
The following short proof uses Chebyshev's association inequality, which states, assuming integrability, $\E^\p[f(X)g(X)]\ge \E^\p[f(X)]\E^\p[g(X)]$ for any increasing functions $f,g$ and random variable $X$.
\begin{proof}[Proof.]
Write $\E^{\theta}$ for $\E^{\p_\theta}$. 
For any  $\theta\le \theta_0<\theta'$, we have 
\begin{align*}
     1=  
\E^{\theta_0} \left[
\frac{ p_{\theta'}(Z)}{ p_{\theta_0}(Z)}
\right]
& =
\E^{\theta} \left[
\frac{ p_{\theta'}(Z)}{ p_{\theta_0}(Z)}
\frac { p_{\theta_0}(Z)}{ p_{\theta}(Z)
}
\right]
\\ &\ge 
\E^{\theta} \left[ \frac{p_{\theta'}(Z)}{ p_{\theta_0}(Z)} \right]
\E^{\theta} \left[
\frac{ p_{\theta_0}(Z)}{ p_{\theta}(Z)
}
\right]
= 
\E^{\theta} \left[ 
\frac{p_{\theta'}(Z)}{ p_{\theta_0}(Z)} 
\right],
\end{align*}
where the inequality follows from Chebyshev's association inequality.
This gives $\E^\theta[E]\le 1$ for all $\theta \le \theta_0$, as required.
\end{proof}

The same conclusion in Proposition~\ref{prop:c2-MLR} holds true for testing 
  $\cP=\{\p_\theta: \theta \ge \theta_0\}$
against $\cQ= \{\p_\theta: \theta \le   \theta_1\}
$
with $\theta_0>\theta_1$; moreover, either inequality (or both) in $\cP$ and $\cQ$ can be replaced by a strict inequality.
The same e-variable can also be used to test $\cP=\{\p_\theta: \theta \le \theta_0\}$
against $\cQ= \{\p_\theta: \theta >  \theta_0\}.
$

The example of normal distributions in Section~\ref{sec:LR-e-variable} satisfies the monotone likelihood ratio assumption.
Therefore, the e-variable in 
\eqref{eq:example-gaussian-e3} for any $\mu>0$
is an e-variable for testing 
$\{\mathrm{N}(\theta,1):\theta\le 0 \}$ against $\{\mathrm{N}(\theta,1):\theta >0\}$.

An important special case is given by one-parameter exponential families, which we discuss further in Section~\ref{sec:exp-fam-numeraire}. There, and in Section~\ref{sec:mlr-numeraire}, we will argue that the above e-variable is actually \emph{log-optimal} for $\q=\p_{\theta_1}$ for any $\theta_1>\theta_0$, simply because $p_{\theta_1}/p_{\theta_0}$ is an e-variable for $\cP$. 

\subsection{Testing symmetry}
\label{sec:c2-symmetry}

Take $\cZ = \R$. We now consider the null hypothesis that $Z$ is symmetrically distributed around $0$, that is,
\[
\cP = \left\{ \p \in \cM_1 \colon Z \text{ and } -Z \text{ have the same distribution under } \p \right\}.
\]
An example in clinical trials  is to set $Z=Y_1-Y_0$, where $Y_1$ is a a post-treatment measurement  and $Y_0$ is pre-treatment measurement,  and we are interested in whether the treatment has any effect. Such an application typically has more than one data point, which will be discussed later, but for now we first consider the case of one summary statistic $Z$.  

Fix any alternative hypothesis $\q$ that admits a Lebesgue density $q$. 
Let $\p^*$ be given by its density
\[
p^*(z) = \frac{1}{2}\left( q(z) + q(-z) \right) \id_{\{q(z) > 0\}},
\]
which is absolutely continuous with respect to $\q$. 
Define     
\begin{equation}
\label{eq:e-test-symmetry}
E = \frac{\d\q}{\d\p^*} = \frac{2q(Z)}{q(Z) + q(-Z)} \id_{\{q(Z) > 0\}}.
\end{equation}
We can check that $E$ is an e-variable for $\cP$ because, using the convention $0/0=0$, 
\begin{align*}
\E^\p[E]   & = \E^\p\left[ \frac{q(Z)}{q(Z) + q(-Z)}  \right]
+\E^\p\left[  \frac{q(Z)}{q(Z) + q(-Z)} \right]
\\&= \E^\p\left[ \frac{  q(Z)}{q(Z) + q(-Z)}  \right]+\E^\p\left[  \frac{q(-Z)}{q(Z) + q(-Z)} \right] 
=1.
\end{align*}
Corollary~\ref{coro:adm-e-1} implies that $E$ is admissible because it is exact.
The optimality of this e-variable against $\q$ is studied in Section~\ref{sec:symmetry}.

The e-variable in \eqref{eq:e-test-symmetry} based on one realization is at most $2$, and in applications we need to rely on multiple data points to accumulate more evidence against the null. 
Suppose that multiple iid data points sampled from $X$ are available, summarized by $Z=(X_1,\dots,X_n)$ with $\cZ=\R^n$. We consider the hypothesis $\cP$ of symmetry, that is, $X$ and $-X$
have the same distribution under  $\p $, and for simplicity, assume further that $X$ is continuously distributed. 
There are many e-variables to test $\cP$, and we give a few.
First, let 
$$
E_\lambda =\prod_{i=1}^n \frac{2\exp(\lambda X_i)}{\exp(\lambda X_i) + \exp(-\lambda X_i)},~~~\lambda \in \R,
$$ 
Corollary~\ref{coro:adm-e-1} implies $E_\lambda$ is admissible since it is exact. 
This is the product of several e-variables in \eqref{eq:e-test-symmetry} with a specifically chosen $q$.
Another class of e-variables is 
$$
E'_\lambda =\prod_{i=1}^n   \exp(\lambda X_i - \lambda^2X_i^2/2) ,~~~\lambda \in \R. 
$$ 
We have $E'_\lambda \le E_\lambda$, because 
$$
\exp(x)+\exp(-x)\le 2\exp(x^2/2)
$$
for all $x\in \R$. To remove the reliance on $\lambda$, we can take $E= \int _\R E_\lambda \d \mu(\lambda) $, where $\mu$ is any distribution on   $\R$. 
Thus $E'_\lambda$ is inadmissible.

 Note that the signs of $X_1,\dots,X_n$ are iid and uniformly distributed on $\{-1,1\}$ under the null hypothesis. 
Therefore, a class of e-variables is given by
$$
E_p = \frac{p^{k} (1-p)^{n-k}}{2^n},~~~p\in(0,1),
$$ 
where $k$ is the number of $1$s in $\{X_1,\dots,X_n\}$. 
Similarly to $E_\lambda$ above, $E_p$ is also a product of several e-variables in \eqref{eq:e-test-symmetry}, but now this class is indexed by $p\in(0,1)$. 

A  different approach is to use the signed-rank e-variables, based on 
Wilcoxon's signed-rank test. 
The idea is to  arrange the magnitudes $|X_i|$ of the observations
in the ascending order and assigning to each its {rank},
which is a number in $[n]$.
Under the null hypothesis, any subset of $[n]$ has the same probability $2^{-n}$ to be the  set of ranks of the positive observations. 
This determines the distribution 
of Wilcoxon's statistic $V_n$
defined as the sum of the ranks of the positive observations.
The class of signed-rank e-variables is
defined by 
$$
E_\beta =\exp(\beta V_n) \prod_{i=1}^n \frac{2}{1+\exp(\beta i)} ,~~~\beta \in \R. 
$$ 
It is a standard exercise to check that $E_\beta$ is an e-variable.

\subsection{Testing the mean and the variance}
\label{sec:c2-bounded}

Let $\cZ \subseteq \R_+$  so that the statistic $Z$ takes nonnegative values. Fix $\mu > 0$ and consider the null hypothesis that the mean of $Z$ is at most $\mu$,
\[
\cP = \{\p \in \cM_1 \colon \E^\p[Z] \le \mu\}.
\] 
The alternative hypothesis $\cQ$ 
is that mean of $Z$ is larger than $\mu$. Two natural e-variables are $Z/\mu$ and the constant one. All convex combinations of these are also e-variables; therefore, 
\begin{equation}
    \label{eq:simple-test-mean}
(1-\lambda)+ \lambda \frac{Z}{\mu},~~~\lambda \in [0,1],
\end{equation}
 is an e-variable for $\cP$. 
This construction works outside $\R_+$ as long as $Z$ is bounded from below, by using $(Z-z)/(\mu-z)$ in place of $Z/\mu$, where $z$ is a lower bound on $Z$.
An optimal choice of $\lambda$ against a specific $\q$ is discussed in Section~\ref{sec:bounded-mean} in the case that $\cZ$ is a bounded interval.
Corollary~\ref{coro:adm-e-1}(b) can be used to deduce that the above e-variable is admissible, despite not being exact.

We can also consider the null hypothesis that the mean of $Z$ is at most $\mu$ and the variance is at most $\sigma^2$,
\[
\cP = \{\p \in \cM_1 \colon \E^\p[Z] \le \mu;~ \var^\p(Z)\le \sigma^2\}.
\] 
In this case, because $\E^\p[Z^2]\le \mu^2+ \sigma^2$ for $\p\in \cP$, each of  
\begin{equation}
    \label{eq:simple-test-variance}
(1-\lambda )+ \lambda  \frac{Z^2}{\mu^2+\sigma^2},~~~\lambda \in [0,1],
\end{equation}
 is an e-variable for $\cP$,  in addition to  those in \eqref{eq:simple-test-mean}, as well as the mixtures of \eqref{eq:simple-test-mean} and \eqref{eq:simple-test-variance}.
 If $\cZ=\R$, meaning that $Z$ may take negative values, then \eqref{eq:simple-test-variance}
 is an e-variable for 
\begin{equation}
  \{\p \in \cM_1 \colon \E^\p[Z] = \mu;~ \var^\p(Z)\le \sigma^2\}.
  \label{eq:c2-mv}
\end{equation}

As a further example, suppose that we observe nonnegative data $Z=(X_1,\dots,X_n)$ with the same marginal distribution $\p$ and possibly dependent components, 
and we  wish to test whether $\p$ has a mean at most $\mu$. Then, $(X_1 +\dots+X_n)/(n\mu)$ is an e-value for any dependence structure of $(X_1,\dots,X_n)$, and does not require making any distributional or independence assumption. 
An e-variable can be constructed similarly for the case of testing the mean and variance.
This example shows the advantage of e-values when it comes to dependent data, and this will be further explored in Chapter~\ref{chap:compound}.

\subsection{Testing sub-Gaussian means}
\label{sec:c2-subG-mean}

Now take $\cZ = \R$ and consider the null hypothesis that the observation $Z$ has a 1-sub-Gaussian distribution with nonpositive mean. This class is a nonparametric generalization of Gaussians with nonpositive mean and unit variance. Formally, we set
\[
\cP = \left\{\p \in \cM_1 \colon \E^\p[ e^{\lambda Z - \lambda^2/2} ] \le 1 \text{ for all } \lambda \in [0,\infty) \right\}.
\]
Here the sub-Gaussian restriction is only placed on the
right tail.  
For $\p\in \cP$ and $\lambda \in \R$, let $M_\p(\lambda) = \E^\p[e^{\lambda Z}]$, which is the moment generating function of $Z$. 
The constraint in $\cP$ means
$$
M_\p(\lambda) \le e^{\lambda^2/2},~~~\lambda \ge0.
$$
A Taylor expansion  in a neighborhood of $\lambda=0$ yields,
\begin{align*} 
M_\p(\lambda)&=1+ \lambda \E^\p[Z]  + \frac{\lambda^2}{2}\E^\p[Z^2]  
+o(\lambda^2);
\\  e^{\lambda^2/2} 
&= 1+ \frac{\lambda^2}{2} + o(\lambda^2). 
\end{align*}
Therefore, for $\lambda\ge0$ small enough, 
\begin{equation}\label{eq:subgaussian-derivative}
\lambda \E^\p[Z]  + \frac{\lambda^2}{2}\E^\p[Z^2]   
\le  \frac{\lambda^2}{2},
\end{equation}
and this implies $\E^\p[Z]\le 0$.



Since the standard Gaussian distribution (zero mean, unit variance) satisfy the above with equality (inferred via their moment generating function) for every real $\lambda$, the constraint expands the class to account for two aspects: the mean can be negative, and the right tails of the random variable only have to be lighter than that of a unit-variance Gaussian. It is well-known that any bounded distribution on $[-1,1]$ with zero mean is 1-sub-Gaussian; so the class $\cP$ clearly does not have any common reference measure, containing discrete distributions with uncountable support as well as continuous unbounded distributions.

A class of e-variables for $\cP$ is given by:  
\[ 
\{e^{\lambda Z - \lambda^2/2} \colon \lambda \in [0,\infty)\}, 
\]
and certainly  any convex combinations (mixtures) of its elements are also e-variables.
The determination of the optimal e-variable against a specific $\q$ is studied in Section~\ref{sec:subG-mean}.
The determination of {all possible} e-variables for $\cP$ 
is a nontrivial task. 

We also consider the closely related null hypothesis that the observation $Z$ has a 1-sub-Gaussian distribution with zero (instead of nonpositive) mean:
\[
\cP^0 = \left\{\p \in \cM_1 \colon \E^\p[ e^{\lambda Z - \lambda^2/2} ] \le 1 \text{ for all } \lambda \in \R \right\}.
\]
Here, the sub-Gaussian restriction is   placed on both the left and the right tails.  
Every $\p \in \cP^0$ has zero mean, 
because now \eqref{eq:subgaussian-derivative} holds for both positive and negative $\lambda$. 
Again, a natural class of e-variables for $\cP^0$ is   \[
\{e^{\lambda Z - \lambda^2/2} \colon \lambda \in \R \},
\]
and certainly  any convex combinations (mixtures) of its elements are also e-variables.

It turns out that every $\p \in \cP^0$ also has variance bounded by one, and this follows from \eqref{eq:subgaussian-derivative}  
with $\E^\p[Z]=0$.
Thus, $Z^2$ is also an e-variable, as is any convex combination with the aforementioned class. 
Indeed, $\cP^0$ is a subset of the hypothesis in \eqref{eq:c2-mv} with $\mu=0$ and $\sigma=1$, and also a subset of the one-sided sub-Gaussian class $\cP$.
However, $Z^2$ is not an e-variable for $\cP$, because those could potentially have very heavy left tails, meaning that their variance is not necessarily bounded.

We note in passing that both the above classes of e-variables are admissible, but the rather technical proof is omitted.

\subsection{Constrained hypotheses}
\label{sec:c2-constrained}

Many of the above examples of $\cP$, and many others we encounter in this book, can be called ``constrained hypotheses'', meaning that they are defined via integral constraints. For such classes, we can sometimes cleanly characterize the set of \emph{all} e-variables (and all admissible ones) in terms of those constraints. 

Let $\cL$ denote the set of all real-valued measurable functions on the underlying probability space. A constraint set is any nonempty set of functions $\Psi \subseteq \cL$.  Define the null hypothesis generated by $\Psi$ as
\begin{equation}
\label{eq:constrained-hyp}
\cP = \left\{\p \in \cM_1 : \int |f|\d\p < \infty \text{ and } \int f \d\p \leq 0 \text{ for all }f \in \Psi \right\}.
\end{equation}
We say that a property holds $\cP$-quasi surely if it holds everywhere outside of a set $A$ such that $\sup_{\p\in\cP}\p(A)=0$. 

We start with the case where $\Psi = \{g_1,\dots,g_K\}$ is finite, as in the case of testing a mean, or mean and variance. One can show that the set of all e-variables $\fE$ for $\cP$ is given by all $[0,\infty]$-valued measurable functions which are $\cP$-quasi surely dominated by $1+\sum_{k =1}^K \pi_k g_k$ for some $\pi \in \R_+^K$. As a corollary, every admissible e-variable for $\cP$ is $\cP$-quasi surely equal to $1 + \sum_{k =1}^K\pi_k g_k$ for some $\pi$ in the set
\[
\Pi^\Psi := \left\{\pi \in \R^K_+ : 1 + \sum_{k =1}^K\pi_k g_k \geq 0, ~ \cP\text{-quasi surely} \right\}.
\]
Clearly, the zero vector $\mathbf 0$ is in $\Pi^\Psi$. But for some choices of $\Psi$, $\Pi^\Psi$ may be equal to $\{\mathbf 0\}$, in which case all  e-variables for $\cP$ are dominated by the constant $1$. 

In settings where $\Psi$ is infinite, such as the sub-Gaussian example in the previous subsection, the corresponding statement is more complex.
Define
\[
\mathcal C = 
\{f \in \cL : f \leq g~ \cP\text{-quasi surely  for some } g \in \text{cone}(\Psi) \},
\]
where $\text{cone}(\Psi)$ is the conic hull of $\Psi$. Then, the set of all e-variables $\fE$ for $\cP$ is given by all $[0,\infty]$-valued measurable functions that are $\cP$-quasi surely equal to $1+f$ for some $f \in \overline{\mathcal C}$, where $\overline{\mathcal C}$ is the weak closure of $\mathcal C$ (which is the closure with respect to a particular weak topology that is omitted here for brevity).

The constraints in \eqref{eq:constrained-hyp} are formulated via integrals of $f\in \Psi$ over $\p$. 
In some applications, the constraints are formulated via other functionals, such as risk measures, and some of them are studied in Chapter~\ref{chap:risk-measures}.

\subsection{Testing likelihood ratio bounds}
\label{sec:c2-LRB}

Let $\mathcal Z $ be a general sample space.
Let     $\p_0$ be a fixed probability measure and $\gamma \in [1,\infty)$, and consider the null hypothesis
$$
\mathcal P =\{\p\in \cM_1: \d \p/\d \p_0 \le \gamma\}
$$
against an alternative $\q \ll \p_0$.
To interpret, 
  testing $\mathcal P$ means to test whether   the true data-generating distribution is close to a given benchmark distribution  $\p_0$ in terms of likelihood ratio being bounded by $\gamma$. 
The case $\gamma=1$ corresponds to the simple null hypothesis $\{\p_0\}$.


Clearly, the likelihood ratio $Z^*=\d \q /\d \p_0$ is not an e-variable for $\cP$ unless $\gamma=1$, but if we modify it to 
$
E=Z^*/\gamma
$, then $\E^\p[E]\le \gamma \E^{\p_0}[Z^*/\gamma] \le 1$. 
Hence, $E=Z^*/\gamma$ is an e-variable for $\cP$. We can actually improve this e-variable, as described below. 

For any random variable $Z$,
denote by  $q_t(Z)=\inf \{x\in \R: \p_0(Z\le x)\ge 1-t\}$ for $t\in (0,1)$; that is, $t\mapsto q_{1-t}(Z)$ is the left quantile function of $Z$ under $\p_0$.   
Let $E^*$ be a transform of the likelihood ratio given by
 $$
E^*= \frac{Z^*\vee z_0}{\gamma},
 $$
 where $z_0\ge 0$ is the largest constant such that 
     \begin{align}\label{eq:z0}
  \int_{0}^{1/\gamma} (q_t(Z^*)  \vee z_0) \d t =1.
  \end{align}
  Clearly $E^*\ge E$.
Let us verify that $E^*$ is an e-variable for $\mathcal P$. This follows by noting that
$$
\sup_{\p\in \mathcal P}\E^\p[E^*] = 
\gamma \int_{0}^{1/\gamma} q_t(E^*)\d t =
\int_0^{1/\gamma} (q_t(Z^*)  \vee z_0) \d t =1,
$$
where the first equality is a classic result on risk measures; see  \eqref{eq:ES-rep} in Chapter~\ref{chap:risk-measures}. 
The optimality of the e-variable $E^*$ is proved in Section~\ref{sec:c5-ex}.

\subsection{Testing exchangeability}\label{subsec:c2-numeraire-exch}

Let $\cZ$ be arbitrary, and suppose $Z = (Z_1,\dots,Z_n)$.
Denote by $\mathcal S_n$ the set of all permutations of the set $[n]$.
 For any $\sigma\in \mathcal S_n$, let $Z_\sigma$ denote the permuted vector $(Z_{\sigma(1)},\dots,Z_{\sigma(n)})$. Consider testing the null hypothesis of exchangeability $$\cP =\{\p: Z \laweq Z_\sigma \text{ under $\p$ for all permutations } \sigma\},$$ where $\laweq$ denotes equality of distribution.
Let $\q$ be an alternative distribution of $Z$ with density $q$, which is not exchangeable, i.e., not in $\cP$.

Let $\pi$ be uniformly distributed on  $\mathcal S_n$ independent of $Z$, and let
$$
E =\frac{q(Z)}{q(Z_\pi)}  = \frac{q(Z)}{   \frac1{n!}\sum_{\sigma\in \mathcal S_n} q(Z_\sigma)}.
$$ By exchangeablity, $\E^\p[E]=\E^\p[q(Z)/q(Z_\pi)]=1$ for all $\p\in \cP$, so $E$ is an exact e-variable, which is admissible by Corollary~\ref{coro:adm-e-1}.
The e-variable $E$ has the form of the soft-rank e-value in Section~\ref{sec:c1-ex2}, where we use $r=0$ 
and $R_b=L_b= q(Z_\sigma)$ for $\sigma\in\mathcal S_n$.
The optimality of the e-variable $E$ is proved in Section~\ref{subsec:numeraire-exch}.

\section{Informal grids for realized e-values}
\label{sec:c2-jeffrey}
 
In testing scientific hypotheses, thresholds for p-values are often chosen as $0.01$ or $0.05$ which correspond  to type-I errors controlled at these levels. The e-to-p calibrator $e\mapsto \min( 1/e,1)$ implies that thresholds of $100$ and $20$ for e-values also have the above type-I error control. However, it is not recommended in practice to directly use these thresholds as  the conversion $e\mapsto \min( 1/e,1)$ is typically wasteful. 
In Section~\ref{sec:betting-interpretation}, and also in Section~\ref{sec:2-posthoc} later, we see that for $a\in (0,1)$, an e-value of $1/a$ carries more information than a p-value of $a$. 
  
In order to   judge how significant results of testing using e-values are, the type-I error, based on which p-values are defined, may not be the desirable metric.
Although there are no universally agreed   thresholds to use for e-values, 
 the rule of thumb of Jeffreys, originally designed for likelihood ratios,  may offer some insight, as e-values are generalizations of likelihood ratios (for example, recall the result presented in Section~\ref{sec:c2-likelihood}). We summarize this rule of thumb in Table~\ref{tab:tab1}. 
 For instance, if $e>3.16$, we recommend to associate it with ``substantial'' (similar to $p<0.05$),
 and  if $e>10$, we recommend to associate it with ``strong'' (similar to $p<0.01$). 
One should keep in mind that p-values and e-values do not have a one-to-one correspondence with each other, and the interpretation of the informal grid for e-values does not rely on converting into p-values.

Because Jeffreys focused on likelihood ratios (or Bayes factors), the above rule of thumb applies well to e-values that are obtained ``naturally'', such as likelihood ratios or their generalizations in the subsequent chapters. At the other extreme, for all-or-nothing e-variables (see Definition~\ref{def:all-or-noth}), such as $E=\id_{\{P\le \alpha \}}/\alpha$ for some p-variable $P$, one should not refer to the grids in Table~\ref{tab:tab1}, because for such an e-variable, $E=10$ is equivalent to $P\le 0.1$.
Finally, Section~\ref{sec:Bayes-factor} draws deeper connections between e-values and Bayes factors, but it builds on some more advanced material.
 
 \begin{table}[ht!]
 \begin{center} 
 \begin{tabular}{lll}
 e-value& level of evidence & Shafer's p-value \\  \hline  
 $[0,1)$ & null hypothesis is supported &$  (0.25,1]$ \\ 
 
$ (1,3.16) $ & no more than a bare mention  &  $(0.0577,0.25)  $\\   
$  (3.16, 10)  $ & substantial  & $  (8.3  {\times}  10^{-3}, 0.0577) $\\  
$  (10, 31.6)  $ & strong & $  (9.4  {\times} 10^{-4}, 8.3  {\times} 10^{-3} )$ \\  
$ (31.6, 100 ) $ & very strong & $  (9.8 {\times} 10^{-5},  9.4  {\times} 10^{-4}) $ \\  
 $ (100,\infty] $ &decisive  & $[0, 9.8 {\times} 10^{-5}     )$ \\  
  \end{tabular}
 \caption{Applying Jeffreys's rule of thumb for likelihood ratios to e-values. For comparison, we also reported Shafer's p-value, which corresponds to the range of $p$ via $e= p^{-1/2}-1$ in \eqref{eq:cal3}. 
 The boundary values can be put in either of the two adjacent  categories.}
 \label{tab:tab1}
 \end{center}
 \end{table} 

 It is important to note that in general, p-values and e-values need not even ``point in the same direction''. We can easily have situations where the p-value appears to provide strong evidence against the null, but the e-value supports the null, or vice versa.
 Let us elaborate with a simple example. Consider testing $\mu = 0$ against $\mu =  \delta$ (for some $\delta >0$) when we see data $Z \sim N(\mu,1)$. Since the classic $z$-test is uniformly most powerful in this setting, we would get the same p-value $P=1-\Phi(Z)$, regardless of the alternative mean $\delta$, while the (arguably natural) likelihood ratio e-value $E=\exp(\delta Z-\delta^2/2)$ depends strongly on $\delta$. 
 For instance, with $\delta=5$ and $X=2$, the p-value $P\approx 0.023$ (suggesting evidence against the null), but
 the e-value $E\approx 0.08$ (suggesting evidence for the null as a likelihood ratio).  
 Thus the realized p-value and e-value could ``point in different directions'' depending on $\delta$ and the realized data. 
 This is another reason that we do not like to in general compare e-value and p-value based methods: while calibration can often preserve the direction (large/small e-values yield small/large p-values and vice versa), natural p-values and e-values that we may directly construct from the data need not have any natural concordance with each other, and so the resulting methods are not always directly comparable.

\section*{Bibliographical note}

This chapter largely resembles results and concepts in various important papers on e-values, but also it contains many results that are derived for the book to make a coherent piece.  

The example in Section~\ref{sec:c2-markov} is based on \cite{vovk2023confidence}. The concept of a calibrator in Section~\ref{sec:2-cali} dates back, at least, to~\cite{Vovk93}; the admissibility results are in~\cite{Vovk/Wang:2021}. \cite{bickel2025small} discussed p-to-e calibrators that depend on the sample size based on Bayes factors. 
Section~\ref{sec:2-rand} is based on~\cite{ignatiadis2024values} and~\cite{ramdas2023randomized}.

The definition of admissibility in Section~\ref{sec:c2-admissible} is derived from~\cite{ramdas2020admissible}, which studies it more generally in sequential settings. Some sufficient conditions in Corollary~\ref{coro:adm-e-1} also follow from the aforementioned paper. 

Among the examples in Section~\ref{sec:c2-ex},  
exponential families and monotone likelihood ratios were studied by~\cite{GrunwaldHK19, grunwald2024optimal}.
The testing symmetry example was first studied by~\cite{ramdas2020admissible}, with follow-up work in \cite{koning2023post,vovk2024nonparametric,larsson2024numeraire}. 
The bounded mean example is a variant of one studied in \cite{waudby2020estimating}, revisited by~\cite{larsson2024numeraire}. 
The mean-and-variance example was studied by \cite{fan2024testing}. 
Sub-Gaussian e-values were studied by~\cite{howard2020time,ramdas2020admissible,larsson2024numeraire}. 
The example on constrained hypotheses was studied by~\cite{clerico2024optimal,larsson2025variables}; our presentation is based on the latter.
The exchangeability example was studied in~\cite{vovk2021testing} and~\cite{ramdas2022testing} in a sequential context. In the batch context, the soft-rank e-value was proposed in the first arXiv preprint of~\cite{wang2022false} (but omitted in the journal version), and revisited by~\cite{koning2023post}. 

The grids in Section~\ref{sec:c2-jeffrey} can be found in \cite[Appendix B]{Jeffreys61}, with a variant proposed in~\cite{kass1995bayes}.

\chapter{Efficiency: E-values under the alternative hypothesis}


\label{chap:alternative}

In Chapter~\ref{chap:markov} we focused only  on the properties of the e-values under the null hypothesis. In this chapter, we describe its desired properties under the alternative hypothesis. This requires the notions of powered e-values, e-power, and log-optimality.

\section{Powered tests, e-values, and p-values}
\label{sec:nontrivial-test}

We first introduce the concept of power for tests and e-values. 
The power of a  binary  test 
$\phi$ against $\q$ is defined as the  probability of rejecting the null hypothesis under the alternative $\q$. 
This concept can be naturally generalized to non-binary tests as we present below.

\begin{definition}[Powered tests and e-values]
\label{def:nontrivial}
    The \emph{power} of a test $\phi$ against an alternative $\q$ is simply $\E^{\q}[\phi]$. The \emph{type-II error (rate)} of $\phi$ equals one minus its power. 
    When testing a null hypothesis $\cP$, a  test $\phi$ is  \emph{powered}   at level $\alpha\in (0,1)$ against $\cQ$ if 
    $$
    \mbox{$\E^{\p}[\phi] \leq \alpha$ for every $\p \in \cP$  and $\E^\q[\phi] > \alpha$ for every $\q \in \cQ$,}
    $$ 
    and it is \emph{uniformly} powered if further $\inf_{\q \in \cQ} \E^\q[\phi] > \alpha$. 
    An e-variable $E$ is  \emph{powered} against $\cQ$ if $\E^\q[E] > 1$ for every $\q \in \cQ$, and it is  \emph{uniformly} powered against $\cQ$ if $\inf_{\q\in \cQ}\E^\q[E] > 1$.
\end{definition}

\begin{remark}
A  closely related concept is  that of level-$\alpha$ \emph{unbiased} tests, i.e., tests $\phi$ satisfying $\E^{\p}[\phi] \leq \alpha$ for every $\p \in \cP$  and $\E^\q[\phi] \ge  \alpha$ for every $\q \in \cQ$. 
For the purpose of this chapter, 
our definition of powered tests is needed, and it excludes  trivial cases such as the constant test $\phi=\alpha$.
\end{remark}

When we mention a powered test $\phi$ without specifying the level, it should be understood as $\alpha=\sup_{\p\in \cP} \E^\p[\phi]$.
  Tests that are not powered can still be practically useful.
For instance,  the existence of
powered tests and e-variables requires that the alternative $\cQ$ does not contain a ``boundary point'' of the null $\cP$, as shown by Example~\ref{ex:nontrivial} below. Nevertheless, although powered tests do not exist for testing $H_0 : \mu< 0$ against $H_1 : \mu \ge  0$ in that example, any powered test for testing 
$H _0 $ against $H'_1 : \mu > 0$ can still be useful to test  $H_0$ against $H_1$ in  practical applications.

\begin{example}
\label{ex:nontrivial}
Suppose that iid data $X_1,\dots,X_n$ are sampled from $N(\mu,1)$ with unknown $\mu$. 
If we test $H_0 : \mu\le 0$ against $H_1 : \mu > 0$, then
$\exp (X_1+\dots+X_n-n/2)$
is a powered e-variable (but not uniformly e-powered).
This e-variable is uniformly powered against  $H_1 : \mu \ge  \epsilon$ for any $\epsilon>0$.

If we test $H_0 : \mu< 0$ against $H_1 : \mu \ge 0$, then 
there does not exist a powered e-variable (or test) because it is impossible to have power against the boundary point $0$. 
\end{example}
\begin{proof}[Proof of the claim in Example~\ref{ex:nontrivial}]  Let $p_m$ be the probability density function of $X^n=(X_1,\dots,X_n)$ sampled from $N(\mu,1)$ with $\mu=-1/m$, and let $p_0$ be that of the case $\mu=0$. For any function $f:\R^n \to \R$, we have $f(\mathbf x)p_m(
\mathbf x) \to f(\mathbf x)p_0(\mathbf x)$ pointwise in $\mathbf x$ as $m\to\infty$. So Fatou's lemma yields that $\int f(\mathbf x)p_0(\mathbf x) \d \mathbf x \leq \liminf_m \int f(\mathbf x)p_m(
\mathbf x) \d \mathbf x $. When $f(X^n)$ is chosen to be an e-variable (resp.~level-$\alpha$ test) for $H_0: \mu < 0$, the right-hand side equals $1$ (resp.~$\alpha$), and the inequality above means that $f$ is also an e-variable (resp.~level-$\alpha$ test) for the boundary point of $\mu=0$. Thus, every e-variable (resp.~level-$\alpha$ test) for $H_0: \mu<0$ is also an e-variable (resp.~level-$\alpha$ test) for $H_0: \mu \leq 0$, or in other words, it is powerless against the point $\mu=0$. 
However, one does exist if we omit $0$ from the alternative --- the aforementioned e-variable is again powered for $H_0 : \mu < 0$ against $H_1 : \mu > 0$.
\end{proof}

\begin{remark}
  Example~\ref{ex:nontrivial} and its proof  demonstrate a more general phenomenon. When testing a parametric model (like the one in the example),  if the null hypothesis is an open subset of the parameter space, the alternative is its complement, and densities exist and are continuous in the parameter, then powered tests and e-variables do not exist.
\end{remark}

Next, we discuss the natural connection between e-values and tests. 
Given any level-$\alpha$ binary test, one can always reproduce its decision by combining an ``all-or-nothing e-value'' with Markov's inequality, as we explain (and generalize) below. 

\begin{definition}[All-or-nothing e-variable]

\label{def:all-or-noth}
    An e-variable $E$ is called an \emph{all-or-nothing e-variable} if  $E = \id_A / c$ for some event $A$ and $c>0$. 
\end{definition}

All-or-nothing e-variables take on only two values: zero or $1/c$ for some $c$. Further, it is clear that $A$ must satisfy  $\p(A) \leq c$ for all $\p \in \cP$. For such e-variables, Markov's inequality at level $c$ holds with equality, so there is no loss in translating between e-values, tests, and p-values.


\begin{fact}[Equivalence between tests and bounded e-variables]
\label{fact:duality}

For any $\alpha \in(0,1)$,
there is a one-to-one correspondence between $[0,1/\alpha$]-valued e-variables and level-$\alpha$ tests. Specifically, for any $\alpha\in (0,1]$, we can use the relationship
$$
E =   \phi/\alpha
$$
to move between  a level-$\alpha$ test $\phi$ and a corresponding e-variable $E$. 

If $\phi$ is binary, then $E$ is  an all-or-nothing e-variable, which takes values in $\{0,1/\alpha\}$. 
In particular, for any  binary test $\phi$ that rejects the null hypothesis 
when a realized p-value is smaller than  or equal to $\alpha$, $\phi$ can be expressed in terms of the all-or-nothing e-variable  $\id_{\{P \leq \alpha\}}/\alpha$.

Thus, there exists a powered level-$\alpha$ binary test  if and only if there exists a powered  e-variable taking values in $\{0,1/\alpha\}$. 
We will expand this observation to include p-variables in Theorem~\ref{th:existence}. 
\end{fact}

As a consequence of Fact~\ref{fact:duality}, 
all binary tests based on p-values --- routinely used in the sciences --- can be  replicated using e-values. A non-binary test $\phi$ is usually provided with the interpretation that it outputs a rejection with probability $\phi$, realized as a number in $[0,1]$. In other words, one can convert a non-binary test $\phi$ to a binary test $\phi'$ as $\phi' = \id_{\{U \leq \phi\}}$, where $U$ is an independent uniform random variable on $[0,1]$. This is closely related to the randomized Markov's inequality discussed in Section~\ref{sec:2-rand}. 

An e-variable $E$ taking values in $[0,1/\alpha]$ is associated with a level-$\alpha$ test $\phi = \alpha E$; in fact, it also yields
 a level-$\beta$ test for any $\beta<\alpha$. A single e-variable yields a family of level-$\alpha$ tests (indexed by $\alpha$), but in the other direction, every level-$\alpha$ test yields a different e-variable.


The above discussion highlights an important conceptual point --- e-values suffice for testing: If a problem is testable (i.e., a powered test exists), then it is testable with e-values. 
This also applies to a notion of powered p-variables, although it is conceptually slightly more complicated, as we define below. 

\begin{definition}[Powered p-values]
For $\beta\in (0,1]$, 
a p-variable $P$ is \emph{powered on $[0,\beta]$} against $\cQ$ 
if  for every $\q\in \cQ$, we have $\q(P\le t) \ge t$ for all $t\in [0,\beta]$ and  strict inequality holds for some $t\in [0,\beta]$. We say $P$ is \emph{powered} if the above holds for $\beta=1$; this means that $P$ is strictly stochastically smaller than a uniform random variable on $[0,1]$ under each $\q\in \cQ$. 
\end{definition}
For practical purposes, it is more than sufficient to consider $\beta=0.25$, for example, as rejections via p-values are always done at small levels. 
In the next section, we will connect powered tests, powered e-variables,
and powered p-variables.

\section{A unified existence result}
\label{sec:exist-results}

We sometimes make the following  simple assumption:
\begin{equation}    \label{assm:U}
       \begin{aligned} &\mbox{There exists  $U\lawis \mathrm{U}[0,1]$ under every $\p \in \cP$ and $\q\in\cQ$,}\\ &\quad \mbox{independent of any other random variables considered.} 
       \end{aligned}\tag{U}
\end{equation}
 This assumption is very weak and is typically satisfied in all statistical applications.

\begin{remark}
To satisfy Assumption \eqref{assm:U}, 
it may be helpful to consider a lifted sample space $\Omega\times [0,1]$ from $\Omega$ and product
  probability measures $\p\times \lambda  $ and $\q\times \lambda  $ from $\p$ and $\q$, where $\lambda$ is the Lebesgue measure. 
Then, choose $U$ with $U(\omega_1,\omega_2)=\omega_2$,
and extend random variables $X$ in the original space to the mappings
$ (\omega_1,\omega_2)\mapsto X(\omega_1)$, so that $U$ is independent of the other random variables that we consider. 
\end{remark}

With this assumption, we can formally connect the existence conditions of powered e-variables, powered p-variables, and powered tests.

\begin{theorem}
\label{th:existence}
Suppose that Assumption~\eqref{assm:U} holds. 
For testing $\cP$ against $\cQ$ (without restriction), the following are equivalent:
\begin{enumerate}[label=(\roman*)]   \item  there exists  a   bounded powered  e-variable;
   \item there exists a powered level-$\alpha$ test  for some (or every) $\alpha\in (0,1)$;   
       \item  there exists  a powered   level-$\alpha$ binary test  for some (or every) $\alpha\in (0,1)$;   
    \item there exists a  powered on $[0,\beta]$ p-variable  for some (or every) $\beta\in (0,1)$;
     \item there exists a   powered p-variable. 
    \end{enumerate}
    In the directions (i)$\Leftrightarrow $(ii),
(iii)$\Rightarrow$(i)
and  (v)$\Rightarrow$(i), we do not need \eqref{assm:U}.
    \end{theorem}    \begin{proof}[Proof.]

In items (ii)--(iv), there is a ``some'' version and an ``every'' version.
When we prove such an item as a consequent, we use the ``every'' version. 
When we take such an item as a premise, we use the ``some'' version.

    The equivalence (i)$\Leftrightarrow$(ii)  follows from Fact~\ref{fact:duality}, noting that from a bounded powered e-variable $E$ we can construct 
    another powered e-variable $\lambda E+(1-\lambda) $ for some $\lambda \in (0,1)$ that  has a smaller upper bound.
    The direction (iii)$\Rightarrow$(i) is also straightforward from Fact~\ref{fact:duality}.

    Next, we will prove 
    in the route 
    that (i)$\Rightarrow$(iv)$\Rightarrow$(v)$\Rightarrow$(i)$\Rightarrow$(iii). 
 
     (i)$\Rightarrow$(iv)
     and  (i)$\Rightarrow$(iii): Let $E_0$ be a bounded powered  e-variable and take $\beta\in (0,1)$. 
    Since $E_0$ is bounded, 
    there exists $\lambda \in (0,1)$ such that 
    $E: =\lambda E_0+(1-\lambda) $
     is bounded above by $1/\beta$.
Clearly, $E $ is an e-variable, and $\E^\q[E]=\E^\q[\lambda E_0+(1-\lambda) ] >1$ for every $\q\in \cQ$.
Let $P=U/E$, where $U$ is in \eqref{assm:U} and independent of $E$. 
By Corollary~\ref{coro:ip-comb}, $P$ is a p-variable. 
To see that it is powered on $[0,\beta]$, we can check that, for $t\in (0,\beta]$ and $\q\in \cQ$,
$$\q(P\le t) =
\q(U\le t E)
= \E^\q[(tE)\wedge 1 ] = \E^\q[tE]>t.
$$
Thus, (iv) holds. 
To see (iii), a powered level-$\beta$ binary test is $\id_{\{P\le \beta\}}$ for every $\beta \in(0,1)$.

  (iv)$\Rightarrow$(v): 
  Let $P$ be a  powered on $[0,\beta]$ p-variable,
  $U$ be in \eqref{assm:U} independent of $P$, and $P'$ be given by 
  $$
  P'= P\id_{\{P\le \beta \}} + (\beta+ (1-\beta) U)  \id_{\{P>\beta\}}. 
  $$
  Let $\s\in \cP\cup\cQ$. For $t\in (0,\beta]$, we have 
$\s(P'\le t ) = \s(P\le t)  $.
For $t\in (\beta,1)$, we have  
  \begin{align*}
    \s(P'\le t)& = \s(P\le t\wedge \beta) + 
  \s(\beta+ (1-\beta) U\le t) \s(P>\beta) 
 \\&  = \s(P\le \beta) + \frac{t-\beta}{1-\beta}(1-\s(P\le \beta)) 
  \\&  = \s(P\le \beta)\left(1-\frac{t-\beta}{1-\beta}\right) + \frac{t-\beta}{1-\beta}.  
  \end{align*}
Note that if $\s(P\le \beta)=\beta $ then 
$\s(P'\le t)=t$ for $t\ge \beta$.   
  Putting these together, we get that  for $t\in (0,1)$,
   $\s(P'\le t)\le t$  if $\s\in \cP$,
 and $\s(P'\le t)\ge t$  if $\s\in \cQ$.
 Moreover, 
 $\s(P'\le t)> t$ for some $t\in (0,\beta]$ if $\s\in \cQ$.
Therefore, $P'$ is a powered p-variable for $\cP$ against $\cQ$.

  (v)$\Rightarrow$(i):
    Let $P$ be a  
    powered p-variable,
    and  
    $$
    f:\R_+\to \R_+, ~~~f(x)= {(2-2x)_+}.
    $$
    Note that $\int_0^1 f (x)\d x = 1$ and $f$ is a calibrator in Section~\ref{sec:2-cali}. 
Define
    $E=f(P)$, which is an e-variable for $\cP$.  
    Moreover, for each $\q\in \cQ$, we have $\E^\q[E]>1$ by the definition of the  powered p-variable $P$ and the fact that $f$ is strictly decreasing on $[0,1]$.
    Therefore, $E$ is bounded and powered against $\cQ$.
\end{proof} 
 
Condition~\eqref{assm:U} is needed for implications involving constructing powered p-variables, because without \eqref{assm:U}, e.g., in the case of discrete sample spaces (see Example~\ref{ex:discrete} below), one may not be able to construct a p-variable $P$ that satisfies 
$\p(P\le \alpha)\le \alpha\le \q(P\le \alpha)$ for all $\alpha \in(0,1)$ and all $\p,\q\in \cQ$ even when powered e-variables exist. 
Basically, finite sample spaces are not suitable for constructing powered p-variables, but there could still be potentially useful (discrete) p-variables  such as the one for the  Bernoulli distributions in Section~\ref{sec:LR-e-variable}. 
\begin{example}
\label{ex:discrete}
Consider 
a finite sample space $\Omega$
and 
testing $\p $
against $\q $
with $\p\ne \q$ and $\p(\{\omega\})>0$ for all $\omega\in \Omega$ (as in the Bernoulli setting in Section~\ref{sec:LR-e-variable}). A powered e-variable is the likelihood ratio $\q/\p$.
Any p-variable $P$ defined on $\Omega$ takes finitely many positive values. Therefore, $\q(P\le \alpha)=0$ for $\alpha>0 $ small enough, showing that $P$ cannot be  $[0,\beta]$ powered against $\q$ for any $\beta \in(0,1]$.
\end{example}

Theorem~\ref{th:existence} offers an equivalence between  the existence conditions, but it is silent about how to check such existence for   given $\cP$ and $\cQ$. This issue will be addressed in  
Section~\ref{sec:existence}.




\section{E-power and powered e-values}
\label{sec:e-power}

Given this high-level equivalence  between tests and bounded e-values discussed in Fact~\ref{fact:duality}, one may be tempted to ask: What is new here? Are we just simply working with tests in disguise? The answer is this: While all-or-nothing e-values are interesting conceptually in order to relate them to tests and p-values, these are not typically how one would construct e-values in practice. Recall the example from Section~\ref{sec:LR-e-variable} on likelihood ratios. They will commonly have support on the entire positive real line (not just on two points), and never take on the value zero if $\p$ and $\q$ are absolutely continuous with respect to each other.  
To explain this, the following concept of e-power is useful.  

\bigskip


\begin{definition}[E-power]\label{def:e-power}
The \emph{e-power} of an e-variable $E$ against an alternative $\q$ is $\E^\q[\log E]$, assuming that the expectation is well-defined. 
An e-variable $E$ has  \emph{positive e-power} against $\cQ$ if $\E^\q[\log E] > 0$ for each $\q\in \cQ$. It has \emph{uniformly} positive e-power against $\cQ$ if $\inf_{\q \in \cQ} \E^\q[\log E] > 0$.
\end{definition}

The e-power of an e-variable could equal $\infty$ or $-\infty$ (these still count as being well-defined). 
Indeed, whenever we make statements involving expected logarithms, we require them to be well-defined.
The e-power (higher is better) of all-or-nothing e-variables is $-\infty$.

The concept of e-power is closely related to the Kelly criterion in betting and investment. The logarithm reflects the fact that in many applications, e-values are multiplied together. For a simplest example, consider a sequence of iid e-variables $E_1,E_2,\dots$ with finite e-power and
 the process $M=(M_t)_{t\in \N}$ defined by \begin{align}\label{eq:c3-M}
 M_t=\prod_{s=1}^ tE_s,~~~~t\in \N.
 \end{align} 
 The e-power of $E_1$ against $\q$ equals the \emph{asymptotic growth rate} of $M$, 
that is, the $\q$-almost sure limit of $(\log M_t)/t$  as $t\to\infty$ (see Definition~\ref{def:asym-GR}),
justified by the strong law of large numbers. 
The asymptotic growth rate  will be important in the methodology of testing by betting, formally treated in Chapter~\ref{chap:eprocess}.  
When converting $M_t$ into the level-$\alpha$ test $\id_{\{M_t\ge 1/\alpha\}}$ for $\alpha\in (0,1)$, 
 the role of e-power is clear from  the following proposition.

\begin{proposition}[Relating e-power to power]
\label{prop:e-power-iid}
Let $E_1,E_2,\dots$ be a sequence of e-variables that are iid under $\q$, and $M$ be given by \eqref{eq:c3-M}.
    If $\E^\q[\log E_1]>0$, then 
    $\q(M_t>1/\alpha)\to 1$ as $t\to\infty$ for every $\alpha>0$.
        If $\E^\q[\log E_1]<0$, then 
    $\q(M_t>1/\alpha)\to 0$ as $t\to\infty$ for every $\alpha>0$.
\end{proposition}
The proof of Proposition~\ref{prop:e-power-iid} is a simple application of the strong law of  large numbers.

Note that for each $t\in \N$, $E_1$ has positive e-power against $\q$ is and only if $M_t$ has positive e-power against $\q$. 
Proposition~\ref{prop:e-power-iid}
implies that in the simple iid setting, a positive e-power leads to asymptotic power-$1$ tests, and a negative e-power leads to asymptotic power-$0$ tests.
An axiomatic justification of e-power is presented in Section~\ref{sec:axioms}, which uses an axiom based on an idea that is similar to Proposition~\ref{prop:e-power-iid}. 

\begin{example}
\label{ex:c3-LRE}
For testing $\p=\mathrm{N}(0,1)$ against  $\q=\mathrm{N}(\mu,1)$ with $n$ observations, the likelihood ratio e-value  in \eqref{eq:example-gaussian} is 
 $ 
E_n=  \prod_{i=1}^n \exp(\mu X_i - \mu^2/2)
 $. 
 Its e-power is thus $\E^\q[\log E_n] =  n\E^\q[\mu X_1-\mu^2/2] = n\mu^2/2.$
This e-variable has positive e-power against $\q$ and the e-power grows linearly in $n$.
Indeed, these two properties hold for all likelihood ratio e-values with iid observations. 

Let $\q_\theta=\mathrm{N}(\theta,1) $ for $\theta \in \R$. 
Noting that $\E^{\q_\theta}[\log E_n] = n(\mu\theta -\mu^2/2)$,
the e-variable $E_n$ has positive e-power against $\{\q_\theta:\theta >\nu\}$ for $\nu\ge \mu/2$ (uniformly if $\nu>\mu/2$),  and it does not have positive e-power against any $ \q_\theta$ with $\theta \le \mu/2$.
\end{example}

The condition $\E^\q[\log E] > 0$ in Definition~\ref{def:e-power} is equivalent to saying that the geometric mean of $X$ is larger than one: $\exp\E^\q[\log E]> 1$, which is a stronger constraint than requiring the (arithmetic) mean of $E$ to be larger than 1, that is, $\E^\q[E] > 1$ in Definition~\ref{def:nontrivial} (this can be seen from Jensen's inequality).  Hence, if $E$ has positive e-power, then $E$ is a powered e-variable, but the reverse implication is not true in general. 
These two conditions are closely related. Indeed, for a simple alternative $\q$,
from any powered e-variable  one can easily find an e-variable with positive e-power, through the transform $E\mapsto 1-\lambda+\lambda E $ for some small $\lambda >0$ that depends on $E$ and $\q$.  This is formalized in the following result.

\begin{theorem}
\label{thm:nontrivial}
Let $X\ge 0$ be a nonnegative random variable. 
  For  any $\q\in \cM_1$, \begin{equation}
    \label{eq:nontrivial}
\E^\q[X]>1 \iff \E^\q[\log (1-\lambda + \lambda X)]>0 \mbox{ for some $\lambda \in [0,1]$}.
\end{equation} 
In fact, either condition in~\eqref{eq:nontrivial} implies the existence of a \emph{bounded} nonnegative $Y$ with $\E^\q[Y]>1$ and $\E^\q[\log (1-\lambda + \lambda Y)]>0$.
For a general $\cQ\subseteq \cM_1$,  
if $X$ is bounded or $\cQ$ is finite, then 
\begin{equation}
    \label{eq:nontrivial2}
\begin{aligned}
   & \inf_{\q\in \cQ} \E^\q[X]>1   \\& \iff \inf_{\q\in \cQ} \E^\q[\log (1-\lambda + \lambda X)]>0   \mbox{ for some $\lambda \in [0,1]$}.
\end{aligned}
\end{equation} 
As a consequence,  for a finite $\cQ$, an e-variable with positive e-power against $\cQ$ exists if and only if a powered e-variable against $\cQ$ exists.      
\end{theorem}
\begin{proof}[Proof.]

The backward direction in all statements does not require any conditions, and it follows from a simple application of Jensen's inequality, which yields $$0<\E^\q[\log (1-\lambda + \lambda X)]\le \log \E^\q[ 1-\lambda + \lambda X] ~~\Longrightarrow ~~\E^\q[X]>1.$$

For the forward direction, we first note that if $\cQ$ is finite, then  $\min_{\q\in \cQ} \E^\q[X]>1 $
implies $\min_{\q\in \cQ}  \E^\q[X\wedge K]>1$ for some $K\ge 1$. Hence, we can without loss of generality assume that $X$ is bounded. 
It then remains to prove \eqref{eq:nontrivial2} in case $X$ is upper bounded by some $K>1$.

Write $x_+=x\vee 0$ and $x_-=(-x)\vee 0$ for $x\in \R$ and denote by $\delta=\inf_{\q\in \cQ} \E^\q[X]-1 >0$. 
Next, fix  any $\q\in \cQ$.
Since $\E^\q [(X-1)_+] - \E^\q [(X-1)_-] = \E^\q [X-1] \ge  \delta>0 $ and $\E^\q [(X-1)_-] \in [0,1]$, 
we have
$$
\E^\q [(X-1)_+]  \ge \E^\q [(X-1)_-]   +  \delta  
\ge (1+\delta) \E^\q [(X-1)_-].
$$
Fix $\epsilon\in (0,1)$ such that $(1+\epsilon)/(1-\epsilon)<1+\delta$,
and let 
$$\eta = \left(
\frac{1}{1+\epsilon}
-\frac{1}{(1+\delta)(1-\epsilon)}
\right) \delta>0.$$
Since $\E^\q [(X-1)_+]\ge \delta >0$,  we have  
\begin{align*}
   & \frac{\E^\q[(X - 1)_+ ]}{1+\epsilon} - \frac{\E^\q[(X  - 1)_-]}{1-\epsilon}
 \\  & \ge  
\left(
\frac{1}{1+\epsilon}
-\frac{1}{(1+\delta)(1-\epsilon)}
\right)
\E^\q[(X - 1)_+ ]  \ge \eta >0.
\end{align*}
Note that
$\log (1+x) \ge x/(1+\epsilon)$ for $x\in [0,\epsilon)$
and $\log (1+x) \ge x/(1-\epsilon)$ for $x\in (-\epsilon,0)$, that is,
$$
\log(1+x) \ge \frac{x_+}{1+\epsilon}
-\frac{x_-}{1-\epsilon} \mbox{~~~for $x\in (-\epsilon,\epsilon).$}
$$
Hence, for $\lambda\in (0,\epsilon/K)$, implying $\lambda (X-1)\in (-\epsilon,\epsilon)$, we have  
\begin{align*}
  \E^\q[\log (1 + \lambda (X-1) )]   
\ge   \frac{\E^\q[ \lambda (X-1)_+ ]}{1+\epsilon}
- \frac{\E^\q[ \lambda (X-1)_- ]}{1-\epsilon}
\ge \lambda \eta >0.
\end{align*}   
Note that the choice of $(\lambda,\epsilon,\eta)$ does not depend on  $\q\in \cQ$.
Hence, by taking an infimum over $\q\in \cQ$ we complete the proof. 
\end{proof}

The forward direction in \eqref{eq:nontrivial2} fails if  $\cQ$ is infinite and  $X$ is unbounded. Below,
we give an example of an unbounded e-variable $E$  satisfying $\inf_{\q\in \cQ}\E^\q[E] >  1$ (that is, uniformly powered against $\cQ$), but  $$\inf_{\q\in \cQ} \E^\q[\log (1-\lambda + \lambda E)]< 0$$ for all $\lambda \in (0,1].$  Moreover, this example also illustrates that a bounded powered e-variable may not exist even if a uniformly powered e-variable exists. 

 \begin{example}
    \label{ex:e-power-infinite-q} 
    Let 
    $\Omega= \N_0$,  $\p=(1/2) \delta_{0} +(1/2) \delta_1$ where $\delta_x$ is the point mass at $x$,
    and 
    $\cQ=\{\q_n:n\in \N\}$,
    where $\q_n=(1/n) \delta_n + (1-1/n) \delta_0$. 
    \begin{enumerate}
        \item[(i)] 
      Let $E$ be given by $E(n)=2n$ for $n\in \N_0$.     
    This gives $\E^\p[E]=1$ and  $\E^{\q}[E]=2$ for all $\q\in \cQ$. 
    Thus, $E$ is a uniformly powered e-variable, and it is unbounded.
    For any $\lambda \in (0,1]$, we have $$\E^{\q_n} [\log (1-\lambda +\lambda E)] 
    = \frac 1 n \log\left( (1-\lambda )^{n-1} (1+(2n-1)\lambda) \right). 
    $$
    Since $(1-\lambda )^{n-1} (1+(2n-1)\lambda)<1$ for $n$ large enough, 
    we know that  $\inf_{\q\in \cQ} \E^\q[\log (1-\lambda + \lambda E)]< 0$
    for all $\lambda \in (0,1]$. 

\item [(ii)]
There does not exist any  bounded e-variable that is powered against $\cQ$. To see this, for an e-variable $E$ to be powered against $\q_1$, we have $E(1)>1$. 
For $E$ to be powered against $\q_n$, we need $(1/n) E(n) + (1-1/n) E(0)>1$ for all $n\in \N$.
Since $E$ is bounded, this means $E(0)\ge 1$. 
Together with $E(1)>1$, we have $\E^\p[E]>1$, conflicting the fact that $E$ is an e-variable for $\p$. 
    \end{enumerate}

    Note that here $\p$ does not lie in $\cQ$, or even in its convex hull, but it does lie within its closed convex hull, where the closure is taken with respect to the total variation (TV) distance (see Theorem~\ref{th:kraft} below). In particular, by sending $n\to\infty$, $\q_n$ approximates $\delta_0$ arbitrarily well in TV. Averaging $\delta_0$ with $\q_1$ yields $\p$. 
\end{example}

As hinted above, the e-power $\E^\q[\log E]$ may not always be well-defined, due to the unbounded nature of the logarithm function. 
Moreover, there are some special (perhaps pathological) examples in which the e-power, although well-defined, may be infinite, and does not have a useful interpretation. We give a stylized, unrealistic, example here. 

\begin{example}
    \label{ex:pathological}
  Let $Y$ be a Pareto(1) random variable under $\q$, that is,  $\q(Y>x)=1/x$ for $x\ge 1$ (this distribution has an infinite mean),
and  the distribution of $Y$ under $\p$ be such that $\E^\p[\exp({6-Y/100})]\le 1$. 
Note that for any two distributions $F,G$, we can always find $\p,\q$ such that a random variable $Y$ is distributed as $F$ under $\p$ and $G$ under $\q$. 

The  e-variable $ E$ for $\p$ given by $$E=\exp({6-Y/100})$$  
satisfies $ \E^\q[\log E]=\E^\q[6-Y/100]=-\infty$. 
We can easily compute $$\q(E>100) >  \q(E>\exp(5)) = \q(Y<100)=0.99.$$
Therefore,  $E$ has a high probability under $\q$ to take large numbers (providing very strong evidence against the null),
although its e-power  is $-\infty$. 
If another e-variable $E'$ is identically distributed to $1/E$ under   $\q$,
then $\q(E' > 1) =1/600$ but $\E^\q[\log E']=\infty$. 
Thus,
$E'$ has more e-power than
$E$ (and  any   e-variables with finite e-power), which may be regarded as counter-intuitive.
Nevertheless, if we have an infinite iid sequence $E'_1,E'_2,\dots$ following the same distribution as $E'$ under $\q$, then the product $M_n = \prod_{i =1}^n E'_i$ goes to $\infty$ $\q$-almost surely as $n\to\infty$, due to the law of large numbers applied to $(\log M_n)_{n \geq 1}$, and this is consistent with Kelly's criterion. 
A similar example in the sequential setting will be discussed in Section~\ref{sec:path-ex-power}.
\end{example}

\section{Conditions  
for the existence of powered e-variables}
\label{sec:existence}

In this section, we address the existence of powered e-variables and p-variables by means of two characterization results.  
These   results connect  the existence of p-variables and e-variables  to the structure of the null and alternative hypotheses, 
different from Theorem~\ref{th:existence}, which does not discuss how to check the existence.
In these results, 
let $L$ and $M$ be two positive integers  and we continue to assume  \eqref{assm:U} 
in Section~\ref{sec:exist-results}.

The following result uses the total variation distance between sets of distributions, which we now define. 
For any two distributions $\p,\q$, define 
\[
d_\text{TV}(\p,\q) = \sup_A |\p(A) - \q(A)|,
\]
where the supremum is taken over all measurable sets $A$. 
Then, given two sets of distributions $\cP,\cQ$, we can define
\[
d_\text{TV}(\cP,\cQ) = \inf_{\p\in\cP,\q\in\cQ} d_\text{TV}(\p,\q).
\]



The following result is closely connected to a classic theorem of \cite{kraft1955some}, proved by Le Cam in private communication. The Kraft--Le Cam theorem states that for general $\cP$ and $\cQ$ (not necessarily finite) that have a common reference measure, the following holds: for any $\epsilon \geq 0$, there exists a test $\phi$ such that $\inf_{\q \in \cQ} \E^\q[\phi] - \sup_{\p \in \cP} \E^{\p}[\phi] > \epsilon$ if and only if the total variation distance between the closed convex hulls of $\cP$ and $\cQ$ is greater than $2\epsilon$.

\begin{theorem}
\label{th:kraft}
  For testing $\cP$ against $\cQ$ with a common reference measure,
 the following are equivalent:
\begin{enumerate}[label=(\roman*)] 
    \item there exists a test $\phi$ such that $\inf_{\q \in \cQ} \E^\q[\phi] > \sup_{\p \in \cP} \E^{\p}[\phi] $; 
    \item  there exists a uniformly powered  bounded e-variable  $E$ for $\cP$ against $\cQ$; 
    \item  $d_{\text{TV}}(\conv(\cP), \conv(\cQ)) >0$. 
    \end{enumerate}
\end{theorem}
\begin{proof}[Proof.]
    The equivalence of (i) and (iii) is by the Kraft-Le Cam theorem. Statement (ii) immediately implies (i) because the definition of an e-variable implies that the right-hand side of (i) is bounded by one; since the e-variable is uniformly bounded, one can transform it to a test by dividing by this bound. So it only remains to establish that (i) implies (ii). First note that the right-hand side of (i) must be strictly larger than zero. This is because if $\phi$ equals zero $\p$ almost surely for any $\p \in \cP$, then it also equals zero almost surely for any $\q \in \cQ$ because of the common reference measure assumption, making the strict inequality in (i) impossible. Thus, using the $\phi$ in (i), one can define
    \[
    E = \frac{\phi}{\sup_{\p \in \cP} \E^{\p}[\phi]},
    \]
    for which we see that (ii) holds, completing the proof.
\end{proof}

If there is no common reference measure, then condition (iii) is not equivalent to the earlier ones. In fact, even when $\cQ=\{\q\}$ is simple, $d_{\text{TV}}(\conv(\cP), \q) >0$ does not imply the existence of a powered e-variable. The appropriate condition actually  involves an object called the \emph{bipolar} of $\cP$, and is formally presented in Theorem~\ref{thm:nontrivial-iff-bipolar}.
Separately, note that when $\cP$ and $\cQ$ are finite, condition (iii) above reduces to  $\conv(\q_1,\dots,\q_M) \cap  \conv(\p_1,\dots,\p_L) = \emptyset$.

The existence of exact p-variables and e-variables is technically more advanced, and studied in Section~\ref{sec:existence-exact}, where we obtain a similar result (Theorem~\ref{th:c3-span}) to Theorem~\ref{th:kraft}.

\section{Log-optimality and likelihood ratios}
\label{sec:3-p-vs-q}

Next, we turn to the key concept  of log-optimality for a simple alternative $\{\q\}$, which is defined essentially via the comparison of e-power. 

\begin{definition}[Log-optimality]
\label{def:log-opt}
When comparing two e-variables $E_1$ and $E_2$ for $\cP$  against $\q$, $E_1$ is said to have greater (resp.~equal) e-power than $E_2$ if $\E^\q[\log(E_2/E_1)] \leq 0$ (resp.~$=0$).  
Moreover, $E$ is \emph{log-optimal} if  $\E^\q[\log(E'/E)] \leq 0$ for any other e-variable $E'$.
\end{definition} 

Definitions~\ref{def:e-power}
and~\ref{def:log-opt} are obviously connected.
 If the maximum e-power, defined as
\[
\max_{E\in \fE} \E^\q[\log E]  ,
\]
is finite, then the log-optimal e-variable maximizes the e-power, and the definition of log-optimality reduces to saying $\E^\q[\log E'] \leq \E^\q[\log E]$. The current definition allows the maximum e-power to be infinite, which could occur even in simple cases such as testing a Gaussian against a Cauchy because the tails of the Cauchy are much heavier than that of the Gaussian; see also Example~\ref{ex:pathological}.
Note that  
$\E^\q[\log(E'/E)] \leq 0$
and  $\E^\q[\log(E') ] \le \E^\q[\log(E) ] $  are equivalent when $E$ and $E'$ have finite e-power, but the two conditions are not equivalent in general.

\begin{proposition}
For testing any $\cP$, a log-optimal e-variable against $\q$ always exists and is $\q$-almost surely unique.
\end{proposition}
We do not prove the first part in this book, since the proof is quite technical, and defer the second part to Chapter~\ref{chap:numeraire}.
However, in the case of a simple null hypothesis $\p$ with $\q\ll \p$,  
existence and uniqueness of the log-optimal e-variable can be easily argued, as we will see in Theorem~\ref{th:KL}. 


In this section, we present an important result on the log-optimal e-values for a simple null hypothesis $\cP = \{\p\}$ and a simple alternative $\cQ = \{\q\}$. 
We will assume that $\q \ll \p$.
Let $q$ and $p$ represent the densities of $\q$ and $\p$ with respect to some reference measure $\leb$, by which we can write $({\d\q}/{\d\p})(x) = {q(x)}/{p(x)}$.

The classic Neyman--Pearson lemma for hypothesis testing states that when testing a simple null $\p$ against a simple alternative $\q$, the most powerful level-$\alpha$ test involves thresholding the likelihood ratio. We now present a result that can be seen as the e-analog of the Neyman--Pearson lemma. Importantly, we are not interested in designing tests (which can be obtained via Markov's inequality as discussed in Chapter~\ref{chap:markov}), but simply in designing e-variable that is log-optimal (as e-powerful as possible).

To prepare for the statement, define the Kullback--Leibler divergence between $\q$ and $\p$ as
\begin{equation}
\label{eq:kl-def1}
\kl(\q,\p) :=  \E^\q\left[\log \frac{\d \q}{\d\p} \right]= \int q(x) \log(q(x)/p(x)) \leb (\d x) ,
\end{equation}   
where it is defined as $\infty$ if $\q$ is not absolutely continuous with respect to $\p$.
Gibbs' inequality, or the Shannon--Kolmogorov information inequality, states that the Kullback--Leibler divergence between any two distributions is nonnegative: $\kl(\q,\p) \geq 0$.  Further, equality holds if and only if $\p=\q$.

\begin{theorem}\label{th:KL}
    For testing $\p$ against $\q$ with  $\q\ll \p$,
    the ($\p$-almost surely) unique log-optimal e-variable is the likelihood ratio $\d\q / \d\p$. Thus, the optimal e-power equals $\kl(\q,\p)$, which could be infinite.
\end{theorem}
 
\begin{proof}[Proof.]
Let $\s \ll \p$ be an arbitrary distribution. By Proposition~\ref{prop:admissible-lr}, $E:=\d \s/\d \p$ covers all exact e-variables. 
For our setting of a point null $\p$,
it suffices to consider exact e-variables 
$E'$ in Definition~\ref{def:e-power},
because any non-exact e-variable can be made strictly larger by multiplying it by a constant larger than $1$ (the inverse of its expected value under $\p$).

We now argue that it suffices to consider $\s \ll \q$. 
Indeed, if $\s$ is not absolutely continuous with respect to $\q$,
then we can take its absolutely continuous part $\s^*$, and define an e-variable $E^*:=\d \s^*/\d \p $. 
Note that $E=E^*$ $\q$-almost surely, and thus they have the same e-power, but $E^*$ is not exact, so we do not need to consider either $E$ or $E^*$. 

Thus it suffices to consider e-variables $E = \d\s/\d\p$, where $\s \ll \q \ll \p$.
By direct substitution, 
\begin{align}
\E^\q\left[\log \frac{E}{\d \q/\d \p} \right] &=
\int q(x) \log\left(\frac{s(x)/p(x)}{q(x)/p(x)} \right) \leb (\d x) \\
&= - \int q(x) \log\left(\frac{q(x) }{s(x) } \right) \leb (\d x)   \le 0 ,
\end{align}
where the inequality follows from the fact that $\kl(\q, \s) \geq 0$.
This shows that $\d \q/\d \p$ is log-optimal.
 Uniqueness follows from the fact that the inequality holds  as an equality only when $\s=\q$,
 so $E=\d \q/\d \p$ $\q$-almost surely,
 and $\E^\p[\d\q/\d\p]\le \E^\p[E]\le 1=\E^\p[\d\q/\d\p]$ forces $E$ to be unique $\p$-almost surely.
 This completes the proof of the first claim. 
 The second claim follows by simply plugging in $q(x)/p(x)$ into the definition of e-power. 
\end{proof}
\begin{example}
Continue the setting of
 testing $\p=\mathrm{N}(0,1)$ against $\q=\mathrm{N}(\mu,1)$  with $n$ observations in Example~\ref{ex:c3-LRE}. 
  Consider the e-variable $ 
E_n(\theta) =  \prod_{i=1}^n \exp(\theta X_i - \theta^2/2)
 $ for some $\theta\in \R$ that is not necessarily equal to $\mu$.
Its e-power is $\E^\q[\log E_n(\theta)] =  n\E^\q[\theta X_1-\theta^2/2] = n (\theta\mu-\theta^2/2).$
  Simple algebra yields that $\theta=\mu$ gives the largest e-power, confirming the implication of Theorem~\ref{th:KL}. 
\end{example}

In Chapter~\ref{chap:numeraire},
we present a more general result encompassing Theorem~\ref{th:KL}, which even relaxes the assumption $\q \ll \p$.

\section{Conditional log-optimality and   log-optimal calibrators}

In statistical applications,   e-variables and p-variables are  constructed based on available data, and thus they belong to a $\sigma$-algebra $\mathcal G\subseteq \cF$ generated by the data. 
With this restriction, we need to consider  the \emph{log-optimal e-variable with respect to $\mathcal G$}, that is,  
an e-variable $E$ measurable with respect to $\mathcal G$ that satisfies 
$\E^\q[\log(E'/E)] \leq 0$
for all e-variables  $E'$ that are measurable with respect to $\mathcal G$.
In what follows,  for $\p\in \cM_1$, 
  $\p|_\cG$ is the restriction of $\p$ to $\cG$, that is, $\p|_\cG(A)=\p(A)$ for $A\in \cG$.
 
\begin{proposition}
\label{prop:cond-log-opt}
Let $\mathcal G$  be  a sub-$\sigma$-algebra of $\mathcal F$. If $\q\ll \p$, then 
the log-optimal e-variable for testing $\p$ against $\q$  with respect to $\mathcal G$ exists and is  given  ($\p$-almost surely uniquely) by
\begin{equation}
\label{eq:cond-log-opt}
E=  \frac{\d \q|_\cG}{\d \p|_\cG} 
 =  \E^\p \left[\frac{\d \q}{\d \p}\mid \mathcal G \right].
 \end{equation} 
The first equality in \eqref{eq:cond-log-opt} remains true 
if we assume $\q|_{\cG} \ll\p|_{\cG}$ instead of $\q\ll \p$.

\end{proposition}
\begin{proof}[Proof.]
The first equality in \eqref{eq:cond-log-opt} follows directly from Theorem~\ref{th:KL}, by restricting the problem to the smaller measurable space $(\Omega,\cG)$.
The second equality in \eqref{eq:cond-log-opt}
follows from, for any $A\in \mathcal G$, 
$$
 \int _A  \frac{\d \q}{\d \p}  \d \p 
=\q(A)
= \q|_\cG(A) 
=\int _A  \frac{\d \q|_\cG}{\d \p|_\cG}  \d \p|_\cG =\int _A  \frac{\d \q|_\cG}{\d \p|_\cG}  \d \p  ,
$$
and therefore
$ \d \q|_\cG/{\d \p|_\cG} $
is the conditional expectation of $\d \q/\d\p$ on $\cG.$
\end{proof}

A direct consequence of Proposition~\ref{prop:cond-log-opt} is the log-optimal calibrator. 
Suppose that we know a p-variable $P$ is uniformly distributed on $[0,1]$ under   $\p $;  thus it is exact. 
For a given simple alternative $\q$,   
 we can define the \emph{log-optimal calibrator}, that is, a calibrator $f$ such that 
$$
\E^{\q}[\log( f(P) / g(P))]\le 0
$$
for any other calibrator $g$.  
In other words, $f(P)$ is the log-optimal e-variable among all possible e-variables that are calibrated from $P$.
\begin{corollary}
    \label{cor:log-opt-cal}
Let $P$ be an exact p-variable  testing $\p$ against $\q\ll\p$.  
Assume that  $
\E^\p[\d\q /\d \p ~|~ P]$ is decreasing in $P$. Then 
the log-optimal calibrator exists and it is given (almost everywhere uniquely)  by
$$
f(P)=  \E^\p \left[\frac{\d \q}{\d \p}\mid P \right] \mbox{ and $f =0$ on $(1,\infty)$.}
$$
\end{corollary}
\begin{proof}[Proof.]




Clearly, $\int_0^1 f(u) \d u= \E^\p[f(P)] = \E^\p[\d\q/\d\p] =1$, and $f$ is decreasing, so it is a calibrator. Its optimality follows from Proposition~\ref{prop:cond-log-opt}.  
\end{proof}

\section{Asymptotic log-optimality: methods of mixtures and plug-in}
\label{sec:3-mix-plugin}
 
Consider now the case of testing a simple null hypothesis $\p$ against a composite alternative $\cQ$.  
 It is easiest to think about the case where all the distributions involved are iid product probability measures, in which case, with only a slight overload of notation, we can let $\p$ and $\q \in \cQ$ refer to the distributions of each individual observation, and let these be associated with densities $p$ and $q$.
 
Let the data be represented by a vector $X^n$. 
A reasonable objective is that of \emph{asymptotic log-optimality}. Here, the goal is to define a sequence of e-variables $E_n = E(X^n)$ such that for any $\q \in \cQ$,
\begin{equation}\label{eq:asymp-log-optimality}
 \lim_{n \to \infty} \frac{\E^{\q^n}[\log E_n]}{n} \to \kl(\q,\p).
\end{equation}
The right-hand side of the above expression is the e-power of an oracle that knows which alternative distribution is generating the data. If the above limit holds, our e-variable is able to asymptotically achieve the same e-power without this knowledge.

There are two broad classes of methods that do achieve these goals under weak assumptions (that we do not expand on here): the \emph{mixture} and the \emph{plug-in} methods.

\paragraph{Mixture method.} 
Pick a mixture distribution $\nu$ over $\cQ$, and define the e-variable
\begin{equation}\label{eq:mixture-LR}
E_n = \int \prod_{i=1}^n \frac{q(X_i)}{p(X_i)} \nu(\d q).
\end{equation}
Note the order of integral and product: we are calculating the mixture over likelihood ratios. 
The key limitation to applicability of this method is that the mixture integral should be analytically or computationally easy to evaluate.

\paragraph{Plug-in method.} Let $\q_{i-1} \in \cQ$ be a distribution of the $i$-th observation, picked based on the first $i-1$ observations, and let $q_{i-1}$ be its density. Define 
\begin{equation}\label{eq:plugin-LR}
E_n = \prod_{i =1}^n \frac{q_{i-1}(X_i)}{p(X_i)}.
\end{equation}
It is not hard to check that $E_n$ is an e-variable.


It is possible to prove the asymptotic log-optimality~\eqref{eq:asymp-log-optimality} of the above methods under some rather weak assumptions, but we omit these details. We only make some informal remarks here, which are still useful for the practitioner. For asymptotic log-optimality of the mixture method, it is necessary that $\nu$ has full support over $\cQ$; for parametric problems with finite dimensional  alternatives, this condition is often sufficient as well. For the plug-in method, a sensible choice of $\q_{i-1}$ is the posterior mean of $\q\in \cQ$ after seeing $X^{i-1}$, when starting with prior $\nu$. If the mixture method with $\nu$ is asymptotically log-optimal, then one would typically expect the aforementioned plug-in method to also have the same property.

The mixture method is also useful for composite nulls. Here we give an explicit example without derivation.
 
\begin{example}[T-test]
\label{ex:t-test}
Consider the standard t-testing situation with $\cP = \{ \mathrm{N}(0,\sigma^2): \sigma > 0 \}$ and $\cQ = \{ \mathrm{N}(\mu,\sigma^2): \sigma > 0,~\mu \neq 0 \}$ (elements of $\cP$ and $\cQ$ are the marginal distributions). Suppose we observe $n$ iid observations $X_1,\dots,X_n$ from some $\p \in \cP \cup \cQ$. Let $S_n = \sum_{i=1}^n X_i$ and $V_n = \sum_{i=1}^n X_i^2$. Then for any $c > 0$,
\[
\sqrt{\frac{c^2}{n+c^2}}\left( \frac{(n+c^2)V_n}{(n+c^2)V_n - S_n^2} \right)^{n/2}
\]
is an e-variable for $\cP$ that satisfies~\eqref{eq:asymp-log-optimality} for $\cQ$. Note that this e-variable is scale-invariant, meaning that rescaling the data by a constant does not change the e-variable, and when $n=1$, it equals 1 and has zero e-power. 

While the following observation is not about the method of mixtures, it feels appropriate to point out here the surprising fact that there do exist non-constant e-variables that have positive e-power against some alternatives, even when $n=1$. For example, for any $a \in \R$, 
\[
E_1 = |X_1| \exp((1-(X_1-a)^2)/2)
\] 
is one such e-variable for $\cP$; note however that it is not scale-invariant. To see that $E_1$ is an e-variable, note that for all $\sigma$, the inequality $|x| \leq \sigma \exp(x^2/(2\sigma^2) - 1/2)$ holds. Thus, we see that under any $\p \in \cP$, which equals $\mathrm N(0,\sigma^2)$ for some unknown $\sigma$, we have 
\begin{align*}
\E^\p[E_1] &\leq \E^\p[\sigma \exp(X_1^2/(2\sigma^2) - 1/2) \exp((1-(X_1-a)^2)/2)] \\& = \E^\p \left[ \frac{\exp(-(X_1-a)^2/2)}{\sigma^{-1}\exp(-X_1^2/(2\sigma^2))} \right] = 1,
\end{align*} 
where the final inequality holds by recognizing that we are evaluating the expectation of a likelihood ratio of $\mathrm{N}(a,1)$ to $\p$. 
For the distribution $\q=\mathrm{N}(a,1)$,
we can see $\E^\q[E_1]= \E^\q[|Z+a| \exp(1-Z^2/2)]$, where $Z$ follows a standard normal distribution under $\q$, and hence $\E^\q[E_1]>1$ when $a$ is large enough.
\end{example}

Next, we demonstrate a setting with non-iid data where the mixture method is helpful and natural.
 
\begin{example}[Changepoint e-value]
\label{ex:changepoint}
    Suppose we observe $n$ ordered observations $X_1,\dots,X_n$. The null hypothesis equals the singleton $\{\p^n\}$ for some given $\p$. The alternative hypothesis equals $\cQ = \{\p^{n-k} \q^{k}\}$ for some unknown $k \in[n]$ and known $\q$, meaning that $X_1,\dots,X_{n-k}$ are iid from $\p$, and the rest are iid from $\q$. It is natural to consider a mixture over the unknown $k$, in particular forming the  e-value  
    \[
        \frac{1}n\sum_{k=1}^{n } \prod_{i=n-k+1}^{n} \frac{q(X_i)}{p(X_i)},
    \]
    based on uniform-mixture likelihood ratios. 
\end{example}

In fact, the above example underlies the much more general, nonparametric, \emph{sequential} change detection methodology of \emph{e-detectors}, which is out of the scope of this book.

\section{Axiomatic justification of the e-power} 
\label{sec:axioms}

In this section we justify the definition of e-power in Definition~\ref{def:e-power} through an axiomatic approach.

 Let $\mathcal X$ be the set of all nonnegative random variables that are bounded from above and away from $0$. 
 We will restrict our consideration of e-power on $\mathcal X$ to avoid issues with infinity,
 and it can be naturally extended to the set of all nonnegative random variables.
As standard in the axiomatic approach in decision theory, axioms become weaker when they are formulated on a smaller set of objects, making the axiomatic characterization results stronger. Therefore, our restriction to $\mathcal X$ is harmless. 

 Let $\mathcal M_1$ be the set of all probability measures on $(\Omega,\mathcal F)$.
 We will fix the alternative hypothesis $\q \in \mathcal M_1$, which is an atomless probability measure. The requirement of $\q$ being atomless is only used to justify the existence of a non-degenerate iid sequence.

In general, e-variables for the null hypothesis $\mathcal P$ may be arbitrarily distributed under the alternative hypothesis $\q$, and the e-power is a concept concerning the performance of an e-variable under the alternative hypothesis.
Therefore, the choice of the null hypothesis $\mathcal P$ is irrelevant 
 for defining e-power.  
 
 For this reason,  without loss of generality, we   directly consider the candidates for the e-power as mappings $\epo^\q:\mathcal X \to  \R$, without specifying that e-variables need to satisfy $\E^\p [E]\le1$ for $\p\in \mathcal P$.

For $\epo^\q$ to represent the e-power under $\q$, the natural interpretation of $\epo^\q(E_1)>\epo^\q(E_2)$ is that $E_1$ is considered as more powerful than $E_2$. We keep this important interpretation in mind in the following discussions.

Suppose that a statistician has designed a two-stage experiment, and observes an e-value of $e_0$ at the first stage. She needs to choose between two e-variables $E_1$ and $E_2$ 
at second stage. If  $E_1$ is more powerful than $E_2$, then she should choose $E_1$ over $E_2$. 
However, the effective e-variables to compare in this procedure are $ e_0 E_1 
$ and $e_0 E_2$.
Hence, there should be a consistency between the pair $(E_1,E_2)$ and the pair $(e_0E_1,e_0E_2)$,
which is formalized by the following axiom.

\begin{axiom}[Homogeneity]
\label{ax:H}
If  $\epo^\q(E_1 )\ge  \epo^\q(E_2)$, then  $\epo^\q(  e_0 E_1)  \ge \epo^\q( e_0 E_2)$ for any constant $e_0>0$.
\end{axiom} 

 The second axiom concerns the design of e-variables in an asymptotic sense.  
 It states that, to choose between two configurations of iid e-variables, if the  product e-process of  less powerful configuration has asymptotic power $1$, then that from  
 the more powerful configuration should also have asymptotic power $1$. 
This axiom reflects the consideration of the Kelly criterion, which 
we discussed  in Section~\ref{sec:e-power}. 

\begin{axiom}[Asymptotic consistency]
\label{ax:AC} 
Suppose that $(E_n)_{n\in \N}$ and $(E'_n)_{n\in \N}$ 
are two iid sequences of nonnegative random variables under $\q$.
If $\epo^\q(E_1 )\ge  \epo^\q(E_1')$,
then for any $\alpha\in (0,1)$, 
$$ \q\left(\prod_{t=1}^n E'_t \ge \frac 1 \alpha \right) \to 1 ~\Longrightarrow~ \q\left(\prod_{t=1}^n E_t\ge \frac 1 \alpha\right) \to 1,$$ both as $n\to\infty$.
\end{axiom} 

It turns out that Axioms 
\ref{ax:H} and~\ref{ax:AC} are sufficient to jointly characterize the e-power, which has the form $\E^Q[\log E]$ as we saw in  Definition~\ref{def:e-power}.

\begin{theorem}
\label{th:axiom-ch}
Given an atomless probability measure $\q\in \mathcal M_1$, 
a mapping $\epo^\q:\mathcal X\to \R$ satisfies Axioms 
\ref{ax:H} and~\ref{ax:AC} if and only if it can be represented by
 $$\epo^\q(E)=f(\E^\q[\log  E]),~~~E\in \mathcal X$$ for some strictly increasing function $f$.
\end{theorem}
\begin{proof}[Proof.]
It is straightforward to verify that $E\mapsto \E^\q[\log E]$ satisfies 
the two axioms; Axiom~\ref{ax:AC} can be checked with the strong law of large numbers.
Clearly, a strictly increasing transform does not matter for the two axioms. Hence,  $E\mapsto f( \E^\q[\log E])$ also satisfies the two axioms.

Now we show the ``only if'' direction. 
We first suppose for the purpose of contradiction that there exist $E,E'$ such that
\begin{equation}
\label{eq:contra}
\epo^\q(E) > \epo^\q(E') \mbox{~and~} \E^\q[\log E] < \E^\q[\log E'].
\end{equation}
Using Axiom~\ref{ax:H}, we can find some $c>0$ such that 
\begin{equation}
\label{eq:contra2}
\epo^\q(c E) > \epo^\q(c E') \mbox{~and~} \E^\q[\log (c E)] <0< \E^\q[\log (c E')].
\end{equation}
Using Axiom~\ref{ax:AC}, for the iid sequence $(E_n)_{n\in \N}$  distributed as $cE$ 
and the iid sequence $(E'_n)_{n\in \N}$  distributed as $cE'$,
 the strong law of large numbers and \eqref{eq:contra2} together imply $\log \prod_{t=1}^n E _t \to -\infty$ and $\log \prod_{t=1}^n E '_t\to \infty$ in $\q$. 
 This violates Axiom~\ref{ax:AC}. Hence, \eqref{eq:contra} does not hold. 
 This implies
 \begin{equation}
\label{eq:contra3}
\epo^\q(E)  > \epo^\q(E') ~\Longrightarrow ~\E^\q[\log E]  \ge  \E^\q[\log E'].
\end{equation}
Equivalently,
 \begin{equation}
\label{eq:contra4}
\epo^\q(E)  \le \epo^\q(E') ~\Longleftarrow ~\E^\q[\log E]  <  \E^\q[\log E'].
\end{equation}
Using \eqref{eq:contra4}, we know that $c\mapsto \epo^\q(cE)$ is increasing on $(0,\infty)$ for each $E\in \mathcal X$.

Suppose for the purpose of contradiction that there exist $E,E' \in \mathcal X$ such that 
 \begin{equation}
\label{eq:contra5}
\epo^\q(E)  > \epo^\q(E') \mbox{~and~}\E^\q[\log E] =  \E^\q[\log E'].
\end{equation}
Using Axiom~\ref{ax:H}, we have from \eqref{eq:contra5}, for any $c>0$,
 $
\epo^\q(cE)  > \epo^\q(c E').
 $ 
Note that for any $\lambda <1$, we have, from \eqref{eq:contra4}, 
$\epo^\q(\lambda cE)  \le \epo^\q(c E')$, because $\E^\q[\log (\lambda c E)]  <  \E^\q[\log (c E')]$.
Hence, $\lambda\mapsto \epo^\q(\lambda cE)$ is discontinuous at $\lambda =1$ for all $c>0$.
This implies that  $c\mapsto \epo^\q(cE)$ is discontinuous everywhere, a contradiction to the fact that this mapping is monotone and hence has bounded variation. 
 Therefore, \eqref{eq:contra5} does not hold,
 and  by \eqref{eq:contra3}, we have 
 $$
\epo^\q(E)  > \epo^\q(E') ~\Longrightarrow ~\E^\q[\log E] >  \E^\q[\log E'].
 $$
This, together with \eqref{eq:contra3}, 
shows that $\epo^\q$ and $X\mapsto \E^\q[\log X]$ generate the same ordinal relation, and hence the representation holds.   
\end{proof}
Theorem~\ref{th:axiom-ch} implies that if Axioms~\ref{ax:H} and~\ref{ax:AC} are reasonable for a notion of e-power, then the only possible way to quantify the e-power  is to use a strictly increasing transform of $\E^\q[\log E]$, which generates the same ordinal structure as using 
$\E^\q[\log E]$ directly.
It also shows what the e-power essentially guarantees: It is the most suitable when considering the product e-value with asymptotic consistency
(also discussed in Section~\ref{sec:e-power}). A simple example of a strictly increasing transform of the e-power is  the geometric mean of the e-variable $E$, given by $\exp\E^\q[\log E]$. 
 It is nonnegative, and  it is larger than $1$ if the  e-power is positive. 
 
In the above result,  the e-power $\E^\q[\log E]$   in Definition~\ref{def:e-power}  is  justified as a comparative notion.  This serves 
 a very similar purpose to the classic notion of power of tests (Definition~\ref{def:nontrivial}) when we compare e-variables to determine which ones are more ``powerful''.  
When it comes to the  meaning of the actual value, its natural interpretation is the asymptotic 
 growth rate in a multiplicative setting, which is different from the classic notion of power, whose interpretation is the probability of rejection under the alternative.

We can easily extend the mapping $\epo^\q$ to the set of all 
nonnegative random variables, provided that $\E^\q[\log E]$ is well-defined, that is, at least one of $\E^\q[(\log E)_+]$
  and $\E^\q[(\log E)_-]$ is finite.

In what follows, we let  
$$\epo^\q: X\mapsto \E^\q[\log X]$$ for any nonnegative random variable $X$.
  We treat $\epo^\q(X)$ as undefined if $\E^\q[\log X]$ is not well-defined.
  The next proposition gives some convenient properties of $\epo^\q$.
  \begin{proposition}\label{prop:e-power-properties}
    Let $E$ be an e-variable. 
    \begin{enumerate}[label=(\roman*)]
  \item If $\epo^\q(E)> 0$, then 
  $\epo^\q(1-\lambda +\lambda E) >0$ for all $\lambda \in (0,1]$.
\item  If $\epo^\q(E) \ge  0$, then 
  $\epo^\q(1-\lambda +\lambda E) \ge 0$ for all $\lambda \in [0,1]$.
  \item   $\epo^\q(1-\lambda +\lambda E)  $ is always well-defined for   $\lambda \in [0,1)$.
\item  If $\epo^\q(E) <  \infty$, then 
  $\epo^\q(1-\lambda +\lambda E) \in \R $ for all $\lambda \in [0,1)$.
  \end{enumerate}
   \end{proposition}
     \begin{proof}[Proof.]
     For (i), it suffices to note that $f: \lambda \mapsto \epo^\q(1 +\lambda (E-1)) $ is concave with $f(0)=0$ and $f(1)>1$. This implies $f(\lambda)>0$ for all $\lambda \in (0,1]$. 
     The case of (ii) is similar. 
     Part (iii) follows by noticing 
     that $ 1-\lambda+\lambda E  $ is bounded away from $0$, and hence 
     $\E^\q [( \log (1-\lambda+\lambda E))_- ]<\infty$.
     Part (iv) follows from $\E^\q [( \log (1-\lambda+\lambda E))_+ ] \le \E^\q [ (\log  E)_+ ]<\infty$, which holds by noting 
     $ \q(  \log (1-\lambda+\lambda E) >x)\le  \q (\log E>x)$ for $x>0$, together with  the  observation  $\E^\q [( \log (1-\lambda+\lambda E))_- ]<\infty$ in (iii). 
  \end{proof}



\section*{Bibliographical note}

The connection between tests and bounded e-variables described in Fact~\ref{fact:duality} was discussed in \cite{koning2024continuous} in detail. The e-power 
in Section~\ref{sec:e-power} 
 is formally defined in \cite{vovk2024nonparametric}. Example~\ref{ex:pathological} was discussed in \cite{wang2024proposer}. The notion of log-optimality was proposed by~\cite{Kelly56} in the context of betting and investment. \cite{larsson2024numeraire} showed the existence and uniqueness of a log-optimal e-variable with no conditions on $\cP,\q$.
The result in Section~\ref{sec:3-p-vs-q} essentially dates back to~\cite{breiman1961optimal}, though its presentation here mirrors that in~\cite{Shafer:2021}.
The mixture and plug-in methods to handle composite alternatives, as presented in Section~\ref{sec:3-mix-plugin}, are very well known; in fact, both ideas were already mentioned by~\cite{Wald47}.
Later, \cite{robbins1974expected} derived some strong connections between these methods. The t-test e-variable is derived and studied in, for example,~\cite{Lai:1976},~\cite{perez2022estatistics}, and~\cite{wang2023anytime}, with the latter paper having the e-variable for $n=1$.  
The asymptotic optimality of the method of mixtures is strongly related to information-theoretic topics such as redundancy and description length, and many more such optimality results can be inferred using a plethora of results on the latter topic in~\cite{grunwald2007minimum}. The changepoint e-value in Example~\ref{ex:changepoint} is borrowed from~\cite{shin2022detectors}.

\chapter{Post-hoc testing and decision making with e-values}

\label{chap:posthoc}

In this chapter, 
we study e-values in the context of post-hoc testing and decision making and compare them with p-values. 
The post-hoc setting allows for either the risk level, the possible set of decisions, or the risk constraint to  depend on the data that are used for the decision.
It also allows for future data collection and experiment design to depend on previously collected data. 
The main message is that e-values play a central role in the post-hoc  testing and decision problems. We also note that several of the ideas presented here have been extended to \emph{multiple} hypothesis testing in Section~\ref{sec:9-closed-ebh}.

\section{Testing at data-dependent error levels}\label{sec:2-posthoc}

Recall that a level-$\alpha$ test for $\cP$ is defined to be a function $\phi$ such that $\E^\p[\phi] \leq \alpha$ for every $\p \in \cP$. 
As presented so far, the type-I error level $\alpha$ has to be a predefined constant. However, e-variables do yield a generalized form of error control when testing at data-dependent error levels $\hat \alpha$ that could be arbitrarily dependent on the e-variable. 
A natural way to  incorporate random $\hat \alpha$ is through  
a uniform generalization of the constraint $\E^\p[\phi/\alpha]\le 1$, although there may be other ways to define post-hoc validity.

\begin{definition}[Post-hoc valid tests and p-values]\label{def:post-hoc}
A \emph{post-hoc valid family of tests} for $\cP$ is  a collection of tests $(\phi_\alpha)_{\alpha \in (0,1)}$ 
such that for every $\p \in \cP$, and for every random variable $\hat \alpha$ taking values in $(0,1)$, we have $\E^\p\left[ \phi_{\hat \alpha} / \hat \alpha \right] \leq 1$. Equivalently, we require that
\[
\E^\p\left[ \sup_{\alpha \in (0,1)} \frac{\phi_\alpha}{\alpha} \right] \leq 1 \mbox{~~~~for all $\p\in \cP$}.
\] 
A \emph{post-hoc p-variable} for $\cP$ is a nonnegative random variable $P$ such that for every $\p \in \cP$, and for every random variable $\hat \alpha$ taking values in $(0,1)$, we have $$\E^\p\left[ \frac{\id_{\{P \leq \hat \alpha\}} }{ \hat \alpha} \right] \leq 1.$$  
Equivalently, we require that
\begin{equation}
    \label{eq:post-hoc-p-value}
\E^\p\left[ \sup_{\alpha \in (0,1)} \frac{\id_{\{P \leq \alpha\}}}{\alpha} \right] \leq 1,
\end{equation}
for all $\p\in \cP$.
\end{definition}

For a p-variable $P$, 
a weaker condition than 
\eqref{eq:post-hoc-p-value} holds:
\begin{equation}
    \label{eq:post-hoc-p-value2}
\sup_{\alpha \in (0,1)} 
 \E^\p\left[ \frac{\id_{\{P \leq \alpha\}}}{\alpha} \right] \leq 1 \mbox{~~~~for all $\p\in \cP$}.
\end{equation}
Indeed, \eqref{eq:post-hoc-p-value2} is precisely the definition of a p-variable $P$. 
\begin{remark}
\label{rem:p-to-e-random-alpha}
Let $P$ be a p-variable  for $\p$. 
The condition 
\begin{equation}
    \label{eq:post-hoc-p-value-3}
\E^\p\left[ \frac{\id_{\{P \leq T\}}}{T} \right] \leq 1 \end{equation}
also holds for some nonnegative random variables $T$, where $0/0$ is understood as $0$. 
A sufficient condition for \eqref{eq:post-hoc-p-value-3} on $T$ is
\begin{equation}
    \label{eq:post-hoc-p-value-4}
p\mapsto \p(T \le t \mid P\le p)  \mbox{~is increasing 
for all $t\in\R$;}
\end{equation}
 this defines a form of negative dependence between $P$ and $T$. Some examples satisfying \eqref{eq:post-hoc-p-value-4} include  (a) $T$ is a constant; (b) $T$  is independent of $P$; (c) $T$ is a decreasing function of $P$. 
 Note that condition \eqref{eq:post-hoc-p-value-3} means that $\id_{\{P\le T\}}/T$ is an e-variable.
Moreover, \eqref{eq:post-hoc-p-value-3} is implied by  
 $$
p\mapsto \p(T \le t \mid P= p)  \mbox{~is (almost everywhere) increasing 
for all $t\in\R$,}
 $$
which is a stronger condition.
\end{remark} 

Clearly, in Definition~\ref{def:post-hoc}, each member $\phi_\alpha$ of a post-hoc valid family of tests is itself a level-$\alpha$ test for any fixed $\alpha \in (0,1)$, and $P$ is itself a p-variable. 
It is clear that for a post-hoc p-variable $P$, $(\id_{\{P\le\alpha\}})_{\alpha\in(0,1)}$ is a post-hoc valid family of binary tests; this is not true for a usual p-variable. 
Moreover, it does not hurt to require a post-hoc valid family $(\phi_\alpha)_{\alpha \in (0,1)}$  of tests to be increasing (that is, $\phi_\alpha$ is increasing in $\alpha$), because replacing $\phi_\alpha $ with $\sup_{\gamma \ge \alpha} \phi_\gamma$ does not affect its post-hoc validity. 

We say that a post-hoc valid family $(\phi_\alpha)_{\alpha \in (0,1)}$  
 of tests
is dominated by 
another family  $(\phi'_\alpha)_{\alpha \in (0,1)}$   
if $\phi_\alpha\le \phi'_\alpha$ for all $\alpha \in (0,1)$.
Intuitively, a dominated family is less useful than a dominating family, because the latter rejects at least as often as the former.  
In the next result, we will see   testing with e-variables   is essentially the only useful way of constructing post-hoc valid families of (binary) tests.

\begin{proposition}
\label{prop:c2-post-hoc-1}
   Let $E$ be an e-variable   for $\cP$. Then,  $((\alpha E) \wedge  1)_{\alpha \in (0,1)}$ is a post-hoc valid family of tests and $(\id_{\{E \geq 1/\alpha\}})_{\alpha \in (0,1)}$ is a post-hoc valid family of binary tests.
    
    Conversely, every post-hoc valid family of 
      tests is dominated by $((\alpha E) \wedge  1)_{\alpha \in (0,1)}$  for some e-variable $E$ for $\cP$,
      and  every post-hoc valid family of
    binary tests is dominated by $(\id_{\{E \geq 1/\alpha\}})_{\alpha \in (0,1)}$ for some e-variable $E$ for $\cP$.        
\end{proposition}
\begin{proof}[Proof.]
It is elementary to verify  
    \begin{equation}
        \label{eq:c2-simple-eq}
      \id_{\{x \geq 1\}} \leq x\wedge 1  \le x \mbox{~~~for~} x \geq 0 , 
    \end{equation}
    which can be used to show Markov's inequality. 
Using \eqref{eq:c2-simple-eq},  for any $\hat \alpha$ taking values in $(0,1)$, we have
    $$\E^\p\left[ \frac{ \id_{\{E \geq 1/\hat \alpha\}}  }{ \hat \alpha } \right] \leq \E^\p\left[ \frac{ (\hat \alpha E )\wedge 1 }{ \hat \alpha }\right]   \leq \E^\p\left[ \frac{ \hat \alpha E  }{ \hat \alpha }\right] = \E^\p\left[ E  \right] \leq 1.$$ 
This shows that both $(\id_{\{E \geq 1/\alpha\}})_{\alpha \in (0,1)}$ and $((\alpha E) \wedge  1)_{\alpha \in (0,1)}$ are post-hoc valid families of tests.

    For the converse claims, let
   $(\phi_\alpha)_{\alpha \in (0,1)} $ be a post-hoc valid family of tests, and let 
    $E=\sup_{\alpha \in (0,1)} \phi_{\alpha}/\alpha$. 
    By the definition of post-hoc validity, $E$ is an e-variable for $\cP$.
    Moreover, $\phi_\alpha \le (\alpha E)\wedge 1$ by definition, showing the first converse claim. 
    If members of $(\phi_\alpha)_{\alpha \in (0,1)} $ are binary, then   $\{\phi_\alpha=1\} \subseteq \{E\ge 1/\alpha\}$ for every $\alpha \in(0,1)$. Therefore,   $
    \phi_\alpha \le \id_{\{E\ge 1/\alpha\}}$, showing the desired dominance in the second converse claim.   
\end{proof}

It is important that in Proposition~\ref{prop:c2-post-hoc-1}, the same e-variable $E$ is used to construct all members in a post-hoc valid family of tests. 
The next result gives an intimate connection between post-hoc  p-variables and  e-variables.

\begin{proposition}
\label{prop:c2-post-hoc}
For a random variable $P$, the following are equivalent:
    \begin{enumerate}[label=(\roman*)]
        \item $P$ is a post-hoc p-variable  for $\cP$; 
        \item $P= 1/E$ on the event $P <  1$ for some e-variable $E$ for $\cP$;
        \item $P\ge (1/E)\wedge 1$  for some e-variable $E$ for $\cP$.
    \end{enumerate}  
\end{proposition}
\begin{proof}[Proof.] 
To see (i)$\Rightarrow$(ii), 
first note that $\p(P=0)=0$ for all $\p\in \cP$, since $P$ is a p-variable for $\cP$. 
Let $E=(1/P)\id_{\{P< 1\}}$ and $\hat \alpha = P\id_{\{P<1\}} + \gamma\id_{\{P\ge 1\}}$ for some $\gamma \in(0,1)$.
We have that $\hat \alpha$ is $(0,1)$-valued ($\p$-almost surely for every $\p\in \cP$) and $\{P<1\}=\{P\le \hat \alpha\}$. Therefore, for all $\p\in \cP$,
$$
\E^\p [E] = \E^\p\left[  \frac{1}{P} \id_{\{P< 1\}}\right] 
= \E^\p\left[  \frac{1}{P\id_{\{P<1\}}} \id_{\{P\le\hat \alpha\}}\right] 
= 
\E^\p\left[  \frac{1}{\hat \alpha} \id_{\{P\le \hat \alpha\}}\right] 
\le 1,
$$
where the last inequality follows  from the definition of the post-hoc   p-variable.   
Hence, $E$ is an e-variable, and clearly $P=1/E$ when $P< 1$. 
The direction (ii)$\Rightarrow$(iii) is immediate. Finally, the direction (iii)$\Rightarrow$(i) follows from 
\begin{align*}
\E^\p\left[ \sup_{\alpha \in (0,1)} \frac{\id_{\{P \leq \alpha\}}}{\alpha} \right]
&\le 
\E^\p\left[ \sup_{\alpha \in (0,1)}  
 \frac{\id_{\{(1/E)\wedge1  \leq \alpha\}}}{\alpha} \right]  \\
 &=\E^\p\left[ \sup_{\alpha \in (0,1)}  
 \frac{\id_{\{ 1/E   \leq \alpha\}}}{\alpha} \right] \le \E^\p[E]  
\leq 1   
\end{align*}
for all $\p\in \cP$, where the penultimate inequality is due to \eqref{eq:c2-simple-eq}.
\end{proof}

There is a small variation of the result in Proposition~\ref{prop:c2-post-hoc}, further connecting post-hoc p-values and e-values. 
If we change the range of $\alpha $ from $(0,1)$ to $(0,\infty)$, that is,  using the stronger requirement
\begin{equation}
    \label{eq:c4-nick}
    \E^\p\left[ \sup_{\alpha \in (0,\infty)} \frac{\id_{\{P \leq \alpha\}}}{\alpha} \right] \leq 1 \mbox{~~~~for all $\p\in \cP$}
\end{equation} 
in the definition of a post-hoc p-variable $P$,  
then a post-hoc p-variable would be precisely $1/E$ for some e-variable $E$, following the same argument in the proof of Proposition~\ref{prop:c2-post-hoc}.

From Proposition~\ref{prop:c2-post-hoc}, it is immediate that $(1/E)\wedge 1$ is a post-hoc p-variable; so is $1/E$, which may take value larger than $1$ (recall that we allow p-variables and post-hoc p-variables to take values larger than $1$). 
Using the same result, an admissible (i.e., not dominated by another) post-hoc valid family  $(\phi_\alpha)_{\alpha \in (0,1)}$ of binary tests can  always be written as 
$$
\phi_\alpha=\id_{\{P\le \alpha\}} \mbox{~~~for $\alpha \in(0,1)$},
$$ 
for some post-hoc p-variable $P$. 
In other words, a post-hoc p-variable is the smallest level that the associated admissible post-hoc family of binary tests can reject the null hypothesis. P-values are sometimes said to equal the smallest level at which one could have rejected the null, \emph{had that level been fixed in advance}.  Our preceding discussion implies that the italicized text can be ignored for post-hoc p-variables, but in some scientific research  it is also ignored for usual p-variables, leading to misuse or misterpretation of p-values.  

The connection between $1/E$ and post-hoc p-variables 
also illustrates that it is unfair to directly compare the value of $1/E$ for an e-variable with the value of a usual p-variable $P$, as the former offers a stronger type of guarantee: post-hoc validity in the sense of Definition~\ref{def:post-hoc}. 

To interpret this stronger validity, 
compare the situations that we observe  an e-value 
 $E=25$ versus a p-value $P=0.04$ on the same hypothesis.
 We can safely say (in the sense of Definition~\ref{def:post-hoc}) that we reject the null hypothesis based on $E$ at level $0.05$ (and even at $0.04$), but it would not be statistically valid to say that we can reject based on $P$ at level $0.05$, unless this level was prespecified before collecting and seeing data. Indeed,  one cannot claim rejection at level $0.05$ even if one observes $P=0.0001$, if the level $0.05$ is not prespecified.
 The procedure of first seeing $P$ and then specifying the significance level is a statistical malpractice, but it appears often in published studies. For instance, in many papers in the social sciences,  first p-values are computed (sometimes on many different configurations of regression analysis) and then they are marked with  asterisk signs to indicate at which levels they can reject the corresponding null hypothesis.
 For such a purpose, e-values (or post-hoc p-values) are the right tool to use, instead of the usual p-values. 
To re-emphasize, an e-value equal to $1/a$ carries more information than a p-value equal to $a$ for $a\in(0,1)$; see Section~\ref{sec:c2-jeffrey} for more discussions on comparing these values.




\section{Post-hoc decisions using p-values and e-values}
\label{sec:post-hoc-decisions}

In this section, we show that e-value-based decision rules naturally appear in the context of post-hoc decisions, in a similar way to Section~\ref{sec:2-posthoc} but with a different motivation and interpretation.
We consider a framework of decision making  that is more general than standard hypothesis testing, which concerns binary decisions of acceptance and rejection. 

Consider a class of  loss functions $L_b: \{0,1\} \times D_b\to \R$ indexed by scenarios $b\in \mathcal B\subseteq \R$, where $D_b\subseteq \N_0$ represents the set of possible decisions. (It may be instructive to first consider $\mathcal B\subseteq\N_0$ for simplicity.)
A scenario $b$ refers to a possible set of decisions  that may be known after seeing the data. For instance, a scientist may have a different set of possible recommendations to the policy makers when they see different data results (see Example~\ref{ex:vaccine} below). For the moment, we focus on a fixed scenario $b$, and the general case is treated in Section~\ref{sec:c4-22}.

The value $L_b(0,d)$ represents the loss to the decision maker 
with decision $d \in D_b$ when the null hypothesis $\cP$ is true;
the value $L_b(1,d)$ represents the  loss to the decision maker 
with decision  $d \in D_b$ when  the alternative hypothesis is true. The decision maker knows the loss functions and the scenario.

Without loss of generality, 
we assume that
\begin{enumerate}
    \item $d \mapsto L_b(0, d)$ is strictly increasing, and
    \item $ d \mapsto L_b(1,d)$ is strictly decreasing.
\end{enumerate}
This assumption is harmless because 
we can discard all decisions $d'$ with $L_b(0,d')\ge L_b(0,d)$ and $L_b(1,d')\ge L_b(1,d)$  
for some $\delta$, because $d'$ will never be better than $d$. After discarding these decisions, we can rearrange the decisions according to an increasing order of $L_b(0,d)$.
(When $b$ varies, we can do this for each $b$.)
Intuitively, one can interpret a larger $d$ as a more aggressive decision,
and the above assumption means that as $d$ increases, there is more loss to the decision maker if the null is true, and there is more gain to the decision maker if the alternative is true (see Example~\ref{ex:p-hack} below).

The dataset is represented by a random vector $X$ taking values in a set $\cZ$ (e.g., $\R^n$), 
and $\delta(X)\in D_b$ is the decision made based on the data $X$.
When we allow $b$ to vary,
$\delta$ will also depend on $b$.
The rest of this section will show:
\begin{itemize}
    \item For binary decision spaces and fixed scenarios, p-values can safely be used to make decisions. 
    \item In non-binary decision problems or in settings where the scenario could be data-dependent, decisions based on p-values do not yield the desired risk guarantee, but those based on e-values do.
\end{itemize}

For the sake of illustration, we assume $\cQ=\{\q\}$ in this section, that is, a simple alternative. 
The problem can be easily generalized using a composite $\cQ$, by considering a more general objective in the optimization, which is not essential to our discussion below.

\subsection{Decision rules based on p-values and e-values}

Given a fixed scenario $b\in \mathcal B$, the decision maker's problem is  
 \begin{align}
\label{eq:PG-1}
\begin{aligned}
\mbox{ to minimize~~~} &\E^\q[L_b(1,\delta(X))]
\\
\mbox{ over~~~} & \delta:\cZ \to D_b
\\
\mbox{ subject to~~~} &\sup_{\p \in \cP} \E^\p[L_b(0,\delta(X))] \le \gamma,
\end{aligned}
\end{align} 
where $\gamma>0$ is a prespecified risk budget, and the constraint in \eqref{eq:PG-1} is called a risk constraint. The scenario $b$ will be allowed to depend on the data later.

\begin{example}[Ternary significance decisions are common scientific practice]
\label{ex:p-hack}
    A practical example in scientific research in the above framework is to assign a label from ``non-significant (NS), significant (S), very significant (VS)'' to the test result  after seeing data. In our notation, it is $D_b=\{0,1,2\}$, with $0$ representing NS, $1$ representing S, and $2$ representing VS.
    If the null hypothesis is true, then making the decision NS has the least loss. If the alternative hypothesis is true, then making the decision VS has the least loss. 
    A common choice is to set $\gamma=1$,  $L_b(0,0)=0$,
    $L_b(0,1)=1/\alpha_1$
    and   $L_b(0,2)=1/\alpha_2$ for some $0<\alpha_2<\alpha_1<1$, 
    and the values   $L_b(1,\delta)$ depend on the actual problem.
    \end{example}

We define two potential approaches based on p-values and e-values. 
Suppose that a p-variable $P(X)$ and an e-variable $E(X)$ for $\cP$ are computed from the data.  
The decision   $\delta_P$ based on the p-value is given by
\begin{align}
\label{eq:PG-p}
\delta_P(X) = \max \{ d \in D_b:  P(X) L_b(0,d) \le\gamma\}.
\end{align}
The decision   $\delta_E$ based on the e-value is given by
\begin{align}
\label{eq:PG-e}
\delta_E(X) = \max \{ d \in D_b:   L_b(0,d) \le\gamma E(X) \}.
\end{align}

\begin{example}[Significance decisions based on p-values and e-values]
\label{ex:p-hack2}
In the context of Example~\ref{ex:p-hack}, the p-value-based decision rule in
\eqref{eq:PG-p} yields $\delta_P(X)=1$ (significant) if $\alpha_2<P(X)\le \alpha_1$
and $\delta_P(X)=2$ if $P(X)\le \alpha_2$ (very significant); otherwise $\delta_P(X)=0$ (non-significant). 
This is a common practice in scientific research, but it does not satisfy the risk control in \eqref{eq:PG-1} as shown in Example~\ref{ex:p-hack3} below. 

The e-value-based decision rule in
\eqref{eq:PG-e} yields $\delta_E(X)=1$ if $1/\alpha_1\le E(X)< 1/\alpha_2$
and $\delta_E(X)=2$ if $E(X)\ge 1/\alpha_2$; otherwise $\delta_E(X)=0$. This decision rule always satisfies the risk control in \eqref{eq:PG-1}; see Theorem~\ref{th:PG}.
\end{example}

\subsection{Binary decisions in fixed scenarios}

For now,  consider the simplest case $D_b=\{0,1\}$, that is, the case of binary decisions, where  $0$ represents the decision of accepting the null hypothesis, 
and $1$ represents the decision of rejecting the null hypothesis. 
Without loss of generality, the loss function can be specified by $L_b(0,0)=L_b(1,1)=0$ and $L_b(0,1)=L_b(1,0)=1$  (this is without loss of generality because of the monotonicity  assumed above). 
In this case, \eqref{eq:PG-1} becomes 
 \begin{align}
\label{eq:PG-2}
\begin{aligned}
\mbox{ to maximize~~~} &\q(\delta(X)=1)
\\
\mbox{ over~~~} & \delta:\cZ \to \{0,1\}
\\
\mbox{ subject to~~~} &\sup_{\p \in \cP} \p(\delta(X)=1) \le \gamma,
\end{aligned}
\end{align} 
which is the standard problem of maximizing power while keeping type-I error control $\gamma$.

In the context of binary decisions in \eqref{eq:PG-2},  the p-value-based decision rule  means choosing $\delta_P(X)=1$ if $P(X)\le \gamma$ and choosing $\delta_P(X)=0$ if $P(X)>\gamma$; thus, it  is the standard way in which p-values are used in decision making (see also Example~\ref{ex:p-hack2}). 
The e-value-based decision rule means choosing $\delta_E(X)=1$ if $E(X)\ge1/ \gamma$ and choosing $\delta_E(X)=0$ if $E(X)<1/\gamma$.

\begin{proposition}
Given a fixed scenario $b$ with $D_b=\{0,1\}$, 
both $\delta_P$ in \eqref{eq:PG-p} and $\delta_E$ in \eqref{eq:PG-e} satisfy the risk constraint in \eqref{eq:PG-2}.
\end{proposition}
\begin{proof}[Proof.]
 To see that $\delta_P$ is feasible for problem \eqref{eq:PG-2}, 
we note that for $\p\in \cP$, 
$$\E^\p[L_b(0,\delta_P(X))]=\p(\delta_P(X)=1) =\p(P(X)\le \gamma) \le \gamma.$$
To see that $\delta_E$ is feasible for problem \eqref{eq:PG-2}, 
we note that for $\p\in \cP$, 
\begin{align}
\label{eq:PG-e2} \E^\p[L_b(0,\delta_E(X))]\le\E^{\p}[ \gamma E(X)] \le \gamma,\end{align}
as desired. 
\end{proof}

Since $1/E(X)$ is a p-variable (Section~\ref{sec:c2-markov}), the e-value-based approach is typically more conservative than the p-value-based approach. 
In Section~\ref{sec:c4-22}, we will see that when $D_b $ is not binary, $\delta_E$ maintains the risk control, whereas $\delta_P$ does not.

\subsection{Non-binary decisions or data-dependent scenarios}
\label{sec:c4-22}

In many real-world applications, the testing problem, specified by $(L_b,D_b)$, depends on the scenario $b$, which itself depends on the data through, for example, the result of a preliminary exploratory experiment, the evidence provided by the first half of the data, or
 the amount of collected data (which the scientist may not know before collecting them). 
A simple example of data-dependent scenarios with $\mathcal B=\{1,2\}$ is given below. 

As the scenario $b$  depends on the data $X$,  we represent it by a random variable $B$ (which is a function of $X$), whose realization is observable to the decision maker. 
The decision rule $\delta$ now takes both realizations of $X$ and $B$ as its inputs, written as $\delta(X,B)$.

\begin{example}
\label{ex:vaccine}
 After seeing data from a clinical trial experiment, if the data is promising, the scientist finds themselves in the situation $B=1$, in which case they can choose not to vaccinate anyone ($\delta(X,1)=0$)
 or to vaccinate the highest-risk population ($\delta(X,1)=1$).  
 If the vaccine turns out (unexpectedly) to be extremely effective (the case $B=2$), then the second option  $(\delta(X,2)=1)$ changes to vaccinating the \emph{full} population.
 Note that the decision sets $D_1$ and $D_2$ have different meanings, and in reality the scientist must choose from the decision set $D_{B}$.
In this example, we may set $L_b(0,0)=L_b(1,1)=0$ for $b\in \{1,2\}$,
$L_1(0,1)=L_1(1,0)=L_2(1,0)=1$
and  $L_2(0,1)=10$.  (To be pedantic, the first argument of the loss function still indicates whether the null is true or not (unknown to the scientist), while the second argument indicates the decision taken in the scenario in the subscript.) The scientist actually receives the loss $L_{B}(0,\delta(X,B))$ or $L_{B}(1,\delta(X,B))$, depending on the unknown state of nature.
\end{example}

The problem \eqref{eq:PG-1} becomes
 \begin{align}
\label{eq:PG-3}
&\begin{aligned}
\mbox{ to minimize~~~} &\E^\q[L_{B}(1,\delta(X,B))]
\\
\mbox{ over~~~} & \delta: \cZ \times \mathcal B    \to \N_0;~ \delta(x,b)\in  D_{b}
\\
\mbox{ subject to~~~} &\sup_{\p \in \cP} \E^\q[L_{B}(0,\delta(X,B))] \le \gamma, \end{aligned}
\end{align}   
The constraint and objective in \eqref{eq:PG-3} are formulated for the  given random variable $B$.
In case how $B$ depends on the data is not known to the decision maker, one may consider a stronger constraint and a worst-case objective, that is,
\begin{align}
    \label{eq:PG-4}
    \begin{aligned}
\mbox{ to minimize~~~} &\E^\q\left[\sup_{b\in \mathcal B}L_{b}(1,\delta(X,b))\right]
\\
\mbox{ over~~~} & \delta: \cZ \times \mathcal B    \to \N_0;~ \delta(x,b)\in  D_{b}
\\
\mbox{ subject to~~~} &
    \sup_{\p \in \cP} \E^\p\left[\sup_{b\in \mathcal B}L_{b}(0,\delta(X,b))\right] \le \gamma.\end{aligned}
\end{align} 
The objectives to minimize in \eqref{eq:PG-3}
and \eqref{eq:PG-4} are not important in our discussions, and they can be replaced by any objective that is decreasing in $\delta(X,b)$. The risk constraint is our focus. 

In the setting of this section, 
the p-value-based decision $\delta_P$ and  the e-value-based decision $\delta_E$ are given by \eqref{eq:PG-p} and \eqref{eq:PG-e} while treating $b$ as an additional input variable; so we write them as $\delta_P(X,b)$ and $\delta_E(X,b)$ going forward.
Note that for each $b\in \mathcal B,$ the definition of $\delta_E$ guarantees
\begin{align}
\label{eq:PG-e-key}
  L_b(0,\delta_E(X,b)) \le\gamma E(X).
\end{align}

 
\begin{theorem}
\label{th:PG}
For the problems \eqref{eq:PG-3} and \eqref{eq:PG-4},  the following hold.
 \begin{enumerate}[label=(\roman*)]
    \item 
The e-value-based decision $\delta_E$ for any $E(X)$    satisfies the risk constraints  in~\eqref{eq:PG-3}
and  \eqref{eq:PG-4}.
\item  
Any decision rule $\delta^*$ satisfying the risk constraint in~\eqref{eq:PG-3} can be written as $\delta^*(X,B)=\delta_E(X,B)$  for some e-variable $E(X)$.
\item Any decision rule $\delta^*$ satisfying the risk constraint in~\eqref{eq:PG-4} satisfies  $\delta^*(X,b)\le \delta_E(X,b)$ for some e-variable $E(X)$ and all $b\in \mathcal B$, meaning that $\delta^*$ is weakly improved by utilizing e-variables. 
\item 
The p-value-based decision $\delta_P$   does not satisfy the risk constraint in~\eqref{eq:PG-3} or  \eqref{eq:PG-4} in general when $D_b$ has at least 3 elements or $\mathcal B$ has at least 2 elements. 
 \end{enumerate}
Implied by (ii) and (iii), if an optimizer $\delta^*$ of \eqref{eq:PG-3}  or \eqref{eq:PG-4} exists, then it is dominated by some $\delta_E$.

\end{theorem}
\begin{proof}[Proof.]
(i): 
This follows directly follows from \eqref{eq:PG-e-key}, which yields, for $\p\in \cP$, 
 $$ \E^\p\left[ \sup_{b\in \mathcal B} L_{b}(0,\delta_E(X,b))\right]\le\E^{\p}[ \gamma E(X)] \le \gamma,
$$
which guarantees
the risk constraints in  \eqref{eq:PG-3} and \eqref{eq:PG-4}.

(ii):  
Take $E (X)= L_{B}(0,\delta^*(X,B))/\gamma$. 
Then $E(X)$ is an e-variable because $\delta^*$ satisfies the risk constraint in \eqref{eq:PG-3}.
We have 
 $$
\delta_E(X, b) = \max \left\{ d \in D_{b }:   L_{b}(0,d) \le L_{B}(0,\delta^*(X,B))  \right\}.
 $$
This gives $\delta_E(X,B) = \delta^*(X,B)$ under the assumption of a strictly increasing $d \mapsto L_b(0,d)$. 

(iii): Take $E (X)= \sup_{b\in\mathcal B} L_{b}(0,\delta^*(X,b))/\gamma$. 
Then $E(X)$ is an e-variable because $\delta^*$ satisfies the risk constraint in \eqref{eq:PG-4}.
We have 
 $$
\delta_E(X, b) = \max \left\{ d \in D_{b }:   L_{b}(0,d) \le \sup_{b'\in\mathcal B} L_{b'}(0,\delta^*(X,b'))  \right\}.
 $$
This gives $\delta_E(X,b) \geq \delta^*(X,b)$ under the assumption of a strictly increasing $d \mapsto L_b(0,d)$.

(iv): Counter-examples are provided in Example~\ref{ex:p-hack3} below.
\end{proof}

\begin{example}[The risk constraint when $|D_b|\ge 3$ or $|\mathcal B|\ge 2$]
\label{ex:p-hack3}
 
We only show how $\delta_P$ breaks the constraint in \eqref{eq:PG-3}, 
as that in \eqref{eq:PG-4} is more strict. 
Fix $\p\in \cP$ and let $P(X)$ be uniformly distributed on $[0,1]$ under $\p$.
First, consider   $\mathcal B=\{1,2\}$ and $D_b=\{0,1\}$ for $b\in \mathcal B$. Set $L_1(0,1)=\ell_1$ 
and $L_2(0,1)=\ell_2$ with $\ell_1<\ell_2$, let $B=2$ when $P(X)\le \gamma /\ell_2$.
We can compute
\begin{align}
\label{eq:PG-p-cal}
\begin{aligned}
&\E^\p[L_{B}(0,\delta_P(X))]
\\&= \ell_1 \p(\delta_P(X)=1;~B=1) 
+\ell_2  \p(\delta_P(X)=1;~B=2) 
\\&
=\ell_1  \p(\gamma/\ell_2 \le P(X)\le \gamma/\ell_1  ) +\ell_2 \p(P(X)\le \gamma/\ell_2 ) 
\\&
=\gamma \left( 1-\frac{\ell_1}{\ell_2}\right)  + \gamma > \gamma.
\end{aligned} 
\end{align}
For an example in case $D_b=\{0,1,2\}$ and $\mathcal B=\{1\}$, set $L_1(0,1)=\ell_1$ 
and $L_1(0,2)=\ell_2$ with $\ell_1<\ell_2$, as in Example~\ref{ex:p-hack2}.
Following the same calculation as in \eqref{eq:PG-p-cal}, we have
\begin{align*}
\E^\p[L_{1}(0,\delta_P(X))]&= \ell_1 \p(\delta_P(X)=1) 
+\ell_2  \p(\delta_P(X)=2) 
\\&
=\ell_1  \p(\gamma/\ell_2 \le P(X)\le \gamma/\ell_1  ) +\ell_2 \p(P(X)\le \gamma/\ell_2 )  \\& > \gamma.
\end{align*} 

\end{example}


Although the p-value-based decision does not control the risk level  in \eqref{eq:PG-3}, it is still often employed (misused) in real-world situations; see  Examples~\ref{ex:p-hack} and~\ref{ex:p-hack2}. Sometimes the misuse is implicit or unconscious, as the decision maker may not realize that the decision problem has changed due to the collected data, or the fact that there are non-binary decisions (or both). 
Therefore, when using p-value-based methods such as the ones described in Examples~\ref{ex:p-hack} and~\ref{ex:p-hack2}, one should be cautious and mindful about what type of statistical guarantee they actually provide. 
This problem does not exist for the e-value-based rule.


To summarize the message from this section, when decisions are non-binary or the decision problem depends on the data, e-value-based decision rules are needed to satisfy the risk control, and  the optimizers are based on e-values. 

In the simple case of binary decisions,
the frameworks of Sections~\ref{sec:2-posthoc} and~\ref{sec:post-hoc-decisions}
can be connected via setting $\gamma=1$, 
and $L_{b}(0,\delta) = \id_{\{\delta\le b\}} /b$ with $D_b=\mathcal B=(0,1)$, where we slightly generalize the assumption $D_b\subseteq \N_0$; the latter assumption was needed to make sure  that the maximum in \eqref{eq:PG-e} is attained. 
A post-hoc p-variable $\delta$ satisfies the risk constraint for all data-dependent $b$. 
  Section~\ref{sec:2-posthoc} does not specify the loss function under the alternative hypothesis.

\section{Optional continuation of experiments} 
\label{sec:c4-optional}

Suppose that a scientist tests whether a coin is fair (the null hypothesis) or not. It is anticipated that if the coin is not fair, the bias could be around $\theta=0.6$ (the alternative hypothesis). A short calculation shows that to guarantee type-I and type-II errors at $0.05$, the optimal (by the Neyman--Pearson lemma) likelihood ratio test would toss the coin 280 times and reject the null if at least 154 heads appear. However, on performing those coin tosses, suppose that they do not reject the null. This could happen by chance under the alternative hypothesis (a type-II error), or because the scientist had overestimated the strength of the alternative. Perhaps, they should have planned for a coin bias of $\theta=0.55$, but it is too late: Such thoughts are futile, because the type-I error budget of the test has been exhausted and the null was not rejected. However, a conceptually different approach using e-values is possible, as explained next. 

The original test can be rewritten as rejecting the null when the p-value from the likelihood ratio test drops below $0.05$.
Suppose they had instead summarized the data as an e-value from the very beginning; the likelihood ratio would be a natural one in this case. 

Based on the realized value of $E_1$, the scientist (or another scientist) is permitted to toss the coin some more times, say collecting 50 more tosses. Suppose they calculate an e-value $E_2$, based only on the fresh batch of data. In this case, note the important property of conditional expectation $\E^\p[E_2 \mid E_1] = 1$ under the null hypothesis. Then, $E^{(2)} = E_2  E_1$ is also a valid e-value for the same null, which can be verified by the law of iterated expectation.
Note that both the option of running the second experiment and the design of that experiment can depend on the first batch of data, and thus it is a post-hoc decision.

The scientist may still not be happy. Perhaps $E^{(2)}$ is still promising and not conclusive. They have the right to continue their experiment and collect one more batch of data, summarize it as an e-value $E_3$, and calculate $E^{(3)} = E_3  E^{(2)} $, which is itself an e-value for the same null.

The scientist is allowed to repeat this process indefinitely. They can continually monitor the evidence process $E^{(1)},E^{(2)},\dots$ and can reject the null as soon as their current level of evidence exceeds $1/\alpha$. This controls type-I error by  Ville's inequality (Chapter~\ref{chap:eprocess}). Alternatively, they may simply stop at any data-dependent stopping time $\tau$ and report their evidence against the null as $E^{(\tau)}$, which is an e-value due to the optional stopping theorem (Chapter~\ref{chap:eprocess}). When the size of each batch is reduced to one, one ends up with
a fully sequential test, in this case Wald's sequential probability ratio test that we revisit in Section~\ref{sec:6-wald}. 

However, despite the clear relationship to optional stopping and sequential testing above, the two settings (the one discussed above, and the one in Chapter~\ref{chap:eprocess}) are motivated differently. Above, the scientists who may be continuing the study may be different from the first scientist that collected the first batch of data. This could happen in settings where data collection is more expensive than a coin toss; for example, the first scientist may simply have run out of money to collect more data, or may have given up if the evidence looked far from promising. Thus, we call the above setting as ``optional continuation'', as opposed to the setting in Chapter~\ref{chap:eprocess}, which we call ``optional stopping''.

Such optional continuation is particularly prevalent in the meta-analysis of clinical trials. Suppose that Country~B is hesitant to run an expensive trial on some drug (e.g., a vaccine for a pandemic). Country~A goes ahead with such a trial, and the findings are impressive, but Country~B deems it still too risky to launch the drug without its own trial. Piggybacking on the prior success of the drug, now Country~B runs is own trial, and appends its evidence to that collected by Country~A. Now Country~C, which perhaps was not even considering the vaccine, may get convinced to run its own trial. The body of evidence against the null hypothesis (the drug being a placebo) is growing as more countries optionally continue this ``living meta-analysis''. At this meta-level across countries, it does seem like a sequential meta-analysis is proceeding, but not one that was planned or anticipated in advance. Indeed, each country may have only intended to collect its own batch of data. As long as the evidence from each trial was summarized as an e-value, the overall evidence can easily be updated and sequentially monitored. For example, while wealthier countries may launch a vaccine program if the e-value exceeds 20, there may be several poorer countries that will wait till the overall e-value exceeds 100 before launching a vaccine program without their own dedicated trial.

There are three important caveats:
\begin{itemize}
    \item No participants can discard or hide their findings, even if their findings lead to a small e-value (less than $1$). After the decision of running the experiment, the resulting e-value needs to be incorporated into the evidence. 
    \item Building on the above point, it is clear that a single scientist with an e-value equal to zero can permanently destroy the evidence accumulated thus far. This can be accounted for in two ways. Perhaps the scientists agree beforehand on using (e.g., log-optimal) e-values that are not equal to zero by design, or they agree to an accumulation of evidence of the form $\prod_{i=1}^n (1+E_i)/2$, so that each individual batch of data cannot lose more than half of the existing evidence; we can of course generalize this to $\prod_{i=1}^n (\lambda_i + (1-\lambda_i) E_i)$, where $\lambda_i\in [0,1]$ must be specified before seeing $E_i$. Either way, both the e-value being used in the $i$-th round and how to combine the $i$-th e-value with the previous ones must be specified before seeing the $i$-th batch of data. 
    \item The decision to continue the collection of data at the second stage cannot be based on the first p-value not being smaller than $0.05$ (because testing at the first stage with this p-value has used up all the type-I error budget). Either one commits to using a p-value with a single batch of data, or one commits to using an e-value with the possibility of extending the analysis.
\end{itemize}

It is worth emphasizing that it appears difficult to merge p-values that have been collected using the above protocol. While it may appear to be the case that second batch of data (say the trial data from Country B) is independent of the first batch of data (the trial data from Country A), it would be a mistake to summarize these as p-values and attempt to merge them using a standard method for combining independent p-values, like Fisher's rule. The reason is that the second p-value is technically not independent of the first, despite being calculated on ``separate'' data. The second p-value would not even have existed had the first p-value not been promising. So the \emph{number} of p-values is itself random, and the p-values are \emph{sequentially dependent} in a complex way.

Such post-hoc decisions to run trials, or to make decisions based on accumulated data from different trials, is a key feature of e-value-based evidence collection and analysis.

\section{Incentive-compatible principal-agent hypothesis testing} 

\label{sec:c4-PAHT}

In this section, we consider
another setup where e-values play a key role in enabling decision making with constraints.  

Consider two parties, the \emph{principal} and the \emph{agent}.
As a running example for the section, the principal can be thought of as a regulator, like a government agency that grants licenses, while the agent can be thought of as the experimenter, like a pharmaceutical company developing and testing various drugs.
This setting has two unique features:
\begin{itemize}
    \item The principal and agent have different amounts of information about the true state of the world/nature. For example, the agent may have a clearer picture of whether a drug is effective or not, or at least how likely it is to be effective.
    \item The principal and agent have different, typically conflicting, incentives. For example, the principal may not want ineffective drugs (placebos) to be falsely marketed as cures for some disease, while the agent may like to be granted licenses to make money (whether or not the drugs are effective).
\end{itemize}

Usually, a pharmaceutical company declares to the regulator that it is performing a clinical trial in order to determine the efficacy of its drug. For simplicity, let us ignore the multi-phase aspects and other details of such trials, and assume that if the p-value is smaller than some threshold (such as $0.05$), then the regulating agency grants permission (a \emph{license}) to the company to market the drug and profit from its sales.

The problem with the usual hypothesis testing paradigm is that it is not \emph{incentive compatible}. A testing paradigm is incentive compatible if it discourages malicious agents --- those who know the null hypothesis is true, for example, those who have a placebo drug and know it --- from attempting to run such a trial. Indeed, suppose that a trial costs 10 million dollars, then consider the following calculation by a malicious agent:
\begin{itemize}
    \item If the malicious agent estimates that they can make 100 million dollars by selling the drug, then they will not run the trial, because they will expect to lose 5 million dollars (100 million times $0.05$, minus 10 million).
    \item If their anticipated profit is one billion dollars, then a malicious agent will certainly run such a trial, since their expected profit is now 40 million dollars. 
\end{itemize}

If there are twenty times as many ineffective drugs as effective ones, and 5 percent of ineffective drugs are approved, then a quick calculation shows that at least half of all approved drugs will be ineffective.  

Another problem with the above protocol is that the profit margin available to a malicious agent is only limited by the market size, and not by the strength of their statistical evidence. It seems reasonable to suggest a modification in which the more (less) evidence there is that a drug is effective, the more (less) an agent should be able to sell its drugs (or profit from its sale). We assume for simplicity that the profit/sales is determined by the regulator,  not by the true effectiveness of the drug (which will be unknown to the market when sales begin, at least for a reasonable period of time, during which the market will assume that the drug is effective because it was approved by the regulator).

Thus, we may be interested in designing other interaction protocols between the principal and the agent, which have two natural properties: 
\begin{itemize}
    \item the protocol should be incentive compatible, in the sense that malicious agents should not have any incentive to participate.
    \item if an agent participates in the protocol, the sales/profit license increases with the strength of the statistical evidence.
\end{itemize}

We now show how this can be easily achieved with the use of a \emph{statistical contract} based on e-values. We assume that the effectiveness of a drug can be captured by a quality parameter $\theta \in \Theta$ which is known to the agent, while the regulator deems that drugs with $\Theta_0 \subsetneq \Theta$ characterizes the placebos and would prefer not to have them be sold in the marketplace. 
The statistical contract is a protocol in which the agent follows the following steps in order:
\begin{enumerate}
    \item Choose a nonnegative license function $f$ from a menu of options $\cF$. 
    \item Incur a fixed cost $C$ to run a trial in order to collect data $Z \lawis \p_\theta$.
    \item Receive a license to profit at most $L = f(Z)$.
\end{enumerate}
Under such a statistical contract, it is clear that the expected profit for the agent is $\E^{\p_\theta}[f(Z)] - C$, and we may for simplicity assume that they will participate in the trial if and only if this is positive. Thus, the regulator would like to ensure that this quantity is negative if $\theta \in \Theta_0$. 

The above setup leads to the following result.
\begin{proposition}
    The statistical contract is incentive compatible if and only if $
    f(Z)/C
    $
    is an e-variable for 
    $\{\p_{\theta}:\theta\in\Theta_0\}$
    for every $f\in \cF$.
\end{proposition}
The proof is immediate, from realizing that the statistical contract is incentive compatible if and only if $\E^{\p_\theta}[f(Z)] \leq C$ for all $\theta \in \Theta_0$.

As a particular example, if the evidence is going to be summarized as a p-value $P$, then it may be tempting to promise a license that profit margin that is inversely proportional to $P$. However, this would lead to an infinite expected profit under the null, and is not incentive compatible. In contrast, it is incentive compatible to award a license of $C/(2\sqrt{P})$, whose expected value is no more than $C$ under the null. Indeed, $p \mapsto 1/(2\sqrt{p})$ is a p-to-e calibrator.

The setting in this section is connected to the one in Section~\ref{sec:post-hoc-decisions} in the following way. There are two decision makers, making different types of decisions. The regulator is choosing the menu $\cF$ to offer, and then the agent must choose whether to participate in the contract (run a trial). 
\begin{enumerate}
    \item Going first, the regulator has no particular objective to maximize, but they would like to meet their constraint of incentive compatibility, meaning that the menu must discourage malicious agents from participating. This is, as we have seen, achieved by using a menu composed of (scaled) e-values. Of course, if they hope to encourage agents with effective drugs to participate, the menu $\cF$ should not be restricted in any other way: it should be the set of all (scaled) e-values. Assume this is how the regulator plays.
    \item Then, the agent is choosing whether to participate in the statistical contract, given that their eventual payout is data-dependent (that is, they face a data-dependent loss function). They do not need to worry about incentive compatibility (that constraint has already been enforced by the regulator's choice of $\cF$), and they have to simply decide to participate or not before running the trial. Assuming that they know the ground truth and  there exists some $f \in \cF$ that is profitable (in expectation under the alternative),   their decision is clear: they will not participate if they are malicious, and otherwise they may  (depending on their assessment of the efficacy of the drug, and thus of their expected profit). 
\end{enumerate}

In summary, the decision problem for both agents considerably simplifies via a smart use of data-dependent loss functions based on the menu $\cF$ of scaled e-values.



\section*{Bibliographical note}
Section~\ref{sec:2-posthoc} combines results that are explicitly or implicitly found in~\cite{wang2022false,xu2022post,koning2023markov}, in the context of post-hoc $\alpha$ for testing or controlling the false discovery rate or false coverage rate. 
A notion of post-hoc p-values is defined in
 \cite{koning2023markov} through \eqref{eq:c4-nick},  
which is slightly different from our definition and coincides with the reciprocal of e-values. 
The sufficient condition in Remark~\ref{rem:p-to-e-random-alpha} was proved by~\cite{blanchard2008two}.

Section~\ref{sec:post-hoc-decisions} is based on \cite{Grunwald22}, where some omitted details, admissibility of the e-value-based decision rules, and many examples can be found, but our results and presentation differ slightly. 

\cite{GrunwaldHK19} coined the term optional continuation in Section~\ref{sec:c4-optional}. This theme was further explored by~\cite{shafer2023improving}, who argues that it is easier to justify optional continuation using game-theoretic probability than using measure-theoretic probability (the latter relying on the Ionescu-Tulcea theorem for justifying extending the probability space from its initial specification). \cite{SchureG22} detailed an application of these ideas to meta-analysis.

The incentive-compatible principal-agent  hypothesis testing problem in Section~\ref{sec:c4-PAHT} was formulated and studied by~\cite{bates2022principal}. 

\part{Core Ideas}
\label{part:core}

\chapter{The universal inference e-value}
\label{chap:ui}

Universal inference is a general method (or set of methods) to construct an e-variable (or an e-process) for a composite null $\cP$ against a composite alternative $\cQ$ under no \emph{regularity conditions}. The simple takeaway message is: if we can efficiently calculate the maximum likelihood under the null, then we can construct an e-variable, and hence a test.

\section{Motivation: Irregular models}

The development of universal inference was motivated by the difficulty of constructing tests for \emph{irregular problems}, for which we often know no (even asymptotically) valid test. 

Without getting into formal definitions, because they will not be relevant later, irregular problems are usually plagued by many issues all at once: the maximum likelihood estimator may not be asymptotically normal, the Fisher information matrix may not always be invertible (especially on the boundary between null and alternative), the bootstrap may not have provable guarantees, Wilks' theorem may not hold (the log-likelihood ratio may not be asymptotically chi-squared), and so on. We provide two motivating examples here, for which, to the best of our knowledge, universal inference provides the first valid test.

We begin with the (conceptually, not technically) simple example of testing whether our data is Gaussian (null) or is drawn from a mixture of two Gaussians (alternative). Defining $p$ as the density of $(1-\lambda) \mathrm{N}(\mu_1,1) + \lambda \mathrm{N}(\mu_2,1)$, suppose we want to test the null hypothesis represented by the set $\{(\mu_1, \mu_2, \lambda): \lambda = 0\}$, which also equals $\{(\mu_1, \mu_2, \lambda): \mu_1 = \mu_2\}$. The alternative set of distributions is given by triplets $(\mu_1, \mu_2, \lambda)$ that are unrestricted (except for excluding the null). For the null, $\mu_2$ is not identified, and for the alternative, $\lambda$ is not identified. This is a prototypical example of an irregular testing problem, with a composite null and alternative, and it is hard to find classical asymptotic procedures that can provably control the type I error. 
In contrast, the methodology presented here does work without any further assumptions, and has straightforward extensions to 
multidimensional settings, mixtures of more than two components, non-Gaussian distributions, and so on. 

A second example involves testing shape constraints such as log-concavity. A density $p$ on $\R$ is log-concave if
$p= e^g$ for some concave function $g$.
Consider testing
$H_0: $ $p$ is log-concave versus
$H_1: $ $p$ is not log-concave. The method presented in this chapter will be applicable to this problem, while we know of no other nontrivial test. Further, our method is also computationally efficient because the maximum likelihood log-concave density can be solved in polynomial time using a (relatively) efficient optimization algorithm. 

As mentioned earlier, universal inference can be summarized as
\begin{quote}
    a \emph{reduction} from hypothesis testing to maximum likelihood estimation under the null.
\end{quote}
Here, the word ``reduction'' means that the original problem of testing is reduced to an arguably simpler problem, that of (maximum likelihood) estimation. 

We make two assumptions in this chapter:
\begin{enumerate}
    \item $\cP,\cQ$ have a common reference measure $\leb$, so we will associate distributions with their densities with respect to  $\leb$. 
    \item The data $X$ are represented by a vector $X^n =(X_1,\dots,X_n)$ with iid entries. We identify distributions $\p \in \cP$ and $\q \in \cQ$ by the corresponding densities of a single entry, denoted $p(x)$ or $q(x)$. 
    For instance, we use $p\in \cP$ and $\p \in \cP$ interchangeably.
\end{enumerate}

\section{Split likelihood ratio e-variable}
\label{sec:4-split}

As suggested by its name, the key idea of the split likelihood ratio method is based on sample splitting. We divide the $n$ data points into two groups $D_0$ and $D_1$ such that $D_0$ and $D_1$ are independent. For example, take $D_0 = \{ 
X_i : i \in I \}$ and $D_1 = \{ X_i : i \notin I  \}$ where $I$ is a random subset of $[n]$ independent from $X^n$.
Even more precisely, for each point $X_i$, one could flip an independent coin (of some prespecified bias) to decide whether to put it in $D_0$ or in $D_1$.
One may think of these as being roughly equally sized sets as a default choice (i.e.\ the coin bias equals half). The validity of the e-variable does not depend on the relative sizes, just on their independence, but the e-power does. We also note that the split does not have to be random; one could define $D_0=\{X_1,\dots,X_m\}$ for some $1 < m < n$ and $D_1 = \{X_{m+1},\dots,X_n\}$.

The two sets have complementary roles. First, one uses the data in $D_1$ to pick an alternative $\hat q_1 \in \cQ$.    This choice is unrestricted: The maximum likelihood estimator, or a Bayes estimator (a posterior mean or mode based on some prior), 
or a robust one; the choice does not affect validity, but it does affect e-power. Our default recommendation in practice is to use a Bayes estimator with a prior that puts mass everywhere in $\cQ$; this comes with certain asymptotic optimality properties, not discussed here. No particular properties of Bayesian inference are important here; one can simply view it as a smoothed maximum likelihood estimator where the smoothing particularly helps for small sample sizes.

The second step is to use $D_0$ to calculate the maximum likelihood estimator under the null, denoted $\hat p_0 = \argmax_{p \in \cP} \prod_{i \in D_0} p(X_i)$.

The final step uses $D_0$ again to calculate the likelihood ratio between $\hat q_1$ and $\hat p_0$:
\begin{equation}\label{eq:split-lrt}
E = \prod_{i \in D_0} \frac{\hat q_1(X_i)}{\hat p_0(X_i)}.
\end{equation}
We call this the \emph{split likelihood ratio e-variable}. When the roles of the two data splits are swapped, the split likelihood ratio $E^\text{swap}$ is recomputed, then $E^\text{cross-fit} = (E + E^\text{swap})/2$ is called the \emph{cross-fit likelihood ratio e-variable}.
When combined with (some variant of) Markov's inequality (Chapter~\ref{chap:markov}) to yield a level-$\alpha$ test, we call it the \emph{split likelihood ratio test} (split LRT) or \emph{cross-fit LRT}. 

\begin{theorem}
\label{thm:split-lrt}
    The split likelihood ratio statistic defined in~\eqref{eq:split-lrt} is an e-variable for $\cP$. 
\end{theorem}
\begin{proof}[Proof.]
By definition of $\hat p_0$, it holds that $\prod_{i \in D_0} \hat p_0(X_i) \geq \prod_{i \in D_0} p(X_i)$ for any $\p \in \cP$,  and thus
\[
\E^{\p}\left[ \prod_{i \in D_0} \frac{\hat q_1(X_i)}{\hat p_0(X_i)}  \mid D_1 \right] \leq \E^{\p}\left[ \prod_{i \in D_0} \frac{\hat q_1(X_i)}{p(X_i)}  \mid D_1 \right] = 1,
\]
where the final equality holds simply because $D_0$ and $D_1$ are independent, and conditional on $D_1$, $\hat q_1$ is a fixed density, and we recall from Section~\ref{sec:LR-e-variable} that a likelihood ratio of any alternative to $p$ is an exact e-variable for $p$. The proof is completed by taking a further expectation with respect to $D_1$ and the tower property of conditional expectation.
\end{proof}

Recalling the discussion from Section~\ref{sec:3-mix-plugin} about mixture and plug-in methods for handling composite alternatives, one may notice that the split likelihood ratio e-variable is based on the plug-in method. However, the techniques for handling composite nulls and composite alternatives are modular, in the sense that one can swap in the mixture method for the plug-in method in the numerator of these statistics, while leaving the denominator identical.

\begin{remark}[Mixture universal inference]
\label{rem:mix-plugin-ui}
A variant of the split LRT uses $D_1$ to pick a distribution $\nu$ over $\cQ$, and then defines
\[
E = \int_{\cQ} \prod_{i \in D_0} \frac{ q(X_i)}{\hat p_0(X_i)} \nu(\d q).
\]
This is also an e-variable, which can be seen with an almost identical proof to above, by noting that integrals are just averages, and averages of e-variables are e-variables (more generally treated in Chapter~\ref{chap:multiple}).
\end{remark}

\section{Subsampled likelihood ratio e-variable}

The split likelihood ratio statistic has an extra source of variability (beyond the randomness of the data) that was introduced algorithmically by sample splitting. It is possible to effectively remove this extra randomness by averaging, as we describe next. We could swap the roles of $D_0$ and $D_1$, recalculate the e-variable using \eqref{eq:split-lrt} (call it $E'$), and use $(E+E')/2$ as the e-variable. We call this the \emph{crossfit likelihood ratio e-variable}. A better idea is to recalculate the split likelihood ratio e-variable $B$ times, each time involving an identical and independent random split of $X^n$ (these can be done in parallel). This yields e-variables $E^{(1)},\dots,E^{(B)}$ that can be combined by averaging:
\begin{equation}\label{eq:subsampled-ui}
\bar E^{(B)} = \frac{E^{(1)} + \dots + E^{(B)}}{B}.
\end{equation}
We call this the \emph{subsampled likelihood ratio e-variable}. One can show that this (usually strictly) improves the e-power of the method, as proved below. This relies on noting that $E^{(1)},\dots,E^{(B)}$ have the same e-power (since they are identically distributed).

\begin{proposition}
[Improving e-power by averaging]
\label{prop:epower-avg}
    Let $E^{(1)},\dots,E^{(B)}$ be arbitrary e-variables, and let $\bar E^{(B)}$ denote their average. For any $\q$, 
    \[
    \E^\q[\log \bar E^{(B)}] \geq \frac{1}{B}\sum_{b=1}^B \E^\q[\log E^{(b)}],
    \]
    assuming that the right-hand side is well-defined.
    In particular, the e-power (against any $\q$) of the subsampled likelihood ratio statistic~\eqref{eq:subsampled-ui}  is at least as large as that of the split likelihood ratio statistic~\eqref{eq:split-lrt}.
\end{proposition}

This proof follows immediately from concavity of the logarithm function.  
Further, as $B\to\infty$, $\bar E^{(B)}$ converges to a fixed random variable (see Fact~\ref{fact:exch}). 

Proposition~\ref{prop:epower-avg} can be  viewed as a form of Rao--Blackwellization. The split likelihood ratio statistic utilized an external source of randomness (``algorithmic randomness'') by data splitting, thus artificially introducing a source of variance, which is then removed by averaging.

Above, $B$ has to be fixed in advance and cannot depend on the data. 
If one wants to obtain a p-value or a test from the subsampled likelihood ratio statistic, then it is possible for $B$ to be data-dependent, as we now describe.


\paragraph{Subsampled LRT.}
For $b \geq 1$, let $\bar E^{(b)}$ denote $(E^{(1)} + \dots + E^{(b)}) / b$. The subsampled LRT rejects the null as soon as any $\bar E^{(b)}$ exceeds $1/\alpha$ for any $b \geq 1$. Thus, the subsampled LRT can be run sequentially to produce $E^{(1)}, E^{(2)}, \dots$ one by one, from which one can calculate their running averages $\bar E^{(1)}, \bar E^{(2)}, \dots$ one be one, rejecting the null at the first time any of these exceeds $1/\alpha$. One can also terminate this procedure at any large, data-dependent $B$, without violating the type-I error guarantee. 

\begin{proposition}[Improving power by averaging]
    The subsampled LRT described above controls type-I error at level $\alpha$ and is at least as powerful as the split LRT procedure that is based on $E^{(1)}$ and Markov's inequality.
\end{proposition}

The second claim about larger power is evident because the first step of the subsampled LRT corresponds exactly to the split LRT. The first claim about type-I error is less obvious, and is a direct consequence of the \emph{exchangeable Markov inequality} introduced next.


\section{Exchangeable Markov's inequality}
\label{sec:4-exch}

We now state the generalization of the Markov inequality under exchangeability.
All statements are made under a probability measure $\p$ that is implicit. 

Recall that $Z_1,\dots,Z_n$ are called exchangeable if the joint distribution of $(Z_1,\dots,Z_n)$ equals that of $(Z_{\pi(1)},\dots,Z_{\pi(n)})$ for any permutation $\pi$ of $\{1,\dots, n\}$. A sequence  $Z_1,Z_2,\dots,$ is called exchangeable if  $Z_1,\dots,Z_n$ are exchangeable for any $n\geq 1$. 

For instance, iid random variables are exchangeable, and so are identical copies of the same random variable. Also, exchangeable random variables have the same marginal distribution and the same expectation.

\begin{fact}
\label{fact:exch}
    For any exchangeable sequence $Z_1,Z_2,\dots$ of integrable random variables,
    their sample mean process $\bar Z=(\bar Z_n)_{n\in\N}$
    given by $ \bar Z_n:=(Z_1+\dots+Z_n)/n$ 
    is a \emph{backward martingale} with respect to the natural backward filtration
    in the sense that 
    $$
    \E^\p[\bar Z_{n} \mid  \mathcal G_{n+1} ]=\bar Z_{n+1}~\mbox{for $n\in \N$}, \mbox{~where $\mathcal G_{n+1}=\sigma(\bar Z_{n+1},\bar Z_{n+2},\dots)$}.
            $$
    Further, $(\bar Z_n)_{n\in\N}$ converges almost surely to an integrable random variable.
\end{fact}

The proof of the martingale property in Fact~\ref{fact:exch} follows by first noting by direct  calculation that $\E^\p[\bar Z_n | \bar Z_{n+1}] = \bar Z_{n+1}$, and the limiting claim is a consequence of the backward martingale convergence theorem of~\cite{doob1953stochastic}. We are now in a position to state a useful nonasymptotic inequality.

\begin{theorem}[Exchangeable Markov inequality (EMI)]
\label{thm:emi}
For any exchangeable sequence $Z_1,Z_2,\dots$ of random variables and any $\alpha > 0$,
\begin{equation}\label{eq:ex-M}
\p\left( \sup_{n\ge 1} \left| \frac1n \sum_{i=1}^n Z_i \right| \geq \frac 1 \alpha\right) 
\leq \alpha \ \E^\p [|Z_1|].
\end{equation}
In particular, if $Z_1,Z_2,\dots$ are e-variables,
\[
\p \left( \exists n \geq 1:  \frac1n \sum_{i=1}^n  Z_i  \geq \frac 1 \alpha \right) \leq \alpha .
\]
The above claims also hold for finite sets of exchangeable random variables. 
\end{theorem}

The proof is short, albeit technical, and we only briefly mention it here without details. The first inequality is trivial. We next note that if $Z_1$ is not integrable, the desired inequality in \eqref{eq:ex-M} is trivially true. If $Z_1$ is integrable, $( \sum_{i=1}^n \left|Z_i \right|/n)_{n \geq 1}$ is a nonnegative backward martingale by Fact~\ref{fact:exch}, and this allows us to infer the desired inequality by invoking a time-reversed variant  of Ville's inequality (Fact~\ref{Ville} in Chapter~\ref{chap:eprocess}). 

A useful corollary of the EMI is as follows. Let $Z_1,\dots,Z_n$ be any set of (potentially nonexchangeable) arbitrarily dependent random variables. Let $\pi$ be a uniformly random permutation of $\{1,\dots,n\}$. Then, for any $\alpha>0$, 
\begin{equation}\label{eq:ex-M2}
\p\left( \sup_{1 \leq m \leq n}  \frac1m \sum_{i=1}^m \left|Z_{\pi(i)} \right| \geq \frac{1}{\alpha} \right)  \leq \alpha \frac{\E^\p[|Z_1| +\dots+|Z_n|]}{n}.
\end{equation}
Here, the original random variables are effectively made exchangeable by the random permutation, thus allowing us to invoke the original inequality.

We state another variant of the inequality. Suppose we take $N$ arbitrarily dependent random variables and put them in a bag. Suppose $Z_{\pi(1)},\dots,Z_{\pi(n)}$ are $n$ samples drawn uniformly at random with or without replacement from this bag.  Then, we have
\begin{equation}\label{eq:ex-M3}
\p\left( \sup_{1 \leq m \leq n}  \frac1m \sum_{i=1}^m \left|Z_{\pi(i)} \right| \geq \frac 1 \alpha  \right)  \leq \alpha \frac{\E^\p[|Z_1|+\dots+|Z_N|]}{N}.
\end{equation}
This holds because the sampling process induces the exchangeability required for~\eqref{eq:ex-M} to be invoked on the otherwise non-exchangeable random variables. 

We end with a result that combines Theorem~\ref{thm:umi} with Theorem~\ref{thm:emi}.

\begin{theorem}[Exchangeable Uniformly Randomized Markov's Inequality (EUMI)]\label{thm:umi+emi}
Let $Z_1, \dots, Z_n$ 
be a set of exchangeable random variables. Then, for any $a \in (0,1)$, 
\begin{align}\label{eq:emi_umi}
    \p\left( Z_1 \geq U/a \text{ or } \exists t \leq n:  \left|\sum_{i=1}^t Z_i/t\right| \geq 1/a \right) \leq a \cdot \E^\p|Z_1|,
\end{align}
where $U$ is a uniform random variable on $[0,1]$ that is independent of $Z_1,\dots,Z_n$. 
\end{theorem}

Its proof is similar to that of EMI, and is omitted for brevity.

\section{Universal confidence sets}
\label{sec:ui-ci}

It is useful, but not necessary, to adopt a parametric framework for what follows. Assume that the set of distributions under consideration can be represented as $\{p_\theta\}_{\theta \in \Theta}$. If we were testing, then the null and alternative hypotheses would correspond to certain subsets $\Theta_0 \subsetneq \Theta$ and $\Theta_1 \subsetneq \Theta$ respectively. But suppose instead that the underlying distribution of $X^n$ is $p_{\theta^*}$ and we want to design a $(1-\alpha)$ confidence set for $\theta^*$, which is a set $C(X^n)$ such that for any $\theta^*\in\Theta$, 
$$\p_{\theta^*}(\theta^* \in C(X^n)) \geq 1-\alpha.$$
Such a confidence set, without requiring any regularity conditions, is immediately given by inverting the universal inference test. To elaborate, we simply test each point null hypothesis $\theta \in \Theta$ at level $\alpha$, and retain those we failed to reject. Let us describe the set more formally by inverting the split LRT.

We first randomly split the data $X^n$ into $D_0$ and $D_1$. Let $\hat \theta$ be any estimator based only on $D_1$. Define
\[
C(X^n) := \left\{ \theta \in \Theta: \prod_{i \in D_0} \frac{p_{\hat \theta}(X_i)}{p_\theta(X_i)} < \frac 1 \alpha \right\}.
\]
Denoting 
\[
\cL_0(\theta) := \prod_{i\in D_0} p_\theta(X_i),
\]
we see that $C(X^n) = \{ \theta \in \Theta: \cL_0(\hat \theta)/\cL_0(\theta) < 1/\alpha \}$.
Then, we have the following guarantee.
\begin{proposition}
    For any $\theta^* \in \Theta$ and any predefined error level $\alpha \in (0,1)$, the set $C(X^n)$ defined above is a $(1-\alpha)$ confidence set for $\theta^*$.
\end{proposition}

The proof is immediate, relying only on the observation that $E = \prod_{i \in D_0} ({p_{\hat \theta}(X_i)}/{p_{\theta^*}(X_i)})$ is an e-variable, or more accurately that $\id_{\{E \geq 1/\alpha\}}$ is a level-$\alpha$ test. Sometimes, for instance for simply estimating a multivariate normal mean with a known covariance matrix, the above set is in closed form. Typically, the set can be made smaller (in terms of expected diameter) by replacing the split LRT by the subsampled LRT.

\subsection*{An example: multivariate Gaussian data with identity covariance}

It is not obvious what the optimal splitting proportion is for split LRT based confidence sets. For iid Gaussian data with identity covariance in $d$ dimensions, these sets are spherical, and one can derive explicit expressions for the proportion $p_0^*$ of the data in $\cD_0$ that minimizes the squared radius $r^2$ of $C(X^n)$.

\begin{theorem} \label{thm:split_p0}
Let $X_1, \ldots, X_n \overset{\mathrm{iid}}{\sim} \mathrm{N}(\theta^*, I_d)$ under $\p$. 
The splitting proportion that minimizes $\E^\p[r^2 ( C(X^n) ) ]$ is 
\begin{align}
p_0^* &= 1 - \frac{\sqrt{4d^2 + 8d\log\left(1/\alpha\right)} - 2d}{4\log\left(1/\alpha\right)}.
\end{align}
\end{theorem}
The proof is omitted.

As $d\to\infty$ for fixed $\alpha$, the optimal split proportion $p_0^*$ converges to 0.5. Alternatively, as $\alpha\to 0$ for fixed $d$, the proportion $p_0^*$ converges to 1, suggesting that one should use nearly all data for the likelihood estimation. This is not an issue for reasonable $\alpha$ levels, though. For instance, at $d=1$, one would need to set $\alpha < \exp(-40)$ to produce an optimal split proportion $p_0^*$ that exceeds 0.90. 

\begin{figure}[t]
\centering
\includegraphics[width=\textwidth]{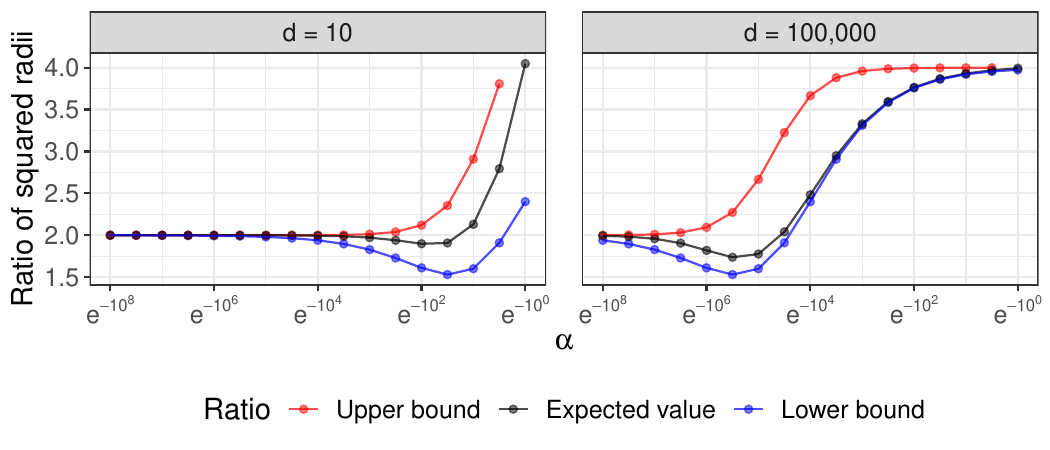}
\caption{Expectation (black), lower bound (blue), and upper bound (red) of $\E\left[r^2(C(X^n)) / r^2(C^\mathrm{LRT}(X^n))\right]$. 
Data points correspond to values at $\alpha = \exp(-10^x)$ for $x$ from 8 to 0 in increments of $-0.5$.}
\label{fig:split_usual}
\end{figure}

Let $C^{\mathrm{LRT}}(X^n)$ denote the confidence set obtained by inverting the classical LRT (whose distribution is known exactly in this Gaussian setting).
With an equal split, $\E[r^2(C(X^n))] = (4/n) \log(1/\alpha) + 4d/n$, and one can easily prove the following two bounds:
\begin{itemize}
    \item As $d \to \infty$, ${\E[r^2(C(X^n))]}/{r^2(C^{\mathrm{LRT}}(X^n))} \to 4$,
    \item As $\alpha \to 0$, ${\E[r^2 (C(X^n))]}/{r^2(C^{\mathrm{LRT}}(X^n))} \to 2$;
\end{itemize}
a more explicit expression is available for any $d$ and $\alpha$, that we visualize in Figure~\ref{fig:split_usual}. One can also derive high probability bounds in the ratio of radii, that deliver the same message as above.

\section{Extensions}

\subsection{Handling nuisance variables} 
Continuing with the theme of estimation,
now consider the case where there is a nuisance parameter in both the null and the alternative. For example, suppose $\theta^* = (\theta_0^*,\theta_1^*)$ but we are only interested in a confidence set for $\theta_0^*$. More generally, suppose we are interested in a confidence set for $\psi^* = g(\theta^*)$ for some given function $g$. Define $g^{-1}(\psi) = \{\theta: g(\theta) = \psi\}$. Now define
\[
B(X^n) = \{\psi: C(X^n) \cap g^{-1}(\psi) \neq \varnothing \}.
\]
Since $C$ is a $(1-\alpha)$ confidence set for $\theta^*$, it is straightforward to check that $B$ is a $(1-\alpha)$ confidence set for $\psi^*$. But the above notation is a bit cumbersome. There is a more straightforward way to write this set. Define the \emph{profile likelihood} on $D_0$ as
\[
\cL^\dagger_0(\psi) := \sup_{\theta: g(\theta)=\psi} \cL_0(\theta),
\]
Then, we can write 
\[
B(X^n) = \left\{ \psi : \frac{\cL^\dagger_0(\psi)}{\cL_0(\hat \theta)} \geq \alpha \right\},
\]
an expression which is possibly easier to work with.

\subsection{Smoothed likelihood}

Sometimes the maximum likelihood under the null may be infinite since the likelihood function is unbounded. One can then use a smoothed likelihood in its place. Consider a kernel $k(x,y)$ such that $\int k(x,y) \d y = 1$ for any $x$. For any density $p_\theta$,
 let its smoothed version be denoted 
\[
\widetilde p_\theta(y) := \int k(x,y) p_\theta(x) \d x,
\]
Note that $\widetilde p_\theta$ is also a probability density.
Let the smoothed empirical density based on $D_0$ be defined as 
\[
\widetilde p_n := \frac1{|D_0|}\sum_{i \in D_0} k(X_i, \cdot ).
\]
Define the smoothed maximum likelihood estimator as the KL projection of $\widetilde p_n$ onto $
\{\widetilde p_\theta\}_{\theta \in \Theta_0}$:
\[
\widetilde \theta_0 := \argmin_{\theta \in \Theta_0} \kl(\widetilde p_n, \widetilde p_\theta).
\]
If we define the smoothed likelihood on the first half of the data $D_0$ as
\[
\widetilde{\cal L}_0(\theta) := \prod_{i\in D_0} \exp \int k(X_i, y) \log \widetilde p_\theta(y) \d y,
\]
then it can be checked that  $\widetilde \theta_0$ maximizes the smoothed likelihood, that is
$
\widetilde \theta_0 = \argmax_{\theta \in \Theta_0} \widetilde{\cal L}_0(\theta).
$
As before, let $\widehat \theta_1 \in \Theta$ be any estimator based on $D_1$. The smoothed split LRT is defined as:
\begin{equation}\label{eq:smoothed-splitLRT}
\text{reject $H_0$ if } \widetilde U_n > 1/\alpha, \text{ where } \widetilde U_n = \frac{\widetilde{\cal L}_0(\widehat\theta_1)}{\widetilde{\cal L}_0(\widetilde\theta_{0})}.
\end{equation}
We now verify that the smoothed split LRT controls type-I error, by simply checking that $\widetilde U_n$ is an e-variable. Indeed, for any fixed $\psi \in \Theta$, we have (where $\E^\theta$ is evaluated with respect to $p_\theta$)
\begin{align*}
\mathbb{E}^{\theta^*}\left[\frac{\widetilde{\cal L}_0(\psi)}{\widetilde{\cal L}_0(\widetilde\theta_{0})}\right] 
&\stackrel{(i)}{\leq} \mathbb{E}^{\theta^*}\left[\frac{\widetilde{\cal L}_0(\psi)}{\widetilde{\cal L}_0(\theta^*)}\right] \\
&=  \prod_{i \in D_0} \int \exp \left( \int k(x,y) \log \frac{\widetilde p_{\psi}(y)}{\widetilde p_{\theta^*}(y) } \d y \right) p_{\theta^*}(x) \d x\\
&\stackrel{(ii)}{\leq}  \prod_{i \in D_0}  \int \left( \int k(x,y) \frac{\widetilde p_{\psi}(y)}{\widetilde p_{\theta^*}(y) } \d y \right) p_{\theta^*}(x) \d x\\
&=  \prod_{i \in D_0}  \int \left(  \frac{\int k(x,y) p_{\theta^*}(x) \d x }{\widetilde p_{\theta^*}(y) } \right) \widetilde p_{\psi}(y) \d y \\
&=  \prod_{i \in D_0}  \int \widetilde p_{\psi}(y) \d y = 1,
\end{align*}
where inequality $(i)$ is because $\widetilde \theta_0$ maximizes the smoothed likelihood, and step $(ii)$ follows by Jensen's inequality.

\subsection{Relaxed, powered and conditional likelihood}

There may be settings
where computing the MLE and/or the maximum likelihood (under the null) is computationally hard.
Suppose one could come up with a relaxation $F_0$ of the null likelihood ${\cal L}_0$ (or of the log-likelihood) in the sense that
\[
\max_\theta F_0(\theta) \geq \max_\theta {\cal L}_0(\theta).
\]
For example, ${\cal L}_0$ may be defined as $-\infty$ outside its domain, but $F$ could extend the domain.
As another example, instead of minimizing the negative log-likelihood which could be nonconvex and hence hard to minimize, 
we could instead minimize a convex relaxation.
In such settings, define 
\[
\widehat\theta^F_0 := \argmax_\theta F_0(\theta).
\]
If we instead define the test statistic
\[
T_n' := \frac{{\cal L}_0(\widehat \theta_1)}{F_0(\widehat \theta_0^F)},
\]
then inference may proceed using $T_n'$ instead of $T_n$ in the split or crossfit LRT. This is simply because of the 
aforementioned property of the relaxation that $F_0(\widehat \theta_0^F) \geq  {\cal L}_0(\widehat \theta_0)$, and hence $T_n' \leq T_n$.

One particular case when this would be useful is the following. Testing the sparsity level
in a high-dimensional linear model corresponds to solving the best subset selection problem, which is NP-hard in the worst case (integer programming). 
There exist well-known quadratic programming relaxations that are far more computationally tractable.

The takeaway message is that \emph{it suffices to upper bound the maximum likelihood under the null in order to perform inference}.

When model misspecification is a concern, inferences can be made robust by replacing the likelihood ${\cal L}$
with the power likelihood ${\cal L}^\eta$ for some $0 < \eta < 1$.
Note that
$$
\mathbb{E}^\theta
\Biggl[
\Biggl(
\frac{ {\cal L}_0(\widehat\theta_1)}
{ {\cal L}_0(\theta)}
\Biggr)^\eta \, \Biggm| \, D_1\Biggr]=
\prod_{i\in D_0}\int
p_{\widehat\theta_1}^\eta (y_i)p_{\theta}^{1-\eta} (y_i) \d y_i \leq 1,
$$
and hence all the aforementioned methods can be used with the power-robustified likelihood as well. (The last inequality follows because the $\eta$-Rényi divergence is nonnegative.)

Our presentation has assumed that the data are drawn iid from some distribution under the null. However, this is not really required (even under the null), and was assumed for expositional simplicity. All that is needed is that we can calculate the likelihood on $D_0$ conditional on $D_1$ (or vice versa). For example, this could be tractable in models involving sampling without replacement from an urn with $M \gg n$ balls. Here $\theta$ could represent the unknown number of balls of different colors. Such hypergeometric sampling schemes obviously result in non-iid sampling, but conditional on the first half of the data (for example how many red, green and blue balls were sampled from the urn in that subset), one can potentially still evaluate the conditional likelihood of the second half of the data and maximize it, rendering it possible to apply our universal tests and confidence sets.

\section{Numerical studies}

\begin{figure}[t]
\centering
\includegraphics[width=\textwidth]{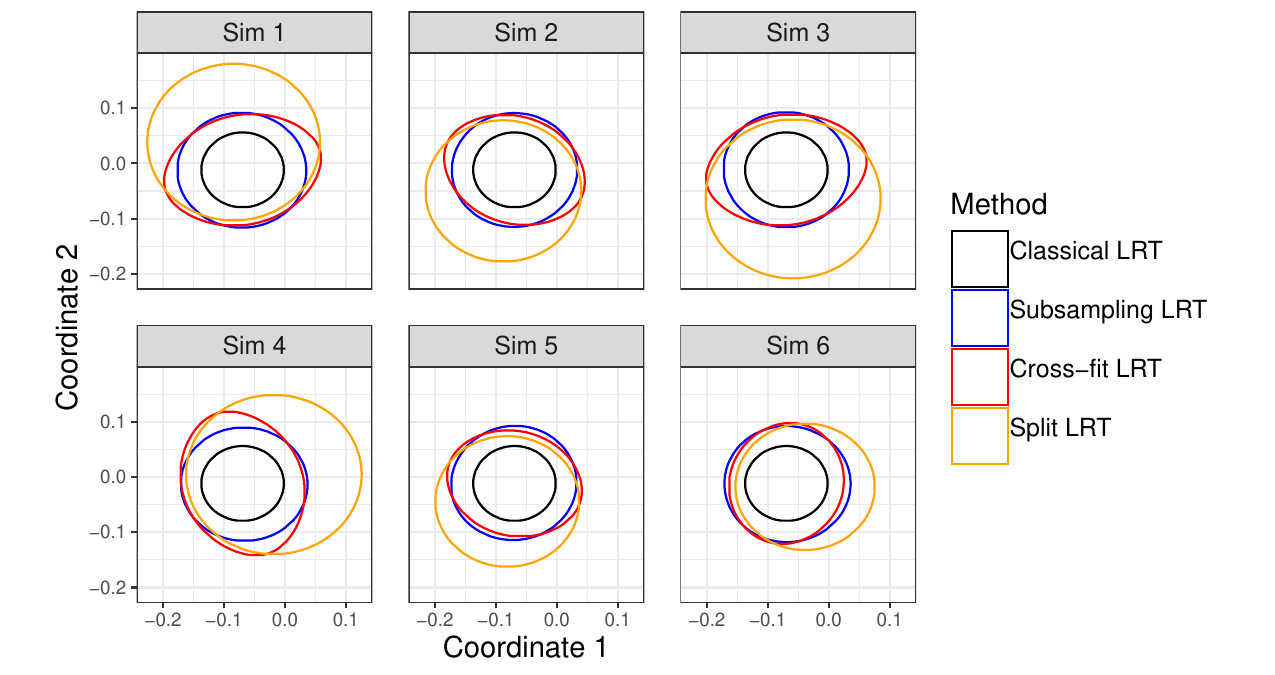}
\caption{Coverage regions of classical LRT (black), subsampling LRT (blue), cross-fit LRT (red), and split LRT (orange) at $\alpha = 0.1$. The six simulations use the same 1000 observations from $\mathrm{N}((0, 0), I_2)$.}
\label{fig:fixed_data_regions}
\end{figure}

\begin{figure}[t]
\centering
\includegraphics[width=\textwidth]{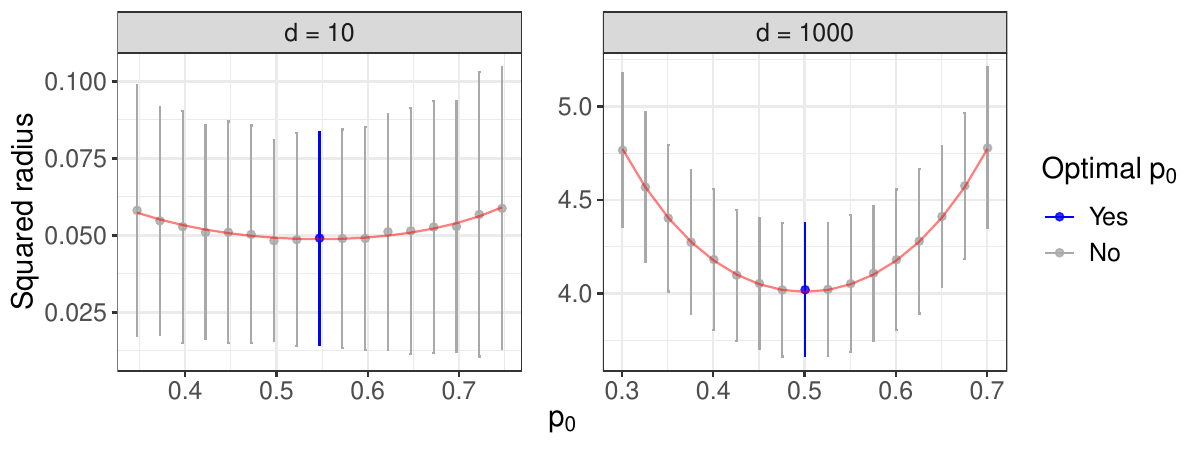}
\caption{Squared radius of multivariate normal split LRT with varying $p_0$. We simulate $X_1, \ldots, X_{1000} \lawis \mathrm{N}(0, I_d)$ and compute the split LRT region at $\alpha = 0.10$ and varying $p_0$. We repeat this simulation 1000 times. At each $p_0$, the circular point is the mean squared radius and the error bar represents the mean squared radius $\pm$ 1.96 standard deviations. Hence, the error bars represent a typical range of squared radius values for each $d$ and $p_0$.  Blue points/lines correspond to $p_0^*$. The red curve is the expected squared radius. See Theorem~\ref{thm:split_p0} proof in the supplement for a derivation of the expected squared radius at $p_0$.}
\label{fig:p0_sq_rad}
\end{figure}

\begin{figure}[t]
\centering
\includegraphics[width=\textwidth]{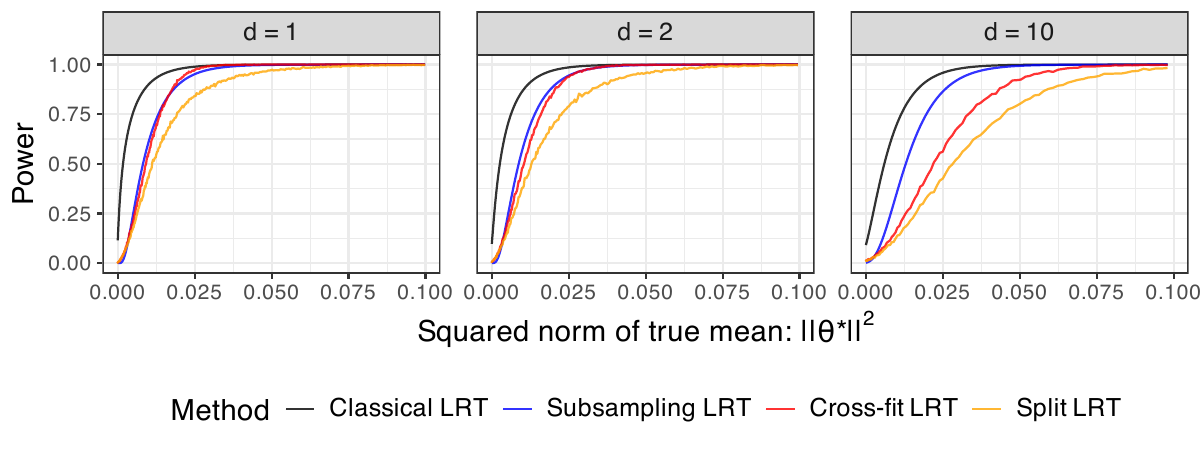}
\caption{Estimated power of classical LRT (black), limiting subsampling LRT (blue), cross-fit LRT (red), and split LRT (orange). We are testing $H_0: \theta^* = 0$ versus $H_1: \theta^* \neq 0$ across varying true $\|\theta^*\|^2$. We use the standard normal cumulative distribution function approximation for the classical and subsampling LRT power calculations, and we use simulations to estimate the cross-fit and split LRT power.}
\label{fig:power_vs_norm_2}
\end{figure}

Since it is easier to visualize confidence sets, Figure~\ref{fig:fixed_data_regions} shows coverage regions of the classical LRT, split LRT, cross-fit LRT, and subsampling LRT at $B = 100$ from six simulations with $\theta^* = (0, 0)$. We generate 1000 observations from $\mathrm{N}(\theta^*, I_2)$ once, and we use this sample for all simulations. Hence, the variation in the split, cross-fit, and subsampling LRTs across simulations is due to algorithmic randomness.

The coverage regions in Figure~\ref{fig:fixed_data_regions} suggest several relationships that can be formally proven. We see that the classical LRT provides the smallest confidence regions. This is not surprising since, even in finite samples, the classical LRT statistic follows a chi-squared distribution under the null in the Gaussian case. It can be proven that the volume of the cross-fit LRT set is less than or equal to the volume of the split LRT set, although the cross-fit set is not entirely contained within the split set. The split and cross-fit approaches both use a single split of the data, but there is a notable improvement from cross-fitting. One can prove that the subsampling set also has less volume than the split LRT set. While any individual split LRT region is guaranteed to be spherical, the subsampling set is not necessarily a spherical region. For large $B$, however, we see that the subsampling region is approximately spherical.

Next, Figure~\ref{fig:p0_sq_rad} shows the average squared radius of the split LRT confidence sets at $p_0^*$ and at surrounding choices of $p_0$. The expected squared radius, given by the red curve, is more sensitive to changes in $p_0$ at higher values of $d$. That is, use of the optimal $p_0^*$ has a greater effect on the split LRT squared radius in higher dimensions, where $p_0^*$ is close to 0.5. 

Figure~\ref{fig:power_vs_norm_2} plots the power of the LRTs against $\|\theta^*\|^2$. Each vector $\theta^*$ has the form $c \mathbf 1.$ To plot the classical and subsampling LRT power, this figure uses the standard normal cumulative distribution function approximation to the non-central $\chi^2$ cumulative distribution function. We use simulations to approximate the power of the split and cross-fit LRTs. For a given value of $\theta^*$, we simulate $n=1000$ observations $Y_1, \ldots, Y_n \lawis \mathrm{N}(\theta^*, I_d)$. We construct split LRT and cross-fit LRT confidence sets from this sample. Then we test whether $\theta = 0$ is in each confidence set. We repeat this procedure 5000 times at each $\theta^*$, and each procedure's  estimated power at $\theta^*$ is the proportion of times that $0 \notin C_n(\alpha)$. 

As we would expect, the power is higher when $\theta^*$ is farther away from $0$. In addition, the classical LRT has the highest power, followed in order by the subsampling LRT, the cross-fit LRT, and the split LRT. Interestingly, at $d=1$ the subsampling and cross-fit LRT have nearly identical approximate power. As $d$ increases, the difference between the subsampling and cross-fit LRT power increases.

\paragraph{Universal Inference for model selection.}

We now compare various methods for the problem
of testing the number
of components in a Gaussian mixture model. 
It is well-known that 
the parametric
family of Gaussian mixtures does not
satisfy the regularity conditions required
for (twice the negative logarithm of) the likelihood ratio statistic to admit
its traditional $\chi^2$ limiting distribution. 
In contrast, the method of universal inference
is valid without any regularity conditions, 
and is therefore a natural candidate for this problem. 

Although the limiting distribution 
of the likelihood ratio statistic is unknown
or intractable for general Gaussian mixtures, 
it admits a simple expression when the underlying
mixing  proportions are known. 
We will assume this to be the case so that
we can use the likelihood ratio test (LRT) as a benchmark, 
but we emphasize that the LRT cannot easily
be used to derive a valid test for more
general Gaussian mixtures, where the universal inference method would remain valid. 

Let
\begin{equation} 
\label{eq:mixture} 
X_1, \dots, X_n \overset{\mathrm{iid}}{\sim} 
0.25\cdot \mathrm{N}(\mu_1,1) + 0.75\cdot \mathrm{N}(\mu_2,1),
\end{equation}
where the only
unknown parameters are $\mu_1,\mu_2 \in \R$, and consider the problem of testing
whether the above mixture has one versus two components, i.e.,
\begin{equation} 
\label{eq:hypothesis_mixture} 
H_0: \mu_1=\mu_2, \quad \text{versus}\quad H_1: \mu_1 \neq \mu_2.
\end{equation}

If $\lambda$
denotes the likelihood-ratio statistic
for these hypotheses, then $-2\log \lambda$ 
admits the limiting distribution $\max(0,Z)^2$,
for $Z \lawis \mathrm{N}(0,1)$, thus a valid level-$\alpha$
test for $H_0$ is to 
$\text{reject if } -2\log \lambda > q_{1-2\alpha},$
the latter being the $1-2\alpha$ quantile
of the $\chi_1^2$ distribution. 
We will refer to this as the LRT test. We compare its numerical
performance to that of the following universal inference tests:
\begin{itemize}
    \item UI:  split LR e-variable with Markov's inequality;
    \item UMI-UI:  split LR e-variable  with the uniformly randomized Markov's inequality; 
    \item SUI: subsampling LR e-variable with Markov's inequality;
    \item UMI-SUI: subsampling LR e-variable  with  uniformly randomized Markov's inequality; 
    \item EMI-SUI: subsampling LR e-variable with exchangeable Markov's inequality;
    \item EUMI-SUI: subsampling LR e-variable with exchangeable, uniformly-randomized Markov's inequality;
\end{itemize}

\begin{figure}[t]
\centering
\includegraphics[width=0.5\textwidth]{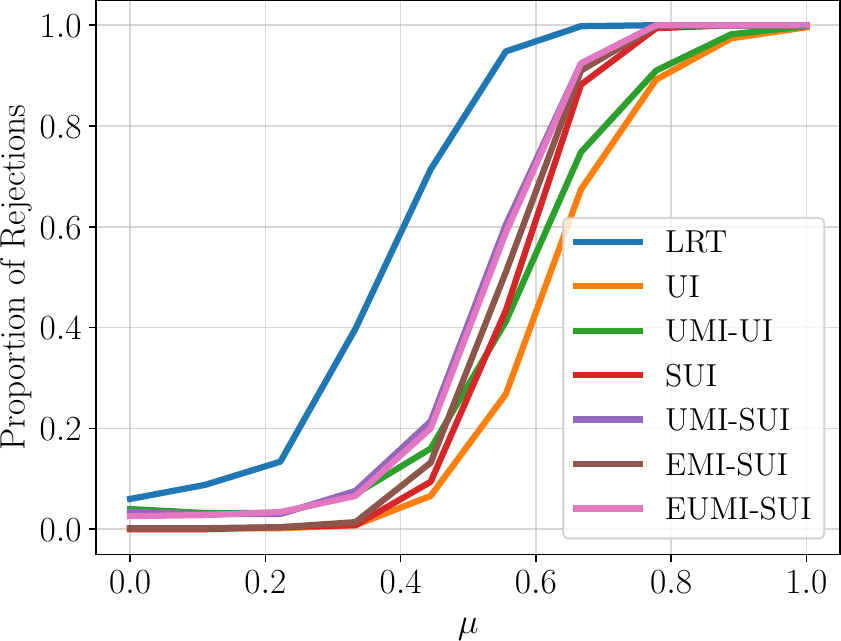}
\caption{\label{fig:ui} 
Empirical power of the seven tests
for the null hypothesis $H_0$ in equation~\eqref{eq:hypothesis_mixture}.
While the benchmark test provided by the LRT
is most powerful, the variants of Universal Inference
based on the EMI, UMI, and EUMI are more powerful
than their MI counterparts. 
}
\end{figure}

Figure~\ref{fig:ui} 
reports the empirical power
of these procedures based on 500 samples
of size $n=500$ from model~\eqref{eq:mixture}, 
under ten equally-spaced values of $\mu := -\mu_1 = \mu_2 \in [0,1]$.  
We observe that the power
of the UMI-UI method
uniformly dominates that of
the original UI method, 
and similarly, the UMI-SUI, EMI-SUI,
and EUMI-SUI methods dominate their SUI counterpart.
The methods UMI-SUI and EUMI-SUI 
exhibit similar performance, and
their absolute increase in power compared
to SUI is on the order
of 15\% for some values of $\mu$. 
Although all methods based on the split
LRT are markedly more conservative
than the LRT, we recall that 
they can be used in arbitrary mixture models,
while the LRT cannot.

\section*{Bibliographical note}

The universal inference method was developed by~\cite{wasserman2020universal}. These authors also showed that for regular models, the universal confidence sets had width scaling as $\sqrt{d/n}$, where $d$ is the dimensionality of the data, as would be expected. 
The exchangeable Markov inequality was first derived in~\cite{manole2021sequential}, but its implications for testing with e-variables was first detailed in~\cite{ramdas2023randomized}.  The smoothed MLE idea is related to the method in~\cite{seo2013universally}. 

Theorem~\ref{thm:split_p0} was proved in \cite{dunn2021gaussian}. The authors find that the universal sets match the standard ones in terms of scaling with respect to $d,n,\alpha$ and the signal to noise ratio, their radii being larger by a small constant factor (smaller than 2). 
\cite{strieder2022choice} also study the optimal splitting ratio for the split LRT, but by examining the power of the test. They derive the (complicated) limiting distribution of the split LRT, and Monte-Carlo simulations show that their proposed split ratio also converges to $1/2$ when testing a fixed number of parameters, but otherwise varies with the dimension of the tested hypothesis.

\cite{dunn2021universal} applied the method to the problem of testing log-concavity, for which no other valid test is known. 
The Gaussian mixture model selection problem was previously studied by
\citep{ghosh1984asymptotic,chen2009hypothesis}, and the limiting distributions we employed can be found in~\cite{goffinet1992testing}. \cite{shi2024universal}  prove that for this problem, the split likelihood ratio statistic
is asymptotically normal with increasing mean and variance, and that it achieves the optimal detection rate of $\sqrt{(\log\log n)/n}$. 
\cite{tse2022note} study some asymptotic approaches to improving the power of universal inference for high-dimensional problems with a large number of nuisance parameters.

The code to reproduce plots is at 
\[
\mbox{\texttt{https://github.com/RobinMDunn/GaussianUniv$\_$sims}}
\]
for the first four figures, and at  
\[
\mbox{\texttt{https://github.com/tmanole/Randomized-Markov} }
\]
for the final figure.

\chapter{The log-optimal e-value and reverse information projection}
\label{chap:numeraire}

While the universal inference e-variable (or e-process) exists in quite some generality and is relatively straightforward to construct, it is only known to be \emph{asymptotically} log-optimal but does not have a finite-sample log-optimality guarantee. It turns out that in even more generality, a log-optimal e-variable exists and is unique. In fact, it always dominates universal inference, rendering the latter method inadmissible. The catch, of course, is that it may not be numerically as easy to construct, but we will be able to characterize it analytically in many cases.

This chapter is dedicated to describing a duality theory that is central to designing log-optimal e-variables. The two main characters in this duality are a special e-variable, called the numeraire, and a special (sub-probability) distribution called the reverse information projection (RIPr). 

\begin{proof}[Setting: $\cP$ versus $\q$.] \renewcommand{\qedsymbol}{}
As elaborated in Chapter~\ref{sec:3-mix-plugin}, there are two central techniques for dealing with composite alternatives --- mixtures and plug-in. For this reason in the current chapter, we consider only a point alternative $\q$ and a composite null $\cP$. 
\end{proof}

Our most general results make no assumptions on $\q$ and $\cP$ whatsoever, not even the existence of a common reference measure or the minimum KL divergence between $\q$ and $\cP$ being finite. However, some of these results have sophisticated proofs that are omitted. We will instead prove some of the results in simpler special cases, making some simplifying assumptions. 

When we use the term e-variable below, we always mean e-variable \emph{for $\cP$}. 
We let $\fE$ denote the set of e-variables:
\[
\fE = \{E \geq 0: \E^\p[E] \leq 1 \text{ for all } \p \in \cP\}.
\]

\section{The numeraire e-variable}
\label{sec:c5-numeraire}

We first give the definition of the numeraire e-variable, which is closely connected to the log-optimal e-variable in Section~\ref{sec:3-p-vs-q}.
\begin{definition} \label{def:numeraire}
A \emph{numeraire e-variable} (or just \emph{numeraire} for short) is a $\q$-almost surely strictly positive e-variable $E^*$ such that $\E^\q[E/E^*] \le 1$ for every e-variable $E \in \fE$.
\end{definition}

Numeraires are unique up to $\q$-nullsets.  Indeed, if $E^*_1$ and $E^*_2$ are numeraires, then the ratio $Y = E^*_2/E^*_1$ satisfies $1 \le 1/\E^\q[Y] \le \E^\q[1/Y] \le 1$ thanks to the numeraire property of $E^*_1$, Jensen's inequality, and the numeraire property of $E^*_2$. Thus Jensen's inequality holds with equality, so $Y$ is $\q$-almost surely equal to a constant which must be one. It follows that $E^*_1$ and $E^*_2$ are $\q$-almost surely equal. In view of this uniqueness, we often speak of `the' numeraire.

In the case of a simple null $\cP = \{\p\}$ with $\q \ll \p$, the numeraire is just the likelihood ratio $E^* = \d\q/\d\p$. Indeed, $E^*$ is $\q$-almost surely strictly positive, it is an e-variable, and for any other e-variable $E \in \fE$ we have by a simple change of measure that $\E^\q[E/E^*] = \E^{\p}[E] \le 1$. 

Even for composite nulls, the numeraire is a likelihood ratio of $\q$ to the `reverse information projection'. Before getting there, we present a few other simple results.

\begin{lemma}[Lifting]
\label{lem:lifting}
Let $E^*$ denote the numeraire for $\cP$, and consider a larger null hypothesis $\cP' \supseteq \cP$. If $E^*$ is still an e-variable for $\cP'$, then it is also the numeraire for $\cP'$.
\end{lemma}
The proof is simple.
If $E$ is an e-variable for $\cP'$, it is also an e-variable for $\cP$ and one has $\E^\q[E/E^*] \leq 1$ by assumption, proving the lemma.

A numeraire $E^*$ is also {log-optimal} (Section~\ref{sec:3-p-vs-q}) in the sense that
\begin{equation} \label{eq_log_optimality}
\text{$\E^\q\left[ \log \frac{E}{E^*} \right] \le 0$ for every e-variable $E \in \fE$,}
\end{equation}
where the left-hand side may be $-\infty$. This follows directly from Jensen's inequality and the numeraire property. The converse is also true, and we record the equivalence in the following proposition. The proof shows that the numeraire property is really the first-order condition for log-optimality. Note also that a numeraire $E^*$ is the $\q$-almost surely unique log-optimal e-variable in the sense of \eqref{eq_log_optimality} even if $\E^\q[\log E^*]$ happens to be infinite.

\begin{proposition}[Log-optimality of numeraire] \label{P_log_optimality}
Let $E^*$ be a $\q$-almost surely strictly positive e-variable. Then $E^*$ is a numeraire if and only if it is log-optimal. In particular, a log-optimal e-variable is unique up to $\q$-nullsets.
\end{proposition}
\begin{proof}[Proof.]
The forward direction was argued above. To prove the converse we assume $E^*$ is log-optimal. For any e-variable $E$ and $t \in (0,1)$, $E(t) = (1-t) E^* + t E$ is an e-variable. Thus by log-optimality, $\E^\q[t^{-1}\log(E(t)/E^*)] \le 0$. The expression inside the expectation equals $t^{-1} \log(1 - t + t E/E^*)$,
which converges to $E/E^* - 1$ as $t$ tends to zero and is bounded below by $t^{-1}\log(1-t)$, hence by $-2 \log 2$ for $t \in (0, 1/2)$. Fatou's lemma thus yields $\E^\q[E/E^* - 1] \le 0$, showing that $E^*$ is a numeraire. Finally, the uniqueness statement follows from the equivalence just proved together with uniqueness of the numeraire up to $\q$-nullsets.
\end{proof}

The following result says that a numeraire always exists without any assumptions whatsoever on $\cP$ or $\q$. Although no assumptions on $\cP$ or $\q$ are needed in general, we do not present this more general theory for reasons of brevity. It turns out that the theory simplifies under the following condition, which, despite being nontrivial, is still very weak.

We say that \emph{$\q$ is absolutely continuous with respect to $\cP$}, written $ \q \ll \cP,$ if $\q$ does not assign positive probability to any event that is null under every $\p \in \cP$. Equivalently,
\begin{equation}\label{def:abs-cont-composite}
    \q \ll \cP \iff  \q(A)>0  \mbox{ implies } \sup_{\p\in \cP} \p(A)>0 \mbox{ for every $A$}.
\end{equation}
 This natural generalization of absolute continuity for pairs of measures is satisfied in most statistically relevant settings, while being significantly weaker than requiring a common dominating reference measure. It is also weaker than requiring that $\q \ll \p$ for every $\p \in \cP$, which would require that $\q(A)=0$ as long as any (as opposed to every) $\p \in \cP$ assigns $\p(A)=0$ (or contrapositively, that $\q(A) > 0$ implies $\p(A) > 0$ for all $\p \in \cP$).

\begin{theorem} \label{T_numeraire}
A numeraire always exists. Moreover, the following conditions are equivalent:
\begin{enumerate}
[label=(\alph*)]
\item\label{T_numeraire_1} The numeraire is $\q$-almost surely finite.
\item\label{T_numeraire_2} Every e-variable is $\q$-almost surely finite.
\item\label{T_numeraire_3} $\q$ is absolutely continuous with respect to $\cP$, that is $\q \ll \cP$.
\end{enumerate}
\end{theorem}
\begin{proof}[Proof]
    The proof that the numeraire exists is sophisticated and technical, thus omitted here. We focus on the equivalence of the other three claims. Clearly~\ref{T_numeraire_2} implies~\ref{T_numeraire_1}. To see that~\ref{T_numeraire_1} implies~\ref{T_numeraire_2}, let $E^*$ be a $\q$-almost surely finite numeraire and note that the numeraire property implies that any e-variable $E$ must be $\q$-almost surely finite too. Next, we prove that~\ref{T_numeraire_2} is equivalent to~\ref{T_numeraire_3}. First, suppose~\ref{T_numeraire_3} fails. Then there is an event $A$ with $\q(A) > 0$ and $\p(A) = 0$ for all $\p \in \cP$, in which case $E = \infty \id_{A}$ is an e-variable that is not $\q$-almost surely finite. Thus~\ref{T_numeraire_2} fails also. The converse direction follows by noting that every e-variable is finite $\p$-almost surely for all $\p \in \cP$. 
\end{proof}

We end with a few additional properties of the numeraire. Since the constant one is always an e-variable, the numeraire $E^*$ satisfies the aesthetically pleasing property that
\[
\E^\p[E^*] \leq 1 \text{ for all }\p \in \cP \text{ and } \E^\q\left[\frac1{E^*}\right] \leq 1.
\]
When $\cP = \{\p\}$ is simple and $\q \ll \p$, both inequalities above hold with equality, attained by the numeraire $\d\q/\d\p$. Another pleasing property is obtained by Jensen's inequality: The numeraire $E^*$ satisfies
\[
\E^\q[E^*/E] \geq 1 \text{ for all } E \in \fE.
\]

Finally, note that the definition of KL divergence in~\eqref{eq:kl-def1} can be rewritten and generalized to nonnegative measures $\p$ in the following manner. Let $\cM_+$ denote the set of all nonnegative measures on $\Omega$. Then, for any $\p \in \cM_+$ and $\q\in \cM_1$, define
\begin{equation}
\label{eq:kl-def2}
\kl(\q,\p) =  -\E^\q\left[\log \frac{\d \p^a}{\d\q} \right] ,
\end{equation}
where $\p^a$ is the absolutely continuous part of $\p$ with respect to $\q$. When $\p\in \cM_1$, this definition coincides with \eqref{eq:kl-def1}.

We first begin with the simple observation that since $\log y \leq y - 1$, we have by a change of measure that for any $E \in \fE$ and $\p \in \conv(\cP)$,
\[
\E^\q\left[ \log \left( E \frac{\d\p^a}{\d\q} \right) \right] \leq \E^\p\left[ E  \right] - 1 \leq 0.
\]
This immediately implies the following fact that we record for later use.


\begin{proposition}[Weak duality]
\label{prop:weak-duality}
For all e-variables $E \in \fE$ and distributions $\p \in \conv(\cP)$, we have $\E^\q[\log E] \leq \kl(\q,\p)$ (where if $\E^\q[\log E]$ is undefined, we treat it as $-\infty$), and thus
\begin{equation}\label{eq:weak-duality}
\sup_{E \in \fE} \E^\q[\log E] \leq \inf_{\p \in \conv(\cP)} \kl(\q,\p).
\end{equation}
\end{proposition}

Note that both sides of \eqref{eq:weak-duality} are nonnegative and well defined, since the constant $1$ is always an e-variable.
It turns out that the numeraire $E^*$ satisfies
\[
\E^\q[\log E^*] = \kl(\q,\p^*),
\]
for an appropriately defined $\p^*$ that we introduce and discuss next.

\section{Strong duality for finite \texorpdfstring{$\cP$}{\cP}}

We start here with the case that $\cP$ is finite. 
Then, it is possible to show the following result. Recall that $\Delta_K$ denotes the probability simplex. 

\begin{proposition}
\label{prop:strong-duality-finite}
Suppose $\cP = \{\p_1,\dots,\p_K\}$ with $\p_k \ll \q $ for all $k \in [K]$, and
$\min_{\p \in \conv(\cP)} \kl(\q,\p) < \infty$. Then,
\begin{equation}\label{eq:duality-simple}
 \E^\q[\log E^*] = \ \max_{E \in \fE} \E^\q[\log E] = 
 \min_{\p \in \conv(\cP)} \kl(\q, \p) 
 = \kl(\q,\p^*),
\end{equation}
where $E^*$ is the numeraire and $\p^*$ achieves the minimum.
\end{proposition}
To elaborate, we know that the numeraire achieves the left hand side ($\q$-almost surely uniquely), so the first equality holds by its definition.
The minimization problem can be rewritten as $$\min_{\lambda \in \Delta_K} \kl\left(\q, \sum_{k=1}^K \lambda_k \p_k\right) .
$$
This is a (strictly) convex minimization problem over the compact convex set $\Delta_K$, so it achieves its minimum uniquely at some $\lambda^*\in\Delta_K$, and we define the reverse information projection (RIPr) to be
\[
\p^* = \sum_{k=1}^K \lambda^*_k \p_k,
\]
so the last equality also holds by definition. Since $\kl(\q,\p^*) < \infty$ by assumption, we must have $\q \ll \p^*$. We now prove the strong duality result in the middle equality. 

\begin{proof}[Proof of Proposition~\ref{prop:strong-duality-finite}]
Recall from the earlier weak duality result of Proposition~\ref{prop:weak-duality} that if we find any e-variable whose e-power equals $\kl(\q,\p^*)$, then the proof is complete, and further, that e-variable must be the numeraire.
By the optimality of $\p^*$ and convexity of $\{\sum_{k=1}^K \lambda_k \p_k: \lambda \in \Delta_K\}$, we that have for any $\p \in \cP$,
\[
\kl(\q,\p^*) = \min_{\alpha \in [0,1]} \kl(\q, \alpha \p^* + (1-\alpha)\p)
\]
The optimality criterion (Karusch--Kuhn--Tucker condition) for the above minimization problem is that the gradient with respect to $\alpha$ is nonpositive at $\alpha=1$ (going from 1 towards 0 increases the objective):
\[
-\E^\q\left[ \frac{\d(\p^* - \p)}{\d\p^*} \right] \leq 0  ~ \Longleftrightarrow  ~ \E^\q\left[ \frac{\d\p}{\d\p^*} \right] \leq 1   ~ \Longleftrightarrow ~ \E^\p\left[ \frac{\d\q}{\d\p^*} \right] \leq 1.
\]
The last inequality means that $\d\q/\d\p^*$ is an e-variable. By definition of e-power, the e-power of this e-variable is $\kl(\q,\p^*)$. 
Thus, the numeraire $E^*$ is $\d\q/\d\p^*$ and the proof is complete.
\end{proof}




\section{Reverse information projection}\label{sec:ripr}

Remarkably, all assumptions in Proposition~\ref{prop:strong-duality-finite} can be dropped. But doing this is very challenging, as explained below.

First, it is clear that in the previous section, $\p^*$ is the closest distribution in $\conv(\cP)$ to $\q$ as measured by KL. However, when the infimum KL divergence equals infinity, it is unclear how do define the ``closest'' distribution to $\q$.

Second, in the previous section, we used the compactness of $\conv(\cP)$ when $\cP$ is finite to argue that the infimum is attained by a distribution in $\conv(\cP)$. However, when $\cP$ is uncountable (as would typically be the case in statistical applications), the set is no longer compact (with respect to the weak topology) and it is not clear that the infimum is achieved inside the set. Indeed, it is typically the case that the minimizer lies outside of $\conv(\cP)$, and in fact may be a \emph{sub}distribution in general (a nonnegative measure that integrates to less than one). In that case, where does it lie, and how do we characterize it? 

Finally, when all absolute continuity requirements are dropped, all arguments involving likelihood ratios (Radon--Nikodym derivatives) and KL divergences become much more subtle, and one cannot proceed as easily as in the previous section.  

Since formally handling the above subtleties requires some technical measure-theoretic care, we only state the results here but do not prove them. The first step is to define the reverse information projection (RIPr) appropriately. 





\begin{definition}[Reverse information projection (RIPr)]
    Let $E^*$ be the numeraire for $\cP$ against $\q$. The reverse information projection of $\q$ onto $\cP$, denoted $\p^*$, is defined via
    \[
    \frac{\d\p^*}{\d\q} = \frac{1}{E^*},
    \]
    where the right-hand side is understood to equal zero on $\{E^*=\infty\}$.
\end{definition}

In words, $\p^*$ is a nonnegative measure  whose likelihood ratio with respect to $\q$ equals $1/E^*$. 
Since $E^*$ exists, and is $\q$-almost surely unique and positive, we see that $\p^*$ is uniquely defined up to $\q$-nullsets. Also, $\p^* \ll \q$ by definition.
Moreover, $\p^*(\Omega)\le 1$ because $\E^\q[1/E^*]\le 1$, and hence $\p^*$ is in general a sub-probability  measure, though it is a probability measure in some special cases.

\subsection{The effective null hypothesis is the bipolar of $\cP$}

Let $\cX_+$  denote the set of all $[0,\infty]$-valued measurable functions on $\Omega$ and recall that $\cM_+$ denotes the set of all nonnegative measures on $\Omega$.  The \emph{polar} of $\cP$ is defined by  
\[
\cP^\circ = \{E \in \cX_+ \colon \E^\p[E] \le 1 \text{ for all } \p \in \cP\}.
\]
This is precisely the set of e-variables, i.e., $\cP^\circ = \fE$. The \emph{bipolar} of $\cP$ is
\[
\cP^{\circ\circ} = \{\p \in \cM_+ \colon \E^\p[E] \le 1 \text{ for all } E \in \cP^\circ\}.
\]
It is easy to see that 
\[
\cP^{\circ\circ} \cap \cM_1 \supseteq \conv(\cP).
\]
Equality holds when $\cP$ is finite,
but the inclusion can be strict.
Because the constant $E = 1$ is always an e-variable, elements of $\cP^{\circ\circ}$ have total mass at most one, and sometimes strictly less than one. For example, if $\p$ belongs to $\cP^{\circ \circ}$, then so does any $\p'\in
\cM_+$ that is set-wise dominated by $\p$, meaning $\p'(A) \le \p(A)$ for all $A$.
Below we give two examples where
$\cP^{\circ\circ} \cap \cM_1$
is a strict superset of $ \conv(\cP).$

\begin{example}
\label{ex:c6-bipol}
    Suppose the null $\cP$ consists of all unit variance Gaussians with  mean less than $0$. Then the standard Gaussian is not in the convex hull of $\cP$, but it is in $\cP^{\circ\circ} \cap \cM_1$, which contains all unit variance Gaussians with nonnegative mean (and their mixtures).
\end{example}

The reader may note that if we had considered the total variation closure of the convex hull of $\cP$ (motivated by Theorem~\ref{th:kraft}) in Example~\ref{ex:c6-bipol}, it would have coincided with $\cP^{\circ\circ} \cap \cM_1$. However, more generally, this is not the case.

\begin{example}
Let $\Omega = [0,1]$, and let 
\[
\cP = \{\delta_x : x \in [0,0.5]\},
\]
where $\delta_x$ is a point mass at $x$. Clearly, $\cP$ has no common reference measure. The total variation closure of the convex hull of $\cP$ still contains only atomic measures. However, $\cP^{\circ\circ} \cap \cM_1$ contains every bounded distribution on $[0,0.5]$.  
\end{example}

In the above example, if,
instead of total variation closure, we had taken a weak closure of the convex hull of $\cP$, then again we would have obtained $\cP^{\circ\circ} \cap \cM_1$. However, one can construct yet other examples to show that this is not true in general, but this gets rather deep into topological measure theory, and is out of the scope of the current book.

The bipolar of $\cP$ can be interpreted as the \emph{effective null hypothesis}: it is the largest family of distributions whose set of e-variables is exactly $\fE$. Thus it consists of those distributions against which e-variables in $\fE$ do not provide any evidence. In particular, we note the following simple yet important consequence, recalling the definition of a powered e-variable from Definition~\ref{def:nontrivial}. 

\begin{theorem}
\label{thm:nontrivial-iff-bipolar}
A powered test (or e-variable) for $\cP$ against $\q$ exists if and only if $\q \notin \cP^{\circ\circ}$.
\end{theorem}
\begin{proof}[Proof.]
The proof follows in two simple steps. It follows from Proposition~\ref{th:existence} that a powered test exists if and only if a powered e-variable exists. The latter condition is possible if and only if $\q \notin \cP^{\circ\circ}$ by definition of $\cP^{\circ\circ}$.
\end{proof}

It turns out that considering only bounded e-variables suffices. Let $\fE^b \subsetneq \cP^\circ = \fE$ be defined as $\fE^b = \{E \in \cX_+ \colon \E^\p[E] \le 1 \text{ for all } \p \in \cP, ~\sup_{\omega \in \Omega} E(\omega) < \infty\}$. Then an application of the monotone convergence theorem yields that
\[
\cP^{\circ\circ} = \{\p \in \cM_+ \colon \E^\p[E] \le 1 \text{ for all } E \in \fE^b\}.
\]

We now claim that even though it is possible that $\p^* \notin \conv(\cP)$ in general, it does lie in $\cP^{\circ\circ}$.

\begin{figure}
    \centering
\includegraphics[width=\textwidth]{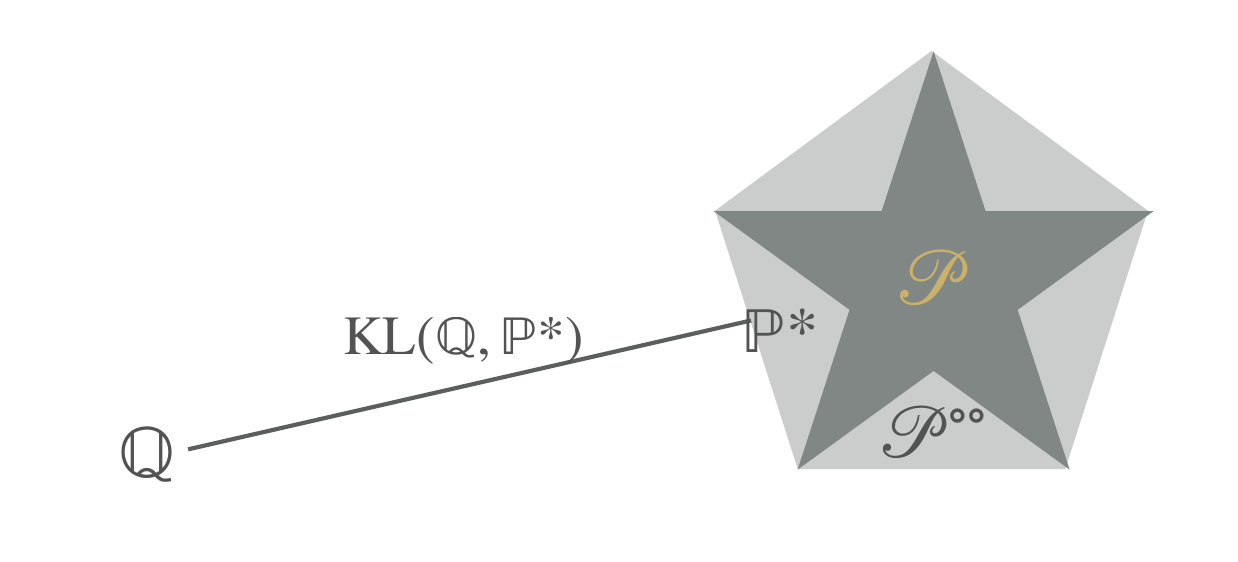}
\caption{A visualization of the RIPr $\p^*$ of $\q$ onto $\cP$, which lies in $\cP^{\circ\circ}$, but not necessarily in $\cM_1$.}
\label{fig:RIPr}
\end{figure}

\begin{proposition}\label{prop:ripr-in-bipolar}
    The reverse information projection of $\q$ onto $\cP$ lies in the bipolar of the null: $\p^* \in \cP^{\circ\circ}$.
\end{proposition}
    The proof follows because for any e-variable $E$, 
    \[
    \E^{\p^*}[E] = \E^{\q}[\id_{\{E^* < \infty\}}E/E^*] \le \E^\q[E/E^*] \le 1,
    \]
    because $E^*$ is the numeraire. The above result can be visualized as Figure~\ref{fig:RIPr}.

\subsection{A general strong duality result}

We now state a more general strong duality result removing almost all measure-theoretic assumptions. As before, $\E^\q[\log E]$ is understood as $-\infty$ whenever $\E^\q[(\log E)_-] = \infty$.  Recall the condition $\q \ll \cP$ from~\eqref{def:abs-cont-composite}.

\begin{theorem}
\label{th:strong-duality-ripr}
Assume $\q \ll \cP$.
    Letting $E^*$ be the numeraire and $\p^*$ the \textnormal{RIPr}, one has the strong duality relation
\begin{equation} \label{eq_strong_duality}
\E^\p[\log E^*] = \sup_{E \in \fE} \E^\q[\log E] = \inf_{\p \in \cP^{\circ\circ}} \kl(\q,\p) = \kl(\q,\p^*),
\end{equation}
where these quantities may equal $\infty$.  
\end{theorem}

We also provide some characterizing properties of the RIPr.

\begin{theorem}
\label{thm:characterize-RIPr}
Assume $\q \ll \cP$ and let $\p^*$ be an element of $\cP^{\circ\circ}$ equivalent to $\q$.
   Let  $\p^a$ denote the absolutely continuous part of $\p$ with respect to $\q$. The following statements imply each other:
   \begin{enumerate}
       \item[(i)] $\p^*$ is the RIPr.
       \item[(ii)] $\displaystyle \E^\q\left[ \frac{\d\p^a}{\d\p^*} \right] \le 1 \text{ for all } \p \in \cP^{\circ\circ}$
       \item[(iii)] $\displaystyle \E^\q\left[ \log \frac{\d\p^a}{\d\p^*} \right] \le 0 \text{ for all } \p \in \cP^{\circ\circ}$.  
   \end{enumerate}
   If we assume further that all elements of $\cP$ are absolutely continuous with respect to a probability measure $\leb$, then all elements of $\cP^{\circ\circ}$ are also absolutely continuous with respect to $\leb$, so that we can identify any $\p \in \cP^{\circ\circ}$ with its density $p = \d \p/\d \leb$. Further, each of the above conditions is equivalent to
   \begin{enumerate}
       \item[(iv)] 
$\E^\q\left[p/p^*\right] \le 1 \text{ for all densities $p$ of $\p \in \cP$.}$
   \end{enumerate}
\end{theorem}

We do not prove the above theorems here. Crucially, note that the fourth condition above only needs to hold over $\cP$ and not over the bipolar, making it easier to verify in settings with a reference measure such as exponential families, as we will see in Section~\ref{sec:exp-fam-numeraire}.

\section{Numeraire is the only e-variable that is a likelihood ratio}\label{sec:numeraire-LR}

We first present a verification theorem that simplifies the task of checking that candidates $E^*$ and $\p^*$ (for the numeraire and RIPr respectively) are in fact optimal.

\begin{theorem} \label{T_verification}
Assume that $\q \ll \cP$. Let $E^*$ be a $\q$-almost surely strictly positive e-variable and let $\p^*$ be the measure given by $\d\p^*/\d\q = 1/E^*$. Then $E^*$ is a numeraire if and only if $\p^*$ belongs to the effective null $\cP^{\circ\circ}$. Thus the numeraire is the only e-variable that is also a likelihood ratio between $\q$ and some equivalent element of $\cP^{\circ\circ}$. Finally, $\E^{\p^*}[E^*] = 1$ and $\p^*$ is maximal in the sense that no other $\p \in \cP^{\circ\circ}$ equivalent to $\q$ can satisfy $\p(A) \ge \p^*(A)$ for all $A$ with strict inequality for some $A$.
\end{theorem}
\begin{proof}[Proof]
Recall that all e-variables are $\q$-almost surely finite thanks to Theorem~\ref{T_numeraire} and the assumption that $\q \ll \cP$; in particular, $\p^*$ is equivalent to $\q$. Now, if $E^*$ is a numeraire, then Proposition~\ref{prop:ripr-in-bipolar} yields that  $\p^* \in \cP^{\circ\circ}$. For the converse, the definition of $\p^*$ and the fact that it belongs to $\cP^{\circ\circ}$ yield
$\E^\q[E/E^*] = \E^{\p^*}[E] \le 1$ for any e-variable $E$. This is the numeraire property. Finally, the definition of $\p^*$ immediately gives $\E^{\p^*}[E^*] = 1$, and $\p^*$ must be maximal because if a `dominating' equivalent $\p \in \cP^{\circ\circ}$ existed we would get the contradiction $1 = \E^{\p^*}[E^*] < \E^\p[E^*] \le 1$.
\end{proof}

We stress that the above result is quite remarkable: the numeraire is clearly a likelihood ratio between $\q$ and some equivalent element of $\cP^{\circ\circ}$, but it is also the unique e-variable with such a representation! This property helps identify the numeraire in several simple examples (eg: Sections~\ref{sec:exp-fam-numeraire} and~\ref{sec:mlr-numeraire}), where we guess the numeraire as being a particular simple likelihood ratio of the alternative to some ``extreme'' point of the null, verify that it is indeed an e-variable, and conclude that it is in fact the log-optimal e-variable.

\begin{corollary} \label{C:240207}
Assume $\q \ll \cP$. Let $\fE_0$ be a family of $[0,\infty]$-valued random variables that generates the null hypothesis in the sense that
\begin{equation} \label{T_verification_1}
\cP = \{\p \in \cM_1 \colon \E^\p[E] \le 1 \text{ for all } E \in \fE_0\}.
\end{equation}
Let $E^*$ be a $\q$-almost surely strictly positive e-variable such that
\begin{equation} \label{T_verification_2}
\E^\q\left[ \frac{1}{E^*} \right] = 1 \text{ and } \E^\q\left[ \frac{E}{E^*} \right] \le 1 \text{ for all } E \in \fE_0.
\end{equation}
Then $E^*$ is the numeraire and the \textnormal{RIPr} belongs to $\cP$.
\end{corollary}
\begin{proof}[Proof]
Let $\p^*$ be given by $\d\p^*/\d\q = 1/E^*$. Then
\eqref{T_verification_2} says that $\p^*$ is a probability measure such that $\E^{\p^*}[E] \le 1$ for all $E \in \fE_0$. Thus by \eqref{T_verification_1}, $\p^*$ belongs to $\cP$, hence to $\cP^{\circ\circ}$. The result now follows from Theorem~\ref{T_verification}.
\end{proof}

Thanks to Lemma~\ref{lem:lifting}, Corollary~\ref{C:240207} also holds if one has `$\supseteq$' instead of `$=$' in \eqref{T_verification_1}. This is intuitive: if one enlarges the generating set $\fE_0$ beyond what is necessary to specify $\cP$, then the right-hand side of~\eqref{T_verification_1} can become smaller than $\cP$. In that case, \eqref{T_verification_2} only gets harder to satisfy, and if we find an e-variable $E^*$ satisfying it, it surely must still be the numeraire.

Another useful corollary is the following \emph{shrinking lemma}, a dual result to the \emph{lifting lemma} (Lemma~\ref{lem:lifting}).
\begin{lemma}[Shrinking]
\label{lem:shrinking}
    Assume that $\q \ll \cP$, and let $\p^*$ denote the RIPr of $\q$ onto $\cP$ and let $X^* = \d\q/\d\p^*$ denote the numeraire for $\cP$ against $\q$. For any $\overline \cP \subseteq \cP$, if $\p^* \in \overline\cP^{\circ\circ}$, then $\p^*$ is the RIPr of $\q$ onto $\overline\cP$, so $X^*$ is still the numeraire for $\overline \cP$ against $\q$.
\end{lemma}
\begin{proof}[Proof.]
   Clearly, $\d\q/\d\p^*$ is an e-variable for $\cP$ and thus also for $\overline \cP \subseteq \cP$. Since $\p^* \in \overline\cP^{\circ\circ}$, and noting that the numeraire is the only e-variable that is a likelihood ratio between $\q$ and some equivalent element of $\overline\cP^{\circ\circ}$, it must be the case that $X^*$ is the numeraire for $\overline \cP$ against $\q$ and hence $\p^*$ is the RIPr of $\q$ onto $\overline \cP$. 
\end{proof}

As an example of how one might apply the shrinking and lifting lemmata referenced above, consider $\cP$ being unit-variance Gaussians with nonpositive mean, and $\q$ being $N(1,1)$. It can be easily checked that the RIPr is the standard Gaussian, and the numeraire is the likelihood ratio of $N(1,1)$ to $N(0,1)$. The shrinking lemma implies that if we consider any subset of $\cP$ that still contains $N(0,1)$, the RIPr and numeraire remain identical. We could have also gone in the other direction (from the singleton null of $N(0,1)$ to the one-sided null $\cP$),
and used the lifting lemma instead to reach the same conclusions.

\section{Numeraire dominates universal inference}
\label{sec:numeraire-beats-ui}

Assume now that there exists a common reference measure $\leb$,  and identify distributions with their respective densities, written in lowercase. For a point alternative $q$ over the data $Z$, the method of universal inference boils down to constructing the e-variable $E^{\mathrm{UI}} = q(Z)/p_\text{max}(Z)$ where $p_\text{max}(Z) = \sup_{p \in \cP} p(Z)$ is the maximum likelihood. (Strictly speaking, the supremum here should be understood as an essential supremum under the reference measure.)

To compare the universal inference e-variable with the numeraire $E^* = q(Z)/p^*(Z)$, we first claim that the RIPr satisfies $p^*(Z) \le p_\text{max}(Z)$ up to $\leb$-nullsets. This is obvious if the RIPr belongs to $\cP$, but in general it only belongs to $\cP^{\circ\circ}$. In this case the claim follows by applying the following lemma with $f = p_\text{max}$.

\begin{lemma}
\label{lem:bipolar-domination}
Let $f$ be a (potentially random) function such that $f \ge p$, $\leb$-almost surely, for all $p \in \cP$. Then $f \ge p$, $\leb$-almost surely, for all $p \in \cP^{\circ\circ}$.
\end{lemma}
We omit the technical proof.
As an immediate consequence of the lemma, we  see that $p^*(Z) \le p_\text{max}(Z)$, and hence $E^* \ge E^{\mathrm{UI}}$, up to $\leb$-nullsets. In nondegenerate situations involving a composite null hypothesis, the inequality will be strict with positive $\leb$-probability. Thus, in such cases, the relatively general method of universal inference is in fact inadmissible. We end by noting that the numeraire imposes weaker assumptions than universal inference (and thus is, in some sense, even more universal), since the latter needs a reference measure to define likelihoods. 
 
\begin{proof}[Composite alternatives]
Recall the composite alternative setup and notation described in Section~\ref{sec:3-mix-plugin}.
A reasonable goal is to derive an e-variable $E_n = E(X^n)$ for $\cP$ such that for any $\q \in \cQ$,
\begin{equation}\label{eq:asymp-log-optimality-RIPr}
 \lim_{n \to \infty} \frac{\E^\q[\log E_n]}{n} \to \kl(\q,\p^*),
\end{equation}
where $\p^*$ is the RIPr of $\q$ onto $\cP$. As we have seen in this chapter, the right-hand side is the growth rate of an oracle that knows the data-generating distribution (when it lies in the alternative).
One can derive appropriate extensions of the mixture and plug-in e-variables described in Section~\ref{sec:3-mix-plugin}. For example, in the iid setting considered there, consider the numeraire for testing $\cP^{(n)} = \{\p^n : \p \in \cP\}$ against the alternative $\widetilde \q^{(n)} = \int \q^n \nu(\d \q)$. We conjecture that this is indeed asymptotically log-optimal under weak assumptions, but such a general result is currently not known. 
\end{proof}

We end with a note on logical coherence. Saying that a composite null $\cP$ is true is equivalent to saying that a simple null $\p$ is true for some $\p \in \cP$. Thus, one may expect to only reject $\cP$ if one can also reject every $\p \in \cP$ separately. A method that satisfies this property is called logically coherent. We can extend this definition to e-variables as well. Let $E^\cP$ denote an e-variable for testing $\cP$. The collection $(E^\cP)_{\cP \subseteq   \cM_1}$ can be called coherent if $E^\cP \leq E^{\cP'}$ whenever $\cP' \subseteq \cP$. Fixing an alternative $\q$ for simplicity, we now make two simple observations, stated here without proof:
\begin{itemize}
    \item If $E^\cP$ equals the universal inference e-variable for $\cP$ against $\q$, the collection $(E^\cP)_{\cP \subseteq  \cM_1}$ is always logically coherent. (This point is easy to verify by definition.)
    \item If $E^\cP$ equals the numeraire e-variable for $\cP$ against $\q$, the collection $(E^\cP)_{\cP\subseteq   \cM_1}$ is generally not logically coherent, except for some special cases.
\end{itemize}
Incoherent families of e-variables may be harder to interpret, because enlarging the null may counterintuitively increase the evidence against the null.
Thus, there appears to be a tradeoff between e-power and coherence/interpretability when comparing universal inference and the numeraire.

\section{Connections between     e-values  and Bayes factors}
\label{sec:Bayes-factor}

For those Bayesians who are willing to consider testing hypotheses, one of the standard quantities to examine is the \emph{Bayes factor}. Suppose we observe $n$ iid data points $X_1,\dots,X_n$, and we are interested in testing whether their marginal distribution lies in a composite null $\cP$ or in a composite alternative $\cQ$. Also suppose that there is a common reference measure, allowing us to identify distributions with their densities. Then the Bayes factor is defined as
\[
B_n = \frac{\int_{\cQ} \prod_{i=1}^n q(X_i)\, \mu(\d q)}{\int_{\cP} \prod_{i=1}^n p(X_i)\,  \nu(\d p)},
\]
where $\mu$ and $\nu$ are (often subjective) \emph{priors} over $\cQ$ and $\cP$. Bayesians often refer to tables like in Section~\ref{sec:c2-jeffrey} in order to interpret the realized value, and typical Bayesians do not threshold $B_n$ to make clear decisions of accepting or rejecting the null, and consequently no effort is made to control the type-I error. The law of iterated expectation yields that $B_n$ is an e-variable \emph{if nature draws a random null from the scientist's prior}---$\E^{p \sim \nu(\d p)}[B_n] \leq 1$---but it remains unclear to many why this is a reasonable assumption. Indeed, this is often the problem with Bayesian methods: the scientist using them is perfectly convinced of its correctness (internal consistency) but their analysis fails to convince all others (lack of external consistency). Thus the main drawbacks of this approach are: (a) it is unclear how to generalize these ideas to nonparametric settings without densities and likelihoods, (b) it is unclear how to formally make decisions with $B_n$ with (frequentist) error guarantees.

In contrast, the typical frequentist approach in the above setting would be to consider the generalized likelihood ratio:
\[
G_n = \frac{\sup_{q \in \cQ} \prod_{i=1}^n q(X_i)\, }{\sup_{p \in \cP} \prod_{i=1}^n p(X_i)}.
\]
Under some regularity conditions,  Wilks' theorem implies that $\log G_n$ has a chi-squared limiting distribution as $n \to \infty$, with degrees of freedom determined by some notion of 
parametric dimensionality difference between $\cQ$ and $\cP$ (for the discussion here we allow $\cP$ and $\cQ$ to be nested, which is different from the main setting of the book). This is then used to construct an asymptotic p-value. However, the main drawback here is once more (a) it is unclear how to generalize these ideas to nonparametric settings without densities and likelihoods, and (b) what to do when the regularity conditions for Wilks' theorem fails, and the asymptotics are unwieldy.

In contrast, we note that universal inference takes a middle ground between the above two approaches, using a mixture (i.e.\ a prior without Bayesian connotations) over the alternative, but a maximum likelihood over the null:
\[
U_n = \frac{\int_{\cQ} \prod_{i=1}^n q(X_i)\, \mu(\d q)}{\sup_{p \in \cP} \prod_{i=1}^n p(X_i)},
\]
which is always an e-variable, under no regularity conditions. However, there are still two drawbacks (a) the previous section shows that this is not optimal, and (b) we still require reference measures and densities to define likelihoods, so nonparametric extensions are not immediately obvious.  

The numeraire resolves these issues. This sections shows that one should calculate the RIPr $\p^*$ of $\q^{\text{eff}} := \int_\cQ \q^n \mu(\d \q)$ onto (the bipolar of) $\{\p^n: \p \in \cP\}$. Since $\p^* \ll \q^{\text{eff}}$ by definition, we can define
\[
E_n  = \frac{\d \q^{\text{eff}}}{\d \p^*}(X_1,\dots,X_n).
\]
This is always well defined, even in nonparametric settings with no reference measures for $\cP,\cQ$, it is always an e-variable (thus guaranteeing the frequentist error control), and it dominates $U_n$ (Section~\ref{sec:numeraire-beats-ui}).

In some situations, the RIPr will actually lie in the convex hull (set of mixtures) of $\{\p^n: \p \in \cP\}$, and in this case, the numeraire can be written in the form of $B_n$. However, for every choice of alternative mixture, there is a unique mixture over the nulls that would yield such a representation (Section~\ref{sec:numeraire-LR}), meaning that one cannot decouple the pair of mixtures as the Bayesian did for $B_n$. In these cases, the numeraire can be seen as the unique Bayes factor which is also an e-variable.

Of course, when the null and alternative are simple, all the above methods coincide. When only the null is simple, the Bayes factor, the universal inference e-variable, and the numeraire coincide.
When neither is simple, the Bayes factor is sometimes (but usually not) an e-variable, the numeraire is sometimes (but usually not) a Bayes factor, and the generalized likelihood ratio is neither an e-variable nor a Bayes factor in general.

\section{Examples}
\label{sec:ripr-examples}

We now turn to examples of numeraire e-variables, revisiting several examples in Section~\ref{sec:c2-ex}. 
Consider a single random observation (statistic) $Z$ in a measurable space $\cZ$, and we take $(\Omega, \cF)$ to be $\cZ$ with its $\sigma$-algebra $\sigma(Z)$. We will always have $\q \ll \cP$. Recall that $\cM_1$ is the set of all probability measures on $\cF$, and denote by $\cX_+$  the set of all $[0,\infty]$-valued measurable functions and  by $\cM_+$  that of nonnegative measures on $\cF$.

\subsection{Exponential family with one-dimensional sufficient statistic}
\label{sec:exp-fam-numeraire}

The example here will be subsumed by the (much shorter) Section~\ref{sec:mlr-numeraire}. Nevertheless, we provide a self-contained treatment below that does not explicitly exploit the monotone likelihood ratio property, which may be less familiar to some.
Here $\cZ$ can be general and is equipped with a reference measure $\leb$. 
Consider an exponential family of densities
\[
p_\theta(z) = e^{\theta T(z) -  A(\theta)} 
\]
with respect to $\leb$, where $A$ is convex and differentiable, the sufficient statistic $T(z)$ is one-dimensional, and the natural parameter $\theta$ ranges in some interval $\Theta \subseteq \R$. The null hypothesis is $\cP = \{p_\theta \colon \theta \in \Theta_0\}$ for some closed subset $\Theta_0 \subseteq \Theta$, and the alternative is $q = p_{\theta_1}$ for some $\theta_1 \in \Theta$. We suppose that $\Theta_0$ has a largest element $\theta_0$ and that $\theta_1 > \theta_0$, as in Section~\ref{sec:mlr-evariable}. It is natural to conjecture that $p_{\theta_0}$ is the RIPr. To confirm this, it suffices to verify the fourth condition in Theorem~\ref{thm:characterize-RIPr}. Using the standard formula for the moment generating function of the sufficient statistic we get
\begin{align}
\E^\q\left[\frac{p_\theta(Z)}{p_{\theta_0}(Z)}\right] &= \E^\q\left[e^{(\theta-\theta_0) T(Z) - A(\theta) + A(\theta_0)}\right] \\
&= e^{A(\theta_1+\theta-\theta_0) - A(\theta_1) - A(\theta) + A(\theta_0)}
\end{align}
for all $\theta \in \Theta_0$.
Any real convex function $f$
satisfies $f(a+b+c)-f(a+c)\ge f(b+c)-f(c)$ for $a,b\ge 0$ and $c\in \R$. Therefore, 
$\E^\q[p_\theta(Z) / p_{\theta_0}(Z)] \le 1$ for all $\theta \in \Theta_0$, so that the fourth condition in Theorem~\ref{thm:characterize-RIPr} holds and $p_{\theta_0}$ is indeed the RIPr. As a result the numeraire is the likelihood ratio
\[
E^* = \frac{p_{\theta_1}(Z)}{p_{\theta_0}(Z)} = e^{(\theta_1-\theta_0) T(Z) - A(\theta_1) + A(\theta_0)}.
\]

\subsection{Monotone likelihood ratio}
\label{sec:mlr-numeraire}

Returning to the example described in Section~\ref{sec:mlr-evariable}, we now point out that for a point alternative $p_{\theta_1}$, the e-variable described there is in fact the numeraire. This generalizes the previous example which is well-known to satisfy the monotone likelihood ratio property (in the sufficient statistic). The proof follows by the same argument as above: since $p_{\theta_1}/p_{\theta_0}$ is an e-variable for $\cP$, then the fourth condition in Theorem~\ref{thm:characterize-RIPr} implies that $p_{\theta_0}$ is the RIPr, and thus $p_{\theta_1}/p_{\theta_0}$  is the numeraire.

\subsection{Testing symmetry}
\label{sec:symmetry}

Take $\cZ = \R$. As in Section~\ref{sec:c2-symmetry}, we consider the null hypothesis that $Z$ is symmetric,
\[
\cP = \left\{ \p \in \cM_1 \colon Z \text{ and } -Z \text{ have the same distribution under } \p \right\},
\]
which is a non-dominated family. We also fix an alternative hypothesis $\q$ that admits a Lebesgue density $q$. It is natural to conjecture that the RIPr is given by the symmetrization $\widetilde \p$ of $\q$, whose density is $\widetilde p(z) = (q(z) + q(-z))/2$. However, this cannot quite be true in general because $\widetilde\p$ need not be equivalent to $\q$. Instead, we claim that the RIPr is the measure $\p^*$ with density
\[
p^*(z) = \frac{1}{2}\left( q(z) + q(-z) \right) \id_{\{q(z) > 0\}}.
\]
This is the absolutely continuous part of $\widetilde\p$ with respect to $\q$. It is a probability measure if $\q$ has symmetric support, and otherwise a proper sub-probability measure. Note that $\widetilde\p$ belongs to $\cP$, which in particular shows that $\q \ll \cP$ since $\q \ll \widetilde\p$. (Alternatively, we again have that $\p(A) = 0$ for all $\p \in \cP$ implies that $A$ is empty, which also yields $\q \ll \cP$.) 

To check that $\p^*$ is the RIPr we will show that the implied candidate numeraire is, in fact, the numeraire. It is given by      
\[
E^* = \frac{\d\q}{\d\p^*} = \frac{2q(Z)}{q(Z) + q(-Z)} \id_{\{q(Z) > 0\}}.
\]
(We recall that the numeraire is only defined up to $\q$-nullsets.) 
We claim that the set of all e-variables is 
\[
\fE = \left\{E \in \cX_+ : E \le 1 + \phi(Z) \text{ for some odd function } \phi \right\},
\]
where we recall that a function $\phi$ is odd if $\phi(z) + \phi(-z) = 0$ for all $z \in \R$.  
Indeed, any $E$ of this form satisfies $\E^\p[E] \le 1 + \E^\p[\phi(Z)] = 1$, where the symmetry of $Z$ and the oddness of $\phi$ were used to get $\E^\p[\phi(Z)] = \E^\p[\phi(-Z)] = -\E^\p[\phi(Z)]$ and hence $\E^\p[\phi(Z)] = 0$. Conversely, for any e-variable $E = f(Z)$ we may write $E = \frac12 (f(Z) + f(-Z)) + \phi(Z)$ where $\phi(z) = \frac12(f(z) - f(-z))$ is the odd part of $f(z)$. For any $z \in \R$ we use the symmetric distribution $\p = \frac12 (\delta_z + \delta_{-z})$ and the fact that $E$ is an e-variable to get $\frac12 (f(z) + f(-z)) = \E^\p[f(Z)] \le 1$. Hence $E \le 1 + \phi(Z)$ as required.

We can now verify the numeraire property. First, $E^*$ is $\q$-almost surely strictly positive and finite, and it is an e-variable because it can be written as $E^* = 1 + \phi^*(Z)$ where
\[
\phi^*(z) = \frac{q(z) - q(-z)}{q(z) + q(-z)}
\]
is odd (we also verified that $E^*$ is an e-variable in Section~\ref{sec:c2-symmetry}).   Next, for any e-variable $E \le 1 + \phi(Z)$ where $\phi$ is odd we get
\begin{align*}
\E^\q\left[ \frac{E}{E^*} \right] &\le \E^\q\left[ (1 + \phi(Z)) \frac{q(Z) + q(-Z)}{2 q(Z)} \right]\\
&= \E^{\widetilde\p}\left[(1 + \phi(Z)) \id_{\{q(Z) > 0\}} \right]
\le   \E^{\widetilde\p}[1 + \phi(Z)] = 1.
\end{align*}
This confirms the numeraire property.

Finally, note that it is not really necessary that $\q$ admit a density. The symmetrization of $\q$ is still well-defined as $\widetilde\p = \frac12 (\q + \widetilde\q)$ where $\widetilde \q$ is the distribution of $-Z$ under $\q$, or equivalently, the push-forward of $\q$ under the reflection map $z \mapsto -z$. The RIPr is then the absolutely continuous part $\p^* = \widetilde \p^a$, and the numeraire is any nonnegative version of the Radon--Nikodym derivative $\d\q/\d\p^*$ such that $\d\q/\d\p^* - 1$ is odd.

\subsection{Testing bounded means}
\label{sec:bounded-mean}

Let $\cZ = [0,1]$ so that the random variable $Z$ takes values in the unit interval.  
Fix $\mu \in (0,1/2)$ and consider (as in Section~\ref{sec:c2-bounded}) the null hypothesis that the mean of $Z$ is at most $\mu$,
\[
\cP = \{\p \in \cM_1 \colon \E^\p[Z] \le \mu\}.
\]
There is no single dominating measure for $\cP$ since it contains the uncountable non-dominated family $\{\delta_z \colon z \in [0,\mu]\}$. However, $\cP$ is generated by $\fE_0 = \{Z/\mu\}$ in the sense of Corollary~\ref{C:240207}, so our strategy will be to locate a candidate numeraire and then apply the corollary to verify that the candidate is, in fact, the numeraire. The alternative hypothesis $\q$ is the uniform distribution on $[0,1]$, and we have $\q \ll \cP$ for the simple reason that $\p(A) = 0$ for all $\p \in \cP$ implies that $A$ must actually be empty.

To find a candidate numeraire, we observe that two natural e-variables are $Z/\mu$ and the constant one. All convex combinations of these are also e-variables; equivalently, $1 + \lambda(Z-\mu)$ is an e-variable for each $\lambda \in [0,1/\mu]$. We now look for a log-optimal e-variable in this class by directly maximizing $f(\lambda) = \E^\q[\log (1 + \lambda(Z-\mu))]$ over $\lambda \in [0,\mu^{-1}]$. This is a strictly concave function whose derivative $f'(\lambda) = \E^\q[(Z-\mu) / (1 + \lambda(Z-\mu))]$ satisfies $f'(0) = \E^\q[Z] - \mu > 0$ and 
$f'(\mu^{-1}) = \E^\q[\mu - \mu^2/Z] = -\infty$. Thus there is a unique interior maximizer $\lambda^* \in (0,\mu^{-1})$, which is characterized by the first-order condition
\begin{equation} \label{eq_bdd_mean_testing_foc}
\E^\q\left[ \frac{Z - \mu}{1 + \lambda^* (Z - \mu)} \right] = 0.
\end{equation}
Since $\q$ is the standard uniform distribution we can be more explicit. Nothing changes if we first multiply both sides by $\lambda^*$, and then the left-hand side becomes
\begin{align*}
\int_0^1 \frac{\lambda^*(z - \mu)}{1 + \lambda^* (z - \mu)} \d z
&= 1 -  \int_0^1 \frac{1}{1 + \lambda^* (z - \mu)} \d z \\
&= 1 - \frac{1}{\lambda^*}\log\left(\frac{1 + \lambda^* (1 - \mu)}{1 - \lambda^*\mu} \right).
\end{align*}
Thus \eqref{eq_bdd_mean_testing_foc} for $\lambda^* \in (0,\mu^{-1})$ is equivalent to
\[
    \frac{1 + \lambda^* (1 - \mu)}{1 - \lambda^*\mu} = e^{\lambda^*},
\]
which is easily solved numerically.
This leads us to the candidate numeraire
\[
E^* = 1 + \lambda^* (Z - \mu),
\]
which is strictly positive and finite. To  verify that this is indeed the numeraire, we use \eqref{eq_bdd_mean_testing_foc}
to get, for any e-variable of the form $X = 1 + \lambda (Z - \mu) = E^* + (\lambda - \lambda^*)(Z-\mu)$, that 
\begin{equation} \label{eq_ex_bdd_mean_num_check}
\E^\q\left[ \frac{X}{E^*} \right] = 1 + (\lambda - \lambda^*) \E^\q\left[ \frac{Z - \mu}{1 + \lambda^* (Z - \mu)} \right] = 1.
\end{equation}
Taking $\lambda = 0$ and $\lambda = 1/\mu$ we see that \eqref{T_verification_2} of Corollary~\ref{C:240207} is satisfied and, hence, that $E^*$ is the numeraire and the RIPr $\p^*$  belongs to $\cP$.

\subsection{Testing sub-Gaussian means}
\label{sec:subG-mean}

Now take $\cZ = \R$ and consider the null hypothesis that the observation $Z$ has a 1-sub-Gaussian distribution with nonpositive mean (as in Section~\ref{sec:c2-subG-mean}). That is, we set
\[
\cP = \left\{\p \in \cM_1 \colon \E^\p[ e^{\lambda Z - \lambda^2/2} ] \le 1 \text{ for all } \lambda \in [0,\infty) \right\}.
\]
The ``one-sided'' restriction $\lambda \in [0,\infty)$ implies that $Z$ has nonpositive (potentially nonzero and even infinite) mean under any $\p \in \cP$. Indeed, monotone convergence and the definition of $\cP$ yield $\E^\p[Z] = \lim_{\lambda \downarrow 0} \E^\p[ (e^{\lambda Z} - 1)/\lambda ] \le \lim_{\lambda \downarrow 0} (e^{\lambda^2/2} - 1)/\lambda = 0$. As in the previous example, $\cP$ does not admit any dominating measure. It is generated by the family $\fE_0 = \{e^{\lambda Z - \lambda^2/2} \colon \lambda \in [0,\infty) \}$, so we will again look for a candidate numeraire and verify it using Corollary~\ref{C:240207}. We let the alternative hypothesis $\q$ be normal with mean $\mu > 0$ and unit variance. As in the previous example, and for the same reason, we have $\q \ll \cP$.

To find a candidate numeraire we maximize $\E^\q[\log X]$ over $X \in \fE_0$. That is, we maximize $\E^\q[ \lambda Z - \lambda^2/2] = \lambda \mu - \lambda^2/2$ over $\lambda \in [0,\infty)$. The maximizer is $\lambda^* = \mu$, which yields the candidate
\[
E^* = e^{\mu Z - \mu^2/2}.
\]
This is finite and strictly positive. Moreover, \eqref{T_verification_2} of Corollary~\ref{C:240207} is satisfied because
\begin{equation} \label{eq_sub_gaussian_testing_foc}
\E^\q\left[ \frac{e^{\lambda Z - \lambda^2/2}}{e^{\mu Z - \mu^2/2}} \right] = \E^\q\left[ e^{(\lambda - \mu) (Z - \mu) - (\lambda - \mu)^2/2} \right] = 1.
\end{equation}
We conclude that $E^*$ is the numeraire and, consequently, that the RIPr $\p^*$ is the standard normal distribution. This example can be easily generalized to $\sigma$-sub-Gaussian distributions for $\sigma \neq 1$, but we omit this for brevity.

\subsection{Testing likelihood ratio bounds}
\label{sec:c5-ex}

Consider the setting in Section~\ref{sec:c2-LRB},
 where for   a fixed probability measure $\p_0$ and $\gamma \in [1,\infty)$,  we test the null hypothesis
$$
\mathcal P =\{\p\in \cM_1: \d \p/\d \p_0 \le \gamma\}
$$
against an alternative $\q \ll \p_0$. 
We have seen in Section~\ref{sec:c2-LRB} that the random variable
$E^*$
  is an e-variable for $\cP$, specified by
$$
E^*= \frac{Z^*\vee z_0}{\gamma},
 $$
 where $Z^*=\d \q /\d \p_0$ and $z_0\ge 0$ is the largest constant such that 
 $ 
  \int_{0}^{1/\gamma} (q_t(Z^*)  \vee z_0) \d t =1. 
$ 
Recall our notation of the left quantile function $q_t(Z)=\inf \{x\in \R: \p_0(Z\le x)\ge 1-t\}$ for $t\in (0,1)$ and any random variable $Z$.


We next show that it is indeed  the numeraire for testing $\mathcal P$ against $\q$.  Let $\p^*\in \cM_+$ be specified by 
 $$ \frac{\d \p^*}{ \d \p_0} =
 \gamma \wedge \frac{\gamma Z^*}{z_0},
  $$
  where $Z^*/z_0$ is set to $0$ if $z_0=0$ and $Z^*=0$,
  and it is set to $\infty$ if $z_0=0$ and $Z^*>0$. 
  By straightforward computations, we can check that   $\q$ is absolutely continuous with respect to $\p^*$, $\p^*(\Omega)\le 1$, and $\d \p^*/\d \p_0\le \gamma$. 
Moreover,
$$\frac{\d \q}{\d \p^*}= \frac{\d \q}{\d \p} \frac{\d \p}{ \d \p^*}\id_{\{ Z^* >0\}} 
=  \frac{Z^* \vee z_0}{\gamma} = E^*.$$ 
Therefore, by Theorem~\ref{th:KL} (it is easy to see that the result holds true even if $\p \in \cM_+$ is not a probability measure),  $E^*$ maximizes $\E^\q[\log X]$ under the constraint $\int X\d \p^* \le 1$.
Since an e-variable for $\mathcal P$ always satisfies the constraint  $\int X\d \p^* \le 1$, we know that $E^*$ is also log-optimal for testing $\mathcal P$ against $\q$, and thus it is the numeraire. 

The special case of $\gamma=1$ and $\q\ll \p_0$ in the above example is precisely Theorem~\ref{th:KL}, with the log-optimal e-variable $E^*= Z^*=\d \q/\d \p_0$ and $z_0$ being the largest constant such that $z_0\le Z^*$ almost surely (typically $z_0=0$).

\subsection{Testing exchangeability}\label{subsec:numeraire-exch}

Consider the problem  of  testing exchangeability in Section~\ref{subsec:c2-numeraire-exch}. 
For the data $Z = (Z_1,\dots,Z_n)$,   
 the null hypothesis of exchangeability  is $$\cP =\{\p: Z \laweq Z_\sigma \text{ under $\p$ for all permutations } \sigma\},$$  and 
 $\q$ is an alternative distribution of $Z$ with density $q$, which is not  in $\cP$. 
Let $\pi$ be uniformly drawn from  $\mathcal S_n$. In Section~\ref{subsec:c2-numeraire-exch}, we have verified that the random variable
$$
E =\frac{q(Z)}{q(Z_\pi)}  = \frac{q(Z)}{   \frac1{n!}\sum_{\sigma\in \mathcal S_n} q(Z_\sigma)}
$$
is an e-variable for $\cP$.
To see that it is 
the numeraire for $\cP$ against $\q$ and further, $q(Z_\pi)$ is the RIPr of $q$ onto $\cP$, we note that the distribution $\p^*$ with density $q(Z_\pi)$ is exchangeable, hence in $\cP$. This then yields the above claim by Theorem~\ref{T_verification}, since $\d\q/\d\p^*$ is the only e-variable which can be represented as a likelihood ratio of $\q$ to some element of the bipolar $\cP^{\circ\circ} \supseteq \cP$.

\section*{Bibliographical note}

This chapter is largely shaped by treatment in~\cite{larsson2024numeraire}.
Historically, the RIPr first appeared in \cite{csiszar1984information}. Foundational work was done by \cite{Li99}, who proved (under some assumptions) its uniqueness and an inequality that (in hindsight) yields the e-variable property of the numeraire, amongst several other facts. \cite{GrunwaldHK19} first described its key role in constructing optimal e-variables. \cite{harremoes2023universal} subsequently extended the original definition to more general settings by relaxing some conditions. Finally, \cite{larsson2024numeraire} eliminated all conditions necessary to define the RIPr, but doing so required defining it via the numeraire e-variable, as done in this chapter. 
 
\cite{GrunwaldHK19} also consider composite alternatives. They first consider a worst-case optimality criterion that they call GROW: finding an e-variable for $\cP$ that maximizes the worst-case e-power over $\cQ$; we find this to be a rather pessimistic objective (such an e-variable tuned for the worst case would not adapt to simpler alternatives that could admit larger growth rates), but it does involve some nice mathematics such as a \emph{joint information projection}. This motivates their REGROW criterion, which is similar to our asymptotic log-optimality objective~\eqref{eq:asymp-log-optimality-RIPr}.

The proof of Lemma~\ref{lem:bipolar-domination} can be found in~\cite{larsson2024numeraire}. 
The numeraire for testing exchangeability in Section~\ref{subsec:numeraire-exch} was obtained by \cite{koning2023post} and~\cite{larsson2024numeraire} using different techniques. The numeraire for null hypotheses that satisfy a different notion of group invariance has been derived in~\cite{perez2022estatistics}. 
The numeraire for one-parameter exponential families was derived in~\cite{grunwald2024optimal}, and then in~\cite{larsson2024numeraire} using different techniques; we follow the latter presentation. Exponential families are treated in much more generality and detail by~\cite{grunwald2024optimal}. \cite{csiszar2003information} has more details on the reverse information projection for exponential families.
The bounded mean example is a variant of one studied in \cite{waudby2020estimating}, with the corresponding duality theory derived explicitly in~\cite{honda2010asymptotically}.

The minimum KL divergence in Theorem~\ref{th:strong-duality-ripr} is closely related to, yet different from, the $\kl_{\inf}$ metric found in the modern multi-armed bandit literature; see~\cite{agrawal2021regret} for a recent example. The exact relationships are yet to be fully worked out. The strong duality in Theorem~\ref{th:strong-duality-ripr} also holds for some other divergences than Kullback--Leibler (corresponding to $\phi(E)$ instead of $\log E$ for some concave $\phi$). See~\cite[Section~6]{larsson2024numeraire} for an example involving R\'enyi divergences, with $\phi(x)=x^{1-\gamma}/(1-\gamma)$ for $\gamma > 1$.

The fact that numeraire e-variables are in general logically incoherent (in contrast to universal inference) was first observed by~\cite{bickel2024discussion}.

\chapter{Sequential anytime-valid inference using e-processes}
\label{chap:eprocess}

Classical sequential testing, as developed by Wald, involves specifying a particular stopping rule that achieves a prespecified bound on the type-I and type-II errors. In contrast, this chapter deals with several central concepts for sequential \emph{anytime-valid} inference (SAVI). The latter is a new paradigm for sequential testing (or estimation) that allows the statistician to stop at any arbitrary stopping time, possibly not anticipated or specified in advance. 

\section{The problematic practice of peeking at p-values} 
\label{sec:c7-optional}

Let us motivate SAVI procedures  by describing the not-so-uncommon practice of \emph{peeking at p-values}, and identify the heart of the problem.

Consider a scientist in a psychology laboratory in a university who  painstakingly collects data, subject by subject, where each subject has to undergo an hour-long experiment. It appears entirely reasonable that they may be interested in examining the data after every few subjects, for example by calculating a p-value. For simplicity, assume that the scientist calculates a p-value $P_n$ after seeing $n$ subjects, so that $P_1,P_2,\dots$ is a sequence of p-values calculated on increasing amounts of data. (Assume for now that the statistical models are well specified so that if the null hypothesis is true, $P_n$ is indeed uniformly distributed.)

Peeking at p-values refers to the practice of stopping data collection when $P_n \leq \alpha$ (such as 0.05, or some other appropriate threshold, e.g., 0.04, if the scientist wants to avoid the appearance of stopping just after the data appears significant at the 0.05 level). Why is this an issue? Because the final sample size of the experiment is random, and in particular it is a stopping time
\[
\tau:=\inf\{n \geq 1: P_n \leq \alpha\}.
\]
The heart of the issue is that even though
\[
 \p(P_n \leq \alpha) \leq \alpha,~\text{for each $n \in \N$,} 
\]
it is to be typically expected that
\[
\p(P_\tau \leq \alpha) > \alpha.
\]
To elaborate, even under the null hypothesis, it would be expected under most typical situations (like a t-test p-value) that $\tau$ is almost surely finite, meaning that one will eventually stop and reject the null, and $\p(P_\tau \leq \alpha) = 1$, meaning that we are predestined to make a type-I error! This is sometimes called ``sampling to a foregone conclusion''.

Of course, the catch is that when the data are reported in a scientific article, it is rarely reported that the data were collected in this fashion (perhaps because the scientist may not even realize that this could create issues). Thus, it may be hard to detect and correct for such issues from the data alone: one really needs to know the statistical protocol for collecting the data to identify this issue. This issue is well known to be a key factor (of many) underlying the ``reproducibility crisis'' in many disciplines, such as psychology and biomedical and social sciences, that was recognized and heavily discussed in the second decade of the twenty-first century~\citep{baker20161}.

Since the desire to peek and analyze the  data as it is being collected is fairly natural (who does not want to potentially stop early and save time and money if that is possible?), it is sensible to provide statistical techniques that can account for such peeking (without really knowing any of its details), rather than seeking to ban it (the current solution, which cannot be enforced). That will be the goal of this chapter: to develop sequential anytime-valid inferential (SAVI) procedures.

We begin by summarizing the dominant paradigm for sequential inference introduced by Wald, which demonstrates the allure of sequential methods and dispels a myth that they are less effective than fixed sample-size procedures because one has to account for the multiple looks at the data.

\section{Wald's sequential probability ratio test}
\label{sec:6-wald}

Before we get into formal definitions, we present a numerical example that can be quite surprising to those unfamiliar with sequential analysis. It convincingly brings out the allure of sequential methods compared to batch methods (those that collect and analyze data as a single batch). 

Suppose we have a coin and wish to test whether it is biased or not. Mathematically, we have access to data $X_1,X_2,\dots$ which are iid $\text{Bernoulli}(p)$ random variables with unknown $p\in (0,1)$. Under the null hypothesis, we suppose $p=1/2$. For the sake of simplicity, under the alternative hypothesis it is 
$p=0.6$. Let the null and alternative densities be represented by $f_0, f_1$ respectively. For a fixed sample size test, the classic Neyman--Pearson lemma dictates that the optimal test is the likelihood ratio test (LRT), which thresholds the likelihood ratio 
\[
L_n = \prod_{i=1}^n  \frac{f_1(X_i)}{f_0(X_i)},
\] 
meaning that we reject the null hypothesis when $L_n$ is larger than some threshold $\gamma$.
Different choices of $n$ and $\gamma$ yield different type-I and type-II errors achieved by the test, but the classic Neyman--Pearson lemma promises that no other test can achieve a better tradeoff of the two errors. If we desire a type-I error of $\alpha=0.05$ and a type-II error of $\beta=0.05$, a numerical calculation shows that one requires a sample of size $n=280$ for our problem setup, and the LRT rejects when $L_{280}$ exceeds $\gamma = 0.962$ (equivalently, when at least 154 heads are observed). Note that $L_n$ is an e-variable, but the Neyman--Pearson threshold is \emph{much} smaller than $1/\alpha = 20$, but is rather a quantile of $L_n$. It is indeed odd that we can reject the null for a threshold smaller than 1, meaning even when the null has greater likelihood than the alternative.

In contrast, Wald's sequential probability ratio test (SPRT) works by sampling the data one by one, continually updating $L_t$. 
Indeed, $(L_t)_{t\in \N}$ is an e-process, which will be formally studied in this chapter, so the SPRT has a lot in common with the SAVI and e-value methodology, but also several differences that are highlighted at the end of this section. 

We describe the SPRT below.
If $L_t$ falls below some threshold $\gamma_0$, we give up and accept the null (more rigorously, we do not reject the null); if it goes above some other threshold $\gamma_1$, we reject the null. Else, we keep sampling. It turns out that $\gamma_0 = \beta$ and $\gamma_1 = 1/\alpha$ suffices to keep the type-I and type-II errors below $\alpha$ and $\beta$ respectively, as a consequence of Ville's inequality, to be introduced later in this chapter. (These choices are more conservative than the approximate thresholds of $\beta/(1-\alpha)$ and $(1-\beta)/\alpha$ often recommended for the SPRT.) It can be shown that this procedure will stop almost surely under the null or alternative. 

In summary, our experiments pitch the Neyman--Pearson LRT rule:
\[
\text{ reject if } L_{280} \geq 0.963,
\]
against a conservative variant of Wald's SPRT rule, which does for $t \geq 1$:
\[
\text{ reject if } L_t \geq 20, \text{ accept if } L_t \leq 0.05, \text{ continue sampling otherwise.}
\]

Numerical experiments show that for the SPRT, the expected sample size (equivalently, expected stopping time) under the null is 137.7, while under the alternative it is 139.3. This is about half the number of samples needed by the LRT, a remarkable improvement over the optimal batch test in such a simple setting!

We can also test the methods out when the data is neither drawn from the null or alternative; this is often called the \emph{misspecified} setting. When the bias of the coin is smaller than anticipated, e.g. $p=0.55$, the LRT has power $0.52$ while the SPRT has power $0.48$, which is quite similar, but the expected sample size of SPRT  is $234.5$. When the bias is larger than anticipated, e.g., $p=0.65$, both tests have power equal to one, but the SPRT expected sample size is just $76.1$.  The numbers are summarized in Table~\ref{tab:c6-wald}.

\begin{table}[t]
    \centering
\resizebox{1\textwidth}{!}{
    \begin{tabular}{lcccc} 
        Data-generating distribution & $p=0.5$  & $p=0.55$  & $p=0.6$ & $p=0.65$ \\
           Setting   &  null & misspecified     & 
                    alternative 
                 & misspecified     \\
        \midrule
        SPRT expected sample size   & 137.7  & 234.5 & 139.3   & 76.1 \\
          LRT  pre-specified sample size   & 280  & 280 & 280 & 280  \\
      Power of SPRT  & 0.04 & 0.52 & 0.96 & 1.00 \\
         Power of LRT & 0.05   & 0.48 & 0.95 & 1.00  \\ 
    \end{tabular}
}
    \caption{A comparison of SPRT and LRT for testing $p=0.5$ against $p=0.6$, with target type-I error rate at $0.05$ and target type-II error rate at $0.05$. All numbers are produced with  Monte Carlo simulations using 10,000 repetitions.}
    \label{tab:c6-wald}
\end{table}

Wald proved that for simple nulls and alternatives, his SPRT provided the optimal tradeoff between type-I error, type-II error and stopping time. He showed that any other stopping rule with lower type-I and type-II errors must have larger expected stopping time. In light of this result, the aforementioned improvement of the SPRT over the LRT is not a coincidence.  We record this below for ease of reference.


\begin{theorem}[Wald's SPRT over Neyman--Pearson's LRT]
    Let $\alpha,\beta$ be the type-I and type-II errors  of the LRT with $n$ samples (using some threshold). Let $\alpha', \beta'$ be the type-I and type-II errors of the SPRT (using some thresholds), and  $\tau$ be its stopping time. 
    If $\alpha \leq \alpha'$ and $\beta\leq \beta'$, then  $\max(\E^\p[\tau],\E^\q[\tau]) \leq n$.
\end{theorem}

Wald's paradigm for sequential testing has been hugely successful: the literature on sequential analysis is large and mature. 
Although using the same process and similar stopping rules, the SAVI goals deviate from Wald's goals in four significant ways:
\begin{enumerate}
    \item \textbf{Decisions versus evidence}: Wald focused on decision theory, that is making a decision quickly with predefined error rates. Our focus here is on designing measures of evidence that can be sequentially updated. Of course, our evidence measures can be used to make decisions, but that is a secondary and not primary goal for us. Said differently,  we do not necessarily make binary decisions (accept or reject). When using SAVI methods to make binary decisions, our rules will resemble Wald's.
    \item \textbf{Single stopping time versus all stopping times:} Wald focused on designing a single stopping time (or rule) at which the scientist would reject or accept the null. We will instead provide a notion of evidence (the e-process) that is anytime-valid in the sense that it has a guarantee (under the null) at any stopping time. When SAVI techniques are applied to the setting of Wald's SPRT, the likelihood ratio e-process gives us a way of continually examining and updating evidence before (and leading up to) Wald's stopping time.
    \item \textbf{Differing optimality criteria:} Our notion of optimality --- log-optimality --- is very different from Wald's notion of optimal tradeoffs of the two errors. For simple nulls and alternatives, it is remarkable that the same statistic (the likelihood ratio martingale) simultaneously optimizes both objectives, but this will in general not be true.
    \item \textbf{Beyond likelihood ratios:} While the SPRT techniques presented above can be extended to certain composite nulls and alternatives, it has typically always dealt with settings where likelihood ratios exist (e.g., parametric settings, with a common reference measure). However, e-processes provide us with a universal language of handling \emph{any} sequential testing problem: we will show that every sequential test (for any composite, possibly nonparametric, null) can be obtained using e-processes.
\end{enumerate}

We now move on to defining the basic concepts in the paradigm of e-processes.

\section{Test supermartingales, sequential tests, and e-processes}

There are four main SAVI tools: e-processes, p-processes, $(1-\alpha)$-confidence sequences and level-$\alpha$ sequential tests; it is immediately apparent that the latter two are associated with a predefined level $\alpha\in [0,1]$ but the former two are not. These are respectively the sequential analogs of e-variables, p-variables, $(1-\alpha)$-confidence intervals and level-$\alpha$ tests. We will primarily focus on e-processes because they are the central object; the others can be easily obtained as a consequence of Ville's inequality applied to one or more e-processes. 

We begin with the technical background necessary for the sequential setting. We assume that the reader is familiar with terminologies in stochastic processes, recapped below.

We work in a fixed filtered probability space, meaning that $\cF$ refers to a filtration $(\cF_t)_{t\in \t}$ (a nested sequence of $\sigma$-algebras), where $\t$ is a finite or infinite set of indices starting from $0$. We only deal with discrete-time processes, so $\t=\N_0$ or $\{0,1,\dots,n\}$ for some $n \in \N$. 
Let $\cF_\infty=\bigcup_{t\in \t}\cF_t.$ 
Further denote by 
 $\T$  the set $\t \setminus \{0\}$.


A sequence of random variables $Y = (Y_t)_{t\in \t}$ is called a \textit{process} if it is adapted to $\cF$---i.e., if $Y_t$ is measurable with respect to $\cF_t$ for every $t$.  A process $Y$ is called \emph{predictable} if $Y_t$ is measurable with respect to $\cF_{t-1}$ for each $t\ge 1$. Inequalities comparing processes are understood to hold at each $t\in \t$.

Often $\cF$ is chosen as the natural filtration of data, that is, $\cF_t := \sigma(X^t)$, where $X^t=(X_1,\dots,X_t)$,  with $\cF_0$ being trivial ($\cF_0=\{\emptyset,\Omega\}$). In this case $Y_t$ being measurable with respect to $\cF_t$ means that $Y_t$ is a measurable function of $X^t$.  
But $\cF$ is sometimes a coarser filtration (we discard information) or a richer one (for example allowing $\cF_0 = \sigma(U_0)$ for a uniformly distributed  random variable $U_0$ that is independent of all other randomness). 

A stopping time (or rule) $\tau$ is a nonnegative integer-valued random variable such that $\{\tau \leq t\} \in \cF_t$ for each $t \geq 0$. In words: we know at each time whether the rule is telling us to stop or keep going. Denote by $\cT$ the set of all stopping times, implicitly depending on $\cF$, including ones that may never stop.
For a stopping time $\tau$, the $\sigma$-algebra $\cF_\tau$ is defined by $\{A\in\cF_\infty: A\cap \{\tau\le t\} \in \cF_t \mbox{ for every $t\in \t$}\}$. 

As in previous chapters, the set of distributions $\cP$ represents our null hypothesis.
A level-$\alpha$ sequential test for $\cP$ is a binary process $\phi$ such that 
\begin{equation}\label{eq:sequential-test}
\sup_{\p \in \cP} \p(\exists t \geq 1: \phi_t = 1) \leq \alpha.
\end{equation}
Interpreting $\phi_t=1$ as rejecting $\cP$,  the above condition means that the probability that we will \emph{ever} falsely reject the null is at most $\alpha$, meaning that with probability $1-\alpha$, $\phi$ equals the sequences of all zeros.
These are also called one-sided or open-ended or power-one tests. 

Without loss of generality, we can assume that if $\phi_t=1$, then $\phi_s = 1$ for all $s \geq t$. Further, we will define $\phi_\infty = \lim_{t \to \infty} \phi_t = \sup_{t \in\N } \phi_t$, meaning that if a test did not reject at any finite time, then it also does not reject at infinity. 

Then, we can more easily identify a sequential test by the stopping time $\tau := \inf\{t \geq 1: \phi_t=1\}$ with the convention $\inf\varnothing =\infty$. Then, a level-$\alpha$ test requires that $\tau < \infty$ with probability at most $\alpha$ under every $\p \in \cP$. We can then rephrase the requirement in~\eqref{eq:sequential-test} as:
\begin{equation}\label{eq:sequential-test-2}
\E^\p[\phi_\tau] \leq \alpha,
\end{equation}
for all $\p \in \cP$ and all $\tau \in \cT$.

In the discrete-time setting of this book, 
an integrable process $M$ is a \emph{martingale} for $\p$ if 
\begin{equation}\label{eq:martingale}
    \E^\p[M_t \mid \cF_{t-1}] = M_{t-1} 
\end{equation}
for all $t \in \t_+$, where all such relations between random variables are understood to hold $\p$-almost surely. $M$ is a \emph{supermartingale} for $\p$ if it satisfies \eqref{eq:martingale} with ``$=$'' relaxed to ``$\leq$''.  
 
\begin{definition}[Test (super)martingales and e-processes]
\label{def:test-super}
\begin{enumerate}
    \item[(i)]  A process $M = (M_t)_{t \in \t}$ is called a \emph{test (super)martingale} for $\cP$ if  for every $\p\in\cP$, it satisfies: (a) $M$ is $\p$-almost surely nonnegative, (b) $M$ is a  (super)martingale under $\p$, and (c) $\E^\p[M_0] \leq 1$. A family of processes $(M^\p)_{\p \in \cP}$ is called a \emph{test (super)martingale family} if each $M^\p$ is a test (super)martingale for $\p$.
   \item[(ii)] A sequence of e-values $E = (E_t)_{t \in \t}$ that is adapted to $\cF$ is called an \emph{e-process} if it satisfies one of two equivalent conditions: (a) $\mathbb{E}^\p[E_\tau] \leq 1$ for any stopping time $\tau \in \cT$ and any $\p\in \cP$; 
(b) there exists a test martingale family $(M^\p)_{\p \in \cP}$ such that $E\le M^\p$, $\p$-almost surely for each $\p \in \cP$. 
\end{enumerate}
\end{definition}

The fact that the above two definitions (ii.a) and (ii.b) for an e-process are equivalent is nontrivial, and its proof is omitted. In item (ii.b), one can require $(M^\p)_{\p \in \cP}$ to be a test supermartingale family instead of a test martingale family, without altering the definition. This is a consequence of Doob's decomposition theorem; we invite the reader to fill in the details. 
It is clear that test (super)martingales are e-processes, but not vice versa. 
We sometimes omit ``for $\cP$'' (as we did above) and it should be clear from context.

\begin{definition}[Asymptotic growth rate and log-optimality]
\label{def:asym-GR}
An e-process $M$ is said to be \emph{consistent} against $\q$ if $M_t \to \infty$, $\q$-almost surely as $t \to \infty$.
The lower limit
 $$
\liminf_{t\to\infty} \frac{\log M_t}{t} 
 $$
   is typically  $\q$-almost surely a constant (which depends on $\q$ and may be infinite), and in that case, we call that constant the \emph{asymptotic growth rate of $M$ against $\q$}. 
\end{definition}

Note that the asymptotic growth rate of any e-process for $\cP$ against $\p \in \cP$ is always nonpositive. A positive asymptotic  growth rate implies consistency, but not vice versa; when an e-process is consistent against $\q$, it need not have a positive asymptotic growth rate against $\q$. But in many simple examples (say involving iid data), the asymptotic growth rate will indeed be positive.

E-processes can be converted to sequential tests using Ville's inequality, which we introduce next.

\section{Optional stopping and Ville's inequality}
\label{sec:OST+Ville}

We first state a cornerstone fact in probability theory, the optional stopping theorem, usually attributed to Doob. A precursor to the optional stopping theorem is the following fact, called the supermartingale convergence theorem:  Any test supermartingale $M$ for $\p$ has a well-defined limit, meaning that the random variable $M_\infty$ exists as the (finite) $\p$-almost sure limit of $M_n$.

\begin{fact}[Optional stopping theorem]
\label{fact:optional-stop}
    If $M$ is a test (super)martingale for $\p$, then for any $\cF$-stopping times $\tau$ and $\sigma \leq \tau$, 
    $$\E^\p[M_\tau \mid \cF_\sigma]  \leq M_\sigma.$$ 
    In particular, by taking $\sigma = 0$, we have $$\E^\p[M_\tau] \leq 1.$$ 
    Consequently, if $M$ is an e-process for $\cP$, then $\E^\p[M_\tau] \leq 1$ for all $\p \in \cP$.
\end{fact}

It is important to note that this particular version of the optional stopping theorem does not have any restrictions placed on the stopping time, and this is primarily due to the nonnegativity of the underlying test supermartingale or e-process. If $M$ were not nonnegative, conditions would have to be placed on the stopping time or boundedness of the process for the same result to hold, but we do not concern ourselves with such variants.

 For our purposes, the relevant version of Ville's inequality can be stated as follows.

 \begin{fact}[Ville's inequality]\label{Ville}
 If $M$ is an e-process for $\cP$, then the following three equivalent statements hold:
\begin{subequations}
\begin{align}
\label{eq:Ville}
\p\left(\exists t \in \N: M_t \geq \frac{1}{\alpha}\right) &\leq \alpha \text{ for every } \p \in \cP \text{ and } \alpha \in [0,1];\\
\label{eq:Ville2}
\Longleftrightarrow \qquad \, \sup_{\p \in \cP} \p\left(\sup_{t\in\N} M_t \geq \frac{1}{\alpha}\right) &\leq \alpha \text{ for every }  \alpha \in [0,1];\\
\label{eq:Ville3}
\Longleftrightarrow \qquad  \sup_{\p \in \cP, \tau \geq 0} \p \left(M_\tau \geq \frac{1}{\alpha}\right) &\leq \alpha \text{ for every } \alpha \in [0,1],
\end{align}
\end{subequations}
where $\tau$ is any $\cF$-stopping time taking values in $[0,\infty]$.
 \end{fact}
 
Clearly, for any fixed $t$, $M_t$ is an e-variable. Thus, $\p(M_t \geq 1/\alpha) \leq \alpha$  holds by Markov's inequality. As such, Ville's inequality is best seen as a time-uniform generalization of Markov's inequality that applies to e-processes rather than e-variables.

Note that \eqref{eq:Ville2} and \eqref{eq:Ville3}  often  only holds with  inequality, 
 but they can sometimes hold with equality.
For example, if we were in continuous rather than discrete time, and $(B_t)_{t \in \R_+}$ is a standard (centered) Brownian motion, then for any nonzero $\lambda$, the process $M = (\exp(\lambda B_t - \lambda^2 t/2))_{t \in \R_+}$ can be checked to be a nonnegative continuous-time martingale.
 (A continuous-time martingale $M=(M_t)_{t\in \R_+}$ satisfies $\E[M_t|\mathcal F_s]=M_s$ for all $0\le s< t$, where $(\mathcal F_t)_{t\in \R_+}$ is an increasing chain of $\sigma$-algebras.)
In this case, Ville's inequality holds with equality.
In general, whether or not equality holds depends on whether $M$ \emph{overshoots} $1/\alpha$ at the stopping time $\tau:=\inf\{t \geq 1: M_t \geq 1/\alpha\}$, where infimum of an empty set equals infinity. A more detailed discussion of the tightness of Ville's inequality is out of the scope of this book. 
 
We also remark that a conditional version of Ville's inequality is also true, though we do not utilize it much. Specifically, if $M$ is a test supermartingale for $\cP$, then \begin{equation}\label{eq:conditional-Ville}
    \sup_{\p \in \cP} \p\left(\left. \exists t \geq s: M_t \geq \frac{M_s}{\alpha} \right| \cF_s \right) \leq \alpha \text{ for every } \alpha \in [0,1].
\end{equation}


\begin{figure}[t]
    \centering
\includegraphics[width=0.95\textwidth]{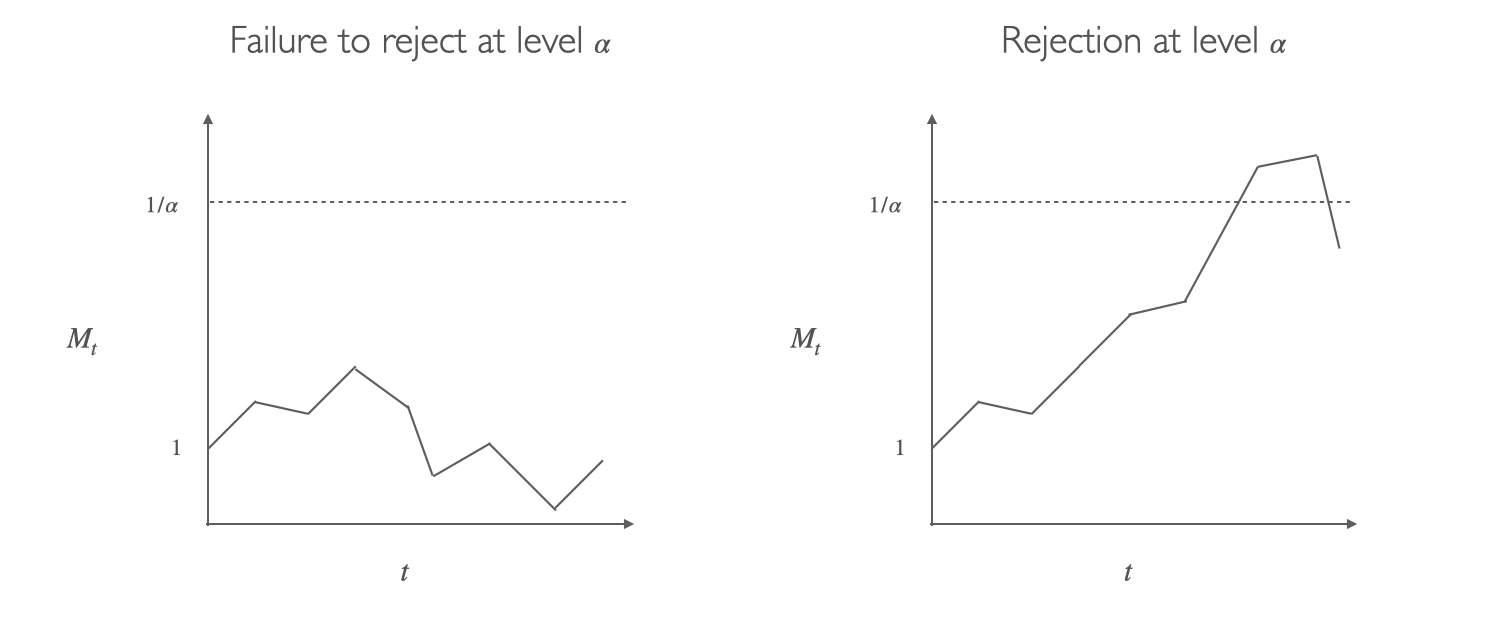}
\caption{A visualization of sequential testing with Ville's inequality applied to an e-process $M$.}
\label{fig:ville}
\end{figure}

 Ville's inequality effectively states that $\phi_t = \id_{\{M_t \geq 1/\alpha\}}$ yields a level-$\alpha$ sequential test; see  
   Figure~\ref{fig:ville} for a simple illustration.
 A converse fact is also true, as presented next.

\begin{proposition}
    Every level-$\alpha$ sequential test $\phi$ for $\cP$ can be obtained by thresholding some e-process $M$ for $\cP$ at $1/\alpha$. 
\end{proposition}
\label{prop:seq-test-via-eprocess}
\begin{proof}[Proof.]
   Define $M_t = \phi_t / \alpha$. Then $M_t$ only takes on two values: 0 and $1/\alpha$. So $M_t \geq 1/\alpha$ if and only if $\phi_t=1$. The only thing left to check is that $M$ is an e-process for $\cP$. Indeed, for any stopping time $\tau$ and any $\p \in \cP$, we have $\E^\p[M_\tau] = \E^\p[\phi_\tau / \alpha] \leq 1$ by~\eqref{eq:sequential-test-2}. This completes the proof.
\end{proof}

Thus, e-processes provide a complete framework for sequential testing. Instead of focusing on the test, we can instead focus on constructing e-processes. However, we do not only think of e-processes as a tool to derive tests. We think of e-processes as measures of evidence against the null: the larger they are (or become), the more evidence we have against the null. We may threshold them to make decisions if we wish to, but they retain interpretability without thresholding as well.

A practical use of  e-processes without thresholding concerns the property of {optional continuation} (discussed in Section~\ref{sec:c4-optional}), meaning that another scientist can continue 
an e-process obtained by a previous scientist (peeking is allowed)  by designing new experiments, and multiplying the new e-process to the value of the previous e-process, to build a new e-process.
Mathematically, this means, if $M^A_s$ is the time-$s$ value of an e-process $M^A$, and $M^B$ is another e-process with $M^B_t=1$ for all $t\le  s$, then $M$ defined by $ M_t=M^A_{s\wedge t} M^B_{t}$ is an e-process. Intuitively, up to time $s$ this process coincides with  $M^A$, 
and after time $s$ the process carries the evidence obtained from $M^A$ and continues with $M^B$.
The same statement holds true if $s$ is not a fixed time, but a stopping time. We put this part formally in the next proposition.

\begin{proposition}
\label{prop:opt-con}
Let $M^A$ and $M^B$ be two e-processes for $\cP$ and $\tau$ a stopping time, such that $M^B_t=1$ on $\{t\le \tau\}$.
Define a process $M$     by $ M_t=M^A_{\tau\wedge t} M^B_{t}$ for $t\in \t$.
Then $M$ is an e-process  for $\cP$.
\end{proposition}
\begin{proof}[Proof.]
For any $\p \in \cP$ 
let $M^{A,\p}$ and $M^{B,\p}$ be two test supermartingales that dominate $M^A$ and $M^B$, respectively, as in Definition~\ref{def:test-super} (ii.b). Clearly, 
$ M_t\le M^{A,\p}_{\tau\wedge t} M^{B,\p}_{t}$.
Below we show that $M^\p:=(M^{A,\p}_{\tau\wedge t} M^{B,\p}_{t})_{t\in\t}$ is a test supermartingale for $\p$, which is sufficient to justify that $M$ is an e-process. 
For $t\ge 1$, noting that $\{\tau<t\}\in \mathcal F_{t-1}$ and   $\id_{\{ \tau <t \}}M^{A,\p}_{\tau } $ is $\mathcal F_{t-1}$-measurable, 
we get \allowdisplaybreaks
\begin{align*}
&\E^\p [M^\p_t\mid \mathcal F_{t-1}] 
\\&= \E^\p [M^{A,\p}_{\tau\wedge t} M^{B,\p}_t \mid \mathcal F_{t-1}] 
\\ &= \E^\p [M^{A,\p}_{\tau\wedge t} M^{B,\p}_t \id_{\{ \tau <t \}}\mid \mathcal F_{t-1}] + \E^\p [M^{A,\p}_{\tau\wedge t} M^{B,\p}_t \id_{\{ \tau  \ge t \}}\mid \mathcal F_{t-1}]  
\\ &=\id_{\{ \tau <t \}}M^{A,\p}_{\tau }   \E^\p [ M^{B,\p}_t \mid \mathcal F_{t-1}] +\id_{\{ \tau  \ge t \}}  \E^\p [M^{A,\p}_{ t}   \mid \mathcal F_{t-1}] 
\\ & \le \id_{\{ \tau <t \}}M^{A,\p}_{\tau }    M^{B,\p}_{t-1}    +\id_{\{ \tau  \ge t \}}   M^{A,\p}_{ t-1}  
\\ &  = \id_{\{ \tau \le t-1 \}}M^{A,\p}_{\tau }    M^{B,\p}_{t-1}    +\id_{\{ \tau > t-1 \}}   M^{A,\p}_{ t-1}  M^{B,\p}_{t-1}  
\\&= M^{A,\p}_{\tau\wedge (t-1)} M^{B,\p}_{t-1} = M^\p_{t-1}.
\end{align*}
Hence, $M^\p$ is a test supermartingale, implying that $M$ is an e-process according to Definition~\ref{def:test-super} (ii.b).
\end{proof}

We make some more observations on the result in Proposition~\ref{prop:opt-con}: First, whether and when the process switches from $M^A$ to $M^B$ depends on data as well as other sources such as personal judgment  after seeing data, which are all included in the filtration $\mathcal F$.
Second, the process $M^B$ does not need to be specified before $\tau$; it only needs to designed after $\tau$ happens.
Third,   switching from one e-process to another can happen multiple times, and the output is always an e-process; this follows by applying Proposition~\ref{prop:opt-con} repeatedly.

We end this section by noting that Ville's inequality~\eqref{eq:Ville2} immediately implies that if $M$ is an e-process, then $(\inf_{s \leq t} 1/M_s)_{t\in\t}$ is a p-process, as defined next. 

\begin{definition}
A nonnegative process
$P=(P_t)_{t\in \t}$ is called a \emph{p-process} if $P_{\tau}$ is a p-variable 
 for any stopping time $\tau \in \cT$ and any $\p\in \cP$.
\end{definition}

Thus an e-process measures evidence based on what we currently have, while a p-process measures evidence based on the best evidence we ever accumulated in the past. Even though all p-processes are not of the form $(\inf_{s \leq t} 1/M_s)_{t\in\t}$ for an e-process $M$, this is certainly a common way to derive p-processes. In this case, it may be clear why p-processes  perhaps provide  an exaggerated sense of evidence: If only the quality, performance and promise of companies could be judged by the maximum wealth they ever had! We present the financial interpretation of e-processes in Section~\ref{sec:6-betting}.

\section{E-processes constructed from   likelihood ratios}
\label{sec:c6-LRP}

\subsection{Log-optimality of likelihood ratio processes}

Suppose that we are testing a simple hypothesis $\p$ versus a simple hypothesis $\q \ll \p$. 
If we observe  iid data $X_1,X_2,\dots$ sequentially, then the likelihood ratio process $M$ given by 
\begin{align}
\label{eq:likelihood-ex}
M_0=1 \mbox{~~~~and~~~~}
M_t=\prod_{i=1}^t \frac{\d \q}{\d \p}(X_i)  \mbox{~for~}  t\in \N
\end{align} is an e-process adapted to the filtration generated by the data, where $\d \q/\d \p(X)$ is the likelihood ratio of $\q$ to $\p$ observing $X$. 
We have encountered this e-process in Section~\ref{sec:6-wald}, used in the SPRT.  Moreover, we can easily see that $M$ is a test martingale for $\p$. 
In fact, even if the data is not iid, the likelihood ratio based on the first $n$ data points, \begin{align}
\label{eq:likelihood-ex-2}
M_n=\frac{\d\q}{\d\p}(X_1,\dots,X_n),
\end{align} 
is still a test martingale for $\p$ (which can be checked by conditioning). Revisiting the example in Section~\ref{sec:LR-e-variable} with the e-variable $E_n= 
  \exp \left( \mu S_n - {n\mu^2}/2\right)$ in \eqref{eq:example-gaussian}, we can check that with $E_0=1$, 
 $$
 \E^\p[E_n|E_{n-1}] = E_{n-1} \E^\p[ \exp (\mu X_{n}-\mu^2/2)]= E_{n-1} \mbox{~~~for $n\in \N$,}
 $$ 
  and hence $E$ is not only an e-process, but also a test martingale. This is true for all likelihood ratio processes in \eqref{eq:likelihood-ex}, including the example in \eqref{eq:example-bernoulli}.

When testing a simple hypothesis $\p$, it is optimal to use a test martingale to construct e-processes: Any e-process for $\p$ can be upper bounded by its Snell envelope under $\p$ (which is a test supermartingale), which in turn can be dominated by the test martingale for $\p$ that is given by its Doob decomposition; we leave the details to the reader.

   In the setting of testing the simple null hypothesis $\p$ against $\q \ll \p$, the likelihood ratio $M_t$ in \eqref{eq:likelihood-ex} is log-optimal among all e-variables using the first $t$ data points  for any  $t\in \N$. This claim follows  from  Theorem~\ref{th:KL}, by noting that $M_t$ is precisely the likelihood ratio of $\q$ to $\p$ based on the first $t$ data points. But a stronger claim can be made about the optimality of the entire \emph{process}, which we do next.

\begin{theorem}[The likelihood ratio process is log-optimal]
\label{th:6-LR-opt}
    Suppose that for testing a simple hypothesis $\p$ versus a simple hypothesis $\q $, the likelihood ratio process $M$ given by \eqref{eq:likelihood-ex-2} exists for every $n$. 
    Then, for any stopping time $\tau$ that is finite $\q$-almost surely, and any e-process $M'$ for $\p$,
   \begin{equation}
    \E^\q[\log M_\tau]\ge  \E^\q[\log M'_\tau].\label{eq:c7-LR-opt-def}
    \end{equation}
\end{theorem} 
\begin{proof}[Proof.]
    Note that $M'_\tau$ is an e-variable because $M'$ is an e-process,
    and $M'_\tau$ is measurable with respect to $\cF_{\tau}$.
Take any event $A \in \cF_\tau$ with $\p(A) = 0$. Then,
\begin{align*}
\q(A) = \E^\q[\id_A] &= \E^\q[\id_A \id_{\{\tau < \infty\}}] 
\\ &= \lim_{n\to\infty}   \E^\q[\id_A \id_{\{\tau \leq n\}}] =\lim_{n\to\infty}  
\E^\p[\id_A \id_{\{\tau \leq n\}} M_n] = 0,
\end{align*}
where the last step follows because $\p(A)=0$ and $M_n$ is finite. 
Thus, $\q \ll \p$ on $\cF_\tau$. Invoking Proposition~\ref{prop:cond-log-opt} completes the proof.
\end{proof}


For testing a simple null $\p$ against a composite alternative $\cQ$, there are two options discussed in Section~\ref{sec:3-mix-plugin}. Both the mixture likelihood ratio in~\eqref{eq:mixture-LR} and the plug-in likelihood ratio in~\eqref{eq:plugin-LR} are actual test martingales and hence e-processes.

For testing composite nulls, the situation is much more complicated, as nontrivial composite martingales may not exist while  nontrivial e-processes may exist. We discuss this case in Section~\ref{sec:univ-eprocess} and later.

 \subsection{Asymptotic log-optimality and the plug-in method}

We now focus on the setting of observing iid data $X_1,X_2,\dots$ sequentially. 
In this subsection, we use
$\p$ for the distribution of
$X_1$ under the null, and 
$\q$ for the distribution of $X_1$ under the alternative ($\q\ne \p$).
Recall that $\mathrm{KL}(\q,\p)$ is the Kullback--Leibler divergence between $\q$ and $\p$
defined in \eqref{eq:kl-def1}.
As in some other parts of the book, we  will slightly abuse the notation here by using $\p$ and $\q$ also for their product measures; 
phrases like ``testing $\p$ against $\q$'' and ``$\q$-almost surely'' refer to the product measures, and it should be clear from the context.

 \begin{proposition}
 \label{prop:c7-largest}
 Consider testing $\p$ against $\q$
in the iid data setting. The likelihood ratio process  
     has  an asymptotic growth rate equal to $\mathrm{KL}(\q, \p)>0$, and hence it is consistent; any other e-process cannot have a larger asymptotic growth rate. 
 \end{proposition}

 \begin{proof}[Proof.]
 Let $M^*$ be the likelihood ratio process 
 and $M$ be another e-process, both for $\p$ against $\q$. 
     The first statement follows from the law of large numbers: as $t\to\infty$, $\q$-almost surely,
      $$
 \frac{\log M^*_t}{t} 
=\frac{1}{t}\sum_{i=1}^t \log \left(\frac{\d \q}{\d \p}(X_i)\right) 
\to \E^\q\left[\log\frac{\d \q}{\d \p}\right] =  \kl(\q,\p).
 $$    
 To show the second statement,
 suppose that  
 $$
   \liminf_{t\to\infty} \frac{\log M_t}{t}
   > \kl(\q,\p)  \mbox{~~~$\q$-almost surely.}
 $$
 Let $\epsilon>0$ and $\tau$ be the first time $t$
 such that 
$  \log M_t  - \log M^*_t >  \epsilon,$
 which is $\q$-almost surely finite by the assumed limits.
 It follows that $\E^\q[\log M_\tau]>\E^\q[\log M^*_\tau] $, which conflicts the conclusion of Theorem~\ref{th:6-LR-opt}.
 \end{proof}

    Since the likelihood ratio process has the largest asymptotic growth rate, 
    we can define a notion of asymptotically log-optimality using the likelihood ratio process as a benchmark. 

\begin{definition}
    \label{def:asym-opt-growth}
      An e-process $M$ is \emph{asymptotically log-optimal} for $\p$ against $\cQ$ if for every $\q \in \cQ$, we have
    $$
    \lim_{t\to \infty} \frac{1}{t}\left(\log M_t- \log M^{\q}_t \right) \ge 0 \mbox{~~ in $L^1$-convergence under $\q$},
    $$
    where $M^{\q}$ is the likelihood ratio process of $\q$ to $\p$.
\end{definition}

The definition of asymptotic log-optimality does not require the iid assumption, although the iid setting is natural because in this case $\log M^\q_t/t$ has a constant and positive $L^1$-limit   $\kl(\q,\p)$. 

Because $\kl(\q,\p)>0$ for every $\q\in \cQ$, any asymptotically log-optimal e-process against $\cQ$ is consistent against every $\q\in \cQ$. Moreover, a mixture of an asymptotically log-optimal e-process
with any other e-process is also 
asymptotically log-optimal.

In  the previous chapters, log-optimality is defined against single $\q$ (the null can be composite), and that is because the log-optimal e-variables (in the sense of Definition~\ref{def:log-opt}) or e-processes (in the sense of \eqref{eq:c7-LR-opt-def})) are determined by $\q$. 
For \emph{asymptotic} log-optimality, we are able to find nontrivial e-processes that work
against a composite alternative $\cQ$, as demonstrated below. 
On the other hand, we present these methods for a simple null in order to define likelihood ratio processes as the benchmark.

    To illustrate asymptotic log-optimality, we now analyze the problem
of testing   parametric models with likelihood ratios discussed in Section~\ref{sec:LR-e-variable}, but with an infinite sequence of iid data $X_1,X_2,\dots$ (we considered finite data in Section~\ref{sec:LR-e-variable}).
Suppose that data follow a distribution  parameterized  by $\theta \in \Theta$ with true unknown value $\theta^\dagger$. 
We assume $\Theta$ is an Euclidean space. 
The null hypothesis $\p$ corresponds to $\theta= \theta_0$
and the alternative hypothesis is $\cQ=\{\q_\theta: \theta \in \Theta_1\}$, where $\Theta_1\subseteq \Theta$ does not include $ \theta_0$. 
Let $\ell$ be the likelihood ratio function given by
 $$
  \ell(x;\theta) = \frac{\d \q_\theta}{\d \p}(x),
 $$
 and let $\theta_i$ be measurable to $\sigma(X_j:j\in[i-1])$ for each $i\in \N$. 
As we see in Section~\ref{sec:LR-e-variable}, the process $M=(M_i)_{i\in\N_0}$  given by
\begin{equation}
\label{eq:c7-asym-log-opt2}
M_0=1 \mbox{~~~and~~~} M_t=\prod_{i=1}^t \ell(X_i;\theta_i) \mbox{ for $t\in \N$} 
\end{equation}
is an e-process for $\p$. 
If $\theta_i=\theta$ for all $i\in \N$, then $M$ is the 
likelihood ratio process, which is log-optimal against  $\q_\theta$ in the sense of Theorem~\ref{th:6-LR-opt}.

For the problem of testing $\p$ against $\cQ$ above,  
it may not be feasible to choose $\theta_i$ deterministically 
because the tester does not know the true value $\theta^\dagger$.
Therefore, it is desirable to design methods that are adaptive,  
hoping for asymptotic log-optimality against the composite hypothesis $\cQ$. 
This is indeed possible, as long as $\theta_i\to \theta$ under $\q_\theta$ in a suitable sense that depends on the parametric model. 
In the next result, we give a sufficient condition based on $L^2$ convergence. 

\begin{proposition}  
\label{prop:c7-asym-opt-param}
For testing $\p $
against $\cQ=\{\q_\theta : \theta\in \Theta_1\}$, 
suppose that for each $x$,
$\theta\mapsto\log \ell(x;\theta) $ 
is concave, and it
has a gradient denoted by $\eta(x;\theta)$. 
If for each $\theta\in\Theta_1$, $\E^{\q_\theta} [\Vert \theta_i- \theta\Vert^2] \to 0$ as $i\to\infty$ 
and $$\sup_{i\in \N}\E^{\q_\theta} \left[\Vert \eta(X_i; \theta_i)\Vert ^2\right] <\infty,$$ then $M$  in \eqref{eq:c7-asym-log-opt2} is asymptotically log-optimal against $\cQ$. 
\end{proposition}


 \begin{proof}[Proof.]
 Take $\theta\in \Theta_1$, and let  $M^*$ be the likelihood ratio process for $\p$ against $\q_\theta$.
 Using the assumed concavity of the log-likelihood function, we get
\begin{align}
     \frac 1 t (\log M_t -   \log M^*_t ) &  = 
\frac 1 t \sum_{i=1}^t 
\left(\log \ell(X_i;\theta_i)- \log \ell(X_i;\theta_i)\right)
\notag
 \\& 
 \ge \frac 1 t \sum_{i=1}^t 
(\theta_i-\theta) \eta(X_i;\theta_i) .
\notag
 \end{align} 
By Cauchy-Schwarz, we have
$$
\left(\E^{\q_\theta}\left [|(\theta_i-\theta) \eta(X_i;\theta_i)| \right]\right)^2
\le \E^{\q_\theta}\left [ \Vert \theta_i-\theta\Vert ^2 \right] \E^{\q_\theta}\left [ \Vert \eta(X_i;\theta_i)\Vert  ^2\right] 
\to 0. 
$$
Using
$$
\E^{\q_\theta}\left[\left|\frac 1 t \sum_{i=1}^t 
(\theta_i-\theta) \eta(X_i;\theta_i) \right|
\right]
\le \frac 1 t \sum_{i=1}^t\E^{\q_\theta}\left[ \left|  
(\theta_i-\theta) \eta(X_i;\theta_i)\right| \right] ,
$$
and 
the fact that the average   of any vanishing real sequence converges to $0$, we conclude 
$$
\liminf_{t\to\infty} \E^{\q_\theta} \left[ \frac{\log (M_t / M^*_t)}{t}  \right] \ge \lim_{t\to\infty} 
\E^{\q_\theta}\left[\left|\frac 1 t \sum_{i=1}^t 
(\theta_i-\theta) \eta(X_i;\theta_i) \right|
\right]
= 0,
$$
which yields the desired asymptotic log-optimality against $\cQ$.
\end{proof}

The concavity of the log-likelihood ratio function  in Proposition~\ref{prop:c7-asym-opt-param} is satisfied by many distributions in the exponential family. 

In the specific  case of normal distributions with mean $\theta\in \R$
as   in Example~\ref{ex:c1-normal}, we have 
$\log \ell(x;\theta)=x\theta - \theta^2/2$, which is concave in $\theta$,
and $\eta (x;\theta)=x-\theta$. 
In this case, the conditions in Proposition~\ref{prop:c7-asym-opt-param} hold as soon as $\theta_i\to \theta$ in $L^2$-convergence under each $\q_\theta$. 
Therefore, the following corollary is obtained.
\begin{corollary}
\label{coro:c7-asym}
    For testing $\p=\mathrm{N}(\theta_0,1)$
against $\cQ=\{\q_\theta= \mathrm{N}(\theta,1): \theta\in \Theta_1\}$, 
if $\E^{\q_\theta} [ (\theta_i- \theta)^2] \to 0$ as $i\to\infty$  for each $\theta\in\Theta_1$, then $M$  in \eqref{eq:c7-asym-log-opt2} is asymptotically log-optimal against $\cQ$.  
\end{corollary}

Corollary~\ref{coro:c7-asym}
justifies some observations made in Figure~\ref{fig:c1-1}.
The condition $\E^{\q_\theta} [(\theta_i- \theta)^2] \to 0$  holds for the MLE, that is, when $\theta_i$ is the sample mean of $X_1,\dots,X_{i-1}$ for the normal model, since $\E^{\q_\theta} [(\theta_i- \theta)^2] =(i-1)^{-1}$ for $i\ge 2$.
Similarly, it also holds for MAP with a normal prior.

 Another sufficient condition for asymptotic log-optimality   
 with consistent estimators
  is when both $\Theta_1$ and the range of the data $X_1$ 
 are compact regions and the derivative $\eta$ of $\theta \mapsto \log \ell(x;\theta)$ is continuous (hence bounded). 
 In this case, if $\theta_i\to \theta$,  $\q_\theta$-almost surely, then 
 we have asymptotic log-optimality. This follows by noting two facts: first,  
$$
 \frac 1 t (\log M_t - \log M^*_t )   = 
\frac 1 t \sum_{i=1}^t (\theta_i-\theta) \eta(X_i;\tilde \theta_i)  $$
for some $\tilde \theta_i$ between  $\theta$ and $\theta_i$ by the mean-value theorem,
and second,
 $
\E^{\q_\theta}[ |(\theta_i-\theta) \eta(X_i;\tilde \theta_i) |]$
vanishes.


\section{Testing by betting}
\label{sec:6-betting}

There is a straightforward way to view test (super)martingales as wealth processes arising from games. The key idea is as follows: 

In order to test a null hypothesis $\cP$, we set up a game of chance, where the statistician (who plays the role of Skeptic from Section~\ref{sec:betting-interpretation}) bets on each observation before observing it. This game must satisfy two properties: 
\begin{itemize}
    \item[(a)] If the null is true, meaning $\p \in \cP$, then no betting strategy should be able to systematically make money (that is, it may make money by chance, but cannot guarantee making money); 
    \item[(b)] if the null is false, then there should exist ``good'' betting strategies that can make money (while one may never be certain of short term profit, there is a way to guarantee making money in the long run).
\end{itemize}
If such a game can be designed, then the wealth of the statistician in the game is directly a measure of evidence against the null: the more the wealth (or more accurately, the larger the ratio of final to initial wealth), the more evidence that the null is likely false.

Let us begin with the simplest example, as discussed in the previous section. We show how to interpret testing $\p$ as a game in Algorithm~\ref{alg:betting-game-simple}. 

\begin{algorithm}
	\caption{Testing by betting: simple null $\p$} 
 \label{alg:betting-game-simple}
	\begin{algorithmic}[1]
            \State Statistician's initial wealth is $W_0 = 1$
		\For {$t=1,2,\ldots$}
				\State Statistician declares bet $S_t: \mathcal X \to [0,\infty]$ such that $$\E^\p[S_t(X) \mid X_1,\dots,X_{t-1}] \leq 1$$
				\State Observe $X_t$
			\State Statistician's wealth is updated as $W_t = W_{t-1} \cdot S_t(X_t)$
		\EndFor
	\end{algorithmic} 
\end{algorithm}

There are three immediate points worthy of comment. First, in the above game, the constraint on $S_t$ is simply that it is an e-variable for $\p$ conditioned on the past. This immediately implies that for any betting strategy, the wealth process $W$ is a test supermartingale for $\p$. 

Second, one may ask: what is the \emph{optimal} bet? Assume for simplicity that the alternative $\q$ is simple, and that $\p \ll \q$. Then, the wealth can grow exponentially fast under $\q$, and so it is natural to optimize the exponent. This means choosing $S_t$ to optimize $\E^\q[\log S_t(X) \mid X_1,\dots,X_{t-1}]$ subject to $S_t$ being an e-variable for $\p$. 
The answer was already derived in Theorem~\ref{th:KL}: it is the likelihood ratio of $\q$ to $\p$, conditioned on $X_1,\dots,X_{t-1}$ (this conditioning is irrelevant if $\p$ and $\q$ are iid product distributions on $\mathcal X^\infty$, but it is relevant if they are not). Thus, the log-optimal wealth process in this game is given precisely by~\eqref{eq:likelihood-ex}.

The last point is that if one wants to test a composite null $\cP$, the only alteration to the protocol is to require that $S_t$ is an e-variable for $\cP$, meaning that the constraint in the above protocol must hold for all $\p \in \cP$. Then, the wealth process $W$ is a test supermartingale for $\cP$.

The reader may now naturally wonder how e-processes can be interpreted in this testing by betting framework, as opposed to simply test supermartingales. The answer is that their game-theoretic interpretation is slightly more involved. The statistician does not play a single game against $\cP$, but instead plays a separate game against each $\p \in \cP$, and accepts its net wealth as the worst wealth amongst all these games. This is made explicit in Algorithm~\ref{alg:e-process}.

\begin{algorithm}[ht!]
	\caption{Testing by betting: composite null $\cP$} 
\label{alg:e-process}
	\begin{algorithmic}[1]
            \For  {$\p \in \cP$}
            \State Statistician's initial wealth in game $\p$ is $W^\p_0 = 1$
            \EndFor
    
		\For {$t=1,2,\ldots$}
				\State Statistician declares bet $S^\p_t: \mathcal X \to [0,\infty]$ such that $$\E^\p[S^\p_t(X) \mid X_1,\dots,X_{t-1}] \leq 1$$
				\State Observe $X_t$
			\State Statistician's wealth in $\p$ is updated as $W^\p_t = W^\p_{t-1} \cdot S^\p_t(X_t)$
                \State Statistician's overall wealth is updated as $W_t = \inf_{\p \in \cP} W^\p_t$
		\EndFor

	\end{algorithmic} 
\end{algorithm}

The above interpretation stems directly from the second of the two equivalent definitions of e-processes. There are some difficulties in the above game: how do we instantiate and play (potentially uncountably) infinitely many games? How do we separately specify a betting strategy in each game? Further, there may be measure-theoretic subtleties associated with the $\inf_{\p \in \cP}$ operation. We sidestep these issues here. For us, the game-theoretic interpretation of e-processes is just that: an interpretation. It helps us understand the relationship to test supermartingales, and indeed recognize why e-processes are a more general concept. However, we will never actually set up and play such a game explicitly in order to derive and use sensible e-processes in practice. Those can be derived more directly, as we shall see next.



\section{The universal inference e-process}
\label{sec:univ-eprocess}

The previous sections have emphasized the fixed-$n$ view, in that the split and subsampled likelihood ratio e-variables were not e-processes. We now note that there is a variant of universal inference that directly yields an e-process. 
Define $E_0 = 1$ and  
\begin{equation}\label{eq:ui-eprocess}
E_t = \prod_{i=1}^t \frac{\hat q_{i-1}(X_i)}{\hat p_{t}(X_i)},
\end{equation}
where $\hat q_{i-1} \in \cQ$ is allowed to depend on $X_1,\dots,X_{i-1}$ and $\hat p_{t} \in \cP$ is the maximum likelihood estimator in $\cP$ based on $X_1,\dots,X_t$.  Note that the numerator can be updated in an online fashion, but all terms in the denominator need to be recalculated after observing each new data point.

\begin{proposition}
   The process $E$ in~\eqref{eq:ui-eprocess} is an e-process for $\cP$.
\end{proposition}
The proof is simple and mimics  the argument for the split likelihood ratio e-variable. By definition of $\hat p_n$ being the maximum likelihood estimator for $\cP$, we see that for each $\p$ in $\cP$, 
\[
E_t \leq \prod_{i=1}^t \frac{\hat q_{i-1}(X_i)}{p(X_i)},
\]
where the right-hand side is a test martingale for $\cP$, thus verifying the defining property of an e-process.

One may view this, and other variants of universal inference, as calculating the ratio  
\[
\frac{\text{out-of-sample alternative likelihood}}{\text{in-sample null likelihood}}.
\] 
To clarify, each $\hat q_{i-1}$ is evaluated on an independent data point  $X_i$, which is outside the sample that was used to choose $\hat q_{i-1}$. In contrast, $\hat p_t$ is evaluated on $X_1,\dots,X_t$, which are all inside the sample that was used to choose it. Using the language of overfitting in machine learning, we overfit under the null by looking at its ``training likelihood'', but we do not overfit under the alternative, by looking at its ``test likelihood''. This is why the threshold of $1/\alpha$ does not depend on $\cP$ and $\cQ$. This is contrast to generalized likelihood ratio methods, which take ratios of maximum likelihoods under alternative and null, but the threshold must then be adjusted (made larger) to take into account the difference in complexities of $\cP$ and $\cQ$.

We remark briefly again that instead of the plug-in method for handling composite $\cQ$, we can use the mixture method instead. In particular, for any distribution $\nu$ over $\cQ$,
\[
E_t = \int_\cQ \prod_{i=1}^t \frac{q(X_i)}{\hat p_{t}(X_i)} \d\nu(q)
\]
is an e-process. Without too much ambiguity, we refer to this also as the universal inference e-process.

\section{Sequential e-values and empirically adaptive e-processes}
\label{sec:EAEP}

In this section, we formally define sequential e-variables (e-values), and a specific form of e-processes built from these  sequential e-variables.
As in the previous section, 
we will use
$\t$ is a finite  or infinite set of indices starting from $0$, and 
 $\T$ is the set $\t \setminus \{0\}$.

\begin{definition}[Sequential e-values]
\label{def:se-values} 
The  e-variables  $E_t$ with $t\in \T$  for $\cP$, adapted to the filtration $\mathcal F$, are \emph{sequential} 
if  $\E^\p[E_t | \mathcal F_{t-1} ]\le1$   for all $t\in \T $ and $\p\in \cP$.  
In case the filtration $\mathcal F$ is not specified, by default we mean the natural filtration; in other words, the  required condition is 
$\E^\p[E_t\mid E_1,\dots,E_{t-1}]\le1$ 
for all $t\in \T$ and $\p\in \cP$.

\end{definition}

A possible interpretation of sequential e-variables  
is that  $E_1,E_2,\dots$ are obtained by laboratories $1,2,\dots$ in this order,
and laboratory $t$ makes sure that its result $E_t$ is a valid e-variable given the previous results $E_1,\dots,E_{t-1}$.  
Note that if sequential e-variables  $E_t$, $t\in \T$ for $\p$  are exact, then they must satisfy $\E^\p[E_t |
 \mathcal F_0]= 1$ for each $t\in \T$ because of the tower property of conditional expectation.

We can build e-processes directly based on sequential e-variables, unspecific to the underlying data or hypothesis.
The first simple observation is that the product process of sequential e-values is always an e-process. 

\begin{proposition}
\label{prop:prod-e-pro}
Let  $(E_t)_{t\in \T } $  be a sequence (finite or infinite) of sequential e-variables for $\cP$. 
Define
$$M_t=\prod_{s=1}^t E_s \mbox{~for $t\in  \T $ and $M_0=1$}.$$
Then $M$ is an e-process for $\cP$.
\end{proposition}
The proof follows by noting that for $\p\in \cP$, $$\E^\p[M_t|\mathcal F_{t-1}]=M_{t-1} \E^\p[E_t|\mathcal F_{t-1}] \le M_{t-1},$$ and  thus $(M_t)_{t\in\t}$ is a  test supermartingale, and thus an e-process.

 More generally, e-processes can be constructed from sequential e-variables through {testing by betting} discussed in Section~\ref{sec:6-betting}. In Chapter~\ref{chap:multiple}, we will see that this is actually the only admissible way of constructing e-processes based only on the sequential e-variables. In what follows, the filtration $\mathcal F$ is the natural filtration of these e-variables.


\begin{definition}[Log-optimal and empirically adaptive e-processes]\label{def:empirical}
Let  $(E_t)_{t\in \T} $  be a sequence (finite or infinite) of sequential e-variables.
    \begin{enumerate}[label=(\roman*)]
    \item  An \emph{e-process  built on} $(E_t)_{t\in 
\T} $ 
is  $M=(M_t)_{t\in\t}$ defined by
\begin{equation}
  M_t =\prod_{s=1}^t \left( (1-\lambda_s) + \lambda_s E_s\right) \mbox{~for $t \in \T$ and }M_0=1,
  \label{eq:e-martingale}  
\end{equation} 
where   for $t \in \T$,
$\lambda_t$ is any measurable function of $E_1,\dots,E_{t-1}$ (i.e., $(\lambda_t)_{t\in \T}$ is predictable with respect to $\mathcal F$) and takes value in $[0,1]$. 
\item For an alternative probability measure $\q$, the \emph{$\q$-log-optimal e-process built on} $(E_t)_{t\in 
\T} $ 
 is  the e-process in \eqref{eq:e-martingale}  such that  
 $\lambda_t$ for $t \in \T$ solves the following equation
  \begin{equation}\label{eq:log-opt-lam}
  \lambda_t \in \argmax_{\lambda \in [0,1]} \E^\q\left[ \log\left(  (1-\lambda) + \lambda  E_t\right)\mid \mathcal F_{t-1}\right].
  \end{equation}
\item For a parameter $\gamma\in (0,1]$, the \emph{empirically adaptive e-process} is  the e-process in \eqref{eq:e-martingale}  such that  $\lambda_t$ for $t \ge 2$ solves the following equation
  $$
  \lambda_t \in \argmax_{\lambda \in [0,\gamma]} \frac 1 {t-1}\sum_{s=1}^{t-1}  \log\left(  (1-\lambda) + \lambda  E_s\right),
  $$ 
  and $\lambda_1=0$.
    \end{enumerate}
\end{definition}

Note that similarly to the product process in Proposition~\ref{prop:prod-e-pro},
the construction \eqref{eq:e-martingale} guarantees that $M_t$ is a supermartingale. 

To interpret the $\q$-log-optimal e-process built on $(E_t)_{t\in\t}$, 
  at every step $t$, $\lambda_t$ is chosen to maximize the e-power in Definition~\ref{def:e-power}, which is allowed to use all past observed e-values.
It
is not necessarily log-optimal among all e-processes testing the alternative $\q$; it is only log-optimal within the class in \eqref{eq:e-martingale}.
If  $E_t$ is independent of $\mathcal F_{t-1}$  under $\q$, then $\lambda_t$ in \eqref{eq:log-opt-lam} is given by maximizing the unconditional expectation,
$$  \lambda_t \in \argmax_{\lambda \in [0,1]} \E^\q\left[ \log\left(  (1-\lambda) + \lambda  E_t\right) \right],$$
which is deterministic. 

The  empirically adaptive martingale does not have a specific alternative, and   at every step $t$, $\lambda_t$ is chosen to maximize the e-power 
  with respect to the empirical distribution of $E_1,\dots,E_{t-1}$; this also explains the choice of the name. The parameter $\gamma$ controls an upper bound for $\lambda_t$, and an uninformative default choice can be $\gamma=1/2$.

  In the construction of the empirically adaptive e-process $M$, for $t\ge 2$,  by Theorem~\ref{thm:nontrivial}, $\lambda_t=0$ if and only if the  empirical mean of $E_1,\dots,E_{t-1}$ is less than or equal to $1$, that is, $\sum_{s=1}^tE_s \le t$.
  It is straightforward to verify that $M$ is a  martingale with respect to the natural filtration of the e-variables if these e-variables are exact; otherwise it is a supermartingale. 
  In particular, $M_\tau$ is an e-variable for any stopping time $\tau$.
  The main advantage of the empirically adaptive martingale is the simple observation that if the e-variables are iid under the alternative hypothesis, then it will have good e-power against the null hypothesis. 
  It generally has good e-power if the iid assumption  holds approximately, and it always has a valid type-I error control at level $\alpha$ for the standard threshold of rejection   $1/\alpha$.
  The following result formalizes the claim about good e-power in an asymptotic setting. 

\begin{theorem}\label{th:adaptive}
Let $(E_t)_{t\in \N}$ be an infinite sequence of  e-variables that is iid under the alternative distribution $\q$ such that $\E^\q[\log E_1] $ is finite.
The empirically adaptive e-process $M$ with $\gamma=1$ satisfies the following:
\begin{enumerate}[label=(\roman*)]
    \item Asymptotic log-optimality of $M$ holds in the following sense:  
    $$
    \lim_{t\to \infty} \frac{1}{t}\left(\log M_t- \log M'_{t}\right) \ge 0 \mbox{~~ in $L^1$-convergence under $\q$}
    $$
    for the $\q$-log-optimal e-process $M'$ built on $(E_t)_{t\in \N}$. 
    \item Consistency  of $M$ holds:  if $\E^\q[E_1]>1$, then $M_t\to \infty$  $\q$-almost surely as $t\to \infty$.
\end{enumerate}
\end{theorem}

Theorem~\ref{th:adaptive} (i) shows that  $M$ has an asymptotic  growth rate that is equal to the $\q$-log-optimal e-process, which is the largest among all choices of e-processes built on $(E_t)_{t\in \N}$. 

We omit a formal proof of Theorem~\ref{th:adaptive}, but its intuition is very simple. Part (i) follows essentially by proving that asymptotically, by the iid assumption, $\lambda_t$ in the empirically adaptive martingale converges in probability to the optimal value $\lambda^*$ that maximizes (over $\lambda$) the e-power of $1-\lambda + \lambda E_1$.  
The e-power of $1-\lambda^* + \lambda^* E_1$ is the asymptotic growth rate of $M_t'$ by the law of large numbers, so the two e-processes share the same  asymptotic growth rate. 
Part (ii) follows by using Theorem~\ref{thm:nontrivial}
to show that the optimal growth rate is positive under the assumption $\E^\q[E_1]>1$, and hence $M_t$, which has the optimal growth rate asymptotically, grows to infinity as $t\to\infty$.  

It should be clear that both statements in Theorem~\ref{th:adaptive} also holds for the 
  $\q$-log-optimal e-process. The advantage of the empirically adaptive e-process is that one does not need to know $\q$.

\section{Obtaining an e-process from e-values using  time-mixture}
\label{sec:c6-time-mix}
 
While the numeraire has been presented as a method for obtaining a single e-variable, one can convert a sequence of numeraires (calculated at different sample sizes $n$) into an e-process using the following simple black-box construction.

For all $n \geq 1$, let $E^{(n)}$ be an e-variable for $\cP$ based on the first $n$ data points $X_1,\dots,X_n$. Take any probability mass function $w$ on the natural numbers ${\mathbb N} = \{1,2, \ldots \}$ and define 
\begin{equation}\label{eq:time-mixture}
M_n = \sum_{j \leq n} w(j) E^{(j)}.  
\end{equation}
with $M_0=1$ by default.
Note that $M$ is an increasing process, which is unusual for e-processes. 

\begin{proposition}
    The time-mixed process $M$ defined in~\eqref{eq:time-mixture} is an e-process for $\cP$.
\end{proposition}
\begin{proof}[Proof.]
    For any $\p \in \cP$ and any random time $\tau$ (possibly, but not necessarily, a stopping time),
  \begin{align*}
    \E^\p[M_\tau]  & = \E^\p\left[\sum_{n=1}^{\infty} {\id}_{\{\tau = n\}} \sum_{j \leq n} w(j) E^{(j)} \right] \\& = \E^\p\left[\sum_{j \geq 1} w(j) E^{(j)} \sum_{n=j}^{\infty} {\id}_{\{\tau = n\}}   \right].
    \end{align*}
    Since $\sum_{n=j}^{\infty} {\id}_{\{\tau = n\}} \leq 1$ and each $E^{(j)}$ is an e-variable, we get
    \[
    \E^\p[M_\tau] \leq \sum_{j \geq 1} w(j) \E^\p\left[E^{(j)} \right] \leq 1,
    \]
    concluding the proof.
\end{proof}

Referring again to the motivation at the start of this section, let $E^{(n)}$ be the numeraire for $\cP$ against some $\q$ based on  $X_1,\dots,X_n$ (more precisely, that means testing $\cP$ against $\q$  conditioned on the $\sigma$-algebra generated by $X_1,\dots,X_n$). In general, $E^{(1)}, E^{(2)}, \dots$ is not an e-process. Now, choose
$$
w(n) = \frac{c}{n (\log n)^2},
$$
for some normalizing constant $c>0$ that ensures $\sum_{n \in \N}w(n)=1$.
we get that 
$$
\log M_n \geq \log E^{(n)} - \log n - 2 \log \log n - O(1).
$$
Thus,
\[\lim_{n \to \infty} \frac1n \log M_n - \frac1n \log E^{(n)} = 0\]
meaning that the asymptotic growth rate of our e-process $M$ matches that of the sequence of numeraire e-variables.

\section{E-processes avoid reasoning about hypothetical worlds}

We end the chapter with a somewhat bizarre example involving a very subtle (mis)use of p-values with optional stopping, whose solution is provided by employing e-processes.

\begin{example}[Quantum-like interference in ordinary experiments]
    Consider a data scientist, Bob, at an IT company. Bob runs an A/B test, which is a two-sample test intended to compare two options A and B of commercial design. Bob fully intends to collect 10,000 observations from the very start. Each observation involves a person choosing between options A and B. The data can thus be seen as a sequence of heads and tails, and let us assume that the null hypothesis is that the coin is fair. At the end of the experiment, Bob finds that his p-value is 0.04, and he is satisfied with that decision, and rejects the null. 
    Later that day, he reports these results to his boss, Alice. Alice is impressed with Bob's statistical rigor, as well as patience. Alice asks Bob how he managed to control his curiosity to look at the data before the experiment was completed. Bob replied that he indeed looked at the data once, after 5000 points were collected, but he simply let the experiment go on. Alice proceeds to ask Bob: ``If you had seen over 4800 tails in the first 5000 tosses, what would you have done?''. Bob replies that he would have of course stopped the experiment. Alice is disappointed, and tells Bob that his reported p-value of 0.04 is incorrect. Bob protests and asks why this  thought experiment in a hypothetical world affects the interpretation of the data what actually happened in the real world. Alice responds that Bob's answer demonstrates that he did not actually have a fixed sample size: he had a data-dependent decision rule, because the sample size across different (hypothetical) runs of the experiment is sometimes 5000 and sometimes 10,000. This changes the p-value calculation, which assumed a fixed sample size of 10,000. A dejected Bob returns to his office to attempt a recalculation of the p-value based on an examination of when he would have stopped in these hypothetical worlds (or to re-run the experiment again, and not peek at the data).

    Similarly, consider an alternate world (in the same setting as above) in which Bob fully intends to collect 5000 samples. He does not peek at the data, and after collecting 5000 samples, he calculates a p-value of 0.04. Alice is once more pleased with Bob's patience and rigor, but she asks him ``Would you have collected more data had the p-value been 0.055?''. If Bob says yes, then Alice immediately says that Bob's actual realized p-value of 0.04 is incorrect. Once more, Bob's answer to a question about a hypothetical world changes the validity of the p-value in the real world. 
\end{example}

While the numbers used in the above example were intentionally extreme in order to make the responses seem believable, what the above example really demonstrates is that reasoning about the validity of p-values requires reasoning about (all possible) ``hypothetical worlds''. Bob's answers to what-if questions in hypothetical worlds affects the validity of his p-value in the real world that did occur. This is clearly not a desirable property, since it is practically infeasible, and indeed one that is often ignored in theory and practice. Despite the small effect on the final p-value that this reasoning may have had in the current experiment (due to the extreme numbers), it points to a fundamental issue that could have a larger effect in other settings, where it may also be more difficult to address.

The effect described above is quantum-like in the sense that merely observing the experiment once before its completion  changes its properties. The mere \emph{possibility} of interference in the experiment seems to have actually interfered in its validity. 

If the final evidence reported was $E_{10,000}$, where $(E_n)_{n \in\N}$ is an e-process, then this issue would have been avoided. Indeed, because 5000 and 10,000 are both valid stopping times, it does not matter that in some other hypothetical world, Bob might have stopped at time 5000. Whenever Bob stops in any world, the reported e-value is a valid e-value.  
Moreover, it is fine to continue sampling if $E_{10,000}$ is not satisfactory to Bob, or to a later researcher.  (Note that we really need $(E_n)_{n \in\N}$ to be an e-process and not just a sequence of e-values.) 

\section{A pathological example for the e-power of e-processes}
\label{sec:path-ex-power}

The e-power has been a central criterion in  the considerations of e-processes so far.
 We close this chapter with two examples, showing that the e-power  has some limitations, similar to the ones discussed in Example~\ref{ex:pathological}, but  in a sequential setting. 
The takeaway message from these examples is 
\begin{enumerate}
    \item[(a)] it is possible to find e-processes with positive  infinite e-power at each fixed time but declines to $0$ with probability $1$ under $\q$;
    \item[(b)]  it is possible to find  e-processes with negative  infinite e-power  at each fixed time but increases to $\infty$ with probability $1$ under $\q$.
\end{enumerate}
Therefore, extra caution should be taken when infinity is involved in the calculation of e-power. 

\begin{example}
\label{ex:patho}
 
 Let $E_1,E_2,\dots,$ be   e-variables for  $\p$ given by
$$
E_t= \exp(t^2-Y_t),~~~t\in \mathbb N,
$$
where $Y_1,Y_2,\dots,$ are iid Pareto(1) random variables under $\q$ (that is,  $\q(Y_1 > x)=1/x$ for $x\ge 1$),
and  they are independent under $\p$
such that $\E^\p[E_t]\le 1$ for each $t\in \mathbb N$. 
Let $M=(M_t)_{t\in \mathbb N}$ be given by  $ 
M_t=\prod_{s=1}^t E_s,
$ which is an e-process.
Simple calculation shows, for any $x>0$,
\begin{align*}
\q(M_t \le x) &= \q\left(\sum_{s=1}^t Y_i  \ge  \frac{t(t+1)(2t+1)}{6} -  \log x\right)
\\ &\le t \, \q\left(Y_1  \ge  \frac{(t+1)(2t+1)}{6} - \frac{ \log x}{t}\right)\to 0 \mbox{~~as $t\to\infty$.}
\end{align*}
Therefore,   $M_t\to \infty$ in probability under $\q$, and the probability of rejecting the null by thresholding $M_t$ approaches $1$.
However,  because  $\E^\q[Y_1]=\infty$, we have $\E^\q[\log E_t]= \E^\q[\log M_t]=-\infty$ for each $t\in \mathbb N$, and  the e-process $M$   has the same e-power $-\infty$ as the zero process.
\end{example}

\begin{example}
\label{ex:patho2}

 Suppose that $E_1,E_2,\dots $ are independent e-variables for  $\p$.
 If $E_1$ satisfies $\E^\q[\log E_1]=-\infty$ like in (a) above  and  $E_t$ has finite e-power for all $t\ge 2$, then the process $M$ given by  $ 
M_t=\prod_{s=1}^t E_s,
$ satisfies $ \E^\q[\log M_t]=-\infty $ for all $t \in \mathbb N$, no matter how fast it grows after $t=1$; in particular, one can easily build $M$ with $M_t\to \infty$ in probability under $\q$.

\end{example}
 
 In either example above, we find an e-process  that is intuitively powerful but has negative infinite e-power. 
 Taking a reciprocal in Examples~\ref{ex:patho}--\ref{ex:patho2},
 we  also find e-processes with positive  infinite e-power but declines to $0$ with probability $1$ under $\q$.

\section*{Bibliographical note}

The sequential probability ratio test described in Section~\ref{sec:6-wald} was introduced in~\cite{Wald45, Wald47}, with its optimality proved in~\cite{wald1948optimum}. The log-optimality of the likelihood ratio process (Theorem~\ref{th:6-LR-opt}) appears in a more restricted form in \cite{KoolenG21}; our result is more general.

Without the term e-process being used, the second definition of an e-process was employed in \cite{howard2020time}, while the first was later proposed in \cite{GrunwaldHK19}. Despite having these new definitions, both papers still primarily focused on test supermartingales obtained as products of e-variables. The equivalence of the two definitions of e-processes was proved in~\cite{ramdas2020admissible} using Snell envelopes and Doob decompositions.  \cite{ramdas2022testing} showed the fundamental difference between e-processes and test supermartingales in their study of testing exchangeability. Certain key geometric concepts like fork-convexity play a key role: test supermartingales for $\cP$ are also test supermartingales for the closed fork-convex hull of $\cP$ but this is not true for e-processes.

The universal inference e-process appears in \citet[Chapter 8]{wasserman2020universal}, but it is a well known statistic in sequential inference. For example it is called the adaptive likelihood ratio in~\citet[\S 5.4.2]{tartakovsky2014sequential}, but they do not appear to define or study e-processes as a more general concept beyond this example. \cite{dixit2023anytime} proved that when using a particular nonparametric mixture technique  (predictive recursion), universal inference is asymptotically log-optimal. 

The modern framework of testing by betting appears in~\cite{Shafer/Vovk:2019}, but its roots can be traced all the way back to~\cite{Ville:1939}. 
The optional stopping theorem was proved by~\cite{doob1953stochastic}. A modern self-contained proof of Ville's inequality~\citep{Ville:1939} can be found in~\cite{howard2020time}. Results on the tightness of Ville's inequality can be found in~\cite{ramdas2020admissible} and~\cite[Proposition 4]{howard2021time}.

The predictable process $(\lambda_t)_{t\in \T}$ in the empirically adaptive e-process is also called the plug-in betting strategy, or the optimal strategy for the growth rate adaptive to the particular alternative (GRAPA), which was  studied by  \cite{waudby2020estimating}.  
Sequential e-variables were first defined in \cite{Vovk/Wang:2021} with respect to their natural filtration.  
Theorem~\ref{th:adaptive} is based on \cite{wang2022backtesting} and its proof can be found in \citet[Theorem 3]{wang2022backtesting}. Example~\ref{ex:patho} was discussed in \cite{wang2024proposer}.

\chapter{Handling multiple e-values}
\label{chap:multiple}

In this chapter, we discuss several methods handling multiple e-values. 
Throughout the chapter, 
let $E_1,\dots,E_K$ be $K$ e-variables, where  $K\ge2$ is a fixed positive integer. 
We write $\cK=[K]$, which will be used in this and the subsequent chapters, wherever we deal with $K$ e-variables.  
The presence of multiple e-values arises in multiple testing  as in Chapter~\ref{chap:compound} but it also can arise in single testing with e-values computed from different parts of the data or data splitting as in Chapters~\ref{chap:ui} and~\ref{chap:eprocess}.
Here, we do not distinguish how they are computed, but focus on handling e-variables for some common hypothesis.  

\section{Merging  e-values under arbitrary dependence}
\label{sec:6-AD}

An important advantage of e-values over p-values is that they are easy to combine.
This is the topic of this section, in which we consider the general case,
without any assumptions on the joint distribution of the input e-variables. Thus, our aim is to design methods to combine arbitrarily dependent e-values. 
The cases of independent e-variables or those  with some specific dependence structures are considered in the next sections. 
In what follows, inequalities between functions are understood to hold everywhere.

\begin{definition}[E-merging functions]
\label{def:e-merg}

\begin{enumerate}
\item [(i)]  An \emph{e-merging function} (of $K$ e-values) is an increasing  Borel function
$F:[0,\infty)^K\to[0,\infty)$
such that for any hypothesis, 
$
  F(E_1,\dots,E_K) $ is an e-variable for any e-variables $E_1,\dots,E_K$.
  \item [(ii)]   An e-merging function  $F$ is \emph{symmetric} if $F(\mathbf e)$ is invariant under any permutation of $\mathbf e$. 
  
\item [(iii)] 
  An e-merging function $F$ \emph{dominates} an e-merging function $G$ if $F\ge G$. 
The domination is \emph{strict}  
if $F\ge G$ and $F(\mathbf{e})>G(\mathbf{e})$ for some $\mathbf{e}\in[0,\infty)^K$.
An e-merging function  is \emph{admissible}
if it is not strictly dominated by any e-merging function.
\item[(iv)] An e-merging function $F$ \emph{essentially dominates} an e-merging function $G$ if,
for all $\textbf{e}\in[0,\infty)^K$,
\[
  G(\mathbf{e}) > 1
  \Longrightarrow
  F(\mathbf{e})\ge G(\mathbf{e}).
\]
\end{enumerate} 
\end{definition}

We will treat $K$ as fixed and will omit the dimension for an e-merging function, which should be clear from the context.   
Similarly to the situation of Section~\ref{sec:2-cali},  the defining  property of $F$ in Definition~\ref{def:e-merg} (i) is required to hold for any hypothesis, but it suffices for $F$ to satisfy the property for one atomless probability measure $\p$  on some $(\Omega,\cF)$; this is formally justified in Appendix~\ref{app:atomless}.

Although an e-variable may take the value $\infty$, this occurs with $0$ probability under the null hypothesis.
Hence, it is safe to set the value of $F$ to $\infty$ if any of its input is $\infty$. 
Therefore, without loss of generality,  
all e-merging functions in this text are defined on $[0,\infty)^K$, and this spares us from specifying terms such as $0\times \infty$ in our analysis.

The increasing monotonicity assumed for e-merging functions reflects the natural assumption that larger individual e-values 
represent stronger  evidence against the null hypothesis, and thus producing a larger e-value. 
\begin{remark}
    An increasing function on $\R^K$ for $K\ge 2$ is not automatically Borel, although it is when $K=1$. An example is given in the discussion after Theorem 4.4 of \cite{Graham/Grimmett:2006}.
\end{remark}

The interpretation of domination is self-evident: We would hope an e-merging function to be large, providing useful e-values. 
The notion of essential domination weakens that of domination in a natural way:
We require that $F$ is not worse than $G$ only in cases where $G$ carries evidence.   A fact about admissibility is that
any e-merging function is dominated by an admissible e-merging function.
For this reason, we are only interested in admissible e-merging functions.

An important example of e-merging function is the 
\emph{arithmetic mean} $\fM_K$, defined by 
\begin{equation}\label{eq:average}
  \fM_K(e_1,\dots,e_K)
  =
  \frac{e_1+\dots+e_K}{K},
  \qquad
  e_1,\dots,e_K \in [0,\infty).
\end{equation}
Another example is the weighted arithmetic mean. 
Let $\Delta_n$ be the standard simplex in $\R^n$  for  $n\in \N$, that is,  $$\Delta_n=\left\{(x_1,\dots,x_n)\in [0,1]^n: \sum_{k=1}^n x_k=1 \right\}.$$ 
It is straightforward to verify that 
for $\boldsymbol \lambda=(\lambda_1,\dots,\lambda_{K+1}) \in \Delta_{K+1}$, the function
\begin{align}
    \label{eq:c7-weighted}
\fM_{\boldsymbol \lambda}:    (e_1,\dots,e_K)\mapsto \sum_{k=1}^K\lambda_k e_k + \lambda_{K+1} ,
\end{align}
is an e-merging function. We will refer to $\fM_{\boldsymbol \lambda}$ as the ${\boldsymbol \lambda}$-weighted arithmetic mean, although strictly speaking it is a weighted average of the input vector $(e_1,\dots,e_K)$ and the constant $1$. 
On the other hand, the geometric mean is also an e-merging function, but it is dominated by the arithmetic average.



We present two main results in this section. First, the arithmetic mean essentially dominates all other symmetric e-merging functions.
Second, all admissible e-merging functions take the form of $\fM_{\boldsymbol\lambda}$.


\begin{proposition}\label{prop:M}
 The arithmetic mean $\fM_K$ essentially dominates any symmetric e-merging function.
\end{proposition}

\begin{proof}[Proof.]
Fix an atomless probability measure $\p $.
  Let $F$ be a symmetric e-merging function.
  Suppose for the purpose of contradiction that there exists $(e_1,\dots,e_K)\in[0,\infty)^K$ such that
  \begin{equation}\label{eq:conj1}
    b:=F(e_1,\dots,e_K)
    > 
    \max\left(\frac{e_1+\dots+e_K}{K},1
    \right)=:a.
  \end{equation} 
  Let $\mathcal S_K$ be the set of all permutations of $\cK$, 
  $\pi$ be randomly and uniformly drawn from $\mathcal S_K$, and $(D_1,\dots,D_K):=(e_{\pi(1)},\dots,e_{\pi(K)})$. 
  Further, let $(D'_1,\dots,D'_K):=(D_1,\dots,D_K)\id_A$,
  where $A$ is an event independent of $\pi$ and satisfying $\p(A) = 1/a$
  (the existence of such random $\pi$ and $A$ is guaranteed for any atomless probability space
  by Lemma~\ref{lem:rich} in the Appendix).
  
  For each $k$, since $D_k$ takes the values $e_1,\dots,e_K$ with equal probability,
  we have $\E^\p[D_k] = (e_1+\dots+e_K)/K$, which implies $\E^\p[D'_k] = (e_1+\dots+e_K)/(Ka) \le 1$.
  Together with the fact that $D'_k$ is nonnegative, we know $D'_k $ is an e-variable for $\p$.
  Moreover, by symmetry,
 \begin{align*}
    \E^\p[F(D'_1,\dots,D'_K)]
   & =
    \p(A)
    F(e_1,\dots,e_K)
    +
    (1-\p(A))
    F(0,\dots,0)
   \\& \ge
    b/a > 1,
  \end{align*}
  a contradiction.
  Therefore, we conclude that there is no $(e_1,\dots,e_K)$ such that \eqref{eq:conj1} holds.
\end{proof}

In particular,  Proposition~\ref{prop:M} implies that if $F$ is an e-merging function that is symmetric
and positively homogeneous
(i.e., $F(\lambda\mathbf{e})=\lambda F(\mathbf{e})$ for all $\lambda>0$; this is satisfied by e.g., $\fM_K$, the geometric mean, the harmonic average, or the minimum),
then $F$ is dominated by $\fM_K$, because  for any $\mathbf e\in [0,\infty)^K$ with $F(\mathbf e)>0$, we have $F(\lambda \mathbf e)\le \fM_K(\lambda e)$ for $\lambda $ large enough, leading to $F(\mathbf e)\le \fM_K(\mathbf e)$.

It is clear that the arithmetic mean $\fM_K$ does not dominate every symmetric e-merging function;
for example, the convex mixtures
of the trivial e-merging function $1$ and $\fM_K$, that is, 
 $
  \lambda + (1-\lambda) \fM_K
$ for $  \lambda\in[0,1]$, 
are pairwise non-comparable with respect to the relation of domination.

In the theorem below, we show that each function $\fM_{\boldsymbol \lambda}$ in \eqref{eq:c7-weighted} is admissible
and that the class \eqref{eq:c7-weighted} is precisely the class of all admissible e-merging functions.
Moreover, every e-merging function is dominated by one of \eqref{eq:c7-weighted}. 

\begin{theorem}\label{th:adm-merg}
For a function $ F:\R_+^K\to\R_+$,
 \begin{enumerate}
 \item[(i)] if $F$ is an  e-merging function, then
  $F\le \fM_{ \boldsymbol \lambda }$ for some  $\boldsymbol \lambda \in \Delta_{K+1}$;  
\item[(ii)] $F$ is an  admissible e-merging function  if and only if     $F = \fM_{ \boldsymbol \lambda }$ for some $\boldsymbol \lambda \in \Delta_{K+1}$. 
\end{enumerate}
\end{theorem}
The proof of Theorem~\ref{th:adm-merg} is put in  Chapter~\ref{chap:combining}, 
as it requires several other advanced technical results.  

\section{Independent e-values and their products}
\label{sec:6-ie}

In this section we consider merging functions for independent e-values.

\begin{definition}[Ie-merging functions]
\label{def:ie-merg}
\begin{enumerate}
\item [(i)]
An \emph{ie-merging function} is  a Borel function
$F:[0,\infty)^K\to[0,\infty)$
such that for any hypothesis, 
$
  F(E_1,\dots,E_K) $ is an e-variable for any independent e-variables $E_1,\dots,E_K$.
  \item[(ii)]An ie-merging function $F$ \emph{weakly dominates} an ie-merging function $G$ if,
for all $e_1,\dots,e_K$,
\[
  (e_1,\dots,e_K)\in[1,\infty)^K
  \Longrightarrow
  F(e_1,\dots,e_K)\ge G(e_1,\dots,e_K).
\]
  \end{enumerate}
\end{definition}

The corresponding definitions of domination, strict domination, and admissibility 
are obtained from  Definition~\ref{def:e-merg} by replacing ``e-merging'' with ``ie-merging''; this also applies to the later definition of se-merging functions.

 One subtle distinction of the above definition from e-merging functions in the previous section is that we do not require increasing monotonicity in all arguments for the ie-merging function. 
This relaxation is useful in the broader class of se-merging functions, which will be explained in Section~\ref{sec:c7-se-merging}.

Similarly to the arithmetic mean for arbitrarily dependent e-values, there is an important ie-merging function: the \emph{product} $\Pi_K$, defined by 
\begin{equation}\label{eq:product}
\Pi_K  (e_1,\dots,e_K)
 =
  \prod_{k=1}^K
  e_k, \qquad
  e_1,\dots,e_K \in [0,\infty).
\end{equation}

For weak dominance, we require that $F$ is not worse than $G$
if all input e-values carry evidence. This requirement is very weak because, especially for a large $K$, it is rarely the case to have all input e-values larger than $1$.

 The following proposition gives a simple condition for admissibility. 

\begin{proposition}\label{prop:IPK}
Fix an atomless probability measure $\p$. 
  For an increasing Borel function $F:[0,\infty)^K\to [0,\infty)$,
  if $\E^\p[F(\mathbf{E})]=1$ for all vectors $\mathbf E$ of exact e-variables
  (resp., for all vectors $\mathbf E$ of  independent exact e-variables),
  then $F$ is an admissible e-merging function
  (resp., an admissible ie-merging function).
  For the statement on the ie-merging function, increasing monotonicity of $F$ is not needed.
\end{proposition}

\begin{proof}[Proof.]
  It is obvious that $F$ is an e-merging function (resp., ie-merging function).
  Next we show that $F$ is admissible.
  Suppose for the purpose of contradiction that there exists an ie-merging function $G$ such that $G\ge F$ and 
  $
    G(e_1,\dots,e_K)>F(e_1,\dots,e_K)
  $
  for some $(e_1,\dots,e_K)\in [0,\infty)^K$.
  Take a vector $(E_1,\dots,E_K) $  of independent exact e-variables 
  such that $\p((E_1,\dots,E_K)=(e_1,\dots,e_K))>0$.
  Such a random vector is easy to construct by considering any distribution with a positive mass on each of $e_1,\dots,e_K$.
  Then we have
  \[
    \p(G(E_1,\dots,E_K) > F(E_1,\dots,E_K)) > 0,
  \]
  which implies
  \[
    \E^\p[G(E_1,\dots,E_K)] > \E^\p[F(E_1,\dots,E_K)] = 1,
  \]
  contradicting the assumption that $G$ is an ie-merging function.
  Therefore, no ie-merging function strictly dominates $F$.
  Noting that an e-merging function is also an ie-merging function,
  admissibility of $F$ is guaranteed under both settings.
\end{proof}

If $E_1,\dots,E_K$ are independent e-variables, their product $E_1\dots E_K$ will also be an e-variable.
This is the analogue of Fisher's  method for p-values, 
according to the rough relation $e\sim1/p$ mentioned in  Section~\ref{sec:2-cali}. 
Fisher's method will be discussed in Chapter~\ref{chap:combining}.
The ie-merging function $\Pi_K$
is admissible by Proposition~\ref{prop:IPK}.
More generally, we can see that the U-statistic functions $U_n$, defined by
\begin{equation}\label{eq:Usf}
  U_n(e_1,\dots,e_K)
  =
  \frac{1}{\binom{K}{n}}
  \sum_{A\subseteq\cK:|A|=n}
\left( \prod_{k\in A}e_k \right),
  \quad
  n\in\{0,1,\dots,K\},
\end{equation}
and their convex mixtures are ie-merging functions.
Notice that this class includes product (for $n=K$),
arithmetic average $\fM_K$ (for $n=1$),
and constant 1 (for $n=0$).
Proposition~\ref{prop:IPK} implies that the U-statistic functions in~\eqref{eq:Usf}
and their convex combinations are admissible ie-merging functions.


Let us now establish a very weak counterpart of Proposition~\ref{prop:M}
for independent e-values; 
on the positive side it will not require the assumption of symmetry needed for Proposition~\ref{prop:M}.

\begin{proposition}\label{prop:M-i}
  The product  ie-merging function $\Pi_K$ weakly dominates any ie-merging function.
\end{proposition}

\begin{proof}[Proof.]
  Indeed, suppose that there exists $(e_1,\dots,e_K)\in[1,\infty)^K$ such that
  \[F(e_1,\dots,e_K)>\Pi_K(e_1,\dots,e_K)=e_1 \cdots e_K.\]
  Fix an atomless $\p\in \Pi$.
  Let $E_1,\dots,E_K$ be independent random variables
  such that each $E_k$ for $k\in\cK$ takes values in the two-element set $\{0,e_k\}$
  and $E_k=e_k$ with probability $1/e_k$.
  Then each $E_k$ is an e-variable but
  \begin{align*}
    \E^\p[F(E_1,\dots,E_K)]
    &\ge
    F(e_1,\dots,e_K)
    \p(E_1=e_1,\dots,E_K=e_K)\\
    &>
    e_1\cdots e_K
    (1/e_1)\cdots(1/e_K)
    =
    1,
  \end{align*}
  which contradicts $F$ being an ie-merging function.
\end{proof}
  A natural question is whether the convex mixtures of \eqref{eq:Usf} form a complete class of admissible ie-merging functions.
  They do not, as shown by the following example.
\begin{example}\label{ex:not-complete}
Define a function $f:[0,\infty)^2\to[0,\infty)$ by
  \[
    f(e_1,e_2)
    =
    \frac12
    \left(
      \frac{e_1}{1 + e_1}
      +
      \frac{e_2}{1 + e_2}
    \right)
    \left(
      1 + e_1 e_2
    \right).
  \]
 Using Proposition~\ref{prop:IPK}, one can check that $f$
  is an admissible ie-merging function.
Nevertheless, it is different from any convex mixture of \eqref{eq:Usf}.
\end{example}

The assumption of the independence of e-variables $E_1,\dots,E_K$
is not necessary for the product $E_1\cdots E_K$ to be an e-variable. 
For instance, if the e-variables $E_1,\dots,E_K$ 
are \emph{negative lower orthant dependent}, meaning that   for all $ (x_1,\dots,x_K) \in \R^K$, 
 \begin{equation}\label{eq:def-nlod}
\p\left(\bigcap_{k\in \cK } \{E_k\le x_k\}\right)\le \prod_{k\in  \cK  } \p\left(  E_k\le x_k\right) ,
 \end{equation} 
then the product $\Pi_K$ and any U-statistic functions  in  \eqref{eq:Usf} of these e-variables are also e-variables. Clearly, if the central inequality in \eqref{eq:def-nlod}
is replaced by an equality,
then we get independence.
In the next section, we consider a particularly useful class of e-variables for which the product and the  U-statistic functions  are valid merging functions. 

\section{Sequential e-values and martingale merging functions}
\label{sec:c7-se-merging}

 The notion of sequential e-variables generalizes that of  independent e-variables, and it appears in several natural contexts that we will explain later. 
 For conditional expectations, since they are defined in an almost-sure sense, we omit ``almost surely'' in their statements.

\begin{definition}[Sequential e-values and se-merging functions]
\label{def:se-merg}~
\begin{enumerate}
\item [(i)]
 The e-variables $E_1,\dots,E_K$ for $\cP$ are \emph{sequential}
if  $$\E^\p[E_k\mid E_1,\dots,E_{k-1}]\le1$$   for all $k\in\cK$ and $\p\in \cP$.
\item [(ii)]
An \emph{se-merging function}  is a Borel function
$F:[0,\infty)^K\to[0,\infty)$
such that for any hypothesis, 
$
  F(E_1,\dots,E_K) $ is an e-variable for any sequential e-variables $E_1,\dots,E_K$. 
  \item[(iii)] 
  An se-merging function $F$ is \emph{exact} if  $F(E_1,\dots,E_K)$ is an exact e-variable for all sequential exact e-variables $E_1,\dots,E_K$.
  \end{enumerate}
\end{definition}

Sequential e-variables can also be defined with a general filtration $(\cF_t)_{t\in \{0,1,\dots,K\}}$ instead of the one generated by $E_1,\dots,E_K$; this is treated in Chapter~\ref{chap:eprocess}.

It is straightforward to check that all convex mixtures of \eqref{eq:Usf} are se-merging functions.
Since independent e-variables are sequential, se-merging functions also produce a valid e-variable for independent e-variables. 
  Example~\ref{ex:not-complete} gives an ie-merging function but not an se-merging function.

Among se-merging functions, the convex mixtures of \eqref{eq:Usf} are admissible. 
They are also admissible in the   class of ie-merging functions by Proposition~\ref{prop:IPK}. 
We also note that it suffices for $E_1,\dots,E_K$ to be sequential in any one of the $K!$ possible orders for these merging methods 
to be valid.

A particular class of se-merging functions, larger than \eqref{eq:Usf}, is defined below. 
\begin{definition}[Martingale merging functions]
\label{def:martingale-merging}
A function $F:[0,\infty)^K\to [0,\infty)$ is a \emph{martingale merging function} if 
\begin{align}
\label{eq:c7-MMF}
F(e_1,\dots,e_K)= \prod_{k=1}^K (1-\lambda_{k} +\lambda_{k} e_k),
\end{align}
where $\lambda_k:=\lambda_k(e_1,\dots,e_{k-1})$ for $k\ge 2$ is a function of $e_1,\dots,e_{k-1}$ taking values in $[0,1]$, and $\lambda_1\in [0,1]$ is a constant. 
\end{definition}

One can verify that, although not completely obvious, the arithmetic mean, the product and other U-statistic functions in \eqref{eq:Usf} are  indeed all martingale merging functions. We leave these as exercises for the reader.

A martingale merging function in \eqref{eq:c7-MMF} has a clear betting interpretation: At each round $k$, one uses $\lambda_k$ proportion of the capital to bet on the next e-value. With this interpretation, the arithmetic mean corresponds to betting a fixed amount $1/K$ (but not a fixed proportion) at each round. 
The e-process $(M_t)_{t\in T}$ built on $(E_t)_{t\in \T}$ in Definition~\ref{def:empirical} is constructed by merging e-values with a martingale merging function in Definition~\ref{def:martingale-merging}  for each fixed $t$ and also in a consistent way across different $t$. 

The above connection to betting  in Chapter~\ref{chap:eprocess} also reveals that 
non-monotonicity is quite natural under the betting interpretation: 
 If we gain evidence in an early round,
 then we may reduce our bet in the next round, 
which leads to non-monotonicity of the resulting merging function (see Example~\ref{ex:hit-and-stop} below).

It turns out that  martingale merging functions play a  special role among all se-merging functions.    But let us first verify that martingale merging functions are exact se-merging functions.

\begin{proposition}
    Any martingale merging function is an exact se-merging function. 
    Moreover, a convex combination of martingale merging functions is again a martingale merging function. 
\end{proposition}
\begin{proof}[Proof.]
Let $F$ be a martingale merging function. 
For sequential exact e-variables $E_1,\dots,E_K$ for $\p$, by writing $L_K=\lambda_K(E_1,\dots,E_{K-1})$
and $\cF_{K-1} =\sigma(E_1,\dots,E_{K-1})$, we have 
\begin{align*}
&\E^\p[F(E_1,\dots,E_K)\mid \cF_{K-1}]
\\& = F(E_1,\dots,E_{K-1},1)  (1-L_K +L_K\E^\p [E_K| \cF_{K-1} ])\\
&
\le F(E_1,\dots,E_{K-1},1)  .
\end{align*}
Therefore, 
$\E^\p[F(E_1,\dots,E_K)]
\le \E^\p[F(E_1,\dots,E_{K-1},1)].$
Using an induction, we get 
$$\E^\p[F(E_1,\dots,E_K)]\le  \E^\p[F(1,\dots,1)]=1.$$
This shows that the martingale merging function $F$ is an se-merging function.
If the e-variables are exact, then the above inequalities are equalities,
which shows that $F$ is exact. 

  The last statement follows by noting    the following equivalent formulation of a martingale merging function:
$F\ge 0$ satisfies, for some measurable functions $t_k:[0,\infty)^{k}\to[0,\infty)$, $k\in \cK$,
that for each $k\in \cK$,
  \begin{align}\label{eq:equivalent}
    &F(e_1,\dots,e_k,1,\dots,1)
    \\&=
  F(e_1,\dots,e_{k-1},1,\dots,1)
    +
    t_k(e_1,\dots,e_{k-1})
    (e_{k} - 1),
  \end{align}
  and $F(1,\dots,1)=1.$
  A convex combination of several choices of $F$ is the same as a convex combination of several choices of $(t_k)_{k\in \cK}$.
\end{proof}

The next result shows that the class of martingale merging functions is
precisely that of all admissible se-merging functions. 
\begin{theorem}\label{th:mart-merg}
  Any se-merging function is dominated by a martingale merging function.
\end{theorem}
The proof of Theorem~\ref{th:mart-merg} is quite technical and not pursued in this book.

Different from the merging functions obtained in Sections~\ref{sec:6-AD} and~\ref{sec:6-ie}, 
a useful martingale merging function need not be increasing in all arguments.
This is because $\lambda_k$ 
is generally not increasing or decreasing in its arguments.
Although monotonicity does not hold,
any martingale merging function $F$ satisfies a property of \emph{sequential monotonicity}:
for fixed $k\in \cK$ and  $(e_1,\dots,e_{k-1})\in [0,\infty)^{k-1}$,
the function $e_{k}\mapsto F(e_1,\dots,e_{k-1}, e_k,1,\dots,1)$ is increasing.
Intuitively it means that, for given $k-1$ e-values, the overall combined e-value is larger if the next observed e-value  $e_k$ is larger, assuming that all future e-values are $1$. Treating future e-values as $1$ is equivalent to using $\lambda_j=0$ for all $j>k$, or not seeing these e-values at all. 

\begin{example}[Hit and stop]\label{ex:hit-and-stop}
  The non-monotonicity of $F$ appears naturally in a ``hit-and-stop'' strategy:
  for a fixed $\alpha\in (0,1)$ and for each $k\in\cK$,
  if $F(e_1,\dots,e_k,1,\dots,1) \ge 1/\alpha$, then we choose $\lambda_{k+1}=0$
  (which implies $\lambda_{j}=0$ for all $j\ge k+1$);
  otherwise we choose $\lambda_{k+1}>0$. 
  It is clear from \eqref{eq:c7-MMF} that $F$ is not increasing.
  This strategy corresponds to stopping an e-process in Chapter~\ref{chap:eprocess} as soon as it can make a rejection at level $\alpha$ by up-shooting $1/\alpha$. 
\end{example}

\begin{example}[Merging via the empirically adaptive e-process]
\label{ex:eMMF}
The empirically adaptive e-process in Definition~\ref{def:empirical}
induces a martingale merging function.  For a parameter $\gamma\in (0,1]$, the martingale merging function   is defined as by \eqref{eq:c7-MMF} with 
$\lambda_1=0$ and $\lambda_k$ for $k\in \{2,\dots,K-1\}$ are given by
  $$
  \lambda_k (e_1,\dots,e_{k-1}) =\argmax_{\lambda \in [0,\gamma]} \frac 1 {k-1}\sum_{s=1}^{k-1}  \log\left(  (1-\lambda) + \lambda  e_s\right).
  $$ 
  This  martingale merging function is uniquely determined by its parameter $\gamma.$  
\end{example} 
Advantages of the 
empirically adaptive martingale merging function in Example~\ref{ex:eMMF}
are justified in Theorem~\ref{th:adaptive} via its corresponding e-process. This function can be used as a default choice without prior information among se-merging functions.

\section{Mean-variance trade-off}

We now pay special attention to the product function $\Pi_K$ in \eqref{eq:product}, which is both an ie-merging function and an se-merging function.
Proposition~\ref{prop:IPK} shows that this function is optimal in a weak sense among all ie-merging functions (hence among all se-merging functions).
However, different from the case of the arithmetic mean for arbitrarily dependent e-values treated in Section~\ref{sec:6-AD},  this does not suggest 
that $\Pi_K$ is the best merging function to use for independent  or sequential e-values. 
Indeed, $\Pi_K$ has an undesirable property, which also highlights the usefulness of other martingale merging functions. 

We first present a simple lemma that is useful in the proof of the next result.
\begin{lemma}\label{lem:product-1}
  For any ie-merging function $F:[0,\infty)^K\to [0,\infty)$ and $a\ge 1$, 
  the function $G:(e_1,\dots,e_K )\mapsto F(a e_1,e_2,\dots,e_K)/a$ 
  is an ie-merging function.
\end{lemma}
 
\begin{proof}[Proof.]
Fix an atomless probability measure $\p$.
  Take an event $A$ with $\p(A)=1/a$.
  For any independent e-variables $E_1',\dots,E_K'$,
  we can take independent e-variables $E_1,\dots,E_K$ independent of $A$  
  such that $(E_1,\dots,E_K)$ is distributed identically to $(E'_1,\dots,E_K')$,
  since the probability space is atomless (guaranteed by Lemma~\ref{lem:rich}).
  Therefore, it suffices to show $\E^\p [G(E_1,\dots,E_K)]\le1$
  for independent e-variables $E_1,\dots,E_K$ independent of $A$.
  Let $E_1' := a E_1\id_A$; this is an e-variable.
  We can compute
  \begin{align*}
    \E^\p[G(E_1,\dots,E_K)]
    &=
    \frac 1 a \E^\p[F(a E_1,E_2,\dots,E_K)]
   \\ & \le
    \E^\p[F(E_1',E_2,\dots,E_K)]
    \le
    1,
  \end{align*}
  where the first inequality follows from $F\ge0$.
  Hence, $G$ is an ie-merging function.
\end{proof}

Now we can compare the product function and other se-merging functions
with respect to the moments of the resulting e-variable.
Recall that $\E^\q[E]>1$ 
is the defining condition  that $E$ is powered against $\q$.

\begin{proposition}\label{prop:trade-off}

Let $F:[0,\infty)^K\to [0,\infty)$ be an se-merging function. 
  For independent nonnegative random variables $E_1,\dots,E_K$ under $\q$ with $\E^\q[E_k]\ge 1$, $k\in\cK$,
  we have, for every $m\in [1,\infty)$, \begin{equation}\label{eq:2nd}
    \E^\q \left[F(E_1,\dots,E_K)^m \right]
    \le
    \prod_{k=1}^K \E^\q [E_k^m]=\E^\q\left[\Pi_K(E_1,\dots,E_K)^m\right].
  \end{equation}
  In particular, 
  if $F$ is exact and $E_1,\dots,E_K$ are independent exact e-variables for $\p$, then \begin{equation}\label{eq:3rd} \var^\p(F(E_1,\dots,E_K ))
    \le
    \var^\p(\Pi_K(E_1,\dots,E_K)).
  \end{equation}
\end{proposition}

\begin{proof}[Proof.] 
  First, we argue that it suffices to show \eqref{eq:2nd}
  for  $E_1,\dots,E_K$ being exact e-variables for $\q$. 
  For independent nonnegative $X_1,\dots,X_K$ with mean larger than or equal to $1$,
  set $a_k:=\E^\p[X_k]\ge1$, $E_k:=X_k/a_k$ for $k\in\cK$,
  and $a:=\prod_{k=1}^K a_k$.
  Let $G:(e_1,\dots,e_K )\mapsto F(a_1 e_1, \dots,a_K e_K)/a$.
  Clearly, $E_1,\dots,E_K$ are independent exact e-variables. 
  Using Lemma~\ref{lem:product-1} repeatedly,
  we deduce that $G$ is an se-merging function.
  For $m\in [1,\infty)$, if \eqref{eq:2nd} holds for all exact e-variables,
  then 
  \begin{align*}
    \E^\p\left[F(X_1,\dots,X_K)^m\right]
    &=
    a^m
    \E^\p\left[\left(\frac{F(a_1E_1,\dots,a_KE_K)}{a}\right)^m\right] \\
    &=
    a^m
    \E^\p\left[G(E_1,\dots,E_K)^m\right] \\
    &\le
    a^m
    \E^\p\left[\Pi_K(E_1,\dots,E_K)^m\right] \\
    &=
    \E^\p \left[\Pi_K(X_1,\dots,X_K)^m\right].
  \end{align*}
  Therefore, the general case of \eqref{eq:2nd}
  follows from the case of exact e-variables.

  Let $E_1,\dots,E_K$ be independent exact e-variables;
  we will show
  \begin{equation}\label{eq:4th}
    \E^\q\left[F(E_1,\dots,E_K)^m \right]
    \le
    \prod_{k=1}^K \E^\q [E_k^m].
  \end{equation}
  Using Theorem~\ref{th:mart-merg}, $F$ is dominated by a martingale merging function.
  Hence, it suffices to show the proposition for a martingale merging function $F$. 
  We show the proposition by induction.
  Note that $\E^\q[E_1^m]\ge \E^\q[E_1]^m\ge  1$. Moreover, this expectation is increasing in $m$ because for $m>t$, $\E^\q[E_1^m] \ge (\E^\q[E_1^{t}])^{m/t} \ge \E^\q[E_1^{t}]$. 
  The inequality \eqref{eq:4th} holds for $K=1$ since,  by  convexity of $x\mapsto x^m$, we have
  \[
    \E^\q[(1-\lambda +\lambda E_1)^m]
    \le (1-\lambda) + \lambda  \E^\q [E_1^m]
    \le \E^\q [E_1^m]
  \]
  for all $\lambda\in[0,1]$.
  To argue by induction, suppose that
  \begin{equation}\label{eq:5th}
    \E^\q[G(E_1,\dots,E_{K-1})^m] \le \prod_{k=1}^{K-1}\E^\q[E_k^m]
  \end{equation}
  for every se-merging function $G:[0,\infty)^{K-1}\to [0,\infty)$.
  Since $F$ is a martingale merging function, 
  we can write $$F(e_1,\dots,e_K) = G(\mathbf e_{[K-1]}) \left(1-\lambda_K(\mathbf e_{[K-1]}) + \lambda_K(\mathbf e_{[K-1]}) e_K
  \right),$$
  for some $\lambda_K:[0,\infty)^{K-1}\to [0,1]$, where $\mathbf e_{[K-1]}=(e_1,\dots,e_{K-1})$.
  Let us   write $\mathbf Y= (E_1,\dots,E_{K-1})$
  and $L=\lambda_K(E_1,\dots,E_{K-1})$. 
We have
  \begin{align*}
    \E^\q[F(E_1,\dots,E_{K})^m \mid \mathbf Y]
    &=
    \E^\q[G(\mathbf Y)^m( 1- L  +  L  E_K )^m \mid \mathbf Y]  
    \\ 
&    \le G(\mathbf Y)^m  \left ( 1- L  +  L  \E^\q [E_K ^m | \mathbf Y]\right)  
  \\
    &\le  G(\mathbf Y)^m \E^\q[E_K^m].
  \end{align*}
  As a consequence,
  \[
    \E^\p[F(E_1,\dots,E_{K})^m]
    \le
    \E^\p[G(\mathbf Y)^m] \E^\p[E_K^m]
    \le
    \prod_{k=1}^K \E^\p[E_k^m],
  \]
  where the last inequality follows by the inductive assumption \eqref{eq:5th}.
  Therefore, we obtain \eqref{eq:2nd}.
  If $F$ is an exact se-merging function,
  we obtain \eqref{eq:3rd} from \eqref{eq:4th} with $m=2$ and $\p=\q$
  since $\E^\p[F(E_1,\dots,E_K)]=\E^\p[\Pi_K(E_1,\dots,E_K)]=1$.
\end{proof}
 
Proposition~\ref{prop:trade-off} has two implications.
First, under the alternative $\q$, if the e-variables are powered, then all moments of the resulting combined e-variable are maximized by the product function. 
Second, under the null,
 the product function $\Pi_K$
results in the largest variance among all exact se-merging functions if the e-variables are exact and independent.

Although having a large mean under $\q$ can be seen as desirable, 
having a large  variance  is generally not a desirable property. Hence the product function can be seen as an extreme choice of ie-merging function.

In particular, if the mean under $\q$, 
 $\E^\q [F(E_1,\dots,E_K)  ]$,  is chosen as an objective to optimize, then 
 the product function $\Pi_K$ dominates all other se-merging functions $F$  under the assumption that $\E^\q[E_k]\ge 1$ for all $k$. 
 Putting this into the context of testing by betting as in Section~\ref{sec:EAEP}, this means the ``all-in'' strategy of choosing $\lambda_t=1$ for all $t$ yields an e-process with the largest expected value at any given time, if each sequential e-value has mean larger than $1$ under $\q$.
 Nevertheless, we remind the reader that e-power is the right objective, instead of the expected value. 
Recall that the e-power is measured by the expected logarithm under $\q$, not any of the moments.
Generally, $\Pi_K(E_1,\dots,E_K)$ does not have the largest e-power. A comparison is provided in Example~\ref{ex:all-in} below. 


 \begin{example}[All-in versus log-optimal] 
 \label{ex:all-in}
We consider the e-processes built on $(E_k)_{k\in \N}$ as in Definition~\ref{def:empirical},   allowing an infinite number of e-variables.
 Suppose that under $\q$, each   $E_k$ for $k\in \N$ is distributed as $\q(E_k=4)=\q(E_k=0)=1/2$ and they are independent. 
 \begin{enumerate}
     \item[(a)]  The all-in strategy is to choose $\lambda_k=1$ for all $k$, leading to
 $M_k^{\text{all-in}}=\Pi_k(E_1,\dots,E_k)$  for any $k\in \N$, which  has two-point distribution with probability $2^{-k}$ taking the value $4^k$ and with probability $1-2^{-k}$ taking the value 0.
     \item[(b)]
On the other hand, we can compute that the $\q$-log-optimal e-process $(M_k^{*})_{k\in \N}$ built on $(E_t)_{t\in [k]}$ in Definition~\ref{def:empirical}  has $\lambda_k=1/3$ for all $k\in \N$, leading to $\E^\q[ \log (1-\lambda_k+\lambda_k E_t)]=\log(4/3)$. 
Then,  by the law of large numbers,
$M_k^*\sim (4/3)^k$ $\q$-almost surely as  $k\to\infty$.
(This asymptotic growth rate is also achieved by the empirically adaptive  e-process by 
Theorem~\ref{th:adaptive}.)  
 \end{enumerate} 
One can see that, although $\E^\q[M_k^{\text{all-in}}] \ge  \E^\q[M_k^{*}]$ for all $k\in\N$ (indeed, $M_k^{\text{all-in}}$ has the largest expectation among all e-processes built on $(E_t)_{t\in [k]}$ by Proposition~\ref{prop:trade-off}), its distribution has a large chance of getting $0$, and moreover $M_k^{\text{all-in}}\to 0$ $\q$-almost surely. These are highly undesirable features of the all-in betting strategy. 
 \end{example}

\subsection*{Summary}

The following points summarize the results in the previous few sections on choosing   merging functions.
\begin{enumerate}
    \item [(i)] For merging arbitrarily dependent e-values, one needs to use the arithmetic mean or a weighted arithmetic mean.  
    \item  [(ii)] For merging sequential e-values, one needs to use martingale merging functions.  In particular, the empirically adaptive martingale merging function in Example~\ref{ex:eMMF} can be used as a default method without prior information on the e-values' behavior under the alternative hypothesis.
        \item [(iii)] For merging independent e-values, one can use, among many choices, the product or the U-statistic functions. Although the product e-merging function has a weak optimality, it also suffers from an undesirable property of having large variance. 
        \item[(iv)] For merging iid e-values,  the empirically adaptive martingale merging function in Example~\ref{ex:eMMF} has asymptotic log-optimality in the sense of Theorem \ref{th:adaptive}.
\end{enumerate}

\section{Family-wise error rate and merging e-values}
\label{sec:FWER}
 
 The merging methods developed in  this chapter can be used to design multiple testing procedures based on e-values. 
In Chapter~\ref{chap:compound}, we will study procedures that control the false discovery rate in great detail. 
In this section, we briefly explain how e-values can be used for testing multiple hypotheses with family-wise error rate (FWER) control, and we omit procedures for other error metrics. 

We give a minimal set of definitions here, and they will be re-elaborated in Chapter~\ref{chap:compound}.
Consider $K$ different null hypotheses, described by sets of probability measures $\mathcal{P}_1,\ldots,\mathcal{P}_K$. For a probability measure $\p\in \cM_1$, the null set of hypotheses   is given by $\cN =\{k\in \cK:\p\in \cP_k\}$ (which depends on $\p$).

Suppose that $E_k$ is an e-variable for $\cP_k$, $k\in \cK$.
For any set $A\subseteq \cK$,
we can define an e-variable for the intersection of $\cP_k$,~$k\in A$, via
$$
E_A=F(E_k:k\in A),
$$
 where $F$ is an e-merging function of $|A|$ e-values (or an ie-merging function if we assume that the e-variables are independent). 
 For instance, taking $F$ as the arithmetic mean $\fM_{|A|}$, we get the e-variable
\begin{align}
\label{eq:c8-mean-e-A}
E_A=\frac{1}{|A|}\sum_{k\in A}E_k.
\end{align}
The collection $(E_A)_{A\subseteq\cK}$ in \eqref{eq:c8-mean-e-A} is called the \emph{mean e-collection}.

With any e-merging function $F$, 
$E_A$ is an e-variable for $\bigcap_{k\in A}\cP_k$.
An e-testing procedure takes $K$ e-values as input and output a set of indices $\mathcal D\subseteq \cK$ (more formal definitions are given in Chapter~\ref{chap:compound}). 
We would like the procedure to control; FWER at level $\alpha\in (0,1)$, meaning the guarantee
$$
\p(\mathcal D\cap \cN  \ne  \varnothing)\le \alpha
$$
 for all $\p\in \cM_1$. 
In other words, the probability of rejecting any null hypothesis is at most $\alpha$, for every possible distribution $\p$.

A simple  procedure to control FWER based on $(E_A)_{A\subseteq\cK}$ is described below. 
Let 
$$E^*_k =\min_{A \subseteq \cK :k\in A} E_A.$$
The e-testing procedure $\cD$ rejects all indices $k\in \cK$ with $E^*_k\ge 1/\alpha$.
To see its FWER control, note that for any $\p\in \cM_1$,
$$
\max_{k\in \cN  } E^*_k
=\max_{k\in \cN  }
\min_{A:k\in A} E_A
\le \max_{k\in \cN  } E_{\cN} =E_{\cN},
$$
and 
$$
\p(\mathcal D\cap \cN  \ne  \varnothing) 
= \p\left(\max_{k\in \cN } E^*_k\ge 1/\alpha
\right)\le 
 \p\left( E_{\cN}\ge 1/\alpha
\right)\le \alpha. 
$$
Note that for the FWER control, we only require that each $E_A$ is an e-variable for $\bigcap_{k\in A}\cP_k$, and it  need not be computed by merging  e-values from  the individual hypotheses. Such $(E_A)_{A\subseteq\cK}$ is called an e-collection; see Section~\ref{sec:9-closed-ebh}.

\begin{algorithm}[t]
  \caption{Adjusting e-values for FWER using the arithmetic mean}
  \label{alg:BH}
  \begin{algorithmic}[1]
    \Require
      A sequence of e-values $e_1,\dots,e_K$ (realizations of $E_1,\dots,E_K$).
    \State Find a permutation $\pi:\cK\to\cK$ such that $e_{\pi(1)}\le\dots\le e_{\pi(K)}$.
    \State Set $e_{(k)}:=e_{\pi(k)}$, $k\in\cK$ (these are the ascending {order statistics}).
    \State $S_0:=0$
    \For{$i=1,\dots,K$}
      \State $S_i := S_{i-1} + e_{(i)}$
    \EndFor
    \For{$k=1,\dots,K$}
      \State $e^*_{\pi(k)}:=e_{\pi(k)}$
      \For{$i=1,\dots,k-1$}
        \State $e := \frac{e_{\pi(k)}+S_i}{i+1}$
        \If{$e < e^*_{\pi(k)}$}
          \State $e^*_{\pi(k)} := e$
        \EndIf
      \EndFor
    \EndFor
  \end{algorithmic}
\end{algorithm}

For the mean e-collection in \eqref{eq:c8-mean-e-A}, a simple $O(K^2)$ algorithm to compute $E_k^*$ from the original e-values $E_1,\dots,E_K$ is described in Algorithm~\ref{alg:BH}.
If we use the product ie-merging function instead, then a simple $O(K)$ algorithm is available,  omitted here.  
Algorithm~\ref{alg:BH} effectively computes
\begin{align*}
  E^*_{\pi(k)} 
  &=
  \min_{i\in[ k-1]}
  \frac{E_{\pi(k)}+E_{(1)}+\dots+E_{(i)}}{i+1}
\end{align*}
for $k\in\cK$,
where $\pi$ is the ordering permutation such that $E_{(j)}=E_{\pi(j)}$ is the $j$th ascending order statistic among $E_1,\dots,E_K$.

\section*{Bibliographical note}

The content in this chapter is mainly based on \cite{Vovk/Wang:2021,vovk2024merging}. 
E-merging functions were introduced by \cite{Vovk/Wang:2021}.
The proof of Theorem~\ref{th:mart-merg} is  in   \cite{vovk2024merging}.
Negative lower orthant dependence was introduced by \cite{BSS82}, and the corresponding results on merging e-values are in \cite{chi2022multiple}.
The FWER controlling procedures based on e-values were studied by \cite{Vovk/Wang:2021} and \cite{hartog2025family}, while recent connections to closed testing were explored in~\cite{xu2025closure} and~\cite{goeman2025epartitioning}.

\chapter{False discovery rate control using compound e-values}
\label{chap:compound}

Our purpose in this chapter is to show how e-values are (a) useful for multiple testing while controlling the false discovery rate (FDR), (b) inherent and central to all FDR controlling procedures. The main setting is explained next, part of which is briefly mentioned in Section~\ref{sec:FWER}.

Suppose we observe data $X$ drawn according to some unknown  distribution $\p^{\dagger}$. 
Here, $X$ denotes the entire dataset available to the researcher.
 We consider testing
$K$ different hypotheses.  The null hypotheses are described by sets of probability measures $\mathcal{P}_1,\ldots,\mathcal{P}_K$. Recall the notation $\cK=[K]$, which avoids a clash with certain decreasing order statistics used in this chapter. 
For $k \in \mathcal K$, the $k$-th null hypothesis $H_k$ is $\p^{\dagger}\in \cP_k$, corresponding to $H_0:\p^{\dagger}\in \cP$ in Chapter~\ref{chap:introduction}. 
Note that it is possible that $\p^{\dagger} \notin \bigcup_{k \in \mathcal{K}} \mathcal{P}_k$, meaning that all null hypotheses could be false.
As in the setting of single hypothesis testing,
the unknown data-generating $\p^{\dagger}$ is used to formally state the multiple hypothesis testing problem and will not appear in the subsequent analysis. 

We write 
 $$\mathcal{N}(\p) = \{k \in \mathcal{K}\,:\,  \p \in \mathcal{P}_k\} \subseteq \mathcal{K},~~~~\p\in \cM_1,$$ 
 which is 
  the index set of  true null hypotheses if $\p$ is the data generating distribution.
Denote by 
$\cN=\cN(\p)$ and  $K_0=|\cN(\p)|$, where we omit the reliance on $\p$.
In the setting of multiple testing, we will let $\cP$ denote a generic set of distributions. Unless otherwise specified, we simply take (explained in Section~\ref{sec:c8-com-e}) \begin{equation*}
\mathcal{P} = \bigcup_{k \in \mathcal{K}} \mathcal{P}_k. 
\end{equation*} 

 Other specifications of $\mathcal P$ may  represent some assumptions or restrictions on the distributions in the testing problem, which will be useful in Section~\ref{sec:c8-pBH} and Chapter~\ref{chap:approximate}.
A central component of this chapter is the e-Benjamini--Hochberg (e-BH) procedure, introduced in Section~\ref{sec:c9-1}.


\section{Compound e-values and  the e-BH procedure}
\label{sec:c9-1}

\subsection{Compound e-variables}

\label{sec:c8-com-e}


A typical setting for our study here is that each hypothesis is associated with an e-variable $E_k$. Sometimes, we also consider and compare with the classic setting, where each hypothesis is associated with a p-value $P_k$. However, more generally, we will only require $(E_1,\dots,E_K)$ to satisfy a more relaxed definition, given below.

\begin{definition}[Compound e-variables]
\label{defi:compound_evalues}

Fix the null hypotheses $(\cP_1,\dots,\cP_K)$ and a set $\cP $ of distributions. 
Let $E_1,\dots,E_K$ be $[0,\infty]$-valued random variables. We say that 
$E_1,\dots,E_K$ are \emph{compound} e-variables for $(\cP_1,\dots,\cP_K)$ under $\cP$ if
$$\sum_{k : \p \in \cP_k} \E^{\p}[E_k] \leq K \qquad \mbox{for all $\p\in \mathcal P$.}$$
They are called \emph{tight} compound e-variables if the supremum of the left hand side over $\p \in \cP$ equals $K$.

\end{definition}

We omit ``under $\cP$'' in case $\cP$ is the set of all distributions.  For all conditions to be checked in the above definition, it suffices to consider  $\p \in   \bigcup_{k \in \mathcal{K}}  \mathcal{P}_k$.
Therefore, it is without loss of generality to assume $\mathcal P \subseteq \bigcup_{k \in \mathcal{K}} \mathcal{P}_k$, and by default (without any restrictions) it is
\begin{equation*}
\mathcal{P} = \bigcup_{k \in \mathcal{K}} \mathcal{P}_k. \end{equation*}

When $K=1$, Definition~\ref{defi:compound_evalues} coincides with the definition of an e-variable. Of course, a vector of e-variables is a trivial but important special case of compound e-variables. As is clear from the definition, no restriction is placed on the dependence structure of the e-variables.

\begin{example}[Weighted e-values]
\label{ex:c8-weighted}
If each $E_k$ is an e-variable for $\cP_k$, then $E_1,\dots,E_K$ are compound e-variables for $(\cP_1,\dots,\cP_K)$. Further, let $w_1,\ldots,w_K \geq 0$ be deterministic nonnegative numbers such that $\sum_{k \in \mathcal{K}} w_k \leq K$ and define $\tilde E_k = {E}_k w_k$. Then $\tilde E_1,\dots,\tilde E_K$ are compound e-variables. 
\end{example}

Note that the class of compound e-values is much richer than the class of weighted e-values in Example~\ref{ex:c8-weighted}, because the later satisfies the stronger constraint 
$$
\sum_{k\in \mathcal K} \sup_{\p\in \cP_k} \E^\p[\tilde E_k] \le K, 
$$
compared to the condition for the compound e-values  
$$
\sup_{\p\in \bigcup_{k\in \mathcal K}\cP_k} \sum_{k\in \mathcal N(\p)} \E^\p[E_k] \le K.
$$

We  now give a more elaborate example of a nonparametric setting where one naturally encounters compound e-variables. More importantly, these e-variables are dependent in a complicated way, making the usual assumptions on the dependence structure, such as independence, inapplicable. 

\begin{example}[Multi-armed bandit testing]  \label{ex:21}
Suppose there are $K$ traders (or machines), and a researcher is interested in knowing which ones are skillful (or useful). This is a classic problem in finance.
For $k\in \cK$,
the null hypothesis $H_k$ is that trader $k$ is not skillful, 
meaning that  they make no profit on average (without loss of generality we can assume the market risk-free return rate is $0$).
The nonnegative random variables $X_{k,1},\dots,X_{k,n}$ are the monthly realized 
performance (i.e., the ratio of payoff to investment; $X_{k,t}>1$ presents a profit and $X_{k,t}<1$ means a loss)
 of agent $k$ from month $1$ to month $n$.
 The no-skill null hypothesis is $\E[X_{k,t} \mid \mathcal F_{t-1}] \le 1$ for $t\in [n]$,
 where the $\sigma$-field
$\mathcal F_t$ represents the available market information up to time $t\in\{0,\dots,n\}$,
and we naturally assume that $(X_{k,t})_t$ is adapted to $(\mathcal F_t)_t$. 
Since the agents are changing investment strategies over time 
and all strategies depend on the   financial market evolution,
there is complicated serial dependence within $(X_{k,1},\dots,X_{k,n})$
for single $k$, as well as cross dependence among agents $k\in \cK$. 
Because of the complicated serial  dependence  and the lack of distributional assumptions of the performance data, it is difficult to obtain useful p-values for these agents.
Nevertheless, we can easily obtain useful e-values: For instance, $E_k=\prod_{t=1}^n X_{k,t}$ is a valid e-value, as well as any mixture of U-statistics of $X_{k,1},\dots,X_{k,n}$, including the mean and the product; indeed, $X_{k,1},\dots,X_{k,n}$ are   sequential e-variables in Definition~\ref{def:se-values}, and they can be combined using merging methods in Chapter~\ref{chap:multiple}.
Moreover, the obtained e-values $E_1,\dots,E_K$ are dependent in a complicated way. Even if these e-values are not very large, they can be useful for other studies on these  traders.   
For this problem, other sophisticated e-values can also be constructed such as the ones from the empirically adaptive e-processes; see Chapter~\ref{chap:eprocess}.

In the above setting, we implicitly assumed that each agent has the same initial wealth $1$. If they have different initial wealth values $c_1,\dots,c_K$ with $\bar c = (c_1+\dots+c_K)/K$, then the wealth of the $k$-th trader is $E'_k = c_k E_k$, and it is easy to check that $E'_1/\bar c,\dots,E'_K/ \bar c$ are compound e-variables, while $E'_1/c_1, \dots, E'_K/c_K$ are e-variables.  
 \end{example}
 
We conclude this section with two general constructions of compound e-variables.

\begin{example}[Convex combinations of compound e-values]
    \label{exam:convex_combi}
    Suppose that $E_1^{(\ell)}, \dots, E_K^{(\ell)}$ are compound e-variables for $\ell=1,\ldots,L$. Let $w_1,\dotsc,w_L \geq 0$ be deterministic nonnegative numbers such that $\sum_{\ell=1}^L w_\ell = 1$. Define $E_k = \sum_{\ell=1}^L w_\ell E_k^{(\ell)}$. Then $E_1,\dots,E_K$ are compound e-variables. 
\end{example}

\begin{example}[Combining evidence from follow-up studies]
Suppose we summarize data $X$ into  compound e-variables $E_1,\dots,E_K$. Based on these, suppose that we select a subset of indices $\mathcal{S} \subseteq \mathcal{K}$ of interest, for which we seek additional evidence. With this goal in mind, suppose we collect additional data $Y$ and summarize it in the form of e-variables $\{E'_j\}_{j \in \mathcal{S}}$ that are conditionally valid given the past data; in particular, $\mathbb E^{\p}[E'_j | E_j] \leq 1$ for all $j \in \mathcal{S}$ with $\p \in \mathcal{P}_j$. Then defining $E'_j = 1$ for $j \notin \mathcal{S}$, we have that $E_1E'_1, \dots, E_KE'_K$ are also compound e-variables.
\end{example}

\subsection{FDR and the e-BH procedure}
\label{sec:c8-eBH}
We first define the concept of FDR, central to modern multiple hypotheses testing. 
\begin{definition}[Multiple testing procedures]
Let $X$ represent the available data used for testing. 
\begin{enumerate}[label=(\roman*)]
\item 
A \emph{multiple testing procedure}  $\cD$ is a Borel function of $X$  that produces a subset of $\cK$
representing the indices of rejected hypotheses. A testing procedure $\cD$ that takes e-values or compound e-values as input is called an \emph{e-testing procedure}. 
\item   The hypotheses that are rejected by $\cD$, given by $\cD(X)\in 2^{\mathcal K}$, are called \emph{discoveries}. 
\end{enumerate}
\end{definition}
We overload notation and write both the procedure (the mapping from $X$ to $\cK$) and the realized set of discoveries (a random subset of $\cK$) as $\cD$, because it is typically clear which one is meant from context.  
We use the convention
$$\frac{0}{0}=0$$
 in the calculation of error metrics in multiple hypothesis testing.

\begin{definition}[FDP and FDR]
\label{def:FDR}
In the following objects, $\cD$ should be interpreted as the random subset $\cD(X)$ of $\cK$. 
\begin{enumerate}[label=(\roman*)]

\item  A \emph{false discovery} is 
a true null hypothesis rejected by $\cD$.
Let 
$$
F_{\cD} = \sum_{k: \p \in \mathcal{P}_k} \id_{\{k\in \cD\}} = |\mathcal N(\p) \cap \cD|
\mbox{~~~and~~~}R_{\cD}=|\cD|,
$$ 
which are the number of  false discoveries and the total number of discoveries, respectively. 

\item The 
\emph{false discovery proportion (FDP)} is ${F_{\cD}}/{R_{\cD} }$.

\item 
The \emph{false discovery rate (FDR)} under $\p$ is the expected value of the FDP,
that is,  $$ 
\mathrm{FDR}_{\cD}=\E^{\p} \left[\frac{F_{\cD}}{R_{\cD} }\right] = \E^{\p} \left[\frac{F_{\cD}}{R_{\cD} \vee 1}\right].$$  
\item 
A multiple testing procedure $\cD$ has \emph{FDR control} (at most) $\alpha\in (0,1)$  under $\cP$  if $\mathrm{FDR}_{\cD}\le \alpha$  
for all $\p \in \cP$. 
We omit ``under $\cP$'' when $\cP=\cM_1$ 
or $\cP=\bigcup_{k \in \mathcal{K}} \mathcal{P}_k$.
\end{enumerate} 
\end{definition}

Clearly, both $\mathrm{FDR}_{\cD}$
and $F_{\cD}$ depend on the unknown $\p$, and we omit this reliance for brevity. 
 Similarly to the case of compound e-values, FDR is trivially $0$ if $\p$ is outside $\bigcup_{k \in \mathcal{K}} \mathcal{P}_k$, so we do not need to check such $\p$.

The FDR has been  the most popular error metric in multiple testing, with wide applications in many scientific disciplines. 
There are many FDR-controlling procedures with p-values as the input.
The Benjamini--Hochberg (BH) procedure of \cite{benjamini1995controlling} is arguably the most popular and successful in this context (which is, to our knowledge, the most cited paper in statistics at the time of writing of this book).
We will study a corresponding e-testing procedure called the e-BH procedure, which as an intimate connection to the BH procedure that will be explained later. 

\begin{definition}[The e-BH procedure]
\label{defi:eBH}
Suppose that $e_1,\dots, e_K$ are realized compound e-values for the hypotheses $(\cP_1,\dots,\cP_K)$. 
For $k\in \mathcal K$, let $e_{[k]}$ be the $k$-th order statistic of $e_1,\ldots,e_K$, from the largest to the smallest.  The \emph{e-BH procedure at level $\alpha \in (0,1)$} rejects all hypotheses with the largest $k^*$ e-values, where 
\begin{equation*} 
\label{eq:e-k-intro} 
k^*:=\max\left\{k\in \mathcal K: \frac{k e_{[k]}}{K} \ge \frac{1}{\alpha}\right\},
\end{equation*}   
with the convention $\max(\varnothing) = 0$.
\end{definition} 

\begin{figure}[htbp]
   \begin{center}      \includegraphics[width=0.7\textwidth, trim={0 0.5cm 0 1.8cm}, clip]{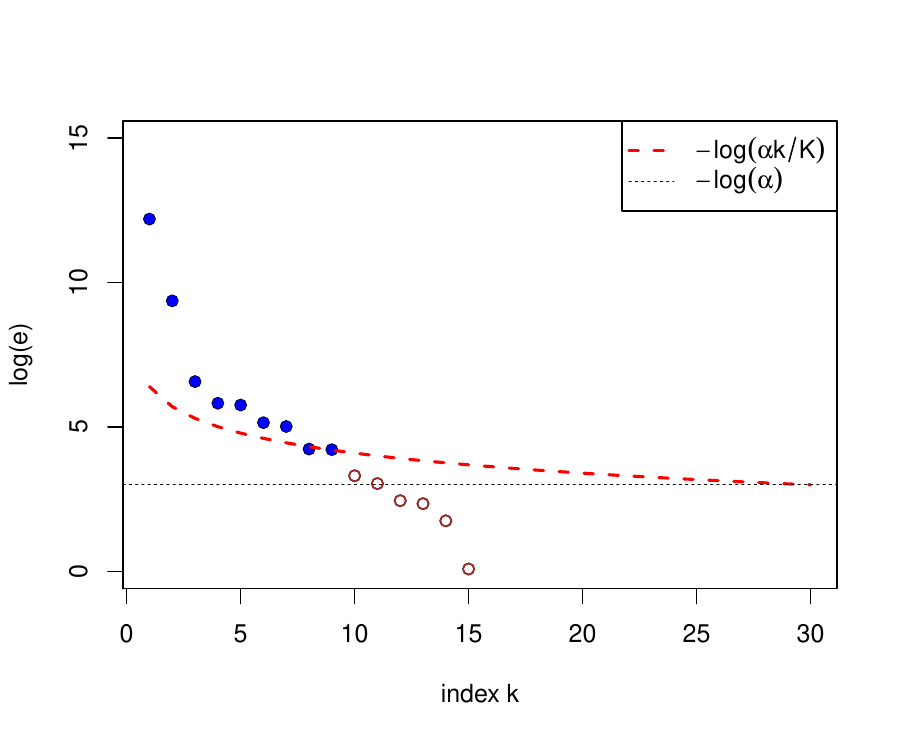}\\ \includegraphics[width=0.7\textwidth, trim={0 0.5cm 0 1.8cm}, clip]{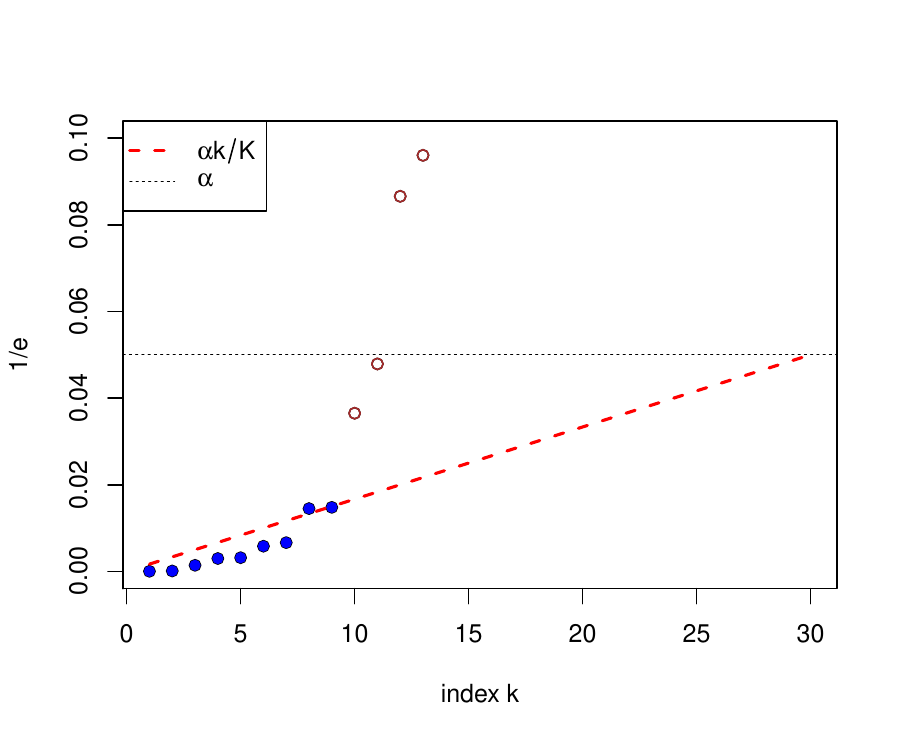}
   \end{center} 
   \caption{An illustration of the e-BH procedure in Definition~\ref{defi:eBH} applied to a simulated sample with $K=30$ and $\alpha=0.05$. 
   We compare $e_{[k]}$ with $K/(\alpha k)$ for $k\in \cK$
   on the log scale in the top panel and on the reciprocal scale in the bottom panel (some small e-values are outside the plotted relevant range). 
   Hypotheses associated with e-values no less than the last point above the red dashed curve in the top panel are rejected, including one that lies below the curve. 
   The bottom panel also illustrates the BH procedure in Definition~\ref{defi:BH}, with the reciprocals of e-values as p-values.
   }
   \label{fig:c9-1}
   \end{figure}

A stylized illustration of the e-BH procedure is depicted in Figure~\ref{fig:c9-1}.

Before showing the FDR guarantee  of the e-BH procedure, 
we slightly enlarge the scope of our methodology. 

\begin{definition}[Self-consistent e-testing procedures]
\label{def:self-cons-e}
An e-testing procedure $\cD$ is  \emph{self-consistent at level $\alpha\in(0,1)$} if  every rejected hypothesis $H_k$ has a realized compound e-value $e_k$ that satisfies
$$
e_k \ge \frac{K}{\alpha R_{\cD}  }.
$$ 
This class includes the e-BH procedure as a special case.
\end{definition} 

 The  e-BH procedure dominates all other self-consistent e-testing procedures by definition, meaning that it maximizes the number of discoveries within this class.
Nevertheless, self-consistent testing procedures that reject a smaller set than that of the e-BH procedure may be useful if we want the set of discoveries to also satisfy some additional structural or logical constraint. For example, if the hypotheses are structured as a graph, then we may to reject a set of hypotheses that form connected subgraphs.

The main result in this section is the FDR control of self-consistent e-testing procedures for any compound e-variables. 

\begin{theorem}[FDR control and post-hoc FDR control]\label{th:compliant}
Any self-consistent e-testing procedure  at level $\alpha \in (0,1)$, including the e-BH procedure, has FDR at most $\alpha$. This claim holds for any $\p \in \cM_1$ and for arbitrary compound e-variables, regardless of the dependence structure that $\p$ induces amongst them. 

In fact, for any compound e-variables, and for any class of 
self-consistent e-testing procedures $(\cD_\alpha)_{\alpha \in (0,1)}$ indexed by their level $\alpha$, the post-hoc FDR guarantee holds for any $\p \in \cM_1$:
\[
\E^{\p} \left[ \sup_{\alpha \in (0,1)} \frac{F_{\cD_\alpha}}{\alpha R_{\cD_\alpha}} \right] \leq 1.
\]
Further, if the compound e-variables are in fact e-variables, then the upper bound above can be improved to $K_0/K$.
\end{theorem}
The last claim above allows for a formal guarantee even when $\alpha$ is itself a function of the data. 
Omitting the supremum in the final statement yields the preceding FDR claim.
\begin{proof}[Proof.]
Let $\mathbf E=(E_1,\dots,E_K)$ be an arbitrary vector of compound e-variables fed to the testing procedure $\cD_\alpha$ and $\cD_{\alpha} (\mathbf E)$ be the set of rejected hypotheses. 
By direct substitution, we have
\begin{align} 
 \frac{F_{\cD_{\alpha}}}{\alpha R_{\cD_\alpha} }   =
 \frac{ |\cD_{\alpha} (\mathbf E)\cap \mathcal N |  }{\alpha (R_{\cD_{\alpha}} \vee 1)}  
 & =\sum_{k\in \mathcal N}  \frac{ \id_{\{k\in \cD_{\alpha}(\mathbf E)\}}  }{\alpha (R_{\cD_{\alpha}} \vee 1)}  
 \\& \le \sum_{k\in \mathcal N}  \frac{ \id_{\{k\in \cD_{\alpha}(\mathbf E)\}} E_k   }{K}   \le \sum_{k\in \mathcal N}   \frac{ E_k  }{K} , 
\label{eq:compliant}
\end{align} 
where the first inequality  is due to self-consistency. 
By the definition of compound e-variables, we immediately obtain for arbitrary $\p \in \cM_1$ that
\begin{align}
\E^{\p}\left[\sup_{\alpha \in (0,1)} \frac{F_{\cD_{\alpha}}}{\alpha R_{\cD_{\alpha}} }  \right]  \le   1.
\end{align}
If $E_1,\dots,E_K$ are e-variables, then we have
\begin{align}
    \label{eq:c8-fdrproof} 
\E^{\p}\left[\sup_{\alpha \in (0,1)} \frac{F_{\cD_{\alpha}}}{\alpha R_{\cD_{\alpha}} }  \right]  \le   \sum_{k\in \mathcal N}  \E^{\p}\left[ \frac{   E_k  }{K}   \right] \le \frac{ K_0}{K},
\end{align}
as claimed. 
\end{proof}

Most notably, the e-BH procedure controls FDR regardless of how the e-variables are dependent. This key feature distinguishes it from the usual BH procedure with input p-values, for which the FDR guarantee requires certain dependence assumptions, a topic we return to later in Section~\ref{sec:c8-pBH}.

In the next proposition we provide an alternative description of the e-BH procedure, which illustrates that the e-BH procedure rejects e-values according to a threshold $t_\alpha$ that depends on all other input e-values. 
This alternative description is useful in Section~\ref{sec:c8-1-2}.
 \begin{proposition}
\label{prop:c8-1-1}
Let $e_1,\dots,e_k\in \R_+$ and $\alpha \in (0,1)$. Define
  \begin{align}
  \label{eq:r1-talpha}
  R(t)= |\{k\in \mathcal K: e_k\ge t\}| \vee 1;~~~
 t_\alpha = \inf \{t\in [0,\infty):t{R(t)} \ge K/\alpha\}.
\end{align}
For each $k$, the e-BH procedure at level $\alpha$ applied to the realizations $e_1,\dots,e_K$ and rejects 
 $H_k$ if and only if $e_k\ge t_\alpha$.  Moreover, $t_\alpha R(t_\alpha)= K/\alpha$, and $t_\alpha$ takes values in $\{K/(k \alpha): k \in \mathcal K\}$.
 \end{proposition}

 \begin{proof}[Proof.]
 For $t\in [0,\infty)$, 
 let $f(t)$ be the number of true null hypotheses with an e-value $e_k$ larger than or equal to $t$.
  Define the quantity
 $
 g(t) =  t{R(t)}/K,
$
 and by definition
 $
 t_\alpha = \inf \{t\in [0,\infty):g(t) \ge 1/\alpha\}.
 $
 Clearly $t_\alpha \in [1/\alpha,K/\alpha]$ since $g(t)\le t$ and $R(t)\ge 1$.
Since $g$ only has downside jumps and $g(0)=0$, we know $g(t_\alpha)=1/\alpha$, and thus
 $
t_\alpha R(t_\alpha)= K/\alpha.
 $

If $e_{k}\ge t_\alpha$, then $H_k$ is rejected by the definition of e-BH.
 If $e_{[k]} < t_\alpha$, then by definition of $g$, we have 
 $$\frac{k e_{[k]}}{K}  \le \frac{R(e_{[k]}) e_{[k]}}{K} =g(e_{[k]}) < 1/\alpha.$$
 Thus, each $H_k$ is rejected by the e-BH procedure if and only if $e_k\ge t_\alpha$.
 \end{proof}


\subsection{Boosting the e-values with distributional information}
\label{sec:c8-1-2}

If we know some information of the null distribution of the e-variables $E_1,\dots,E_K$ fed to the e-BH procedure, then we can enhance the power of the e-BH procedure by a mechanism called \emph{boosting e-values}. 

Let $$K/\mathcal K :=\{K/k:k\in \mathcal K\},$$ and  define a truncation function $T:[0,\infty]\to [0,K]$ by
letting $T(x)$ be the largest number in $K/\mathcal K\cup\{0\}$ that is no larger than $x$. In other words, 
\begin{equation} \label{eq:def-T}
T(x)  = \frac{K}{\lceil  K/x\rceil}\id_{\{x\ge 1\}}\mbox{~~with $T(\infty)=K$}.
\end{equation} 
Note that $T$ truncates $x$ to take only values in $K/\mathcal K\cup\{0\}$.  
For each $k\in \mathcal K$, 
take a \emph{boosting factor} $b_k\ge 1$   such that 
 \begin{align}
 \label{eq:enhance-2}  \E^\p[T(\alpha b_k E_k)]\le \alpha~~ 
 \mbox{for $\p \in \cP_k$,}
 \end{align} 
 and let $E_k'=b_k E_k$. We call $E_1',\dots,E_K'$ the boosted e-values. 
Note that choosing $b_k=1$ always satisfies \eqref{eq:enhance-2} because $T(x)\le x$ for all $x\in \R_+$. 
One would try to choose the largest $b_k$ possible subject to \eqref{eq:enhance-2}. 
Computing the precise value of $b_k$   depends on our knowledge of the distribution of $E_k$ under $\p$. In some situations, one may have partial but not full distributional information of $E_k$ under $\p$, leading to a smaller $b_k$ than the one that attains equality in \eqref{eq:enhance-2}.
If no distributional information is available, i.e., one only knows that $E_k$ is an e-variable for $\cP_k$ but nothing more, then one has to set $b_k=1$.

\begin{theorem}\label{th:c8-boost}
  Applied to arbitrary nonnegative random variables $ E'_1,\dots,E'_K $, the  e-BH procedure $\cD$ at level $\alpha\in (0,1)$ satisfies  
  $$
\E^{\p}\left[\frac{F_{\cD}}{R_{\cD} }\right] \le \frac{  1}{K}  \sum_{k\in \mathcal N}\E^{\p}\left[   T(\alpha E'_k)\right],
$$
for any $\p\in\cM_1$.
In particular, if $E'_1,\dots,E'_K$ are  e-variables for $\cP_1,\dots,\cP_K$ or the boosted e-values via \eqref{eq:enhance-2}, then the FDR is controlled at $K_0\alpha/K$.
\end{theorem}
The proof of Theorem~\ref{th:c8-boost} can be described with a very simple intuition: If e-values $e_1,\dots,e_K$ are replaced by 
$T( \alpha e_1)/\alpha,\dots,T(\alpha e_K)/\alpha$, then they will lead to the same set of rejected hypotheses,
because by Proposition~\ref{prop:c8-1-1},  e-values are rejected only at thresholds in $(\alpha^{-1} K ) /\mathcal K$. 
Therefore, for the desired FDR control, it suffices to verify that $T(\alpha E_k) /\alpha$ are e-values. This leads to condition \eqref{eq:enhance-2}, justifying the use of $b_kE_k$.

We can also boost compound e-values by using $b_1,\dots,b_K$ satisfying 
 \begin{align}
 \label{eq:enhance-2-ce} \sum_{k:\p\in \cP_k}  \E^\p[T(\alpha b_k E_k)]\le K \alpha~~ 
 \mbox{for $\p \in \cP$,}
 \end{align} 
 following the same argument. 

\subsection{Connection between the e-BH and the BH procedures}
\label{sec:c8-pBH}


There are two main connections we wish to emphasize in this subsection. First,  the BH procedure and its guarantees under PRDS can be recovered as an instance of the e-BH procedure via a certain \emph{boosting} operation, and second,  the BH procedure with the Benjamini--Yekutieli correction (which we call the BHY procedure for brevity)  can be recovered by applying the e-BH procedure to some e-values. We do not expand on the first point here, but do prove the second in Proposition \ref{prop:c8-recover-BY}. 
In sum,  BH is a special case of a \emph{boosted} e-BH, and BHY is a special case of e-BH.

\begin{definition}[BH and BHY procedures]
\label{defi:BH} 

Suppose that $p_1,\dots, p_K$ are realized p-values for the hypotheses $(\cP_1,\dots,\cP_K)$. 
For $k\in \mathcal K$, let $p_{(k)}$ be the $k$-th order statistics of $p_1,\dots,p_K$, from the smallest to the largest. 
The \emph{BH procedure} at level $\alpha\in (0,1)$ rejects all hypotheses with the smallest $k^*$ p-values,
where 
\begin{equation}
\label{eq:p-k}
k^*=\max\left\{k\in \mathcal K: \frac{K p_{(k)}}{k} \le \alpha\right\},
\end{equation}
with the convention $\max(\varnothing)=0$.
The \emph{BHY procedure}  is
to apply the BH procedure to $\ell_K p_1,\dots,\ell_K p_K$, where $\ell_K=\sum_{k=1}^K k^{-1}$. 
\end{definition}

The following are classic results:
\begin{enumerate}
    \item[(i)] The BHY procedure at level $\alpha$ controls FDR at level $(K_0/K)\alpha$.
    \item[(ii)] The BH procedure at level $\alpha$ controls FDR at level $(K_0/K)\alpha$ if the p-variables are independent or satisfy the PRDS condition defined below (that is, under the set $\mathcal P$ of distributions  satisfying the specified conditions); see Remark~\ref{rem:proof-bh-evalue}. A generalization of this is Theorem~\ref{th:ep-BH}.

\end{enumerate}

\begin{definition}[Positive regression dependence on a subset]
\label{definition:prds}

A set $A\subseteq \R^K$ is said to be \emph{increasing}
if $\mathbf x\in A$ implies $\mathbf y\in A$ for all $\mathbf y\ge \mathbf x$.

A vector $(P_1,\dots,P_K)$ of p-variables satisfies positive regression dependence on a subset (PRDS) if for any null index $k\in \mathcal N$ and increasing set $A  \subseteq  \R^K$, the
function $x\mapsto \p\{(P_{\ell})_{\ell \in \mathcal{K}} \in A\mid P_k\le x\}$ is increasing on $ [0,1]$.  
\end{definition}
A caveat of PRDS in  Definition~\ref{definition:prds} is that it enforces certain positive dependence between the nulls and non-nulls, and it is often difficult to verify in many applications, such as genome-wise association studies and finance (Example~\ref{ex:21}).

We can see that the price to pay to get validity under arbitrary dependence for the BH procedure is a factor of $\ell_K$, which is close to $\log K$. 
Since p-values are less useful when they are large, this  correction makes the BH procedure to reject less, and typically much less, hypotheses.

With e-variables $E_1,\dots,E_K$, it is immediate that the e-BH procedure is precisely the BH procedure applied to $1/E_1,\dots,1/E_K$.
Since $P_k:= 1/E_k$ is a p-variable for $\cP_k$ (Proposition~\ref{prop:e-to-p}), 
if the BH procedure controls FDR with p-values $P_1,\dots,P_K$, then the e-BH procedure controls FDR with e-values $E_1,\dots,E_K$. 
However, the e-BH procedure guarantees something more: 
the FDR control is valid under arbitrary dependence among the e-values, whereas the BH procedure needs the $\ell_K$ correction under such a setting.

\begin{proposition}
\label{prop:c8-recover-BY}
\begin{enumerate}
    \item [(i)] With $T$ in \eqref{eq:def-T} and $\alpha\in (0,1)$, the function $f:p\mapsto T(\alpha /(\ell_K p))/\alpha$ is a calibrator. 
    \item [(ii)] Let $P_1,\dots,P_K$ be p-variables for $(\cP_1,\dots,\cP_K)$. 
The BHY procedure applied to $P_1,\dots,P_k$
is identical to 
the e-BH procedure to the e-variables $f(P_1),\dots,f(P_K)$.
\end{enumerate}
\end{proposition}
\begin{proof}[Proof.]
(i) Clearly $f$ is nonnegative and decreasing, and it takes value $0$ on $(1,\infty)$. Let $U$ be a uniformly distributed random variable on $[0,1]$ (under some $\p$). We can compute
\begin{align*} \E^\p[f(U)]&=   \E^\p\left[ \frac{1}{\alpha} T\left(\frac{\alpha}{\ell_K U}\right)\right]\\
&= \frac{1}{\alpha}   \sum_{k=1}^K \frac{K}{k } \p\left(  \frac{\alpha (k-1)}{K\ell_K} \le U < \frac{\alpha k}{K\ell_K}\right)    =\frac{1}{\alpha } \sum_{k=1}^K\frac \alpha {\ell_K k} =1.
\end{align*}
This shows that $f$ is a calibrator. 

(ii) As we explained in Section~\ref{sec:c8-1-2}, applying the e-BH procedure to $$
\frac 1\alpha T\left(\frac{\alpha }{\ell_K P_1}\right),\dots, \frac 1\alpha T\left(\frac{\alpha }{\ell_K P_K}\right)$$
is the same as applying the e-BH procedure to $ 1/(\ell_K P_1),\dots, 1/(\ell_K P_K)$, which is equivalent to applying the BH procedure to $\ell_K P_1,\dots, \ell_K P_K$, the latter being the BHY procedure. 
\end{proof}

As a direct consequence of Proposition~\ref{prop:c8-recover-BY}, the FDR guarantee  of the BHY procedure  directly follows from Theorem~\ref{th:compliant}. 


\begin{remark}\label{rem:proof-bh-evalue}
There is a short proof of FDR control of the BH procedure under PRDS that utilizes e-values. 
First, note that $\widehat\alpha := \alpha k^*/K$ is a decreasing function of all $P_1,\dots,P_K$, where $k^*$ is the number of discoveries of the BH procedure.  
Using this,
the PRDS condition implies that $T=\widehat \alpha$
and $P=P_k$  satisfy \eqref{eq:post-hoc-p-value-4} in Remark~\ref{rem:p-to-e-random-alpha} for each null p-value $P_k$, and therefore (as explained in Remark~\ref{rem:p-to-e-random-alpha}) 
\[
\frac{\id_{\{P_k \leq \widehat \alpha\}}}{\widehat \alpha }
\]
is an e-variable. Then, for any $\p \in \cM_1$, we have by the definition of BH and FDR that
\begin{align*}
    \FDR = \E^\p \left[ \frac{\sum_{k \in \cN} \id_{\{P_k \leq \widehat \alpha\}} }{k^* } \right] =\frac{\alpha}{K} \sum_{k \in \cN} \E^\p \left[ \frac{ \id_{\{P_k \leq \widehat \alpha\}} }{\widehat\alpha } \right] \leq \frac{\alpha |\cN|}{K},
\end{align*}
as required.
\end{remark}

\section{Universality of compound e-values and e-BH}

In this section, we present the (perhaps surprising) fact that all FDR controlling procedures can be written as the e-BH procedure applied to some compound e-values, through two results that are simple to prove. Then, we discuss their implications in combining and derandomizing multiple testing procedures. 


The first result enables the construction of   compound e-values from any   FDR controlling procedure.  

\begin{theorem}[From FDR control to compound e-values]
\label{prop:fdr_control_to_compound_evalues}
Let $\cD$ be any procedure that controls the FDR at a given level $\alpha \in (0,1)$. Let $V_k \in \{0,1\}$ be an indicator of whether hypothesis $H_k$ is rejected by $\cD$, so that $F_{\cD} = \sum_{k \in \mathcal N} V_k$ and $R_{\cD} = \sum_{k \in \mathcal K} V_k$. Define
\begin{align}\label{eq:universal-e} E_k = \frac{K}{\alpha}\frac{V_k}{R_{\cD} \lor 1}.
\end{align}
Then, $E_1,\dots, E_K$ are compound e-variables for $(\cP_1,\dots,\cP_K)$. 
\end{theorem}
\begin{proof}[Proof.]
The result can be seen from, for every $\p\in \cP$,
$$ \sum_{k: \p\in \cP_k} \E^{\p}[ E_k] = \sum_{k: \p\in \cP_k} \E^{\p }\left[ \frac{K}{\alpha}\frac{V_k}{R_{\cD} \lor 1}\right] = \frac{K}{\alpha}\E^{\p }\left[\frac{F_{\cD}}{R_{\cD} \lor 1}\right] \leq \frac{K}{\alpha}\alpha = K,$$
where the inequality follows from the definition of FDR guarantee. 
\end{proof}

In \eqref{eq:universal-e}, the random variables $E_1,\dots,E_K$ are well-defined, noting that $V_k$  and $R_\mathcal D$ do  not depend on $\p$, although $F_{\mathcal D}$ depends on $\p$. 


The second result shows that every FDR procedure is e-BH.
In what follows, for a given set $\cP$ of hypotheses, an FDR procedure $\mathcal D$ at level $\alpha $ is admissible if 
for any FDR procedure $\mathcal D'$ at level $\alpha $ such that $\mathcal D \subseteq \mathcal D'$, we have $
\mathcal D=\mathcal D'$.

\begin{theorem}[Universality of e-BH with compound e-values]
\label{theo:universality_eBH}
Let $\cD$ be any procedure that controls the FDR at level $\alpha$. Then, there exists a choice of compound e-values such that the e-BH procedure yields identical discoveries as $\cD$. 
Further, if $\cD$ is admissible, then these compound e-values can be chosen to be tight. 
\end{theorem} 
\begin{proof}[Proof.]
    The first part directly follows from Theorem~\ref{prop:fdr_control_to_compound_evalues}. 
  To show the last statement, consider an FDR procedure $\mathcal D$ at level $\alpha$ and take the compound e-values $E_1,\dots,E_K$ from  \eqref{eq:universal-e}. Define 
    $$
K^*:= \sup_{\p \in   \cP} \sum_{k : \p \in \cP_k} \E^{\p}[E_k].
$$
If $K^*$ is equal to $K$, then there is nothing to show as $E_1,\dots,E_K$ are tight compound e-values. Otherwise, $K^*<K$. The case $K^*=0$ means one never rejects any hypotheses, for which choosing $E_1'=\dots=E_K'=1$ would suffice as tight compound e-values that produce $\mathcal D$ via e-BH.
For $K^*>0$, we let 
    $$
    E_k'=\frac{K}{K^*}E_k \mbox{~~~for $k\in \mathcal K$},
    $$
    and apply the e-BH procedure to $E'_1,\dots,E'_K$ (which are tight compound e-values). 
    This new procedure controls FDR at level $\alpha$ and  produces at least as many discoveries as $\mathcal D$.
    Since $\mathcal D$ is admissible,
    this procedure must coincide with $\mathcal D$, and hence $\mathcal D$ is  e-BH applied to the tight compound e-values  $E'_1,\dots,E'_K$.
\end{proof}

It is worth noting that the above result holds for every fixed $(\cP_1,\dots,\cP_K)$, meaning that
the procedure $\mathcal D$ that we consider could be tuned to this particular set of hypotheses and does not need to control FDR for any other settings. If this tuned procedure $\mathcal D$ is admissible, it must still be recoverable via e-BH with tight compound e-values. 

For the classic BH procedure in Definition~\ref{defi:BH}, 
it controls FDR  if we assume  that the probability measures in $\cP_1,\dots,\cP_K$ guarantee that the p-values are independent or PRDS. 
It can be identified with the e-BH procedure with e-variables given in 
\eqref{eq:universal-e}, but these e-variables are not necessarily independent or PRDS.

The connections between FDR controlling procedures, e-BH, and compound e-values laid out above motivate the following general and practical mechanism for combining discoveries across multiple testing procedures. To be concrete, in order to test the hypotheses $cP_1,\dots,\cP_K$, suppose we run $L$ different multiple testing procedures $\cD_{\mathrm{GeBH}},\dots,\cD_L$ (this can be easily generalized to the case where each procedure tests a different subset of hypotheses). Assume that the $\ell$-th procedure 
controls the FDR at level $\alpha_{\ell}$.  
Then we may proceed as follows:
\begin{enumerate}
\item For the $\ell$-th multiple testing procedure, form $E_1^{(\ell)},\dots E_K^{(\ell)}$,  which are compound e-variables for $\cP_1,\dots,\cP_K$ (which always exist by Theorem~\ref{prop:fdr_control_to_compound_evalues}). 
They could be   
--- but need not necessarily be --- the implied ones formed in Theorem~\ref{prop:fdr_control_to_compound_evalues}. 
\item Fix weights $w_1,\dots,w_L \geq 0$ with $\sum_{\ell=1}^L w_{\ell} = 1$. Then construct $E_1,\dots, E_K$ by convex combination as in Example~\ref{exam:convex_combi}, i.e.,  $E_{k} = \sum_{\ell=1}^L w_{\ell} E_k^{(\ell)}$ for all $k \in \mathcal{K}$.
These will be compound e-values for $\cP_1,\dots,\cP_K$.
We also allow $w_1,\dots,w_L\ge 0 $ to be random, as long as they are independent of all e-values used in the procedures, and in that case it suffices to require $\sum_{\ell=1}^L \E^{\p} [w_{\ell}] \le  1$ for any $\p \in \cP$ (this condition is similar to the condition for compound e-values). 
\item Apply the e-BH procedure  at level $\alpha$ to the new compound e-values $E_1,\dots, E_K$.
\end{enumerate}
The above construction is guaranteed to control the FDR at level $\alpha$.
Note that the value of $\alpha_{\ell}$ does not matter, as it is only used in the construction of $E_1^{(\ell)},\dots,E_K^{(\ell)}$ , possibly implicitly (e.g., in Theorem~\ref{prop:fdr_control_to_compound_evalues}).   

One important application of the above recipe is derandomization. Suppose that $\cD_{\ell}$ is a randomized multiple testing procedure, that is, it is a function of both the data $X$ as well as a random variable $U_{\ell}$ generated during the analysis. Such randomness may not be desirable, since different random number generation seeds will lead to different sets of discoveries. In such cases, the above recipe can be used to construct a new derandomized procedure $\cD$ that is less sensitive to $U_1,\dots,U_L$ (formally, full derandomization would occur if $U_{\ell}$ are iid and $L\to \infty$).

\section{Log-optimal simple separable compound e-values}


In this section we focus on sequence models, wherein the data $X$ may be written as $X=(X_k : k \in \mathcal{K})$ and in principle we may test each hypothesis $H_k$ using only $X_k$. 
We record the following definition, which we will call upon throughout this section.
\begin{definition}[Simple separable e-variables]
\label{def:simple_separable}
$E_1,\dots,E_K$ are called separable if for all $k$, $E_k$ is $X_k$-measurable, that is $E_k = E_k(X_k)$. They are called simple separable if $E_k = f(X_k)$ for some function $f$ (which is the same for all $k$).
\end{definition}
In what follows, we provide a brief summary of compound decision theory, which motivates the nomenclature ``compound e-values'' (Section~\ref{subsec:compound}) and then we provide a construction  of compound e-values motivated by compound decision theory in Section~\ref{subsec:best_simple_separable}  
(see also Section~\ref{subsec:t_tests} for another related construction).

\subsection{Background: Robbins' compound decision theory}
\label{subsec:compound}

Let $\p_{\bmu}$ be the probability measure governing a Gaussian sequence model:
\begin{equation}
\label{eq:gaussian_sequence}
{\bmu} = (\mu_1,\dotsc,\mu_K) \text{ fixed},\quad X_k \sim \mathrm{N}(\mu_k, 1) \text{ for } k \in \mathcal{K}.
\end{equation}
Suppose we are interested in constructing estimators $\hat{\mu}_k$ of $\mu_k$ such that the following expected compound loss would be small:
 \begin{equation}
\label{eq:compound_risk}
\frac{1}{K}\sum_{k \in \mathcal{K}} \E^{\bmu}[(\hat{\mu}_k - \mu_k)^2],
\end{equation}
where $\E^{\bmu}$ means $\E^{\p_{\bmu}}$. 
If $\mu_k$ has a prior distribution $B$, then, by denoting $\p_B$ the joint distribution of $\bmu$ and $\mathbf X=(X_1,\dots,X_K)$,  the best estimator is the Bayes estimator $\hat{\mu}_k^{B} =   \E^{\p_B}{[\mu_k | \mathbf X]}  = \E^{\p_B}{[\mu_k |  X_k]}$, which is simple separable in the sense of Definition~\ref{def:simple_separable}.  

Now suppose $\bmu$ is deterministic as in~\eqref{eq:gaussian_sequence}. 
Given knowledge of $\bmu$,
which function $s_{\bmu}: \R \to \R$ leads to a simple separable estimator $\hat{\mu}_k^{s_{\bmu}} = s_{\bmu}(X_k)$ that minimizes the risk in~\eqref{eq:compound_risk}? The answer lies in the fundamental theorem of compound decisions, which formally connects~\eqref{eq:gaussian_sequence} to the univariate Bayesian problem with prior $M = \sum_{k \in \cK} \delta_{\mu_k}/K$ equal to the empirical distribution of $\bmu$ (with $\delta_{\mu_k}$ denoting the Dirac measure at $\mu_k$):
\begin{equation}
\label{eq:bayes_compound}
\mu' \sim M ,\;\; X' \mid \mu' \sim \mathrm{N}(\mu', 1). 
\end{equation}

The fundamental theorem of compound decision theory states the following. Given any fixed  $s: \R \to \R$, we have that
\begin{align}
\frac{1}{K}\sum_{k \in \mathcal{K}} \E^{\bmu}[(s(X_k) - \mu_k)^2]  &= \sum_{k \in \mathcal{K}}\frac{1}{K}\E^{\mu_k}[(s(X_k) - \mu_k)^2] \\
&= \E^{\p_{M}}[(s(X') - \mu')^2]. \label{eq:fundamental_compound_decisions}
\end{align}
In the left-hand side display, $\bmu \in \R^{K}$ is treated as fixed (as in~\eqref{eq:gaussian_sequence}), while in the right-hand side display, we only have a single random $\mu' \in \R$ that is randomly drawn from $M$ as in~\eqref{eq:bayes_compound}. From~\eqref{eq:fundamental_compound_decisions} it immediately follows that the optimal simple separable estimator is the Bayes estimator under $M$, that is $s(x) = \E^{\p_{M}}[\mu' \mid X'=x]$.

This construction motivates the construction of feasible non-separable estimators $\hat{\mu}_k = \hat{\mu}_k(X)$ that have risk close to the optimal simple separable estimator. Thus the optimal simple separable estimator defines an oracle benchmark. Just as in classical decision theory, one often restricts attention to subclasses of estimators, e.g., equivariant or unbiased, one sets a benchmark defined by simple separable estimators. What is slightly unconventional here is a fundamental asymmetry: the oracle estimator must be simple separable but has access to the true $\bmu$, while the feasible estimator is non-separable but must work without knowledge of $\bmu$.

\subsection{The log-optimal simple separable compound e-values}
\label{subsec:best_simple_separable}
The connection of compound e-values to the fundamental theorem of compound decisions is born from the relaxation of the requirement of having a vector of e-variables to that of requiring only compound e-variables.

To further clarify the connection, we provide a construction of optimal simple separable compound e-values in sequence models with dominated null marginals via the fundamental theorem of compound decisions.

Suppose that $X=(X_k : k \in \mathcal{K})$, where $X_k$ takes values in a space $\mathcal{Y}$. For any distribution $\s \in \cM_1$, let $\s_k$ be the $k$-th marginal of $\s$, that is, the distribution of $X_k$ when $X \sim \s$.

Suppose that for all $k \in \mathcal{K}$, the $k$-th null hypothesis 
 is specified as a simple point null hypothesis on the $k$-th marginal, i.e.,  $$\cP_k = \{\s \,:\, \s_k = \p_{k}\}$$ for some prespecified distribution $\p_{k}$ with $\d \nu$-density equal to $p_{k}$, where $\nu$ is a common dominating measure. 
 We assume that there exists at least one possible true data-generating distribution $\p$ such that $\mathcal{N} = \mathcal{K}$ (as would be the case, e.g., if the $X_k$ are independent and we take $\p = \bigtimes_{k \in \mathcal{K}}  \p_{k}$).

For $k \in \mathcal{K}$, we let $\q_k$ be distributions on $\mathcal{Y}$ with $\d\nu$-densities $q_k$, representing the alternative distributions. We seek to solve the following optimization problem: 
\begin{equation}
\label{eq:optim_compound}
    \begin{aligned}
    & \underset{s: \mathcal{Y} \to [0,\infty] }{\text{maximize}} 
    & & \frac{1}{K}\sum_{k \in \mathcal{K}} \E^{\q_k}[\log(s(X_k))] \\
    & \text{subject to}
    & &  s(X_1),\dots,s(X_K) \text{~are compound } \\ 
    &
    & & \text{e-variables  for } \cP_1,\dots,\cP_K.
    \end{aligned}
\end{equation}

\begin{theorem}[Log-optimal simple separable compound e-values]
\label{th:compound-optimal}
The optimal solution to optimization problem~\eqref{eq:optim_compound} is given by the likelihood ratio of mixtures over the null and alternative:
\begin{equation}
\label{eq:optimal_evalue}
s(x) = \frac{\sum_{j \in \mathcal{K}} q_j(x)}{\sum_{ j \in \mathcal{K}} p_j(x)}.
\end{equation}
\end{theorem}
 
\begin{proof}[Proof.]
Let $\p$ be any probability measure with $k$-th marginal given by $\p_k$ for all $k \in \mathcal{K}$. 
We will first solve the optimization problem subject to the weaker constraint that $s(X_1),\dotsc,s(X_K)$ are compound e-values for $\p_1,\dots,\p_K$.

Define the mixtures  $\overline{\p} = \sum_{k \in \mathcal{K}} \p_k/K$ and $\overline{\q} = \sum_{k \in \mathcal{K}} \q_k/K$. Then, for any fixed $s : \mathcal{Y} \to [0,\infty]$, we have the following equalities:
\begin{align*}
\frac{1}{K}\sum_{k \in \mathcal{K}} \E^{\q_k}[\log(s(X_k))] &= \E ^{ \overline{\q}} [\log(s(X))],\; \\
\frac{1}{K}\sum_{k \in \mathcal{K}} \E^{\p_k}[s(X_k)] &= \E^ {\overline{\p}} [ s(X)],
\end{align*}
where $X$ is a generic random variable following the distribution specified with the expectation.
Thus, it suffices to solve the following optimization problem:
\begin{equation*}
    \begin{aligned}
    & \underset{s : \mathcal{Y} \to [0,\infty]}{\text{maximize}} 
    & & \E ^{\overline{\q}}[\log(s(X))] \\
    & \text{subject to} 
    & &  \E ^{\overline{\p}}[s(X)] \leq 1.
\end{aligned}
\end{equation*}
As seen in Chapter~\ref{chap:alternative}, this has a unique solution given by~\eqref{eq:optimal_evalue}. Finally we may verify that this $s$ indeed yields compound e-values for $\p_1,\dots,\p_K$. To this end, write $\overline{q}$ for the $\d \nu$-density of $\overline{\q}$. Then, we have the following:
$$\sum_{k: \p \in \cP_k} \mathbb E^{\p}[E_k] = \sum_{k: \p \in \cP_k} \int \frac{K \overline{q}(x)}{ \sum_{ j \in \mathcal{K}} p_j(x)}p_k(x)\nu(\d x) \leq \int K \overline{q}(x) \nu(\d x) = K,$$
where the first equality holds because for any $\p \in \cP_k$, its $k$-th marginal equals $\p_k$, and the inequality holds by dropping the constraint on $k$ in the sum.
This completes the proof.
\end{proof}

Note that any tight simple separable compound e-value for testing the nulls in the above problem must be of the form~\eqref{eq:optimal_evalue} for some densities $q_1,\dots,q_K$.

\begin{proposition}
Suppose that $E_1,\dots,E_K$ are tight simple separable compound e-values 
(for the same nulls defined as in the previous parts of this section). 
Then there exist $\d \nu$-densities $q_1,\dots,q_K$ such that $E_k$ may be represented as in~\eqref{eq:optimal_evalue} on the support of $\sum_{k \in \mathcal{K}} \p_k/K$ and it is without loss of generality to take $q_1=\dots=q_K$. 
\end{proposition}
\begin{proof}[Proof.]
First, by definition of simple and separable, there exists a function $s$ such that $E_k = s(X_k)$. The supremum in the definition of tightness must be attained when $\mathcal{N} = K$, and so:
$$\sum_{k \in \mathcal{K}} \E^{\p_k}[s(X_k)]=K.$$
Now define $\overline{p}$ as the $\d \nu$ density of $\overline{\p}$, where the latter object is defined as in~Theorem~\ref{th:compound-optimal}. Then the above may be written as:
$$\int s(x) \overline{p}(x) \d \nu(x) = 1.$$
Thus, if we define $q_k(x) = s(x)\overline{p}(x)$ for all $k \in \mathcal{K}$, we find that $q_k$ is indeed a $\d \nu$-density and $s(x)$ may be represented as in~\eqref{eq:optimal_evalue}.
\end{proof}

The objective value of the optimization problem~\eqref{eq:optim_compound} is given by $\kl(\overline{\q},{\overline{\p}})$. By convexity, $\kl({\overline{\q}},{\overline{\p}}) \leq \sum_{k \in \mathcal{K}} \kl({\q_k},{\p_k})$, which means that the optimal simple separable compound e-value has a worse objective than using the optimal separable e-value for each individual testing problem of $\p_k$ vs $\q_k$ given by $E_k^* = q_k(X_k)/p_k(X_k)$.

In view of the above, why is the optimal simple separable compound e-value of interest? The reason is that if $K$ is large (but the individual statistical problems have limited sample sizes), then it may be possible to mimic the log-optimal simple separable compound e-value without exact knowledge of the null $\p_k$, e.g., via empirical Bayes methods.

To elaborate, suppose we have a ``nuisance parameter'' $\theta \in \Theta$ that is common to all hypotheses, meaning that all null hypotheses are actually indexed by some $\theta$: $\cP^\theta_1,\dots,\cP^\theta_K$, but $\theta$ is unknown. If $\theta$ were known, we could simply find an e-value separately for each hypothesis. Since $\theta$ is unknown, we cannot do that, but since $\theta$ is common across all hypotheses, we could potentially learn $\theta$ by pooling data from all hypotheses, and then form (approximate) compound e-values from the data that asymptotically match the log-optimal simple separable compound e-values. 

Let us look at an explicit example that is in the spirit of the above discussion, but the resulting compound e-values are not of the above form.


\section{A minimally adaptive e-BH procedure}
\label{sec:c8-minimal}
 
We describe a  minimally adaptive procedure by proposing a tiny but uniform improvement of the e-BH procedure. This improvement is negligible for large values of $K$ and it may only be practically interesting for small $K$ such as $K\le 10$. 
We will focus on e-values and boosted e-values. 

First, choose an e-merging function $F:[0,\infty)^K \to [0,\infty)$ in Chapter~\ref{chap:multiple}. We allow for a general choice of $F$ other than the arithmetic mean $\fM_K$ as it will be useful for the case of boosted e-values.

With a chosen e-merging function $F$ and a level $\alpha\in (0,1)$, the improved e-BH procedure, denoted by $\cD^F_\alpha$, is designed as follows. We 
  first test the global null $\bigcap_{k=1}^K \cP_k$ via the rejection condition $F(e_1,\ldots,e_K)\ge 1/\alpha$, which has a type-I error of at most $\alpha$, and if the global null is rejected,  we then apply the e-BH procedure at level $\alpha'=K\alpha/(K-1)$.
In other words, 
\begin{enumerate}
\item if $F(e_1,\ldots,e_K) < 1/\alpha$, then $\cD^F_\alpha =\varnothing$;
\item  if $F(e_1,\ldots,e_K) \ge 1/\alpha$, then $\cD^F_\alpha=\cD_{\alpha'}$ where $\alpha'=K\alpha/(K-1)$ and $\cD_{\alpha'}$ is the e-BH procedure at level $\alpha'$.
\end{enumerate}
The next proposition shows that by choosing $F=\fM_K$, 
the resulting improved BH procedure dominates the base BH procedure.

\begin{proposition}\label{prop:tiny}
For $\alpha \in (0,1)$, the improved e-BH procedure $\cD^F_\alpha$ applied to arbitrary e-values has false discovery rate at most $\alpha$.
In case $F=\fM_K$, $\cD^{\fM_K}_\alpha$  dominates the  e-BH procedure $\cD_\alpha$, that is, $\cD_\alpha\subseteq \cD^{\fM_K}_\alpha$. 
\end{proposition}

\begin{proof}[Proof.]
Let
$A$ be the event that $F(E_1,\dots,E_K)\ge 1/\alpha$. 
If $K_0<K$, then  by using Theorem~\ref{th:compliant}, \begin{align*}  \E^{\p} \left [ \frac {F_{\cD^F_\alpha }}{R_{\cD^F_\alpha} } \right ]    &=  \E^{\p} \left [   \frac {F_{\cD_{\alpha'} }}{R_{\cD_{\alpha'}} } \id_A\right ]   + \E^{\p} \left[ \frac {F_{\varnothing}}{R_{\varnothing} } (1-\id_A) \right]    \\&  = \E^{\p} \left[\frac {F_{\cD_{\alpha'} }}{R_{\cD_{\alpha'}} } \id_A\right] \le \E^{\p} \left[ \frac {F_{\cD_{\alpha'} }}{R_{\cD_{\alpha'}} } \right] \le \frac{K_0}{K}\alpha' \le \alpha.
\end{align*}  

 If $K_0=K$, then the false discovery rate of  $\cD^F_\alpha$ is at most the probability $\p(A)$ of rejecting the global null via $F(E_1,\dots,E_K)\ge 1/\alpha$. 
 In this case, $\p(A)\le \alpha$ by Markov's inequality and the fact that $F$ is an e-merging function.
Hence, in either case, the FDR of  $\cD^F_\alpha$ is at most $\alpha$.

To show the   statement on dominance, let \begin{equation}\label{eq:e-simes}
S:(e_1,\dots,e_K)\mapsto \max_{k\in \cK} \frac { k e_{[k]}} K.\end{equation} The function
$S$    is an e-merging function  and it is dominated by $\fM_K$ because 
$\sum_{k=1}^Ke_k \ge ke_{[k]}$.
By definition of the e-BH procedure, we have $S(e_1,\dots,e_K)<1/\alpha$ implies  $\cD_\alpha=\varnothing$. 
Therefore, if $\fM_K(e_1,\dots,e_K)<1/\alpha$, then $\cD_\alpha=\varnothing = \cD^{\fM_K}_\alpha$.
Moreover, since $\alpha <\alpha'$, we always have $ \cD_\alpha \subseteq \cD_{\alpha'}$.
Hence, $ \cD_\alpha \subseteq \cD^{\fM_K}_\alpha$.   
\end{proof}

Next, we briefly discuss the case of boosted e-values. The arithmetic average of boosted e-values is not necessarily a valid e-value, so one must be a bit more careful.
Nevertheless, it turns out that
we can use  the function $S$ in \eqref{eq:e-simes} on the boosted e-values in Section~\ref{sec:c8-1-2}. The new procedure can be described as the following steps.
\begin{enumerate}
\item Boost the raw e-values with level $\alpha$.
\item If $S(e'_1,\ldots,e'_K)<1/\alpha$ where $e'_1,\ldots,e'_K$ are the boosted e-values in step 1, then return $\varnothing$.
\item Else: boost  the raw  e-values with level $\alpha'=K\alpha/(K-1)$.
\item Return the discoveries by applying the base e-BH procedure to the boosted e-values in step 3.
\end{enumerate}

This new procedure dominates the e-BH procedure, and it has FDR at most $\alpha$. The proof is   similar to that of Proposition~\ref{prop:tiny}.

\section{Stochastic rounding of e-values and randomized e-BH}
\newcommand{\alphath}{\widehat \alpha}
\label{sec:c8-stoch}

We next discuss how to improve the e-BH procedure by utilizing randomization.  
The randomized Markov's inequality in Theorem~\ref{thm:umi} gives rise to a randomized test based on e-values, and a randomized e-to-p calibrator. Here, we show how to use external randomization to convert e-values to other (stochastically rounded) e-values. 

Consider any closed set $\cG \subseteq [0,\infty]$ (where $\cG$ stands for ``grid'' since $\cG$ will often be countable or finite below).  Denote $g_* \coloneqq \inf\{x: x\in \cG\}$ and $g^* \coloneqq \sup\{x: x \in \cG\}$, and note that $g_*,g^* \in \cG$ since $\cG$ is closed. Further, for any $x \in [g_*, g^*]$, let 
\begin{align}
x^+ \coloneqq \inf\{y \in \cG: y \geq x\}, \qquad x_- \coloneqq \sup\{y \in \cG: y \leq x\},
\end{align} and note that $x^+, x_- \in \cG$ with $x^+ \geq x$ and $x_- \leq x$, and if $x \notin \cG$, then $x_- < x < x^+$.

Now define the stochastic rounding of $x \in [0,\infty]$ onto $\cG$, denoted by $S_\cG(x)$, as follows. If $x < g_*$, $x > g^*$, $x \in \cG$, or $x^+ = \infty$ then define $S_\cG(x)=x$.
Otherwise, define
\begin{align}
    S_\cG(x) =
\begin{cases}
 x_- \text{ with probability } \frac{x^+ - x}{x^+ - x_-}, \\
 x^+ \text{ with probability } \frac{x- x_-}{x^+ - x_-}.
\end{cases}
\label{eq:StochasticRounding}
\end{align}
Note that $S_\cG(x)$ need not lie in $\cG$, because if $x$ lies outside the range of $\cG$ then it is left unchanged. Also note that when $\cG = [0,\infty]$, we have $S_\cG(x)=x$ for all $x\in[0,\infty]$.
In what follows, $\p$ denotes the joint distribution of the random variable $X$ and the stochastic rounding $S_\cG$. 

One can also use a generalized version of rounding, where one stochastically rounds an input $x$ by sampling from a mixture distribution over all values in $\cG$ that are larger than $x$, instead of only $x^+$, but we do not discuss this idea further here.

\begin{proposition}\label{prop:stoch-round-e}
For any grid $\cG$, and any integrable random variable $X$, $\E^\p[X] = \E^\p[S_\cG(X)]$.
In particular, if $X$ is an e-variable, then $S_\cG(X)$ is also an e-variable.
\end{proposition}
The proof is simple: by design, $\E^\p[S_\cG(x)] = x$ for any real $x$, and thus when applied to any random variable, it leaves the expectation unchanged.

A key property is that $S_\cG(X)$ can be larger than $X$ with (usually) positive probability, since it can get rounded up to $X^+$, and is at least $X_-$, even when rounded down.

Remarkably, elements of the grid can also be random. For a concrete demonstration, consider the grid $$\cG(\widehat \alpha) = \{0,\widehat \alpha^{-1},\infty\}$$ and let $S_{\widehat \alpha}$ be shorthand for $S_{\cG(\widehat \alpha)}$.
 
\begin{proposition}
\label{prop:dynamic-round-e}
     If $X$ is an e-value for $\p$, and  $\widehat{\alpha}$ taking values in $ [0, 1]$ possibly depending on $X$, $S_{\widehat{\alpha}}(X)$ is also an e-value and satisfies $\E^\p[S_{\widehat{\alpha}}(X)] = \E^\p[X]$.
\end{proposition}
    To prove the above result, notice that another way of writing $S_{\widehat{\alpha}}$ is as follows:
\begin{equation}\label{eq:rand-e}
    S_{\widehat{\alpha}}(X) = (X \cdot \id_{\{X \geq \widehat{\alpha}^{-1}\}}) \vee (\widehat\alpha^{-1} \cdot \id_{\{X \geq U\widehat\alpha^{-1}\}}),
\end{equation}
where $U$ is a uniform random variable over $[0, 1]$ that is independent of $X$ and $\widehat\alpha$, and we treat $0 / 0 = 0$. Indeed, we can check that $\E^\p[S_{\widehat{\alpha}}(X)] = \E^\p[X]$ by first taking expectation with respect to $U$ (while conditioning on $X$) and then with respect to $X$.
\begin{proof}[Proof.]
We rewrite $S_{\widehat\alpha}(X)$ into the following equivalent form:
\begin{align}
    S_{\widehat\alpha}(X) = X \cdot \id_{\{X \geq \widehat{\alpha}^{-1}\}} + \underset{A}{\underbrace{\widehat{\alpha}^{-1} \cdot \id_{\{\widehat{\alpha}^{-1}
    > X \geq U\widehat{\alpha}^{-1}\}}}}, \label{eq:RandomERewrite}
\end{align}
where $U \lawis \mathrm{U}[0, 1]$ is independent of $X$ and $\widehat{\alpha}$.
First, we derive the following equality for the expectation of $A$ as indicated in \eqref{eq:RandomERewrite}:
\begin{align}
    \E^\p[A]
    &=\E^\p[\E^\p[\widehat{\alpha}^{-1} \cdot \id_{\{\widehat{\alpha}^{-1} > X \geq U\widehat{\alpha}^{-1}\}} \mid X, \widehat{\alpha}]]\\
    &= \E^\p[\widehat{\alpha}^{-1} \cdot \id_{\{X < \widehat{\alpha}^{-1}\}} \p(U \leq \widehat{\alpha} X \mid X, \widehat{\alpha})]\\
    &= \E^\p[\widehat{\alpha}^{-1} \cdot \id_{\{X < \widehat{\alpha}^{-1}\}} \cdot (\widehat{\alpha} X \wedge 1)] \\
    &= \E^\p[X \cdot \id_{\{X < \widehat{\alpha}^{-1}\}}]
    \label{eq:AUpperBound}.
\end{align} 
Then, we can simply upper bound the expectation of $S_{\widehat{\alpha}}(X)$ as follows.
\begin{align}
    \E^\p[S_{\widehat\alpha}(X)] &= \E^\p[X\cdot\id_{\{X \geq \widehat{\alpha}^{-1}\}}] + \E^\p[X \cdot \id_{\{X < \widehat{\alpha}^{-1}\}}] \\
    &= \E^\p[X] \leq 1, \label{eq:add-parts}
\end{align} 
where the first equality is by application of \eqref{eq:RandomERewrite} and \eqref{eq:AUpperBound}. Hence, we have shown our desired result. Note that if $U$ were to be a superuniform (stochastically larger than uniform) random variable, the last equality in \eqref{eq:add-parts} would be an inequality.
\end{proof}

We make a final important observation.

\begin{proposition}
    If $X \geq \widehat \alpha^{-1}$, then $S_{\widehat \alpha}(X) \geq \widehat \alpha^{-1}$.
\end{proposition}

The proof is omitted, since it is similar to that of the preceding result.

Above, we introduced stochastic rounding above for any random variable $X$.
Next, we apply the stochastic rounding to e-values and the e-BH procedure.
In this context,
    stochastic rounding trades off e-power for power. Indeed, it is not hard to check, by Jensen's inequality, that for an e-variable $E$, $$\E^\p[\log(S_{\widehat \alpha}(E))] \leq \E^\p[\log E].$$ Thus, the e-power of a stochastically rounded e-value is always at most that of the original e-value, while its power (at level $\widehat \alpha$) is at least that of the original e-value.

Fix $\alpha \in (0,1)$, and let 
\begin{gather}
    \alpha_i = \alpha i / K, \qquad \cG_\alpha = \{\alpha_i^{-1}: i \in \mathcal K\} \cup \{0, \infty\}
\end{gather}
denote the set of possible levels that e-BH may reject e-values at and in addition to $0$ and $\infty$. If $E \geq K / (\alpha k)$ for some $k \in \cK$, then $S_{\cG_\alpha}(E) \geq K / (\alpha k)$ as well. Thus, an e-value that is stochastically rounded to $\cG_\alpha$ can only improve power when used in conjunction with e-BH.  



One can view the e-BH procedure as selecting a data-dependent threshold
\begin{align}
\label{eq:alpha-hat-star}
    \alphath^* \coloneqq \alpha(k^* + 1) / K
\end{align}
and rejecting the $i$th hypothesis if and only if $E_i \geq 1 / \alphath^*$ (indeed, if $k^*=K$, the claim is trivially true, and if $k^*<K$, then there can only be $k^*$ hypotheses with e-values larger than $K/(\alpha(k^*+1))$, otherwise we violate the maximality of $k^*$ in its definition). Now, observe that $1 / \alphath^*$ can only take values from the grid $\{K / (k\alpha): k \in \cK\}$.
If we ``rounded'' down each e-value $E_i$ to the closest value in the grid that is less than $E_i$ (or 0 if $E_i$ is smaller than any value in the grid), the discovery set that is output by e-BH would be identical to the one where no rounding had occurred. But if we could round \emph{up} the e-value, we would potentially gain power; however, this would inflate its expectation and it would no longer be an e-value. Our key insight is that if we \emph{randomly} rounded every e-value up or down --- appropriately so that its expectation is unchanged --- then we could increase power with positive probability.

The fact that e-values are often continuous and will typically lie between grid points provides a broad opportunity to significantly increase the power. In what follows, we use independent external randomness to stochastically round e-values, as introduced previously, and increase the number of discoveries made by e-BH.



Let $k^*$ be the number of discoveries made by e-BH applied to e-values $E_1,\dots,E_K$. 
The \emph{grid-rounded e-BH (Ge-BH)  procedure}  simply applies the e-BH procedure to the set of e-values $\{S_{\cG_\alpha}(E_{k})\}_{k \in \cK}$.
Let $\cD_{\mathrm{GeBH}}$  be the set of rejections made by Ge-BH,
and $\cD_{\mathrm{eBH}}$ be the set of rejections made by e-BH at level $\alpha.$
%
\begin{theorem}
\label{thm:grid-ebh}
For any arbitrarily dependent e-values $(E_1, \dots, E_K)$, the Ge-BH procedure ensures $\FDR \leq \alpha$ and $\cD_{\mathrm{GeBH}}\supseteq \cD_{\rm eBH}$. Further, for any $\p \in \cM_1$, $\p(\cD_{\mathrm{GeBH}} \supsetneq \cD_{\rm eBH}) > 0$, i.e., the probability that Ge-BH makes extra discoveries over e-BH is positive, if and only if 
\begin{equation}
    \p\left(\exists k \in [K - k^*]: E_{[k^* + j]} > \frac{K}{\alpha(k^* + k + 1)}, ~\forall j \in [k]\right) > 0. \label{eq:kcal-cond}
\end{equation}
\end{theorem}

The proof follows from the fact that if $E_i$ ever takes on a value that is between levels in $\cG_\alpha$, the e-BH procedure will  make the same rejections as  if $(E_i)_{-}$ were substituted in its place. Stochastic rounding guarantees that $S_{\cG_\alpha}(E_i) \geq (E_i)_{-}$ almost surely, so it can only increase the number of rejections. Further, when $E_i$ is between two levels in $\cG_\alpha$, then $S_{\cG_\alpha}(E_i) =  (E_i)_{+} > E_i$ with positive probability, which leads to rejecting hypotheses that e-BH did not reject. This intuition leads us to the condition in \eqref{eq:kcal-cond} for which Ge-BH has strictly more power than e-BH. We omit the proof.

It is worth remarking that \eqref{eq:kcal-cond} is an extremely weak condition that would be very frequently satisfied. For example, if the e-values are independent and continuously distributed over $[0, K / \alpha ]$ (or a larger interval), then \eqref{eq:kcal-cond} will hold.  As an explicit example, if the data $Z_i$ for testing the $i$-th hypothesis are Gaussian with variance $\sigma^2$, and we are testing whether the mean of $Z_i$ is nonpositive (or equal to zero) against the alternative that the mean is positive, all admissible e-values are  mixtures of likelihood ratios between a positive mean Gaussian and a zero mean Gaussian: these likelihood ratio e-values take the form $\exp(\lambda Z_i - \lambda^2\sigma^2/2)$ for $\lambda > 0$, which are clearly continuous and unbounded, as are many other e-values for testing parametric and nonparametric hypotheses. The independence mentioned above is far from necessary, but it is sufficient to ensure that the probability in \eqref{eq:kcal-cond} is not pathologically equal to zero due to some awkward worst-case dependence structure.

One can further improve the above procedure with more randomization. We now describe a more powerful procedure below, as applied to input e-values $(E_1,\dots,E_K)$.
This procedures makes two iterations of rounding, first with the grid $\mathcal G_\alpha$, and then  with  
the data-dependent grid $\{0, (\widehat \alpha^*)^{-1},\infty\}$,  where $\widehat \alpha^*$ is in in~\eqref{eq:alpha-hat-star}. This procedure is explained in Algorithm~\ref{alg:sebh}, and it is called the \emph{doubly-rounded e-BH (De-BH) procedure}, with 
 $\cD_{\rm DeBH}$ denoting its discovery set resulting.

\begin{algorithm}[t]
    \caption{The De-BH procedure on input e-values $(E_1, \dots, E_K)$.\label{alg:rboth}}
\label{alg:sebh}
    Apply the e-BH procedure to $(S_{\cG_\alpha}(E_1), \dots, S_{\cG_\alpha}(E_K))$.\\
    Let $k^*$ denote the number of rejections above.
    Define $\widehat \alpha^*=\alpha(k^*+1)/K$.\\
    Reject the $i$th hypothesis if $S_{\widehat\alpha^*}(S_{\cG_\alpha}(E_i)) \geq (\widehat \alpha^*)^{-1}$. \\
\end{algorithm}

Finally, we mention the \emph{uniformly-randomized 
e-BH
(Ue-BH) procedure}: this runs the e-BH procedure on $(E_1/U,\dots,E_K/U)$, where $U$ is an independent $\mathrm{U}[0,1]$ random variable. Let $\cD_{\rm UeBH}$ be the resulting discovery set.

\begin{theorem}
\label{thm:rboth-ebh}
    For any arbitrarily dependent e-values $(E_1, \dots, E_K)$, the De-BH and Ue-BH procedures both control $\FDR$ at level $\alpha$. Further, the discoveries of De-BH is a superset of that of Ge-BH and e-BH, and the discoveries of Ue-BH is a superset of that of e-BH.
    Further, $\p(\cD_{\rm DeBH} \supsetneq \cD_{\mathrm{GeBH}}) > 0$ and $\p(\cD_{\rm UeBH} \supsetneq \cD_{\rm eBH})>0$ if and only if 
    \begin{equation}
    \label{eq:cond-higher-power}
    \p(\exists k \in \cK: E_k \in (0,\widehat{\alpha}^{*-1})) > 0.
    \end{equation}
\end{theorem}

The proof is omitted, but follows from earlier observations on stochastic rounding.

As an important special case, since the BHY procedure for FDR control with arbitrarily dependent p-values (Definition~\ref{defi:BH}) can be reinterpreted as the e-BH procedure applied to e-values obtained via a particular p-to-e calibrator, the above results also mean that the BHY procedure can be dominated by a randomized variant, which always rejects a superset of discoveries and yet controls the FDR at the same level.

\paragraph{A note on randomization.} While the set of randomized multiple testing procedures contains all deterministic ones, it is far from clear when the most powerful randomized procedure is strictly more powerful (in at least some situations) than the most powerful deterministic one.
It could be the case that randomization confers no additional benefit in any situation. 
As an important counterpoint, we do not know of any way for randomization to improve the power of other procedures such as the BH procedure under positive dependence, and indeed we conjecture that it does not. Both the BH and BHY procedures have their FDR upper bounds being achieved with equality in a particular setting; yet it appears as if the power of BHY can be essentially improved in almost all other situations, while we conjecture that without further knowledge, the BH procedure cannot be strictly improved in any situation without worsening its power in other situations. Thus, it is not the case that randomization should indeed always or ``naturally'' help in multiple testing.

\begin{figure}[t]
    \centering
    \begin{subfigure}{\textwidth}
        \centering
    \includegraphics[width=0.6\textwidth]{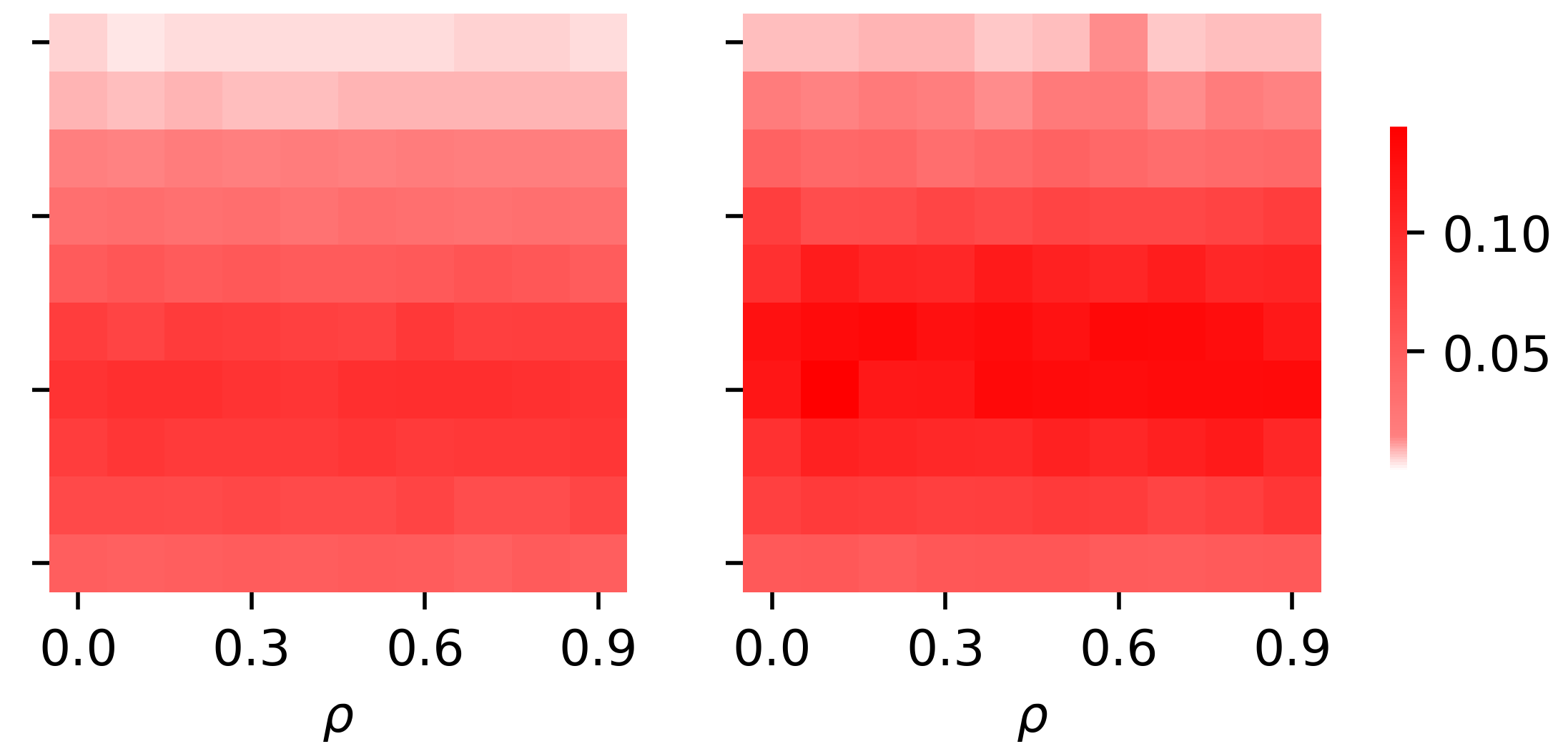}
    \caption{Positive dependence.}
    \end{subfigure}
    \begin{subfigure}{\textwidth}
        \centering
    \includegraphics[width=0.6\textwidth]{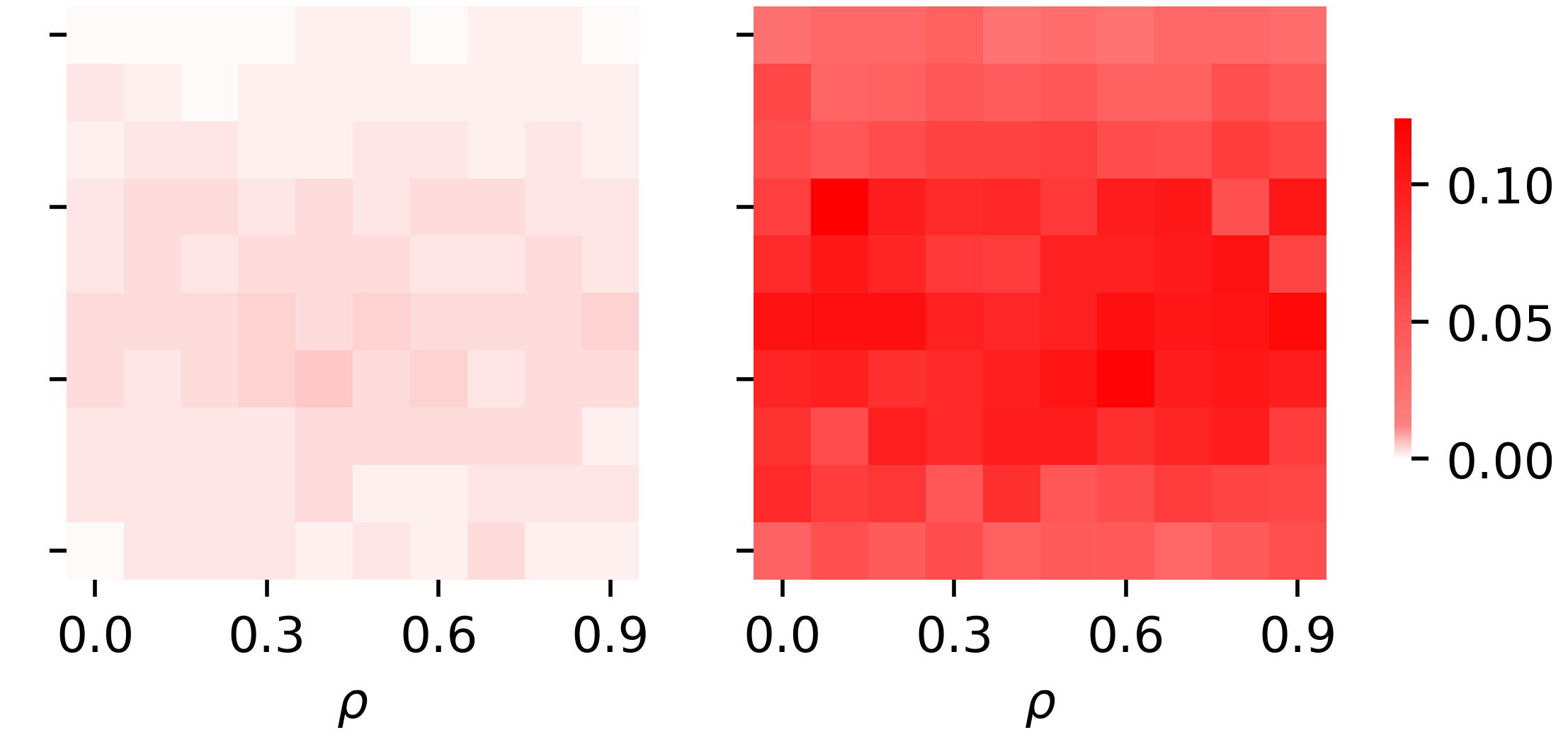}
    \caption{Negative dependence.}
    \end{subfigure}
    \caption{Heatmaps of the power improvement of De-BH over e-BH (left) and Ue-BH over e-BH (right),  across various values of $\mu$ (y-axis) and $\rho$ (x-axis).
    }
    \label{fig:ebh-heatmap}
\end{figure}

\paragraph{Numerical experiments.}
In order to gauge the amount of improvement that randomization delivers, we run simulations of the e-BH procedure, where we test $K = 100$ hypotheses. Let $\pi_0 = (K - |\cN|)/ K = 0.3$, i.e., the proportion of hypotheses where the null is false. For each hypothesis, we sample $Z_i \lawis \mathrm{N}(\mu_i, 1)$, and perform the following one-sided test:
\begin{align}
    H_0:&\ \mu_i = 0\ \text{vs.}\ H_1:\ \mu_i \geq 0,
\end{align} for each $i \in \cK$. We consider two dependence settings:
\begin{enumerate}
    \item Positive dependence: $\cov(Z_i, Z_j) = \rho^{|i - j|}$ for each $i, j \in \cK$, i.e., the covariance matrix is Toeplitz.
    \item Negative dependence: $\cov(Z_i, Z_j) = -\rho / (K - 1)$, i.e., the covariance is equal among all pairs of samples.
\end{enumerate}
We let $\mu$ range in $[1, 4]$, and $\rho$ range in $[0, 0.9]$.
For each $i \in \cK$, where the null hypothesis is false, we let $\mu_i = \mu$, i.e., all non-null hypotheses have a single mean. Power is defined as the expected proportion of non-null hypotheses that are discovered, i.e., $\E[|\cD \backslash \cN| / (K - |\cN|)]$. We calculate the power averaged over 500 trials. The e-variable used is $E_i = \exp(\lambda Z_i - \lambda^2/2)$, with $\lambda = \sqrt{2 \log(1/\alpha)}$ and $\alpha=0.05$. The results are summarized in Figure~\ref{fig:ebh-heatmap}.

\section{The closed e-BH procedure}
\label{sec:9-closed-ebh}

We end this chapter by describing a procedure that (usually strictly) dominates the e-BH procedure. Due to strong connections with the classical closure principle in family-wise error rate control, we call it the \emph{closed} e-BH procedure. This section has four goals:
\begin{itemize}
    \item To improve the e-BH procedure, when starting with $K$ e-values, one for each hypothesis.
    \item To recover any given FDR procedure for any problem as an instance of the closed e-BH procedure.
    \item To recover and/or improve any given post-hoc FDR procedure using the closed e-BH procedure.
    \item To point out that the techniques developed are applicable to many error metrics beyond FDR.
\end{itemize}

First, let us show that the closed e-BH procedure is both universal and unimprovable for FDR control. 

\subsection{Improving the e-BH procedure}

For any $A \subseteq \cK$, define the intersection hypothesis $\cP_A := \bigcap_{k \in A} \cP_k$. 
The set
$\{E_A\}_{A \subseteq \cK}$ is called an \emph{e-collection} if $E_A$ is an e-variable for $\cP_A$ for every $A \subseteq \cK$ (see Section~\ref{sec:FWER}). We take $E_\emptyset = 1$.

If the e-collection is not explicitly specified, then we assume that we only have access to one e-variable $E_k$ for each base hypothesis $\cP_k$, $k \in \cK$, which  induces  an e-collection by defining   
\[
E_A = \frac1{|A|}\sum_{k \in A}E_k,~~~~A\subseteq \cK.
\]
We refer to the above as \emph{the mean e-collection}. For any sets  $A \subseteq \cK$ and  $R \subseteq \cK$, define
\[
\FDP_A(R) = \frac{|A \cap R|}{|R| \vee 1},
\]
which is the false discovery proportion of $R$ if (hypothetically) the set of all true nulls were $A$. Using this notation, the actual FDP of a procedure that outputs $R$ is $\FDP_\cN(R)$, which is  unknown since $\cN$ is unknown.

\begin{definition}[Closed e-BH procedure]
Let $\{E_A\}_{A \subseteq \cK}$ be an e-collection. Define the \emph{candidate discovery sets} as
\begin{equation}
\label{eq:candidates-closed}
\mathcal C_\alpha = \left\{ R \subseteq \cK: E_A \geq \frac{\FDP_A(R)}{\alpha} \text{ for all }A \subseteq  \cK \right\}.
\end{equation}
Define an \emph{e-closed} procedure as a procedure that outputs some set $R \in \mathcal C_\alpha$. 
Finally, define the \emph{closed e-BH procedure} ($\overline{\text{e-BH}}$) at level $\alpha$ as the procedure that rejects a largest set in $\mathcal C_\alpha$ (if there are many largest sets of equal size, we may pick any of them). 
\end{definition}

When working with the mean e-collection derived from $\cK$ base e-values, let $R_{[1:k]}$ denote the indices corresponding to the largest $k$ base e-values, with $R_0 = \emptyset$. One can easily check that in this case, the largest set in $\mathcal C_\alpha$ equals $R_{[1:\bar k]}$, where
\[
\bar{k} = \max\{k \geq 0: R_{[1:k]} \in \mathcal C_\alpha\}.
\]

For a comparison to e-BH, recall that the set of all self-consistent discovery sets is given by:
\[
\mathcal C_\alpha' = \left\{ R \subseteq \cK: \min_{k \in R} E_k \geq \frac{K}{\alpha|R|}  \right\},
\]
where we insist that $\emptyset \in \mathcal C_\alpha, \mathcal C_\alpha'$ by convention, and note that the e-BH procedure rejects the largest self-consistent set $R_{[1:k^*]}$.
Also, let $k^\sharp$ represent the number of rejections of the minimally adaptive e-BH procedure in Section~\ref{sec:c8-minimal}.

\begin{theorem}
\label{th:closed-ebh}
   Any procedure that reports some $R \in \mathcal C_\alpha$, and in particular the closed e-BH procedure,  controls the FDR at level $\alpha$.
 Moreover, the level post-hoc guarantee holds for any $\p\in \cM_1$:
    \[ \E^\p\left[\sup_{\alpha \in (0,1)} \sup_{R \in \mathcal C_\alpha} \frac{\FDP_\cN (R)}{\alpha}\right] \leq 1.
    \]
    Further, for the mean e-collection, $\mathcal C_\alpha' \subseteq \mathcal C_\alpha$; thus, $\overline{\text{e-BH}}$ always rejects a superset of e-BH. In fact, $\overline{\text{e-BH}}$ also rejects a superset of the minimally adaptive e-BH procedure, and thus $\bar k \geq k^\sharp \geq k^*$.
\end{theorem}
\begin{proof}[Proof.]
    The proof is straightforward. First, fix $\p$ and thus $\cN=\cN(\p)$, which we assume is nonempty; otherwise the proof is trivial. Choosing $A=\cN$, we see that every set $R \in \mathcal C_\alpha$ satisfies $\FDP_{\cN}(R) \leq \alpha E_{\cN}$, and thus 
    \[
         \E^{\p}\left[\sup_{\alpha \in (0,1)} \sup_{R \in \mathcal C_\alpha} \frac{\FDP_\cN(R)}{\alpha}\right] \leq  \E^{\p}[E_{\cN}] \leq 1,
    \]
    completing the proof of the first claim. 
    
    To prove  $\mathcal C_\alpha' \subseteq \mathcal C_\alpha$, note the deterministic inequality for any non-empty $A,R \subseteq \cK$:
    \begin{equation}
    \label{eq:determinstic-eq}
     E_A = \frac{1}{|A|}\sum_{k \in A}E_k \geq \frac{ |A \cap R|}{|A|}  \min_{k \in A\cap R} E_k \geq \frac{ |A \cap R| }{K} \min_{k \in  R} E_k.
    \end{equation}
    For any non-empty self-consistent set $R \in \mathcal C_\alpha'$, we have $\min_{k \in  R} E_k \geq K/(\alpha|R|)$, and so the above inequality simplifies to
    \[
    E_A \geq \frac{|A \cap R|}{\alpha|R|},
    \]
    which implies that $R \in \mathcal C_\alpha$, completing the proof of the second claim.

    We now argue that $R_{[1:k^\sharp]}$, the rejection set of minimally adaptive e-BH, is also a candidate set, that is $R_{[1:k^\sharp]} \in \mathcal C_\alpha$, so that  $R_{[1:\bar k]} \supseteq R_{[1:k^\sharp]}$. We naturally consider the case where $k^\sharp \geq 1$, implying that $E_{\cK} \geq 1/\alpha = \FDP_{\cK}(R_{[1:k^\sharp]})/\alpha$. For any other $A$ such that $|A| \leq K-1$. Altering the last step of~\eqref{eq:determinstic-eq} to take $|A| \leq K-1$ into account yields
    \[
     E_A  \geq \frac{1}{K-1}  |A \cap R_{[1:k^\sharp]}| \min_{k \in  R_{[1:k^\sharp]}} E_k \geq \frac{|A \cap R_{[1:k^\sharp]}|}{\alpha |R_{[1:k^\sharp]}|} = \FDP_{A}(R_{[1:k^\sharp]})/\alpha,
    \]
    where the second inequality follows because the minimally adaptive e-BH satisfies $\min_{k \in  R_{[1:k^\sharp]}} E_k \geq (K-1)/(|R_{[1:k^\sharp]}| \alpha)$. This completes the proof.
\end{proof}

Though we do not present it here, the closed e-BH procedure can be run in $O(K^3)$ time when using the mean e-collection.  For other e-collections, as with the usual closed testing procedure, the rejection set of closed e-BH may take exponential time to compute. 

To better illustrate the features of the closed e-BH procedure, we take a close look at the thresholds on the individual e-values for rejection. 
Choose $A = \{j\}$ for each $j \in R_{[1:\bar k]}$, we see that~\eqref{eq:candidates-closed}  enforces the individual constraint that
\[
E_j \geq \frac{1}{\alpha \bar k} \text{ for every } j \in R_{[1:\bar k]}.
\]
(Of course, it imposes other ``joint constraints'' on averages of these e-values as well: for instance, the average of the smallest $m$  rejected e-values must exceed $m/(\alpha \bar k)$.)
In contrast, the e-BH procedure satisfies
\[
E_j \geq \frac{K}{\alpha k^*} \geq \frac1\alpha , \text{ for every } j \in R_{[1:k^*]},
\]
which is a significantly harsher threshold.
Meanwhile, the minimally adaptive e-BH procedure  that makes $k^\sharp \geq k^*$ rejections satisfies
\[
E_j \geq \frac{K-1}{k^\sharp \alpha } \text{ for every } j \in R_{[1:k^\sharp]},
\]
which again is a much more stringent condition.

\begin{example}[Closed e-BH improves minimally adaptive e-BH]
    We give a simple example with $K=3$ to show that the closed e-BH rule can reject strictly more than the minimally adaptive e-BH rule. For the set of input e-values (60, 39, 11), the e-BH procedure only makes one rejection at $\alpha = 0.05$, narrowly missing out on the second rejection. The minimally adaptive e-BH procedure effectively runs e-BH at level $3\alpha/2$, and makes two rejections.
    The closed e-BH procedure rejects all 3 hypotheses, which can be checked by verifying that the set $\{1,2,3\}$ lies in $\mathcal C$ in~\eqref{eq:candidates-closed}, since each e-value is larger than 20/3, all pairwise averages exceed 40/3 and the overall average exceeds 20.

    Why is it even possible to reject e-values below $1/\alpha$, indeed as small as $1/(3\alpha)$ in the above example, while still controlling FDR at $\alpha$? The informal explanation is as follows. Assume for the moment that the first two e-values are so large that it is almost without doubt that they are non-nulls, while the third hypothesis is indeed null. Then, the FDP can only be 0 or 1/3; it is 1/3 only when the third e-value exceeds $1/(3\alpha)$, which happens with probability at most $3\alpha$ by Markov's inequality. Thus, the FDR is still bounded by $\alpha$ in this case. This intuition may be helpful to understand the more abstract proof of Theorem~\ref{th:closed-ebh}.
\end{example}

Our results also have implications for the 
BHY procedure  
in Definition~\ref{defi:BH}.  Recall from Proposition~\ref{prop:c8-recover-BY}  that the BHY procedure is an instance of the e-BH procedure with e-values formed via calibrating the input p-values using a particular p-to-e calibrator. Define the \emph{closed BHY} procedure as the one obtained by using the closed e-BH procedure with the mean e-collection formed by the aforementioned e-values. Both procedures output a set of the form $R_{[1:k]}$.

\begin{corollary}[Improving the BH procedure under arbitrary dependence]
\label{cor:closedby-beats-by}
The discoveries made by the closed BHY procedure are a (usually strict) superset of those made by the BHY procedure, while still controlling the FDR at level $\alpha$. Thus, the BHY procedure is not admissible in the sense that it can be improved while maintaining the same FDR control. 
\end{corollary}

\subsection{Recovering any FDR procedure}

This subsection will show that that every FDR controlling procedure can be written as the closed e-BH procedure applied to an appropriate e-collection.



\begin{theorem}
\label{th:closed-ebh-universal}
Consider an arbitrary FDR controlling procedure $\cD$. 
Then, for each $\alpha$, there exists an e-collection $\{E_A\}_{A \subseteq \cK}$ for which the corresponding closed e-BH procedure recovers $\cD$ at level $\alpha$.
\end{theorem}
\begin{proof}[Proof.]
Suppose that running $\cD$ at level $\alpha$ on some dataset produces a rejection set $R = R(\alpha)$. Now, define the e-collection $\{E_A\}_{A \subseteq \cK}$ given by 
\begin{equation}
\label{eq:e-collection-recovery}
E_A = \FDP_A(R)/\alpha.
\end{equation}
    First note that each $E_A$ as defined above is an e-value for $\cP_A$ because $\cD$ is an FDR controlling procedure. 
    Defining the candidates $\mathcal C_\alpha$ by~\eqref{eq:candidates-closed}, we observe that $R \in \mathcal C_\alpha$ by construction. 
    Indeed, one can check that $\mathcal C_\alpha = \{\varnothing, R\}$: no other nonempty set $R' \neq R$ can lie in $\mathcal C_\alpha$ because the inequality $E_A \geq \FDP_A(R')/\alpha$ cannot hold for all $A$, noting that $|A\cap R'|/|R'|>|A\cap R|/|R|$ for $A=R'$ when $R'\subseteq R$ and $A=R'\setminus R$ otherwise.
    Thus the only set in $\mathcal C_\alpha$ is $R$ itself, and the closed e-BH procedure rejects it.
\end{proof}

The e-collection in the above proof is worth further examination. If $R = \emptyset$, then all e-values equal zero, so we will assume for simplicity that $R$ is nonempty. For every $j \in R$, we have $E_j = 1/(\alpha|R|)$, and for every $j \notin R$, we have $E_j=0$. Meanwhile, $E_{\cK} = E_R =  1/\alpha$ and $E_{\cK \backslash R} = 0$. Consider the directed acyclic graph formed by $2^{\cK}-1$ nodes, indexed by nonempty sets $A \subseteq \cK$, where the $K$ leaves are the singletons, the root is the set $\cK$, and a node $A$ is a parent of node $B$ if $A \supsetneq B$ and $|A \backslash B|=1$. Then these e-values are increasing from leaves to the root. Adding nonrejected elements to any $A$ does not change its assigned e-value, but adding rejected elements to $A$ increases its e-value (from $1/(|R|\alpha)$ in the leaf to $1/\alpha$ in the root).

We end by noting that there is no contradiction that every FDR procedure can be written both as the e-BH procedure with compound e-values and as the closed e-BH procedure with an e-collection.

In fact, we can also write down a version of the closed e-BH procedure that works with compound e-values. Indeed, given an input vector $\widetilde{\mathbf{E}}$ of compound e-values, we can define an e-collection using $E_A = K^{-1} \sum_{k \in A} \widetilde E_k$. Then, one can show that the corresponding closed e-BH procedure makes more discoveries than the e-BH procedure, using a proof nearly identical to that of Theorem~\ref{th:closed-ebh}. We omit the details.

\subsection{Recovering or improving any post-hoc FDR procedure}

Recall that a post-hoc FDR procedure $\cD$ is one that produces a family of rejection sets $\{R_\alpha\}_{\alpha \in (0,1)}$ that satisfies 
\[ 
\sup_{\p \in \cM_1}
\E^\p\left[\sup_{\alpha \in (0,1)} \frac{\FDP_\cN (R_\alpha)}{\alpha}\right] \leq 1.
\]
(By default, take $R_0=\emptyset$ and $R_1=\cK$ without loss of generality.)
We now claim that we can recover any such method using the closed e-BH procedure.

\begin{theorem}
    Consider an arbitrary post-hoc FDR procedure $\cD$. Then, there exists an e-collection $\{E_A\}_{A \subseteq \cK}$ for which the corresponding closed e-BH procedure (at any $\alpha$) recovers or improves $\cD$ (at that $\alpha$).
\end{theorem}
\begin{proof}[Proof.]
For each $A \subseteq \cK$, define the e-value
\begin{equation}
\label{eq:e-collection-posthoc}
\bar E_A = \sup_{\beta \in (0,1]} \frac{\FDP_A(R_\beta)}{\beta}.
\end{equation}
Then $(\bar E_A)_{A \subseteq \cK}$ form an e-collection because $\cD$ is a post-hoc FDR procedure. 

At any (potentially data-dependent) $\alpha \in (0,1)$, we form $\mathcal C_\alpha$ using~\eqref{eq:candidates-closed}. Observe that $R_\alpha \in \mathcal C_\alpha$ by construction. Thus, the closed e-BH procedure either rejects $R_\alpha$ or a larger set.
\end{proof}

Let us briefly compare the construction of the   e-collection in~\eqref{eq:e-collection-posthoc} to~\eqref{eq:e-collection-recovery} in the preceding subsection. First, in the former and unlike the latter, the e-collection does not depend on $\alpha$ ($E_A$ depends on $\alpha$, but $\bar E_A$ does not). 
Second, if the e-BH procedure (that takes $\cK$ base e-values as input) is plugged into Theorem~\ref{th:closed-ebh-universal}, we get back the e-BH procedure itself (as Theorem~\ref{th:closed-ebh-universal} claims). Every e-value in the above e-collection $\{\bar E_A\}$ is at least as large as the preceding e-collection $\{E_A\}$, typically strictly larger, and so the above construction will yield a procedure strictly improving e-BH, at least at some level $\alpha$.
Note also that the previous construction sets $\bar E_j = 0$ every $j \notin R_{k^*}$. However, the above construction has no e-values equal to 0 because every index is eventually included in $R_\beta$ for some $\beta$.

%

\subsection{A general e-closure principle for arbitrary error metrics}

    What we have implicitly developed in this section is a much more far-reaching ``e-closure'' principle that applies to any error metric, not just the FDR. To elaborate, let $L_A(R)$ represent an arbitrary error metric that quantifies the quality of the selected set $R$, when the true set of null hypotheses happens to equal $A$. Suppose we would like to design procedures $\cD$ whose rejected sets $R$ satisfy
    \begin{equation}
    \label{eq:closure-general-loss}
    \sup_{\p \in \cM_1} \E^\p[L_\cN(R)] \leq \alpha,
    \end{equation}
    for any $\alpha$ in some range; usually the losses $L$ are nonnegative and bounded, and so the range can often be taken to be $[0,1]$, as we will do below for simplicity.
    In the case of FDR control, we had $L_A = \FDP_A$. But we could consider alternatives such as 
    \begin{itemize}
        \item $k$-familywise error rate ($k$-FWER):  $L_A(R) = \id_{\{|A\cap R| \geq k\}}$,
        \item per-family error rate (PFER):  $L_A(R)=|A\cap R|/K$ (the division by $K$ is to rescale the range to $[0,1]$).
        \item false discovery exceedence (FDX): $\id_{\{\FDP_A(R)>\gamma\}}$ for some $\gamma \in (0,1)$.
    \end{itemize}
    Let $\cL$ denote the set of all such loss functions.
    
    Analogous to what was done in this section, given an e-collection $\{E_A\}_{A \subseteq \cK}$, we could define the candidate set for a particular $L$ as
    \[
    \mathcal C_\alpha(L) = \left\{R \subseteq \cK : E_A \geq \frac{L_A(R)}{\alpha} \text{ for all } A \subseteq \cK \right\}.
    \]
    Then, any set $R \in \mathcal C_\alpha(L)$ satisfies~\eqref{eq:closure-general-loss}, and we could, in particular, choose the largest such set to reject. For example, for FWER control (i.e.\ $k$-FWER control with $k=1$), e-closure with the mean e-collection recovers the e-Holm procedure, whose details are omitted here for brevity.

    Further, if we begin with any procedure $\cD$ that guarantees \eqref{eq:closure-general-loss}, it can be recovered or improved by applying our e-closure technique on the e-collection defined by $E_A = L_A(R)/\alpha$. In fact, we would get the stronger claim that
    \[
    \sup_{\p \in \cM_1} \E^\p \left[ 
    \sup_{L \in \cL}
    \sup_{\alpha \in [0,1]} \sup_{R \in \mathcal C_\alpha(L)} \frac{L_\cN(R)}{\alpha}\right] \leq 1.
    \]
    In particular, it means that both the loss function $L$ and the level $\alpha$ can be chosen post-hoc. For example, if the scientist does not like the rejection set achieved with FWER control at level $0.01$, they can report the results when requiring FDR control at level $0.05$ (and indeed, can do through all such loss functions and error levels before deciding what to report).

    Further developments of e-closure are omitted here.

\subsection{Numerical experiments}

\begin{figure}[t]
\centering
    \begin{subfigure}{0.33\textwidth}
        \includegraphics[trim=0 0 3.2cm 0,clip,width=\textwidth]{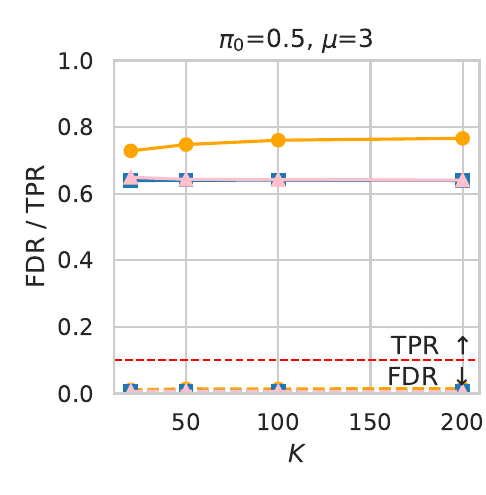}
    \end{subfigure}\begin{subfigure}{0.33\textwidth}
        \includegraphics[trim=0 0 3.2cm 0,clip,width=\textwidth]{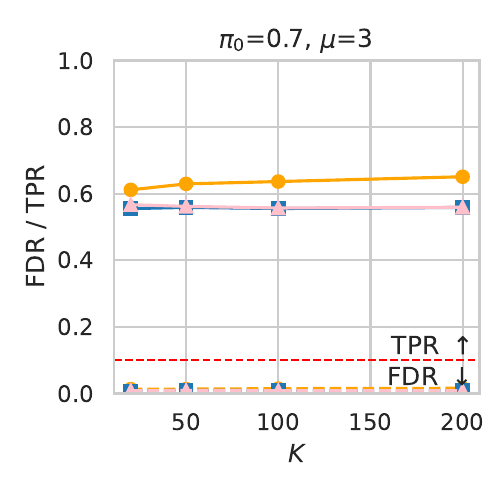}
    \end{subfigure}\begin{subfigure}{0.33\textwidth}
        \includegraphics[trim=0 0 3.2cm 0,clip,width=\textwidth]{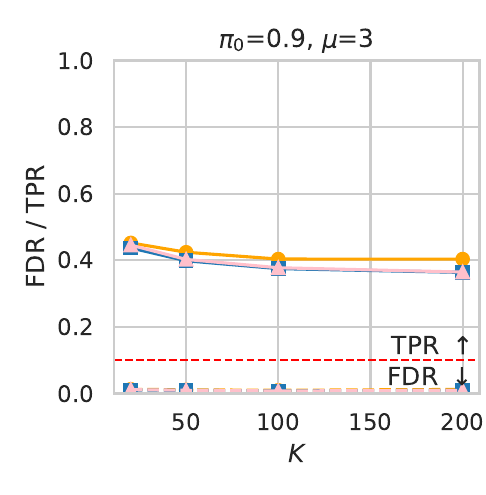}
    \end{subfigure}
    \caption{Plot of FDR and TPR for different null proportions $\pi_0$ and alternative mean $\mu = 3$. All differences are significant; the closed e-BH (orange circle, upper curve in all plots) dominates both e-BH and minimally adaptive e-BH (eBHm; almost identical lower curves in all plots).}
    \label{fig:gaussian-sim}
\end{figure}
\begin{figure}[t]
\centering
    \begin{subfigure}{0.33\textwidth}
        \includegraphics[trim=0 0 3.2cm 0,clip,width=\textwidth]{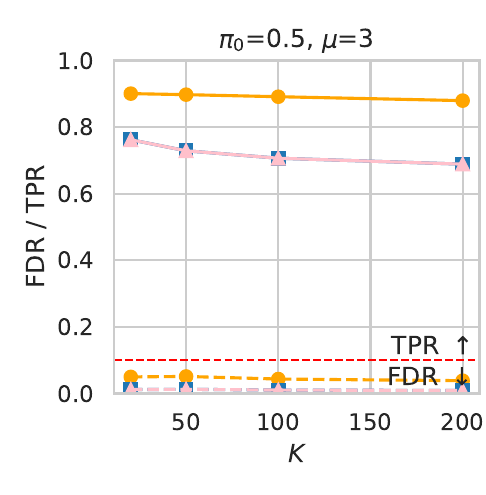}
    \end{subfigure}\begin{subfigure}{0.33\textwidth}
        \includegraphics[trim=0 0 3.2cm 0,clip,width=\textwidth]{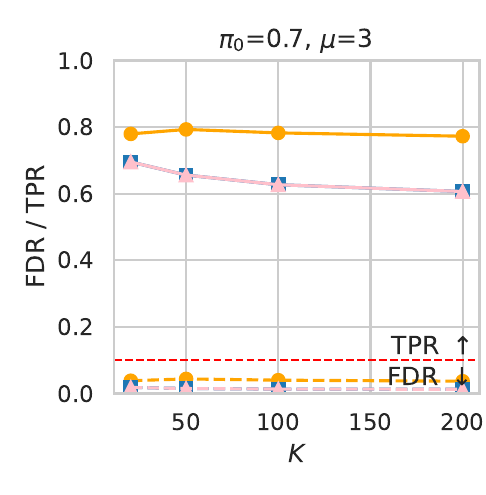}
    \end{subfigure}\begin{subfigure}{0.33\textwidth}
        \includegraphics[trim=0 0 3.2cm 0,clip,width=\textwidth]{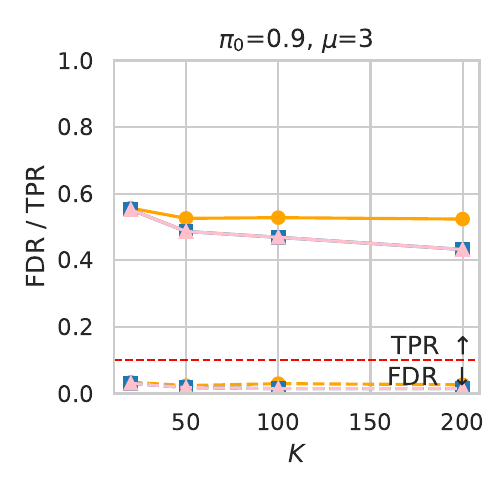}
    \end{subfigure}
    \caption{Plot of FDR and TPR for different null proportions $\pi_0$ and alternative mean $\mu = 3$. All differences are significant; the closed BHY procedure (orange circle, upper curve in all plots) dominates the BHY procedure and its minimally adaptive variant (almost identical lower curves in all plots).}
    \label{fig:by-sim}
\end{figure}

We perform some simple simulations for testing one-sided  Gaussian means. Let $X_k \sim N(\mu_k,1),~ k \in \cK$, and let the $k$-th null hypothesis be given by $H_k: \mu_k \leq 0$ against the alternative that $\mu_k>0$. For our first simulation, we let the data be independent across $k$. Let $\pi_0 \coloneqq |\cN| / K$ be the null proportion, $\mu$ be the signal of all alternatives; so let $X_k \sim \mathrm N(0, 1)$ for $k \in \cN$ and $X_k \sim \mathrm  N(\mu, 1)$ otherwise.
Define $E_k \coloneqq \exp(\lambda X_i - \lambda^2 / 2)$ for $\lambda = \mu =3$, which is the log-optimal choice for $H_k$.
We let $\alpha = 0.1$ in our simulations. We plot the empirical average FDR and true positive rate (TPR) $\coloneqq \E[|\cN^c \cap R| / |\cN^c|]$ of e-BH, minimally adaptive e-BH, and closed e-BH over $1000$ trials. We plot the results in Figure~\ref{fig:gaussian-sim}, and can see that the closed e-BH procedure improves noticeably over the standard e-BH and the minimally adaptive e-BH procedures across all settings, while controlling the FDR below $\alpha$.

We also evaluate the BHY procedure, the application of minimally adaptive BHY procedure (which applies minimally adaptive e-BH to the calibrated p-values from Proposition~\ref{prop:c8-recover-BY}), and the closed BHY procedure from Corollary~\ref{cor:closedby-beats-by}, where we set $P_i \coloneqq 1 - \Phi(X_i)$. Here, we draw $(X_1, \dots, X_K)$ from a multivariate normal with the same component means as described above and with Toeplitz-like covariance matrix where $|\text{Cov}(X_i, X_j)| = \exp(-|i - j| / 10) / 5$ with the covariance being positive if $|i - j|$ is even and negative otherwise. Similarly, our results in Figure~\ref{fig:by-sim} show that our closed BHY procedure has substantial power increase over the original BHY procedure.

Comparing Figures~\ref{fig:gaussian-sim} and~\ref{fig:by-sim} also presents a frequently observed empirical phenomenon: the e-value-based procedures tend to be less powerful than p-value-based procedures. But this observation comes with three caveats: 
\begin{itemize}
\item
 e-value-based procedures offer a much stronger type of post-hoc type-I error control; 
 \item 
 as the BHY example shows, it is sometimes possible to improve p-value-based procedures by employing e-values;
\item
 the empirical observations depend heavily on which e-values are used, especially since it is possible to recover all p-value-based procedures using e-value-based procedures.
\end{itemize}
We remind the reader of the more detailed discussion in Section~\ref{c1:p-vs-e}.

\section*{Bibliographical notes}

\citet{benjamini1995controlling} introduced the FDR error metric and developed the BH procedure for FDR control for independent p-values. \cite{benjamini2001control} extended the result to  p-values under the PRDS condition. That paper also showed that under arbitrary dependence, the FDR control achieved by the BHY procedure equals $\alpha \ell_K$ (the Benjamini--Yekutieli correction), a fact mentioned in Section~\ref{sec:c8-pBH}.

The e-BH procedure was proposed by \cite{wang2022false}, with some simple codes for the e-BH procedure available at 
\begin{quote}
\texttt{https://github.com/ruoduwang/e-BH}.
\end{quote}
The proof of Theorem~\ref{th:c8-boost} is in \cite{wang2022false}, where other boosting methods for e-values than the one in Section~\ref{sec:c8-1-2} are also studied, including the boosting methods under the PRDS condition.

Compound e-values were first used (without being given a name) in~\citet[Theorem 3]{wang2022false}. The term compound e-values first appeared in~\citet{ignatiadis2024values}.
Other authors have used the term generalized e-values~\citep{banerjee2023harnessing, bashari2023derandomized, zhao2024false, lee2024boosting} or relaxed e-values~\citep{ren2024derandomised, gablenz2024catch}. A unified treatment of compound e-values was provided in~\cite{ignatiadis2024compound}, from which this chapter draws many results. The term ``self-consistent procedures'' originated from the literature on multiple testing procedures based on p-values such as  \cite{blanchard2008two}, who also developed the proof technique in Remark~\ref{rem:proof-bh-evalue} (despite not using the terminology of e-values).

The term compound e-values pays tribute   to Robbins' compound decision theory~\citep{robbins1951asymptotically} in which multiple statistical problems are connected through a (compound) loss function that averages over individual losses. 
We refer the reader to e.g.,~\citet{copas1969compound, zhang2003compound, jiang2009general} for more comprehensive accounts and present  a brief overview here. 
As a very concrete example (which was studied in detail already in~\citet{robbins1951asymptotically}), consider the Gaussian sequence model. Empirical Bayes~\citep{robbins1956empirical} ideas work well for this task, for example,~\citet{jiang2009general} provide sharp guarantees for the performance of a nonparametric empirical Bayes method in mimicking the best simple separable estimator in the Gaussian sequence model in~\eqref{eq:gaussian_sequence} as $K \to \infty$.

The combination and derandomization recipe was used by
\citet{ren2024derandomised} to derandomize the model-X knockoff filter~\citep{Candes2018}, a flexible set of methods for variable selection in regression with finite-sample FDR control
 that previously relied on additional randomness. Some further applications include the following: \citet{banerjee2023harnessing} use the same recipe for meta-analysis in which the $\ell$-th study only reports the set of tested hypotheses, the set of discoveries, as well as the targeted FDR level. \citet{li2024note} also discuss designing multiple testing procedures in various contexts by aggregating e-values from different procedures. The current chapter unifies and extends several of these observations.

The minimally adaptive e-BH procedure in Section~\ref{sec:c8-minimal} was proposed by \cite{ignatiadis2024values}, and it was inspired by a similar improvement of the BH procedure in \cite{SG17}.
The stochastic rounding of e-values in Section~\ref{sec:c8-stoch} was developed in \cite{xu2023more}. The intersection of closure and e-values has been tackled in many papers including \cite{Vovk/Wang:2021, vovk2023confidence,fischer2024online,hartog2025family}, and the closed e-BH procedure in Section~\ref{sec:9-closed-ebh} is based on~\cite{xu2025closure} and~\cite{goeman2025epartitioning}. 
Online multiple testing with e-values was studied in~\cite{xu2024online}. Recently,~\cite{fischer2024anonline} proposed an online analog of the BH and e-BH procedures.

 The code to reproduce the randomized e-BH figures is at 
\begin{quote}\texttt{https://github.com/neilzxu/rand-ebh},
 \end{quote} and the code to reproduce the closed e-BH figures is at 
 \begin{quote}
\texttt{https://github.com/neilzxu/closed-ebh}.
 \end{quote}

\part{Advanced Topics}
\label{part:advanced}

\chapter{Approximate and asymptotic e-values}
\label{chap:approximate}

We have been studying e-values and p-values that are valid under the null. In many statistical contexts, one may have approximate or asymptotic p-values and e-values, that satisfy relaxed conditions, and offer relaxed statistical guarantees.  
This chapter formalizes the notions of approximate and asymptotic e-values, compound e-values, and asymptotic compound e-values, show their versatility and usefulness for multiple testing, as well as provide certain optimal constructions and examples. 

\section{Approximate p-values and e-values}

\label{sec:approx}

We first define approximate p-variables and e-variables. The basic idea is that, for these approximate objects, the conditions for p-variables and e-variables in Definition~\ref{def:e-value} should be satisfied approximately. It turns out that   two different forms of approximation error   often exist in statistical applications: 
\begin{enumerate}
    \item  One, represented by $\epsilon$ below,  measures how much the value of the random variable deviates from a p-variable or e-variable.
    For instance, the random variable $(1+\epsilon) P$ for  $P\lawis \mathrm{U}(0,1)$ can be seen as ``$\epsilon$-away'' from a p-variable.
    \item The other one, represented by $\delta$ below measures with how much probability  the random variable deviates  from a p-variable or e-variable.
      For instance, the random variable $ P\id_A $  for  $P\lawis \mathrm{U}[0,1]$  and an event $A$ with probability $\p(A)=1-\delta$ can be seen as ``$\delta$-away'' from a p-variable.
      
\end{enumerate}
This consideration  leads to the following formal definition of approximate p-variables and e-variables, and their connections to the above two forms of approximation error is justified in Proposition~\ref{prop:alter-p}.

\begin{definition}[Approximate p-variables and e-variables]
\label{defi:approx}
Fix the null hypothesis $\cP$ and two functions $\varepsilon:\cP \to \R_+ $ and $\delta:\cP \to [0,1]$, and write $\varepsilon_\p=\varepsilon(\p)$
and $\delta_\p=\delta(\p)$ for $\p \in \cP$. 
\begin{enumerate}[label=(\roman*)]
    \item 
An  $\R_+$-valued  random variable $P$ is an \emph{$(\varepsilon, \delta)$-approximate} p-variable for $\cP$ if 
$$\p(P\le t) \le  (1+\varepsilon_\p )t +\delta_\p \qquad \mbox{for all $t\in (0,1)$ and all $\p\in \mathcal P$.} $$
\item 
A $[0,\infty]$-valued random variable $E$ is an \emph{$(\varepsilon, \delta)$-approximate} e-variable for $\cP$ if  
 $$\E^\p [E\wedge t] \le  1+\varepsilon_\p  +\delta_\p t  \qquad \mbox{for all $t\in \R_+$ and all $\p\in \mathcal P$.}$$ 
\end{enumerate}
\end{definition}

For  two constants $\epsilon \in \R_+$ and $\delta\in [0,1] $, 
  $(\epsilon,\delta)$-approximation is understood as those in Definition~\ref{defi:approx} with the corresponding constant functions. Similarly to ``e-values'' and ``p-values'', we also speak of ``approximate e-values'' and ``approximate p-values'' when we do not emphasize them as random variables.

We make some quick observations. Clearly,  $(0,0)$-approximate p-variables (e-variables) are just p-variables (e-variables).
For a constant $\epsilon\in \R_+$, $P$ is an $(\epsilon,0)$-approximate p-variable if and only if $(1+\epsilon)P$ is a p-variable,
and 
$E$ is an $(\epsilon,0)$-approximate e-variable if and only if $E/(1+\epsilon) $ is an e-variable. 
The $\delta$-parameter is more lenient than the $\epsilon$-parameter; for instance, $(\epsilon,\delta)$-approximation implies $(0,\epsilon+\delta)$-approximation in both classes (see Proposition~\ref{prop:varepsilon_to_delta} below for a sharper result).
If $\delta=1$, then $(\epsilon,\delta)$-approximation puts no requirements on $P$ or $E$.

As mentioned above, the definitions of approximate p-variables and approximate e-variables can be interpreted as allowing for a  wiggle room  in the definitions of p-variables and e-variables in two ways: (a) the probability and expectation requirements   have a multiplicative factor $1+\epsilon_\p$; (b) the requirements are only required to hold on an event of probability $1-\delta_\p$. 
This will be made  transparent  through the following proposition. 
 
\begin{proposition}
\label{prop:alter-p}
Let $\varepsilon:\cP \to \R_+ $ and $\delta:\cP \to [0,1]$ and $P$ and $E$ be nonnegative random variables ($[0,\infty]$-valued for $E$). For the following  statements,
\begin{enumerate} 
    \item[(P1)] $P$ is an $(\varepsilon, \delta)$-approximate p-variable for $\cP$;
\item[(P2)]  
 for every $\p \in \cP$,  there exists an event $A$ such that
$$\p(P\le t , A) \le  (1+\varepsilon_\p )t ~\mbox{for all $t\in (0,1)$ ~~and~~ }\p(A) \geq 1-\delta_\p;$$
  \item[(E1)]   $E$ is an $(\varepsilon, \delta)$-approximate e-variable for $\cP$;
\item[(E2)]   for every $\p \in \cP$, 
there exists an event $A$ such that
$$ \E^{\p}[E\id_{A}] \leq 1+\varepsilon_\p \mbox{  ~~and~~ }\p(A) \geq 1-\delta_\p;$$
\end{enumerate}
we have that (P2) implies (P1) and  (E2) implies (E1).
Moreover, if  each $\p\in \cP$ is atomless, then (P1)--(P2)
  are equivalent,
  and (E1)--(E2)
  are equivalent.
\end{proposition}

\begin{proof}[Proof.]
When $\p$ is atomless, for any random variable $X$, there exists 
  an event $A$ with $\p(A)=1-\delta_\p$ such that $\{X < t \} \subseteq  A\subseteq  \{X\le  t \}$ for some $t \in \R $; a formal statement is presented in Lemma~\ref{lem:atomless-p}.
  This fact will be used in the proof below. 
  If $\delta_\p=1$ there is nothing to prove, so we assume $\delta_\p<1$ below.

(P2) $\Rightarrow$ (P1): It suffices to notice  $$\p(P\le t) \le \p(P\le t,A)+\p(A^c) \le (1+\epsilon_\p) t +\delta_\p$$ for $\p\in \cP$ and $t\in [0,1]$. 

(P1) $\Rightarrow$ (P2) in case each $\p$ is atomless:  
Let $A$ be an event with $\p(A)=1-\delta_\p$ such that $\{P > s \} \subseteq  A\subseteq  \{P\ge  s \}$ for some $s \in [0,1]$.   
Note that $\p(P\le t,A^c)=\min\{\p(P\le t),\p(A^c)\}$.
Therefore, 
$$\p(P\le t, A) = 
( \p(P\le t) - \p(A^c))_+ \le ((1+\varepsilon_\p )t +\delta_\p-\delta_\p)_+ = (1+\varepsilon_\p)t$$
for all $t\in (0,1)$,
as desired.

(E2) $\Rightarrow$ (E1): It suffices to notice  
\begin{align}
\E^\p[E\wedge t] =
\E^\p[(E\wedge t) \id_A] + 
\E^\p[(E\wedge t) \id_{A^c}]
&\le \E^\p[E \id_A] + 
\E^\p[t \id_{A^c}] \\
&\le 
 (1+\varepsilon_\p ) +\delta_\p t 
\end{align}
 for $\p\in \cP$ and $t\in \R_+$. 

(E1) $\Rightarrow$ (E2) in case each $\p$ is atomless:  
Let $A$ be an event with $\p(A)=1-\delta_\p$ such that $\{E < t \} \subseteq  A\subseteq  \{E\le  t \}$ for some $t \in \R_+$. 
Note that $E\ge t$ on $A^c$.
We have, 
\begin{align}
\E^\p[E \id_A]   
= \E^\p[(E\wedge t)\id_A]
&=  \E^\p[ E\wedge t ]
-\E^\p[ (E\wedge t) \id_{A^c}]\\
&\leq 1+\epsilon_\p + t \delta_\p
-\E^\p[ t\id_{A^c}]
= 1+\epsilon_\p ,
\end{align}
as desired. 
\end{proof}

Without the atomless assumption on each $\p\in \cP$, the implications (P1)$\Rightarrow$(P2)
and (E1)$\Rightarrow$(E2)
in  Proposition~\ref{prop:alter-p}
do not hold, because the event $A$ with desired probability may not exist (counterexamples can be found on discrete spaces). We emphasize that the atomless assumption does not pertain to the distribution of $P$ or $E$ (i.e., $P$ or $E$ can be   discrete). For instance, this atomless assumption is automatic if external randomization is permitted.

The two formulations (P1) and (P2) in the above proposition are obviously equivalent when $\delta=0$ without requiring any condition; the same holds for (E1) and (E2).

In Chapter~\ref{chap:markov}, we described calibrators from p-values to e-values and those from e-values to p-values. Recall from Chapter~\ref{chap:markov}  that a (p-to-e) calibrator is a decreasing function $f:[0,\infty)\to [0,\infty]$
satisfying $\int_0^1 f(x)\d x \le 1$ and $f(x)=0$ for $x\in (1,\infty)$.
An e-to-p calibrator is a decreasing function $g:[0,\infty]\to [0,\infty)$
such that $g(x)\ge (1/x)\wedge 1$ for all $x\in [0,\infty]$, because the only admissible e-to-p calibrator is $x\mapsto (1/x)\wedge 1$.  
The next theorem demonstrates that we can calibrate approximate p-values and approximate e-values in exactly the same way as in calibrating p-values and e-values. 
\begin{theorem}
\label{th:calibration}
 Let $\varepsilon:\cP \to \R_+ $ and $\delta:\cP \to [0,1]$. 
Using any calibrators, 
    calibrating an $(\varepsilon, \delta)$-approximate  p-variable  yields an $(\varepsilon, \delta)$-approximate e-variable, and vice versa.
\end{theorem}

\begin{proof}[Proof.]
Fix any $\p\in \cP$.   If $\delta_\p=1$ there is nothing to prove, so we assume $\delta_\p<1$ below. 
For the first statement, 
let $P$ be an $(\varepsilon, \delta)$-approximate  p-variable   for $\p$ and $f$ be a (p-to-e) calibrator and $b=(1+\epsilon_\p)/(1-\delta_\p) \ge 1$. 
Let $F$ be a cdf on $[0,1/b]$ given by
$F(s)=(1+\epsilon_\p)s + \delta_\p$ for $s\in [0,1/b]$.
We have $\p(P\le s   ) \le  F(s) $ for $s \in [0,1]$ and $F(0)=\delta_\p$.
Hence, noting that  the calibrator  $f$ is decreasing, we have, for all $t\in \R_+$, 
\begin{align*}
    \E^\p[f(P)\wedge t ] & 
    \le  \int_0^{1/b}( f(s)\wedge t) \d F(s)     
    \\ &\le  (1+\epsilon_\p) \int_0^{1/b}   f(s)\d s + t F(0) 
\\&\le  (1+\epsilon_\p) \int_0^1  f(s)\d s  +  \delta_\p t \le  1+\epsilon_\p + \delta_\p t. 
\end{align*}
For the second statement, it suffices to consider the calibrator $x\mapsto (1/x)\wedge 1$, and we define $P = (1/E)\wedge 1$. 
We have, for  $t\in (0,1)$,  
\begin{align*}
\p(P \leq t) 
&=\p(E\ge 1/t) 
\\ &=  \p(E\wedge(1/t) \ge 1/t) 
\\ &\le  t \E^\p[E\wedge(1/t) ] \leq t \left( 1+\epsilon_\p   +    \frac{\delta_\p}{t} \right)= (1+\epsilon_\p)  t +  \delta_\p,
\end{align*} 
as desired.
\end{proof}

\begin{remark}
Following a similar proof, the statements in Theorem~\ref{th:calibration} hold true also if we use the alternative formulations (P2) and (E2) in Proposition~\ref{prop:alter-p} for both approximate p-variables and approximate e-variables.
 \end{remark}

 The roles of $\varepsilon$ and $\delta$ are asymmetric and, as mentioned earlier,  the wiggle room created by  $\delta$ is more lenient. The following result shows that an $(\varepsilon, \delta)$-approximate p-variable (or e-variable) with $\varepsilon >0$ is always also an $(0, \delta')$-approximate p-variable (or e-variable) where $\delta'$ can be chosen strictly smaller than $\varepsilon + \delta$.

\begin{proposition}
\label{prop:varepsilon_to_delta}
Let  $\varepsilon:\cP \to \R_+ $, $\delta:\cP \to [0,1]$, $\delta' = (\varepsilon + \delta)/(1+\varepsilon)$.
\begin{enumerate}[label=(\roman*)]
\item An ($\varepsilon, \delta$)-approximate p-variable for $\mathcal{P}$ is also an $(0, \delta')$-approximate p-variable for $\mathcal{P}$.
\item An ($\varepsilon, \delta$)-approximate e-variable for $\mathcal{P}$ is also an $(0, \delta')$-approximate e-variable for $\mathcal{P}$.
\end{enumerate}
\end{proposition}

\begin{proof}[Proof.]
Let $P$ be an ($\varepsilon, \delta$)-approximate p-variable. Fix $\p \in \cP$. Then, for any $t \in (0, (1-\delta_{\p})/(1+\varepsilon_{\p})),$ we have that:
$$ \p(P \leq t) \leq (1+\varepsilon_{\p})t + \delta \leq t + \varepsilon_{\p}\frac{1-\delta_{\p}}{1+\varepsilon_{\p}} + \delta_{\p} = t + \delta'_{\p}.$$
Meanwhile the inequality $\p(P \leq t) \leq t + \delta'_{\p}$ is clearly true for $t \geq (1-\delta_{\p})/(1+\varepsilon_{\p})$ (since the right-hand side is $\geq 1$). Thus $P$ is an $(0, \delta')$-approximate p-variable for $\cP$.

The argument for the approximate e-variable is analogous (noting that only $t\ge 1$ needs to be checked in the condition for approximate e-variables) and omitted. 
\end{proof}

Proposition~\ref{prop:varepsilon_to_delta} also suggests that controlling $\delta$ is a stronger requirement, leading to two different notions of asymptotic e-values and p-values in Section~\ref{sec:asymp-e}.

\section{Asymptotic p-values and e-values}
\label{sec:asymp-e}

\subsection{Definitions and conditions}\label{sec:c9-conditions}

We consider a sequence of datasets $X=X^{(n)}$ of growing size as $n \to \infty$. Based on $X^{(n)}$ we compute a $[0,1]$-valued statistic $P^{(n)} = P^{(n)}(X^{(n)})$.  
Similarly, we can compute a nonnegative extended random variable $E^{(n)}= E^{(n)}(X^{(n)})$. The following definitions give some intuitive conditions on $P^{(n)}$ and $E^{(n)}$ for them to be considered as p-variables and e-variables in an asymptotic sense.

\begin{definition}[Asymptotic p-variables and e-variables]
\label{defi:asymp_evariables_pvariables}
A sequence $(E^{(n)})_{n\in\N}$ of $[0,\infty]$-valued random variables is called a sequence of:
\begin{enumerate}[label=(\roman*)]
    \item  \emph{asymptotic} e-variables for $\mathcal P$,  if  for each $n \in \N$, $E^{(n)}$ is $(\varepsilon_n, \delta_n)$-approximate for $\mathcal P$ 
    for some $(\varepsilon_n, \delta_n)$, satisfying  
    \begin{equation}\label{eq:asymp-e} \mbox{for each $\p\in \mathcal P$,}\quad  \lim_{n\to \infty}  \varepsilon_n(\p) = \lim_{n\to \infty} \delta_n(\p) = 0;
    \end{equation}
    
    \item  \emph{strongly asymptotic} e-variables for $\mathcal P$, if \eqref{eq:asymp-e} is replaced by   
$$    \mbox{for each $\p\in \mathcal P$,}  \lim_{n\to \infty}   \varepsilon_n(\p) = 0 \mbox{~and~} 
 \delta_n(\p) = 0 \mbox{~for $n$ large enough;}$$ 
 
    \item \emph{uniformly asymptotic} e-variables for $\mathcal P$, if 
     \eqref{eq:asymp-e} is replaced by  $$ 
    \lim_{n\to \infty} \sup_{\p \in \cP} \varepsilon_n(\p) = \lim_{n\to \infty} \sup_{\p \in \cP} \delta_n(\p) = 0;$$     
    \item  \emph{uniformly strongly asymptotic} e-variables for $\mathcal P$, if condition \eqref{eq:asymp-e} is replaced by $$\lim_{n\to \infty} \sup_{\p \in \cP} \varepsilon_n(\p) = 0 \quad \mbox{~and~} \quad  \sup_{\p \in \cP} 
 \delta_n(\p) = 0 \mbox{~for $n$ large enough;}$$

\end{enumerate}
 Sequences $(P^{(n)})_{n\in\N}$ of  (uniformly, strongly) asymptotic p-variables for $\mathcal P$ are defined analogously.

\end{definition}

We sometimes   simply say that a nonnegative sequence is asymptotic for $\cP$, and it should be clear from context whether we mean  a sequence of asymptotic e-variables (allowed to be $[0,\infty]$-valued) or that of asymptotic p-variables  ($\R_+$-valued). 
Clearly, if $\mathcal P$ is finite, then    asymptotic e-variables or p-variables are also uniformly asymptotic. 


Several conditions for asymptotic p-values and e-values are  often encountered in typical statistical contexts. 
We summarize them together in the following proposition.

\begin{proposition}
\label{prop:asym-e}
\begin{enumerate}[label=(\roman*)]
      \item     If under each $\p\in\cP$, a nonnegative sequence  
converges to an e-variable (resp.~a p-variable) in distribution, then it
is a sequence of asymptotic e-variables  (resp.~p-variables) for~$\cP$.  
\item The nonnegative sequence  $(E^{(n)})_{n\in\N}$ is uniformly strongly asymptotic  for $\mathcal P$ if and only if $$\limsup_{n \to \infty} \sup_{\p\in \cP} \E^\p[E^{(n)}]\le 1 .$$ Similarly, it is   strongly asymptotic  for $\mathcal{P}$ if and only if $$\limsup_{n \to \infty} \E^\p[E^{(n)}]\le 1 \mbox{ for all $\p\in  \cP$}.$$
    \item If a sequence of asymptotic e-variables is  uniformly integrable for each $\p\in \cP$, then it is  strongly asymptotic.

\item The nonnegative sequence $(P^{(n)})_{n\in\N}$  
is uniformly asymptotic  for $\mathcal P$ if and only if  
$$\limsup_{n\to\infty}\sup_{\p\in \cP} \p(P^{(n)}\le t)\le t $$ for all $t\in (0,1)$.
Similarly, it
is asymptotic  for $\mathcal P$ if and only if  
$\limsup_{n\to\infty}\p(P^{(n)}\le t)\le t $ for all $\p\in \mathcal P$ and $t\in (0,1)$.

\end{enumerate}

\end{proposition}
\begin{proof}[Proof.]
\begin{enumerate}[label=(\roman*)]

\item  Fix $\p\in \cP$. 
Suppose that $(P^{(n)})_{n\in\N}$ converges in distribution to a p-variable $P$. Note that 
$$ \limsup_{n\to\infty} \p(P_n\le t) \le \p(P\le t) \le t$$ for each $t\in [0,1]$.
This implies 
$ \lim_{n\to\infty} \p(P_n\le t)\vee t =t$  for each $t\in [0,1]$. 
Since point-wise convergence of increasing functions to a continuous function on a compact interval implies uniform convergence,  we have 
\begin{align}
& \lim_{n\to\infty} \sup_{t\in [0,1]} (\p(P_n\le t)\vee t-t) = 0 \\
& \implies \lim_{n\to\infty} \sup_{t\in [0,1]}(\p(P_n\le t)-t)_+ = 0. 
\end{align}
This shows that $(P^{(n)})_{n\in\N}$ 
 is a sequence of asymptotic p-variables.

Next we show the statement on asymptotic e-variables. For $\delta \in [0,1]$,
and a nonnegative random variable $E$, 
consider the function
$$f_\delta(E)=\sup_{t\in \R_+} (\E^\p [E\wedge t]  -\delta  t) -1. $$ 
By definition, if $\epsilon :=
 f _\delta (E)
\vee 0 $ is finite, then $E$ is an $(\epsilon ,\delta )$-approximate e-variable for $\{\p\}$. 
Let $E$ be the limit of the nonnegative sequence $(E^{(n)})_{n\in\N}$ in distribution, and take $m\in \R_+$ be such that $\p(E\ge m)<\delta$. 
Since $E^{(n)}\to E$ in distribution, there exists $N\in \N$ such that $\p(E^{(n)}\ge m)<\delta$ for $n>N$ by the portmanteau theorem.
It follows that, for all $t\in \R_+$, 
\begin{align*}
\E^\p[E^{(n)}\wedge t] - \E^\p[E^{(n)}\wedge m] - \delta t 
&\le  
\p(E^{(n)} \ge m) (t-m)_+ - \delta t 
\\
&\le \delta (t-m)_+-\delta t
\le 0,
\end{align*}
and hence  $f_\delta(E^{(n)})\le \E^\p[E^{(n)}\wedge m]-1\to \E^\p[E\wedge m]-1\le 0$ as $n\to\infty$.  
This implies, in particular, that for any $\delta>0$,  
there exists $N'\in \N$ such that for $n>N'$, $E^{(n)}$ is an $(\delta,\delta)$-approximate e-variable. 
Since $\p\in \cP$ is arbitrary, this shows that $(E^{(n)})_{n\in\N}$ is a sequence of
 asymptotic e-variables for $\cP$.  
 
    \item  We prove this result for the non-uniform case. The other argument is analogous. First suppose that $(E^{(n)})_{n \in \N}$ is a sequence of strongly   asymptotic e-variables for $\mathcal{P}$. Fix $\p \in \cP$. Then:
$$ \E^\p[E^{(n)}] \leq 1+ \varepsilon_n(\p) \to 1 \text{ as } n \to \infty,$$
where the first inequality holds for large enough $n$.
For the other direction, note that any nonnegative extended random variable $E^{(n)}$ is an $(\varepsilon_n,0)$-approximate e-variable for $\mathcal{P}$ with the choice $\varepsilon_n(\p) = (\E^\p[E^{(n)}]-1)_+$ as long as $\E^\p[E^{(n)}]$ is finite.
If  $\limsup_{n \to \infty} \E^\p[E^{(n)}]\le 1$, this implies that $\E^\p[E^{(n)}] < \infty$ for $n$ large enough and the above $\varepsilon_n(\p)$ satisfies $\lim_{n \to \infty} \varepsilon_n(\p) = 0$, that is, $(E^{(n)})_{n \in \N}$ is a sequence of strongly   asymptotic e-variables.

\item
Fix $\p\in \cP$. Note that for $t\in \R_+$,
$$
\E^{\p}[E^{(n)}] \leq \E^{\p}[E^{(n)}\wedge t] + \E^{\p}[E^{(n)}\id_{\{E^{(n)} \geq t\}}].
$$
By uniform integrability, for any $\eta>0$ there exists $T>0$ such that 
$$\sup_{n \in \mathbb N}\E^{\p}\left[E^{(n)}\id_{\{E^{(n)} \geq T\}}\right]  \leq \eta.$$
Thus,
$$
\E^{\p}[E^{(n)}] \leq 1 + \varepsilon_n(\p) + \delta_n(\p) T + \eta,
$$
and hence,
$$ \limsup_{n \to \infty} \E^{\p}[E^{(n)}]  \leq 1+\eta.$$
Since $\eta>0$ is arbitrary, we conclude by using part (ii).

\item The proof is similar to that of part (i) and omitted here. \qedhere

\end{enumerate}

\end{proof}

In Proposition~\ref{prop:asym-e}, 
part (i)  gives a sufficient condition that best describes the intuition behind  asymptotic e-variables and p-variables, via convergence in distribution.
Converging to an e-variable or p-variable in distribution is a sufficient condition for a sequence to be asymptotic e-variables or asymptotic p-variables, but this condition is not necessary.  For instance, any sequence that is  smaller than an e-variable  is a sequence of asymptotic e-variables, and it does not need to converge. 
Part (ii) gives 
an intuitive condition on the expectation for strongly asymptotic e-variables.   
Part (iii) connects asymptotic e-variables and strongly asymptotic e-variables 
via uniform integrability, and part (iv) gives equivalent conditions for uniformly (but not strongly) asymptotic p-variables.

We can naturally calibrate asymptotic p-variables and e-variables by applying a (possibly different) calibrator to each element in the sequence, in the sense of Theorem~\ref{th:calibration}.

\begin{proposition}
\label{prop:asymptotic_calibration}
    Calibrating a sequence of (uniformly, strongly) asymptotic p-variables  yields a sequence of (uniformly, strongly) asymptotic e-variables, and vice versa. 
\end{proposition} 
The above result follows directly from Theorem~\ref{th:calibration}. The most common way of constructing asymptotic p-variables and e-variables is via limit theorems such as the central limit theorem. 
\begin{proposition}
\label{proposition:weak_convergence_to_asymptotic}
Let $Z^{(n)}$ be a statistic computed based on $X^{(n)}$ for $n\in \N$. Suppose that
$$
Z^{(n)} \to Z \text{~in distribution as }\, n \to \infty \,\text{ under every }\, \p \in \cP,
$$
where $Z $ is a random variable with a continuous cumulative distribution function $F$.  
\begin{enumerate}[label=(\roman*)]
\item Let $P^{(n)} = F(Z^{(n)})$ or $P^{(n)}=1-F(Z^{(n)})$ for $n\in \N$. Then $(P^{(n)})_{n\in\N}$ is a sequence of asymptotic p-variables for $\cP$.
\item Let $h: \mathbb R \to [0,\infty]$ be a function with $\E[h(Z)] \leq 1$ that is increasing, decreasing, or upper semicontinuous, and  $E^{(n)} = h(Z^{(n)})$ for $n\in \N$. Then $(E^{(n)})_{n\in\N}$ is a sequence of asymptotic e-variables for $\cP$.  
If $h$ is further bounded, then $(E^{(n)})_{n\in\N}$ is strongly asymptotic  for $\cP$.  
\end{enumerate}
\end{proposition}

\begin{proof}[Proof.]
By the continuous mapping theorem, $P^{(n)}$ converges weakly to the uniform distribution, and hence item (i) directly follows from Proposition~\ref{prop:asym-e}, part (i). 
 
For (ii) we may argue as follows. 
If  $h$ would have been continuous, then the result would follow 
from the continuous mapping theorem and  Proposition~\ref{prop:asym-e}, part (i).
For upper semicontinuity,
note that $\limsup_{n\to\infty} \p(E^{(n)}\ge  x) \le \p(h(Z)\ge x)$ for all $x\in \R$ by convergence in distribution. This gives 
$\limsup_{n\to \infty} f_\delta(E^{(n)}) \le f_\delta(h(Z))\le 0 $ using the notation and argument in the proof of Proposition~\ref{prop:asym-e}, part (i). 
The last statement follows from 
Proposition~\ref{prop:asym-e}, part (iii).
When $h$ is increasing or decreasing, replacing it with its upper semicontinuous version, which does not change $\E[h(Z)]$, yields the desired statements. 
\end{proof}


Some common cases where we can apply Proposition~\ref{proposition:weak_convergence_to_asymptotic} include: 
\begin{enumerate}
\item $Z \lawis  \mathrm{N}(0,1)$ as would follow from the central limit theorem.
\item $Z \lawis  |\mathrm{N}(0,1)|$, where  $ |\mathrm{N}(0,1)|$ is the folded normal distribution (i.e., the distribution of $|X|$ where $X\lawis \mathrm{N}(0,1)$).  
\item $Z \lawis  \chi^2_{\nu}$, where $\chi^2_{\nu}$ denotes the chi-square distribution with $\nu$ degrees of freedom. Such a limit occurs e.g., for likelihood ratio or chi-square tests.
\end{enumerate}

\section{An asymptotic e-value via the central limit theorem}
\label{sec:c9-ex-asym}

We next provide a prototypical example of asymptotic e-variables 
along the lines of 
the central limit theorem.

\begin{example}[Testing the mean with unknown variance]
    \label{ex:asym-variance}

Suppose we observe the dataset $X^{(n)} = (X^1,\dots,X^n)$ where $X^i \simiid \p$ are real-valued. We assume that their variance, $\var^{\p}(X^i)$, is positive and finite, but is unknown. We seek to test the hypothesis 
\begin{equation}
\label{eq:asym-variance-p}\cP=\{\p\in \cM_1\,:\, \E^{\p}[X^i] = 0,\, \var^{\p}({X^i}) \in (0,\infty) \},
\end{equation}
where $\cM_1$ denotes the set of all probability measures on $\mathbb R$. 
Define 
\begin{equation}
\begin{aligned}
 &\bar{X}^{(n)} = \frac{1}{n}\sum_{i=1}^n X^i,\;\;\;\, S^{(n)} = \sqrt{\frac{1}{n}\sum_{i=1}^n (X^i)^2},  \\
 &\hat{\sigma}^{(n)} = \sqrt{\frac{1}{n-1}\sum_{i=1}^n \left(X^i - \bar{X}^{(n)} \right)^2},
\end{aligned}
\label{eq:one_sample_summary}
\end{equation}
and for fixed $\lambda \in \R$, let
\begin{equation}
\label{eq:asymptotic_evalue}
E^{(n)}= \exp\left(\lambda \frac{ \sqrt{n} \bar{X}^{(n)}}{S^{(n)}} - \frac{\lambda^2}{2}\right),
\end{equation}
and
\begin{equation}
\label{eq:asymptotic_evalue2}
\tilde{E}^{(n)}= \frac12 \exp\left(\lambda \frac{\sqrt{n}\bar{X}^{(n)}}{S^{(n)}}- \frac{\lambda^2}{2}\right) + \frac12 \exp\left(-\lambda \frac{\sqrt{n}\bar{X}^{(n)}}{S^{(n)}} - \frac{\lambda^2}{2}\right).
\end{equation}
The next theorem justifies that \eqref{eq:asymptotic_evalue} defines a sequence of strongly asymptotic e-variables. 
\end{example}

\begin{theorem}
\label{prop:c9-asym}
The sequences of random variables $E^{(n)}$ in \eqref{eq:asymptotic_evalue}
and $\tilde{E}^{(n)}$ in \eqref{eq:asymptotic_evalue2} 
are strongly asymptotic e-variables for $\cP$ in \eqref{eq:asym-variance-p}. 
If $\hat{\sigma}^{(n)}$ is used in place of $S^{(n)}$ in \eqref{eq:asymptotic_evalue} and \eqref{eq:asymptotic_evalue2}, then $E^{(n)}$ and $\tilde{E}^{(n)}$ are asymptotic e-variables for $\cP$.
In either case, for $\lambda>0$, $E^{(n)}$ (resp.~$\tilde{E}^{(n)}$) grows to infinity with probability $1$ under the alternative that $X^1,\dots,X^n$ are drawn iid from $\mathbb Q$ with $\E^{\mathbb Q}[X^i]> 
0$ (resp.~$\E^{\mathbb Q}[X^i]\ne 
0$) and $\var^{\mathbb Q}(X^i) < \infty$. 
\end{theorem}

\begin{proof}[Proof.] 
Take any $\p \in \cP$. By the central limit theorem, the law of large numbers and Slutsky's theorem, $\sqrt{n}\bar{X}^{(n)}/S^{(n)}$ converges in distribution to $\mathrm{N}(0,1)$ as $n \to \infty$. 
We may verify that
$$ \mathbb E^{\mathbb P} \left[\exp\left(\lambda Z - \frac{\lambda^2}{2}\right)\right] =1,$$
for $Z \sim \mathrm{N}(0,1)$. 
By Proposition~\ref{proposition:weak_convergence_to_asymptotic} it follows that $E^{(n)}$ are asymptotic e-variables. To upgrade this result to strongly asymptotic, we can
apply Theorem 2.5 of~\citet{gine1997when} to show that the sequence $(E^{(n)})_{n \in \N}$ is uniformly integrable and then conclude the proof by Proposition~\ref{prop:asym-e}, part (iii).

Next take an alternative $\q$ with $\E^{\mathbb Q}[X^i]> 0$, $\var^{\mathbb Q}(X^i) < \infty$, and suppose that $\lambda >0 $. Observe the following. By the strong law of large numbers
$$ 
\bar{X}^{(n)} \to \E^{\mathbb Q}[X^i]>0\;\; \text{ almost surely as } \;\; n \to \infty,
$$
and
$$
S^{(n)} \to \sqrt{  \var^{\mathbb Q}({X^i}) + \E^{\mathbb Q}[X^i]^2 }>0  \;\; \text{ almost surely as } \;\; n \to \infty.
$$
The above imply that
$$
\lambda \frac{\sqrt{n} \bar{X}^{(n)}}{S^{(n)}} \to \infty\;\; \text{ almost surely as } \;\; n \to \infty,
$$
which also means that $E^{(n)}$ goes to $\infty$ with probability 1.
The argument for $\tilde{E}^{(n)}$ is analogous and omitted. The argument for $\hat{\sigma}^{(n)}$ follows from Proposition~\ref{proposition:weak_convergence_to_asymptotic}.
\end{proof}

The above conclusion extends to a broader class of asymptotic e-variables  obtained by mixing in~\eqref{eq:asymptotic_evalue} over $\lambda$ with respect to a distribution with bounded support. Concretely, let $G$ be any distribution on $\mathbb R$ with compact support. Then, under the assumptions of this section,
$$
\int \exp\left(\lambda \frac{\sqrt{n}\bar{X}^{(n)}}{S^{(n)}} - \frac{\lambda^2}{2}\right) \d G(\lambda),
$$
is a strongly asymptotic e-variable. 

One central purpose of defining asymptotic e-variables is to handle cases where one does not know how to construct a nonasymptotic e-variable. Indeed, for the above $\cP$, we do not know of nonasymptotic e-variables (other than constants). Had a bound on the variance been known, such nonasymptotic e-variables do exist, as shown next.

\begin{example}[Testing the mean with bounded variance]
\label{ex:c2-mv}
Given some $\sigma>0$, consider the null hypothesis
\[
\cP_\sigma=\left\{\p\in \cM_1: \E^{\p}[X] = 0,\; \var^\p({X}) \le  \sigma^2 \right\}.
\]
Then $X^2/\sigma^2$ is an e-variable for $\cP_\sigma$ (see Section~\ref{sec:c2-bounded}), and for any $\lambda \in \R$ the following are e-variables for $\cP_\sigma$:
\[
\exp\left(\lambda X - \frac{\lambda^2}{6}(X^2 + 2 \sigma^2)\right)
\]
and
\[
\exp\left(\phi(\lambda X) - \frac{\lambda^2}{2}\sigma^2 \right),
\]
where $\phi(x) = \log(1+x+x^2/2)$ for $x \geq 0$ and $\phi(x) = -\log(1-x+x^2/2)$ for $x<0$.
The different choices of e-variables have large e-power against different alternatives.
\end{example}  

\section{Approximate/asymptotic compound p-values and e-values}
\label{sec:app-asym-com}

In this section we extend the setting of Sections~\ref{sec:approx} and~\ref{sec:asymp-e} to the multiple hypothesis testing setting, introduced in Chapter~\ref{chap:compound}.

\subsection{Compound  p-values and e-values}
\label{sec:further_defi}


 Compound e-variables are defined in Definition~\ref{defi:compound_evalues}, and they satisfy 
 $$\sum_{k : \p \in \cP_k} \E^{\p}[E_k] \leq K \qquad \mbox{for all $\p\in \mathcal P$.}$$
 We next define their counterpart, compound p-variables. 



\begin{definition}[Compound p-variables]
\label{defi:compound_pvalues}
Fix the null hypotheses $(\cP_1,\dots,\cP_K)$ and a set $\cP $ of distributions.  Let $P_1,\dots,P_K$ be nonnegative random variables. We say that 
$P_1,\dots,P_K$ are \emph{compound} p-variables for $(\cP_1,\dots,\cP_K)$ under $\cP$ if 
$$
\sum_{k : \p \in \cP_k} \p(P_k\leq t) \leq Kt 
$$ 
for all $t\in (0,1)$ and all $\p\in \mathcal P$.
\end{definition}


Similar to Theorem~\ref{th:calibration}, compound p-variables 
and  compound e-variables
can be transformed  via calibration.

\begin{theorem}
\label{prop:compound-pe1}
Let $h $ be any (p-to-e) calibrator.
  If  $P_1,\dots,P_K$ are compound p-variables for $(\cP_1,\dots,\cP_K)$ under $\cP$,   then
$ 
E_k =   h  ( P_k  ),\, k \in \mathcal{K},
$ 
are compound e-variables for $(\cP_1,\dots,\cP_K)$ under $\cP$.
Similarly, calibration in the e-to-p direction is also valid.  
\end{theorem}
\begin{proof}[Proof.]
Fix an atomless $\p $ and it suffices to verify the statement for this $\p$. Let $K_0 = |\mathcal{N}(\p)|$ and take 
$L \lawis \mathrm{U}(\mathcal{N}(\p))$ and 
$U\lawis \mathrm{U}[0,1]$, 
mutually independent and independent of everything else. Then $$\tilde P:= P_{L} \id_{\{U\le K_0/K\}} +  \id_{\{ U>K_0/K\} }$$ is a bona fide p-value because for $t\in (0,1)$,
$$ \p( \tilde P  \leq t) =  \p(P_{L} \leq t) \p( U\leq K_0/K)  = \frac{K_0}{K} \frac{1}{K_0}\sum_{k: \p \in \mathcal{P}_k} \p(P_k \leq  t) \leq t.$$
Therefore, 
\begin{align*} \frac{1}{K} \sum_{k: \p \in \mathcal{P}_k }\E^{\p}\left[   h( P_k) \right]& = \frac{K_0}{K} 
 \E^{\p}\left[   h( P_{L} ) \right] 
\\& \le   \E^{\p}\left[   h( P_L) \id_{\{U\le K_0/K\}}+ h(1)  \id_{\{ U>K_0/K\} }  \right]  \\& =  \E^{\p}\left[ h(\tilde P) \right] \leq 1,
\end{align*}
where the last inequality holds since $h$ is a calibrator.
The other direction of calibration can be checked directly.
\end{proof}

Similarly to Example~\ref{ex:c8-weighted}, 
an example of compound p-variables   is given by weighted p-variables.

\begin{example}[Weighted p-values]
    \label{ex:c9-weighted}
    Let $w_1,\ldots,w_K \geq 0$ be deterministic nonnegative numbers such that $\sum_{k \in \mathcal{K}} w_k \leq K$. 
 Let $P_1^\prime,\dots,P_K^\prime$ be p-variables and define $\tilde{P}_k = P_k^\prime/w_k$ for $k \in \mathcal{K}$ with the convention $0/0 = 0$.  Then, $\tilde{P}_1,\dots,\tilde{P}_K$ are compound p-variables. 
\end{example}

\begin{example}[Combining p-variables and compound e-variables]
Suppose that $P_k$ is a p-variable for $\cP_k$ for all $k \in \mathcal{K}$; $P_1,\dots,P_K$ can be arbitrarily dependent. Moreover, let $E_1,\dotsc,E_K$ be compound e-variables for $(\cP_1,\dots,\cP_K)$ under $\cP$. If $(E_1,\dots,E_K)$ is independent of $(P_1,\dots,P_K)$, then $P_1/E_1,\dots,P_K/E_K$ are compound p-variables for $(\cP_1,\dots,\cP_K)$ under $\cP$. 
This can be seen from Theorem~\ref{thm:umi}. 
\end{example}

\subsection{Approximate compound p/e-values}

We next generalize the approximate p-variables and e-variables in Definition~\ref{defi:approx} to introduce  approximate compound p-variables and e-variables. 

  Throughout this subsection, we assume that each $\p$ is atomless; this means that we can simulate a uniformly distributed random variable under $\p$. This enables us to formulate the definition of approximate compound p-variables and e-variables following the equivalent formulation for approximate p-variables/e-variables in Proposition~\ref{prop:alter-p}.
  This formulation is helpful for our later results on the approximate and asymptotic FDR control (Theorem~\ref{theo:eBH_controls_the_FDR}).

\begin{definition}[Approximate compound p-variables and e-variables]
\label{defi:approx_compound}
Fix the null hypotheses $\cP_1,\dots,\cP_K$, a set of distributions $\cP$, and two functions $\varepsilon:\cP \to \R_+ $ and $\delta:\cP \to [0,1]$, and write $\varepsilon_\p=\varepsilon(\p)$
and $\delta_\p=\delta(\p)$ for $\p \in \cP$.  
\begin{enumerate}[label=(\roman*)]
    \item Nonnegative random variables $P_1,\dots,P_K$ are \emph{$(\varepsilon, \delta)$-approximate compound} p-variables for $(\cP_1,\dots,\cP_K)$ under $\cP$ if for every $\p \in \cP$,  there exists an event $A$ such that $ \p(A) \geq 1-\delta_\p$ and
$$\sum_{k: \p \in \cP_k}  \p(P_k\le t , A) \le  K(1+\varepsilon_\p )t ~\mbox{~for all $t\in (0,1)$};$$
\item  $[0,\infty]$-valued  random variables   $E_1,\dots,E_K$ are \emph{$(\varepsilon, \delta)$-approximate compound} e-variables for $(\cP_1,\dots,\cP_K)$ under $\cP$ if  for every $\p \in \cP$, 
there exists an event $A$ such that $ \p(A) \geq 1-\delta_\p$ and
$$ \sum_{k: \p \in \cP_k}  \E^{\p}[E_k\id_{A}] \leq K(1+\varepsilon_\p).$$
\end{enumerate}
We sometimes omit ``under $\cP$'' when it is clear.
\end{definition}
From per-hypothesis approximate p-variables (or e-variables),   we can recover approximate compound p-variables (or e-variables) for the multiple testing problem.
\begin{proposition}
Suppose that for all $k \in \mathcal{K}$, $P_k$ (resp.~$E_k$) is an $(\varepsilon_k,\delta_k)$-approximate p-variable (resp.~e-variable) for $\cP_k$. Then, $P_1,\dots,P_K$ (resp.~$E_1,\dots,E_K$) are $(\varepsilon,\delta)$-approximate compound p-variables (resp.~e-variables) for the choice:
$$
\varepsilon(\p) = \frac{1}{K} \sum_{k:\p \in \cP_k} \varepsilon_k(\p),\qquad \delta(\p) = \sum_{k:\p \in \cP_k} \delta_k(\p).
$$

\end{proposition}
\begin{proof}[Proof.]
We carry out the argument only for compound e-variables (the argument for p-variables is analogous).
Take any $\p \in \cP$ and let $A_k$ be the event provided by Proposition~\ref{prop:alter-p} for the $k$-th ($\varepsilon_k, \delta_k)$-approximate e-variable. Then, letting $A := \bigcap_{k: \p \in \cP_k} A_k$, it follows that,
$$\p(A^c) = \p\left(\bigcup_{k:\p \in \cP_k} A_k^c\right)\leq \sum_{k:\p \in \cP_k} \delta_k(\p) = \delta(\p).$$ 
We get
$$ \sum_{k: \p \in \cP_k}  \E^{\p}[E_k\id_{A}] \leq \sum_{k: \p \in \cP_k}  \E^{\p}[E_k\id_{A_k}] \leq \sum_{k: \p \in \cP_k} (1+\varepsilon_k(\p)) \leq K(1 + \varepsilon(\p)).$$
This proves the desired statement for compound e-variables.
\end{proof}


We next provide an analogous statement to that of Proposition~\ref{prop:alter-p} for approximate compound e-variables. 
A corresponding representation for 
approximate compound p-variables
seems to be more complicated, and that is the reason why we chose to formulate  Definition~\ref{defi:approx_compound} through the event $A$ instead of using (E1) below.

\begin{proposition}
\label{prop:equiv_approx_compound_e}
Let $\varepsilon:\cP \to \R_+ $ and $\delta:\cP \to [0,1]$ and $E_1,\dots,E_K$ be $[0,\infty]$-valued random variables. The following two statements are equivalent (assuming that $\p$ is atomless for all $\p\in \cP$):
\begin{enumerate}
\item[(E1)] For all $t\in \R_+$ and $\p\in \mathcal P$, 
$$
\E^\p \left[\left(\sum_{k: \p \in \cP_k}E_k \right)\wedge t \right] \le  K(1+\varepsilon_\p)  +\delta_\p t;
$$
\item[(E2)]  $E_1,\dots,E_K$ are $(\varepsilon, \delta)$-approximate compound e-variables for the hypotheses $\cP_1,\dots,\cP_K$ under $\cP$.
\end{enumerate}
\end{proposition}
\begin{proof}[Proof.]
We start with (E1). By direct manipulation we see that (E1) is equivalent to the following statement:
$$
\E^\p \left[\left(\frac{1}{K}\sum_{k: \p \in \cP_k}E_k \right)\wedge t \right] \le  1+\varepsilon_\p  +\delta_\p t    \mbox{~for all $t\in \R_+$ and all $\p\in \mathcal P$}.
$$
Now fix $\p \in \cP$. According to the above statement,  
$$E:=\frac{1}{K}\sum_{k: \p \in \cP_k}E_k \qquad \mbox{is an $(\varepsilon,\delta)$-approximate e-variable for $\{\p\}$}.$$
By Proposition~\ref{prop:alter-p} there exists an event $A$ with $\p(A) \geq 1-\delta$ such that $\E^{\p}[E\id_{A}] \leq 1+\varepsilon_\p$, that is,
$$ \sum_{k: \p \in \cP_k}  \E^{\p}[E_k\id_{A}] \leq K(1+\varepsilon_\p).
$$
Since $\p \in \cP$ was arbitrary, the above establishes the implication (E1) $\Rightarrow$ (E2). 

For the other direction,
fix an arbitrary $\p\in \cP$. Suppose (E2) holds, that is, $\sum_{k: \p \in \cP_k}  \E^{\p}[E_k\id_{A}] \leq K(1+\varepsilon_\p)$ for some event $A$ with $\p(A) \geq 1-\delta$. Define $E := (1/K) \sum_{k: \p \in \cP_k}E_k$, and  
by applying the direction (E2)$\Rightarrow $(E1) in Proposition~\ref{prop:alter-p}  with the hypothesis $\{\p\}$, we obtain the condition in (E1), completing the proof. 
\end{proof}

We can calibrate approximate compound p-variables to approximate compound e-variables and vice versa as per usual if we use \emph{the same} calibrator for all hypotheses. For instance, given a p-to-e calibrator $h$, if $P_1,\dots,P_K$ are approximate compound p-variables for $\cP_1,\dots,\cP_K$ under $\cP$, then $h(P_1),\dots,h(P_k)$ are approximate compound e-variables for the same hypotheses.

\begin{proposition}
\label{th:calibration_approximate_compound}
 Let $\varepsilon:\cP \to \R_+ $ and $\delta:\cP \to [0,1]$. 
    Calibrating $(\varepsilon, \delta)$-approximate compound p-variables  yields  $(\varepsilon, \delta)$-approximate compound e-variables, and vice versa.
\end{proposition}

\begin{proof}[Proof.]
For the first statement, 
let $P_1,\dotsc,P_K$ be $(\varepsilon, \delta)$-approximate compound p-variables. Fix any $\p\in \cP$ and let $h$ be a (p-to-e) calibrator. Let $A$ be an event with $\p(A) \ge  1-\delta_\p$ such that
$\sum_{k: \p \in \cP_k} \p(P\le t,A )\le K(1+\epsilon_\p) t$ for $t \in (0,1)$, 
and let $b=(1+\epsilon_\p)/\p(A)-1\ge 0.$
Let $\p_A(\cdot) = \p(\cdot \mid A)$. We have
$$\frac{1}{K}\sum_{k: \p \in \cP_k} \p_A(P_k\le t) \le (1+b)t \mbox{~~for $t \in (0,1)$}.$$
Define $U$,$L$, $P_L$, and $\tilde{P}$ as in the proof of Theorem~\ref{prop:compound-pe1}. Then,  $\tilde{P}$ satisfies
\begin{align*}
\p_A( \tilde P  \leq t) &=  \p_A(P_{L} \leq t) \p_A( U\leq K_0/K)  
\\&= \frac{K_0}{K} \frac{1}{K_0}\sum_{k: \p \in \mathcal{P}_k} \p_A(P_k \leq  t) \leq (1+b)t.
\end{align*}
This means that $\tilde P$ is a $(b, 0)$-approximate p-variable for $\{\p_A\}$. By Theorem~\ref{th:calibration}, $h(\tilde{P})$ is a $(b,0)$-approximate e-variable for $\{\p_A\}$. Thus,
\begin{align*}
\frac{1}{K} \sum_{k: \p \in \mathcal{P}_k }\E^{\p_A}\left[   h( P_k) \right] &= \frac{K_0}{K} 
 \E^{\p_A}\left[   h( P_{L} ) \right] \\
 &\le   \E^{\p_A}\left[ h(\tilde P) \right] \leq 1 + b =\frac{1+\epsilon_\p}{\p(A)}.
 \end{align*}
It follows that $\sum_{k: \p \in \mathcal{P}_k }\E^{\p}[h( P_k) \id_A] \leq K(1+\epsilon_\p)$, and therefore, we  have that $h(P_1),\dots,h(P_K)$ are $(\varepsilon,\delta)$-approximate compound e-variables.

The second claim follows by Markov's inequality.
Fix $\p \in \cP$ and let $A$ be the event given by the definition of approximate compound p-variables. Let $P_k = (1/E_k)\wedge 1$. Then, for any $t\in (0,1)$, 
\begin{align*}
 \sum_{k: \p \in \cP_k} \p(P_k \leq t,A) 
 &=  \sum_{k: \p \in \cP_k}  \p(E_k\ge 1/t,A) 
 \\ &=  \sum_{k: \p \in \cP_k}  \p(E_k\id_{A}  \ge 1/t) 
 \\ &\le   \sum_{k: \p \in \cP_k}  \E^\p[E_k\id_{A}] t \leq  K(1+\epsilon_\p)  t,
\end{align*} 
as desired.
\end{proof}

\begin{remark}
We may define an alternative notion of $(\varepsilon, \delta)^*$-approximate compound e-variables and p-variables
as follows. 
    Fix the null hypotheses $\cP_1,\dots,\cP_K$, a set of distributions $\cP$, and two functions $\varepsilon:\cP \to \R_+ $ and $\delta:\cP \to [0,1]$, and write $\varepsilon_\p=\varepsilon(\p)$
and $\delta_\p=\delta(\p)$ for $\p \in \cP$. 
\begin{enumerate}[label=(\roman*)]
    \item 
Nonnegative random variables $P_1,\dots,P_K$ are \emph{$(\varepsilon, \delta)^*$-approximate compound} p-variables for $(\cP_1,\dots,\cP_K)$ under $\cP$ if 
$$\sum_{k: \p \in \cP_k} \p(P_k\le t) \le  K ( (1+\varepsilon_\p )t +\delta_\p ) ~\  \mbox{for all $t\in (0,1)$ and   $\p\in \mathcal P$.} $$
\item $[0,\infty]$-valued random variables $E_1,\dots,E_K$ are \emph{$(\varepsilon, \delta)^*$-approximate compound} e-variables for $(\cP_1,\dots,\cP_K)$ under $\cP$ if  
 $$\sum_{k: \p \in \cP_k} 
 \E^\p [E_k\wedge t] \le  K(1+\varepsilon_\p  +\delta_\p t)   \mbox{~for all $t\in \R_+$ and $\p\in \mathcal P$.}$$ 
\end{enumerate}
This definition is also a natural generalization of Definition~\ref{defi:approx}, but it is different from Definition 
\ref{defi:approx_compound}. 
The difference is transparent from Proposition~\ref{prop:equiv_approx_compound_e}.
We stick to Definition 
\ref{defi:approx_compound} because it fits better with the approximate and asymptotic FDR control that we study below. 
\end{remark}

\subsection{Asymptotic compound p/e-values}
\label{subsec:asymp_compound_evalues}

In the context of multiple testing procedures in Chapter~\ref{chap:compound}, compound e-values are the essential quantities for their FDR control. 
The next definition introduces   compound e-variables  in the approximate or asymptotic sense. 
 The reader may directly anticipate the definitions that follow. However, what is most important, is to clarify the asymptotic setup.
As before, we seek to test hypotheses $\mathcal{P}_1,\dots,\mathcal{P}_K$ based on data $X^{(n)}$. 
We are interested in asymptotics in which either $n \to \infty$ (with $K$ fixed) or $K \to \infty$ (with $n$ fixed) or both $K,n \to \infty$. 

Formally, we consider triangular-array type of asymptotics indexed by a single parameter $m \in \mathbb N$ that determines the dataset size, $n=n(m)$, and the number of hypotheses, $K=K(m)$. We require that $m \to \infty$ and that $\max\{K, n\} \to \infty$ as $m \to \infty$. 

All three asymptotic regimes are frequently encountered in the multiple testing literature. When using Gaussian approximations to construct Z-test or t-test p-values, one is relying on asymptotics as the size of $X$ grows to infinity. When using empirical Bayes arguments to share strength and learn nuisance parameters across hypotheses, one is relying on asymptotics as $K$ grows to infinity.


\begin{definition}[Asymptotic compound p-variables and e-variables]
\label{def:asym-comp-e}
A sequence 
$$\left( E^{(n(m))}_1,\dots, E^{(n(m))}_{K(m)}\right)_{m\in\N}$$ 
of tuples of $[0,\infty]$-valued random variables is called a sequence of:
\begin{enumerate}[label=(\alph*)]
    \item  \emph{asymptotic compound}  e-variables for $\cP_1,\dots,\cP_{K(m)}$ under $\cP$,  if  for each $m$, $E^{(n(m))}_1,\dots, E^{(n(m))}_{K(m)}$ are $(\varepsilon_m, \delta_m)$-approximate compound e-variables for $\cP_1,\dots,\cP_{K(m)}$ under $\cP$
    for some $(\varepsilon_m, \delta_m)$, satisfying  
    \begin{equation}\label{eq:asymp-compound-e} \mbox{for each $\p\in \mathcal P$,}\quad  \lim_{m\to \infty}  \varepsilon_m(\p) = \lim_{m\to \infty} \delta_m(\p) = 0;
    \end{equation}
    
    \item  \emph{strongly asymptotic compound} e-variables for $\cP_1,\dots,\cP_{K(m)}$ under $\cP$, if condition \eqref{eq:asymp-compound-e}  is replaced by   
$$   
\begin{aligned}
    \mbox{for each $\p\in \mathcal P$:} \quad & \lim_{m\to \infty}   \varepsilon_m(\p) = 0   \\&\mbox{~and~}    
 \delta_m(\p) = 0 \mbox{~for $m$ large enough.}
 \end{aligned}$$ 
\end{enumerate}
Uniformly (strongly) asymptotic compound e-variables are defined analogously. 
\end{definition}
Sequences 
of (uniformly, strongly) asymptotic compound p-variables can be defined analogously.

For the remainder of this section, we assume that each $\p$ is atomless; this means that we can simulate a uniformly distributed random variable under $\p$. This allows us to use the equivalent definition in Proposition~\ref{prop:alter-p}  of approximate e-variables through an event $A$ with a constraint on its probability.
As in Definition~\ref{def:FDR}, we denote by $F_{\cD} $ and $R_{\cD}$ the number of false discoveries and the total number of  discoveries, that is,
$$
F_{\cD} = \sum_{k: \p \in \mathcal{P}_k} \id_{\{k\in \cD\}} = |\mathcal N (\p) \cap \cD| \mbox{~~~and~~~} R_{\cD}=|\cD|.
$$
  
\begin{definition}[Approximate and asymptotic FDR control]
\label{def:asym-FDR}
Let  $\cP $ be a set of distributions,  $\mathcal{D}$ be a multiple testing procedure, and let $\varepsilon:\cP \to \R_+ $ and $\delta:\cP \to [0,1]$. The procedure $\mathcal{D}$ is said to have
\begin{enumerate}[label=(\alph*)]
\item \emph{$(\varepsilon,\delta)$-approximate} FDR \emph{control} at level $\alpha$ for   $(\mathcal{P}_1,\dotsc,\mathcal{P}_K)$ under $\cP$, if for every $\p \in \cP$, 
there exists an event $A$ such that
$$\E^{\p}\left[\frac{F_{\cD}}{R_{\cD}} \id_{A}\right] \le  \alpha(1+\varepsilon_\p ) ~\mbox{~~and~~ }\p(A) \geq 1-\delta_\p;$$
\item \emph{asymptotic} FDR \emph{control} at level $\alpha$
for   $(\cP_1,\dots,\cP_{K(m)})$ under $\cP$ if for each $m \in \mathbb N$  it has $(\varepsilon_m, \delta_m)$-approximate FDR control at level $\alpha$
for some $(\varepsilon_m, \delta_m)$, satisfying  
$$
\mbox{for each $\p\in \mathcal P$,}\quad  \lim_{m\to \infty}  \varepsilon_m(\p) = \lim_{m\to \infty} \delta_m(\p) = 0.
$$
\end{enumerate}
\end{definition}

It is clear that 
the $(0,0)$-approximate FDR control is  equivalent to  the 
FDR control defined in Definition~\ref{def:FDR}. 
  A procedure with $(\varepsilon, \delta)$-approximate FDR control also satisfies $$\mathrm{FDR}_{\cD} = \E^{\p}\left[\frac{F_{\cD}}{R_{\cD}}  \right]  \leq \alpha(1+\varepsilon_{\p}) + \delta_{\p} \text{ for any } \p \in \cP.$$
Conversely, any procedure satisfying the above inequality 
has $(\epsilon,\delta/\alpha)$-approximate FDR control, by invoking the equivalence in Proposition~\ref{prop:alter-p}; more precisely,
the above condition implies, for all $t\ge 1$,
$$\E^{\p}\left[\frac{F_{\cD}}{\alpha {R_{\cD}}} \wedge t\right] \le\E^{\p}\left[\frac{F_{\cD}}{\alpha {R_{\cD}}} \right] 
\le (1+\epsilon_\p) + \frac{\delta_\p}{\alpha} \le(1+\epsilon_\p) + \frac{\delta_\p}{\alpha} t . $$


 
With Definitions~\ref{def:asym-comp-e} and~\ref{def:asym-FDR}, we can formally state that the e-BH procedure applied to asymptotic compound e-values has asymptotic FDR control.  

\begin{theorem}
\label{theo:eBH_controls_the_FDR}
Suppose that we apply the e-BH procedure at level $\alpha$ to $E_1,\dots,E_K$. Then:
\begin{enumerate}[label=(\roman*)]
\item If $E_1,\dots,E_K$ are compound e-variables, then e-BH controls the FDR at level $\alpha$.
\item If $E_1,\dots,E_K$ are $(\varepsilon, \delta)$-approximate compound e-variables, then e-BH satisfies $(\varepsilon, \delta)$-approximate FDR control, and so it controls
the FDR at level $\alpha(1+\varepsilon) + \delta$. 
\item If $E_1,\dots,E_K$ are asymptotic compound e-variables, then e-BH has asymptotic FDR control at level $\alpha $.
\end{enumerate}
The same conclusions also hold for any  self-consistent e-testing procedures at level $\alpha$ (Definition~\ref{def:self-cons-e}).
\end{theorem}
 
Theorem~\ref{theo:eBH_controls_the_FDR}
follows directly from the definitions and the previous results on the e-BH procedure.

\subsection{The simultaneous t-test sequence model}
\label{subsec:t_tests}

We next describe a prototypical model of constructing asymptotic compound e-variables in a setting with composite null hypotheses.

In our setting, we observe $K$ groups of data of size $n$, that is, we observe $X=(X_k^1,\ldots,X_k^n)_{k \in \mathcal{K}}$. We consider distributions $\p$ according to which all observations are jointly independent and the data within each group are iid. The universe of distributions $\cP$ is understood to obey the above assumptions. Given any $\p \in \cP$, we write $\p_k$ for its $k$-th group-wise marginal so that $X \sim \p$ entails that 
$$
X_k^i \simiid \p_k\;\; \text{ for }\;\; i=1,\dotsc,n.
$$
We seek to test the hypotheses,
$$
\cP_k = \{ \p \in \cP\,:\,  \E^{\p_k}[X_k^i] = 0,\,\var^{\p_k}(X_k^i) \in (0,\infty)\}.
$$
To this end, we start by constructing summary statistics in a group-wise manner, similar to~\eqref{eq:one_sample_summary}:
\begin{equation}
\begin{aligned}
&
\bar{X}^{(n)}_k = \frac{1}{n}\sum_{i=1}^n X^i_k,\;\;\;\, S^{(n)}_k = \sqrt{\frac{1}{n}\sum_{i=1}^n (X^i_k)^2},\\& \hat{\sigma}_k^{(n)} = \sqrt{\frac{1}{n-1}\sum_{i=1}^n \left(X^i_k - \bar{X}_k^{(n)} \right)^2}.
\end{aligned}
\label{eq:one_sample_summary_mtp}
\end{equation}
To motivate the construction that follows, fix any $\p \in \cP$ and write  $\mu_k(\p) =\E^{\p_k}[X_k^i]$ and $\sigma_k^2(\p) = \var^{\p_k}(X_k^i)$ (supposing these exist). 
Observe that for any $k \in \mathcal{K}$, 
$$\E^{\p_k}[(S_k^{(n)})^2] = \sigma_k^2(\p) + \mu_k^2(\p) \, \text{ for all }\,  k \in \mathcal{K},
$$ and so $$ \E^{\p_k}[(S_k^{(n)})^2] = \sigma_k^2(\p)\, \text{ for }\, \p \in \cP_k.
$$
One can check that for any $\p \in (\bigcup_{k \in \mathcal{K}} \cP_k)\bigcap \cP$, 
$$
\sum_{k: \p \in \cP_k} \E^{\p}\left[\tilde{E}_k^{(n)}(\p)\right]  = K, \quad \text{where} \quad \tilde{E}_k^{(n)}(\p) = \frac{K \big(S_k^{(n)}\big)^2}{ \sum_{j: \p \in \cP_j} \sigma_j^2(\p)}.
$$
The above construction is not directly usable in practice because of the dependence of $\tilde{E}_k^{(n)}(\p)$ on $\p$. However, it motivates a construction in which we conservatively estimate the $\p$-dependent quantity $\sum_{j: \p \in \cP_j} \sigma_j^2(\p)$ by $\sum_{j \in \mathcal{K}}(\hat{\sigma}^{(n)}_{j})^2$ and then let
\begin{equation*}
E^{(n)}_k = \frac{K \big(S_k^{(n)}\big)^2}{\sum_{j \in \mathcal{K}}\big(\hat{\sigma}^{(n)}_{j}\big)^2},\quad k \in \mathcal{K}.
\end{equation*}
We will show that the above are asymptotic compound e-variables under an asymptotic setup indexed by $m \in \mathbb N$ such that $n=n(m)$, $K=K(m)$ and $\max\{n, K\} \to \infty$ as $m \to \infty$. In case $K \to \infty$, we further restrict the class $\cP$ to require that for any $\p \in \cP$, $|\{k \in \mathbb N: \p \in \cP_k\}| = \infty$. In addition to this restriction (which remains implicit in our notation), our formal result will be stated in terms of the following two sets of distributions. The first set imposes some restrictions on moments of the marginals for null groups,
\begin{align*}
\cP^{\mathrm{B}}& =\{\p \in \cP: \E^{\p_k}[|X_k^i|^{2(1+\delta)}] \leq C,\;\var^{\p_k}(X_k^i) \geq c, ~\forall k \text{ with }\p \in \cP_k\},
\end{align*}
where $ c, C \in (0, \infty)$ 
and $\delta \in (0,1]$ are fixed constants. The second set also imposes normality of the marginals for null groups,
\begin{align*}
\cP^{\mathrm{N}} &=\{\p \in \cP: \p_k = \mathrm{N}(0, \sigma_k^2(\p) ),  c \leq \sigma_k^2(\p) \leq C, ~\forall k\text{ with }\p \in \cP_k\},
\end{align*}
for some constants $c,C \in (0, \infty)$.

\begin{proposition}
Suppose that $n=n(m) \geq 2$ and $K=K(m) \geq 1$ are such that $\max\{n, K\} \to \infty$ as $m \to \infty$. 
Then:
\begin{enumerate}[label=(\roman*)]
\item $(E_1^{(n)},\ldots, E_K^{(n)})_{m\in \N}$ is a sequence of asymptotic compound e-variables for $(\cP_1,\ldots,\cP_K)$ under $\cP^{\mathrm{B}}$. 
\item $(E_1^{(n)},\ldots, E_K^{(n)})_{m\in \N}$ is a sequence of  strongly asymptotic compound e-variables for $(\cP_1,\ldots,\cP_K)$ under $\cP^{\mathrm{N}}$. 
\end{enumerate}
\end{proposition}

\begin{proof}[Proof.]
Fix an arbitrary distribution $\p \in (\bigcup_k \cP_k) \bigcap \cP^{\mathrm{B}}$ and write $\sigma_k^2 = \sigma_k^2(\p)$, omitting the dependence on $\p$.
We write $K_0(\p) = |\mathcal{N}(\p)| = \{k \in \mathcal{K}: \p \in \cP_k\}$. Note that $\max\{K_0(\p), n\} \to \infty$ and that
\begin{equation}
\label{eq:asymp_comp_evalue_renormalized}
\frac{1}{K}\sum_{k: \p \in \cP_k} E_k^{(n)} \leq \frac{\sum_{k:\p \in \cP_k } \big(S_k^{(n)}\big)^2}{\sum_{k: \p \in \cP_k}\big(\hat{\sigma}^{(n)}_{k}\big)^2}.
\end{equation}
We will argue that the right-hand side of the above expression converges in probability (and thus also in distribution) to the constant $1$. Arguing as in Proposition~\ref{prop:asym-e}, part  (i), it will then follow that $E_1^{(n)},\ldots, E_K^{(n)}$ are asymptotic compound e-variables. We claim that
$$
\frac{n\sum_{k:\p \in \cP_k } \big(S_k^{(n)}\big)^2}{ n\sum_{k:\p \in \cP_k } \sigma_k^2} \stackrel{\p}{\to} 1 \;\;\text{ as }\;\; m \to \infty.
$$
Writing $U^{(m)} = n\sum_{k:\p \in \cP_k } \big(S_k^{(n)}\big)^2 = \sum_{k: \p \in \cP_k} \sum_{i=1}^n (X_{k}^i)^2$, we see that $\E^{\p}[U^{(m)}] = n\sum_{k:\p \in \cP_k } \sigma_k^2$. Moreover, by our assumptions and using the Marcinkiewicz--Zygmund inequality, there exists another constant $C'$ that depends only on $\delta>0$ such that:
$$
\E^{\p}[U^{(m)}] \geq cnK_0(\p),\quad \E^{\p}[|U^{(m)} - \E^{\p}[U^{(m)}]|^{1+\delta}] \leq C'  n K_0(\p).
$$
Thus, since $n K_0(\p) \to \infty$ and using Markov's inequality, it follows that $U^{(m)}/ \E^{\p}[U^{(m)}]$ converges in probability to $1$ as desired. We can analogously argue that:
$$
\frac{n\sum_{k:\p \in \cP_k } \big(\hat{\sigma}^{(n)}_{k}\big)^2}{ n\sum_{k:\p \in \cP_k } \sigma_k^2} \stackrel{\p}{\to} 1 \;\;\text{ as }\;\; m \to \infty.
$$
Thus implies that $E_1^{(n)},\ldots, E_K^{(n)}$ are asymptotic compound e-variables.

Next fix $\p \in (\bigcup_k \cP_k) \bigcap \cP^{\mathrm{N}}$. It suffices to prove uniform integrability of the right-hand side in~\eqref{eq:asymp_comp_evalue_renormalized} for all $m$ such that $(n-1)K_0(\p) \geq 5$ (see Proposition~\ref{prop:asym-e}, part (iii)). We prove uniform integrability by controlling the following $\eta$-th moment for $\eta \in (1, 1.25)$:
\begin{align*}
&\E^{\p}\left[ \left( \frac{\sum_{k:\p \in \cP_k } \big(S_k^{(n)}\big)^2}{\sum_{k: \p \in \cP_k}\big(\hat{\sigma}^{(n)}_{k}\big)^2} \right)^{\eta} \right] \\&\leq
\E^{\p}\left[ \left(\sum_{k:\p \in \cP_k } \big(S_k^{(n)}\big)^2\right)^{2\eta}\right]^{1/2} \E^{\p}\left[ \left( \frac{1}{\sum_{k: \p \in \cP_k}\big(\hat{\sigma}^{(n)}_{k}\big)^2} \right)^{2\eta} \right]^{1/2}.
\end{align*}
Notice that by our assumptions, $\sigma_k^2=\sigma_k^2(\p)$ is both uniformly lower bounded and upper bounded for all $k$ such that $\p \in \cP_k$. Ignoring multiplicative constants, it suffices to bound:
$$
\E^{\p}\left[ \left(\sum_{k:\p \in \cP_k } n\left(\frac{S_k^{(n)}}{\sigma_k}\right)^2\right)^{2\eta}\right]^{1/2} \E^{\p}\left[ \left( \frac{(n-1)^{-1}}{\sum_{k: \p \in \cP_k}\big({\hat{\sigma}^{(n)}_{k}}/{\sigma_k^2}  \big)^2} \right)^{2\eta}\right]^{1/2}. 
$$
This can be accomplished by noting that,
$$\sum_{k:\p \in \cP_k } n\big(S_k^{(n)}/\sigma_k\big)^2 \lawis \chi^2_{nK_0(\p)},~\sum_{k: \p \in \cP_k}(n-1)\big(\hat{\sigma}^{(n)}_{k}/\sigma_k \big)^2 \lawis \chi^2_{(n-1)K_0(\p)} ,$$
where $\chi^2_d$ denotes the chi-square distribution with $d$ degrees of freedom. The remainder of the calculation amounts to direct evaluation of the $(2\eta)$-th moment (resp.~$(-2\eta)$-th moment) of chi-squared random variables.
\end{proof}

\section*{Bibliographical note}

Materials in this chapter are largely based on~\cite{ignatiadis2024compound}.  
 The notion of $(\varepsilon, 0)$-approximate e-variable (not using this terminology) was considered by~\citet[Section 6.2]{wang2022false}.
The standard examples following Proposition~\ref{proposition:weak_convergence_to_asymptotic} in Section~\ref{sec:c9-conditions} can be found in standard textbooks 
 \cite{vandervaart1998asymptotic, lehmann2005testing}.
 Example~\ref{ex:c2-mv} is found in \citet{wang2023catoni}.  The notion of compound p-variables was introduced by~\citet{armstrong2022false} under a different name (average significance controlling p-values and tests) and further studied in an empirical Bayes setting by~\citet{ignatiadis2024empirical}. 
The notion of
approximate compound e-variables appears implicitly in~\citet{ignatiadis2024values}.
 A notion related to approximate FDR control was used as proof technique (without being given that name) by~\citet{blanchard2008two}.
  The construction in Section~\ref{subsec:t_tests} appeared in~\citet{ignatiadis2024compound,ignatiadis2024values} and pertains to the problem of conducting multiple t-tests.

\chapter{Combining e-values and using e-values as weights}
\label{chap:combining}

In this chapter, we 
first present a classic result in optimal transport. Using this result, we 
prove  Theorem~\ref{th:adm-merg}, which characterizes the class of all e-merging functions. Then, we discuss how to combine an e-value and a p-value. Based on this combination, we illustrate how e-values can be used as weights in multiple testing procedures.  

Throughout this chapter, $K\ge 2$ is a positive integer. We continue to fix one atomless probability measure $\p$  on some $(\Omega,\cF)$ as in Chapter~\ref{chap:multiple}.

\section{Optimal transport duality}

A technical tool that will be useful in two main proofs in Sections~\ref{sec:10-emerg} and~\ref{sec:10-2} is optimal transport duality. 
We state a classic version of multi-marginal optimal transport duality without proving it. 

Let $F:[0,\infty)^K\to \R$ be a bounded  Borel function.
  Denote by $\mathcal B$ the set of   Borel functions on $[0,\infty)$,
  and define the operator $\bigoplus$ 
  on $(\phi_1,\dots,\phi_K)\in\mathcal B^K$ 
  by  
  \[
    \left(
      \bigoplus_{k=1}^K \phi_k
    \right)
    (x_1,\dots,x_K)
    :=
    \sum_{k=1}^K \phi_k(x_k),
 \quad(x_1,\dots,x_K)\in[0,\infty)^K.
  \]  
 Define 
$$D_F=\left\{(\phi_1,\dots,\phi_K) \in \mathcal B^K: \bigoplus_{k=1}^K   \phi_k   \ge F\right\}.$$   
Let $\Gamma(\mu_1,\dots,\mu_K)$ be the set of all Borel measures on $\R^K$ with marginal distributions $\mu_1,\dots,\mu_K$ on $\R$.   In what follows, we always write $\boldsymbol \phi=(\phi_1,\dots,\phi_K)$
and 
$\boldsymbol \mu=(\mu_1,\dots,\mu_K)$. 

  \begin{lemma}[Optimal transport duality]
  \label{lem:duality}For any bounded  Borel function $F:[0,\infty)^K\to \R$ 
  and distributions $\mu_1,\dots,\mu_K$ on $\R$,
   \begin{align}
\label{eq:duality-eq}
  \sup_{\pi\in \Gamma(\boldsymbol \mu) } \int F \d \pi = \inf_{\boldsymbol \phi \in D_F}  \sum_{k=1}^K\int \phi_k \d \mu_k,
\end{align}
where in the infimum we only consider those with $  \sum_{k=1}^K\int \phi_k \d \mu_k$ well-defined. 
  Moreover, if $F$ is upper semicontinuous, then the infimum is attainable, and the functions $\phi_1,\dots,\phi_K$ 
 in \eqref{eq:duality-eq} can be chosen as upper semicontinuous. 
  If $F$ is decreasing (nonnegative), then $\phi_1,\dots,\phi_K$ can be chosen as decreasing (nonnegative). 
   If $F$ is bounded in $[0,1]$, then $\phi _1,\dots,\phi_K$ can be chosen as bounded in $[0,1]$. 
   
  \end{lemma}

  At an abstract level, the reason why optimal transport duality appears in this chapter and the next one is that merging functions both for  p-values and for e-values need to produce an output under arbitrary dependence, and optimization under arbitrary dependence is precisely the topic of optimal transport theory. 
  In our context, it concerns multi-marginal optimal transport, to be precise.

\section{Admissible e-merging functions}

\label{sec:10-emerg}

The main content of this section is to prove Theorem~\ref{th:adm-merg}, recalled here.
\begin{proof}[Recalling Theorem~\ref{th:adm-merg}]
\renewcommand{\qedsymbol}{}
For a function $ F:\R_+^K\to\R_+$,
 \begin{enumerate}
 \item[(i)] if $F$ is an  e-merging function, then
  $F\le \fM_{ \boldsymbol \lambda }$ for some  $\boldsymbol \lambda \in \Delta_{K+1}$;  
\item[(ii)] $F$ is an  admissible e-merging function  if and only if     $F = \fM_{ \boldsymbol \lambda }$ for some $\boldsymbol \lambda \in \Delta_{K+1}$. 
\end{enumerate}
\end{proof}
It suffices to show statement (i). Statement (ii) follows directly from (i) and the simple fact that $\fM_{\boldsymbol \lambda }$ is an e-merging function.

First, we characterize the simplest case $K=1$, which will be used in the proof of Theorem~\ref{th:adm-merg}.   
All e-variables in this section are for $\p$.

\begin{lemma}\label{lem:mean2}
Let $\theta\in[ 1,\infty]$ and $r\in \R_+$. 
If a function $g:[0,\theta] \to \R_+$ satisfies 
$\E^\p [g(E)]\le r  $ for all   e-variables $E$ taking     values in $[0,\theta]$, then there exists $\lambda \in [0,1]$ such that  $g(x)\le r (1-\lambda + \lambda x )$ for $x\in [0,\theta]\cap \R$.
\end{lemma}


As a particular case of Lemma~\ref{lem:mean2} with $r=1$, if $g$ is an e-merging function of dimension $1$, then  there exists $\lambda \in [0,1]$ such that 
$g(x) \le (1-\lambda) + \lambda x  $ for $x\in \R_+$.
Hence, the desired statement in Theorem~\ref{th:adm-merg}  holds for $K=1$. 

Lemma~\ref{lem:mean2} has the following useful implication. 
If one is supplied with an e-variable $E$, but wants another e-value (hoping that it would work better than $E$ in some context), then one has to use $(1-\lambda)+ \lambda E$ for some $\lambda \in[0,1]$, as all other ways are dominated by this class.
 For instance, $\sqrt{E}$ is a valid e-variable for any $E\in \fE$, but Lemma~\ref{lem:mean2} says that it is dominated by $(1-\lambda)+ \lambda E$ for some $\lambda$, which turns out to be $1/2$ in this case.
 This is the reason why there is no  way of betting in each round other than using $(1-\lambda_t+\lambda_t E_t)$ in Section~\ref{sec:EAEP} to build an e-process from given sequential e-variables $(E_t)_{t\in T_+}$.

\begin{proof}[Proof of Lemma~\ref{lem:mean2}.]
 
The case $r=0$ is trivial as $g$ is the constant function $0$. 
Otherwise, $r>0$, and without loss of generality we can assume $r=1$.

If $\theta=1$, then   $g(x)\le 1$ for all $x$, and taking $\lambda=0$ gives the desired inequality. 
In what follows, we assume $\theta \in (1,\infty]$,
and all points $x,y$ that appear below are in $[0,\theta]\cap \R$.

 First, it is easy to note that $g(y)\le y$ for $y>1$; indeed, if $g(y)>y$ then taking a random variable $X$ with $\p(X=y)=1/y$ and $0$ otherwise  gives $\E^\p[g(X)]> 1$ and breaks the assumption.
 Moreover, $g(y)\le 1$ for $y\le 1$ is also clear, which in particular implies $g(1)\le 1$.

Suppose for the purpose of contradiction that the statement in the lemma does not hold. 
This means that for each $\lambda\in [0,1]$, either (a) $g(x) >(1-\lambda) + \lambda x$ for some $x<1$ or (b) $g(y) >(1-\lambda) + \lambda y$ for some $y>1$ (or both).
Since $g(y)\le y$ for $y>1$ and $g(x)\le 1$ for $x<1$, we know that
$\lambda =1$ implies (a)
and  $\lambda =0$ 
implies (b). 

We claim that there exists $\lambda_0\in (0,1)$ for which both (a) and (b) happen. 
To show this claim, let 
\begin{align*} 
\Lambda_0&= \{ \lambda \in [0,1]: g(y) >(1-\lambda) + \lambda y \mbox{ for some $y>1$}\};\\\Lambda_1&= \{ \lambda \in [0,1]: g(x) >(1-\lambda) + \lambda x \mbox{ for some $x<1$}\}.\end{align*}
Clearly, the above arguments show $\Lambda_0\cup \Lambda_1=[0,1]$, $0\in \Lambda_0$, and $1\in \Lambda_1$. Moreover, since the function $\lambda 
\mapsto  (1-\lambda) + \lambda x$ is monotone  for either $x<1$ or $x>1$, we know that both $\Lambda_0$ and $\Lambda_1$ are intervals. 
Let  
$ \lambda_*=\sup \Lambda_0$ and  
$ \lambda^*=\inf \Lambda_1.$
We will argue  $\lambda_*\not \in \Lambda_0$ 
and $ \lambda^* \not \in   \Lambda_1.$ 
If $\lambda_* \in \Lambda_0$, then there exists $y>1$ such that   $g(y) >(1-\lambda_*) + \lambda_* y $. By continuity, there exists $\hat \lambda_*>\lambda_*$ such that  $g(y) >(1-\hat \lambda_*) + \hat\lambda_* y $, showing that $\hat \lambda_*\in \Lambda_0$, contradicting the definition of $\lambda_*$. Therefore, $\lambda_* \not \in \Lambda_0$.
Similarly,  $\lambda^* \not \in \Lambda_1$, following the same argument.
If $\lambda_*  = \lambda^*$, then this point is not contained in $\Lambda_0\cup \Lambda_1$, a contradiction to $\Lambda_0\cup \Lambda_1=[0,1]$. Hence, it must be  $\lambda_* > \lambda^*$, which implies that $\Lambda_0\cap \Lambda_1 $ is not empty.

Let $x_0<1$ and $y_0>1$ be such that 
  $$ g(x_0)> 1-\lambda_0 +  \lambda_0 x_0  \mbox{~~~and~~~} g(y_0)> 1-\lambda_0 +  \lambda_0 y_0 .$$    
Let $X$ be distributed as 
$\p(X=y_0) = (1-x_0)/(y_0-x_0)$
and 
 $\p(X=x_0) = ( y_0-1)/(y_0-x_0)$, which clearly satisfies 
 $\E^\p[X]=1$ and is binary.   
Moreover, 
\begin{align*}
\E^\p[g(X)] &=
\frac{ 1-x_0   }{y_0-x_0 }  g(y_0 ) + 
\frac{ y_0-1   }{y_0-x_0 } g(x_0)
\\&> 
\frac{ 1-x_0   }{y_0-x_0 } \left(1-\lambda_0 +  \lambda_0 y_0 \right) + 
\frac{ y_0-1   }{y_0-x_0 }\left(1-\lambda_0 +  \lambda_0 x_0\right)=1.
\end{align*}
This yields a contradiction.  
\end{proof}

   In what follows,  $\mathbf 0$ and $\mathbf 1$ represent a vector  of zeros and a vector of ones of the appropriate dimension, respectively. 
 The next lemma allows us to only consider upper semicontinuous e-merging functions. 

\begin{lemma}\label{lem:usc}
  If $F$ is an e-merging function, then its upper semicontinuous version  $F^*$ is given by
\begin{align}\label{eq:usc}
  F^*(\mathbf{e})
  =
  \lim_{\epsilon\downarrow 0}
  F(\mathbf{e} + \epsilon\mathbf 1),
  \quad
  \mathbf{e}\in[0,\infty)^K;
\end{align}  is also an e-merging function, and $F^*\ge F$.
\end{lemma}
\begin{proof}[Proof.]
  Take any vector $\mathbf{E} $ of e-variables.
  For every rational $\epsilon \in (0,1)$, let $A_\epsilon$ be an event independent of $\mathbf{E}$ with $\p(A_\epsilon)=1-\epsilon$,
  and $\mathbf{E}_{\epsilon} = (\mathbf{E}+\epsilon\mathbf 1)\id_{A_\epsilon}$
  (here we use the convention that $\mathbf{E}_\epsilon=\mathbf{0} $
  if the event $A_{\epsilon}$ does not occur).
  For each $\epsilon$, $\E[\mathbf{E}_\epsilon] \le (1-\epsilon) (\mathbf 1+\epsilon \mathbf{1}) \le \mathbf{1}$. 
  Therefore, $\mathbf{E}_\epsilon $ is a vector of e-variables,  and hence
  \begin{align*}
    1
    \ge
    \E[F(\mathbf{E}_\epsilon)]
    =
    (1-\epsilon)
    \E\left[F(\mathbf{E}+\epsilon\mathbf{1})\right]
    +
    \epsilon F(\mathbf{0}),
  \end{align*}
  which implies
  \[
    \E
    \left[
      F(\mathbf{E}+\epsilon\mathbf{1})
    \right]
    \le
    \frac{1 - \epsilon F(\mathbf{0})}{1-\epsilon}.
  \]
  Fatou's lemma yields
  \begin{align*}
    \E[F^*(\mathbf{E})]
    &=
    \E
    \left[
      \lim_{\epsilon \downarrow 0}
      F(\mathbf{E}+\epsilon\mathbf{1})
    \right] \\
    &\le
    \lim_{\epsilon \downarrow 0}
    \E
    \left[
      F(\mathbf{E}+\epsilon\mathbf{1})
    \right]
    \le
    \lim_{\epsilon \downarrow 0}
    \frac{1 - \epsilon F(\mathbf 0)}{1-\epsilon}
    =
    1.
  \end{align*}
  Therefore, $F^*$ is an e-merging function. The domination $F^*\ge F$ follows from the fact that $F$ is increasing. 
\end{proof}

\begin{proof}[Proof of Theorem~\ref{th:adm-merg}]
By Lemma~\ref{lem:usc}, it suffices to consider an upper semicontinuous e-merging function $F$.  
Let $\mathcal M_{\mathcal E}$ be the set of all distributions on $\R_+$ with mean no larger than $1$, i.e., the set of all distributions of e-variables,
and   $\mathcal M_{\mathcal E}^\theta $ be the subset of $\mathcal M_{\mathcal E} $ containing all distributions on $[0,\theta ]$ for $\theta \ge 1$.

Fix $\theta \ge 1$. 
Since $F$ is bounded and nonnegative on $[0,\theta]^K$ and upper semicontinuous, using Lemma~\ref{lem:duality} we get
\begin{align}
\label{eq:c7-duality1}
\begin{aligned}
  \sup_{ \boldsymbol\mu \in  (\mathcal M ^\theta_{\mathcal E})^K  } \sup_{\pi\in \Gamma(\boldsymbol \mu) } \int F \d \pi 
  &=   \sup_{ \boldsymbol\mu \in  (\mathcal M ^\theta_{\mathcal E})^K  } \inf_{\boldsymbol \phi \in D_F}  \sum_{k=1}^K\int \phi_k \d \mu_k  \\
  &=
    \sup_{ \boldsymbol\mu \in  (\mathcal M ^\theta_{\mathcal E})^K  } \inf_{\boldsymbol \phi \in D_F^+}  \sum_{k=1}^K\int \phi_k \d \mu_k,   
    \end{aligned}
\end{align}
where $D_F^+$ is the subset of $D_F$ containing all $\boldsymbol \phi$ with nonnegative and upper semicontinuous components.  
Define the mapping $$J:( \boldsymbol \mu,\boldsymbol \phi)\mapsto \sum_{k=1}^K  
 \int \phi_k \d \mu_k .$$
 We will verify a few conditions, which allows us to apply a minimax theorem.
 
\begin{enumerate}[label=(\roman*)] 
 \item  
 The mapping $J$ is bilinear, and therefore both convex and concave. 
\item
Since $[0,\theta ]$ is compact, the set $(\mathcal M ^\theta_{\mathcal E})^K$ equipped with the weak topology is tight.
   By Prokhorov's theorem, it is   sequentially compact, and hence compact. 
\item Since each component of $ \boldsymbol \phi $ is upper semicontinuous,
$\boldsymbol \mu \mapsto J( \boldsymbol \mu,\boldsymbol \phi)$ is upper semicontinuous with respect to 
weak convergence in $\boldsymbol \mu$ for each $\boldsymbol \phi$. 
\end{enumerate} 
The above conditions allow us to use 
\citet[Theorem 2]{fan1953minimax} to conclude   
\begin{equation}
\label{eq:c7-duality2}
  \sup_{ \boldsymbol\mu \in  (\mathcal M ^\theta_{\mathcal E})^K   } \inf_{\boldsymbol \phi \in D^+_F}  \sum_{k=1}^K\int \phi_k \d \mu_k  
  = \inf_{\boldsymbol \phi \in D^+_F}    \sup_{ \boldsymbol\mu \in   (\mathcal M ^\theta_{\mathcal E})^K   }  \sum_{k=1}^K\int \phi_k \d \mu_k. 
\end{equation} 
Since $F$ is an e-merging function, 
each term in \eqref{eq:c7-duality1} is bounded by $1$. Using \eqref{eq:c7-duality2}, we get  
\begin{equation}
1 \ge \inf_{\boldsymbol \phi \in D^+_F}  \sup_{\boldsymbol \mu \in (\mathcal M^\theta _{\mathcal E})^K}  \sum_{k=1}^K\int \phi_k \d \mu_k  = \inf_{ \boldsymbol \phi  \in D^+_F} \sum_{k=1}^K  \sup_{ \mu_k  \in \mathcal M^\theta_{\mathcal E}}   \int \phi_k \d \mu_k . \label{eq:value}
 \end{equation}
 Denote by $T_{\phi}= \sup_{ \mu \in \mathcal M_{\mathcal E}^\theta}   \int \phi  \d \mu  $ for $\phi\in \mathcal B$. 
 For any $\epsilon>0$, by   \eqref{eq:value},  we can find  $(\phi_1,\dots,\phi_K)\in D^+ _F$
 such that $\sum_{k=1}^K T_{\phi_k}\le 1+\epsilon$.
 Using Lemma~\ref{lem:mean2}, 
 for each $k\in  \cK$, there exists a constant $h_{\phi_k}\in [0,1]$ such that 
 $$ \phi_k(x) \le T_{\phi_k}\left(1-h_{\phi_k} + h_{\phi_k} x\right) \mbox{~for all $x\in [0,\theta]$.}$$ 
Since $(\phi_1,\dots,\phi_K)\in D^+_F$, we have
$$
F(x_1,\dots,x_K) \le\sum_{k=1}^K  T_{\phi_k}\left(1-h_{\phi_k} + h_{\phi_k} x_k\right)   \mbox{~for all  $x_1,\dots,x_K \in [0,\theta]$}.
$$
This means that there exists $\boldsymbol \lambda_{\epsilon,\theta} \in \Delta_{K+1}$
such that   $F \le (1+\epsilon)\fM_{ \boldsymbol \lambda_{\epsilon,\theta} }$ on $ [0,\theta]^K$. 
Since $\Delta_{K+1}$ is compact, 
we can find a limit $ \boldsymbol \lambda_0 \in \Delta_{K+1}$ of some subsequence of $\boldsymbol \lambda_{\epsilon,\theta} $  as $\epsilon \downarrow 0$ and $\theta\to \infty$.
Continuity of  $ \boldsymbol \lambda 
\mapsto\fM_{ \boldsymbol \lambda  }$   
yields $ F \le \fM_{ \boldsymbol \lambda_0  } $ on $ \R_+^K$, 
thus the desired statement. \end{proof}

\section{Cross-merging: merging an e-value and a p-value}

\label{sec:10-cross}

Here we consider the setting where we have one e-variable $E$ and one p-variable $P$ for the same hypothesis.
This will help to design and understand procedures when both p-values and e-values are available, treated in Section~\ref{sec:10-weighted}.

Propositions~\ref{prop:e-to-p} and~\ref{prop:p-to-e} in Chapter~\ref{chap:markov} identify all admissible calibrators between p-values and e-values.   
  We consider four cases of merging a p-value $P$ and an e-value $E$: whether $P$ and $E$ are independent, and whether the output is a p-value or an e-value.  
  
To discuss these cases, we define combiners, in a similar way as calibrators and merging functions.  
\begin{definition}
    \label{def:combiners}
A function $f: [0,1]\times [0,\infty] \to  [0,\infty]$
is called an i-pe/e combiner 
if  $f(P,E)$ is an e-value for any independent p-value $P$ and e-value $E$,
and $(p,e)\mapsto f(p,e)$ is decreasing in $p$ and increasing in $e$.
Similarly, we define i-pe/p, pe/p, and pe/e combiners, where i indicates independence, and p and e are self-explanatory.
If the output is a p-value,
the combiner is increasing in $p$ and decreasing in $e$.

\end{definition}

In what follows we omit ``pe'' in ``i-pe/e'' and other terms to be concise. 
We provide four natural combiners to the above four cases, some relying on   an admissible calibrator $h $.
\begin{enumerate} 
    \item[{[ie]}] Return $h(P)  E$ by using the function $C^{\mathrm{ie}}_h:(p,e)\mapsto h(p)e$. The convention here is $0\times \infty=\infty$. 
    \item[{[ip]}]  Return $P/E$, capped at $1$, by  using the function $C^{\mathrm{ip}}:(p,e)\mapsto (p/e)\wedge 1$.
    \item[{[e]}] Return $\lambda h(P)+(1-\lambda) E$  by   using the function $C^{\mathrm{e}}_{\lambda,h} : (p,e)\mapsto \lambda h(p)+(1-\lambda) e $ for some $\lambda \in (0,1)$.
    \item[{[p]}]  Return $2\min(P,1/E)$, capped at $1$,  by  using the function $C^{\mathrm{p}}:(p,e)\mapsto  ( 2 (  p\wedge  e^{-1}) )\wedge 1$.
\end{enumerate}

The superscript in the notation $C^{\rm ie}$ suggests that the combiner assumes independence and outputs an e-value; the other cases are similar.   The combiners  
$C^{\mathrm{ie}}_h$ and $C^{\mathrm{e}}_{\lambda,h} $ depend on $h$ whereas $C^{\mathrm{ip}} $ and $C^{\mathrm{p}}$ do not, similarly to the situation of calibrators: There are many more choices when the output is an e-value, compared to the case when the output is a p-value.  
For the combiner $C^{\mathrm{e}}_{\lambda,h}$, it may be convenient to choose $\lambda=1/2$, so that $C^{\mathrm{e}}_{\lambda,h}(P,E)$ is the arithmetic average of two e-values $h(P)$ and $E$, although this choice has no special advantage, since the positions of $h(P)$ and $E$ are not symmetric.

Validity of these combiners is straightforward to verify. 
Admissibility of the above combiners is obtained in the following result. 

 \begin{theorem} \label{th:combiners}
For  an admissible calibrator $h $ and $\lambda \in(0,1)$,
 $C_h^{\rm ie}$ is an admissible i-pe/e combiner, 
 $C^{\rm ip}$ is an admissible i-pe/p combiner, 
 $C_{\lambda,h}^{\rm e}$ is an admissible pe/e combiner,   
 and 
 $C^{\rm p}$ is an admissible pe/p combiner.
 \end{theorem}  
 
The most useful combiner is  the i-pe/p combiner $C^{\rm ip}$ that produces a p-value based on independent $P$ and $E$. 
We provide a proof for the statement on $C^{\rm ip}$ in  Theorem~\ref{th:combiners} and 
omit the rest of the proof.

\begin{proof}[Proof of  the {[ip]} case  in Theorem~\ref{th:combiners}.]

Checking that  $C^{\rm ip}$ is a valid i-pe/p combiner is simple (seen already in  Corollary~\ref{coro:ip-comb} in Section~\ref{sec:2-rand}): For $t\in (0,1)$ and any  independent pair $(P,E)$ of p-variable and e-variable for $\p$, we have 
\begin{align*}
\p ( C^{\rm ip}(P,E)\le t) &= \p(P/E\le t)
\\&= \p(P\le tE) =\E^\p[\p(P\le tE  \,|\, E) ] \le  \E^\p [tE]\le t. 
\end{align*}
  To show its admissibility,  
suppose for the purpose of contradiction that an i-pe/p combiner $f$ satisfies 
$f\le C^{\rm ip}$ and  
$f (p,e)< (p/e)\wedge 1$ for some $(p,e)\in [0,1]\times [0,\infty)$. 
Since $a\mapsto f(p,a)$ is decreasing,  
we can assume  $e\in [p,\infty)$ by replacing $e$ with $p$ if $e<p$. Take $P\lawis \mathrm {U}[0,1]$ and  $E$ distributed as 
 $\p( E  =e)=1/e =1- \p(E =0)$ if $e\ge 1$
and $\p(E   =e)=1/2 =  \p(E =2-e)$ if $e< 1$.
It is clear that $\E^\p[E ]=1$. 

Since $q\mapsto f(q,e)$ is increasing, 
there exists $p'<p$ such that 
$f (q,e) \le f(p,e) < p'/e$ for all $q\in [0,p]$. 
This gives
$$\p(f(P,e)  \le p'/e)\ge  \p(P\le p)=p.$$ 
For $\alpha = p'/e \in (0,1)$,  if $e\ge 1$, then 
\begin{align*}
\p(f(P,E)  \le \alpha )& = 
\p(f(P,e)\le p'/e )e^{-1}   + \p(f(P,0) \le \alpha)(1-e^{-1}) \\
  &\ge p /e >\alpha.
\end{align*}
If $e<1$, then 
\begin{align*}
\p(f(P,E)  \le \alpha )& = \frac 12   
\p(f(P,e)\le p'/e ) + \frac 12   \p(f(P,2-e) \le \alpha) 
\\ & \ge  \frac 12 p +   \frac 12   \p( P  \le \alpha (2-e))  
\\& =  \frac{1}{2} \left(p + (2-e) p'/e\right)  =\frac{1} {2} (p-p') +   p'/e 
>\alpha.
\end{align*}
Hence, $f(P,E)$  is not a p-value, and this contradicts the fact that $f$ is an i-pe/p combiner. 
This contradiction shows that $C^{\rm ip}$ is admissible.
\end{proof}

We will call  $C^{\rm ip}$ the \emph{quotient combiner} because it outputs the quotient of $P$ and $E$ (capped at $1$). 
This combiner typically leads to more powerful procedures compared to the other combiners and provides the foundation for our insight that e-values can act as unnormalized weights in multiple testing  in Section~\ref{sec:10-weighted}. 

\begin{example}[Randomization]

Consider an e-variable  $E$ for $\p$ and generate an independent uniform variable $U \sim \mathrm{U}[0,1]$. Then, using the quotient combiner, $U/E$ is a p-variable that satisfies $   U/E  \le 1/E$.
This implies that the unique admissible e-to-p calibrator $f:e \mapsto (1/e)\wedge 1$  can be improved if we allow for randomization, as studied already in 
Section~\ref{sec:2-rand}.
Although $ U/E $  may not be  practical in general due to external randomization, it becomes practical when applied  to a p-value $P$ (computed from data) that is independent of $E$.
\end{example}

The quotient combiner is a useful general-purpose method for meta-analysis from two independent datasets. The quotient combiner is immediately applicable when the researcher summarizes the first dataset as a single p-value, and the second dataset as a single e-value. 
Such a situation could occur when the second dataset is collected in such a way, e.g., with optional stopping and continuation, that inference is more natural with e-values. It could also be the case that one dataset comprises of a large sample size, allowing for asymptotic approximations to compute p-values, while the second dataset is smaller and may require finite-sample inference methods.

The quotient combiner is also useful when the above data constraints are not in place and the researcher can in principle compute both a p-value $P'$ and an e-value $E$ on the second dataset, both of which are independent of the p-value $P$ computed on the first dataset. In that case, the researcher could apply a p-value combination method based on $P$ and $P'$, e.g., Fisher's combination  $P_{\rm F}:=1- \chi_4(-2\log ( P P'))$, 
where $\chi_4$ is the cdf of the chi-square distribution with $4$ degrees of freedom. However, the researcher may still prefer to proceed with the quotient combiner $P/E$ (we omit the cap at $1$ for simplicity, which does not affect testing results). We suggest the following rule of thumb:

\emph{The Fisher combination is preferable to the quotient combiner under dataset exchangeability:} Suppose that the analyst considers the two datasets as a-priori exchangeable. In that case, it may be undesirable to use an asymmetric combination rule such as $P/E$, and Fisher's combination $P_{\rm F}$ is preferable on conceptual grounds. If the two datasets are also exchangeable in terms of their statistical properties (i.e., they have similar power), then $P_{\rm F}$ will typically have higher power than $P/E$. 

\emph{The quotient combiner is preferable to the Fisher combination for imbalanced datasets:} When one dataset (the ``primary'' dataset) is substantially more well-powered (larger anticipated signal or sample size) than the secondary dataset, and the investigator knows which dataset is more well-powered, then the $P/E$ combiner can often outperform Fisher's combination test in terms of power. A proviso is that the p-value is computed on the primary (more well-powered) dataset and the e-value on the secondary dataset.

\begin{figure}
\centering
\includegraphics[width=\linewidth]{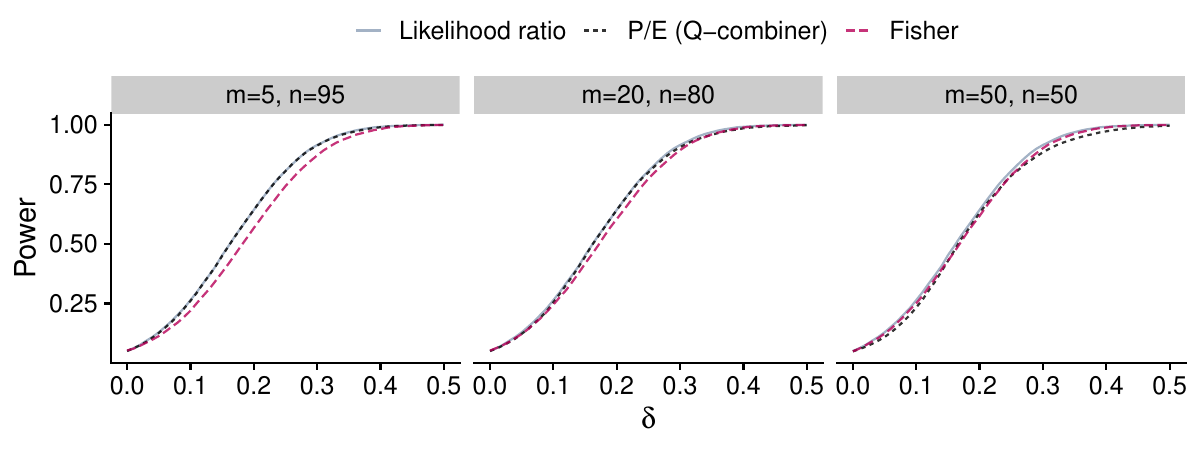}
\caption{\textbf{Simulation study for a meta-analysis combining two samples}: We compare the likelihood ratio test, the quotient combiner (Q-combiner), and the Fisher combination test, plotting power against signal strength $\delta$. The panels correspond to different choices of the two sample sizes $m$ and $n$. The quotient combiner is visibly more powerful on the left (matching the likelihood ratio), and Fisher's combination is marginally better on the right.}
\label{fig:0}
\end{figure}

 \begin{example}[A stylized example illustrating the rule of thumb]
 \label{ex:ep-or-pe}
Suppose we have two independent samples of iid  data points, $\mathbf X=(X_1,\dots,X_m)$  and $\mathbf Y=(Y_1,\dots,Y_n)$, both from a distribution $\p$, where  $n\ge m \ge 1$. 
We seek to test $H_0: \p = \mathrm{N}(0,1)$ 
against $H_1: \p=\mathrm{N}(  \delta,1)$,
where $\delta >  0$ is known.
The optimal (Neyman--Pearson) p-value based on $\mathbf X$ is given by $P_{\mathbf{X}} := 1-\Phi(T_{\mathbf X})$,
 where $\Phi$  is the standard normal distribution function and $T_{\mathbf X}:=m^{-1/2} \sum_{i=1}^m X_i$. Analogously we may compute p-values $P_{\mathbf{Y}}$ based on $\mathbf {Y}$, as well as $P_{\mathbf{Z}}$, where $\mathbf Z=(\mathbf X,\mathbf Y)$ is the full dataset. The optimal e-value $E_{\mathbf{X}}$ based on $\mathbf{X}$ is the likelihood ratio of $\mathrm{N}(\delta,1)^m$ over $\mathrm{N}(0,1)^m$.  

By the Neyman--Pearson Lemma, the most powerful test is the likelihood ratio test, which coincides with a test based on the p-value $P_{\mathbf{Z}}$. We compare $P_{\mathbf{Z}}$ against the quotient combiner $P_E :=P_{\mathbf{Y}}/E_{\mathbf{X}}$ 
and $P_{\mathrm F}$
by considering the hypothesis tests that reject $H_0$ when the p-value is smaller than or equal to some $\alpha > 0$.
Omitting details, we report the following conclusions.

 (i) Theoretically, using the notion of asymptotic relative efficiency, we can prove that $P_E$ is  as efficient as $P_{\mathbf Z}$ in two asymptotically settings: (a) $\alpha \to 0$, that is, when the type-I error is very stringent, and (b) $m/n \to 0$, that is, when $\mathbf{Y}$ is substantially more informative than $\mathbf{X}.$ 

(ii)  A small simulation study (Fig.~\ref{fig:0}) yields the following observations: (a)
 if $m/n$ is small, then $P_E$ has almost the same power as $P_{\mathbf Z}$, and both outperform the $P_{\mathrm F}$; (b) when $m=n$, $P_{\mathrm F}$ has slightly more power than $P_E$, and both have slightly less power than the $P_{\mathbf Z}$.

 These observations also confirm that, when using one p-value and one e-value, the primary dataset should  be used to generate the p-value. 
 \end{example}
 


\section{E-values as weights in p-value-based multiple testing}
\label{sec:10-weighted}

In this section, we study how  the cross-merging methods in Section~\ref{sec:10-cross} can be useful in the context of multiple testing in Chapter~\ref{chap:compound}.

The basic setting is that a vector $\mathbf P=(P_1,\dots,P_K)$ of p-values and a vector $\mathbf E =(E_1,\dots,E_K)$ of e-values are available for the hypotheses $H_1,\dots,H_K$. More precisely, $P_k$ is a p-variable for $H_k$, and $E_k$ is an e-variable for $H_k$ for each $k\in \mathcal K$. 
A useful context could be that e-values are obtained from preliminary data analysis and p-values are 
obtained from follow-up confirmatory research on these hypotheses (recall from Example~\ref{ex:ep-or-pe} that, in the context of testing a single hypothesis, we would like to choose the primary dataset to generate the p-value). 
We allow for some hypotheses with missing p-value or e-value, and in these cases $P_k$ and $E_k$ are set to $1$. The interpretation is that these hypotheses miss either a preliminary or a follow-up analysis.

\begin{definition}[E-weighted p-value procedure (ep-$\mathcal{D}$)]
\label{definition:eweightedpvalue}
Let p-$\cD$ be a multiple testing procedure based on p-values. We define the e-weighted p-value procedure ep-$\cD$ which proceeds as follows: for $k\in \mathcal K$, compute the quotient combiner $P^*_k := (P_k/E_k)\wedge 1$,
and then supply $\mathbf P^*$ to p-$\mathcal{D}$. 
\end{definition}

We mainly focus on the e-weighted version of the BH procedure in Section~\ref{sec:c8-pBH}, which we call the ep-BH procedure.   
Recall that the BH procedure at level $\alpha$ has valid FDR control in case the p-values are PRDS, as mentioned in Section~\ref{sec:c8-pBH}. 
In the e-BH procedure, we interpret the e-values as weights for the p-values. Intuitively, if $E_k>1$, then there is evidence against $H_k$ being a null, and we have $P_k/E_k<P_k$ (assuming $P_k\ne 0$), that is, the weight strengthens the signal of $P_k$. Conversely, if $E_k<1$, then there is no evidence against $H_k$ being a null, and  we have  $P_k/E_k>P_k$. 
The above interpretation of e-values as weights is quite natural, and the perspective is useful in deriving guarantees. Below we prove that ep-BH controls the FDR under the assumption that $\mathbf{P}$ is PRDS
and $\mathbf P$ and $\mathbf E$ are independent. Note that these two conditions do not imply PRDS of $\mathbf{P}^*$,
and hence some arguments are needed to establish the FDR control of the ep-BH procedure. The following result integrates over the randomness in the weights, i.e., the e-values.

\begin{theorem}\label{th:ep-BH}
Suppose that  $\mathbf{P}$ is independent of $\mathbf{E}$ and that $\mathbf{P}$ satisfies PRDS (Definition~\ref{definition:prds}). Then, the ep-BH procedure has FDR at most $K_0\alpha/K$. 
\end{theorem}

\begin{proof}[Proof.]
Let ep-$\cD$ be the ep-BH procedure at level $\alpha$. 
Since $\mathbf{E}$ is independent of $\mathbf{P}$, 
conditional on $\mathbf{E}$, the ep-BH procedure becomes a weighted p-BH procedure with weight vector 
$\mathbf{E}$ applied to the p-values $\mathbf{P}$ that are PRDS. Using well-known existing results on the false discovery rate of the weighted p-BH procedure (e.g., \citealp[Theorem 1]{RBWJ19}),  
we get  
 $$
\E\left [\frac{F_{\text{ep-}\cD}}{R_{\text{ep-}\cD} }\;\Big |\; \mathbf{E}\right ] \le \frac{1}{K} \sum_{k\in \mathcal N} E_k \alpha.
 $$ 
Hence, by iterated expectation,
the FDR of $ \text{ep-}\cD $
is guaranteed to be at most $\E[\sum_{k\in \mathcal N} E_k \alpha/K]\leq K_0 \alpha/K.$  
\end{proof}

Compared with the weighted BH procedure using constant weights, the ep-BH procedure uses weights that do not necessarily add up to $K$, i.e., they are not normalized. This allows us to save evidence obtained from the preliminary analysis using e-values. 
Moreover, the above result does not require any assumption about the dependence within $\mathbf{E}$.

 There are many other types of procedures that take p-values as input, and they can be generalized with e-values as weights.

\section*{Bibliographical note}

The optimal transport duality statements in Lemma~\ref{lem:duality} can be verified with several results in \cite{kellerer1984duality}; the formula \eqref{eq:duality-eq}
can be found in \citet[Theorem 2.1.1]{rachev2006mass} and \citet[Theorem 2.3]{R13}.
A standard textbook on optimal transport theory is \cite{villani2009optimal}.
\cite{Vovk/Wang:2021} identified all symmetric  admissible e-merging functions, but not the non-symmetric ones. 
Theorem~\ref{th:adm-merg} as presented in Section~\ref{sec:10-emerg} was obtained by \cite{wang2024only}.

The content of Sections~\ref{sec:10-cross}--\ref{sec:10-weighted} is based on \cite{ignatiadis2024values}.
The full proof of Theorem~\ref{th:combiners}  can be found in \citet[Appendix S1]{ignatiadis2024values}.  Using deterministic weights for the BH procedure to control FDR was studied by \cite{benjamini1997multiple}. The e-weighted  BH procedure in Section~\ref{sec:10-weighted} was introduced by \cite{ignatiadis2024values}, along with results on many other  procedures   weighted by e-values.

\chapter{Using e-values to combine dependent p-values}
\label{chap:pmerge}

In this chapter we study merging functions of p-values, and show how e-values appear as an essential step in such merging functions. 
The simplest form of merging is to take an average, which we study first. We turn to the more general deterministic methods for arbitrarily dependent p-values, followed by  methods based on the assumption of exchangeability among p-values, and methods based on randomization. 

As in Chapter~\ref{chap:combining}, in most parts of this chapter, we consider a fixed positive integer $K\ge 2$ and an atomless probability measure $\p$  on $(\Omega,\cF)$ without further specification. 
In a few instances, we mention ``for all $\p\in \cP$'' or ``for all $K\ge 2$'',
and thus temporarily allowing $\p$ and $K$ to vary, but it should be clear from context.

\section{The arithmetic mean of p-values}
\label{sec:c12-average-p}


\subsection{Average p-values}

Intuitively, an average p-value   should be the arithmetic mean of several p-values.
The distribution of such an arithmetic mean may depend on how the underlying null hypothesis $\cP$ is specified, for instance, simple or composite. We would like to define an average p-value   based on its distribution    under each $\p$, as in the cases of p-values and e-values.
This issue is clarified in the next proposition.  We will use the notation $P$ for average p-variables, as they share many properties with p-variables.




\begin{proposition}
\label{prop:c12-avp-eq}
 For a random variable $P$ and an atomless probability measure $\p$, the following are equivalent:
    \begin{enumerate}[label=(\roman*)] 
        \item $P$ is the arithmetic mean of finitely many p-variables for $\p$; 
        \item $P$ is a (possibly continuous) mixture of p-variables for $\p$; 
        \item $P\ge \E^\p[U|f(U)]$ for some measurable function $f$ and a random variable $U$ whose distribution is standard uniform under $\p$;
        \item  $\E^\p[g(P)] \ge \int_0^1 g(u)\d u  $ 
        for all increasing concave functions $g:\R\to \R$.  
        \item $\E^\p[(v-P)_+] \le v^2/2 $ for all $v\in (0,1)$.
    \end{enumerate}
Without assuming that $\p$ is atomless, we have (i)$\Rightarrow$(ii)$\Rightarrow$(iv)$\Leftrightarrow$(v).
\end{proposition}

\begin{proof}[Proof.]
We will 
show (i)$\Rightarrow$(ii)$\Rightarrow$(iv)$\Leftrightarrow$(v) without assuming that $\p$ is atomless
and (iii)$\Rightarrow$(iv).
The two remaining directions (v)$\Rightarrow$(i) and 
(v)$\Rightarrow$(iii) involve advanced results and are briefly explained.

(i)$\Rightarrow$(ii) follows by definition.
(ii)$\Rightarrow$(iv) follows from Jensen's inequality: For a probability measure $\nu$ on a space $\Theta$ and a set $\{P_{\theta}: \theta \in \Theta\}$ of p-variables, 
we have  
$$
\E^\p\left[g\left(\int_\Theta P_\theta \nu(\d \theta)\right)\right]
\ge\E^\p\left[ \int_\Theta  g(P_\theta)  \nu(\d \theta)\right]
\ge \E^\p[g(U)],
$$ 
where we interchanged the expectation and the integral by Fubini's theorem, noting that $g$ is bounded from below on $[0,\infty)$. 
(iv)$\Rightarrow$(v) is straightforward by taking $g(x)=-(v-x)_+$, and (v)$\Rightarrow$(iv) is due to the well-known fact that any increasing concave functions on $(0,1)$ belong to  the convex hull of $-(v-x)_++b$ for $v\in \R$ and $b\in \R$. 
(iii)$\Rightarrow$(iv)  
follows from conditional Jensen's inequality:
$$
\E^\p[g(P)]
\ge \E^\p[g(\E^\p[U|f(U)])] 
\ge \E^\p[g(U)]. 
$$  

The implication (v)$\Rightarrow$(i) relies on the following nontrivial fact: 
the average of three standard uniform random variables can attain any distribution satisfying (v) with mean $1/2$. The implication
(v)$\Rightarrow$(iii)  relies on a refinement of Strassen's theorem. 
For these two directions, see the bibliographic note.
\end{proof}

Proposition~\ref{prop:c12-avp-eq} leads to the following definition of average p-values purely based on their distributions under the null.

\begin{definition}[Average p-values]
\label{def:avg-p2}
A random variable $P$ is an 
\emph{average p-variable} for $\cP$ if $\E^\p[(v-P)_+] \le v^2/2 $ for all $v\in (0,1)$ for all $\p\in \cP$. 
\end{definition}

It is straightforward to check that a standard uniform random variable satisfies each  condition  in Proposition~\ref{prop:c12-avp-eq}, and hence a p-variable is also an average p-variable.  
Moreover, for any p-variables $P_1,\dots,P_K$ for $\cP$, 
their arithmetic mean is an average p-variable, regardless of whether $\cP$ is simple or composite. 
Similarly to the case of p-variables, 
a truncation of an average p-variable at $1$ is still an average p-variable.

\begin{remark}
Condition (iv) of  Proposition~\ref{prop:c12-avp-eq}  is known as second-order stochastic dominance in economics. Indeed, p-variables can be defined using the so-called first-order stochastic dominance; that is,  $P$ is a p-variable for $\p$ if and only if $\E^\p[g(P)]\ge \E^\p[g(U)]$ for all increasing functions $g:\R\to \R$ (this was briefly mentioned in Section~\ref{sec:1-def-E-P-test}). Clearly, this condition is stronger than condition (iv).
\end{remark}

All results in this and the next few subsections hold for both simple hypotheses and composite hypotheses, and we will focus on the case of a fixed atomless $\p$ from now on.

As a simple consequence of Proposition~\ref{prop:c12-avp-eq},  the set of average p-variables for  $\p$ is convex and closed under convergence in distribution, whereas the set of p-variables is not convex. 
Indeed, the set of  average p-variables  is the convex hull of 
 the set of p-variables.


Although a p-variable is also an average p-variable,  the converse is not true. Indeed, even  the arithmetic mean of iid   standard uniform p-variables
is not a p-variable, as it is
asymptotically concentrated around $1/2$ (by the law of large numbers), conflicting the definition of a p-variable.
It turns out that a factor of $2$ 
will convert an average p-variable into a p-variable. 
\begin{proposition}\label{prop:2avg}
    An average p-variable is  not necessarily a p-variable, 
     but twice an average p-variable is a p-variable. 
     Moreover, 
the factor $2$ cannot be improved; that is, for any $\gamma\in [0,2)$, 
$\gamma P$ is not a p-variable for some average p-variable $P$.
\end{proposition}


\begin{proof}[Proof.] 
The first statement is explained above. 
For the second statement,
using Proposition~\ref{prop:c12-avp-eq}, 
we can write an average p-variable $P$ for $\p$ as $P=(P_1+\dots+P_K)/K$, where  $P_1,\dots,P_K$  are  p-variables  for $\p$ and $K\ge 2$.
Take any $\alpha \in(0,1)$. 
Note that $p\mapsto (2-2p)_+$  is a calibrator.
By Proposition~\ref{prop:rescale-cal}, $p\mapsto \gamma (2-2 \gamma p)_+$
is a calibrator for any $\gamma \ge 1$.
This implies that $E_k:=(2-2P_k/\alpha)_+ /\alpha $
is an e-variable for each $k\in \cK$ and $\alpha \in (0,1)$.
We have 
\begin{align*}
\p\left(\frac{2(P_1+\dots+P_K)}{K} \le \alpha \right)  &= 
\p\left(\frac{1}{K}\sum_{k=1}^K \left(2- \frac{2P_k}{\alpha}\right)  \ge 1\right)
\\ &=
\p\left(\frac{1}{K}\sum_{k=1}^K \frac{1}{\alpha} \left(2- \frac{2P_k}{\alpha}\right)  \ge \frac{1}{\alpha}\right)
\\ &\le 
\p\left(\frac{1}{K}\sum_{k=1}^K   E_k \ge \frac{1}{\alpha}\right)\le \alpha,
\end{align*}
where the last inequality is due to Markov's inequality.

For the last statement,  take  a random variable $U$ that is uniformly distributed on $[0,1]$ under $\p$. 
Clearly, $U$ and $1-U$ are p-variables, and hence $1/2$ is an average p-variable. The constant $\gamma/2$ is not a p-variable for any $\gamma<2$. 
\end{proof}

  Proposition~\ref{prop:2avg} implies that for any collection of p-variables,
twice their arithmetic mean is a p-variable.
Similarly, twice their median is a p-variable, which we will study later (see Example~\ref{ex:order-fam}).

\subsection{Special cases: mid p-values and posterior predictive p-values}
 
In addition to the natural appearance as the arithmetic mean of p-values, another useful class of average p-values 
is that of \emph{mid p-values}, which are computed based on discrete test statistics. 
We first recall the usual practice to obtain p-values for a simple hypothesis $\p$, mentioned in Fact~\ref{fact:basic} in Section~\ref{sec:1-def-E-P-test}. 
Let $T$ be a test statistic which is a function of the observed data. 
Here, a smaller value of $T$ represents stronger evidence against the null hypothesis. 
The p-variable $P$ is usually computed from the conditional probability  \begin{equation}\label{eq:exceedance1} 
P= \p(T'\le T  |T)=F(T),\end{equation} 
where $F$ is the distribution of $T$, and $T'$ is an independent copy of $T$; here and below, a copy of  $T$ is a random variable identically distributed as $T$.  

If $T$ has a continuous distribution under $\p$, then the p-variable $P$ defined by \eqref{eq:exceedance1} has a standard uniform distribution. 
If the test statistic $T$ is discrete, e.g., when testing a binomial model, 
$P$ is strictly stochastically larger than a uniform random variable on $[0,1]$. 

The discreteness  of $T$   leads to a   conservative p-value; in particular, $\E[P]>1/2$. One way to address this is to randomize the  p-value to make it uniform on $[0,1]$; however, randomization is generally undesirable in testing.  
As the most natural alternative, \emph{mid p-values}  arise in the presence of discrete test statistics. 
For the test statistic $T$,  the mid p-value is given by 
 \begin{equation}\label{eq:mid}
P_T=  \frac{1} 2\p(T'\le T|T)+\frac{1}2 \p(T'<T|T)  =  \frac 12 F(T-)+ \frac 12 F(T),
\end{equation}
where $F(t-)=\lim_{s\uparrow t} F(t)$ for $t\in \R$.
 Clearly, $  P_T\le  F(T)$ and $\E[  P_T ] =  1/2$.  
If $T$ is continuously distributed, then $P_T=F(T)$ is uniform on $[0,1]$. 
In case $T$ is discrete, $P_T$ is not a p-variable. 
In Figure~\ref{fig:r1-2} we present some examples of quantile functions of p-, mid p- and average p-variables.

\begin{figure}[t]    
 \begin{subfigure}[b]{0.31\textwidth}
         \centering
\begin{tikzpicture}
\draw[<->] (0,3.4) -- (0,0) -- (3.4,0);
\draw[gray,dotted] (0,0) -- (3,3);
\draw[gray,dotted] (3,0) -- (3,3);
\draw[gray,dotted] (0,3) -- (3,3); 
\node[below] at (3,0) {$1$}; 
\node[left] at (0,3) {$1$}; 
\draw (0,0.12) -- (0.24,0.12);
\draw (0.24,0.42) -- (0.6,0.42);
\draw (0.6,0.8) -- (1.0,0.8);
\draw (1.0,1.2) -- (1.4,1.2);
\draw (1.4,1.8) -- (2.2,1.8);
\draw (2.2,2.6) -- (3,2.6); 
\draw[red, thick, densely dotted] (0,0.24) -- (0.24,0.24);
\draw[red, thick, densely dotted] (0.24,0.6) -- (0.6,0.6);
\draw[red, thick, densely dotted] (0.6,1.0) -- (1.0,1.0);
\draw[red, thick, densely dotted] (1.0,1.4) -- (1.4,1.4);
\draw[red, thick, densely dotted] (1.4,2.2) -- (2.2,2.2);
\draw[red, thick, densely dotted] (2.2,3) -- (3,3); 
\end{tikzpicture} 
         \caption{mid p-variable (black) and p-variable (red dotted) for a discrete statistic}
         \label{fig:r1-2-1}
     \end{subfigure}
     ~~
      \begin{subfigure}[b]{0.31\textwidth}
         \centering
\begin{tikzpicture}
\draw[<->] (0,3.4) -- (0,0) -- (3.4,0);
\draw[gray,dotted] (0,0) -- (3,3);
\draw[gray,dotted] (3,0) -- (3,3);
\draw[gray,dotted] (0,3) -- (3,3); 
\node[below] at (3,0) {$1$}; 
\node[left] at (0,3) {$1$};  
\draw (0,0.2) -- (0.4,0.2);
\draw (0.4,0.4) --  (1.2,1.2);
\draw (1.2,1.5) -- (1.8,1.5);
\draw (1.8,1.8) -- (2.2,2.2);
\draw (2.2,2.6) -- (3,2.6); 
\draw[red, thick, densely dotted] (0,0.4) -- (0.4,0.4); 
\draw[red, thick, densely dotted] (0.4,0.4) --  (1.2,1.2);
\draw[red, thick, densely dotted] (1.2,1.8) -- (1.8,1.8); 
\draw[red, thick, densely dotted] (1.8,1.8) -- (2.2,2.2);
\draw[red, thick, densely dotted] (2.2,3) -- (3,3); 
\end{tikzpicture} 
         \caption{mid p-variable (black) and p-variable (red dotted) for a  hybrid statistic}
        \label{fig:r1-2-2}
     \end{subfigure}
     ~~
      \begin{subfigure}[b]{0.31\textwidth}
         \centering
\begin{tikzpicture}
\draw[<->] (0,3.4) -- (0,0) -- (3.4,0);
\draw[gray,dotted] (0,0) -- (3,3);
\draw[gray,dotted] (3,0) -- (3,3);
\draw[gray,dotted] (0,3) -- (3,3); 
\node[below] at (3,0) {$1$}; 
\node[left] at (0,3) {$1$}; 
\draw [black] (0,0) .. controls (1.2,2.6) and (2.3,0.8).. (3,3);
\end{tikzpicture} 
         \caption{average p-variable (not a mid p-variable)\\~}
          \label{fig:r1-2-3}
     \end{subfigure}
 
\caption{Examples of quantile functions of p-, mid p- and average p-variables}
\label{fig:r1-2}
\end{figure}

A mid p-value is an average p-value for any atomless probability measure $\p$.
Suppose that there exists  a standard uniform random variable $V$   independent of $T$. 
Let 
$U=F(T-)+ V(F(T)-F(T-))
$, which is uniformly distributed on $[0,1]$. 
It follows that 
\begin{align*}
    \frac 12 F(T-)+ \frac 12 F(T)
    = \E^\p[U |T],
\end{align*}
and then we conclude by Jensen's inequality that a mid p-value is an average p-value.
(The above $U$ exists even without assuming the existence of $V$ independent of $T$; see Lemma~\ref{lem:quantile2}.)

Another example of average p-values is the class of posterior predictive p-values in   Bayesian analysis, which we explain briefly below.

Let $X$ be a vector of data. The null hypothesis $H_0$ is given by 
$\{\psi\in \Psi_0\}$ where $\Psi_0$ is a subset of the parameter space $\Psi$ on which a prior distribution is specified. 
The posterior predictive p-value   is defined as the realization of the random variable
$$
P_B:= \p(D( X',\psi) \ge D( X,\psi) \mid  X),
$$
where $D$ is a function (taking a similar role as  test statistics), $ X'$ and $ X$  are iid conditional on $\psi$, and the probability is computed under the joint posterior distribution of $( X',\psi)$. 
Note that $P_B$ can be rewritten as 
$$
P_B=\int  \p(D( X',y) \ge D( X,y) \mid  X,y)  \Pi(\d y| X),
$$
where $\Pi$ is the posterior distribution of $\psi$ given the data $  X$. 
The fact that $P_B$ is an average p-variable follows from Proposition~\ref{prop:c12-avp-eq}, condition (ii).  In this context, average p-variables are obtained by integrating p-variables over the posterior distribution of some unobservable parameter.  We omit a detailed discussion of the Bayesian interpretation, which is not the focus on the book, but simply make the point that average p-values have a natural appearance in the Bayesian context.

\subsection{Calibration between average p-values and e-values}

In the next result, we show that average p-values can be calibrated into an e-value through p-to-e calibrators in Section~\ref{sec:2-cali}, under the additional assumption that the calibrator is convex.

\begin{theorem}
\label{th:cal-avg}
For an average p-variable $P$, 
$f(P)$ is an e-variable for  any calibrator $f$ that is convex on $[0,1]$.
For any e-variable $E$, 
$(2E)^{-1}\wedge 1$ is an average p-variable. 
\end{theorem}
\begin{proof}[Proof.]
To show the first statement,  note that it suffices to consider p-variables $P_1,\dots,P_K$ taking values in $[0,1]$. Let $P=(P_1+\dots+P_K)/K$ and  $E=f(P)$. By the convexity of the calibrator $f$, $$\E^\p[E]  \le 
   \E^\p\left[\frac 1 K \sum_{k=1}^K f (P_k)\right]
   \le 1,
    $$
    and  hence $E$ is an e-variable.

Next, we show the second  statement. 
Clearly, it suffices to show that $1/(2E)$ is an average p-variable for any e-variable $E$ with mean $1$, since any e-variable with mean less than $1$ is dominated by an e-variable with mean $1$. 
Let $\delta_x$ be the point mass at $x\in \R$. 

Assume that $E$ has a two-point distribution (including the  point mass $\delta_1$). With $\E^\p[E]=1$, the distribution of $E$ can be characterized with two parameters $p\in (0,1) $ and $a\in [0,1/p]$ as
$$
 p \delta_{1+(1-p)a}  + (1-p)  \delta_{1-pa}.
$$ Write $x_1= 1/(2+2(1-p)a)\in (0,1/2]$
and $x_2= 1/(2-2pa)\in [1/2,\infty)$.
The distribution   of $P=1/(2E)$ (we allow $P$ to take the value $\infty$ in case $a=1/p$, which is harmless) is given by
$$  p \delta_{x_1}  + (1-p)  \delta_{x_2}.$$
To show that $P$ is an average p-variable, let  $U \lawis \mathrm{U}[0,1]$ and $V=\E^\p[U|\id_{\{U\le x_1 \}}]$. We will show $\E^\p[g(P)] 
\ge \E^\p[g(V)] \ge \E^\p[g(U)]  
$ for all increasing concave functions $g$. 
The second inequality follows from Jensen's inequality, and it suffices to prove the first one. 
The distribution of $V $
is given by
$$
 2x_1 \delta_{x_1}  + (1-2x_1)  \delta_{x_3},
$$
where $x_3=x_1+1/2$. 
Note that  $2 x_1 \ge 1/(2-p) \ge p$.  
Moreover, by Jensen's inequality,
$$\E^\p[P]= \E^\p\left[\frac{1}{2E}\right] 
\ge \frac{1}{2\E^\p[E]}=\frac 12.
$$
Note that
\begin{align*}
\E^\p[g(P)] - \E^\p[g(V)]
= (1-p) g(x_2) - (2x_1-p) g(x_1) - (1-2x_1) g(x_3). 
\end{align*}
Since $\E^\p[P]\ge \E^\p[V]=1/2$, we know that $$\frac{2x_1-p}{1-p} \delta_{x_1} - \frac{ 1-2x_1  }{1-p}\delta_{x_3}$$ has a mean at most $x_2$.
Jensen's inequality then gives $\E^\p[g(P)] - \E^\p[g(V)]\ge 0.$ 
Hence, $P$ is an average p-variable. 

We use two other known facts that can be checked directly:
(a) For a general e-variable $E$ with mean $1$, its distribution can be rewritten as a mixture of two-point distributions with mean $1$. 
(b) The set of average p-variables is closed under  mixtures of distributions.
Using these two facts, we conclude that $f(E)$ is an average p-variable. 
\end{proof}

Convex calibrators are quite common. 
Indeed, all calibrators   \eqref{eq:cal1}--\eqref{eq:cal4} in Section~\ref{sec:2-cali} are   convex on $[0,1]$.
Practically, average p-values can be calibrated into e-values using any calibrators introduced before, except for the all-or-nothing one $p\mapsto \id_{\{p\le \alpha\}}/\alpha$ for some $\alpha \in (0,1)$.

 Figure~\ref{fig:1} summarizes the calibrators among p-values, average p-values and e-values, including the ones  in Section~\ref{sec:2-cali}.
 As implied by Theorem~\ref{th:cal-avg}, 
 the calibrator $e\mapsto (2e)^{-1} $ to an average p-value has an extra factor of $1/2$ compared to the standard calibrator $e\mapsto e^{-1}$ to a p-value. 
A composition of $e\mapsto (2e)^{-1}\wedge 1$ (from an e-value to an average p-value) 
and   
$p\mapsto (2 p)\wedge 1$
(from an average p-value to a p-value)
leads to  the unique admissible e-to-p calibration $e\mapsto e^{-1}\wedge1$ in Section~\ref{sec:2-cali}, thus showing that average p-values serve as an intermediate object between e-values and p-values.

\begin{figure}[t]
\begin{center}
\tikzstyle{bag} = [text width=5em, text centered]
\tikzstyle{end} = []
  \begin{tikzpicture}[
      mycircle/.style={
         circle,
         draw=black,
         fill=gray,
         fill opacity = 0.1,
         text opacity=1,
         inner sep=0pt,
         minimum size=40pt,
         font=\small},
               myinvisiblecircle/.style={
         circle,
         draw=black ,
         fill=gray,
         fill opacity = 0.1,
         text opacity=1,
         inner sep=0pt,
         minimum size=55pt,
         font=\small},
      myarrow/.style={-Stealth},
      node distance=1.5cm and 3.5cm
      ] 
      \node[myinvisiblecircle]  (c1) {}; 
      \node at (0,-0.26) {\small average};
           \node at (0,-0.53) {\small p-value};
                  \draw[   fill=gray,
         fill opacity = 0.04] (0,0.31) ellipse (0.85cm and 0.32cm);
       \node at (0,0.31) (c4) {\scriptsize mid p-value};
      \node[mycircle, below right=of c1] (c2) {e-value}; 
      \node[mycircle,below left=of c1] (c3) {p-value}; 

    \foreach \i/\j/\txt/\p in {c1.320/c2.140/{$e= f(p^*)$, $f$ convex}/below,
       c2.130/c1.330/{$p^*=  (2e)^{-1}\wedge 1$}/above,
      c1.200/c3.50/{$p=  (2p^*)\wedge 1$}/above,     
      c3.40/c1.210/$p^*=  p$/below,    
      c2.175/c3.5/{$p =  e^{-1}\wedge 1 $}/above,
      c3.355/c2.185/$e= f(p)$/below}
       \draw [myarrow] (\i) -- node[sloped,font=\small,\p] {\txt} (\j); 
    \end{tikzpicture}
    \end{center}
    \caption{Calibration among p-values, average p-values and e-values, where the calibrator $f:[0,1]\to [0,\infty]$ is   left-continuous and decreasing with $f(0)=\infty$ and $\int_0^1 f(t)\d t=1$. All these calibrators are admissible in the sense of Section~\ref{sec:2-cali}.}
    \label{fig:1}
    \end{figure}

\subsection{Tests based on average p-values}

For an average p-variable $P$,
the test $\id_{\{P\le \alpha\}}$
is not a level-$\alpha$ test in general, because $P$ is not necessarily a p-variable.
Nevertheless,  $\id_{\{P\le V\}}$ for a random variable $V$ with mean $\alpha$ independent of $P$ often gives a  level-$\alpha$ test,
as shown in the next result.
In what follows, independence is understood to hold for all $\p\in \cP$.

\begin{theorem}
\label{th:avg-p-random}
Let $\alpha \in (0,1)$, $V\ge 0$ be independent of $P$ such that $\E^\p[V\wedge 1] \le \alpha $ for each $\p\in \cP$.
\begin{enumerate}[label=(\roman*)]
\item  If $P$ is a p-variable, then 
$\id_{\{P \le V\}}$ is a level-$\alpha$ test  for $\cP$. 
    \item  If $P$ is an average p-variable and $V$ has a decreasing density function on $(0,\infty)$  under each $\p\in \cP$, then 
$\id_{\{P \le V\}}$ is a level-$\alpha$ test  for $\cP$.
\end{enumerate}
In particular, $\id_{\{P \le 2\alpha U\}}$ is a level-$\alpha$ test  for $\cP$, where $U $ is independent of $P$ and uniformly distributed on $[0,1]$  under each $\p\in \cP$.
\end{theorem}
\begin{proof}[Proof.]
Part (i) follows from $\E^\p[\id_{\{P \le V\}}]
=\E^\p[V\wedge 1] \le \alpha$ for $\p\in \cP$.

Part (ii) follows from the following argument. For $\p\in \cP$, since $V$ has a decreasing density,  the function $g:x\mapsto \p(V\ge x)$ is decreasing convex. Hence, using Proposition~\ref{prop:c12-avp-eq}, we have 
$$\E^\p[\id_{\{P \le V\}}]
= \E^\p[g(P)] \le \int_0^1 g(u)\d u = \E^\p[V\wedge 1]\le \alpha,
$$
as desired. 
\end{proof}

Theorem~\ref{th:avg-p-random} gives rise to randomized tests based on the arithmetic mean of several p-values regardless of their dependence structure; i.e.\
\begin{align}
    \label{eq:randomized-avg1}
\id_{\{(P_1+\dots+P_K)/K\le 2\alpha U\} }
\end{align} 
is a level-$\alpha$ test for any p-variables $P_1,\dots,P_K$ and $U\lawis \mathrm{U}[0,1]$ independent of them.
If randomization is not allowed, then we can use   $$
\id_{\{(P_1+\dots+P_K)/K\le 2\alpha  \}},
$$
which is weaker  than \eqref{eq:randomized-avg1} in terms of power.
This is a method of testing by combining arbitrarily dependent p-values. 
The next few sections are dedicated to a formal theory of general methods for combining p-values, including an exchangeable improvement  and another  
 randomized test based on the average p-value.

\section{Merging p-values under arbitrary  dependence}
\label{sec:10-2}

The structure of merging functions for p-values is much richer than that of e-merging functions. 
It turns out that admissible ways of merging p-values under arbitrary dependence have to go through merging e-values. 
At a high level, the reason why merging p-values relies on e-values is due to optimal transport duality in  Lemma~\ref{lem:duality}, which converts   constraints on probability into constraints on expectation. 

\subsection{P-merging functions and examples}
\label{sec:p-merging-examples}

We first give the definition of p-merging functions, which is analogous to e-merging functions in Definition~\ref{def:e-merg}.
 
 \begin{definition}[P-merging functions]
\label{def:p-merg}
\begin{enumerate}[label=(\roman*)]
\item 
 A  \emph{p-merging function} (of $K$ p-values) is  an increasing Borel function $F:[0,\infty)^K\to[0,\infty)$
such that for any hypothesis, 
$
  F(P_1,\dots,P_K) $ is a p-variable for any p-variables $P_1,\dots,P_K$. 
\item [(ii)] A p-merging function $F$ is \emph{symmetric} if $F(\mathbf p)$ is invariant under any permutation of $\mathbf p$. 
\item 
A p-merging function $F$ is \emph{homogeneous} if $F(\lambda\mathbf p)=\lambda F(\mathbf p)$ for all $\lambda \in (0,1]$ and $\mathbf p$ with $F(\mathbf p) \le 1$. 
\item  
  A p-merging function $F$ \emph{dominates} a p-merging function $G$ if $F\le G$.
The domination is \emph{strict}  
if $F\le  G$ and $F(\mathbf{p})<G(\mathbf{p})$ for some $\mathbf{p}\in[0,\infty)^K$.
A p-merging function   is \emph{admissible}
if it is not strictly dominated by any p-merging function. 
\end{enumerate} 
 \end{definition} 
 
 Similarly to the situation in Chapters~\ref{chap:multiple} and~\ref{chap:combining},  it suffices to consider p-variables for the fixed atomless probability measure $\p$. Let $\fU$ be the set of all p-variables for  $\p$. 
   A p-merging function $F$ is  \emph{precise} if
\[
  \sup_{\mathbf P\in \fU^K} \p(F (\mathbf P)\le \epsilon) = \epsilon \text{ for all }\epsilon \in (0,1).
\]
In other words, $\epsilon$ by $\epsilon$, $F$ attains the largest possible probability allowed for $F(\mathbf P)$ to be a p-value. 
 
 Note that p-values are more useful to be small, and hence domination between p-merging functions is defined via an opposite direction to that between e-merging functions.  
We will also speak of admissibility within smaller classes of p-merging functions, such as the class of symmetric p-merging functions.

All p-merging functions that we encounter in this chapter are homogeneous and symmetric.
Although we allow the domain of $F$ to be $[0,\infty)^K$ in order to simplify presentation,
the informative part of $F$ is its restriction to $[0,1]^K$.

 We have some specific notation for this section. Denote by
 $\mathbf 0$   the $K$-vector of zeros, and $\mathbf 1$   the $K$-vector of ones.
All vector inequalities and the operation $\wedge$ of taking the minimum of two vectors are component-wise.
For  $a,b,x,y\in \R$, $ax \wedge by $ should be understood as $ (ax)\wedge (by)$.

We collect some basic properties of admissible p-merging functions,
which will be useful in our analysis later.
We skip the proofs.

\begin{proposition}\label{prop:p-properties}
\begin{enumerate}[label=(\roman*)] 
  \item   Any admissible p-merging function is precise and lower semicontinuous,
  takes value $0$ on $[0,\infty)^K\setminus(0,\infty)^K$,
  and satisfies $F(\mathbf{p})=F(\mathbf{p}\wedge\mathbf{1})\wedge1$
  for all $\mathbf{p}\in[0,\infty)^K$.
  \item  The point-wise limit of a sequence of p-merging functions is a p-merging function.
  \item   Any p-merging function is dominated by an admissible p-merging function.
    \end{enumerate}
\end{proposition}

In what follows, we will always write $\mathbf p=(p_1,\dots,p_K)$. 
Two natural families of  p-merging functions are  the \emph{order-family}  based on order statistics 
and the \emph{mean-family} based on generalized mean. 
All functions that we see below are homogeneous and symmetric. 
  \begin{example}[Order-family]
  \label{ex:order-fam}
The order-family is parameterized by $k\in\cK$,
and its $k$-th element is the function 
\begin{equation}\label{eq:Ruger}
  G_{k,K}: \mathbf p
  \mapsto
  \frac{K}{k} 
  p_{(k)}
  \wedge
  1,
\end{equation}
where $p_{(k)}$ is the $k$-th ascending order statistic of $p_1,\dots,p_K$. One can verify that the functions $G_{k,K}$, $k\in \cK$ are precise p-merging functions; see also Example~\ref{ex:ex-th:e}.
  \end{example}
    \begin{example}[Mean-family]
The mean-family is parameterized by $r\in[-\infty,\infty]$,
and its element with index $r$ has the form
\begin{equation}\label{eq:merge1}
  F_{r,K}: \mathbf p
  \mapsto
  b_{r,K}\fM_{r,K}(\mathbf p)\wedge1,
\end{equation} 
where $\fM_{r,K}$ is given by
 \[
  \fM_{r,K}(\mathbf p)
  =
  \left(
    \frac{p_1^r + \dots + p_K^r }{K}
  \right)^{1/r}
 \]
and $b_{r,K}\ge1$ is a suitable constant making $F_{r,K}$ a p-merging function.
The arithmetic mean $\fM_{1,K}$ of p-values has been studied in Section~\ref{sec:c12-average-p}.

The average $\fM_{r,K}$ is also defined for $r\in\{0,\infty,-\infty\}$ as the limiting cases of \eqref{eq:merge1},
which correspond to the geometric mean, the maximum, and the minimum, respectively. 
Another useful member of the mean-family is the multiple $F_{-1,K}$ of the harmonic mean $\fM_{-1,K}$.  
  \end{example} 
    \begin{example}[Bonferroni correction and maximum] 
      \label{ex:c10-Bonf}
The initial and final elements of the M- and O-families coincide:
the initial element is the Bonferroni p-merging function
\begin{equation}\label{eq:Bonferroni}
  G_{1,K} = F_{-\infty,K}: \mathbf p
  \mapsto
  K\min(\mathbf p)
  \wedge
  1,
\end{equation}
and the final element is the maximum p-merging function
\[
  G_{K,K} = F_{\infty,K}:
  \mathbf p
  \mapsto
  \max(\mathbf p).
\] 
  \end{example}
  
  \begin{example}[Simes and Hommel functions]
  \label{ex:c10-Simes}
By minimizing over $k$ in the order-family, we get  the \emph{Simes function}
\[
  S_K := \bigwedge_{k=1}^K G_{k,K}  
\]
which 
is not a p-merging function, but it produces a p-variable if the input p-variables are independent or PRDS.  Moreover, one can show that the Simes function is the minimum of all symmetric p-merging functions. 
Multiplied by a constant 
$\ell_K := \sum_{k=1}^K k^{-1}$ representing a factor of penalization, 
we get a precise p-merging function
\[
  H_K := \left(\sum_{k=1}^K \frac 1k \right) \bigwedge_{k=1}^K G_{k,K} = \ell_K S_K,
\]    
 which we call the \emph{Hommel function}.
\end{example}

\subsection{Using e-values to merge p-values}
\label{sec:10-2-OT}

We first explain a simple way of merging p-values through e-values. 
Recall that a calibrator $f:[0,\infty)\to [0,\infty]$ transforms a p-variable to an e-variable, and they satisfy $\int_0^1 f(p)\d p\le 1$  and $f(p)=0$ for $p>1$.  
Let $\Delta_K$ be the standard $K$-simplex, that is,
$$\Delta_K := \{(\lambda_1,\dots,\lambda_K)\in[0,1]^K:\lambda_1+\dots+\lambda_K=1\}.$$ 
For any calibrators $f_1,\dots,f_K$ and any  $(\lambda_1,\dots,\lambda_K)\in\Delta_K$,
define the function 
\begin{equation}\label{eq:general}
  G(\mathbf{p}):
   =
  \lambda_1 f_1(p_1)
  + \dots +
  \lambda_K f_K(p_K).
\end{equation} 
For any vector $\mathbf P$ of p-variables,
the function $G$ produces an e-variable  $G(\mathbf P)$, since it is a convex combination of e-variables $f_1(P_1),\dots,f_K(P_K)$. 
Markov's inequality implies
$$
\p\left(G(\mathbf P)\ge \frac{1}{\epsilon} \right) \le \epsilon~~~\mbox{
for all $\epsilon \in (0,1)$.}
$$
Thus, $1/G$ is a p-merging function, and this can also be seen as applying the e-to-p calibrator $t\mapsto (1/t)\wedge 1$ in Proposition~\ref{prop:e-to-p} to the e-variable $G(\mathbf P)$.

Such a ``naive'' (but valid)  procedure  of merging p-values through merging e-values, producing $1/G$, is generally not admissible, as the last conversion step is wasteful. 
However,  perhaps surprisingly,
all admissible and homogeneous p-merging functions can be obtained via the above procedure,  choosing different functions $G$ for each $\epsilon$. 
This is shown in Theorem~\ref{th:e-merge-p} below, the main result of this section.

 To explain this result, we reformulate p-merging functions by their rejection regions.
The \emph{rejection region} of a p-merging function $F$ at level $\epsilon\in (0,1)$ is defined as
\begin{equation}\label{eq:region}
  R_\epsilon(F) := \left\{\mathbf p \in [0,\infty)^K: F(\mathbf p)\le \epsilon \right\}.
\end{equation}
If $F$ is homogeneous, then $R_\epsilon(F)$, $\epsilon\in(0,1)$,
takes the form $R_\epsilon(F)=\epsilon A$ for some $A\subseteq [0,\infty)^K$. 
Conversely, any increasing collection of Borel lower sets
$\{R_\epsilon \subseteq [0,\infty)^K: \epsilon \in (0,1)\}$
determines an increasing Borel function $F: [0,\infty)^K\to[0,1]$ by the equation
\begin{equation}\label{eq:10-inf}
  F(\mathbf p) = \inf\{\epsilon\in(0,1): \mathbf p\in R_{\epsilon}\},
\end{equation}
with the convention $\inf\varnothing = 1$.
It is immediate that $F$ is a p-merging function if and only if
$\p(\mathbf P \in R_\epsilon) \le \epsilon$ for all $\epsilon\in (0,1)$ and $\mathbf P\in \fU^K$.

\begin{theorem}\label{th:e-merge-p}

\begin{enumerate}[label=(\roman*)]
\item For any  p-merging function $F$ and any $\epsilon \in(0,1)$,  there exists a homogeneous p-merging function $G$
such that $R_\epsilon(F) \subseteq R_\epsilon(G)$. 
\item For any admissible p-merging function $F$ and any $\epsilon \in(0,1)$, 
  there exist $(\lambda_1,\dots,\lambda_K)\in\Delta_K$ and admissible calibrators $g_1  ,\dots,g_K  $
  such that
  \begin{equation}\label{eq:calibrator}
    R_\epsilon(F)
    \subseteq 
    \left\{
      \mathbf p \in [0,\infty)^K:
      \sum_{k=1}^K \lambda_k g_k (p_k) \ge \frac{1}{\epsilon}
    \right\}.
  \end{equation}
\item  If $F$ is an admissible homogeneous p-merging function, 
then   there exist $(\lambda_1,\dots,\lambda_K)\in\Delta_K$ and admissible calibrators $f_1,\dots,f_K$
  such that for each $\epsilon \in (0,1)$,
  \begin{equation}\label{eq:calibrator-h}
    R_\epsilon(F)
    = 
    \left\{
      \mathbf p \in [0,\infty)^K:
      \sum_{k=1}^K \lambda_k  f_k\left(\frac{p_k}{   \epsilon}\right) \ge 1
    \right\}.
  \end{equation}
  \item  For any $(\lambda_1,\dots,\lambda_K)\in\Delta_K$ and calibrators $f_1,\dots,f_K$,
  \eqref{eq:calibrator-h} determines a homogeneous p-merging function.
\end{enumerate} 

\end{theorem}

 As a consequence of Theorem~\ref{th:e-merge-p} part (i), if the level $\alpha$ of type-I error control is determined before choosing the p-merging function, then it suffices to consider homogeneous ones, since their rejection sets are at least as larger as  those of other p-merging functions.  
Note that this does not imply that there exists a homogeneous p-merging function $G$ dominating $F$ in general, because the construction of $G$ depends on the given $\alpha.$

The calibrators in (ii) and (iii) of  Theorem~\ref{th:e-merge-p} are connected via 
$g_k(x) = f_k(x/\epsilon)/\epsilon$, $x\in[0,\infty)$, for each $k\in \cK$. This choice yields admissible calibrators by Proposition~\ref{prop:rescale-cal}. 
Therefore, when $F$ is an admissible homogeneous p-merging function, 
the calibrators for different $\epsilon$-levels are generated by one tuple of calibrators through  rescaling in the sense of Proposition~\ref{prop:rescale-cal}.

The p-merging function in \eqref{eq:calibrator-h} can be explicitly expressed via \eqref{eq:10-inf} as
\begin{equation}\label{eq:calibrator-h-exp}
   F(\mathbf p)
    = \inf
    \left\{\epsilon\in (0,1]:
      \sum_{k=1}^K \lambda_k  f_k\left(\frac{p_k}{   \epsilon}\right) \ge 1
    \right\} ~~
    \text{for  $\mathbf p \in [0,\infty)^K$}.
  \end{equation}

To prove Theorem~\ref{th:e-merge-p}, we need an additional lemma.   
We say that a set $R\subseteq [0,\infty)^K$ is a decreasing set if
$\mathbf x \in R $ implies $\mathbf y\in R $ for all
$ \mathbf y\in [0,\infty)^K$ with $\mathbf y
\le  \mathbf x$ (componentwise).
Let $\fU_0$ be the set of uniformly distributed random variables on $[0,1]$ under $\p$ (i.e., they are exact p-variables).
Clearly, $\fU_0\subseteq \fU$. 
For any decreasing set $L$,  we have  
\begin{equation}\label{eq:r1-larger}
\sup_{\mathbf P\in \fU_0^K} \p(\mathbf P\in   L) =\sup_{\mathbf P\in \fU^K} \p(\mathbf P\in   L), 
\end{equation}
because replacing a p-variable with a smaller exact p-variable does not reduce the above probability.
 The fact \eqref{eq:r1-larger} will be repeatedly used in the proof below.

\begin{lemma}
\label{lem:OTD}
Let $R\subseteq [0,\infty)^K$ be a decreasing Borel set. 
For any $\beta\in (0,1)$, we have 
$$
\sup_{\mathbf P\in \fU_0^K} \p(\mathbf P\in \beta R) \ge \beta \iff
\sup_{\mathbf P\in \fU_0^K} \p(\mathbf P\in   R)  = 1. 
$$
\end{lemma}
\begin{proof}[Proof.]

We first prove the $\Leftarrow$ direction  
  by   contraposition. 
Suppose $$\gamma: =\sup_{\mathbf P\in \fU_0^K}  \p(\mathbf P\in \beta R)< \beta.$$
Take an  event $A$ with probability $\beta$ and any   $\mathbf P \in \fU_0^K$  independent of $A$. Define $\mathbf P^* $ by $$\mathbf P^*= 
\beta  \mathbf P \times  \id_A+ \mathbf 1 \times \id_{A^c}.$$  
It is straightforward to check $\mathbf P^*\in \fU^K$. Hence,  by \eqref{eq:r1-larger},
$$\beta 
\p\left (    \mathbf P\in     R \right)
= \p(A) \p\left (    \mathbf P\in    R  \right)
\le  \p( \mathbf P^*    \in  \beta  R  )  \le \gamma ,$$
and  thus $\p(\mathbf P\in   R)  \le \gamma/\beta$. 
Since  $\gamma/\beta<1$, this then yields that $\sup_{\mathbf P\in \fU_0^K}  \p(\mathbf P\in   R) < 1$  and completes the $\Leftarrow$ direction. 

Next  we show the $\Rightarrow$ direction. 
Suppose   $\sup_{\mathbf P\in \fU_0^K} \p(\mathbf P\in \beta R) \ge \beta$.
For any $\epsilon\in (0,\beta)$, there exists $\mathbf P=(P_1,\dots,P_K)\in \fU_0^K$ such that $  \p(\mathbf P\in \beta R) > \beta-\epsilon.$
Let $A=\{\mathbf P\in \beta R\}$, $\gamma=\p(A)$,
and  $B$ be an event containing $A$ with $\p(B)=\beta\vee \gamma$. Let $\mathbf P^*=(P_1^*,\dots,P_K^*)$ follow the conditional distribution of $\mathbf P/\beta$ given $B$. We have
$$
\p( \mathbf P^* \in R) = \p(\mathbf P \in \beta R\mid B) =\p(A \mid B) = \frac{\gamma}{\beta\vee \gamma}.
$$ 
Note that  for $k\in \{1,\dots,K\}$, 
$$
\p(P_k^* \le  p ) = \p(P_k /\beta \le p \mid B) \le \frac{\p(P_k\le \beta p )}{\p(B)} 
=  \frac{\beta p }{\beta\vee \gamma}  \le p,
$$
and hence   $\mathbf P^*\in \fU^K$.
Since $\gamma>\beta-\epsilon $ and $\epsilon\in (0,\beta)$ is arbitrary, we can conclude      
$ 
\sup_{\mathbf P\in \fU^K} \p(\mathbf P\in   R) = 1$, yielding  $\sup_{\mathbf P\in \fU_0^K} \p(\mathbf P\in   R)  = 1
$   via \eqref{eq:r1-larger}.
\end{proof}

\begin{proof}[Proof of Theorem~\ref{th:e-merge-p}.]

For part (i), since any p-merging function has an admissible p-merging function dominating it (Proposition~\ref{prop:p-properties}), it suffices to prove (i) for admissible p-merging functions. 

Let $F$ be an admissible p-merging function and fix  $\epsilon \in (0,1)$. 
Note that the set $R_\epsilon(F)$ is a lower set, and it is closed due to Proposition~\ref{prop:p-properties}. Hence $\id_{R_\epsilon(F)}$ is upper semicontinuous. 
Using optimal transport duality,  Lemma~\ref{lem:duality},
\begin{align*} & \min_{(h_1,\dots,h_K)\in\mathcal B^K}
  \left\{
    \sum_{k=1}^K \int_0^1 h_k(x)\d x:  \bigoplus_{k=1}^K h_k \ge \id_{R_\epsilon(F)}
  \right\}
  \\ & =
  \max_{\mathbf P\in\fU_0^K}
  \p(\mathbf P\in R_\epsilon(F))
  =
  \epsilon,
\end{align*}
where the last equality holds because $F$ is precise (Proposition~\ref{prop:p-properties}).
Take $(h^\epsilon_1,\dots,h^\epsilon_K)\in \mathcal B^K$ such that
$$\bigoplus_{k=1}^K h^\epsilon_k \ge \id_{R_\epsilon(F)}
\mbox{~~and~~}
\sum_{k=1}^K \int_0^1 h^\epsilon_k(x)\d x = \epsilon.$$
 By the corresponding conditions on $\id_{R_\epsilon(F)}$ in Lemma~\ref{lem:duality}, we can choose each $h^\epsilon_k$ to be nonnegative, decreasing and left-continuous. Moreover, one can also require $h^\epsilon_k(0)=\infty$ for each $k$.

 Lemma~\ref{lem:OTD} gives 
\begin{equation}\label{eq:jm}
  \max_{\mathbf P\in \fU_0^K} \p( \mathbf P\in R_\epsilon (F)) = \epsilon
  ~~\Longrightarrow~~
  \max_{\mathbf P\in \fU_0^K} \p(\epsilon \mathbf P\in  R_\epsilon ( F)  ) = 1.
\end{equation}
Therefore, using duality in Lemma~\ref{lem:duality} again,
\[
  \min_{(h_1,\dots,h_K)\in\mathcal B^K}
  \left\{
    \sum_{k=1}^K \frac1\epsilon \int_0^\epsilon h_k(x)\d x:\bigoplus_{k=1}^K h_k \ge \id_{R_\epsilon(F)}
  \right\}
  =
  1,
\]
implying $\sum_{k=1}^K \int_0^\epsilon h^\epsilon_k(x) \d x\ge\epsilon$.
As  $\sum_{k=1}^K \int_0^1 h^\epsilon_k(x)\d x = \epsilon$ and each $h^\epsilon_k$ is nonnegative,
we know $h^\epsilon_k(x)=0$ for $x>\epsilon$.

Define the set $A_\epsilon := \{\mathbf p \in [0,\infty)^K:  \sum_{k=1}^K h^\epsilon_k(p_k) \ge  1\}$.
Since $\bigoplus_{k=1}^K h^\epsilon_k \ge \id_{R_\epsilon (F)}$,
we have $R_\epsilon ( F)  \subseteq  A_\epsilon$.
Note that $A_\epsilon $ is a closed lower set.
By Markov's inequality, 
\[
  \sup_{\mathbf P\in \fU_0^K}
  \p \left( \bigoplus_{k=1}^K h^\epsilon_k(\mathbf P) \ge 1 \right)
  \le
  \sup_{P\in \fU_0}
  \sum_{k=1}^K\E^\p[h^\epsilon_k(P)]
  =
  \epsilon. 
\]
Hence, we can define a function $G:[0,\infty)^K\to \R$ via  $R_\delta ( G) = \delta\epsilon^{-1} A_\epsilon  $ 
for all $\delta \in(0,1)$. 
By the above properties of $A_\epsilon$ and Lemma~\ref{lem:OTD}, we can check that $G$ is a valid homogeneous p-merging function, and $R_\epsilon(F)\subseteq R_\epsilon(G)$.
This proves  part (i).
To see the result in part (ii), it suffices to take $\lambda_k = \epsilon^{-1} \int_0^\epsilon g^\epsilon(x) \d x$
and
$g_k: [0,\infty)\to \R$, $x\mapsto  \epsilon^{-1} h^\epsilon_k( x)/\lambda_k $ for each $k=1,\dots,K$, with $f_k=1$ if $\lambda_k=0$, which gives that $g_1,\dots,g_K$ are admissible calibrators, and 
$$
A_\epsilon =   \left\{ 
  \mathbf p \in [0,\infty)^K:  \sum_{k=1}^K \lambda_k g _k( p_k) \ge  \frac{1}{\epsilon}  \right\}.
$$

 Before proving part (iii), we first prove part (iv).
 For a function $F$ defined by \eqref{eq:calibrator-h-exp},
 the condition $\p(F(\mathbf P)\le \epsilon)\le \epsilon$ for $\mathbf P\in \fU_0^K$ follows
  from Proposition~\ref{prop:rescale-cal} and Markov's inequality.
  Its increasing monotonicity and homogeneity are straightforward to check. Therefore, it is a homogeneous p-merging function. 

Now we can prove part (iii).  
Let $f_k(x) = \epsilon g_k(\epsilon x)$ for $x\in [0,\infty)$,
and define $H$ via the rejection regions in \eqref{eq:calibrator-h}, which is a homogeneous p-merging function by part (iv). 
Then we have $R_\epsilon(F) \subseteq A_\epsilon = R_\epsilon(H)$.
Since $F$ and $H$ are both homogeneous, this means 
$R_\delta(F) \subseteq R_\delta(H)   $ for all $\delta\in (0,1)$.
The admissibility of $F$ now gives  $F=H$, and this proves part (iii). 
\end{proof}

If a homogeneous p-merging function $F$ is symmetric,
then $f_1,\dots,f_K$ in Theorem~\ref{th:e-merge-p} can be chosen identical, and the same applies to $\lambda_1,\dots,\lambda_K$.

\begin{theorem}\label{th:e-merge-p2}
  For any $F$ that is admissible within the family of homogeneous symmetric p-merging functions,
  there exists an admissible calibrator $f$ such that for each $\epsilon \in (0,1)$,
  \begin{equation}\label{eq:calibrator2}
    R_\epsilon (F)
    =
    \left\{\mathbf p \in [0,\infty)^K: \frac 1 K \sum_{k=1}^K f\left(\frac{p_k}{\epsilon}\right) \ge 1 \right\}.
  \end{equation}  
  Conversely, for any  calibrator $f$, \eqref{eq:calibrator2} determines a homogeneous symmetric p-merging function.
\end{theorem}

The p-merging function in \eqref{eq:calibrator2} can be explicitly expressed as
\begin{equation}\label{eq:calibrator-2-exp}
   F(\mathbf p)
    = \inf
    \left\{\epsilon\in (0,1]:
      \frac 1K\sum_{k=1}^K  f  \left(\frac{p_k}{   \epsilon}\right) \ge 1
    \right\}~
    \text{for  $\mathbf p \in [0,\infty)^K$}.
  \end{equation}

\begin{proof}[Proof.]
  The proof is similar to that of Theorem~\ref{th:e-merge-p} and we only mention the differences.
  For the first statement, it suffices to notice two facts.
  First, if $R_{\epsilon}$ is symmetric, then $h^\epsilon_1,\dots,h^\epsilon_K$ in the proof of Theorem~\ref{th:e-merge-p} can be chosen as identical;
  for instance, one can choose the average of them.
  Second, the symmetry of $R_\epsilon (F)$ guarantees that $G$ in the proof of Theorem~\ref{th:e-merge-p} is symmetric,
  and hence it is sufficient to require the admissibility of $F$  within homogeneous symmetric p-merging functions in this proposition.
  The last statement in the proposition follows from Theorem~\ref{th:e-merge-p}
  by noting that \eqref{eq:calibrator2} defines a symmetric rejection region.
\end{proof}

The requirement $f(0)=\infty$ for an admissible calibrator $f$
implies that the combined test \eqref{eq:calibrator2} gives a rejection as soon as one of the input p-values is $0$,
which is obviously necessary for admissibility (Proposition~\ref{prop:p-properties}).
Although many examples in the M- and O-families, in particular $F_{r,K}$ for $r>0$ and $G_{k,K}$ for $k>1$, do not satisfy this,
we can make the zero-one adjustment 
\begin{equation}\label{eq:zero}
  \widetilde F(\mathbf{p})
  :=
  \begin{cases}
    F(\mathbf{p}\wedge\mathbf{1}) \wedge 1 & \text{if $\mathbf{p}\in(0,\infty)^K$}\\
    0 & \text{otherwise},
  \end{cases}
\end{equation}
which does not affect the validity of the p-merging function.

 For a decreasing function $f:[0,\infty)\to [0,\infty]$ and a p-merging function $F$ taking values in $[0,1]$,
we say that $f$ \emph{induces} $F$ if \eqref{eq:calibrator2} holds;
similarly, we may say that $\lambda_1,\dots,\lambda_K$ and $f_1,\dots,f_K$ \emph{induce} $F$ if \eqref{eq:calibrator} holds.
Theorems~\ref{th:e-merge-p} and~\ref{th:e-merge-p2} imply that any admissible p-merging function $F$ is induced by some  admissible calibrators (but they need not be uniquely determined by $F$).

The converse direction, constructing an admissible p-merging function, is more delicate. In general,
  a p-merging function induced by admissible calibrators is not necessarily admissible.
  Using \eqref{eq:jm} and a compactness argument,
  a necessary and sufficient condition for a calibrator $f$ to induce a precise p-merging function
  (a weaker requirement than admissibility)
  via \eqref{eq:calibrator2} is \begin{equation}\label{eq:existp}
    \p\left( \frac 1K  \sum_{k=1}^K f(P_k) \ge 1\right)=1
    \quad
    \text{for some $(P_1,\dots,P_K)\in \fU_0^K$}.
  \end{equation}
  Condition \eqref{eq:existp} may be difficult to check for a given $f$ in general.
 A special known case is that $f$ is convex.
 In this case,
  \eqref{eq:existp} holds if and only if $f \le K$ on $(0,1]$, which is a highly nontrivial mathematical result.  
  This condition turns out to  be also sufficient for admissibility of $F$ in Theorem~\ref{th:e-merge-p2}. The proof of Theorem~\ref{th:admissible} is advanced and omitted here. 
 
\begin{theorem}\label{th:admissible}
If $f$ is an admissible calibrator 
 that is strictly convex on $(0,1]$, $f\le K$ on $(0,1]$, and $f(1)=0$, then the p-merging function induced by $f$ via \eqref{eq:calibrator2} is admissible.
 
 More generally, the same holds true if $f$ is an admissible calibrator satisfying the following condition: for some $\eta\in[0,1/K)$ and $\tau=1-(K-1)\eta$, $f=K$ on $(0,\eta]$, $f(\eta+)\le K$ (where $f(\eta+)$ is the right limit of $f$ at $\eta$),  $f$ is strictly convex on $(\eta,\tau]$, and $f(\tau)=0$.
\end{theorem}
\begin{remark}
There is a corresponding  statement for the case of  concavity instead of convexity in Theorem \ref{th:admissible}. 
 The p-merging function induced by $f$ via \eqref{eq:calibrator2} is admissible if $f$ satisfies the following condition: for some $\eta\in[0,1/K)$ and $\tau=1-(K-1)\eta$, $f=K$ on $(0,\eta]$, $f(\eta+)\ge K/(K-1)$,  $f$ is  strictly concave on $(\eta,\tau]$, and $f(\tau)=0$.
\end{remark}

\subsection{Examples}
\label{sec:10-2-EX}

We present a few examples of p-merging functions, continuing those in  Section~\ref{sec:p-merging-examples}, through the constructions Theorem~\ref{th:e-merge-p}  using calibrators.

\begin{example}[Order-family]\label{ex:ex-th:e}
  The p-merging function $F:= G_{k,K}$, $k\in \cK$,
  is induced by  the calibrator $f(p)=(K/k) \id_{p\in [0,k/K]}$.
  Its zero-one adjusted version is admissible for $k\in [K-1]$.
\end{example}

\begin{example}[Twice the arithmetic mean]\label{ex:ex-th:e3}
  The  function
 $
    2\fM_{1,K},
 $ 
  which is 
  twice the arithmetic mean,
  is a precise p-merging function but not admissible. 
We have verified that it is a p-merging function in Proposition~\ref{prop:2avg} in Section~\ref{sec:c12-average-p}.
  
  This p-merging function can be written 
  in the form   \eqref{eq:calibrator-2-exp} with $f(p)=2-2p$ for $p\ge 0$, but 
  this $f$ is not a calibrator because $f(p)$ is not $0$ for $p> 1$. Using the calibrator $p\mapsto (2-2p)_+$ in \eqref{eq:calibrator-2-exp} gives a p-merging function that dominates $2\fM_{1,K}$, a slight improvement. 
\end{example}

\begin{example}[$\e$ times the geometric mean]\label{ex:ex-th:e4}
   The  function
   $
     \e \fM_{0,K},
  $
  that is, $\e$ times  the geometric mean,
  is a p-merging function that is not precise, but approximately precise for large $K$.
  
    This p-merging function can be written 
  in the form   \eqref{eq:calibrator-2-exp} with $f(p)=-\log(p)$ for $p\ge  0$. 
  Using the calibrator $p\mapsto (-\log(p))_+$ in \eqref{eq:calibrator-2-exp} gives a slight improvement.
\end{example}

\begin{example}[$\log K$ times the harmonic mean] \label{ex:ex-th:e5}
Let $T_K=\log K + \log \log K + 1$.
  The  function
 $
(T_K +1)\fM_{-1,K} ,
$  
 that is, $T_K+1$ times the harmonic mean,
  is a p-merging function that is not precise, but approximately precise for large $K$.
  
    This p-merging function can be written 
  in the form   \eqref{eq:calibrator-2-exp} with $$f(p)= \left(\frac{1}{T_K p}-\frac{1}{T_K}\right) \wedge K$$ for $p\ge  0$, and it is an elementary exercise in calculus to check that $p\mapsto (f(p))_+$ is a calibrator, and using it in \eqref{eq:calibrator-2-exp} gives a slight improvement.
\end{example}

\begin{example}[Hommel function]\label{ex:ex-th:e6}
The Hommel function $H_K$ is a precise p-merging function but not admissible. 
For $K\ge 4$, it is strictly dominated by 
the p-merging function $H_K^*$ induced by
the calibrator $f$ given by $$f(p)= \frac{K\id_{\{\ell_K p\le 1\}}}{\lceil K\ell_Kp \rceil},~~p\ge 0,$$
and $H_K^*$ is admissible when $K$ is not a prime number. The latter point is an intriguing fact; for example, if $K=2$ or $3$, $H_K^*$ is shown to be not admissible. 
\end{example}

\section{Combining exchangeable p-values}
\label{sec:c11-2}

In this section, we additionally assume that the input p-variables are exchangeable. 
  Exchangeability among p-values is  encountered, for example, in statistical testing via sample splitting. 
  Reasons for  sample splitting include 
  to relax the assumptions needed to obtain theoretical guarantees and to reduce computational costs. The drawback of methods based on sample-splitting is that the obtained p-values are affected by the randomness of the split, and thus repeated re-sampling can be performed, resulting in exchangeable p-values. 

One can of course combine the obtained p-values by using the  rules obtained in Section~\ref{sec:10-2} (such as twice the average in Example~\ref{ex:ex-th:e3}) for arbitrarily dependent p-values. But we can do better by exploiting   exchangeability.  
However, the constant $2$ in front of the arithmetic average cannot be directly improved. Indeed, we show that if we stick to symmetric merging functions, there is no hope to improve even assuming exchangeability.

\begin{definition}
An  {\emph{ex-p-merging function}} is an increasing Borel function $F:[0,\infty)^{K}\to [0,\infty)$
such that 
$\p(F(\mathbf P) \le \alpha)\le \alpha$ for all $\alpha \in (0,1)$ and $\mathbf P \in \fU^K$ that is exchangeable. 
It is \emph{homogeneous} if $F(\lambda \mathbf p)=\lambda  F(\mathbf p)$ for all $ \lambda  \in (0,1]$ and $\mathbf p $ with $F(\mathbf p)\le 1$.
\end{definition}

\begin{proposition}\label{prop:ex-p-m}
A symmetric ex-p-merging function is necessarily a p-merging function. Hence, for an ex-p-merging function to strictly dominate an admissible p-merging function, it cannot be symmetric. 
\end{proposition}

\begin{proof}[Proof.]
Let $\mathbf P\in \fU^K$,
and let $\pi$ be a random permutation of $\{1,\dots,K\}$, uniformly drawn from all permutations of $\{1,\dots,K\}$ and independent of $\mathbf P$.
Let $\mathbf P^\pi= (P_{\pi(1)}, \dots, P_{\pi (K)})$. 
Note that $\mathbf P^\pi$  is exchangeable by construction.
If $F$ is a symmetric {ex-p-merging function},
it must satisfy $F(\mathbf P^{\pi})=F(\mathbf P)$.
Because $F(\mathbf P^{\pi})$ is a p-variable, 
so is $ F(\mathbf P)$, showing that $F$ is a p-merging function. 
\end{proof}

In many practical settings, the exchangeable p-values can be generated one by one by repeating the same randomized procedure many times, generating a stream of p-values. We introduce combination rules that would simply process these p-values in the order that they are generated.

The next result follows by combining the exchangeable Markov inequality in Theorem~\ref{thm:emi}, the p-merging method from Theorem~\ref{th:e-merge-p2} and the formula \eqref{eq:calibrator-h-exp}.

\begin{theorem}[Combination of exchangeable p-values]\label{th:ex-comb-p}
Define the function $F:[0,\infty)^K\to [0,1]$ by
    $$F(\mathbf p)
    = \inf
    \left\{\epsilon\in (0,1]:
      \bigvee_{\ell=1}^K \left( \frac1{\ell }\sum_{k=1}^\ell  f  \left(\frac{p_k}{   \epsilon}\right)\right) \ge 1
    \right\},$$
where $f$ is a  calibrator. 
Then, for any vector $\mathbf P$ of  p-variables that is exchangeable,
$F(\mathbf P)$ is a p-variable. 
That is, $F$ is an ex-p-merging function. 
\end{theorem}

In fact, for the method in Theorem~\ref{th:ex-comb-p},  we do not need to fix the number of p-values ahead of time, they can just be processed online, yielding a p-value whenever this procedure is stopped. This makes our merging rules particularly simple and practical. However, note that the online procedure requires that the calibrator  $f$ does not depend on $K$.

\section{Randomized combinations of p-values}
\label{sec:c11-3}

In Section~\ref{sec:10-2}, we obtained methods to merge p-values 
using e-values  in Theorems~\ref{th:e-merge-p} and~\ref{th:e-merge-p2},
and they are admissible with the condition in Theorem~\ref{th:admissible}. When randomization is allowed, the Markov inequality can be enhanced to the randomized Markov inequality, presented in Theorem~\ref{thm:umi}.
This allows us to design more powerful way of merging p-values under randomization.
We continue to work under the assumption that the p-values are arbitrarily dependent, as in Section~\ref{sec:10-2}.
The following result describes a method that uses randomization to improve p-merging functions.  This result can be seen as a combination of Theorem ~\ref{thm:umi} and~\ref{th:e-merge-p} and the formula \eqref{eq:calibrator-h-exp}.

\begin{theorem}[Randomized combination of p-values]\label{th:rand-comb-p}
Define the function $F:[0,\infty)^K\times [0,\infty) \to [0,1]$ by
    \begin{equation}
  \label{eq:c11-random}  F(\mathbf p, u)
    = \inf
    \left\{\epsilon\in (0,1]:
     \sum_{k=1}^K  \lambda_kf_k  \left(\frac{p_k}{   \epsilon}\right) \ge u
    \right\},
    \end{equation}
where $(\lambda_1,\dots,\lambda_K)\in\Delta_K$ and $f_1,\dots,f_K$ are calibrators. 
Then, for any vector $\mathbf P$ of  p-variables and another p-variable $V$ independent of $\mathbf P$,
$F(\mathbf P, V)$ is a p-variable.
    
\end{theorem} 
\begin{proof}[Proof.]
Let $U\lawis \mathrm{U}[0,1]$. The 
condition
$\p(F(\mathbf P, U)\le \alpha)\le\alpha $ for all $\alpha \in (0,1)$ 
follows from the randomized Markov inequality in Theorem~\ref{thm:umi}. For a general p-variable $V$, it suffices to notice that $F(\mathbf p, u)$ is increasing in $u$, and hence $\p(F(\mathbf P, V)\le \alpha)\le\p(F(\mathbf P, U)\le \alpha) \le \alpha $. 
\end{proof}

A simple subclass of \eqref{eq:c11-random} is  $$F(\mathbf p, u)
    = \inf
    \left\{\epsilon\in (0,1]:
      \frac 1K\sum_{k=1}^K  f  \left(\frac{p_k}{   \epsilon}\right) \ge u
    \right\},$$
    for a calibrator $f$, similarly to \eqref{eq:calibrator-h-exp}.

For the deterministic choice $V=1$, Theorem~\ref{th:rand-comb-p} gives the validity statement in Theorem~\ref{th:e-merge-p} (iv). 
In practice, when randomization is allowed, $V$ should be chosen as uniformly distributed on $[0,1]$.
This merging method produces a smaller p-value than the deterministic p-merging function in \eqref{eq:calibrator-h-exp}.

Randomization is undesirable in many statistical applications. One deterministic remedy is that, when one of the p-variables $P_k$ in $\mathbf P$ is known to be independent of the others, one can use $V=P_k$ and apply the merging function to the rest of the p-variables with randomization through $V$. 

If randomization is permitted, one can also improve the existing rules for combining arbitrarily dependent p-values by using the exchangeable combination rule in Section~\ref{sec:c11-3} applied to a random permutation of the p-values.

Table~\ref{tab:c11-1} summarizes some examples of p-merging methods under arbitrary dependence, randomization, or exchangeability. 
They are derived from Theorems~\ref{th:e-merge-p2},~\ref{th:ex-comb-p}  and~\ref{th:rand-comb-p} with different calibrators $f$, and correspond to Examples~\ref{ex:ex-th:e}--\ref{ex:ex-th:e5}.

\renewcommand{\arraystretch}{2.5}
\begin{table}[t]
    \centering
    \resizebox{\textwidth}{!}{
    \begin{tabular}{l|c|c|c}
    Context
    & {\begin{minipage}{0.22\textwidth}\centering 
 Arbitrary dependence \\ (Section~\ref{sec:10-2})
    \end{minipage}}
    & {\begin{minipage}{0.22\textwidth}\centering 
 Exchangeability\\ (Section~\ref{sec:c11-2})
    \end{minipage}}
    & 
    {\begin{minipage}{0.24\textwidth}\centering 
 Arbitrary  dependence,\\  randomized (Section~\ref{sec:c11-3})
    \end{minipage}}
    \\  
Result & Theorem~\ref{th:e-merge-p2} & Theorem~\ref{th:ex-comb-p} & Theorem~\ref{th:rand-comb-p}
      \\ \hline
Order statistics & $\frac{K}{k} p_{(k)}$ & $\frac{K}{k} \bigwedge_{m=1}^K p^m_{(\lambda_m)}$ & $\frac{K}{k} p_{(\lceil Uk \rceil)}$ 
     \\ Arithmetic mean &  $2 \fM_1(\mathbf{p})$ & $2 \bigwedge_{m=1}^K \fM_1(\mathbf{p}_m)$ &   $\frac{2}{2-U} \fM_1(\mathbf{p})$
    \\ Geometric mean &  $\e \fM_0(\mathbf{p})$   &   $\e \bigwedge_{m=1}^K \fM_0(\mathbf{p}_m)$ & $\e^U \fM_0(\mathbf{p})$ 
    \\ Harmonic mean & $(T_K+1) \fM_{-1}(\mathbf{p})$   & $(T_K+1) \bigwedge_{m=1}^K \fM_{-1} (\mathbf{p}_m)$   &  $(T_KU + 1)\fM_{-1}(\mathbf{p})$ \\
    \end{tabular}} 
    \caption{Some combination rules corresponding to  Examples~\ref{ex:ex-th:e}--\ref{ex:ex-th:e5}.
    The order-statistics family is indexed by $k\in \{1,\dots,K\}$.
    Here, $\mathbf{p}=(p_1, \dots, p_K)$ denotes the vector of p-values, and $\mathbf{p}_m$ represents the vector containing the first $m$ values of $\mathbf{p}$. In the table, $p_{(k)}$ is the $k$-th smallest value of $\mathbf{p}$, while $p_{(\lambda_m)}^m$ is the $\lambda_m=\lceil m k / K  \rceil$ ordered value of $\mathbf{p}_m$. 
    The random variable $U$ is uniformly distributed on the interval $[0,1]$ independent of the p-values.
The functions $\fM_1$, $\fM_{0}$ and $\fM_{-1}$ respectively denote the arithmetic mean, the geometric mean, and the harmonic mean of the suitable dimension (we omit the dimension from the subscript). The value $T_K$ is given by $T_K=\log K + \log \log K +1$   for $K \geq 2$. 
Another randomized combination rule  $\frac{1}{2-2U} \fM_1(\mathbf{p})$ can be derived from \eqref{eq:randomized-avg1}.
}
   
    \label{tab:c11-1}
\end{table}

\section*{Bibliographical note}

 Merging p-values has a long history, and some early works include
 \cite{Tippett:1931}, \cite{Pearson:1933} and \cite{Fisher:1948-short}. 
 The Simes and Hommel functions in Example~\ref{ex:c10-Simes} were proposed by \cite{Simes:1986} and \cite{Hommel:1983}, respectively. 
 The order-family was proposed by 
 \cite{Ruger:1978}. 
 The harmonic mean was proposed by \cite{wilson2019harmonic}.
The mean-family was formally studied by 
 \cite{vovk2020combining}.
In some of the above works, often the p-values are assumed to be independent or follow a certain dependence structure; key exceptions are \cite{Ruger:1978}, \cite{ruschendorf1982random} and \cite{vovk2020combining}.

The content of Section~\ref{sec:c12-average-p} is  based on \cite{wang2024testing}, where average p-values are called p*-values. 
Proposition~\ref{prop:2avg} was first shown by \cite{ruschendorf1982random} based on  optimal transport duality and then by \cite{meng1994posterior} based on convex ordering inequalities; our proof is shorter and simpler. 
The full proof of Proposition~\ref{prop:c12-avp-eq}  can be found in \cite{wang2024testing}, where
  (iv)$\Leftrightarrow$(v) follows from
\citet[Theorem 4.A.2]{shaked2007stochastic},
(iv)$\Rightarrow$(i) is based on \citet[Theorem 5] {mao2019sums} and (iv)$\Rightarrow$(iii) is based on the refinement of Strassen's theorem in \citet[Theorem 3.1]{nutz2024martingale}. 
Mid p-values were studied by \cite{lancaster1952statistical}, and some methods of combining independent mid p-values were proposed by \cite{rubin2019meta}.  
The posterior predictive p-values were proposed by \cite{meng1994posterior}. 
The content of Section~\ref{sec:10-2} 
 is mainly based on   \cite{vovk2022admissible},
 and the content of Sections~\ref{sec:c11-2}--\ref{sec:c11-3} is based on \cite{gasparin2024combining}.

The necessary and sufficient condition for \eqref{eq:existp} to hold for convex $f$ is given in Theorem 3.2 of \cite{WW16}. This result is also needed to prove Theorem~\ref{th:admissible}.
The proof of Theorem~\ref{th:admissible} and justifications for statements in Section~\ref{sec:10-2-EX} are  found in Section 6 of \cite{vovk2022admissible}.  

\chapter{False coverage rate control using e-confidence intervals}
\label{chap:eci}

This chapter studies confidence intervals formulated by e-values, that we call e-confidence intervals. These will often not be intervals, but confidence sets (regions), but we still call them confidence intervals for simplicity (and to enable the common abbreviation CI).

We will see that these e-confidence intervals have additional properties compared to the usually constructed confidence intervals.
In this chapter, we let $\cP$ denote the set of all possible data distributions. The data will be drawn from some $\p^{\dagger} \in \cP$, and we will be interested in estimating $\vartheta(\p^{\dagger})$, where $\vartheta:\cP\to \Theta$ is a predefined functional of interest (such as the mean, median, and so on).
Since $\p^\dagger$ is unknown, we will consider generic distributions $\p\in \cP$, as in Chapter~\ref{chap:compound}.
Let $$\cP_{\theta} := \{\p \in \cP: \vartheta(\p) = \theta \}$$ denote the set of all distributions whose functional equals $\theta$.
We continue to write $\cK=[K]$ as in the previous chapters.

\section{Defining and constructing e-CIs}


We first define e-confidence intervals as a generalization of the usual confidence intervals. 

\begin{definition}[Families of e-variables and e-confidence intervals (e-CIs)]
We say that  $\{E(\theta)\}_{\theta \in \Theta}$ is a \emph{family of e-variables} for $\{\cP_{\theta}\}_{\theta \in \Theta}$ if for each $\theta \in \Theta$, $E(\theta)$ is an e-variable for $\cP_{\theta}$.

For a fixed $\alpha\in[0,1]$, a set $C(\alpha) \subseteq \Theta$ is called a $(1-\alpha)$-confidence interval (CI) for a functional $\vartheta$ if $\p(\vartheta(\p) \in C(\alpha)) \geq 1-\alpha$ for all $\p \in \cP$. We say that $\{C(\alpha)\}_{\alpha \in [0,1]}$ is a \emph{family of confidence intervals} for $\vartheta$ if for every $\alpha \in [0,1]$, $C(\alpha)$ is a $(1-\alpha)$-CI for $\vartheta$.

For a fixed $\alpha \in [0,1]$, we say that $C(\alpha)$ is an $(1-\alpha)$-e-CI if there exists a family of e-variables $\{E(\theta,\alpha)\}_{\theta \in \Theta}$ for $\{\cP_{\theta}\}_{\theta \in \Theta}$ such that
\begin{align}
    C(\alpha) = \left\{\theta \in \Theta: E(\theta,\alpha) < \frac{1}{\alpha}\right\}.
    \label{eqn:ECI}
\end{align}
We often omit the level $(1-\alpha)$ for CIs and e-CIs, which should be clear from the notation $C(\alpha)$ or from the context.

We say that $\{C(\alpha)\}_{\alpha \in [0,1]}$ is a family of \emph{e-confidence intervals (e-CIs)} if $C(\alpha)$ is an e-CI for every   $\alpha \in [0,1]$.
Finally, a family of e-CIs   called \emph{level-free} if the above e-variable $E(\theta,\alpha)$ does not depend on $\alpha$; in this case we simply write it $E(\theta)$.
\end{definition}

Without loss of generality, we can take $C(0)=\Theta$, $C(1)=\emptyset$.
Moreover, we usually would assume $C(\alpha) \supseteq C(\alpha')$ if $\alpha \leq \alpha'$; this is automatically true for a level-free family of e-CIs.

It is straightforward to observe that every e-CI is a CI and every family of e-CIs is a family of CIs: for the latter case, observe that for any $\p \in \cP_\theta$, we have $\p(\theta \not\in C(\alpha)) = \p(E(\theta,\alpha) \geq  1 / \alpha) \leq \alpha$ by Markov's inequality, since $\E^\p[E(\theta,\alpha)] \leq 1$ for any $\p \in \cP_\theta$. 
(This also explains why we formulated $C(\alpha)$ using $<1/\alpha$ instead of $\le 1/\alpha$; our choice leads to a smaller CI.)

Moreover, every CI is actually an e-CI by using
$$E(\theta,\alpha)=\id_{\{\theta \in C(\alpha)\}},$$
which is an  all-or-nothing e-variable in Section~\ref{sec:nontrivial-test}.
Therefore, the more interesting object to us is the level-free e-CI families.
 
\begin{example}[e-CI for Gaussian mean with unknown variance]
    We continue from Example~\ref{ex:t-test} to provide an intriguing e-CI for a Gaussian mean $\mu$, with an unknown variance $\sigma^2$. The following is an e-CI that can be constructed with a single observation, i.e. $n=1$:
    \[
    X_1 \pm \alpha^{-1} \exp((X_1^2 -1 )/2) .
    \]
    It is formed via the family of e-variables $\{|X_1-\mu|\exp((1-X_1^2)/2)\}_{\mu \in \R}$ for testing $\{\bigcup_{\sigma > 0} N(\mu,\sigma^2)\}_{\mu \in \R}$.
    In contrast, the classical t-CI for $\mu$ needs at least two observations to yield a nontrivial interval. 
\end{example}

We have already seen one general construction of e-CIs: those obtained via universal inference in Section~\ref{sec:ui-ci}. Here, we develop two more examples:  via stopped confidence sequences and calibrated CIs.

\subsection{E-CIs from stopped confidence sequences}

Let us first define the central concept of a \emph{confidence sequence}. The setup here is inherently sequential, as in Chapter~\ref{chap:eprocess}. One observes an increasing amount of data from an unknown distribution $\p$, and desires to estimate its functional $\vartheta$. Rather than a single $\sigma$-algebra $\cF$, one must now consider a filtration $\{\cF_t\}_{t \geq 0}$, which is a nested sequence of $\sigma$-algebras, indicating obtaining an increasing amount of information with the passage of time. 
A stopping time $\tau$ is a random variable such that $\{\tau \leq t\}$ is $\cF_t$-measurable for all $t\ge 0 $.

\begin{definition}
    Given  $\alpha \in [0, 1]$, a $(1 - \alpha)$-\emph{confidence sequence} for a functional $\vartheta$ is a sequence of sets $(C^t(\alpha))_{t \in \N}$ such that $C^\tau(\alpha) \subseteq \Theta$ is an $(1 - \alpha)$-CI for any stopping time $\tau$. It also has the following equivalent definition: for any $\p \in \cP$, we have
\begin{align}
    \p(\vartheta(\p) \in C^t(\alpha) \text{ for all }t \in \N) \geq 1 - \alpha .
\end{align}
\end{definition}

A universal way to construct such an object is using a family of  {e-processes}, which are sequential versions of e-values, as introduced in Chapter~\ref{chap:eprocess}.

For any fixed $\alpha$, let $E(\theta,\alpha) \equiv \{E_t(\theta,\alpha)\}_{t \geq 0}$ denote an e-process for $\cP_\theta$. Then, Ville's inequality implies that the set 
\[
 C^t(\alpha) = \left\{\theta \in \Theta: E_t(\theta,\alpha) < \frac{1}{\alpha}\right\}
\]
is a confidence sequence for $\vartheta$.

\begin{proposition}
\label{prop:c11-stop-eci}
    If $\{C^t(\alpha)\}_{t \in \N}$ is a confidence sequence as constructed above, then for any stopping $\tau$, $C^\tau(\alpha)$ is an e-CI. Further, if $E(\theta,\alpha)$ does not depend on $\alpha$, then $\{C^\tau(\alpha)\}_{\alpha \in [0,1]}$ is a level-free family of e-CIs.
\end{proposition}

Proposition~\ref{prop:c11-stop-eci} directly follows from Ville's inequality stated in Fact~\ref{Ville}.

\subsection{E-CIs by calibrating CIs}

Recall the notion of a calibrator from Definition~\ref{def:calib}. Define $f^{-1}$ to be the right inverse of the calibrator $f$, i.e.,\ $f^{-1}(x) := \sup \{p: f(p) \geq x\}$.  When $f$ is invertible, $f^{-1}$ is the inverse of $f$. Using $f^{-1}$, we can convert any CI to an e-CI. Before we do so, define two properties about any set-valued CI function $C: [0, 1] \mapsto  2^{\Theta}$. Let $C$ be \textit{decreasing} if, for any $\alpha, \beta \in [0, 1]$, $\alpha \leq \beta$ implies that $C(\alpha)\supseteq C(\beta)$. Further, define \textit{continuous from below} to be the property that for any $\alpha \in [0, 1]$, $C(\alpha) = \bigcup_{\beta > \alpha} C(\beta)$.

\begin{theorem}
\label{th:calibrate-eci}
Let $C: [0, 1] \mapsto  2^{\Theta}$ be a decreasing function that is continuous from below such that $C(\alpha)$ produces a $(1 - \alpha)$-CI, and $f$ be a calibrator with right inverse $f^{-1}$. Then, the following set $ C^{\mathrm{cal}}(\alpha)$ is a $(1 - \alpha)$-e-CI:
\begin{align}
    C^{\mathrm{cal}}(\alpha) &= C\left(f^{-1}\left(\frac{1}{\alpha}\right)\right)\supseteq C(\alpha).
\end{align}\label{thm:Calibration}
Moreover,  $\{C^{\mathrm{cal}}(\alpha)\}_{\alpha \in [0,1]}$ is a level-free family of e-CIs.
\end{theorem} 

Theorem~\ref{th:calibrate-eci} exploits the duality between CIs and families of tests (or families of p-values). For each $\theta \in \Theta$, we can test the null hypothesis $\cP_\theta$. The CI then corresponds to the values of $\theta$ for which the aforementioned test does not reject at level $\alpha$. We can invert this process to use the CI to define a p-value for every null hypothesis $\cP_\theta$. We can then calibrate the implicit p-value to yield an e-value, from which we can produce our e-CI. This is formalized in the proof below.


\begin{proof}[Proof of Theorem~\ref{th:calibrate-eci}.]

The inclusion
$
 C^{\mathrm{cal}}(\alpha)  \supseteq C(\alpha)
 $ follows directly from the fact that $f(p)\le 1/p$ for all $p\in [0,1]$ (explained in Section~\ref{sec:2-cali}), and hence $f^{-1}(1/\alpha) \ge \alpha$, and thus $C( f^{-1}(1/\alpha)  ) \supseteq C(\alpha)$. We next prove that  $\{C^{\mathrm{cal}}(\alpha)\}_{\alpha \in [0,1]}$ is a level-free family of e-CIs.

For any $\theta \in \Theta$, note that the following is a p-variable for $\cP_\theta$:
\begin{align}
    P^{\mathrm{dual}}(\theta) \coloneqq \inf \left\{\alpha \in [0, 1]: \theta \not \in C(\alpha)\right\}.
\end{align}
Consequently, $E^{\mathrm{cal}}(\theta) \coloneqq f(P^{\mathrm{dual}}(\theta))$ is an e-variable. Hence, \[
\left\{\theta \in \Theta: E^{\mathrm{cal}}(\theta) < \frac{1}{\alpha}\right\}\] yields a $(1 - \alpha)$-e-CI. The theorem now follows because
\begin{align}
    \left\{\theta \in \Theta: E^{\mathrm{cal}}(\theta) < \frac{1}{\alpha}\right\} &= \left\{\theta \in \Theta: f(P^{\mathrm{dual}}(\theta)) < \frac{1}{\alpha}\right\} \\
    &= \left\{ \theta \in \Theta : P^{\mathrm{dual}}(\theta) >  \max \left\{p: f(p) \geq \frac{1}{\alpha}\right\} \right\}\\
    &= \left\{ \theta \in \Theta : P^{\mathrm{dual}}(\theta) >  f^{-1}\left(\frac{1}{\alpha}\right) \right\} \\
    &=  C\left(f^{-1}\left(\frac{1}{\alpha}\right)\right).
\end{align} 
Above, the third equality is a result of $f$ being decreasing and upper semicontinuous at $1 / \alpha$; hence the supremum is achieved and the equality holds. The final equality is because $C$ is decreasing and continuous from below at $f^{-1}\left(1 / \alpha\right)$; if $P^{\mathrm{dual}}(\theta) = f^{-1}\left(1 / \alpha\right)$, then $\theta \not\in C_i(P^{\mathrm{dual}}(\theta))$ by $C$ being continuous from below.
\end{proof}

\section{The e-BY procedure for false coverage rate}\label{sec:eby}

\subsection{The goal of FCR control}


Consider a scientist who observes some data, $\mathbf{X} = (X_1, \dots, X_K)$ drawn from some   distribution $\p\in \cP$. 
The scientist is potentially interested in the values of $K$ of its functionals $\vartheta_1,\dots,\vartheta_K$ through the   parameter values $\boldsymbol\theta^*  := (\theta_1^*, \dots, \theta_K^*)$, where $\theta^*_i := \vartheta_i(\p)$. 
The  functionals  $\vartheta_1,\dots,\vartheta_K$ could each take values in different spaces, but this complicates notation, and we will not explicitly need these sets later on. Each $\vartheta_i$ need not be bijective;
for example, it could capture the median of a distribution. 

The scientist may be interested in identifying the $L \geq 1$ indices corresponding to the largest $L$ values of $\boldsymbol\theta^*$, assuming they are real-valued. Or they may be interested in any index $k$ such that $\theta^*_k$ is larger than some prespecified (or data dependent) threshold. In short, the scientist may not know which indices are of interest to them before observing the data, and it is quite likely that after observing the data, only a small fraction of indices are actually of interest.  

For each $i \in \mathcal K$, we assume that from $X_i$, the scientist can construct a ($1 - \alpha$)-confidence interval for $\theta^*_i$.

    The scientist uses the data $\mathbf{X}$ to select a subset of ``interesting'' parameters, $S \subseteq \mathcal K$, using some potentially complex data-dependent selection rule $\mathcal S: \mathbf{X} \mapsto S $. The scientist must then devise confidence levels for the CI of each selected parameter, $\{\alpha_i\}_{i \in S}$, that \emph{can depend on the data }$\mathbf{X}$. The  \emph{false coverage proportion} (FCP) and \emph{false coverage rate} (FCR) of such a procedure are defined by
\begin{align}
    \FCP = \FCP(S,\{\alpha_i\}_{i \in S}) := \frac{\sum_{i \in S}\id_{\{\theta_i^* \not\in C_i(\alpha_i)\}}}{|S| \vee 1},~~~  \FCR = \E^{\p}\left[\FCP\right].
\end{align} 
Note the obvious similarity between the concepts of FCR and FCP
and between the concepts of FDR and FDP in Chapter~\ref{chap:compound}; all of these terms depend on $\p\in \cP$ through $\theta^*_i=\vartheta_i(\p)$ and $\E^\p$. 
 If every $\alpha_i$ equals the same constant (say $\gamma$), we use the more succinct notation $\FCP(S,\gamma)$ to abbreviate $\FCP(S,\{\alpha_i\}_{i \in S})$.

Our goal is to design a method for choosing $\{\alpha_i\}_{i \in S}$ that guarantees $\FCR \leq \delta$ (for all $\p\in \cP$) at a predefined level $\delta \in [0, 1]$ provided by the scientists in advance,
\emph{regardless of what the selection rule $\mathcal S$ is, and in particular even if the rule is unknown and we only observe the selected set $S = \mathcal S (\mathbf X)$.}

\subsection{The e-BY procedure}

Our primary point of comparison is the so-called Benjamini--Yekutieli (BY) procedure. The BY procedure's choice of $\{\alpha_i\}_{i \in S}$ and resulting guarantees depend upon  assumptions (or knowledge) of the dependence structure in $\mathbf{X}$ and the selection algorithm $\mathcal S$.
Under certain restrictions (omitted here for brevity) on $\mathcal S$, the BY procedure sets $\alpha_i  = \delta |S| / K$ to ensure that the FCR is controlled at level $\delta$ for mutually independent $X_1, \dots X_K$.
However, when no such assumptions can be made (i.e., under arbitrary dependence and an unknown selection rule) the BY procedure sets $\alpha_i = \delta |S| / (K\ell_K)$, where $\ell_K := \sum_{i = 1}^K i^{-1} \approx \log K$ is the $K$-th harmonic number. Clearly, the BY procedure produces much more conservative CIs when no assumptions can be made about dependence or selection.

The above facts are reminiscent of the BH procedure discussed in Section~\ref{sec:c8-pBH}, and one may, analogously to the situation of e-BH versus BH procedures in Chapter~\ref{chap:compound}, hope that an e-CI based procedure has FCR control without the additional assumptions and the multiplicative penalty $\ell_K$. This will be addressed next. 

Now, we formally define the e-BY procedure as follows.

\begin{definition}
    The \emph{e-BY procedure} at level $\delta \in [0, 1]$ sets $\alpha_i = \delta |S| / K$ for each $i \in S$.
    \label{def:eBY}
\end{definition}
We show that a FCR bound can be proven quite simply given the fact that e-CIs are constructed for each selected parameter.

\begin{theorem}
    Let $\{C_i(\alpha)\}_{\alpha \in [0,1]}$ be a level-free family of e-CIs  for each $i \in \mathcal K$. Then, the e-BY procedure ensures $\FCR \leq \delta$ for any $\delta \in (0, 1)$ under any dependence structure between $X_1, \dots, X_K$, and for any selection rule $\mathcal S$. In fact, the e-BY procedure satisfies the stronger guarantee:
    \[
    \E^{\p}\left[ \sup_{S \in 2^{\mathcal K}} \sup_{\delta \in (0,1)} \frac{\FCP(S, \delta |S|/K )}{\delta } \right] \leq 1~~~~~\mbox{for all $\p\in \cP.$}
    \]
\label{thm:eBY}
\end{theorem}
\begin{proof}[Proof of Theorem~\ref{thm:eBY}]
We directly show an upper bound for the FCR as follows:
\begin{align}
    \FCR &= \mathbb{E}^\p\left[ \frac{\sum_{i \in S} \id_{\{\theta_i^* \notin C_i(\delta |S|/K)\}} }{|S| \vee 1} \right]\\
               &= \mathbb{E}^\p\left[ \frac{\sum_{i \in \mathcal K} \id_{\{E_i(\theta_i^*)|S|\delta / K > 1\}} \cdot \id_{\{i \in S\}} }{|S| \vee 1} \right]\\
    &\leq \sum_{i \in \mathcal K} \mathbb{E}^\p\left[ \frac{E_i(\theta_i^*) |S| \delta} {K(|S| \vee 1)} \right] =  \sum_{i \in \mathcal K} \frac{\delta}{K} \mathbb{E}^\p\bigg[E_i(\theta_i^*)  \cdot \frac{ |S|}{|S| \vee 1} \bigg] \leq \delta,
\end{align} 
where the first inequality is because $\id_{\{x>1\}}\le x$ for all $x\ge 0$. The second inequality is a result of the definition of the e-value for $\theta^*_i$ having its expectation under $\p$ be upper bounded by 1. This achieves our desired bound.

The proof of the more general claim follows by an easy amendment of the above proof.
\end{proof}

    We note in passing that FCR control of the e-BY procedure implies FDR control of the e-BH procedure, while the converse is not true. 

\section{Combining CIs via majority vote}

Given $K$ uncertainty sets that are arbitrarily dependent --- for example, confidence intervals for an unknown parameter obtained with $K$ different estimators, or prediction sets obtained via conformal prediction based on $K$ different algorithms on shared data --- we now address the question of how to efficiently combine them in a black-box manner to produce a single uncertainty set. We present a simple and broadly applicable majority vote procedure that produces a merged set with nearly the same error guarantee as the input sets. We then extend this core idea in a few ways: we show that weighted averaging can be a powerful way to incorporate prior information, and a simple randomization trick produces strictly smaller merged sets without altering the coverage guarantee. Further improvements can be obtained if the sets are exchangeable. Underlying all of these methods are e-values and various versions of Markov's inequality.

Formally, we start with a collection of $K$ different sets $C_k$ (one from each ``agent'') for the same target quantity $c$, each having a confidence level $1-\alpha$ for some $\alpha \in (0,1)$:
\begin{equation}
\label{eq:coverage}
\p\big(c \in C_k\big) \geq 1-\alpha, \quad k\in \mathcal K.
\end{equation}
The probability measure $\p$ in this section is the one that generates the data. 
We say that $C_k$ has \emph{exact} coverage if $\p(c \in C_k) = 1-\alpha$.
Since the sets $C_k$ are based on data, they are random quantities by definition, but $c$ can be either fixed or random; for example, in the case of confidence sets for a target functional of a distribution it is fixed, but it is random in the case of prediction sets for an outcome (e.g., conformal prediction). Our method will be agnostic to such details.

Our objective as the ``aggregator'' of uncertainty is to combine the sets in a black-box manner in order to create a new set that exhibits favorable properties in both coverage and size. A first (trivial) solution is to define the set $C^J$ as the union of the others:
\begin{equation*}
    C^J = \bigcup_{k=1}^K C_k.
\end{equation*}
Clearly, $C^J$ respects the property defined in \eqref{eq:coverage}, but the resulting set is typically too large and has significantly inflated coverage. On the other hand, the set resulting from the intersection $C^I = \bigcap_{k=1}^K C_k$ is narrower, but typically has inadequate coverage --- it guarantees at least $1-K\alpha$ coverage by the Bonferroni inequality, but this is uninformative when $K$ is large.

    If the aggregator knows the $(1-\alpha)$-confidence intervals not just for a single $\alpha$ but for every $\alpha \in (0,1)$, using the duality between CIs and tests, one can calculate a p-value for each $\theta \in \Theta$, each agent $k \in \mathcal K$. We have already seen many ways to combine dependent p-values, for example, by averaging them and multiplying by two, and these can be used to combine these p-values across agents into a single one and then obtain a $(1-\alpha)$-confidence interval for any $\alpha$ of the aggregator's choice.
    The current section addresses the setting where only a single interval is known from each agent, ruling out the above averaging schemes. 

\subsection{The majority vote procedure}
\label{sec:gen_maj_vote}

Let the observed data  $Z$ lie in a sample space $\mathcal{Z}$, while our target $c$ is a point in a measure space $(\cS, \mathcal{A}, \nu)$, where $\mathcal{A}$ is a $\sigma$-algebra on $\cS$ and $\nu$ is a measure on $\cS$. As mentioned earlier, it is important to note that $c$ can itself be a random variable. The sets $C_k = C_k(z) \subseteq \cS$, $k\in\mathcal K$, based on the observed data, follow the property \eqref{eq:coverage}, where the probability refers to the joint distribution $(Z,c)$ (or only $Z$ if $c$ is fixed). Let us define a new set $C^M$ that includes all points \emph{voted} by at least half of the sets:
\begin{equation}
\label{eq:cm}
    C^M:=\bigg\{s \in \cS: \frac{1}{K} \sum_{k=1}^K \id_{\{s \in C_k\}} > \frac{1}{2}  \bigg\}.
\end{equation}

\begin{theorem}
\label{th:CM}
    Let $C_1, \dots, C_K$ be $K \geq 2$ different sets based on the observed data $z$, satisfying property \eqref{eq:coverage}. Then, the set $C^M$  defined in \eqref{eq:cm} has coverage of at least $1-2\alpha$:
    \begin{equation}
    \label{eq: coverage_cm}
    \p\big(c \in C^M\big) \geq 1-2\alpha.
    \end{equation}
\end{theorem}
\begin{proof}[Proof]
    Let $\phi_k=\phi_k(Z,c)=\id_{\{c \notin C_k\}}$ be a Bernoulli random variable such that $\E^\p[\phi_k] \leq \alpha$, $k\in\mathcal K$.    
    Thus $E_k = \phi_k/\alpha$ is an e-value, implying that $\bar E = (E_1+\dots+E_K)/K$ is an e-value.
    Thus,
    \[
    \begin{split}
    \p(c \notin C^M) &= \p\left(\frac{1}{K} \sum_{k=1}^K \phi_k \geq \frac{1}{2} \right) = \p\left(\bar E \geq \frac1{2\alpha} \right) \leq 2\alpha,
    \end{split}
    \]
    by Markov's inequality, concluding the proof.
\end{proof}

Actually, a slightly tighter bound can be obtained if $K$ is odd.
    In this case, for a point to be contained in $C^M$, it must be voted for by at least $\lceil K/2 \rceil$ of the other sets. This implies that, with the same arguments as used in Theorem~\ref{th:CM}, the probability of miscoverage is equal to $$\alpha K/\lceil K/2 \rceil = 2\alpha K/(K+1),$$ approaching \eqref{eq: coverage_cm} for large $K$.

The conservative threshold of $1/2$ is needed to handle arbitrary dependence between the input sets. If they were independent, the right-hand side of~\eqref{eq:cm} can be replaced by the $\alpha$-quantile of a Binomial distribution with $K$ trials of probability $1-\alpha$.

        Theorem~\ref{th:CM} is known to be tight in a worst-case sense; a simple example shows that if $K$ is odd and if the sets have a particular joint distribution, then the error will equal $(\alpha K)/\lceil K/2 \rceil$. This worst-case distribution allows for only two types of cases: either all agents provide the same set that contains $c$ (so majority vote is correct), or $\lfloor K/2 \rfloor$ sets contain $c$ but the others do not (so majority vote is incorrect). Each of the latter cases happens with some probability $p$, so the probability that majority vote makes an error is $\binom{K}{\lfloor K/2 \rfloor +1} p$. The probability that any particular agent makes an error is $\binom{K-1}{\lfloor K/2 \rfloor}p$, which we set as our choice of $\alpha$, and then we see that the probability of error for majority vote simplifies to $\alpha K/\lceil K/2 \rceil$.  Despite the apparent tightness of majority vote in the worst-case, we will develop several ways to improve this procedure in non-worst-case instances, while retaining the same worst-case performance.

        \begin{remark}[When does majority vote overcover and when does it undercover?]
    While the worst-case theoretical guarantee for majority vote is a coverage level of $1-2\alpha$, sometimes it will get close to the desired $1-\alpha$ coverage, and sometimes it may even overcover, achieving coverage closer to one. Here, we provide some intuition for when to expect each type of behavior in practice assuming $\alpha < 1/2$, foreshadowing many results to come. If the sets are actually independent (or nearly so), we should expect the method to have coverage more than $1-\alpha$. This can be seen via an application of Hoeffding's inequality in place of Markov's inequality in the proof of Theorem~\ref{th:CM}: since each $\phi_k$ has expectation (at most) $\alpha$, we should expect $\frac{1}{K} \sum_{k=1}^K \phi_k$ to concentrate around $\alpha$, and the probability that this average exceeds $1/2$ is exponentially small (as opposed to $2\alpha$), being at most $\exp(-2K(1/2-\alpha)^2)$ by Hoeffding's inequality. In contrast, if the sets are identical (the opposite extreme of independence), clearly the method has coverage $1-\alpha$. As argued in the previous paragraph, there is a worst-case dependence structure that forces majority to vote to have an error of (essentially) $2\alpha$. Finally, if the sets are exchangeable, it appears more likely that the coverage will be closer to $1-\alpha$ than $1-2\alpha$. While one informal reason may be that exchangeability connects the two extremes of independence and being identical (with coverages close to $1$ and $1-\alpha$), a slightly more formal  reason is that under exchangeability, we will later in this section actually devise a strictly tighter set $C^E$ than $C^M$ which also achieves the same coverage guarantee of $1-2\alpha$, thus making $C^M$ itself likely to have a substantially higher coverage.
\end{remark}


The above method and result can be easily generalized beyond the threshold value of 1/2. We record it as a result for easier reference. For any $\tau \in [0,1)$, let
\begin{equation}\label{eq:c^tau}
    C^\tau:=\bigg\{s \in \cS: \frac{1}{K} \sum_{k=1}^K \id_{\{s \in C_k\}} > \tau \bigg\}.
\end{equation}

\begin{theorem}\label{th:ctau}
     Let $C_1, \dots, C_K$ be $K \geq 2$ different sets satisfying property
     \eqref{eq:coverage}. Then, 
    \begin{equation*}
    \p\big(c \in C^\tau \big) \geq 1-\alpha/(1-\tau).
    \end{equation*}
\end{theorem}
The proof follows the same lines as the original result outlined in Theorem~\ref{th:CM} and is thus omitted. 
As expected, it can be noted that the obtained bounds decrease as $\tau$ increases. In fact, for larger values of $\tau$, smaller sets will be obtained. One can check that this result also yields the right bound for the intersection ($\tau=1 - 1/K$) and the union ($\tau=0$). In certain situations, it is possible to identify an upper bound to the coverage of the set resulting from the majority vote.

Even if the input sets are intervals, the majority vote set may be a union of intervals. One can easily construct a simple aggregation algorithm to find this set quickly by sorting the endpoints of the input intervals and checking some simple conditions. 

\begin{figure}
    \centering
    \includegraphics[width=0.95\textwidth]{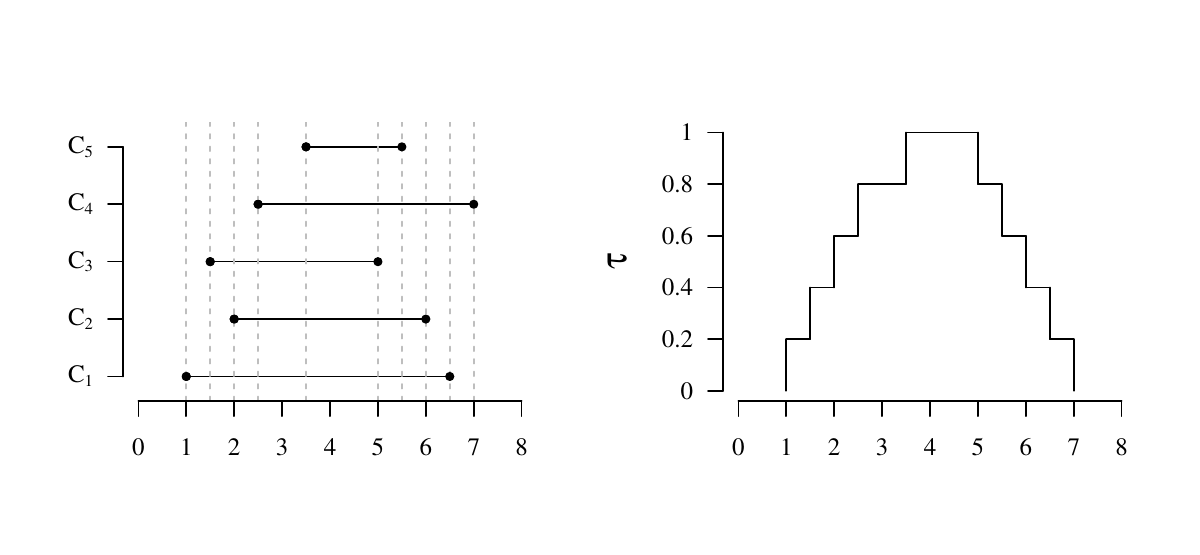}
    \caption{Visualization of the majority vote procedure when $\bigcap_{k=1}^K C_k \neq \emptyset$.}
    \label{fig:hist}
\end{figure}

Frequently, the majority vote set is indeed an interval. In fact, it is easy to check that if $C_1,\dots,C_K$ are one-dimensional intervals and $\bigcap_{k=1}^K C_k \neq \varnothing$, then $C^\tau$ is an interval for any $\tau$. See Figure~\ref{fig:hist}.

How large can the majority vote set be?
One naive way to combine the $K$ sets is to randomly select one of them as the final set; this method clearly has coverage $1-\alpha$, and its length is in between their union and intersection, so it seems reasonable to ask how it compares to majority vote. Surprisingly, majority vote is not always strictly better than this approach in terms of the expected length of the set: consider, for example, three nested intervals $C_1, C_2, C_3$ of width $10, 8$ and $3$, respectively. The majority vote set is $C_2$, with a length of $8$, but randomly selecting an interval results in an average length of $7$. However, we show next that the majority vote set cannot be more than twice as large. In addition, when the input sets are intervals, it is never wider than the largest interval. 

\begin{theorem}
    \label{thm:length_mv}
    Let $\nu(C^\tau)$ be the measure (size) of $C^\tau$ in \eqref{eq:c^tau}. Then, for all $\tau \in [0,1)$,
    \begin{equation}\label{eq:m^tau}
        \nu(C^\tau) \leq \frac{1}{K\tau}\sum_{k=1}^K \nu(C_k).
    \end{equation}
    If the input sets are $K$ one-dimensional intervals, then for all $\tau \in [\frac{1}{2}, 1)$, 
    \begin{equation}\label{eq:m^tau_max}
        \nu(C^\tau) \leq \max_k \nu(C_k).
    \end{equation}
\end{theorem}
Above, $\nu$ could be the Lebesgue measure (for intervals), or the counting measure (for discrete, categorical sets), for example. The proof uses the fact that the majority vote set is elementwise monotonic in its input sets, meaning that if any of the input sets gets larger, the majority vote set can never get smaller. From \eqref{eq:m^tau} we have that if $\tau=1/2$, then the  measure of the majority vote set is never larger than 2 times the average of the measure of the initial sets. This result is essentially tight as can be seen in the following example involving one-dimensional intervals. For odd $K$, let $(K+1)/2$ intervals have a length $L$, while the rest have a length of nearly 0. The average length is then $(K+1)L/(2K)$, and the majority vote has length $L$, whose ratio approaches $1/2$ for large $K$. In addition, \eqref{eq:m^tau} gives the right bound for the intersection ($\tau \uparrow 1$) and for the union ($\tau \uparrow 1/K$).


Consider a scenario similar to that of the last example. Suppose we have $K > 2$ confidence intervals, with odd $K$, such that $C_1= \dots = C_{\lceil K/2 \rceil} = (-\epsilon,1)$ and $C_{\lceil K/2 \rceil + 1} = \dots = C_K = (0,1+\epsilon)$, for some $\epsilon>0$. By definition, if $\tau=(K-2)/(2K)<1/2$, then the set $C^\tau$ coincides with the union of the initial sets, which is larger than the input intervals and has width $1+2\epsilon$. On the other hand, choosing $\tau=1/2$ we obtain a set whose size is $1+\epsilon$, matching the width of the input sets.

\subsection{Exchangeable, randomized and weighted variants}

Surprisingly, when $C_1, \dots, C_K$ are not independent, but are exchangeable, something better than a naive majority vote can be accomplished. To describe the method, let $C^M(1:k)$ denote the majority vote of sets $C_1, \dots, C_k$. Now define
\[
C^E := \bigcap_{k=1}^K C^M(1:k),
\]
which can be equivalently represented as 
\begin{equation}\label{eq:setce}
   C^E = \left\{s \in \cS: \frac{1}{k} \sum_{j=1}^k \id_{\{s \in C_j\}}>\frac{1}{2} \text{ holds for all }k\in \mathcal K \right\}. 
\end{equation}
Essentially, $C^E$ is formed by the intersection of sets obtained through sequential processing of the sets derived from the majority vote.

\begin{theorem}\label{thm:exch}
    If $C_1, \dots, C_K$ are $K \geq 2$ exchangeable sets having coverage $1-\alpha$, then $C^E$ is a $1-2\alpha$ uncertainty set, and it is never worse than majority vote ($C^E \subseteq C^M$).
\end{theorem}

The proof mimics that of Theorem~\ref{th:CM}, but uses the exchangeable Markov inequality from Chapter~\ref{chap:ui} in place of Markov's inequality. This method shares a same idea as in designing the  ex-p-merging functions in Section~\ref{sec:c11-2}.

This result points at a simple way of improving majority vote for arbitrarily dependent sets: process them in a random order. To elaborate, let $\pi$ be a uniformly random permutation of $\{1, \dots,K\}$ that is independent of the $K$ sets, and define
\begin{equation}\label{eq:c^pi}
    C^\pi := \bigcap_{k=1}^K C^M(\pi(1):\pi(k)).
\end{equation}

Since $C^M(\pi(1):\pi(K)) = C^M(1:K)$, $C^\pi$ is also never worse than majority vote despite satisfying the same coverage guarantee:

\begin{corollary}\label{cor:arbit-exch}
    If $C_1, \dots, C_K$ are $K \geq 2$ arbitrarily dependent uncertainty sets having coverage $1-\alpha$, and $\pi$ is a uniformly random permutation independent of them, then $C^\pi$ is a $1-2\alpha$ uncertainty set, and it is never worse than majority vote ($C^\pi \subseteq C^M$).
\end{corollary}

The proof follows as a direct corollary of Theorem~\ref{thm:exch} by noting that the random permutation $\pi$ induces exchangeability of the sets (the joint distribution of every permutation of sets is the same, due to the random permutation). Of course, if the sets were already ``randomly labeled'' 1 to $K$ (for example, to make sure there was no special significance to the labels), then the aggregator does not need to perform an extra random permutation.

Moving in a different direction below, we demonstrate that the majority vote can be improved with the aim of achieving a tighter set through the use of independent randomization, while maintaining the same coverage level.

Let $U$ be an independent random variable that is distributed uniformly on $[0,1]$, and let $u$ be a realization.
We then define a new set $C^R$ as
\begin{equation}
\label{eq:cr}
    C^R := \left\{s \in \cS: \frac{1}{K}\sum_{k=1}^K \id_{\{s \in C_k\}} > \frac{1}{2} + \frac u 2\right\}.
\end{equation}
As a small variant, define
\begin{equation}
\label{eq:cu}
    C^U := \left\{s \in \cS: \frac{1}{K}\sum_{k=1}^K \id_{\{s \in C_k\}} > u \right\}.
\end{equation}

\begin{theorem}
Let $C_1, \dots, C_K$ be $K \geq 2$ different sets satisfying the property \eqref{eq:coverage}. Then, the set $C^R$ has coverage at least $1-2\alpha$ and is never larger than $C^M$, while the set $C^U$ has coverage at least $1-\alpha$ and is never smaller than $C^R$.
\end{theorem}

The proof follows that of Theorem~\ref{th:CM}, but uses the randomized Markov inequality in Theorem~\ref{thm:umi} from Chapter~\ref{chap:markov} in place of Markov's inequality. This is related to randomized p-merging functions in Section~\ref{sec:c11-3}, in particular to the randomized order-family combination.
Even though $C^U$ does not improve on $C^M$, we include it since it involves random thresholding and delivers the same coverage as the input sets, a feature that we do not know how to obtain without randomization.

It is not unusual for each interval to be assigned distinct ``weights'' (importances) in the voting procedure. This can occur, for instance, when prior studies empirically demonstrate that specific methods for constructing uncertainty sets consistently outperform others. Alternatively, a researcher might assign varying weights to the sets based on their own prior insights.  

Since it is desirable to have as small an interval as possible if coverage \eqref{eq:coverage} is respected, we would like to assign a higher weight to intervals of smaller size. However, the weights must be assigned before seeing the sets.

Assume as before that the sets $C_1, \dots, C_K$ based on the observed data follow the property \eqref{eq:coverage}. In addition, let $w = (w_1, \dots, w_K)$ be a set of weights, such that 
\begin{gather}
\label{eq:weights}
    w_k \in [0, 1] \quad  \mathrm{and} \quad \sum_{k=1}^K w_k = 1.
\end{gather}
These weights can be interpreted as the aggregator's prior belief in the quality of the received sets. A higher weight signifies that we attribute greater importance to that specific interval. 
As before, let $U$ be an independent random variable that is distributed uniformly on $[0,1]$, and let $u$ be a realization.
We then define a new set $C^W$ as:
\begin{equation}
\label{eq:cw}
    C^W := \left\{s \in \cS: \sum_{k=1}^K w_k \id_{\{s \in C_k\}} > \frac{1}{2} + \frac u 2\right\}.
\end{equation}

\begin{theorem}
\label{th:CR}
    Let $C_1, \dots, C_K$ be $K \geq 2$ different sets satisfying property \eqref{eq:coverage}. Then, the set $C^W$ defined in \eqref{eq:cw} has coverage of at least $1-2\alpha$:
    \begin{equation}
    \label{eq: coverage_cr}
    \p\big(c \in C^W\big) \geq 1-2\alpha.
    \end{equation}
    In addition, let $\nu(C^W)$ be the measure associated with the set $C^W$, then
    \begin{equation}
        \label{eq:length_cr}
        \nu(C^W) \leq 2 \sum_{k=1}^K w_k \nu(C_k).
    \end{equation}
\end{theorem}

If the weights are equal to $w_k = 1/K$ for all $k\in\mathcal K$, then the set $C^W$ coincides with the set $C^R$ defined in \eqref{eq:cr} and it is a subset of that in \eqref{eq:cm}. This means that in the case of a democratic (equal-weighted) vote $C^R \subseteq C^M$, since $C^M$ is obtained by choosing $u=0$.
Furthermore, \eqref{eq:length_cr} says that the measure of the set obtained using the weighted majority method cannot be more than twice the average measure obtained by randomly selecting one of the intervals with probabilities proportional to $w$. When there are only two sets and $u=0$, the weighted majority vote set will correspond to the set with the greater weight.

\begin{proposition}
    If $C_1, \dots, C_K$ are $K \geq 2$ different sets having coverage $1-\alpha_1, \dots, 1-\alpha_K$ (possibly unknown), then the set $C^W$  defined  in \eqref{eq:cw} has coverage
    \begin{equation*}
        \p(c \in C^W) \geq 1-2\sum_{k=1}^K w_k \alpha_k.
    \end{equation*}
    In particular, this implies that the majority vote of asymptotic $(1-\alpha)$ intervals has asymptotic coverage at least $(1-2\alpha)$.
\end{proposition}

The proof is identical to that of Theorem~\ref{th:CR}, with the exception that the expected value for the variables $\phi_k$ is equal to $\alpha_k$, and is thus omitted. If the $\alpha_k$ levels are known (which they may not be, unless the agents report it and are accurate), and if one in particular wishes to achieve a target level $1-\alpha$, then it is always possible to find weights $(w_1, \dots, w_K)$ that achieve this as long as $\alpha/2$ is in the convex hull of $(\alpha_1, \dots, \alpha_K)$.

Suppose now that we have $K$ different \emph{confidence sequences} for a parameter that need to be combined into a single confidence sequence.
For this setting we show a simple result:

\begin{proposition}
    Given $K$ different $1-\alpha$ confidence sequences for the same parameter that are being tracked in parallel, their majority vote set is a $1-2\alpha$ confidence sequence.
\end{proposition}

It may not be initially apparent how to deal with the time-uniformity in the definitions. The proof proceeds by first observing that an equivalent definition of a confidence sequence is a confidence interval that is valid at any arbitrary stopping time $\tau$ (here the underlying filtration is implicitly that generated by the data itself).
The sequence $(C_k^{(t)})_{t \geq 1}$ is a confidence sequence if and only if for every stopping time $\tau$,
$
\p(c \in C_k^{(\tau)}) \geq 1-\alpha,
$
for all $k\in\mathcal K$. Now, the proof follows by applying our earlier results.

Now, let  the data be fixed, but consider the setting where an unknown number of confidence sets arrive one at a time in a random order, and need to be combined on the fly. Now, we propose to simply take a majority vote of the sequences we have seen thus far. Borrowing terminology from earlier, denote by
\begin{equation}\label{eq:sequential_mer}
C^E(1:t) := \bigcap_{i=1}^t C^M(1:i).    
\end{equation}

We claim that the above sequence of sets  is actually a $(1-2\alpha)$ confidence {sequence} for $c$.
\begin{theorem}
\label{thm:exch2}
    Given an exchangeable sequence of confidence sets $C_1, C_2,\dots$ (or confidence sets arriving in a uniformly random order), the sequence of sets formed by their ``running majority vote'' $(C^E(1:t))_{t \geq 1}$ is a $(1-2\alpha)$ confidence sequence:
    \[
    \p\left(\exists t \geq 1: c \notin C^{E}(1:t)\right)\leq 2\alpha.
    \]
\end{theorem}

Such a result is useful when derandomizing a statistical procedure, by repeating it many times one by one, and attempting to combine the results of these repetitions on the fly.

\subsection{Median of Midpoints and Median of Median of Means}

A special case for the majority vote set in \eqref{eq:cm} occurs if the input sets are (or can be bounded by) intervals of the same width, and can thus each can be represented as their midpoint plus/minus a half-width, then the following rule results in a more computationally efficient procedure. We call this method the ``Median of Midpoints''.

\begin{theorem}\label{thm:median-of-midpoints}
    Suppose the input sets $C_1,\dots,C_K$ are intervals having the same width, and let $C^M$ be their majority vote set. Let the midpoint of $C_k$ be denoted by $c_k$. 
    If $K$ is odd, let $c_{(\lceil K/2 \rceil)}$ denote their median, and let $C^{(K/2)}$ denote the interval whose midpoint is $c_{(\lceil K/2 \rceil)}$. 
    If $K$ is even, let $C^{(K/2)}$ be defined by the intersection of the sets whose midpoints are $c_{(K/2)}$ and $c_{(1+ K/2)}$.
    Then, $C^{(K/2)} \supseteq C^M$ and is hence also a $1-2\alpha$ uncertainty set. Further, $C^{(K/2)}$ has at most the same width as the input sets. If $\bigcap_{k=1}^K C_k \neq \emptyset$, then $C^{(K/2)}=C^M$.
 \end{theorem}

Since the paper has so far focused on uncertainty sets, we first describe a general result that derandomizes point estimators ``directly''.

\begin{corollary}
\label{thm:point}
    Suppose $\hat \theta_1,\dots,\hat \theta_K$ are $K$ univariate point estimators of $\theta$ that are each built using $n$ data points and satisfy a high probability concentration bound:
    \[ 
    \p( |\hat \theta_k - \theta| \leq w(n,\alpha)) \geq 1-\alpha,
    \] 
    for some function $w$. Then, their median $\theta_{(\lceil K/2 \rceil)}$ satisfies
    \[ 
    \p(|\hat \theta_{(\lceil K/2 \rceil)} - \theta| \leq w(n,\alpha)) \geq 1-2\alpha. 
    \] 
    Further, if $\hat \theta_1,\hat \theta_2,\dots,\hat \theta_K, \hat\theta_{K+1} \dots$ are exchangeable, then
    \[
    \p(\forall K \geq 1: |\hat \theta_{(\lceil K/2 \rceil)} - \theta| \leq w(n,\alpha)) \geq 1-2\alpha. 
    \]
\end{corollary}

We present an application of this idea to derandomizing the median of means estimator.
The aim of the celebrated median-of-means procedure is to produce estimators of the mean that have sub-Gaussian tails, despite the underlying unbounded data having only a finite variance. The setup assumes $n$ iid data points from an unknown univariate distribution $P$ with unknown finite variance, whose mean $\mu$ we seek to estimate. The method works by randomly dividing the $n$ points into $B$ buckets of roughly $n/B$ points each. One takes the mean within each bucket, and then calculates the median $\hat \mu^{\mathrm{MoM}}$ of those numbers. 

The method does not produce confidence intervals for the mean; indeed, one needs to know a bound on the variance $\sigma^2$ for any nonasymptotic CI to exist. But the point estimator $\hat \mu^{\mathrm{MoM}}$ satisfies the following sub-Gaussian tail bound: for any $t \geq 0$,
\[
\p(|\hat \mu^{\mathrm{MoM}} - \mu| \geq C \sigma \sqrt{t/n}) \leq 2 e^{-t},
\]
for a constant $C = \sqrt{\pi} + o(1)$, where the $o(1)$ term vanishes as $B, n/B \to \infty$.
This is a much faster rate than the $1/t^2$ on the right-hand side obtained for the sample mean via Chebyshev's inequality. However, one drawback of the method is that it is randomized, meaning that it depends on the random split of the data into buckets.

Here, we point out that Theorem~\ref{thm:point} yields a simple way to derandomize $\hat \mu^{\mathrm{MoM}}$. We simply repeat the median-of-means procedure independently $K$ times to obtain exchangeable estimators $\{\hat \mu^{\mathrm{MoM}}_k\}_{k=1}^K$, and then report the ``median of median of means'' (MoMoM):
\[
\hat \mu^{\mathrm{MoMoM}}_K := \hat \mu^{\mathrm{MoM}}_{\lceil K/2 \rceil}.
\] 

Clearly, $\{\hat \mu^{\mathrm{MoM}}_k\}_{k=1}^K$ form an exchangeable set, whose empirical distribution will stabilize as $K \to \infty$, and thus whose median will also be stabilize for large $K$. 

\begin{corollary}
    Under the setup described above, the \emph{median of median of means} satisfies a nearly identical sub-Gaussian behavior as the original:
\[
\p\left(\exists K \geq 1: |\hat \mu^{\mathrm{MoMoM}}_K - \mu| \geq C \sigma \sqrt{t/n}\right) \leq 4 e^{-t}.
\]
\end{corollary}

The simultaneity over $K$ allows us to choose any sequence of constants $k_1 \leq k_2 \leq \dots$ and produce the sequence of estimators $\hat \mu^{\mathrm{MoMoM}}_{k_1}, \hat \mu^{\mathrm{MoMoM}}_{k_2}, \dots$, and plot this sequence as it is being generated. We can then stop the derandomization  whenever the plot appears to stabilize, with the guarantee now holding for whatever data-dependent random $K$ we stopped at. Clearly, one can use the same technique to derandomize other random estimators by taking their median while ``importing'' their nonasymptotic bounds. We will see an empirical example in the next section.

\subsection{Empirical examples}

We now present a simulation to investigate the relative performance of our proposed methods. We apply the majority vote procedure on simulated high-dimensional data with $n=100$ observations and $p=120$ regressors. Specifically, we simulate the design matrix $\mathbf{X} \in \R^{n \times p}$, where each column is independent of the others and contains standard normal entries. The outcome vector is equal to $\mathbf{y} = \mathbf{X}\boldsymbol{\beta} + \boldsymbol{\varepsilon}$, where $\boldsymbol{\beta}$ is a sparse vector with only the first $m=10$ elements different from 0 (generated independently from   $\mathrm{N}(0, 4)$) while $\boldsymbol{\varepsilon} \lawis \mathrm{N}(0, I_n)$. A test point $(\mathbf{x}_{n+1}, y_{n+1})$ is generated with the same data-generating mechanism. At each iteration we estimate the regression function using the Lasso algorithm with penalty parameter $\lambda$ varying over a fixed sequence of values $K=20$ and then construct a conformal prediction interval for each $\lambda$ in $\mathbf{x}_{n+1}$ using the split conformal method presented in package R \texttt{ConformalInference}. We run  $10,000$ iterations at $\alpha=0.05$. 
We then merge the $K$ different sets using the method described in \eqref{eq:cw} with $w_k=1/K$ for  $k\in \cK$. 

\begin{figure}[h!]
    \centering
    \includegraphics[width=0.95\textwidth]{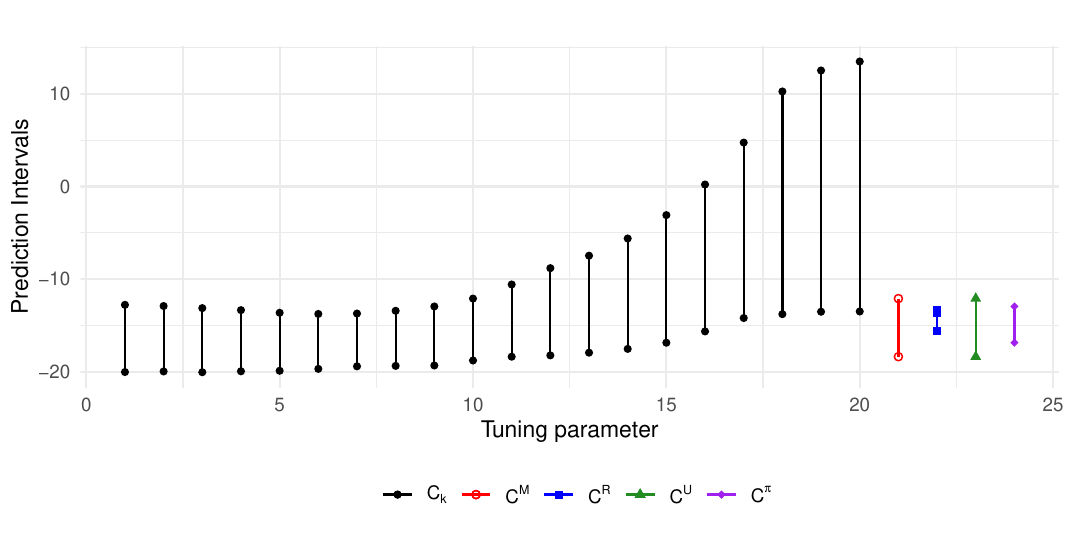}
    \caption{Intervals obtained using different values of $\lambda$, $C^M(x)$, $C^R(x), C^U(x)$ and $
    C^\pi(x)$. For standardization, the value of $U$ in the randomized thresholds is set to $1/2$. The smallest set $C^R$. Since $U=1/2$, the sets $C^M$ and $C^U$ coincides.}
    \label{fig:lasso_sims}
\end{figure}

An example of the result is shown in Figure~\ref{fig:lasso_sims}. The empirical coverages of the intervals $C^M$ and $C^R$ are $$\frac{1}{B}\sum_{b=1}^B \id_{\{y_{n+1}^b \in C^M_b(\mathbf{x}_{n+1})\} }=0.97 \mbox{ ~and~ }\frac{1}{B}\sum_{b=1}^B \id_{\{y_{n+1}^b \in C^R_b(\mathbf{x}_{n+1})\}}=0.92.$$ By definition, the second method produces narrower intervals while maintaining the coverage level $1-2\alpha$. As explained in the previous sections, by inducing exchangeability through permutation, it may be possible to enhance the majority vote results. In fact, the empirical coverage of the sets $C^\pi$ is $0.93$ while the sets are smaller than the ones produced by the simple majority vote. Furthermore, we tested the sets $C^U$ defined in \eqref{eq:cu} and obtained an empirical coverage equal to $0.96$, which is very close to the nominal level $1-\alpha$. In all five cases, the occurrence of obtaining an union of intervals as output is very low, specifically less than 1\% of the iterations.

We now present the results on a real dataset on Parkinson's disease \citep{parkinson_data}. The goal is to predict the total UPDRS (Unified Parkinson's Disease Rating Scale) score using a range of biomedical voice measurements from people suffering early-stage Parkinson's disease. We used split conformal prediction and $K=4$ different algorithms (linear model, lasso, random forest, neural net) to obtain the conformal prediction sets. In particular, we choose $n=5000$ random observations to construct our intervals and the others $n_0=875$ observations as test points. Prior weights also in this case are uniform over the $K$ models that represent the different agents, so a priori all methods are of the same importance. If previous studies had been carried out, one could, for example, put more weight on methods with better performance. Otherwise, one can assign a higher weight to more flexible algorithms such as random forest or neural net. 

The results are reported in Table~\ref{tab:parkinson} where it is possible to note that all merging procedures obtain good results in terms of length and coverage. In addition, also the randomized vote obtains good results in terms of coverage, with an empirical length that is slightly larger than the one obtained by the neural net. The percentage of times that a union of intervals or an empty set is outputted is nearly zero for all methods. In this situation, the intervals produced by the random forest outperform the others in terms of size of the sets; as a consequence, one may wish to put more weight into the method, which results in smaller intervals on average.

\begin{table*}[h!]
\centering
\resizebox{\textwidth}{!}{
\begin{tabular}{ccccccccc}
  \hline
  Methods & LM & Lasso & RF & NN & $C^M$ & $C^R$ & $C^U$ & $C^\pi$ \\ 
  \hline
  Coverage & 0.958 & 0.960 & 0.949 & 0.961 & 0.951 & 0.923 & 0.961 & 0.918  \\ 
  Width & 40.143 & 40.150 & 13.286 & 32.533 & 29.508 & 20.620 & 32.544 & 20.710 \\ 
   \hline
\end{tabular}
}
\caption{Coverage and length of the methods for the Parkinson’s dataset.}
\label{tab:parkinson}
\end{table*}

As previously mentioned, a potential method to derandomize the median of means estimator is to repeat the procedure $K$ times and subsequently take the median of the estimators. The natural question that arises is: how large should $K$ be in practice to achieve a derandomized result? Indeed, considering computational and time resources, one would prefer a small value for $K$. We conducted a simulation study to study the impact of $K$. 

We simulated observations from a standard Student's t-distribution with 3 degrees of freedom and from a skew-t distribution with 3 degrees of freedom and 1 as a skewness parameter. The number of batches is equal to $B=21$ and the number of data points generated are $n=210$. For each iteration, we computed $\hat\mu^\mathrm{MoMoM}_1, \dots, \hat\mu^\mathrm{MoMoM}_K$ with $K = 70$ and considered the absolute value of the difference $\hat\mu^\mathrm{MoMoM}_k-\hat\mu^\mathrm{MoMoM}_{k-1}$ as a measure of stability. 

\begin{figure}
    \centering
\includegraphics[width=\textwidth]{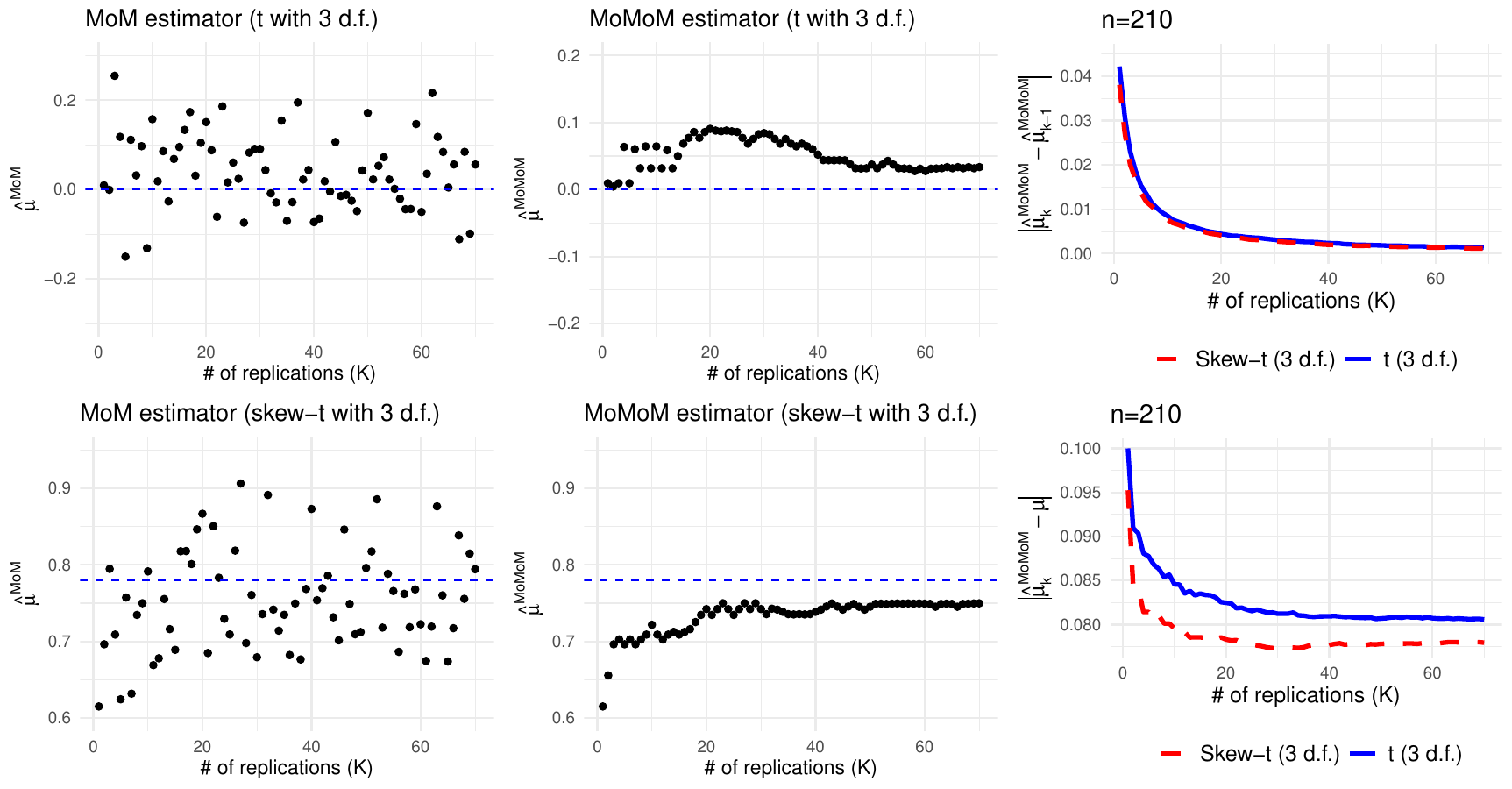}
    \caption{First column: MoM estimator obtained during various replications in a \emph{a single run of the procedure} using data generated from the t-distribution (above) and the skew-t distribution (below). The dashed line is the true mean $\mu$. Second column: MoMoM estimator obtained during the replications; both series stabilize for $k>40$. Third column: average over 1000 runs of the absolute difference $|\hat\mu^\mathrm{MoMoM}_k-\hat\mu^\mathrm{MoMoM}_{k-1}|$ against $k$ (above) and average over 1000 runs of the absolute difference $|\hat\mu^\mathrm{MoMoM}_k-\mu|$ against $K$. Note that the third column is only available to us in the simulation, but plots like the second column can be constructed on the fly as $K$ increases, and it can be adaptively tracked and stopped, while retaining the same statistical guarantee at the stopped $K$.}
    \label{fig:momom_diff}
\end{figure}

In the first two columns of Figure~\ref{fig:momom_diff}, we present two examples of the procedure, where it can be observed that in the first case, $\hat\mu^\mathrm{MoMoM}$ tends to stabilize after 50 iterations, while in the second case, where the data are generated from a skew-t distribution, 25 iterations are sufficient.
In the third column of Figure~\ref{fig:momom_diff}, we report the empirical average over 1000 replications of $|\hat\mu^\mathrm{MoMoM}_k-\hat\mu^\mathrm{MoMoM}_{k-1}|$. We observe an inflection point around 25 for both distributions, and it seems that once this threshold is reached, the gain becomes negligible.

\section*{Bibliographical note}

E-confidence intervals were studied by \cite{vovk2023confidence} and \cite{xu2022post}. 
The t-CI example is from~\cite{wang2023anytime}.
Confidence sequences were introduced by \cite{Darling/Robbins:1967a}.
The original BY procedure for false coverage rate control was proposed by~\cite{benjamini2005false}, while the e-BY procedure was proposed by in~\cite{xu2022post}. 

\cite{kuncheva2003limits} studied majority voting for classifiers and derived its miscoverage rate. \cite{solari2022multi} used majority vote to derandomize and stabilize split conformal prediction~\citep{vovk2005algorithmic}. 
Both their results can in turn be improved using the randomized and exchangeable extensions of majority vote, which were proposed by~\cite{gasparin2024merging}. The latter paper also presents several applications to derandomization of other statistical procedures based on sample splitting such as median of means~\citep{lugosi2019mean} and HulC~\citep{kuchibhotla2024hulc}. The weighted variants of majority vote were exploited by~\cite{gasparin2024conformal} in the development of an online conformal prediction algorithm that aggregates the outputs of many models in an online manner, upweighting the good models over time.

The median of means estimator stems back to, at least, a book by~\cite{nemirovskij1983problem}, but some modern references include~\cite{devroye2016sub}, \cite{lugosi2019mean}. \cite{minsker2022u} proposed a derandomized variant that is computationally intensive: it involves taking the median of all possible sets of size $n/B$. The code to reproduce experiments on derandomizing median of means (and conformal prediction) is at 
\begin{equation}
\mbox{\texttt{https://github.com/matteogaspa/MergingUncertaintySets}}.
\end{equation}

\chapter{Further contrasts between e-values and p-values}
\label{chap:theoretic}

In this chapter, we present a few  results  strengthening the theoretical foundation of e-values. This includes a duality between e-values and p-values,
  representations of e-variables and p-variables as conditional expectations and probabilities,
  and equivalent conditions for the existence of powered   exact p-values and e-values.

\section{A duality between e-values and p-values}
\label{sec:c2-como}

We present a duality between e-values and p-values, further clarifying their intimate connection.

We consider a simple setting where a real-valued test statistic $X$ is available, and a larger value of $X$ represents evidence against the null hypothesis.
Simple examples are \eqref{eq:example-gaussian} and \eqref{eq:example-bernoulli} in
Section~\ref{sec:LR-e-variable}, where the test statistic $X$ is the sample sum.  
In this context, an e-variable should be an increasing function of $X$, and a p-variable should be a decreasing function of $X$. 
Denote by $$S(x;\p)=\p(X \ge x)$$ for $x\in \R$,
which is the left-continuous survival function of $X$ under $\p$.

The next proposition explains that under the monotonicity restriction,  $S(X;\p)$ provides natural bounds on both p-variables and e-variables.

\begin{proposition}
\label{prop:ep-duality1}
Suppose that $E$ is an e-variable for $\p$ that is an increasing function of $X$,
and $P$ is a p-variable for $\p$ that is a decreasing function of $X$. 
Then, $$P\ge S(X;\p) \mbox{~~and~~} E\le 1/S(X;\p) \mbox{~~~$\p$-almost surely}.$$
\end{proposition}
\begin{proof}[Proof.]
For the statement on the p-variable $P$,
first note that for any decreasing function $T$ and an iid copy $X'$ of $X$, we have 
$\{T(X')\le  T(X)\}\supseteq 
\{X'\ge X\}$, which implies
$$
\p(T(X')\le  T(X) \mid X ) \ge \p(X'\ge X\mid X)=S(X;\p) \mbox{~~~$\p$-almost surely.} 
$$
Using this and the fact that the cdf $F_P$ of $P$ under $\p$ is smaller than or equal to the identity on $[0,1]$ (this is the defining property of a p-variable), we get
$P\ge F_P(P) \ge   S(X;\p)$.

The statement on the e-variable $E$ follows directly by observing that $1/E$ is a p-variable, guaranteed by Proposition~\ref{prop:e-to-p}, and then applying the result in the first part.
\end{proof}

In fact, in Proposition~\ref{prop:ep-duality1}, $S(X;\p)$ is a p-variable for $\p$, but $1/S(X;\p)$ is not an e-variable for $\p$.
Nevertheless, in the next result we will see that, $1/S(X;\p)$ is the supremum of all e-variables under the monotonicity condition.  
 
\begin{theorem}[Duality between e-values and p-values]\label{th:ep-duality2}
Suppose that    probability measures in $\mathcal P$ are absolutely continuous with respect to a reference measure $\leb$.  Let $\mathcal E^X$ be the set of all e-variables for $\cP$ that are increasing functions of  $X$,
 and $\mathcal U^X$ be the set of all p-variables for $\cP$ that are decreasing functions of  $X$. 
	Then, $$ \sup_{E\in \mathcal E^X}E  = \inf_{\p\in \mathcal P} \frac{1}{S(X;\p)}
  = \frac{1}{\inf_{P\in \mathcal U^X} P} 
 ~~~\mbox{$\leb$-almost surely}.$$
 In particular, if $\mathcal P=\{\p\}$, then $\sup_{E\in \mathcal E^X}E = 1/S(X;\p)$, $\p$-almost surely.
\end{theorem}
\begin{proof}[Proof.]
Proposition~\ref{prop:ep-duality1}  and the fact that $1/E$ for an e-variable $E$ is a  p-variable guarantee 
\begin{equation}\label{eq:ep-duality3}
\sup_{E\in \mathcal E^X}E  \le \sup_{P\in \mathcal U^X} \frac{1}{ P} \le  \inf_{\p\in \mathcal P} \frac{1}{S(X;\p)}   ~~~\mbox{$\leb$-almost surely}.
\end{equation}    
Below we will show the equalities. 
	For $s \in \R$, let $$E_s = \frac1 {\sup_{\p\in \mathcal P} S(s;\p)}\id_{\{X\ge s\}} \mbox{~ where we set $E_s=1$ in case of $0/0$},$$
	which is  increasing in $X$.
 Moreover, $E_s$ is an e-variable, because $\E^\p[E_s]= S(s;\p)/{\sup_{\p'\in \mathcal P} S(s;\p')} \le 1$ for all $\p\in \mathcal P$. 
	It is straightforward to see that the supremum over $s\in \R$ for $E_s$ is achieved at $s=X$, that is,
	$$\sup_{s\in \R}E_s  = \frac{1}{\sup_{\p\in \mathcal P}  S(X ; \p) }$$ 
	and  hence  the equalities in \eqref{eq:ep-duality3} hold.
\end{proof}

The quantity 
$\sup_{\p\in \mathcal P} {S(X;\p)}$
in Theorem~\ref{th:ep-duality2} is a p-variable for $\mathcal P$, again the best one under the monotonicity restriction, as it is equal to  $\inf_{P\in \mathcal U^X} P$.
Theorem~\ref{th:ep-duality2} suggests that  the supremum of e-values over a set $\mathcal E\subseteq \mathcal E^X$ can be seen as a conceptual middle ground between an e-value ($\mathcal E$ is a singleton)  and the reciprocal of a p-value ($\mathcal E=\mathcal E^X$ in Theorem~\ref{th:ep-duality2}, which is the maximal set).

The next corollary strengthens the Markov inequality for e-variables, which
is immediate from Theorem~\ref{th:ep-duality2}.

\begin{corollary}\label{coro:c2-sup}
For any set $\mathcal E$ of e-variables for $\mathcal P$ that are increasing functions of $X$, we have $$\sup_{\p\in\mathcal P} \p\left(\sup_{E\in \mathcal E} E\ge \frac{1}{\alpha}\right)\le \alpha~~~\mbox{ for all $\alpha \in (0,1)$}.$$
\end{corollary}

If increasing monotonicity of $E$ in $X$ is replaced by decreasing monotonicity (and analogously for $P$), similar conclusions in this section hold true, with the left-continuous survival function of $X$ replaced by its cdf. 
Indeed, what matters for the type-I error control in Corollary~\ref{coro:c2-sup} is that e-variables are \emph{comonotonic}, meaning that they are increasing functions of the same random variable (in our context, it is $X$ or $-X$, but it can be anything). 
 Section~\ref{sec:c12-como} has a more general treatment of this property.

\section{Representations as conditional expectations and probabilities}
\label{sec:c2-cond-exp}

In this section, we present two representation results on p-variables and e-variables, which helps to further explain ``e'' in ``e-values''.

In what follows, $X$ represents a generic data sample, which takes values in $\R^n$; formally $X$ is a measurable mapping from $\Omega$ to $\R^n$. In fact, the space $\R^n$ does not matter in this section. We say that a p-variable or an e-variable $Y$ is built on $X$ if it is a measurable function of $X$. 
 In the language of measure theory, this means that $Y$ is measurable with respect to the $\sigma$-algebra of $X$.
\begin{theorem}[Representation theorem for e-variables]
\label{th:rep-e}
For a random variable $E$, the following are equivalent:
\begin{enumerate}[label=(\roman*)]
    \item $E$ is an e-variable for $\mathcal P$ built on $X$; 
    \item For every $\p\in \cP$, 
 there exists $\q \in \cM_1$ with $ \q \ll \p$ such that  
    \begin{equation}\label{eq:c2-rep1}
    E\le  \E^\p\left[ \frac{\d \q }{\d \p}\mid X\right]\mbox{~~~~~$\p$-almost surely}.
    \end{equation}
\end{enumerate}
If all  $\p\in \mathcal P$ are equivalent to some measure $\leb$, then every e-variable built on $X$  is dominated by  (meaning $\le$)  
$$
 \inf_{\p\in \mathcal P} \E^\p\left[ \frac{\d \q_\p }{\d \p}\mid X\right],~~~\mbox{which is an e-variable for $\mathcal P$},
$$
for some class $(\q_\p)_{\p\in \cP}$ in $\cM_1$ absolutely continuous with respect to $\leb$.
\end{theorem}

 For given $\p,\q\in \cM_1$,
the conditional expectation in \eqref{eq:c2-rep1} often takes the form
$$\E^\p\left[ \frac{\d \q }{\d \p}\mid X\right]=\frac{p(X)}{q(X)},$$
where $p$ and $q$ are the density functions of the data $X$ under $\p$ and $\q$, when they exist, as in Chapter~\ref{chap:ui}.

\begin{proof}[Proof.]
The direction (ii)$\Rightarrow$(i) is straightforward, because for all $\p\in \mathcal P$, the tower property of conditional expectation gives
$$
\E^\p[E] \le \E^\p \left[   \E^\p  \left[ \frac{\d \q}{\d \p}\mid X\right] \right]  =  \E^\p  \left[ \frac{\d \q }{\d \p} \right] = 1.
$$
Next we show (i)$\Rightarrow$(ii). 
For each $\p\in \mathcal P$, let    $\q$ be a probability measure defined by 
$\d \q/\d \p=E/\E^\p[E]$ if $\E^\p[E]>0$,
and $  \q =  \p$ if $\E^\p[E]=0$.
We have, for $\p$ with $\E^\p[E]>0$, 
$$
  \E^\p  \left[ \frac{\d \q}{\d \p }\mid X\right]=  \E^\p  \left[ \frac{E}{\E^\p[E]}\mid X\right]    =   \frac{E}{\E^\p[E]} \ge E \mbox{~~~~~$\p$-almost surely}.
$$
For  $\p$ with $\E^\p[E]=0$, we have $E=0$ $\p$-almost surely, which trivially satisfies \eqref{eq:c2-rep1}.

The last claim follows by taking an infimum over $\p$ in \eqref{eq:c2-rep1}. 
\end{proof}

 Theorem~\ref{th:rep-e} generalizes
Proposition~\ref{prop:admissible-lr}, which   connects e-variables to likelihood ratios in a simple null setting. 
Theorem~\ref{th:rep-e} implies that all e-variables built from a dataset are indeed the conditional expectation of some likelihood ratio. 
In this sense, e-values are not only defined by constraints on expectations, but \emph{they are indeed conditional expectations on the data}. Moreover, Proposition~\ref{prop:cond-log-opt} gives the optimality of this construction when testing a simple null $\p$ against $\q$, where $\q$ in \eqref{eq:c2-rep1} is the alternative hypothesis.

The next result shows that \emph{p-variables are conditional probabilities on the data}.

\begin{theorem}[Representation theorem for p-variables]
\label{th:rep-p}
For a random variable $P$, the following are equivalent:
\begin{enumerate}[label=(\roman*)]
    \item $P$ is a p-variable for $\mathcal P$ built on $X$; 
    \item There exists a measurable function $T: \R^n \to \R$ such that for every $\p\in \cP$,   \begin{equation} P\ge  \p\left(T(X')\le T(X)\mid X\right) \mbox{~~~~~$\p$-almost surely},\label{eq:c2-rep2}
\end{equation} 
    where $X'$ is an iid copy of $X$ under $\p$.
\end{enumerate} 
If all  $\p\in \mathcal P$ are equivalent  to some measure $\leb$, then every p-variable built on $X$  is dominated by  (meaning $\ge$) 
 $$
  \sup_{\p\in \mathcal P} \p\left(T(X'_\p)\le T(X)\mid X\right) \mbox{~~~which is a p-variable for $\mathcal P$}, $$
for some measurable function $T: \R^n \to \R$,  where $X'_\p$ is an iid copy of $X$ under $\p$.
\end{theorem}

\begin{proof}[Proof.]
We first show (ii)$\Rightarrow$(i).
Fix any $\p\in \mathcal P$, and denote by $F_T$ the cdf of $T(X)$ under $\p$. 
Note that 
$\p(T(X')\le T(X)\mid X) = F_T(T(X))$. 
Hence, $\p(P\le \alpha) \le \p(F_T(T(X))\le \alpha ) \le \alpha, $ a well-known fact in probability. This shows that $P$ is an e-variable for $\mathcal P$.

Next, we show (i)$\Rightarrow$(ii). 
Take $T(X)=P$, treating $P$ as a function of $X$. Fix any $\p\in \mathcal P$, and denote by $F_P$ the cdf of $P$ under $\p$.  
Since $P$ is a p-variable, $F_P$ is by definition smaller than or equal to the identity on $[0,1]$. 
It follows that $P\ge F_P(P) = \p(T(X')\le T(X)\mid X).$

The last claim follows by taking a supremum over $\p$ in \eqref{eq:c2-rep2}. 
\end{proof}

\section{Existence of powered exact e-values and p-values}
\label{sec:existence-exact}

By considering finite hypotheses, 
we can derive another existence result, similar to Section~\ref{sec:existence}, in the case that we require the e-variable to be \emph{exact} (meaning that its expectation equals one under every $\p \in \cP$), which turns out to need the same geometric condition as requiring that the p-variable is exact (meaning it is uniformly distributed for every $\p \in \cP$).  
For a given  $\cP$, 
 a random variable
$X$ is \emph{pivotal} 
if it has the same distribution under all $\p\in \cP$.

\begin{theorem}\label{th:c3-span}
Consider testing a finite $\cP=\{\p_1,\dots,\p_L\}$ against a finite $\cQ=\{\q_1,\dots,\q_M\}$. If \eqref{assm:U} holds, the following are equivalent:
    \begin{enumerate}[label=(\roman*)]
    \item there exists an exact (hence pivotal)   p-variable that is powered against $\cQ$;
        \item there exists a pivotal, exact,  bounded e-variable that has positive e-power against $\cQ$;
\item there exists an exact  e-variable that is powered against $\cQ$;
        \item it holds that $\conv(\q_1,\dots,\q_M)\cap \mathrm{Span}(\p_1,\dots,\p_L) = \emptyset$.
    \end{enumerate}
\end{theorem}

Compared to the Theorem~\ref{th:existence}, the convex hull of the null has been replaced by the span, defined as 
$$\mathrm{Span}(\p_1,\dots,\p_L)
=
\left\{\sum_{\ell=1}^L \lambda_\ell\p_\ell: \lambda_\ell \in \R \mbox{ for each $\ell\in [L]$}\right\},
$$
which is a set of signed measures on the sample space (the convex hull would only allow each $\lambda_\ell\ge 0$ and $\sum_{\ell=1}^L\lambda_\ell =1$).

Though we do not provide a full proof of the theorem, which requires modern tools from optimal transport, the appearance of the span can be intuitively justified. If an exact e-variable $E$ is powered against $\cQ$, then for every $\q \in \conv(\cQ)$, $\E^{\q}[E] > 1$. But if such a $\q$ can be written as a linear combination  $\sum_{\ell=1}^L \lambda_\ell\p_\ell$ of elements of $\cP$ (which implies 
 $\sum_{\ell=1}^L\lambda_\ell =1$ because $\cQ$ is a probability), we get a contradiction, because we know that $\E^{\p}[E] = 1$ for every $\p \in \cP$. This proves the direction (iii)$\Rightarrow$(iv) above. 


\section*{Bibliographical note}

The duality in Section~\ref{sec:c2-como} is based on some results in \cite{blier2024improved}. 
The observations in 
Section~\ref{sec:c2-cond-exp} were made in
\cite{wang2021ruodu}. 
 Section~\ref{sec:existence-exact} summarizes results from~\cite{zhang2023exact}. In particular, Theorem~\ref{th:c3-span}  was proved  in \cite{zhang2023exact} using tools from optimal transport theory, and with a weaker assumption than \eqref{assm:U}, called joint non-atomicity, which does not require randomization.

\chapter{Improved  e-value thresholds under additional conditions}
\label{chap:shape}

Fix an atomless probability measure $\p$ for which all e-variables in this chapter are defined. 
Similarly to the case of Chapters~\ref{chap:multiple} and~\ref{chap:combining}, it suffices to consider the null $\{\p\}$, and all results  can be easily translated into the case of general composite nulls. 

The Markov inequality  in Proposition~\ref{prop:markov} gives  
$
 \p(E\ge 1/\alpha) \le \alpha $ for any $E$ in the set  $\fE$  of all e-variables, and it cannot be improved without further assumptions.  Nevertheless, the probability $\p(E>1/\alpha)$ can be improved if the e-variable $E$ is constrained in a subset $\mathcal E\subseteq \fE$ of e-variables, which we will study in this section. 
The set $\mathcal E$ will be chosen to satisfy some distributional conditions to be specified later. 

The application context of results in this section is when the tester does not know the 
distribution of the e-value (e.g., due to its complicated or black-box design), but have some shape information about the e-value. 

We consider a simple hypothesis in this section, but if a composite null hypothesis is considered, then it suffices to assume the corresponding shape information for all distributions in the null hypothesis. 

\section{Conditional e-to-p calibrators on a subset of e-values}
\label{sec:c12-1}
 
We are interested in the quantity  $R_{\gamma}(\mathcal E)$ for $\gamma >0 $, defined by
$$
R_{\gamma}(\mathcal E)= \sup_{E\in \mathcal E} \p(E\ge 1/\gamma ),
$$  
that is, the largest probability that $\p(E\ge 1/\gamma )$ can attain for $E\in \mathcal E$.
 Hence, for any set $\mathcal E$ of e-variables, it holds that $R_{\gamma}(\mathcal E) \le  R_{\gamma}(\fE)=\gamma$ for  $\gamma\in (0,1]$.

 We are interested only in the case  $\gamma\in (0,1]$. For $\gamma>1$,
 it is usually the case that $R_{\gamma}(\mathcal E)=1$ because for most classes $\mathcal E$ that we consider, either $1\in \mathcal E$ or $1$ is the limit of elements of $\mathcal E$. 
 

To build a level-$\alpha$ test using the e-variable $E$,   
one needs to find a threshold $t>0$,   better to be smaller, such that 
 $
 \p(E \ge t) \le \alpha.
  $
  For this, we intuitively should use  the smallest $t$ such that 
 $R_{1/t}(\mathcal E)\le \alpha$, which satisfies  $R_{1/t}(\mathcal E)=\alpha$
in case   $\gamma\mapsto R_{\gamma}(\mathcal E)$  is continuous.
Because of the usual interpretation of an e-variable  $E$   that $E\le 1$ carries no evidence against the null hypothesis, we consider $t\ge 1$.


The following result computes the smallest threshold  $t$   for the e-test.  Denote by $q_\beta(X)$ the left $(1-\beta)$-quantile function of $X$ (as in Sections~\ref{sec:c2-LRB} and~\ref{sec:c5-ex}), that is, 
 $$q_\beta(X) =\inf \{x\in \R: \p(X\le x)\ge 1-\beta\} \mbox{~~~for $\beta \in (0,1)$}.$$

Conditions on the quantile function   can be translated into those on the cdf. A summary of these translations is provided in Lemma~\ref{lem:quantile1} in Appendix~\ref{app:quantile}.
 \begin{lemma}\label{lem:1}
 For $\alpha \in (0,1)$, the quantity
 $ 
T_{\alpha}(\mathcal E) : = \inf\{t\ge 1: R_{1/t}(\mathcal E) \le \alpha\}
 $ satisfies 
$$T_{\alpha}(\mathcal E)= \left( \sup_{E\in \mathcal E} q_{\alpha}(E)\right)\vee 1.$$
If $\gamma\mapsto R_{\gamma}(\mathcal E)$  is continuous, then 
 $T_{\alpha}(\mathcal E)$
 is the smallest real number $t\ge 1$ such that  $ \p(E \ge t) \le \alpha$ for all $E\in \mathcal E$.  
 \end{lemma}
   \begin{proof}[Proof.] 
   Since $ R_{\gamma}(\mathcal E)\le \gamma$ for $\gamma \in (0,1]$, 
   the set $\{t\ge 1: R_{1/t}(\mathcal E) \le \alpha\}$ is not empty.
By using  Lemma~\ref{lem:quantile1}, we have  
\begin{align*}
T_\alpha(\mathcal E)  &  \ge \inf \{ t\ge 1 :   \p(E > t) \le \alpha \mbox{ for all $E\in \mathcal E$} \} \\&  =
\inf \{ t\ge 1 :   q_\alpha(E)\le t \mbox{ for all $E\in \mathcal E$} \} 
  = \left( \sup_{E\in \mathcal E} q_{\alpha}(E)\right)\vee 1.
\end{align*}
Take any $\epsilon\in (0,1)$. We have
 \begin{align*}
T_{\alpha}(\mathcal E) -\epsilon  &  =  \inf \{t\ge 1-\epsilon : R_{1/(t+\epsilon)}(\mathcal E) \le \alpha\}   
\\& = \inf \left\{ t\ge 1 -\epsilon :    \sup_{E\in \mathcal E} \p(E \ge  t +\epsilon) \le \alpha  \right\}  
\\& = \inf   \{ t \ge 1  -\epsilon :     \p(E \ge  t+\epsilon) \le \alpha \mbox{ for all $E\in \mathcal E$} \}  
 \\&\le \inf \{ t\ge 1  -\epsilon :   \p(E > t ) \le \alpha \mbox{ for all $E\in \mathcal E$} \}  
\\ &  = \left(\sup_{E\in \mathcal E} q_{\alpha}(E)\right)\vee (1-\epsilon) \le  \left( \sup_{E\in \mathcal E} q_{\alpha}(E)\right)\vee 1\le T_\alpha(\mathcal E) .
\end{align*}
Since $\epsilon\in (0,1)$ is arbitrary, we have $T_\alpha(\mathcal E) =\left( \sup_{E\in \mathcal E} q_{\alpha}(E)\right)\vee 1$, showing  the first statement. 
If  $\gamma\mapsto R_{\gamma}(\mathcal E)$  is continuous, then 
 \begin{align*}
T_{\alpha}(\mathcal E)  &  = \min\{t\ge 1 : R_{1/t}(\mathcal E) \le \alpha\} \\
  & = \min \{ t \ge 1 :     \p(E \ge  t) \le \alpha \mbox{ for all $E\in \mathcal E$} \},
\end{align*}  
showing the second statement.
 \end{proof}

In Section~\ref{sec:2-cali} we have seen  e-to-p calibrators, which are mappings that convert any e-value into a p-value.
  The function $\gamma\mapsto R_{1/\gamma}(\mathcal E)$ serves to refine this concept.  
For a subset $\mathcal E$
  of e-variables, we say that a function $f:[0,\infty]\to[0,\infty)$ is a \emph{conditional e-to-p calibrator on $\mathcal E$} if $f$ is decreasing, and 
  $f(E)$ is a p-variable for all $E\in\mathcal E$.
 Different from Proposition~\ref{prop:e-to-p}, which gives that  
  $x\mapsto (1/x)\wedge 1$
  is the only useful e-to-p 
 calibrator on $\fE$,  we can find better conditional e-to-p 
 calibrators based on $R_{1/\gamma}(\mathcal E)$ than $x\mapsto (1/x)\wedge 1$ for   various subsets  $\mathcal E$ of  $\fE$.
  Moreover,
  the following result implies that any class $\mathcal E$ admits a smallest conditional e-to-p 
 calibrator. This is in sharp contrast to the set of calibrators from p-values to e-values, which does not admit a smallest element.

\begin{proposition}\label{th:calibrator}
    The function $x \mapsto R_{1/x }\left(\mathcal{E}\right)$ on $[0,\infty]$ is a conditional e-to-p  calibrator on $\mathcal{E}$, and it is the smallest such calibrator. 
\end{proposition}
\begin{proof}[Proof.]
    We first show that the function $g:x \mapsto R_{1/x}\left(\mathcal{E}\right) $ is a conditional e-to-p calibrator.
    This means
    $
    \p(g(E)\le \alpha) \le\alpha 
    $
    for all $\alpha \in (0,1)$ and $E\in \mathcal E$.
    By definition, for $x\in [0,\infty]$ and $E\in\mathcal E$, 
    \begin{align*}
        g( x) 
       \ge  \p(E\ge x ) 
    = \p(-E\le -x
        ) = \  F_{-E}(-x)
        ,
    \end{align*}
where $F_{X}$ is the cdf of a random variable $X$. Since $F_{X}(X)$ is a p-variable for any random variable $X$, we have 
$\p(F_{-E}(-E)\le \alpha )\le \alpha$ for all $E\in \mathcal E$. 
Therefore,  for all $\alpha \in (0,1)$, $$\p(g(E)\le \alpha) \le \p(F_{-E}(-E)\le\alpha)\le \alpha.$$

    We now prove that $g$ is the smallest calibrator on $\mathcal{E}$. Suppose 
 that $f$ is another conditional  e-to-p calibrator on $\mathcal{E}$ such that  $f(x)< g(x)$ for at least one $x\in [1,\infty]$.  By definition, there exists   $E\in \mathcal E$
 such that  $
        f(x) < \p(E \geq x).
 $
    This implies  $$\p(f(E) \leq f(x)) \geq \p(E \geq x)> f(x),$$ which violates  the definition of $f(E)$ as a p-variable.
\end{proof}

Proposition~\ref{th:calibrator} implies that by computing $R_{\gamma}(\mathcal E)$ or an upper bound on it, we can find conditional e-to-p  calibrators to convert e-values  realized by elements of $\mathcal E$ into p-values, which work better than the e-to-p 
 calibrator $x\mapsto (1/x)\wedge 1$ on $\fE$. This can be useful in procedures that take p-values as input, such as the Benjamini--Hochberg procedure.

We make a further observation on the probability bound  $R_\gamma(\cE)$ and the threshold $T_\alpha(\cE)$. 
For a class $(E_\theta)_{\theta\in \Theta}$ of e-variables in $\mathcal E$ and a probability measure $\pi$ on $\Theta$, 
let  $E_\pi$ follow the mixture  distribution, that is, 
$\p(E_\pi\le x) = \int_\Theta \p(E_\theta\le x) \pi (\d \theta)$ for $x\in \R$. It is straightforward to see $\p(E_\pi \ge 1/\gamma) \le R_\gamma(\mathcal E)$. 
In other words,  
the probability bounds and thresholds obtained for $\cE$ also apply to the convex hull of $\cE$ with respect to distribution mixtures.

\section{Comonotonic e-variables}

\label{sec:c12-como}

In this section, we discuss $R_{\gamma}(\cE)$ for a set $\cE$ of comonotonic e-variables.
Our results here can be seen as a generalized version of Proposition~\ref{prop:ep-duality1}. 

Consider the simple example in Section~\ref{sec:c1-ex2}, that is, to test $\mathrm{N}(0,1)$ against $\mathrm{N}(\mu,1)$ with $\mu> 0$ for iid observations $X_1,\dots,X_n$. 
A natural e-variable $E_\mu$   is   the likelihood ratio 
\begin{equation}
    \label{eq:c12-LR}
E_\mu =\exp( {\mu S_n -n\mu^2/2}),
\end{equation} 
where $S_n=\sum_{i=1}^n X_i$.  An interesting observation is that $E_\mu$, $\mu>0$ are \emph{comonotonic} random variables.  

\begin{definition}[Comonotonicity]
A set $\cX_C$ of random variables  
is \emph{comonotonic} if there exists a common random variable $Z$ such that each $X\in \cX_C$ is an increasing function of $Z$. We also say that the random variables in $\cX_C$ are comonotonic.
\end{definition}

 The following lemma allows us to analyze   $R_{\gamma}$ for a set of comonotonic e-variables. 

\begin{lemma}\label{lemma:comonotonic-survival}
	Let $x\in \R$ and $\gamma \in [0,1]$. 
	Suppose that $\mathcal X_C$ is a collection of comonotonic random variables
	satisfying $\p(X \ge  x)\le \gamma $ for each $X\in \mathcal X_C$.
	Then $\p (\sup_{X\in \mathcal X_C} X \ge x )\le \gamma.$
\end{lemma}
\begin{proof}[Proof.]

Suppose that all elements of $\cX_C$ are increasing functions of $Z$. 
Note that the sets
$\{X\ge x\}$ for $X\in \cX$ are nested since $\cX_C$ is a set of comonotonic random variables, i.e., each of the set is of the form $\{Z\ge z\}$ or $\{Z> z\}$ for different $z$ and a common $Z$. 
Therefore, we can move the supremum outside  the probability and get 
$$\p\left(\sup_{X\in \mathcal X_C} X \ge x \right)=\sup_{X\in \mathcal X_C} \p( X \ge x )\le \gamma,$$
as claimed.
\end{proof}

The next result follows immediately.
\begin{proposition} \label{prop:sup-comonotonic}
For any collection  $
\mathcal E \subseteq \fE$ of comonotonic e-variables for $\p$, we have $  \p(\sup_{E\in \cE} E  \ge T_\alpha(\mathcal E))\le \alpha$.
\end{proposition} 

In particular, Proposition~\ref{prop:sup-comonotonic} implies  $  \p(\sup_{E\in \cE} E  \ge 1/\alpha)\le \alpha$ for any comonotonic set $\mathcal E$, which we have seen in Corollary~\ref{coro:c2-sup} with the help of a statistic $X$.
As a consequence, $(\sup_{E\in \cE} E)^{-1} $ is a p-variable.

In the setting where e-variables are computed from likelihood ratios, 
we discussed using the mixture of likelihood ratios  in 
Section~\ref{sec:3-mix-plugin}. 
Proposition~\ref{prop:sup-comonotonic} suggests that, under comonotonicity, we can use the supremum of e-variables instead of their mixture, and the supremum has a higher power. 

\begin{example}[Testing normal mean]
In the example of testing $\mathrm{N}(0,1)$ against $\{\mathrm{N}(\mu,1): \mu>0\}$ with $n$ data points, instead of using $E_\mu$ in \eqref{eq:c12-LR}  
for a fixed $\mu>0$ in  or its mixture, we can use
$$
E=\sup_{\mu>0}  E_\mu = \sup_{\mu>0} \exp( {\mu S_n  - n \mu^2/2}) =  \exp\left(\frac{(S_n)_+^2}{2n}\right),
$$ 
which, although not being an e-variable, gives $\p(E \ge 1/\alpha)\le \alpha$ under the null.   
If we are testing $\{\mathrm{N}(\mu,1):\mu\le 0\}$ against $\{\mathrm{N}(\mu,1): \mu>0\}$, the same type-I error guarantee holds true, following from Corollary~\ref{coro:c2-sup}.
\end{example}

The above observation on the normal distributions can be generalized to other families. 
Suppose that we test
$\p=\q_{\theta_0}$ against $\{\q_{\theta} : \theta \in \Theta\}$ with a class of likelihood ratio e-variables
 $$
E_\theta =  \prod_{i = 1}^{n} \frac{  q_{\theta}(X_i)}{p(X_i) }.
$$  
As long the e-variables are increasing or decreasing functions of the same test statistic, 
we can take the supremum over these e-variables, and particular this yields 
 $E=E_{\hat \theta}$, where $\hat \theta$ is  the maximum likelihood estimator.
Although $E$ is not an e-variable, 
the test using $\id_{\{E\ge 1/\alpha\}}$ has level $\alpha$. Note that here the maximum likelihood estimator $\hat \theta$ appears in the numerator, instead of the denominator as in the universal inference in Chapter~\ref{chap:ui} (which yields e-variables).

In addition to the normal distributions, 
testing several other members of the exponential family yields a class of likelihood ratio functions that are monotonic in a common test statistic, and we omit a detailed discussion here.

Other than likelihood ratios, several examples in Section~\ref{sec:c2-ex}   also involve parametric classes $( E_\theta)_{\theta \in \Theta}$ of comonotonic e-variables, and $  \p(\sup_{\theta \in \Theta} E_\theta  \ge 1/\alpha)\le \alpha$ holds for these examples.

\section{Density conditions}

We next study  some common conditions on the density function of the e-variables in $\mathcal E$ to compute $R_{\gamma}(\cE)$ or an upper bound on it.

Three conditions, motivated by different applications, are   modeled by the following sets of e-variables: 
\begin{align*}
\mathcal E_{\rm D}& =\{E\in\fE: \mbox{$E$  has a decreasing density on its support}\},
\\
\mathcal E_{\rm D>1}& =\{E\in\fE: \mbox{$E$ has a decreasing density over $[1, \infty)$}\},
\\
\mathcal E_{\rm U}& =\{E\in\fE: \mbox{$E$  has a unimodal density on $[0,\infty)$}\}.
\end{align*} 
A random variable on $[0,\infty)$ has a unimodal density, or a unimodal distribution, if it has a density function 
that is increasing on $[0,x]$ and decreasing on $[x,\infty)$ for some constant $x\in [0,\infty)$, 
and it can possibly have a point mass at $x$. Such $x$ is called a mode of the random variable. 

The assumption of decreasing densities in $\mathcal E_{\rm D}$ is satisfied by, for instance,  exponential and Pareto distributions. Note that here assumptions are made on the densities of the e-variables themselves, not on those in the hypotheses of the underlying testing problem.

Whenever we consider a density, we include its limit case; that is, we allow $E\in \mathcal E_{\rm D}$ to have a point mass at the left end-point of its support. This convention does not change any results but simplifies many arguments in the proofs. 
The assumption of  decreasing density over $[1, \infty)$ but arbitrary elsewhere   in    $\mathcal E_{\rm D>1}$ is  useful and can be observed from some e-variables obtained from universal inference discussed in Chapter~\ref{chap:ui}. 
The assumption of unimodal density in  $\mathcal E_{\rm U}$  is satisfied by, for instance, log-normal and gamma distributions.  
Moreover, as explained in Section~\ref{sec:c12-1}, our obtained bounds also apply to the mixtures of the   distributions within each class.  
The set of all modes of $E\in\mathcal E_{\rm U}$ may be a singleton or an interval.
Note that $\mathcal E_{\rm D}\subseteq \mathcal E_{\rm U} $,
and hence $R_\gamma (\mathcal E_{\rm U}   ) \ge  R_\gamma (\mathcal E_{\rm D}   )$ for all $\gamma \in (0,1]$.

\begin{theorem}[Decreasing and unimodal densities]\label{th:dec-uni}
For $\gamma \in (0,1]$,
\begin{enumerate}
\item [(i)]
$R_\gamma ({\mathcal E_{\rm D}}  ) =\gamma /2$ if $\gamma \ne 1$ and $R_1 ({\mathcal E_{\rm D}}  ) = 1$;
\item [(ii)]    
    $R_\gamma(\mathcal{E}_{\rm D>1}) =  {\gamma }/({1+ \sqrt{1-\gamma^2}})$; 
    \item[(iii)]  
$R_\gamma (\mathcal E_{\rm U}   ) =(\gamma /2)\vee (2\gamma -1)$.
\end{enumerate}
\end{theorem}

\begin{figure}[t]
    \centering
     \includegraphics[width=0.75\textwidth]{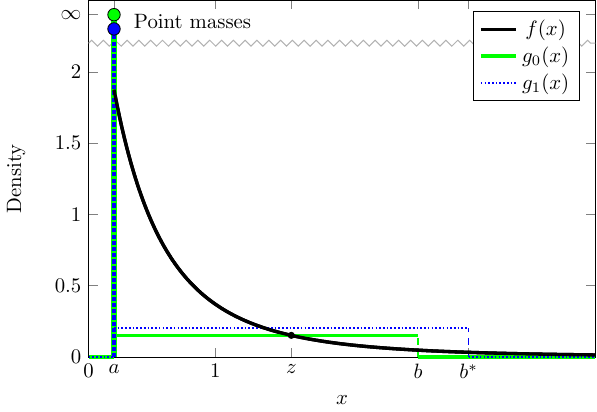}
    \caption{Proof sketch in part (i) of Theorem~\ref{th:dec-uni}. 
    Here, $f$ and $g_0$ have the same area (probability) exceeding $z$, with $g_0$ having a smaller mean and a point mass at $0$; $g_1$ has mean 1 and larger probability of exceeding $z$ than $g_0$.}
    \label{fig:pf-sketch-decreasing}
\end{figure}

\begin{proof}[Proof.]
We only present the proof of part (i) here, which is the simplest, but 
it illustrates the basic techniques of proving more complicated statements. The proof of parts (ii) and (iii) are omitted. 
We first show $\p(E\ge 1/\gamma )\le \gamma/2$ for each $E\in \mathcal E_{\rm D}$ and $\gamma \in (0, 1)$.
Let $f$ be the density function of $E$ and let $a$ be the left end-point of the support of $E$. For a fixed $z > 1$, set $b := z + \p(E \ge z)/f(z)$.
We construct an e-variable $E_0$ with density $g_0$ such that 
(a) $g_0(x) = f(z)$ for $a < x < b$;
(b) $E_0$ has a point mass at $a$ with probability 
$1 - \int_a^\infty g_0(x) \d x$.

We illustrate in Figure~\ref{fig:pf-sketch-decreasing} how to construct $g_0$ from $f$.
Since $f$ has a decreasing density over $[a, \infty)$, we will have $f(x) \ge g_0(x)$ for $x \in [a, z)$ and $f(x) \le g_0(x)$ for $x > z$.
We can construct $g_0$ by shifting the area between $f(x)$ and $f(z)$, for $x \in (a, z)$ to a point mass at $a$ and shifting the area between $f(x)$ and $0$, for $x > b$ to the area between $f(x)$ and $f(z)$, for $x \in (z, b)$. 
Note that $b$ is chosen so that $\p(E_0 \in [z, b)) = \p(E \ge z)$, 
so $\p(E_0 \ge z) = \p(E \ge z)$. 
Further, since we construct $g_0$ by shifting the density of $f$ to the left, we have $\E^\p[E_0] \le \E^\p[E]$.  

Therefore, for any $E \in \mathcal{E}_{\rm D}$ with left end-point of support $a$, there exists an e-variable $E' \in \mathcal{E}_{\rm D}$ such that $E'$ has a point mass at $a$ and a uniform density on $[a, b]$ for some constant $b$, and $\p(E'\ge z) = \p(E \ge z)$. To show $\p(E \ge z) \leq 1/(2z)$, it suffices to look at the collection of e-variables that have a point mass at $a$ and a uniform density on $[a, b]$ for any $b > 1$.
Moreover, it suffices to consider the case $\E^\p[E']=1$.

Note that $\E^\p[E'] = (1 - w)a + w (a + b)/2$ for some $w\in [0, 1]$ and $b \geq 1$. The requirement $\E^\p[E'] = 1$ implies that $w = 2(1 - a)/(b - a)$. Now, we seek $b$ that maximizes $\p(E'\ge z)$. 
We have 
$\p(E' \ge z) = w (b-z)/(b-a) = 2(1-a)(b-z)/(b-a)^2$. Taking the derivative with respect to $b$ and setting the result equal to zero, we find that $b^* := 2z - a$ for $z > 1$. Substituting this value of $b^*$ back into $\p(E'\ge z)$, we find that $\p(E'\ge z) = (1-a)/(2(z-a))$.  
We illustrate the density of $E'$ that maximizes $\p(E'\ge z)$, denoted by $g_1$, in Figure~\ref{fig:pf-sketch-decreasing}.
Since $z > a$, the latter probability is a decreasing function of $a$, hence its supremum occurs for $a = 0$ and we have $\p(E\ge z) \le 1/(2z)$. 

It remains to show $ \p(E\ge 1/\gamma ) = \gamma/2$  for some $E\in \mathcal E_{\rm D}$. Take $E$ following a mixture distribution of a uniform distribution on $[0,2/\gamma]$ with weight $\gamma$ and a point mass at $0$. It is easy to see that $\p(E\ge 1/\gamma) =\gamma/2$ and $\E^\p[E]=1$, showing the desired result. The last statement of part (i) is trivial by checking with $E=1$.  
\end{proof}

The main message from Theorem~\ref{th:dec-uni} is that when an e-variable has a decreasing density, its threshold can be boosted by a factor of $2$. That is, $T_\alpha(\mathcal{E}_{\rm D}) = 1/(2\alpha)$. 
Further, $T_\alpha(\mathcal{E}_{\rm D>1}) = 1/(2\alpha) + \alpha/2 = T_\alpha(\mathcal{E}_{\rm D}) + \alpha/2$.
Therefore, for small values of $\alpha$, we can boost the threshold by a factor slightly less than 2 when $E \in \mathcal{E}_{\rm D>1}$. This validates our message that, for small $\alpha$, the distribution of uninformative e-values (between 0 and 1) has little impact on the rejection region; the shape of the distribution of large e-values has the most influence on the rejection threshold.

\begin{figure}[t]
\centering \includegraphics[width=0.95\textwidth]{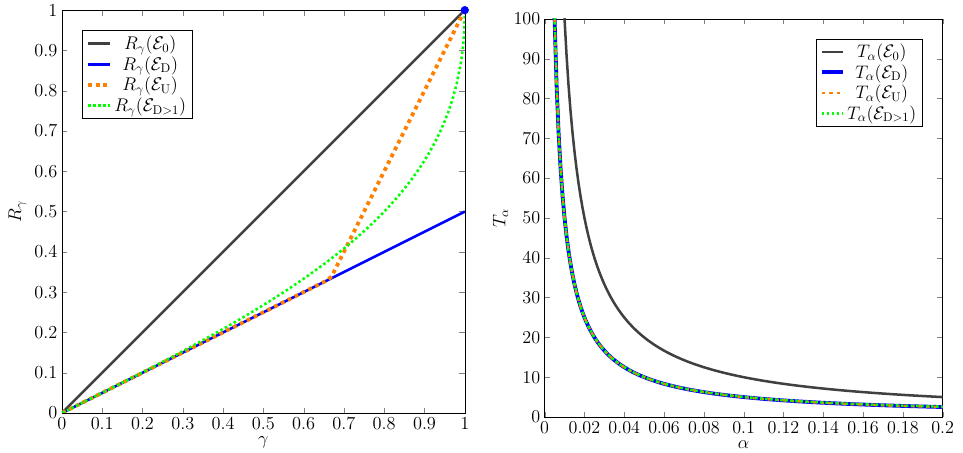}
\caption{Comparison of worst-case type-I errors and improved thresholds for decreasing and unimodal densities. In the left panel, the diagonal line is $R_{\gamma}(\fE)$, corresponding to Markov's inequality, and all other curves have a slope of approximately $1/2$ for small $\gamma$. }    \label{fig:threshold-comparison-1}
\end{figure}

In part (iii) of Theorem~\ref{th:dec-uni}, 
$R_\gamma (\mathcal E_{\rm U}   ) = \gamma /2 $,
for $\gamma \in (0,2/3]$   and $R_\gamma (\mathcal E_{\rm U}   ) = 2\gamma -1$ for $\gamma \in (2/3,1]$.
Since $\gamma\le 2/3$ is the most practical situation (meaning a type-I error control of $1/3$), the main message from part (iii) is similar to that from part (i): when an e-variable has a unimodal density,
its threshold can be boosted by a factor of $2$ in practice. 
We provide in Figure~\ref{fig:threshold-comparison-1} a comparison of different worst-case type-I errors and improved thresholds for $\mathcal{E}_{\rm 0}, \mathcal{E}_{\rm D},  \mathcal{E}_{\rm D>1}$ and $\mathcal{E}_{\rm U}$, where  we use $\cE_0$ for the full set $\fE$.

\section{Conditions on log-transformed e-variables}

In many applications of e-values, the final e-variable $E$ is the product of many e-variables, especially in the context of e-processes in Chapter~\ref{chap:eprocess}. In such a setting, assuming $\log E$ has some simple distributional properties may be reasonable, due to effects of central limit theorems. 
We say that a random variable $X$ has a symmetric distribution
if there exists $c\in \R$ such that $X$ and $c-X$ are identically distributed. 
We consider six different sets in this section:
\begin{align*}
\mathcal E_{\rm LS}& =\{E\in\fE: \mbox{$\log E$ has a symmetric distribution}\},
 \\
 \mathcal E_{\rm LU}& =\{E\in\fE: \mbox{$\log E$ has a unimodal density}\},
 \\ \mathcal E_{\rm LD>0}& =\{E\in\fE: \mbox{$\log E$ has a decreasing density over $[0, \infty)$}\},
\\ \mathcal E_{\rm LD}& =\{E\in\fE: \mbox{$\log E$ has a decreasing density}\},
\\
 \mathcal E_{\rm LUS}& =\{E\in\fE: \mbox{$\log E$ has a unimodal and symmetric distribution}\},
 \\
\mathcal E_{\rm LN}& =\{E\in\fE: \mbox{$E$  has a log-normal distribution}\}.
 \end{align*} 
 In all sets above, we require $\p(E=0)=0$, so that $\log E$ is a real-valued random variable. 
Note that $$\mathcal E_{\rm LN}\subseteq \mathcal E_{\rm LUS} \subseteq \mathcal E_{\rm LS} \mbox{~~~and~~~} \mathcal E_{\rm LD}\subseteq \mathcal E_{\rm LD>0}.$$
The point mass distributions $x\in (0,1]$ are included in all sets above, which are degenerate cases of log-normal distributions. 
Note that $E\in \fE$ has a log-normal distribution if and only if $\log E\lawis \mathrm{N}(\mu,\sigma^2)$ with $\mu+\sigma^2/2\le 0$.

\begin{theorem} 
\label{th:c12-log}
For $\gamma \in (0,1]$, 
\begin{enumerate}
\item [(i)]
     $R_\gamma (\mathcal E_{\rm LS}   ) =\gamma \wedge (1/2)$  if $\gamma\ne 1$ and $R_1(\mathcal E_{\rm LS}  ) =1$;
     \item[(ii)] $R_\gamma (\mathcal E_{\rm LU}   ) =\gamma$;
     
     \item[(iii)] $\mathcal{E}_{\rm LD>0} \subseteq \mathcal{E}_{\rm D>1}$ and 
$$ 
\frac{\gamma}{\e - \gamma} \le 
R_\gamma(\mathcal{E}_{\rm LD>0}) = 
\e^{-a_\gamma} \le  
\frac{\gamma }{1+ \sqrt{1-\gamma^2}} ,$$ 
where $a_\gamma$ is the unique solution of $\e^a(1 - a - \log \gamma) = 1$ for  $a \in (-\log \gamma, \infty)$;  
     
     \item[(iv)] 
     $R_\gamma (\mathcal E_{\rm LD}   ) =R_\gamma (\mathcal E_{\rm LUS}   ) $
 and     $$\frac{\gamma}e\le R_\gamma (\mathcal E_{\rm LD}   ) =R_\gamma (\mathcal E_{\rm LUS}   ) 
  \le \frac{\gamma}{\e(1-\gamma^2)} \wedge \left(   
\frac{\gamma }{1+ \sqrt{1-\gamma^2}} \right);$$
     
     \item[(v)]
     if $\gamma\ne1$ then $$R_\gamma (\mathcal E_{\rm LN}   ) =\Phi\left(-\sqrt{-2\log\gamma}\right) \leq  \frac{\gamma}{2\sqrt{-\pi \log \gamma}},$$
     where $\Phi$ is the standard normal cdf, and $R_1(\mathcal E_{\rm LN}   )=1$.
     \end{enumerate}
 \end{theorem}

 We omit the proof of Theorem~\ref{th:c12-log}.
To summarize the results, 
for log-transformed symmetric distributions 
and log-transformed unimodal distributions,
the standard Markov's bound cannot be improved for $\gamma\le 1/2$.
An improvement of roughly the order of $1/\e$ for small $\gamma$ is possible 
for log-transformed   decreasing densities on the positive real line, log-transformed  decreasing densities
and log-transformed  unimodal-symmetric densities.

\begin{figure}[t]
    \centering
  \includegraphics[width=0.95\textwidth]{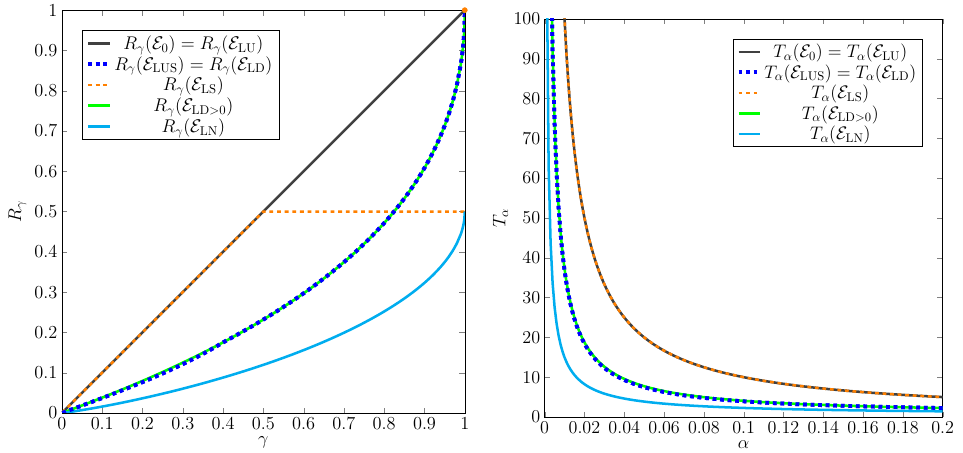}
    \caption{Comparison of worst-case type-I errors and improved thresholds for log-transformed e-variables.  
     In the left panel, the diagonal line is $R_{\gamma}(\fE)$, corresponding two the Markov inequality, and the curves for $\mathrm{LUS}$, $\mathrm{LD}$ and $\mathrm{LD}{>}0$ have a slope of approximately $\e^{-1}$ for small $\gamma$. The curve for $\mathrm{LN}$ is the smallest. 
    Results for $\mathrm{LUS}$ are conservative bounds since we only find an upper bound for $R_\gamma(\mathcal{E}_{\rm LUS})$ and $T_\alpha(\mathcal{E}_{\rm LUS})$.}
    \label{fig:threshold-comparison-3}
\end{figure}

We plot the bounds on worst-case type-I errors and improved thresholds from Theorem~\ref{th:c12-log} in Figure~\ref{fig:threshold-comparison-3}. Table~\ref{tab:improved-threshold} summarizes some numerical values for results in  Theorems~\ref{th:dec-uni} and~\ref{th:c12-log} (with $\cE_0=\fE$). 

\begin{table}[ht!]
	\centering
\centering\renewcommand{\arraystretch}{1.3}
	\begin{tabular}{rrrrrrr}
         &         \multicolumn{6}{c}{$\alpha$}         \\
         &  0.001 &  0.01 &  0.02 &  0.05 &  0.1 &  0.2 \\ \cmidrule(lr){2-7}
		       $\mathcal{E}_{0}, \mathcal{E}_{\rm LS}$  &   1000 &   100 &    50 &    20 &   10 &    5 \\
		    $\mathcal{E}_{\rm D}, \mathcal{E}_{\rm U}$ &    500 &    50 &    25 &    10 &    5 &  2.5 \\
		                       $\mathcal{E}_{\rm D>1}$  &    500 & 50.01 & 25.01 & 10.03 & 5.05 & 2.60 \\
		 $\mathcal{E}_{\rm LUS}, \mathcal{E}_{\rm LD}$    & 367.88 & 36.82 & 18.45 &  7.49 & 3.93 & 2.25 \\
		                      $\mathcal{E}_{\rm LD>0}$   & 368.25 & 37.16 & 18.77 &  7.73 & 4.07 & 2.25 \\
		                        $\mathcal{E}_{\rm LN}$ &    118 & 14.97 &  8.24 &  3.87 & 2.27 & 1.42
	\end{tabular}
    	\caption{Improved thresholds $T_{\alpha}(\mathcal{E})$ for different distributional hypothesis.}\label{tab:improved-threshold}
\end{table}

\section*{Bibliographical note}
The content of this chapter is largely based on \cite{blier2024improved}, where proofs of all results can be found, as well as some improved thresholds for combined e-values and the e-BH procedure. 
Computing probability bounds with distributional information is a well-studied topic in operations research and risk management; 
some general references on related techniques are \cite{shaked2007stochastic} and \cite{ruschendorf2024model}.
Comonotonicity in Section 
\ref{sec:c12-como} is a central concept in dependence modeling. It is also central to  decision theory under ambiguity via its connection to Choquet integrals \citep{schmeidler1989subjective}. A good reference on comonotonicity is \cite{denneberg1994non}.


\chapter{E-values and risk measures}
\label{chap:risk-measures}

This chapter studies the connection between e-values and risk measures. 
By risk measures, we mean general mappings from random variables or distributions to real numbers, which include the usual statistical functions such as the mean, the variance, and the quantiles. 

\section{Risk measures and statistical functions}
\label{sec:RM}

Risk measures are quantitative tools  used to quantify the riskiness level associated with financial positions, such as asset returns and losses, both at the individual level and at the portfolio level, in a fixed time period.  
Recall that $\cX$ is the set of all random variables on a given measurable space $(\Omega,\mathcal F)$,
and $\cM_1(\R)$ is the set of all distributions on $\R$.
There are two common formulations of risk measures (just like the expected value, which can be seen as either a mapping from   random variables or one from distributions, to the extended real line).

\begin{enumerate}
    \item[(a)] A risk measure can be formulated as  a mapping from a subset of $\cX$ to $(-\infty,\infty]$.
    We will use $\cR$ for such mappings.
    In this setting, 
    $\cR(X)$ represents the evaluation of the riskiness of a random variable $X$, representing a financial loss or gain. 
    \item[(b)]  A risk measure can be formulated  as a mapping from a subset of $\cM_1(\R)$ to $(-\infty,\infty]$.
    We will use $\rho, \phi$ for such mappings. 
        In this setting, 
    $\rho(F)$ represents the evaluation of the riskiness of the distribution $F$ of some financial loss or gain.      
\end{enumerate}

Mappings in both (a) and (b) will be called risk measures. 
When a reference probability measure $\p_0$ is fixed, the two formulations above can be connected via
$\cR(X)=\rho(F)$ where $F\in \cM_1(\R)$ and $X\lawis F$ (under $\p_0$). 
This correspondence is one-to-one  if $\cR$ satisfies \emph{law invariance}: $\cR(X)=\cR(Y)$ if $X\laweq Y$. 
In this case, a risk measure $\cR$ or $\rho$ is also called a \emph{statistical function}.

Certainly, common statistical functions, such as the mean, the variance, or the median, can also be naturally defined for distributions on $\R^d$ for $d\in \N$, and in this chapter we focus on the case of distributions on $\R$. We keep our discussions on risk measures minimal and omit the  generalizations such as set-valued, vector-valued, systemic, and dynamic risk measures, and many other variants, that are actively studied in the finance literature.

A main interpretation of a risk measure in finance is that its value $\cR(X)$ for $X\in \cX$ or $\rho(F)$ for $F\in \cM_1$ represents the regulatory capital requirement of the random loss $X$ or distribution $F$. When $X$ is treated as a random asset price  instead of loss, $\cR(X)$ can be interpreted as the price to purchase $X$ (the ``ask price'') in an incomplete financial market. The price $\cR$ would have been linear if the financial market were complete, and the nonlinearity is due to market frictions such as hedging limitations, bid-ask spreads, and transaction fees.

With the capital requirement interpretation, the most commonly used risk measures in the banking and insurance industries are the \emph{Value-at-Risk} (VaR) and the \emph{Expected Shortfall} (ES), defined as the following two classes of risk measures. Since they are law invariant, we will use the same notation for their version in both (a) and (b).  
We identify all distributions in $\cM_1(\R)$ with their cdf in this chapter. 
 
At level $\beta\in(0,1)$, the VaR is defined as the left $\beta$-quantile:
$$\VaR_\beta(F)=\inf\{x\in \R: F(x)\ge \beta\},~~F\in\cM_1(\R),$$
and the  ES  is defined as
$$
\ES_\beta(F)=\frac{1}{1-\beta}\int_\beta^1 \VaR_\gamma (F) \d \gamma,~~F\in\cM_1(\R).
$$ 
ES is also called 
 Conditional VaR (CVaR) and Tail VaR (TVaR) in different parts of the literature. 
Here, $\ES_\beta(F)$ may take the value $\infty$ if $F$ does not have a finite mean.
We also write $\VaR_\beta(X)$  and $\ES_\beta(X)$   for $X\in \cX$  to convert between versions in (a) and (b).
Note that $\VaR_\beta(X)=q_{1-\beta}(X)$ for the quantile function defined in Chapters~\ref{chap:markov},~\ref{chap:numeraire} and~\ref{chap:shape}. 
A well-known connection between ES and VaR is, for $X$ with finite mean,
\begin{align}
\label{eq:c13-VaR-formula}       \VaR _\beta(X) & \in   \argmin_{z\in \R} \left\{z + \frac{1}{1-\beta}\E^{\p_0} [(X-z)_+]\right\},  \\
  \ES_\beta(X) &=     \min_{z\in \R} \left\{z + \frac{1}{1-\beta}\E^{\p_0} [(X-z)_+]\right\} .
  \label{eq:c13-ES-formula} 
\end{align}
As of 2025,
ES at level $\beta=0.975$ is the standard measure for market risk in the current global banking regulation provided by the Basel Committee on Banking Supervision. 
VaR at various levels is the standard measure in insurance regulation and operational risk, and it is closely connected to the probability of default criterion used to measure credit risk.

A popular class of risk measures is that of \emph{coherent risk measures}.
Let $\cX_\cR\subseteq \cX$ be a convex cone in $\cX$ containing 
all constant random variables (identified with constants in $\R$).
A coherent risk measure $\cR$  is a mapping from $\cX_\cR $ to $(-\infty,\infty]$ satisfying four economic axioms: 
\begin{enumerate}[label=(\roman*)]
    \item Monotonicity: $\cR(X)\leq \cR(Y)$ for all $X,Y\in \cX_\cR$ with $X\leq Y$;
  \item Cash invariance: $\cR(X+c)=\cR(X)+c$ for all $X\in \cX_\cR$ and $c\in\mathbb{R}$;
    \item Subadditivity: $\cR(  X+ Y)\leq  \cR(X) + \cR(Y)$ for all $X,Y\in \cX_\cR$; 
    \item Positive homogeneity: $\cR(\lambda X)=\lambda \cR(X)$ for all $X\in \cX_\cR$ and $\lambda > 0$.    
\end{enumerate}
The four axioms have clear financial meanings. For instance, monotonicity means that a larger loss has bigger risk; subadditivity means that a portfolio is not riskier  than the individual risks summed due to possible diversification effects. 
We omit the detailed interpretations here.  Note also that subadditivity and positive homogeneity together imply convexity. 
For $\beta \in(0,1)$,
$\ES_\beta :\cX \to (-\infty,\infty]$ satisfies all of the above, and $\VaR_\beta  :\cX \to \R $ satisfies all but subadditivity. Hence, $\ES_\beta$ is coherent, but $\VaR_\beta$ is not.

Finally,  a risk measure $\cR:\cX_\cR\to (-\infty,\infty]$ is said to be \emph{continuous from above ($\downarrow$-continuous)} if 
$$
 \cR(X_n)\to\cR(X) \mbox{ for any bounded sequence 
 $(X_n)_{n\in \N} \downarrow X$,}
$$ 
where $(X_n)_{n\in \N} \downarrow X$ is understood point-wise.
In the next section, we will consider the set $\cX_+$ of all nonnegative random variables, and connect ($\downarrow$-continuous) coherent risk measures on $\cX_+$ to e-values. 
We summarize the domains of $\cR$ that we will encounter in Table~\ref{tab:c13-domains}.

\begin{table}[ht!]
    \centering
\centering\renewcommand{\arraystretch}{1.2}
    \begin{tabular}{ll}
      $\cX$   &  the set of all random variables \\
      
       $\cX_B$  & the set of  bounded random variables \\
       $\cX_{B+}$  & the set of  bounded nonnegative random variables \\
       $\cX_+$ & the set of nonnegative random variables  \\ 
       $\cX_\cR$  & a generic domain of  $\cR$ that is a convex cone containing $\R$
    \end{tabular}
    \caption{Domains of risk measures $\cR$.} 
    \label{tab:c13-domains}
\end{table}

\begin{remark}
    \label{rem:pf-2times-mean}
    The convexity of $\ES_\beta$   on $\cX$ yields a simple alternative proof of  the result that twice an average p-variable is a p-variable, presented in 
    Proposition~\ref{prop:2avg}.
  Let $P=(P_1+\dots+P_K)/K$, where $P_1,\dots,P_K$ are (arbitrarily dependent) p-variables. Without loss of generality we can assume  that each of $P_1,\dots,P_K$ is uniformly distributed on $[0,1]$.
  It suffices to show  
 $
\VaR_{\alpha}(2P) \ge {\alpha} $ for all $\alpha \in(0,1)$ (see Lemma~\ref{lem:quantile1} for conversion between probability and quantile conditions), which is equivalent to
 $
\VaR_{1-\alpha}(-2P) \le - {\alpha}  $ for all $\alpha \in(0,1)$. This follows from
\begin{align*}
\VaR_{1-\alpha} (- 2 P)     & \le \ES_{1-\alpha}(-2P )   
  \le \frac {2}{K} \sum_{k=1}^K \ES_{1-\alpha}(-P_k)  = - \alpha,
\end{align*}  
where we used the convexity of $\ES_{1-\alpha}$.
\end{remark}

\section{Connecting e-values to coherent risk measures}

\label{sec:super-ex}

Let us first notice that, for any  null hypothesis $\cP$, the requirement of a random variable $E\ge 0$ to be an e-variable is simply a risk measure constraint
 $\cR(E)\le 1$.
 This risk measure $\cR$ is defined as
\begin{align}\label{eq:c13-super-ex}
\cR(X)=\sup_{\p\in\cP} \E^\p[X]
\end{align}
for $X$ in a set such that the above expectations are well-defined (such as the set of nonnegative random variables).

The mapping $\cR$ in \eqref{eq:c13-super-ex}, sometimes called a super-expectation, has a long history as an important objects in different fields (and with different names), such as robust statistics, decision theory, optimization, finance, imprecise probability, and game-theoretic probability.  
As we see next, $\cR$ is a coherent risk measure, and it has essentially the only form of  coherent risk measures. 
\begin{theorem}
\label{th:coh-RM}
Let $\cX_{B+}$ be the set of all bounded nonnegative random variables. 
A mapping $\cR:\cX_{B+}\to\R$  
is a $\downarrow$-continuous coherent risk measure  if and only if 
\begin{align} 
\label{eq:c13-coherent}
\cR(X)=\sup_{\p\in \cP}\E^\p[X] ,~~~~X\in \cX_{B+}
\end{align} 
for some set $\cP\subseteq \cM_1$.
\end{theorem}
\begin{proof}[Proof.]
It is straightforward to verify the ``if'' statement.
For the ``only if'' statement,
 let $\cX_{B}$ be the set of all bounded  random variables. By Theorem 4.22 of \cite{Follmer/Schied:2011}, any $\cR':\cX_B\to \R$ 
    is a $\downarrow$-continuous coherent risk measure if and only if it has a representation in 
    \eqref{eq:c13-coherent} on $\cX_B$. 
    For $\cR $ on $\cX_{B+}$, it suffices to extend it to $\cX_B$ by letting $\cR(X)=\cR(X-t)+t$ for $X\in \cX_B\setminus \cX_{B+}$ where $t$ is the infimum value of $X$. One can easily check  that with this extension, $\cR$ is a coherent risk measure on $\cX_B$, and hence the aforementioned result guarantees \eqref{eq:c13-coherent} on $\cX_{B+}$.
\end{proof}
The set $\cX_{B+}$ of bounded and nonnegative random variables in Theorem~\ref{th:coh-RM} can be replaced by larger sets, but usually characterization results such as Theorem~\ref{th:coh-RM} are formulated with small domains for the greatest strength. 
A useful remark is that if $\cR: \cX_B\to \R$ 
is law invariant, then the representation \eqref{eq:c13-coherent}  holds without requiring $\downarrow$-continuity.  
Moreover, on any convex cone of random variables, \eqref{eq:c13-coherent}  always defines a coherent risk measure, possibly taking the value $\infty$. This in particular holds on $\cX_+$, leading to the next observation. 

\begin{corollary}
    For any $\cP$, 
    there exists a coherent risk measure $\cR:\cX_+\to[0,\infty]$ such that 
    all e-variables for $\cP$ are precisely those $E
    \ge 0 $ that satisfy $\cR(E)\le 1$.
\end{corollary}

For any $\cP\subseteq \cM_1$, the log-optimal e-variable for $\cQ $
can be reformulated as the solution to the following problem  
\begin{align}
\label{eq:c13-max}
\begin{aligned}
\mbox{ to maximize~~~} & \E^\q[\log X ]
\\
\mbox{ subject to~~~} & \cR(X)\le 1;~ X\ge 0,
\end{aligned}
\end{align} 
where $\cR$ is a coherent risk measure with representation  \eqref{eq:c13-super-ex} on $\cX_+$.
The problem \eqref{eq:c13-max} has a clear financial interpretation: 
We aim to maximize our expected utility of the random wealth $X$ under $\q$ with a log utility function, and the constraint is that our wealth $X$ cannot be negative, and its price $\cR(X)$ is bounded by our initial budget $1$. Recall that when $X$ represents a random financial wealth,  the risk measure value $\cR(X)$ can be interpreted as the price of $X$ in an incomplete  financial market. 
 This interpretation is consistent with testing by betting treated in Chapter~\ref{chap:eprocess}.

We have seen one example of \eqref{eq:c13-max} in Section~\ref{sec:c5-ex}. 
The dual representation of $\ES_\beta $ in \eqref{eq:c13-coherent} is  given by 
\begin{equation}\label{eq:ES-rep}
\ES_\beta(X) = \sup\left\{\E^\p[X] : \p\in \cM_1,~ \frac{\d \p}{\d\p_0}\le \frac{1}{1-\beta} \right\},~~~X\in \cX.
\end{equation}
In other words, the constraint of $\ES_\beta$ formulates all e-variables for the null hypothesis $\mathcal P$ that the true data-generating distribution is close to $\p_0$ in terms of likelihood ratio being bounded by $1/(1-\beta)$.
Similarly, many other families of coherent risk measures correspond to other testing problems.

Finally, when the alternative hypothesis $\cQ$ is composite,  one has  a few options to pick. The first is to consider  the problem of maximizing the mixture e-power for some probability measure $\nu$   over $\cQ$:
 \begin{align}
\label{eq:c13-mix-max}
\begin{aligned}
\mbox{ to maximize~~~} & \int_\cQ \E^\q[\log X ] \d \nu 
\\
\mbox{ subject to~~~} & \cR(X)\le 1;~ X\ge 0,
\end{aligned}
\end{align} 
where  and $\cR$ is again the coherent risk measure in  \eqref{eq:c13-super-ex}.
The second is to  maximize the worst-case e-power
 \begin{align}
\label{eq:c13-worst-max}
\begin{aligned}
\mbox{ to maximize~~~} & \inf_{\q \in \cQ} \E^\q[\log X ] 
\\
\mbox{ subject to~~~} & \cR(X)\le 1;~ X\ge 0.
\end{aligned}
\end{align}  
The third is to maximize the relative e-power compared with an oracle e-variable,
 \begin{align}
\label{eq:c13-relative-max}
\begin{aligned}
\mbox{ to maximize~~~} & \inf_{\q \in \cQ} \left(\E^\q[\log X ]-\E^\q[\log E^*_\q ] \right)
\\
\mbox{ subject to~~~} & \cR(X)\le 1;~ X\ge 0,
\end{aligned}
\end{align}  
where $E^*_\q$ is the numeraire for $\cP$ against $\q$ (from Chapter~\ref{chap:numeraire}).

When there are multiple observations available for the testing problem, 
the solution to the third problem \eqref{eq:c13-relative-max} 
is able to learn the true data-generating distribution with data 
as discussed in Section~\ref{sec:3-mix-plugin}, whereas the solutions to the  first two problems do not.


We end this section by noting that the problems \eqref{eq:c13-mix-max}, \eqref{eq:c13-worst-max},  
and \eqref{eq:c13-relative-max} all have   concave objectives to maximize, subject to (possibly infinitely many) linear constraints, and therefore it is a convex program and can be solved efficiently in many cases, in particular when the state space, $\cQ$ and $\cP$ are all finite.

\section{E-values for testing risk measures}

\label{sec:e-statistics}

In this section, we discuss how to design e-values to test the value of a given law-invariant risk measure.
We will mainly use the formulation $\rho$ that maps a subset of $\cM_1(\R)$ to $\R$.

 \subsection{E-statistics}

The general goal is to test the value of $\rho(X_t)$ through a   sequence  $(X_t)_{t\in \T}$ of data points, which are not assumed to be iid. Here, the index set $\T$ may be finite or infinite. 
For instance, the null hypothesis may be
\begin{equation}
    \label{eq:c13-H0-s}
H_0: \mbox{$X_t $ has distribution $F_t$ with $\rho(F_t) \le r_t$ for $t\in \T$.} 
\end{equation} 
Here, the values $r_t$ can be interpreted as the risk forecast for $X_t$ evaluated with $\rho$.

Our general testing procedure is to follow the idea of testing by betting described in Chapter~\ref{chap:eprocess}, in the following two steps.   
\begin{enumerate}
    \item Compute an e-variable $E_t$ from each data point $X_t$, such that $E_1,E_2,\dots $ are sequential e-variables. 
    \item Build an e-process via testing by betting,  such as the empirically adaptive e-process in Section~\ref{sec:EAEP}. 
\end{enumerate}

To follow the above procedure, we need to build an e-variable from each data point. Below, 
let $X$ be a random variable representing one data point. We will first study how we can build an e-variable from this data point.   
Let $\cM_*\subseteq \cM_\rho$ represent the set of distributions of interest (it may be smaller than the domain of $\rho$). 
The simplest problems are to test 
\begin{equation}
    \label{eq:c13-H01}
H_0: \mbox{$X $ follows any $F\in \cM_*$ with $\rho(F) = r$},
\end{equation}
and to test
\begin{equation}
    \label{eq:c13-H012}
H_0: \mbox{$X $ follows any $F\in \cM_*$ with $\rho(F) \le  r$.}
\end{equation}
In finance, we are mainly interested in the one-sided hypothesis in \eqref{eq:c13-H012} and  its sequential version \eqref{eq:c13-H0-s}, as they have a clear interpretation: We are testing whether the value $r$ is sufficiently conservative for the risk measure $\rho$ of financial risks. In banking, $r$ usually represents regulatory capital charge for a risky investment, and $r\ge \rho(F)$ means that the capital reserve is safe and passing the regulatory requirement.
Two-sided tests for \eqref{eq:c13-H01} can be easily constructed by combining one e-variable for testing $\rho(F)\le r$ and one e-variable for testing $\rho(F)\ge r$ via taking an average. 

We also consider the situation in which we have another  statistical function  $\phi$, for which we have information such as $\phi(F)=z$, but we are not interested in testing it. We assume the domain of $\phi$ contains $\cM_\rho$. An example of $(\rho,\phi)$ is $(\var,\mean)$, where we are only interested in testing variance, but the information on the mean is available. The corresponding hypotheses are 
\begin{equation}
    \label{eq:c13-H01-z} 
H_0: \mbox{$X $ follows any $F\in \cM_*$ with $\rho(F) = r$ and $\phi(F)=z$}, 
\end{equation} 
and
\begin{equation}
    \label{eq:c13-H01-z2}
H_0: \mbox{$X $ follows any $F\in \cM_*$ with $\rho(F) \le r$ and $\phi(F)=z$.}
\end{equation} 
We can easily allow $\phi$ to take values in $\R^d$ instead of $\R$, and we omit such  a generalization.

To test the above hypotheses, we  are interested in finding a
 function $e:(x,r,z)\mapsto \R$ such that 
 $e(X,r,z)$ is an e-variable for
the hypotheses in \eqref{eq:c13-H01}--\eqref{eq:c13-H01-z2}.
 Below, $\rho(\cM_*)$ is the set of values taken by $\rho(F)$ for $F\in \cM_*$. 

\begin{definition}[E-statistics for $\rho$]\label{def:e-stat-1}
Fix $\cM_*\subseteq\cM_\rho$ and $\rho:\cM_\rho\to \R  $, and let $e:\R^2\to [0,\infty]$ be a measurable function. 
\begin{enumerate}[label=(\roman*)]
\item The function $e $ 
is a  \emph{point e-statistic for $\rho$} if $e(X,r)$ is an e-variable for $H_0$ in \eqref{eq:c13-H01},
and it  is a \emph{one-sided e-statistic for $\rho$}
 if $e(X,r)$ is an e-variable for $H_0$ in \eqref{eq:c13-H012}.
\item    A one-sided e-statistic $e $ for $\rho$  is a \emph{backtest e-statistic for $\rho$}   if  
$\int_\R e(x, r ) \d F(x)> 1$ for all $r \in \rho(\cM_*)$ and $F\in \cM_*$ with $\rho(F)> r$.  

 \item  A backtest e-statistic  $e$ is   \emph{monotone} if  $r\mapsto e(x,r)$ is decreasing for each $x$. 
\end{enumerate}
\end{definition}
Note that $\cM_*$ is implicit in the above definition.  
Among the above concepts, it is clear that 
$$ 
\mbox{backtest e-statistic} ~\Longrightarrow~
\mbox{one-sided e-statistic}
~\Longrightarrow~
\mbox{point e-statistic}.
$$
 For example, if $\psi$  
is the mean and $\cM_*$ is the set of distributions on $(0,\infty)$,  then $e:(x,r)\mapsto x/r$   is a  monotone backtest e-statistic, satisfying all requirements in Definition~\ref{def:e-stat-1} (see Example~\ref{ex:e-stat-1} below). 
The property of having a mean larger than $1$ under the alternative is crucial for the e-statistic.
By Theorem~\ref{thm:nontrivial}, if $e$ is a backtest e-statistic and $ \rho(F)>r$, then  there exists $\lambda\in (0,1)$ such that $(1-\lambda) +\lambda e(X,r)$ has positive e-power. In the context of an infinite sequence of observations, this condition is sufficient for establishing consistency of the empirically adaptive e-process by Theorem~\ref{th:adaptive}.

We   analogously define e-statistics for $(\rho,\phi)$, which has one more variable.
\begin{definition}[E-statistics for $(\rho,\phi)$]\label{def:e-stat-2}
Fix $\cM_*\subseteq\cM_\rho$ and $(\rho,\psi):\cM_\rho\to \R^2 $,  and let $e:\R^3\to [0,\infty]$ be a measurable function.
\begin{enumerate}[label=(\roman*)] 
 \item The function $e $ 
is a  \emph{point e-statistic for $(\rho,\psi)$} if $e(X,r)$ is an e-variable for $H_0$ in \eqref{eq:c13-H01-z},  and it  is a \emph{one-sided e-statistic for $(\rho,\psi)$}
 if $e(X,r)$ is an e-variable for $H_0$ in \eqref{eq:c13-H01-z2}.  
\item  A one-sided e-statistic $e$ for  $(\rho,\phi)$ is a \emph{backtest e-statistic for $(\rho,\phi)$}   if  
$\int_\R e(x, r,z) \d F(x)> 1$ for all $(r,z) \in (\rho,\phi)(\cM_*)$ and $F\in \cM_*$ with $\rho(F)> r$.    

 \item  A backtest e-statistic  $e$ is  \emph{monotone} if $r\mapsto e(x,r,z)$ is decreasing for each $(x,z)$. 
\end{enumerate}
\end{definition}

The crucial interpretation of a backtest e-statistic for $(\rho,\phi)$ is that, if the risk measure  $\rho$ is underestimated, 
then a backtest e-statistic has power against the null hypothesis, regardless of whether the prediction of the auxiliary risk measure $\phi$ is truthful.

To explain the motivation behind monotone backtest e-statistics, we consider a financial context. In banking regulation, the risk measure $\rho$ is used to compute regulatory capital, and $\phi$ has no direct financial consequence. Therefore, a firm may have incentive to forecast $\phi$ arbitrarily (correctly or incorrectly), but forecasting a large $\rho$ would result in financial cost.
This motivation explains the term ``backtest'', which is a testing procedure in finance for risk models.
If a backtest e-statistic is monotone, then an overestimation of the risk is rewarded: A  firm being scrutinized by the regulator can deliberately report a higher risk value (which   means a higher capital reserve) to pass to the regulatory test, thus rewarding prudence.   
  We omit a detailed discussion on the financial interpretation here.
  In a non-financial   context, a statistician may be only interested in point or one-sided e-statistics.

As mentioned above, after choosing a suitable e-statistic $e$, we can test $H_0$ in \eqref{eq:c13-H0-s}, as well as other similar hypotheses, by constructing the 
e-variables  $E_t= e(X_t,r_t)$ for $t\in \T$,
 where $r_t$ is the forecast for the risk of $X_t$ conditional on $\mathcal F_{t-1}$, which is the typical situation in financial risk management (e.g., firms provide forecast for  the next-day risk of their portfolio).
 These e-variables are sequential, because $\E^\p[e(X_t,r_t)|\mathcal F_{t-1}] \le 1$ under any distribution $\p$ in the null hypothesis.   
 Therefore, 
even if $X_1,X_2,\dots$ are not iid and possibly dependent, 
we can   build an e-process based on $(E_t)_{t\in \T}$  as in Section~\ref{sec:EAEP} to test the null hypothesis, that is, 
\begin{align}
  M_t =\prod_{s=1}^t \left( (1-\lambda_s) + \lambda_s E_s\right) \mbox{~for $t \in \T$ and }M_0=1,
  \label{eq:c16-e-martingale}  
\end{align}  
The predictable process $(\lambda_t)_{t\in \T}$ in   \eqref{eq:c16-e-martingale} can be chosen in various ways. For instance, it may be computed through  the empirically adaptive e-process (Definition~\ref{def:empirical}), or computed based  on the empirical distribution of the sample $(X_t)_{t\in \T}$  without using the forecasts $(r_t)_{t\in \T}$. 
 In what follows, we will focus on e-statistics for some commonly used statistical functions.

\subsection{Mean, variance,   quantiles, and expected losses}
\label{sec:c13-mean}
 
 We derive backtest e-statistics for some common  statistical functions. In this chapter, we use the convention 
that $0/0=1$ and $1/0=\infty$. 
Let $\cM^k\subseteq \cM_1(\R)$ be the set of all distributions with finite $k$-th-moment for $k>0$. 
All expectations below only concern the distribution of $X$, and we omit the probability $\p$ in $\E[\cdot ]$.

\begin{example}[Backtest e-statistic for the mean]
\label{ex:e-stat-1}
Let $\cM_*$ be the set of distributions on $\R_+$ in $\cM^1$.
Define the function
$e(x,r) = x/r$ for $x,r\ge 0$.
In this case, we have $\E[e(X,r)] \le 1$ for all random variables $X$ with distribution in $\cM_*$ and $ r \ge \E[X]$. Moreover, for any such $X$, $\E[X]>r\ge 0$ implies
$\E[e(X,r)]>1$. 
Therefore,  $e $ is a monotone $\cM_*$-backtest e-statistic for the mean. 
\end{example}

\begin{example}[Backtest e-statistic for $(\var,\mean)$]
\label{ex:e-stat-2}
Consider $(\var,\mean):\cM^2\to\R^2$, which is the pair of the variance and the mean.
The function
$e(x,r,z) = (x-z)^2/r$ for $x,z\in \R$ and $r\ge 0$ is a monotone backtest e-statistic  for $(\var,\mean)$. To see this, for all random variables $X$ with distribution in $\cM^2$, we have
$$\E[e(X, r, \E[X])] =\frac{\E[(X-\E[X])^2]}{r}\le 1$$
for $r \ge \var(X)$.
Moreover, since $z=\E[X]$ minimizes $\E[(X-z)^2]$ over $z\in \R$ and $\mathrm{var}(X)=\E[(X-\E[X])^2]$, 
$\var(X)>r \ge 0$ implies
 $$\E[e(X, r,z)] =\frac{\E[(X-z)^2]}{r} \ge \frac{\var(X)}{r}  >1.$$ 

\end{example}

\begin{example}[Backtest e-statistic for a quantile]
\label{ex:e-stat-5} 
Take $\beta \in (0,1)$. 
Define the function
\begin{equation}\label{eq:eqp}
e_{\beta }^q(x,r) = \frac{1}{1-\beta }\id_{\{x> r\}},~~x,r \in \R.
\end{equation}
We have $\E[e_\beta^q(X,r)] \le 1$ for any random variable $X$ with distribution $F$ and $r \ge \VaR_\beta (F)$. Moreover, $\VaR_\beta (F)>r$ implies $\p(X>r)>1-\beta$, and hence
$\E[e_\beta^q(X,r)]>1$. 
Therefore,  $e_\beta^q $ is a monotone backtest e-statistic for the $\beta$-quantile. 
\end{example}

\begin{example}[Backtest e-statistic for an expected loss]
\label{ex:3}
For some $a\in \R$, let $\ell :\R\to [a,\infty)$ be a function that is interpreted as a loss.
Define the function
$e(x,r) = (\ell(x)-a)/{(r-a)}$ for $x\in \R$ and $r\ge a$.
Analogously to Example~\ref{ex:e-stat-1},   $e $ is a monotone backtest e-statistic for the expected loss $F\mapsto \int \ell \d F$ on its natural domain.   
\end{example}

For fixed $\rho$ or $(\rho,\phi)$, the choice of a backtest e-statistic $e$  is not  unique. 
For instance, a linear combination of $e$ and $1$ with  weights straightly between $0$ and $1$ is also a backtest e-statistic for $\rho$ or $(\rho,\phi)$. 
Depending on the specific situation, either e-statistic may be useful in practice.

\subsection{E-statistics for testing ES}
\label{sec:c13-ES}

The functional $ (\rho,\phi) = (\var,\mean)$ in Example~\ref{ex:e-stat-2} is an example of a \emph{Bayes pair}; that is, there exists a measurable function $L:\R^{2}\to \R$  such that
\begin{align}\label{eq:bayes}
\begin{aligned}
 \phi(F) &\in \argmin_{z \in \R } \int L(z,x) \d F(x) \text{ and }\\ 
  \rho(F)&=\min_{z \in \R}  \int L(z,x) \d F(x)\quad \mbox{ for } F\in \cM_*,
  \end{aligned}
 \end{align}
where $\int L(z,x) \d F(x)$ is assumed to be well-defined for each $z\in \R$ and $F \in \cM_*$. 
The function $L$ is the square loss $(z-x)^2$ in the case of $ (\var,\mean)$.
Bayes pairs often admit backtest e-statistics. A typical example commonly used in risk management practice is $(\ES_\beta,\VaR_\beta):\cM^1\to \R^2$ treated below.
Note that we restrict the domain to $\cM^1$ so that $\ES_\beta$ takes real values.

For $\beta \in (0,1)$, define the function 
\begin{equation}\label{eq:ep}
e^{\mathrm{ES}}_\beta  (x,r,z) = \frac{(x-z)_+}{(1-\beta )(r-z)},~~~x\in \R,~ z \le r.
\end{equation}
This
defines a backtest e-statistic for $(\ES_\beta ,\VaR_\beta )$. Recall the convention that $0/0=1$ and $1/0=\infty$, and set $e^{\mathrm{ES}}_\beta (x,r,z)=\infty$ if $r<z$, which is a case of no relevance since $\ES_\beta (F)\ge \VaR_\beta (F)$ for any $F \in \cM^1$. 

\begin{theorem}\label{th:ep}
For $\beta \in (0,1)$, the function $e^{\mathrm{ES}}_\beta $ is a monotone backtest e-statistic for $( \ES_\beta ,\VaR_\beta )$.
\end{theorem}
\begin{proof}[Proof.]
By  \eqref{eq:c13-VaR-formula}  and \eqref{eq:c13-ES-formula},
the pair $(\ES_\beta,\VaR_\beta)$ is a Bayes pair in \eqref{eq:bayes} with     $L:(z,x)\mapsto z + (x-z)_+/(1-\beta)$. Since $L(z,x) \ge z$, for $r \ge z$, we have that $e^{\mathrm{ES}}_\beta(x,r,z)=(L(z,x)-z)/(r-z) \ge 0$, and it is decreasing in $r$. We have for $r \ge \ES_\beta(X)$ that 
 \[
 \E\left[\frac{L(\VaR_\beta(X),X)-\VaR_\beta(X)}{r-\VaR_\beta(X)}\right] = \frac{\ES_\beta(X)-\VaR_\beta(X)}{r - \VaR_\beta(X)} \le 1.
 \]
 Furthermore, for $z < r \le \ES_\beta(X)$, 
\[
\E\left[\frac{L(z,X)-z}{r-z}\right] \ge \frac{\ES_\beta(X)-z}{r-z} \ge 1
\] 
with equality if and only if $r = \ES_\beta(X)$. 
\end{proof}

As a consequence of Theorem~\ref{th:ep}, 
$e_\beta^\ES(X,r,z)$ is   an e-variable for  the null hypothesis 
$$H_0 :  \ES_\beta(X) \le r \mbox{~and~}\VaR_\beta(X) = z.$$
A generalization of this result is that  $e_\beta^\ES(X,r,z)$ is also an e-variable for  the larger null hypothesis 
$$H_0 :  \ES_\beta(X)-\VaR_\beta(X)\le r-z \mbox{~and~}\VaR_\beta(X) \le z.$$  
Nevertheless, we   note that 
$e_\beta^\ES(X,r,z)$ is not an e-variable for  
$H_0$: $\ES_\beta(X) \le r $ and $\VaR_\beta(X) \le z$.
   
While Examples~\ref{ex:e-stat-1}-\ref{ex:3}  and Theorem~\ref{th:ep} show that interesting backtest e-statistics exist, much more can be said about their general structure, in particular, they are shown to be essentially the only choices, up to some transforms.

\section*
{Bibliographical note}

 Coherent risk measures were introduced by \cite{artzner1999coherent}, and the representation in a form similar to Theorem~\ref{th:coh-RM} was obtained by \cite{delbaen2002coherent}.
 The representation under law invariance, without assuming continuity, was obtained by \cite{jouini2006law}, and generalized by \cite{delbaen2012monetary}. 
 Coherent risk measures were used as pricing formulas by \cite{wang2000class}.
There are many more general  classes of risk measures than the coherent risk measures; a good reference is  \cite{Follmer/Schied:2011}, which also discusses pricing in incomplete financial markets.
A comprehensive treatment of 
the use of risk measures in risk management and financial regulation is  \cite{mcneil2015quantitative}.
In the finance literature, many authors  formulate risk measures on random gains instead of random losses, and their $X$ is our $-X$ in that case (for instance, monotonicity would be formulated in an opposite direction).  
The formulas \eqref{eq:c13-VaR-formula} and \eqref{eq:c13-ES-formula} were obtained by  \citet[Theorem 10]{rockafellar2002conditional}, where ES is called a CVaR,   common in the optimization literature. 
For the dual representation of ES in \eqref{eq:ES-rep}, see e.g., \citet[Theorem 8.14]{mcneil2015quantitative}.
Seven pedagogical proofs of the subadditivity of ES are provided by \cite{embrechts2015seven}.

Super-expectations 
have a very long history. For instance, they were axiomatized by \citet[Chapter 10]{huber1981robust} in the context of robust statistics. 
In addition to statistics and finance, super-expectations are also fundamentally important objects in the areas of decision theory under ambiguity (\cite{schmeidler1989subjective,gilboa1989maxmin}), 
imprecise probabilities (\cite{walley1991statistical}), 
non-linear expectations (\cite{peng2010nonlinear}),
and game-theoretic probability (\cite{Shafer/Vovk:2001,Shafer/Vovk:2019}). 

The content of Section~\ref{sec:e-statistics} is largely based on \cite{wang2022backtesting}, where they showed that most of the backtest e-statistics 
in Sections~\ref{sec:c13-mean} and~\ref{sec:c13-ES} 
are essentially the only possible choices in a suitable sense.
They also studied backtesting methods for risk measures based on e-statistics with financial data. The use of e-values and e-processes for forecast comparison was first  studied by \cite{henzi2021valid} on probability forecasts.  
Sequential testing and estimation of quantiles were also studied in \cite{howard2022sequential}. 
Bayes pairs were introduced by \cite{embrechts2021bayes}. 
The backtest e-statistic for the pair of ES and VaR  in Theorem~\ref{th:ep} is closely connected to the joint elicitability of ES and VaR, studied by \cite{fissler2016higher}.
Building e-statistics for the mean and variance and testing them via e-processes was also studied  by \cite{fan2024testing}.

\backmatter
\appendix

\renewcommand*{\thechapter}{A} 
\setcounter{theorem}{0}

\chapter{Appendix}


\section{Atomless probability spaces}
\label{app:atomless}

In several of our definitions, such as those of a calibrator or a merging function,
we have a universal quantifier over probability spaces.
In this section, we justify our claim in the main text that, when studying concepts like calibrators, e-merging function and p-merging functions, we can safely work with only one atomless probability measure.

Formally, a probability measure $\p:\mathcal F\to [0,1]$  is \emph{atomless}
if it has no \emph{atoms}; an atom is a set $A\in\mathcal F$
such that $ \p (A)>0$ and $\p (B)\in\{0,\p (A)\}$ for any $B\in\mathcal F$ with $B\subseteq A$.
With the same idea, we also say that the probability space $(\Omega,\mathcal F,\p)$ is atomless.
We used an alternative definition (point (ii) below) of atomless probability measures in the main text, which is justified by the following well-known lemma. 

\begin{lemma}\label{lem:rich}
  The following three statements are equivalent for any probability space $(\Omega,\cF,\p)$:
  \begin{enumerate}
  \item[(i)]
    $(\Omega,\cF,\p)$ is atomless;
  \item[(ii)]
    there is a random variable on $(\Omega,\cF,\p)$ that is uniformly distributed on $[0,1]$;
  \item[(iii)]
    for any Polish space $S$ and any probability measure $R$ on $S$,
    there is a random element on $(\Omega,\cF,\p)$ with values in $S$
    that is distributed as $R$.
  \end{enumerate}
\end{lemma}

Typical examples of a Polish space in item (iii) that are useful for us  
are $\R^K$ and countable sets.

\begin{proof}[Proof.]
  The equivalence between (i) and (ii) is stated in \citet[Proposition~A.31]{Follmer/Schied:2011}.
  It remains to prove that (ii) implies (iii).
  According to Kuratowski's isomorphism theorem \citep[Theorem 15.6]{Kechris:1995},
  $S$ is Borel isomorphic to $\R$, $\N$, or a finite set
  (the last two equipped with the discrete topology).
  The only nontrivial case is where $S$ is Borel isomorphic to $\R$,
  in which case we can assume $S=\R$.
  It remains to apply \citet[Proposition~A.31]{Follmer/Schied:2011} again.
\end{proof}

If $(\Omega,\mathcal F)$ is a measurable space and $\mathcal P$ is a collection of probability measures on $(\Omega,\cF)$,
we refer to $(\Omega,\mathcal F,\mathcal P)$ as a \emph{statistical model}. 
We say that it is \emph{rich}
if there exists a random variable on $(\Omega,\mathcal F)$ that is uniformly distributed on $[0,1]$
under any $\p\in\mathcal P$.

\begin{remark}\label{rem:U}
  Intuitively, any statistical model $(\Omega,\mathcal F,\mathcal P)$ can be made rich
  by complementing it with a random number generator
  producing a uniform random value in $[0,1]$:
  we replace $\Omega$ by $\Omega\times[0,1]$,
  $\mathcal F$ by $\mathcal F\times\mathcal B$,
  and each $\p\in\mathcal P$ by $\p\times \leb$,
  where $([0,1],\mathcal B,\leb)$ is the standard Borel space $([0,1],\mathcal B)$ equipped with the uniform probability measure $\leb$.
  If $\mathcal P=\{\p\}$ contains a single probability measure $\p$,
  being rich is equivalent to being atomless
  by Lemma~\ref{lem:rich}.
\end{remark}

For a statistical model $(\Omega,\mathcal F,\mathcal P)$,
an \emph{e-variable} is a random variable $E:\Omega\to[0,\infty]$ satisfying
\[
  \sup_{\p\in\mathcal P}
  \E^\p[E]
  \le
  1,
\]
as in Chapter~\ref{chap:introduction}.
The set of e-variables is denoted by $\fE_{\mathcal P}$.
The set $\fE_{\mathcal P}^K$ is the $K$-fold product set $\fE_{\mathcal P}\times 
\dots \times \fE_{\mathcal P}$. 

An \emph{e-merging function for $(\Omega,\mathcal F,\mathcal P)$}, 
is an increasing Borel function $F:[0,\infty)^K\to[0,\infty)$ such that,
for all $E_1,\dots,E_K$,
\begin{equation*}
  (E_1,\dots,E_K)\in\fE_{\mathcal P}^K
  ~\Longrightarrow~
  F(E_1,\dots,E_K) \in \fE_{\mathcal P}.
\end{equation*}
This definition requires that $K$ e-values for $(\Omega,\mathcal F,\mathcal P)$
be transformed into an e-value for $(\Omega,\mathcal F,\mathcal P)$.  

\begin{proposition}\label{prop:model-space}
  Let $F: [0,\infty)^K \to [0,\infty)$ be an increasing  Borel function.
  The following statements are equivalent:
  \begin{enumerate} 
  \item[(i)] $F$ is an e-merging function for some rich statistical model;
  \item[(ii)] $F$ is an e-merging function for all statistical models (i.e., an e-merging function in Definition~\ref{def:e-merg});
  \item[(iii)] $F$ is an e-merging function for all statistical models  $(\Omega,\mathcal F,\mathcal P)$ with a singleton $\mathcal P$.
  \end{enumerate} 
\end{proposition}

\begin{proof}[Proof.]
  Let us first check that, for any two rich statistical models
  $(\Omega,\mathcal F,\mathcal P)$ and $(\Omega',\mathcal F',\mathcal P')$,
  we always have
  \begin{equation}\label{eq:equal}
  \begin{aligned}
   &    \sup
    \left\{
      \E^\p[F(\mathbf{E})]
     :
      \p\in\mathcal P, \enspace \mathbf{E}\in\fE_{\mathcal P}^K
    \right\}
   \\& =
    \sup
    \left\{
      \E^{\p'}[F(\mathbf{E}')]
  :
      \p'\in\mathcal P', \enspace \mathbf{E}'\in\fE_{\mathcal P'}^K
    \right\}.
  \end{aligned} 
  \end{equation}
  Suppose $
    \sup
     \{
      \E^\p[F(\mathbf{E})] : \p\in\mathcal P, \enspace \mathbf{E}\in\fE_{\mathcal P}^K
     \}
    >
    c $
  for some constant $c$.
  Then there exist $\mathbf{E}\in\fE^K_{\mathcal P}$ and $\p\in\mathcal P$ such that $\E^\p[F(\mathbf{E})] > c$.
  Take a random vector $\mathbf{E}'=(E_1',\dots,E_K')$ on $(\Omega',\mathcal F')$
  such that $\mathbf{E}'$ is distributed under each $\p'\in\mathcal P'$ identically to the distribution of $\mathbf{E}$ under $\p$.
  This is possible as $\mathcal P'$ is rich
  by Lemma~\ref{lem:rich} applied to the probability space $([0,1],\mathcal B,\leb)$, where $\leb$ is the uniform (Lebesgue) measure and $\mathcal B $ is the Borel $\sigma$-algebra on $[0,1]$.
  By construction, $\mathbf{E}'\in\fE_{\mathcal P'}^K$ and $\E^{\p'}[F(\mathbf{E}')] > c$ for all $\p'\in\mathcal P'$.
  This shows that \eqref{eq:equal} holds with ``$=$'' replaced by ``$\le$'',  
  and we obtain the other direction by symmetry.

  The implications $\text{(ii)}\Rightarrow\text{(iii)}$ and $\text{(iii)}\Rightarrow\text{(i)}$ are obvious.
  To check $\text{(i)}\Rightarrow\text{(ii)}$, suppose  that $F$ is an e-merging function for some rich statistical model.
  Consider any statistical model.
  Its product with the uniform probability measure on $[0,1]$ will be a rich statistical model
  (see Remark~\ref{rem:U}).
  It follows from \eqref{eq:equal} that $F$ will be an e-merging function for the product.
  Therefore, it will be an e-merging function for the original statistical model.
\end{proof}

\begin{remark}
  The assumption of being rich is essential
  in item (i) of Proposition~\ref{prop:model-space}.
  For instance, if we take $\mathcal P := \{\delta_\omega\mid \omega\in\Omega\}$,
  where $\delta_{\omega}$ is the point mass at $\omega$,
  then $\fE_{\mathcal P}$ is the set of all random variables taking values in $[0,1]$.
  In this case, the maximum of e-variables is still an e-variable,
  but the maximum function is not a valid e-merging function
  as seen from Theorem~\ref{th:adm-merg}.
\end{remark}




Proposition~\ref{prop:model-space} shows that in the definition of an e-merging function
it suffices to consider a fixed atomless probability space $(\Omega,\mathcal F,\p)$.
Following the same idea, 
similar statements can be made for ie-merging functions and se-merging functions in Chapter~\ref{chap:multiple}.  

We can state a similar proposition  for calibrators.
A \emph{p-variable for a statistical model $(\Omega,\mathcal F,\mathcal P)$}
is a random variable $P:\Omega\to[0,\infty)$ satisfying
\[  \p(P\le\epsilon)
  \le
  \epsilon 
  \mbox{~for all }
  \epsilon\in(0,1) \mbox{~and~}
  \p\in\mathcal P.
\]
The set of p-variables for $(\Omega,\mathcal F,\mathcal P)$ is denoted by $\fU_{\mathcal P}$.
A decreasing function $f:[0,1]\to[0,\infty]$ is a \emph{(p-to-e) calibrator for $(\Omega,\mathcal F,\mathcal P)$}
if    $f(P)\in\fE_{\mathcal P}$ for any p-variable $P\in\fU_{\mathcal P}$.

\begin{proposition}\label{prop:model-calibrator}
  Let $f: [0,1] \to [0,\infty]$ be a decreasing  Borel function.
  The following statements are equivalent:
  \begin{enumerate}
  \item[(i)] $f$ is a calibrator for some rich statistical model;
  \item[(ii)] $f$ is a calibrator for all statistical models (i.e., a p-to-e calibrator in Definition~\ref{def:calib});
  \item[(iii)] $f$ is a calibrator for all statistical models  $(\Omega,\mathcal F,\mathcal P)$ with a singleton $\mathcal P$.
  \end{enumerate} 
\end{proposition}

The proof is similar to that of Proposition~\ref{prop:model-space}. 
There is an obvious analogue of Proposition~\ref{prop:model-calibrator}
for e-to-p calibrators  in Definition~\ref{def:calib} and combiners in Chapter~\ref{chap:combining}, which we omit. 
The content of this section is based on the appendices of \cite{Vovk/Wang:2021}.

The next lemma is useful for some proofs in Chapter~\ref{chap:theoretic}, and it is a restatement of  \citet[Lemma A.3]{wang2021axiomatic}.

\begin{lemma}\label{lem:atomless-p}
Suppose that $(\Omega,\cF,\p)$ is an atomless probability space. 
  For any random variable $X$ and $p\in (0,1)$,
  there exists an event $A\in \mathcal F$ with $\p(A)=p$  such that
$  \{  X >x\} \subseteq A\subseteq   \{ X \ge x\} $ for some $x\in \R$, which can be chosen as  $x=\inf \{t\in \R: \p(X\le t)\ge  1-p \}$.
\end{lemma}

\begin{proof}[Proof.]
Note that the chosen $x$ is the left   $(1-p)$-quantile of $X$, which satisfies  $\p(X> x) \le  p \le \p(X \ge  x) $; see Lemma~\ref{lem:quantile1}.
If either of the above inequalities is an equality, then we can take $A=\{X>x\}$ or $A=\{X\ge x\}$. 
 By Lemma~\ref{lem:rich},
 there exists a   random variable $U$ uniformly distributed on $[0,1]$.  
 For $s\in [0,1]$, 
 let $A_s= \{X>x\} \cup( \{X\ge x\}\cap \{U\le  s\})$.
 Clearly, $\p(A_0)=\p( X>x)\le p$
 and $\p(A_1)= \p(X\ge x)\ge p.$
 Moreover, for any $s,t\in [0,1]$ with $s\le t$, we have $\p(A_t)-\p(A_s)= \p(A_t\setminus A_s) \le \p(s\le U\le t) = t-s$.
 Hence,
 $s\mapsto \p(A_s)$ is continuous.
As a consequence, there exists $s
\in [0,1] $ such that  $\p(A_s)=p$. The condition $  \{  X >x\} \subseteq A_s\subseteq   \{ X \ge x\} $ holds by construction.
\end{proof}

\setcounter{chapter}{0} 
\section{Quantile functions}
\label{app:quantile}
In this section, we provide a few  fundamental facts about quantile functions, which are helpful to understand some results in the book using quantiles  in   Chapters~\ref{chap:markov},~\ref{chap:numeraire},~\ref{chap:shape} and~\ref{chap:risk-measures}.

Fix a probability measure $\p$. 
We define two versions of the quantile function under $\p$. 
For a random variable $X$,  its left quantile function 
is defined as
\begin{equation}
    \label{eq:app-left-q}
    Q_\alpha ^-(X)=\inf \{x\in \R: \p(X\le x)\ge  \alpha \}~~~~\mbox{for $ \alpha  \in (0,1]$}
\end{equation} 
and  its right quantile function  is defined as
\begin{equation}
    \label{eq:app-right-q}
        Q_\alpha ^+(X)=\inf \{x\in \R: \p(X\le x)> \alpha 
 \}~~~~\mbox{for $\alpha  \in [0,1)$}.
\end{equation}
Note that $Q_\alpha ^- =q_{1-\alpha }$ for $\alpha \in (0,1)$
with  $q_{1-\alpha }$   in Sections~\ref{sec:c2-LRB},~\ref{sec:c5-ex} and~\ref{sec:c12-1}, but we slightly generalize the definition here by  including the possibility of $\alpha =1$.

The first result is a list of translation rules between probability statements and quantile statements.

\begin{lemma}
    \label{lem:quantile1} 
For $\alpha \in [0,1]$, $t\in \R$ and  any random variable $X $, the following equivalence statements hold:
\begin{enumerate}[label=(\roman*)]
\item    
 $ Q_{\alpha}^{-}(X) > t$ $\iff$  $ \p (X \leq t)<\alpha $;
  \item $ Q_{\alpha}^{-}(X) \le  t$  $\iff$ $  \p(X \leq t)\ge \alpha $;  
        \item $ Q_{\alpha}^{-}(X)  \ge   t$ $\iff$   $  \p(X \leq s)< \alpha $ for all $s<t$; 
        
    \item  $ Q_{\alpha}^{-}(X)  <  t$  $\iff$ $  \p(X \leq s)\ge \alpha $ for some $s<t$; 
    
\item  $ Q_{\alpha}^{+}(X) < t$ $\iff$ $   \p(X < t) >\alpha $; 

 \item  $ Q_{\alpha}^{+}(X) \ge t$ $\iff$  $   \p(X < t)\le \alpha $; 

    \item   $ Q_{\alpha}^{+}(X)  \le   t$  $\iff$ $  \p(X < s)> \alpha $ for all $s>t$; 
        \item  $ Q_{\alpha}^{+}(X)   >   t$ $\iff$ $  \p(X < s)\le  \alpha $ for some $s>t$.
\end{enumerate} 
\end{lemma}
 \begin{proof}[Proof.]
To show (i), denote by
$ A_{\alpha} = \left\{ t\in \R: \p(X \leq t) \geq \alpha \right\}. $
Note that $A_\alpha$ is  closed   in $\R$ since $t\mapsto \p(X\le t)$ is upper semicontinuous. This gives $Q_{\alpha}^{-}(X) = \min A_{\alpha}.$
Hence,
$\alpha>\p(X\le t) \Longleftrightarrow t\not \in A_\alpha  \Longleftrightarrow  Q_{\alpha}^{-}(X) >t.$ 
Part
(ii) follows from (i) directly;
(iii) follows from (i) and (iv) follows from (ii). 

To show (v), denote by
$ B_{\alpha} = \left\{ t\in \R: \p(X < t) \leq \alpha \right\}$ which is  closed   in $\R$ since $t\mapsto \p(X< t)$ is lower semicontinuous. This gives $Q_{\alpha}^{+}(X) = \max B_{\alpha}.$
It follows that 
$\alpha<\p(X< t) \Longleftrightarrow t\not \in B_\alpha  \Longleftrightarrow  Q_{\alpha}^{+}(X) <t.$  Part (vi) follows from (v) directly; 
(vii) follows from (v) and (viii) follows from (vi). 
\end{proof}

The next lemma clarifies that the two quantile functions are almost everywhere equal. 
\begin{lemma}
    \label{lem:quantile1-a}
For any random variable $X$, $Q_\alpha^-(X)=Q_\alpha^+(X)$ for almost every $\alpha\in (0,1)$.
\end{lemma}
\begin{proof}[Proof.]
    Note that, for any $\beta <\alpha $, we have 
    $
    Q^-_\alpha(X) \ge Q^+_\beta(X) \ge Q^-_\beta(X).
    $
    Hence, by sending $\alpha \downarrow \beta$, we see that the two values $Q^+_\beta(X)$ and $Q^-_\beta(X)$ can only differ at discontinuous points of the curve $\alpha \mapsto Q^-_\alpha(X)$. Since $\alpha \mapsto Q^-_\alpha(X)$ is increasing, its jump points are countable, and thus the statement holds true. 
\end{proof}

The next lemma shows that if the quantile level $\alpha$ is replaced by a uniformly distributed random variable $U$ on $[0,1]$, then 
$Q_U^-(X)$  is identically distributed as $X$, where $Q_U^-(X)$ is understood as the random variable given by $\omega\mapsto Q_{U(\omega)}(X)$. This is true also for $Q_U^+(X)$.
Note that $Q_0^-(X)$ an $Q_1^+(X)$ are undefined, but this does not matter since $U$ takes the value $0$ or $1$ with zero probability.

\begin{lemma}
    \label{lem:quantile1-b}
For any random variable $X$ and any uniformly distributed random variable $U$ on $[0,1]$, 
the random variables $X $, $Q_U^-(X) $ and $Q_U^+(X)$ 
are identically distributed. 
\end{lemma}
\begin{proof}[Proof.]
By Lemma~\ref{lem:quantile1-a}, $Q_U^-(X) = Q_U^+(X)$  almost surely. Therefore, it suffices to consider $Q_U^-(X)$.
By Lemma~\ref{lem:quantile1}, we have, for any $\alpha\in (0,1)$ and $t\in \R$,
\begin{equation}\label{eq:app-quantile-1}
\begin{aligned}
& \p(X\le t)<\alpha \iff Q_\alpha^-(X)>t \\&~\Longrightarrow~ \p(Q_U^-(X)>t)\ge 1-\alpha 
 \iff  \p(Q_U^-(X)\le t)\le \alpha,
 \end{aligned}
\end{equation}
and 
\begin{equation}\label{eq:app-quantile-2}
\p(X\le t)\ge \alpha \iff Q_\alpha^-(X)\le t ~\Longrightarrow~ \p(Q_U^-(X)\le t)\ge \alpha. 
\end{equation}
Note that \eqref{eq:app-quantile-1} implies 
\begin{align*}
\p(X\le t)\le \alpha 
& \iff \p(X\le t)< \beta \mbox{~for all $\beta >\alpha$} \\ & ~\Longrightarrow~   \p(Q_U^-(X)\le t)\le \beta \mbox{~for all $\beta >\alpha$}  \\& \iff \p(Q_U^-(X)\le t)\le \alpha.
\end{align*}
Putting this and \eqref{eq:app-quantile-2} together we get $X\laweq Q_U^-(X).$
\end{proof}
Lemma~\ref{lem:quantile1-b} implies,  in particular, the well known formula
$$
\E^\p[X] =\int_0^1 Q^-_\alpha(X) \d \alpha = \int_0^1 Q^+_\alpha(X) \d \alpha.
$$

The next result is due to \citet[Theorem 1]{ryff1965orbits}, although presented in a different form.
 
\begin{lemma}
    \label{lem:quantile2}
    Suppose that $(\Omega,\mathcal F,\p)$ is an atomless probability space. 
For any random variable $X$, there exists a uniformly distributed random variable $U$ on $[0,1]$ such that $X=Q_U^-(X)=Q_U^+(X)$ almost surely. 
\end{lemma}

We first explain the simplest and most intuitive  case.  
Let $F$ be the cdf of $X$. 
If $X$ is continuously distributed, then it suffices to choose $U=F(X)$,
which satisfies the conditions in Lemma~\ref{lem:quantile2} by a standard probabilistic exercise. The complicated case is when $X$ is not continuously distributed and there does not exist a $\mathrm{U}[0,1]$ random variable independent of $X$.

 \begin{proof}[Proof.] 
By Lemma~\ref{lem:quantile1-a}, $Q_U^-(X) = Q_U^+(X)$  almost surely, so we will only show the first equality. 
First, suppose that there exists a standard uniform random variable $V$ independent of $X$.
 Define 
\begin{align}
\label{eq:app-quantile-3}
U= F(X)-V (F(X)-F(X-)),
\end{align}
where $F(x-)=\lim_{y\uparrow x}F(y)$ for any $x\in \R$.
It is a standard exercise in probability theory to check that $U$ is uniformly distributed on $[0,1]$ and 
$X=Q_U^-(X)$ almost surely. 

Next, we consider the more difficult case that there does not exist any uniformly distributed random variable independent of $X$.
In this case, 
let $B$ be the set of discontinuity points of $F$, which is a countable set. If $B$ is empty then $X$ is continuously distributed, and $U$ can be chosen as $F(X)$. Next, suppose that $B$ is not empty.
Denote by $A_b=\{X=b\}$ for $b\in B$. We have $\p(A_b)>0$, because $\p(A_b)=F(b)-F(b-).$ 
Since  the probability space $(\Omega,\mathcal F,\p)$  is atomless and $\p(A_b)>0$,  the space $ (A_b, \mathcal F|_{A_b}, \p|_{A_b})$ is also atomless. 
Hence, there exists a $\mathrm U[0,1]$-distributed random variable $V_b$ on this space. 
We extend $V_b$ to $(\Omega,\mathcal F,\p)$ by setting $V_b=0$ outside $A_b$. Finally,  define 
$$
U^*= F(X)-V_X (F(X)-F(X-)),
$$
where $V_X$ is understood as the random variable $V_X(\omega)= V_{b}(\omega)$ if $X(\omega)=b\in B$,
and otherwise  $V_X(\omega)=0$. 
Note that $V_X$ is $\mathcal F$-measurable since $B$ is countable. 
Moreover, $U^*$ has the same distribution as $U$ in \eqref{eq:app-quantile-3}, because conditional on $A_b$ for each $b\in B$, $V_b$ and $V$ are identically distributed,
and outside  $\bigcup_{b\in B} A_b$
the random variables $V$ and  $V_X$ do  not matter.
Following the same statement for $U$ in \eqref{eq:app-quantile-3}, we obtain that $U^*$ satisfies   the desired conditions in the lemma. 
\end{proof}


\bibliographystyle{plainnat}
\bibliography{bibliography}

\printindex

\end{document}